\documentclass[11 pt, twoside, openright]{amsbook}
\usepackage{forest}
\usepackage{hyperref}
\usepackage{amssymb}
\usepackage{microtype}
\usepackage{amsmath, enumerate}
\usepackage{mathabx}
\usepackage{amsthm}
\usepackage{color}
\usepackage{mathrsfs}
\usepackage{txfonts}
\usepackage{tikz}
\usepackage{adjustbox}
\usepackage{placeins}
\usetikzlibrary{positioning}
\usetikzlibrary{patterns}
\usetikzlibrary{decorations.pathreplacing}
\usetikzlibrary{positioning,arrows.meta, backgrounds}
\usetikzlibrary{calc}
\usetikzlibrary{angles,quotes}
\usepackage{graphicx} % Required for inserting images
\numberwithin{equation}{chapter}
\numberwithin{section}{chapter}

\newtheorem*{theorem*}{Theorem}
\usepackage{chngcntr}

\counterwithout{figure}{chapter}

\newtheorem{theorem}{Theorem}[chapter]

\newtheorem{definition}[theorem]{Definition}
\newtheorem{proposition}[theorem]{Proposition}
\newtheorem{corollary}[theorem]{Corollary}
\newtheorem{lemma}[theorem]{Lemma}

\newtheorem{letteredtheorem}{Theorem}

\usepackage{cleveref}

\title{Directional maximal operators in the plane}

\author{Edward Kroc} 
\address{University of British Columbia, Vancouver, Canada.}
\email{ed.kroc@ubc.ca}

\author{Juyoung Lee} 
\address{Korea Institute of Advanced Study (KIAS), Seoul, Korea.}
\email{juyounglee@kias.re.kr}  

\author{Malabika Pramanik}
\address{University of British Columbia, Vancouver, Canada.}
\email{malabika@math.ubc.ca}  
\date{\today}

	\hypersetup{hidelinks}
\begin{document}
	
	\UseRawInputEncoding
	\arraycolsep=1pt

	\maketitle 
	%\begin{abstract}
	
	\chapter*{Abstract} 
	%\addcontentsline{toc}{chapter}{Abstract} % Adds it to the Table of Contents	
		This monograph investigates the Lebesgue boundedness of planar directional maximal operators $D_{\Omega}$. These are maximal averages of functions over line segments in $\mathbb R^2$ whose slopes lie in a specified set $\Omega \subseteq \mathbb R$. A large body of multi-authored work has identified a geometric property of $\Omega$, called finite-order lacunarity, as a key factor in ensuring that $D_{\Omega}$ is Lebesgue bounded. While several variations of this notion exist in the literature, they all centre on the distribution of gaps in the slope set $\Omega$.
%		The precise formulation of this notion has a few variations in the literature, though they all remain closely aligned in spirit; the central theme in the study of gaps in the slope set $\Omega$.     
		\vskip0.1in
		\noindent Building on earlier work, an article of Bateman (2009) asserted a dichotomy for such operators. Namely,  $D_{\Omega}$ is bounded on 
		$L^p(\mathbb{R}^2)$ for all $p \in (1,\infty)$ precisely when the slope set $\Omega$ 
		is finite-order lacunary, or equivalently, when $\Omega$ does not 
		admit Kakeya-type sets. Conversely, sublacunary direction sets $\Omega$ admit Kakeya-like phenomena, implying in turn that  $D_{\Omega}$ is unbounded on $L^p(\mathbb R^2)$ for all $p \in [1, \infty)$.
		\vskip0.1in
		\noindent Recent work of Hagelstein, Radillo-Murguia, and 
		Stokolos (2024) identified a gap in the proof of this assertion and produced explicit examples for which the separation mechanism underlying that proof fails, demonstrating the need for a corrected framework. 
		%counterexamples, revealing that the previously accepted characterization cannot hold as stated. 
		%We provide a rigorous resolution under an adjusted notion of admissible finite-order lacunarity.   
%		we prove that $D_\Omega$ and $M_\Omega$ 
%		are bounded on $L^p$ if and only if $\Omega$ is admissible lacunary of 
%		finite order, equivalently if it does not admit Kakeya-type sets, thereby completing the characterization in the planar case. 
We establish the corrected characterization by introducing a new notion of admissible finite-order lacunarity that faithfully reflects the combinatorial structure of the underlying direction set. This leads to a tree-theoretic characterization in terms of finite splitting number and provides the foundation for new geometric and probabilistic constructions establishing the equivalence between admissible finite-order lacunarity, the absence of Kakeya-type sets, and the boundedness of directional maximal operators. The proof develops a pruning construction that yields a weaker
Euclidean separation condition and a corresponding sticky-map
framework sufficient to construct Kakeya-type sets for every
sublacunary slope set.

	\tableofcontents
	
{\allowdisplaybreaks	

\setcounter{chapter}{0}
	\chapter{Introduction} \label{chapter: intro}
	
	The study of maximal operators is a central theme in harmonic analysis, with deep connections to partial differential equations, geometric measure theory and incidence geometry. Among the fundamental examples are directional maximal operators - maximal averages of functions in Euclidean space along lines or rectangles in specified directions. This article investigates planar directional maximal operators and characterizes the geometric and combinatorial properties of the direction set that govern their Lebesgue-space boundedness.
	%	This article is devoted to the study of directional maximal operators in the plane. Specifically, we are interested in identifying  the geometric properties of a direction set that govern the Lebesgue mapping properties of the corresponding directional maximal operator. 
	 
	\vskip0.1in
	\noindent Let us introduce the operators of primary interest.  Given a set of slopes $\Omega \subseteq \mathbb R$,
	the \emph{planar directional maximal operator} $D_{\Omega}$ is defined by
	\begin{equation} \label{dir-max-op-def}
	D_\Omega f(x) := \sup_{\omega \in \Omega} \sup_{h>0}
	\frac{1}{2h} \int_{-h}^{h} | f(x + t(1,\omega)) | \, dt, \qquad x \in \mathbb R^2,
	\end{equation}
	for a function $f: \mathbb R^2 \rightarrow \mathbb C$ that is locally integrable along every line in $\mathbb{R}^2$. In other words, $D_{\Omega}f(x)$ is the supremum of the averages of $|f|$ along line segments centred at $x$ whose slopes lie in $\Omega$. 
	\vskip0.1in
	\noindent Closely related to $D_{\Omega}$ is the {\em{Nikodym-like maximal function}}
	 \begin{equation} \label{Nikodym-like-op-def}
	M_\Omega f(x) := \sup_{\omega \in \Omega} \sup_{R \in \mathscr{R}_\omega}
	\frac{1}{|R|} \int_R | f(x+y)| \, dy, \qquad x \in \mathbb R^2
	\end{equation} 
	where $\mathscr{R}_\omega$ denotes the family of rectangles in $\mathbb{R}^2$ that are 
	centred at the origin, with longer side parallel to the vector $(1,\omega)$. Here and throughout, $| \cdot |$ denotes Lebesgue measure. The ambient dimension $d \in \{1, 2\}$ will be clear from the context. 
	\vskip0.1in
	\noindent Our focus is the finiteness of operator norm bounds for $D_{\Omega}$ and $M_{\Omega}$ on Lebesgue spaces. It is well-known that the operators $D_{\Omega}$ and $M_{\Omega}$ have the same $L^p$ boundedness behaviour  for each $p \in (1, \infty)$. Indeed, there is an absolute constant $C_0 > 0$ for which the pointwise inequality 
	\begin{equation} \label{DM-pointwise-ineq}
		D_{\Omega} f(x) \leq M_{\Omega} f(x) \leq C_0 M_{\text{HL}} D_{\Omega} f(x) 
		\end{equation}
		holds for all smooth functions $f$ with compact support and all $x \in \mathbb R^2$. 
Here $M_{\text{HL}}$ is the classical two-dimensional Hardy-Littlewood maximal function, which is $L^p$ bounded for all $p \in (1, \infty)$. A proof of \eqref{DM-pointwise-ineq} is included in Appendix A, Section \ref{DM-ineq-proof-section}. Roughly speaking, the inequality \eqref{DM-pointwise-ineq} reflects the fact that a line segment average arises as a limit over thin rectangular averages. Conversely, rectangular averages are controlled by line averages through the Hardy-Littlewood maximal function. 
	 \vskip0.1in
	 \noindent A preliminary reduction of the slope set $\Omega$ is convenient. Although $\Omega$ may a priori be unbounded, we may express
	\[ \Omega = \{ \tan \theta : \theta \in \Theta \}, \quad \Theta \subseteq [0, \pi]\setminus\{ \frac{\pi}{2} \}   \] 
	as a union of four subsets $\Omega_n$ corresponding to four angular sectors: 
	\[ \Omega_n := \{ \tan \theta : \theta \in \Theta_n \}, \quad \Theta_n = \Theta \cap \Bigl[\frac{n \pi}{4} + \bigl[0, \frac{\pi}{4} \bigr] \Bigr], \quad n = 0, 1, 2, 3. \] 
	Since
\begin{equation} \label{Omega bounded justification} 
D_\Omega=\max_{0\leq n\leq 3}D_{\Omega_n} \leq \sum_{n=0}^{3}D_{\Omega_n}
\quad\text{and}\quad
M_\Omega=\max_{0\leq n\leq 3}M_{\Omega_n} \leq \sum_{n=0}^{3}M_{\Omega_n},
\end{equation} 
it suffices to consider one angular sector at a time.
	After a rotation of coordinates, which leaves the $L^p$ operator norms unchanged, each $\Omega_n$ can be reduced to a slope set contained in the unit interval:
	\begin{equation} \label{Omega bounded} 
	\Omega \subseteq [0,1].
	\end{equation}  
	This normalization allows us to reduce an unbounded slope set to a finite union of bounded ones without any substantive change to the operators $D_{\Omega}$ and $M_{\Omega}$.     
	\section{Formulation of the problem} 
	The guiding question of this paper is the following:
	\vskip0.1in
	\begin{center}
	\begin{minipage}{4in}
		{\em{Question: Given $p \in (1, \infty)$, for which slope sets $\Omega \subseteq [0,1]$ is $D_\Omega$, or equivalently $M_\Omega$, bounded on $L^p(\mathbb{R}^2)$?}}
	\end{minipage}
	\end{center}
	\vskip0.1in  
	\noindent Two extremal examples clarify the scope of the problem. 
	\vskip0.1in	
	\begin{itemize} 
	\item When $\Omega$ is a singleton, one-dimensional Hardy-Littlewood theory immediately implies that both $D_\Omega$ and $M_\Omega$ are bounded on $L^p(\mathbb{R}^2)$ for every $p\in (1, \infty)$. 
	\vskip0.1in 
	\item At the opposite extreme, if $\Omega = [0,1]$, neither operator can be bounded on $L^p(\mathbb{R}^2)$ for any $p\in (1, \infty)$. This follows from the existence of planar Kakeya-Besicovitch sets, which have Lebesgue measure zero yet contain line segments in every direction \cite[Chapter X, \S2.2]{SteinHA}. 
	\end{itemize} 
	\vskip0.1in
	The heart of the problem is therefore to quantify how ``thin'' an infinite 
	direction set must be to ensure Lebesgue boundedness of $D_{\Omega}$ and $M_{\Omega}$. The difficulty lies in determining which direction sets are sparse enough to prevent Kakeya-type phenomena -- existence of small sets with lines in many directions -- and which
possess enough combinatorial richness to support them . It is natural that the ideas and methods required to resolve this question lie at the intersection of harmonic analysis, geometric measure theory and combinatorics.

	\section{Finite-order lacunarity: an intuitive picture}
	
	A central concept in the study of directional maximal operators  is \emph{admissible finite-order lacunarity}. This is a quantitative measure of how gaps, or lacunas, are distributed within a set of slopes. As Theorem \ref{thm:main} will show, admissible finite-order lacunarity of $\Omega$ is not merely sufficient for Lebesgue boundedness of $D_{\Omega}$ and $M_{\Omega}$, but gives an exact characterization of it. 
	\vskip0.1in
	\noindent The formal, recursive definition of admissible finite-order lacunarity
	is given in Definition \ref{defn: Admissible finite order lacunarity} in  Chapter \ref{section: finite-order lacunarity}. Before turning to it, we describe the geometric intuition underlying the notion.
	%but let us offer the geometric intuition behind the definition here as motivation for the main result. 
	\vskip0.1in
	\noindent Roughly speaking, a sequence of slopes $\Omega \subseteq \mathbb R$ is called
	\emph{lacunary}, or lacunary of first order, if its elements cluster toward a limit point at a
	geometric (or faster) rate. For example, the sequence
	\begin{equation} \label{example-Omega} 
	\mathcal L = \left\{ 1, \tfrac{1}{2}, \tfrac{1}{4}, \tfrac{1}{8}, \ldots \right\}
	\end{equation} 
	is lacunary: each term is obtained by multiplying the previous one
	by a fixed constant $\lambda = \tfrac{1}{2}$. The slopes
	accumulate at $0$ in a highly controlled manner, with successive gaps shrinking proportionally to the distance from the limit. 
	\vskip0.1in
	\noindent Going a step further, one can envision a set consisting of a countable collection of lacunary sequences; for example
	\begin{equation}
		\mathcal L^{[2]} = \bigcup \mathcal L_j = \bigl\{2^{-j} + 2^{-k} : k \geq j \geq 1\bigr\}. 
		\end{equation} 
		Here the $j^{\text{th}}$ lacunary sequence $\mathcal L_j = \{2^{-j} + 2^{-k} : k \geq j\}$ converging to $2^{-j}$ lies in the interval $[2^{-j}, 2^{-j+1}]$, between two consecutive elements of a ``parent'' lacunary sequence $\mathcal L$, given by \eqref{example-Omega}. The set $\mathcal L^{[2]}$ is an example of a lacunary set of order 2. In other words, a lacunary set of order 2 is obtained by inserting lacunary sequences into the {\em{gaps}} of a first-order lacunary sequence. Iterating this construction finitely many times produces lacunary sets of arbitrary finite order.  For instance, 
		\[\mathcal L^{[N]} = \left\{ 2^{-j_1} + 2^{-j_2} + \ldots + 2^{-j_{N}} : 1 \leq j_1 \leq j_2 \leq \ldots \leq j_{N}\right\}  \] is
	lacunary of order $N$.
	\vskip0.1in
	\noindent In contrast, a set of slopes is
	\emph{sublacunary} if it cannot be captured by any finite iteration of lacunary clustering, nor by a finite union of sets obtained in this way. Sublacunary sets exhibit a form of combinatorial ``thickness''. They lack the rapid clustering of lacunary sets, instead accumulating slowly or along fractal patterns. Standard examples of a sublacunary set include intervals or dense subsets of intervals, such as $\mathbb Q \cap [0,1]$. Sequences that converge to their limit more slowly than geometrically, such as $\{\frac{1}{n} : n\geq 1\}$, or sets of positive Hausdorff dimension, like the standard Cantor middle-third set
	\begin{equation} \label{Cantor set}  \mathcal C:= \Bigl\{ \sum_{j=1}^{\infty} \varepsilon_j 3^{-j}: \varepsilon_j \in \{0, 2\} \text{ for all } j \geq 1 \Bigr\} \end{equation} 
	are also sublacunary. 
	\vskip0.1in
	\noindent On a heuristic level, therefore, admissible finite-order lacunarity measures structured sparsity. As Theorem \ref{thm:main} will show, {\em{sublacunary slope sets allow for Kakeya-type configurations}}: rectangles oriented along such directions can be arranged to overlap heavily, producing a union of small total area, yet separate almost completely when translated along their long axes by a fixed factor of their length \footnote{Definition \ref{Kakeya-type set} contains a precise description of what it means for a slope set to admit Kakeya-type sets.}. This geometric instability forces unboundedness of the directional maximal operators, and in fact characterizes sublacunarity. By comparison, admissible finite-order lacunary sets provide just enough directional separation to rule out Kakeya-type phenomena. Making this intuition rigorous is a primary objective of this article. The precise statement appears in Theorem \ref{thm:main}.
\vskip0.1in
\noindent A subtle but crucial issue arises when attempting to formalize this intuition. Finite recursive clustering alone does not suffice to rule out Kakeya-type configurations. Without uniform quantitative control on clustering rates, and without preventing unintended accumulation away from designated limit points, a recursively constructed set may still exhibit sublacunary behaviour. \footnote{See Section \ref{section: comparison with Bateman lacunarity} for relevant counterexamples.} Moreover, compatibility with known results and with the tree-based and
probabilistic arguments developed later requires finite-order lacunarity to be stable under finite unions and well-behaved transformations. Definition \ref{defn: Admissible finite order lacunarity}  is designed precisely to enforce these structural requirements.
%It is not enough to require that a slope set be generated by a finite recursive clustering procedure; one must also impose uniform quantitative control on the rate of clustering and prevent unintended accumulation away from designated limit points. The stability of the theory under taking subsets and finite unions also turns out to be essential. The formulation given later in Definition \ref{defn: Admissible finite order lacunarity} incorporates this hereditary requirement and ensures compatibility with the tree-based and probabilistic arguments that follow.	
	%there exist collections of rectangles oriented along sublacunary directions, which overlap heavily at one scale but separate rapidly under translation. More precisely, the  total area is small but if translated along their respective long axes by a fixed factor of their length, they spring apart and become essentially disjoint. This phenomenon, in turn, forces unboundedness of the directional maximal operators. This feature, in fact, {\em{characterizes}} sublacunarity. Finite-order lacunary sets provide just enough separation among directions to prevent Kakeya-like phenomena. 
	%whereas 
	%sublacunary slope sets retain enough density to reproduce substantial geometric overlap that drives unboundedness. 
	
	\section{History and the main result}
	
	The study of directional maximal operators has a long and rich history, spanning almost a century. Originating in classical problems related to Lebesgue differentiation, density bases and covering theorems \cite{{Besicovitch}, {Buseman-Feller}, {JMZ}, {Nikodym}}, the field now connects to diverse topics in  Fourier analysis \cite{{Cordoba}, {Cordoba-Fefferman1}, {Cordoba-Fefferman2}, {deGuzman1}, {deGuzman2}, {Fefferman}, {SteinHA}, {Grafakos}}, geometric measure theory \cite{{Katz1}, {Katz2}, {KatzTao}, {KatzTao2}, {KatzLabaTao}, {KatzZhal}, {Wolff}},  probability and stochastic processes \cite{{Grimmett}, {Lyons1}, {Lyons2}, {LyonsPeres}} and combinatorics \cite{{Gauvan-thesis}, {Gauvan1}, {HareRonning1}, {HareRonning2}}. Of direct relevance to this paper are the foundational results of C\'ordoba \cite{Cordoba}, Str\"omberg \cite{{Stromberg1}, {Stromberg2}}, and Nagel--Stein--Wainger \cite{NagelSteinWainger}, 
	later developed by Sj\"ogren--Sj\"olin \cite{SjogrenSjolin}, Alfonseca \cite{Alfonseca}, Katz \cite{{Katz1}, {Katz2}}, Bateman-Katz \cite{BatemanKatz}, Parcet--Rogers \cite{{ParcetRogers1}, {ParcetRogers2}}, and many others \cite{{AlfonsecaSoriaVargas1}, {AlfonsecaSoriaVargas2}, {Carbery}, {DanielloGauvanMoonens}, {DanielloGauvanMoonensRosenblatt}, {DanielloMoonens}, {DanielloMoonensRosenblatt}, {DuoandikoetxeaVargas}, {Hagelstein2013}, {HagelsteinParissis}, {HagelsteinStokolos2008}, {HagelsteinStokolos2009}, {Hare2000}, {KrocThesis}, {KrocPramanik}, {Moonens2016}, {Oniani}, {Stokolos1}, {Vargas}}. 
	\vskip0.1in
	\noindent To place our work in context, let us briefly recall some of the landmark results from the extensive literature on directional maximal operators in the plane. For slope sets $\Omega$ of infinite cardinality, the results most
relevant here point to the following boundedness dichotomy:
	\begin{enumerate}[A.] 
	\item {\em{Positive results, where $D_{\Omega}$ is shown to be $L^p$-bounded for all $p>1$}}. This happens, for instance, when $\Omega$ obeys one of these conditions:
	\begin{itemize}
	\item \cite{NagelSteinWainger} $\Omega = \{a_j: j \geq 1\}$ is a lacunary sequence, i.e. $0 < a_{j+1} \leq \lambda a_j$ for some $\lambda \in (0, 1)$ and all $j$.
	\item \cite{{SjogrenSjolin}, {Alfonseca}} The set of angles $\Theta$ corresponding to the slope set $\Omega = \{ \tan \theta : \theta \in \Theta \}$ is finite-order lacunary in the sense of the referenced articles.
	\end{itemize} 
	\vskip0.05in
	\item \label{negative item} {\em{Negative results, where $D_{\Omega}$ is unbounded on $L^p(\mathbb R^2)$ for all $p \in (1, \infty)$.}} Examples in this category include
	\begin{itemize} 
	\item \cite[Chapter X, \S2.2]{SteinHA} Any $\Omega$ that is dense in some interval.
	\item \cite{{Hare2000}, {BatemanKatz}} $\Omega = $ the standard middle-third Cantor set \eqref{Cantor set}.
	\end{itemize} 
	\end{enumerate} 
	These examples suggest that an appropriate formulation of finite-order lacunarity governs the boundary between boundedness and unboundedness of $D_{\Omega}$.
	\vskip0.1in
	\noindent A particularly influential development in the progression of the subject, building on this body of examples, was the work of Bateman
	\cite{Bateman}. It identified, for the first time, a connection between finite-order lacunarity of a set $\Omega$ and the branching structure of a binary tree encoding $\Omega$. 
	A key concept in this framework is the {\em{splitting number of the slope tree}}, which quantifies how often branching occurs along rays. Informally, finite splitting number corresponds to controlled clustering (finite-order lacunarity), while unbounded splitting reflects combinatorial richness sufficient to generate Kakeya-type configurations. This tree-theoretic perspective, originating in \cite{{BatemanKatz}, {Bateman}}, will play a central role in our argument.
	\vskip0.1in
	\noindent In \cite{{BatemanKatz}, {Bateman}}, this connection between clustering and trees was
used to construct Kakeya-type sets as unions of tubes oriented along a sublacunary set of slopes. We will revisit this connection in Chapters \ref{trees-section} and \ref{section: split implies lacunarity}.  The main result of \cite{Bateman} was a striking three-point characterization: namely,
	boundedness of $D_{\Omega}$ and $M_{\Omega}$ on $L^p(\mathbb{R}^2)$ 
	for some $p\in (1, \infty)$ is equivalent to the absence of Kakeya-type sets, which in turn is equivalent to the slope set $\Omega$ being finite-order lacunary in the sense of that article. 
	\vskip0.1in
	\noindent The characterization proposed in \cite{Bateman} was correct in spirit, but not as stated. Subsequent work revealed two distinct issues: the
definition of finite-order lacunarity did not capture the intended structural class, and the proof relied on an additional geometric separation property that was not guaranteed by that definition. Specifically: 
	\vskip0.1in 
\begin{itemize} 
	\item Kroc's thesis 
	\cite[Example (e), p 27]{KrocThesis} noted that the notion of finite-order lacunarity in \cite{Bateman} is insufficient for the main theorem. The recursive definition presented there does not guarantee the splitting structure of the slope tree required by the proof; see Section \ref{section: comparison with Bateman lacunarity}. 
	\vskip0.1in 
	\item Additionally, the geometric consideration at the heart of the argument was later shown to be incomplete. The recent work of Hagelstein, Radillo-Murguia, and Stokolos \cite{HRS2024a,HRS2024b} identified a gap in the Euclidean separation
argument underlying the Kakeya-type construction in \cite{Bateman}. They produced explicit counterexamples to the assumptions implicit in the proof strategy, demonstrating the need for a corrected framework. 
	\end{itemize} 
	%\noindent However, a few gaps were subsequently found:
	%\begin{itemize} 
	%	\item According to \cite{Bateman}, a set is lacunary set of order $N$ if it is covered by a lacunary set of order $N-1$ with lacunary sequences converging to every point of the latter. It has been pointed out in the Ph.D. thesis of Kroc \cite[Example (e), p 27]{KrocThesis} that this definition of finite-order lacunarity does not produce the main result of \cite{Bateman}. See also Section \ref{} below. 
	%	\item Further, the proof of the main result in \cite{Bateman} relied on a probabilistic construction that was later shown to be incomplete. Recent work of 
	%Hagelstein, Radillo-Murguia, and Stokolos 
	%\cite{{HRS2024a},{HRS2024b}} have identified the gap and produced explicit counterexamples, 
	%demonstrating that the argument as written in \cite{Bateman} cannot be carried through, as stated.  
	%\end{itemize}
	%Thus, the precise characterization of boundedness for planar 
	%directional maximal operators is considered unresolved at the moment, though a statement along similar lines is widely believed to be true.
	%\vskip0.1in

	\noindent In this paper we provide such a framework, giving a rigorous resolution of the planar problem.  
	\vskip0.1in
	\noindent Although the statement of our main result, Theorem \ref{thm:main}, formally resembles the three-fold characterization asserted in \cite{Bateman}, the equivalence established here rests on a corrected and coherent formulation of admissible finite-order lacunarity and a substantially different proof. The new ideas are discussed in Section \ref{section: main contributions}.
%	The distinction is not cosmetic; the admissible notion introduced in Definition \ref{defn: Admissible finite order lacunarity} is precisely aligned with the slope-tree structure and enforces key Euclidean separation properties required for a probabilistic construction to succeed. A detailed analysis of these differences and their impact on the arguments appears in Sections \ref{section: main contributions} and \ref{section: proof organization}.   
	%Our main theorem confirms the equivalence originally stated in 
	%\cite{Bateman}, but under a refined notion of admissible finite-order 
	%lacunarity and with a complete and self-contained proof. 
	
%	The aim of this paper is to provide a rigorous resolution of this question. 
%	Our main theorem confirms the equivalence originally stated in 
%	\cite{Bateman2009}, with an appropriate notion of admissible finite-order lacunarity and a complete and self-contained proof that addresses the issues raised in \cite{{HRS2024a}, {HRS2024b}}. 
	
 \begin{theorem}\label{thm:main}
		Let $\Omega \subseteq \mathbb R$ be a set of slopes, and let $p\in (1, \infty)$ be a Lebesgue exponent. 
		Then the following statements are equivalent:
		\begin{enumerate}[1.]
			\item The set $\Omega$ admits Kakeya-type sets, in the terminology of Definition \ref{Kakeya-type set}.  \label{condition 1: Kakeya-type sets}
			\item The directional maximal operators $D_\Omega$ and $M_\Omega$ are unbounded on $L^p(\mathbb{R}^2)$. \label{condition 2: operators unbdd}
			\item The set $\Omega$ is sublacunary, i.e. not admissible lacunary of any finite order in the sense of Definition \ref{defn: Admissible finite order lacunarity}.  \label{condition 3: slopes sublacunary}
			
		\end{enumerate}
	\end{theorem}
\noindent Equivalently, $D_\Omega$ and $M_\Omega$ are bounded on
$L^p(\mathbb R^2)$ for some $p\in(1,\infty)$ if and only if
$\Omega$ is admissible lacunary of finite order. In that case, both
operators are bounded on $L^p(\mathbb R^2)$ for every
$p\in(1,\infty)$.
	\vskip0.1in
	\noindent Theorem \ref{thm:main} connects three statements situated in distinct mathematical domains.  Condition \ref{condition 1: Kakeya-type sets} is geometric, tied to the existence of Kakeya-type configurations. Condition \ref{condition 2: operators unbdd} is analytic, centred on directional maximal operators $D_{\Omega}$ and $M_{\Omega}$. Condition \ref{condition 3: slopes sublacunary} is structural and combinatorial, offering a set-theoretic and an equivalent tree-theoretic description of $\Omega$. Thus Theorem \ref{thm:main} identifies admissible finite-order lacunarity as the key structural condition mediating between Kakeya geometry and analytic
boundedness of directional maximal operators.
%the exact structural invariant mediating between these geometric and analytic phenomena.
	%Thus, admissible finite-order lacunarity is exactly the threshold separating bounded from unbounded directional maximal behaviour in the plane, and demarcates the occurrence of Kakeya-type phenomena. 
	
	\section{Main innovations} \label{section: main contributions} 
	%Upon first glance, Theorem \ref{thm:main} may appear formally identical to the three-fold equivalence stated in \cite{Bateman}, with the adjustments necessary to address the points raised in \cite{{KrocThesis},{HRS2024a},{HRS2024b}}. However, there are essential distinctions underpinning this work: 
	%\vskip0.1in
%	While relying heavily on the insights gained from previous literature, the present work introduces a set of conceptual and methodological advances that are key to understanding the combinatorial and geometric dichotomies underpinning lacunarity and Kakeya-type behaviour. We list the significant points of distinction from earlier ideas and methods. 
	
	Theorem \ref{thm:main} prompts a natural question: what new ideas are required to establish the corrected characterization? While the statement resembles that of Bateman \cite{Bateman}, the proof proceeds along substantially different lines. The principal conceptual and methodological advances are summarized below.
	\subsection{Structural and tree-based corrections of finite-order lacunarity}
	A key contribution of this article is the {\em{formulation of finite-order
lacunarity in both set-theoretic and tree-theoretic terms}}, together
with a proof that the two formulations are equivalent.  The relevant statements appear in Definition \ref{defn: Admissible finite order lacunarity}, Propositions \ref{TREE-LACUNARY-TRADITIONAL} and \ref{TREE-LACUNARY-TRADITIONAL-CONVERSE}. In \cite{Bateman}, the recursive description of lacunarity was not fully aligned with the tree structure used in the proof. Consequently,
finite-order lacunarity as defined there was not equivalent to
finiteness of the splitting number, the tree invariant driving the argument. The admissible formulation introduced in Definition \ref{defn: Admissible finite order lacunarity} restores this alignment through 
	\vskip0.1in
	\begin{itemize} 
	\item {\em{precise positioning of lower order building blocks}} within lacunary gaps of a central lacunary sequence to construct sets of higher lacunarity, and  
	\vskip0.1in
	\item {\em{uniform control of the lacunarity constant}} across the potentially infinitely many lower order components located within these gaps.  
	\end{itemize} 
	\vskip0.1in 
	\noindent Together, these conditions yield a notion of lacunary clustering that is equivalent to the corresponding tree-theoretic formulation.
	%This produces a definition involving lacunary clustering that is demonstrably identical to the tree-theoretic structure it is meant to capture. 
	
%	The admissible formulation introduced in Definition \ref{defn: Admissible finite order lacunarity} clarifies the connection between the set-theoretic description of finite order lacunarity and distinctive features reflected in their tree representation.
	
	\subsection{Universal Euclidean separation via pruning} The article identifies a universally attainable form of {\em{Euclidean separation}}, obtained through pruning. As pointed out in \cite{{HRS2024a},{HRS2024b}}, some form of Euclidean separation among the slopes of $\Omega$ is a critical requirement for the proof strategy in \cite{Bateman} to work. Loosely speaking, the separation property is a technical condition imposed on the binary tree representing $\Omega$; it demands a specific alignment between the dyadic and Euclidean distances among certain vertices of the tree. The separation condition identified in \cite{{HRS2024a},{HRS2024b}} requires the existence of a constant $\eta > 0$ obeying the following: 
	\begin{equation} \label{Bateman-separation} 
	\left\{ 
	\begin{minipage}{0.8\textwidth}
	Suppose that $I_1, I_2$ are two dyadic intervals such that each half of $I_i$ contains a slope in $\Omega$, for $i = 1, 2$. If $I$ is the smallest dyadic interval containing both $I_1$ and $I_2$, then dist$(I_1, I_2) \geq \eta |I|$. 
	\end{minipage}
	\right\}
	\end{equation}
	Not all slope sets $\Omega$ enjoy this property. 
	\vskip0.1in
	\begin{itemize} 
	\item Nonetheless, a key observation is that {\em{the tree of every sublacunary slope set contains subtrees of arbitrarily large splitting number exhibiting a controlled Euclidean separation property.}} The pruning mechanism developed in Chapter \ref{Chapter: Pruning of the slope tree}  extracts such a subtree from the slope tree of any sublacunary set; see Proposition \ref{PRUNING STAGE 1}. 
	\vskip0.1in
	\item {\em{The resulting separation condition in the pruned slope tree, while universally attainable, is quantitatively weaker than the condition \eqref{Bateman-separation}}}. Indeed the examples constructed in \cite{{HRS2024a},{HRS2024b}} show that the stronger Euclidean separation criterion \eqref{Bateman-separation} cannot, in general, be achieved even after pruning; see Section \ref{section: HRS pruning}. 
	\end{itemize} 
	\vskip0.1in
	\noindent A key finding in this paper is that the weaker Euclidean separation is sufficient to sustain the Kakeya-type construction.
	% This universal Euclidean separation that we are able to guarantee is weaker than the one specified in \cite{{HRS2024a},{HRS2024b}}. This subset is constructed via a pruning mechanism, described in Chapter \ref{Chapter: Pruning of the slope tree}.  
	
	\subsection{Replacement of the previous definition of stickiness}  
	The weaker separation property necessitates a re-engineering of the argument \eqref{condition 3: slopes sublacunary} $\implies$ \eqref{condition 1: Kakeya-type sets} in Theorem \ref{thm:main}. The success of the sticky-map mechanism used in \cite{Bateman} for constructing Kakeya-type families of rectangles relied on a stronger Euclidean separation assumption \eqref{Bateman-separation} formulated in terms of dyadic intervals, or vertices on dyadic trees. This distance property is no longer available for general sublacunary slope sets. Accordingly, the revised proof must address two significant issues:  
	\vskip0.1in
	\begin{itemize} 
	\item {\em{the identification of the precise tree properties that are preserved under weak separation}}, and 
	\vskip0.1in
	\item {\em{the development of an alternative Kakeya-type construction that retains enough ``geometric stickiness'' at the level of rectangle families}} -- producing significant overlap in one region of the plane, while ensuring sufficient separation in another. 
	\end{itemize} 
	\vskip0.1in
	\noindent A previously successful ingredient for the construction of such families is the notion of a sticky tree-map, which is a prescription for slope allocation to rectangles rooted on a line. In \cite{{Bateman},{BatemanKatz}}, stickiness was enforced at {\em{every}} dyadic scale, without a corresponding verification of Euclidean separation among
the slopes. The counterexamples in \cite{{HRS2024a},{HRS2024b}} show that this step is not
valid in general.  By contrast, our root-to-slope map is required to be sticky only at carefully chosen scales dictated by the Euclidean separation properties of the pruned slope tree.
	\vskip0.1in
	\noindent We term this scale-selective mechanism \emph{inductive stickiness}.
It is formulated on a compressed representation of the pruned tree,
using the theory of compressed trees introduced in Chapter \ref{trees-section}; its
precise definition appears in Chapter \ref{chapter: compressions}. This mechanism preserves the
geometric core of the earlier argument while placing it within the corrected structural framework.

	\subsection{New probabilistic construction}  
The pruned slope tree and the inductive notion of stickiness require a corresponding redesign of the probabilistic part of the construction. In \cite{Bateman}, the Lebesgue smallness of the Kakeya-type configuration was established through a probabilistic argument based on standard results from Bernoulli percolation. The altered combinatorial geometry
of the pruned tree no longer permits a direct application of that framework. We therefore introduce a new randomization scheme adapted to the pruned tree. This modification is essential for producing Kakeya-type configurations under the weaker separation condition.

	\subsection{Intrinsic nature of the new definition}  
	Beyond its role in the proof of Theorem~\ref{thm:main}, the notion of admissible finite-order lacunarity is invariant under bi-Lipschitz changes of parametrization.  Earlier work often passed interchangeably between angular and slope parametrizations without explicitly justifying the preservation of finite-order lacunarity. 
%	 Using bi-Lipschitz invariance of admissible finite-order lacunarity, the distinction between describing directions by slopes or by angles disappears. We show in 
Sections \ref{section: bi-Lipschitz invariance} and \ref{section: slopes and angles} show that the two formulations are equivalent in the admissible setting, and track the corresponding changes in the lacunarity parameter and covering constant.
	%%We formalize this equivalence based on our definition via bi-Lipchitz invariance of finite order lacunarity. Ee clarify the role of the lacunarity and covering constants and emphasize the necessity of uniform quantitative bounds within a fixed lacunarity order.
	% We also clarify the role of the lacunarity constant and the covering constant, and the importance  of keeping these parameters uniformly bounded within a set of fixed lacunarity order.     
	\vskip0.1in
	\noindent Together, these results establish the planar boundedness dichotomy
using admissible finite-order lacunarity. The proof combines its
characterization by finite splitting number with pruning, sticky
maps, and a probabilistic construction adapted to the weaker
separation condition. Bi-Lipschitz invariance ensures that the
characterization applies equally to slope and angular
parametrizations.	

	\section{Applications} 
	\subsection{The counterexample of \cite{HRS2024b}} \label{HRS-application-section}
	In \cite{HRS2024a}, Hagelstein, Radillo-Murguia, and Stokolos  identified a Euclidean separation condition among the slopes in $\Omega$ as a fundamental requirement for the proof in \cite{Bateman} to succeed. In subsequent work \cite[Section 4]{HRS2024b}, they
				gave an example of a set that fails this separation property, proposing it as a critical test case for a corrected framework. We describe this set below and denote it by $\Omega_{\text{HRS}}$. 
				\vskip0.1in
				\noindent Let $\{N_j : j \geq 1\}$ be a sequence of positive integers $\geq 2$ such that 
				\begin{equation} \sum_{j=1}^{\infty} 2^{-N_j} < \infty. \label{N_j summability} \end{equation}  
				For each $j \geq 1$, let us define two binary strings $\pmb{\eta}_j$ and $\pmb{\zeta}_j$ in $\{0, 1\}^{N_j}$:
				\begin{equation}  \pmb{\eta}_j = (0, 1, 1, \ldots, 1), \qquad \pmb{\zeta}_j = (1, 0, 0, \ldots, 0) \label{eta-zeta} \end{equation} 
				The set in \cite{HRS2024b} contains a specially chosen class of dyadic rationals; the digit expansion in base 2 of a member of $\Omega_{\text{HRS}}$ consists of a finite sequence of binary blocks, the $j^{\text{th}}$ block (of length $N_j$) being either $\pmb{\eta}_j$ or $\pmb{\zeta}_j$. In other words,  
				\begin{align}  \Omega_{\text{HRS}} &:= \bigcup_{R=1}^{\infty} \Omega_{\text{HRS}}(R), \text{ where } \label{HRS-example} \\ 
				\Omega_{\text{HRS}}(R) &:= \left\{ \sum_{k=1}^{\infty} \frac{\varepsilon_k}{2^k} \; \Biggl| \;  
				\begin{aligned} &{\pmb{\varepsilon}} = (\varepsilon_1, \varepsilon_2, \ldots ) = \bigl({\pmb{\kappa}}_1, {\pmb{\kappa}}_2, \ldots, {\pmb{\kappa}}_R, 0 , 0 , \ldots \bigr), \\ & \pmb{\kappa}_j \in \{0, 1\}^{N_j}, {\pmb{\kappa}}_j = \text{ either } {\pmb{\eta}}_j \text{ or } {\pmb{\zeta}}_j, \; 1 \leq j \leq R \end{aligned} 
				\right\}. \label{HRS-example-unit} 
				\end{align} 
Here $(\boldsymbol{\kappa}_1,\ldots,\boldsymbol{\kappa}_R,
0,0,\ldots)$ denotes the concatenation of the finite binary blocks
$\boldsymbol{\kappa}_1,\ldots,\boldsymbol{\kappa}_R$, followed by an
infinite string of zeros.				
		\vskip0.1in		
\noindent The set $\Omega_{\text{HRS}}$ is not finite-order lacunary in the sense of  Sj\"ogren--Sj\"olin \cite{SjogrenSjolin}, so the corresponding positive
boundedness results do not apply. On the other hand, it fails the Euclidean separation condition required by Bateman's argument, so that argument does not establish unboundedness either.
The article \cite{HRS2024b} therefore asked whether $D_{\Omega}$ is $L^p$-bounded on the nontrivial Lebesgue spaces for $\Omega = \Omega_{\text{HRS}}$. Theorem \ref{thm:main} gives a negative answer. 				
				\begin{corollary} \label{Corollary: HRS-example}
					The set $\Omega_{\text{HRS}}$ in \eqref{HRS-example} is sublacunary in the sense of Definition \ref{defn: Admissible finite order lacunarity}. As a result, the maximal operators  $D_{\Omega_{\text{HRS}}}$ and $M_{\Omega_{\text{HRS}}}$ are unbounded on $L^p(\mathbb R^2)$  for all $p \in (1, \infty)$.  
					\end{corollary}  
					\noindent 
				The sublacunarity of $\Omega_{\mathrm{HRS}}$ is proved in Section \ref{HRS-application-proof-section} directly from the definitions. Combined with Theorem \ref{thm:main}, this resolves the question left open in \cite{HRS2024b} by placing $\Omega_{\mathrm{HRS}}$ on the sublacunary side of the dichotomy.	
					
	\subsection{Differentiation theorems} In classical analysis,  maximal function estimates provide a standard route to almost everywhere differentiation theorems. The corresponding theory of differentiation bases, summarized in \cite[Chapter X, \S 2, Proposition 1; also \S 2.2, Corollary 1]{SteinHA}, leads to the following $L^p$ differentiation statements as a consequence of Theorem \ref{thm:main}.
	\begin{corollary}
	Given a slope set $\Omega \subseteq \mathbb R$, let 
	\begin{equation} \mathscr{R}(\Omega) := \bigcup_{\omega \in \Omega} \mathscr{R}_{\omega} \label{R-Omega} \end{equation} denote the family of origin-centred rectangles in $\mathbb R^2$ whose orientations are determined by $\Omega$; the long side of every $R \in \mathscr{R}_{\omega}$ has slope $\omega \in \Omega$. 
	\vskip0.1in  
	\begin{enumerate}[(a)]
	\item Suppose that $\Omega$ is admissible finite-order lacunary in the sense of Definition  \ref{defn: Admissible finite order lacunarity}. Then for every $p \in (1, \infty)$ and every $f \in L^p(\mathbb R^2)$, 
	\[ \lim_{\begin{subarray}{c} \text{diam}(R) \rightarrow 0 \\ R \in \mathscr{R}(\Omega)\end{subarray}} \frac{1}{|R|} \int_R f(x+y) \, dy = f(x) \quad \text{ for Lebesgue almost every $x \in \mathbb R^2$.} \] 
	Here $\text{diam}(R)$ represents the diameter of the rectangle $R$. 
	\vskip0.1in
	\item Conversely,  suppose that $\Omega$ is sublacunary in the sense of Definition  \ref{defn: Admissible finite order lacunarity}. Then for every $p \in (1, \infty)$, there exists a non-negative function $f \in L^p(\mathbb R^2)$ for which the corresponding differentiation averages diverge almost everywhere. Specifically, 
	\[ \limsup_{\begin{subarray}{c} \text{diam}(R) \rightarrow 0 \\ R \in \mathscr{R}(\Omega)\end{subarray}} \frac{1}{|R|} \int_R f(x+y) \, dy = \infty \quad \text{ for Lebesgue almost every $x \in \mathbb R^2$.} \] 
	The quantity $\limsup_R$ in the display above is interpreted as follows: 
\[\limsup_{\begin{subarray}{c} \text{diam}(R) \rightarrow 0 \\ R \in \mathscr{R}(\Omega) \end{subarray}} = \lim_{\delta \rightarrow 0+} \sup_{\begin{subarray}{c}\text{diam}(R) < \delta \\ R \in \mathscr{R}(\Omega) \end{subarray}}. \]
	\end{enumerate} 
	\end{corollary}
	
	\section{Roadmap of the proofs} \label{section: proof organization}
Theorem \ref{thm:main} asserts the equivalence of three statements, one geometric, one analytic, and one combinatorial. Accordingly, the proof is organized around establishing the implications between these three formulations. While the implication from Kakeya-type sets to unboundedness of directional maximal operators is short and direct, the remaining two directions require substantially different ideas and constitute the main body of the manuscript.
\subsection{Existence of Kakeya-type sets implies unboundedness of $D_{\Omega}$} 
Chapter \ref{Kakeya-type-section} shows the implication \eqref{condition 1: Kakeya-type sets} $\implies$ \eqref{condition 2: operators unbdd}, namely that Kakeya-type configurations force unboundedness of the associated maximal operators. The argument is geometric: rectangles exhibiting Kakeya-type overlap produce test functions whose maximal averages remain large while their $L^p$-norms become arbitrarily small.

\subsection{Sublacunarity of $\Omega$ implies existence of Kakeya-type sets} 
The implication \eqref{condition 3: slopes sublacunary} $\implies$ \eqref{condition 1: Kakeya-type sets} is the principal focus of this paper and occupies Chapters~\ref{section: finite-order lacunarity}--\ref{chapter: reference tree geometry}. At a conceptual level, the argument transforms the combinatorial richness of a sublacunary slope set into a geometric configuration of tubes exhibiting Kakeya-type behaviour. 
\vskip0.1in
\noindent The proof of \eqref{condition 3: slopes sublacunary} $\implies$ \eqref{condition 1: Kakeya-type sets} begins by introducing the notion of admissible finite-order lacunarity in Chapter~\ref{section: finite-order lacunarity}, together with the examples and structural properties that will be needed later. The first substantive step is to replace the geometric problem by a combinatorial one. Chapters~\ref{trees-section} and~\ref{section: split implies lacunarity} encode the slope set by an $M$-adic tree and show that sublacunarity forces its splitting number to be infinite. Thus a sublacunary slope tree contains subtrees with arbitrarily large splitting number; these subtrees are used to identify the branching phenomenon that ultimately drives the Kakeya construction. The converse structural implication, namely admissible finite order lacunarity implies finite splitting number, is proved in Chapters \ref{chapter: lac-split-appendix}--\ref{section: proof of gap lemma}.
\vskip0.1in
\noindent Large branching alone, however, is insufficient for the geometric argument. The next stage therefore reorganizes the slope tree while preserving its essential combinatorial complexity. Chapter ~\ref{Chapter: Pruning of the slope tree} implements this reorganization through a pruning procedure that extracts a subtree with the structural properties required for the subsequent construction. Chapters~\ref{Chapter: Fundamental heights} and~\ref{chapter: compressions} then identify the distinguished scales (termed {\em{fundamental heights}}) of this pruned tree and recast it in a compressed form that makes its geometry more transparent.
\vskip0.1in
\noindent With this framework in place, the argument turns from combinatorics to geometry. Chapter~\ref{chapter: sticky rectangles chapter} associates with the compressed tree a family of thin tubes, denoted by $\mathbb T_{\mathbb X}$. Chapters~\ref{chapter: near the roots}--\ref{chapter: reference tree geometry} combine probabilistic methods with incidence estimates to show that the union $\mathtt K_{\mathbb X}$ of the tubes in $\mathbb T_{\mathbb X}$ satisfies the defining property of a Kakeya-type set with positive probability. 
\vskip0.1in 
\noindent More precisely, Chapter~\ref{chapter: near the roots} establishes that $\mathtt K_{\mathbb X}$ is large near the root line (where the tubes originate), while Chapters ~\ref{random construction section}--\ref{chapter: reference tree geometry} establish that, for a suitable choice of $\mathbb X$, it is small far from the root line. The final synthesis, presented in Section \ref{section: Sublac-Kakeya}, combines these two estimates to produce the required Kakeya-type sets and complete the implication \eqref{condition 3: slopes sublacunary} $\implies$ \eqref{condition 1: Kakeya-type sets}. The architecture of the argument is summarized in Figure \ref{fig:proof-roadmap-sublacunary-kakeya}. 
\begin{figure}[htbp]
\centering
\begin{tikzpicture}[
  font=\small,
  node distance=5mm,
  flowbox/.style={
    draw,
    rounded corners=3pt,
    align=center,
    inner xsep=7pt,
    inner ysep=5pt,
    text width=.64\textwidth
  },
  branchbox/.style={
    draw,
    rounded corners=3pt,
    align=center,
    inner xsep=5pt,
    inner ysep=5pt,
    text width=.35\textwidth
  },
  flowarrow/.style={->,>=latex}
]

\node[flowbox] (sub)
  {\textbf{Sublacunarity of $\Omega$}};

\node[flowbox, below=of sub] (split)
  {\textbf{Arbitrarily large splitting number}\\
   of the $M$-adic slope tree\\[-1pt]
   {\footnotesize Chapters \ref{trees-section} --\ref{section: split implies lacunarity}}};

\node[flowbox, below=of split] (prune)
  {\textbf{Pruning to a model slope tree}\\[-1pt]
   {\footnotesize Chapter \ref{Chapter: Pruning of the slope tree}}};

\node[flowbox, below=of prune] (compress)
  {\textbf{Fundamental heights, compressed trees,}\\
   \textbf{and inductively sticky maps}\\[-1pt]
   {\footnotesize Chapters \ref{Chapter: Fundamental heights}--\ref{chapter: compressions}}};

\node[flowbox, below=of compress] (tubes)
  {\textbf{Construction of the sticky tube family $\mathbb T_{\mathbb X}$}\\[-1pt]
   {\footnotesize Chapter \ref{chapter: sticky rectangles chapter} }};

\node[branchbox, below left=8mm and 3mm of tubes.south] (near)
  {\textbf{$\mathtt K_{\mathbb X}$ is large near the root line}\\
   Proposition \ref{prop: near the roots}\\[-1pt]
   {\footnotesize Chapter 10}};

\node[branchbox, below right=8mm and 3mm of tubes.south] (far)
  {\textbf{For some $\mathbb X$, $\mathtt K_{\mathbb X}$ is small}\\
   \textbf{far from the root line}\\
   Proposition \ref{prop: away from the roots}\\[-1pt]
   {\footnotesize Chapters \ref{random construction section}--\ref{chapter: reference tree geometry}}};

\coordinate (merge)
  at ($(near.south)!0.5!(far.south)+(0,-7mm)$);

\node[flowbox, below=12mm of merge] (combine)
  {\textbf{Combine the near- and far-root estimates}\\[-1pt]
   {\footnotesize Section \ref{section: Sublac-Kakeya}}};

\node[flowbox, below=of combine] (kakeya)
  {\textbf{Existence of Kakeya-type sets with slopes in $\Omega$}};

\draw[flowarrow] (sub) -- (split);
\draw[flowarrow] (split) -- (prune);
\draw[flowarrow] (prune) -- (compress);
\draw[flowarrow] (compress) -- (tubes);

\draw[flowarrow]
  (tubes.south) -- ++(0,-3mm) -| (near.north);
\draw[flowarrow]
  (tubes.south) -- ++(0,-3mm) -| (far.north);

\draw (near.south) |- (merge);
\draw (far.south) |- (merge);
\draw[flowarrow] (merge) -- (combine.north);

\draw[flowarrow] (combine) -- (kakeya);

\end{tikzpicture}
\caption{\small{Proof structure of the implication
\eqref{condition 3: slopes sublacunary} $\implies$ \eqref{condition 1: Kakeya-type sets} in Theorem \ref{thm:main}.}}
\label{fig:proof-roadmap-sublacunary-kakeya}
\end{figure}
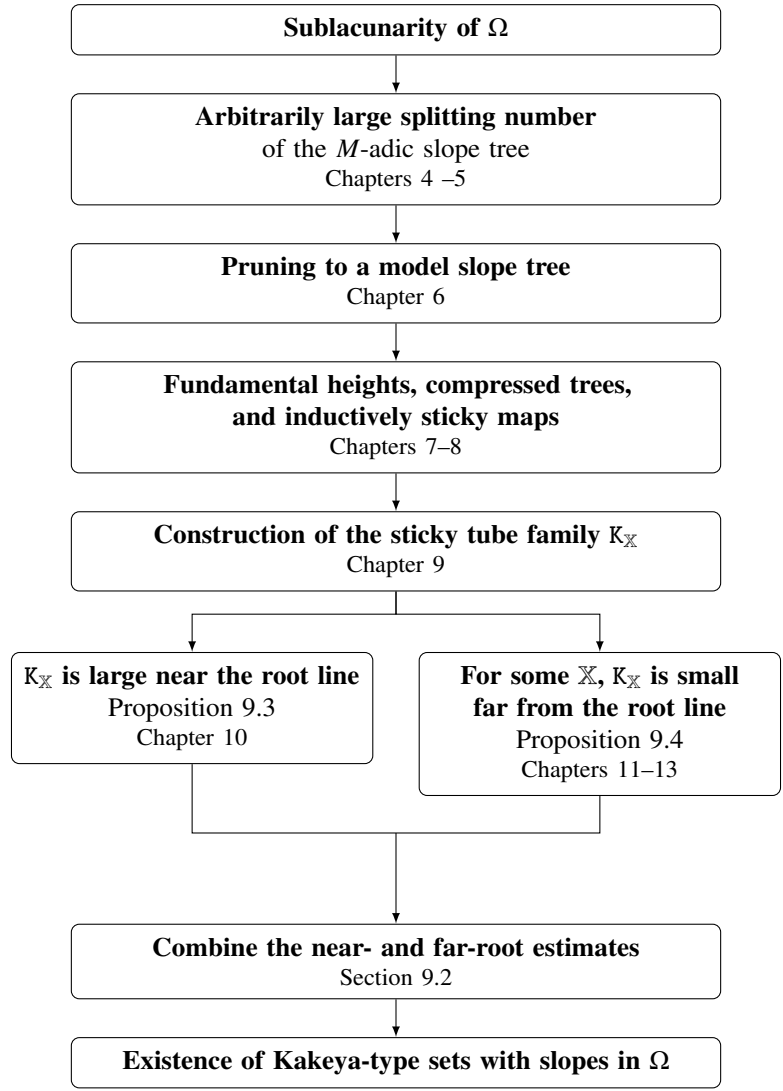

\subsection{Admissible finite-order lacunarity of $\Omega$  implies $D_{\Omega}$ is $L^p$-bounded} 
The remaining implication \eqref{condition 2: operators unbdd} $\implies$ \eqref{condition 3: slopes sublacunary}, proved via the contrapositive, is established in Chapter~\ref{positive direction chapter}. Unlike the previous implication, whose proof proceeds through a geometric construction of Kakeya-type sets, this direction is entirely analytic in nature. Its objective is to show that the hierarchical structure encoded by admissible finite-order lacunarity provides sufficient directional separation to guarantee $L^p$-boundedness of the associated maximal operators.
\vskip0.1in
\noindent The proof begins by replacing the directional averages with suitable smoothed averaging operators, allowing the problem to be analysed using Fourier-analytic techniques. The resulting operators are decomposed into low- and high-frequency contributions, each of which is estimated separately. The principal analytic tool is Christ's almost orthogonality principle, which exploits the multiscale structure of admissible finite-order lacunarity to control interactions between different frequency regimes. These estimates together establish $L^p$-boundedness of the directional maximal operators $D_\Omega$ and $M_\Omega$, completing the proof of Theorem~\ref{thm:main}.
\vskip0.1in
\noindent To preserve the flow of the main exposition, several auxiliary proofs and complementary results have been deferred to the appendices.

\begin{figure}[ht]
\centering

\begin{tikzpicture}[
    maintext/.style={
        align=center,
        font=\small
    },
    connection/.style={
        <->,
        line width=0.7pt
    }
]

% Structural viewpoint
\node[maintext] (lacunarity) at (0,4.2)
{
    \textbf{Admissible finite-order lacunarity of $\Omega$}
};

\node[maintext] (equiv) at (0,3.55)
{
    $\Updownarrow$
};

\node[maintext] (splitting) at (0,2.9)
{
    \textbf{Finite splitting number of the $M$-adic tree for $\Omega$}
};

% Geometric viewpoint
\node[maintext, text width=4.2cm] (kakeya) at (-4.0,-2.1)
{
    \textbf{Non-existence of}\\
    \textbf{Kakeya-type sets with slopes from $\Omega$}
};

% Analytic viewpoint
\node[maintext, text width=4.2cm] (operators) at (4.0,-2.1)
{
    \textbf{$L^p$-boundedness of the directional}\\
    \textbf{maximal operators $D_{\Omega}, M_{\Omega}$}
};

% Equivalences between the structural, geometric,
% and analytic formulations
\draw[connection]
    ([xshift=-1.1cm]splitting.south)
    --
    (kakeya.north);

\draw[connection]
    ([xshift=1.1cm]splitting.south)
    --
    (operators.north);

\end{tikzpicture}

\caption{\small{Theorem~\ref{thm:main} establishes the equivalence
of the structural, geometric, and analytic viewpoints shown above.}}

\label{fig:conceptual-map}
\end{figure}
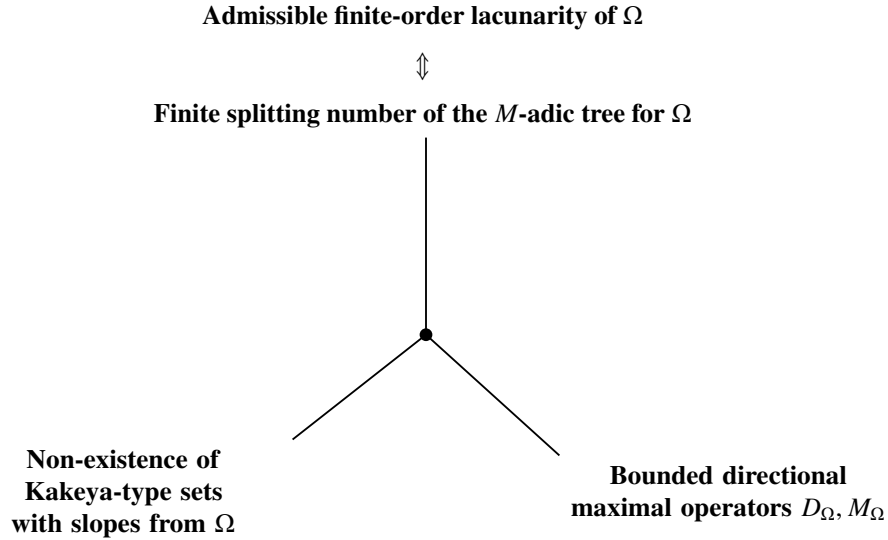
\clearpage
\section{Acknowledgements}
This work was initiated in 2025, during JL's visit to the University of British Columbia supported by a scholarship from Korea Institute for Advanced Study (KIAS). He would like to thank both organizations for making the visit possible. EK and MP were partially supported by Discovery Grants from the Natural Sciences and Engineering Research Council of Canada (NSERC). The authors thank Rafik Zeraoulia for their careful reading of the first arXiv version of this manuscript and for identifying several technical points; these observations led to corrections and clarifications in the present version.
\vskip0.1in
\noindent The authors gratefully acknowledge the assistance of OpenAI's ChatGPT (Version 5.5) with exposition, organization, and presentation of this manuscript, including the preparation of TikZ diagrams. All mathematical content, results, and any remaining errors are the sole responsibility of the authors. 

	\chapter{Directions admitting Kakeya-type sets} \label{Kakeya-type-section}
%Given $\Omega\subseteq [0,1]$, a rectangle in $\mathbb R^2$ is said to have orientation $\omega \in \Omega$ if its longer side has slope $\omega$.  
Condition \ref{condition 1: Kakeya-type sets} in Theorem \ref{thm:main} mentions slope sets $\Omega$ that {\em{admit Kakeya-type sets}}. In this chapter, we explain the meaning of this term, and prove the implication 
 \[ \text{condition \ref{condition 1: Kakeya-type sets}  implies condition \ref{condition 2: operators unbdd} in Theorem \ref{thm:main}}. \]
Namely, a slope set $\Omega$ that allows Kakeya-type configurations prevents the maximal operators $D_{\Omega}, M_{\Omega}$ from being bounded on any non-trivial Lebesgue space $L^p$. 

\section{Background and motivation}
Let us take a moment to review a key geometric property involving oriented rectangles in the plane that inspires the definition to follow. 
The standard Kakeya-Besicovitch set in $\mathbb R^2$ contains, by definition, a unit line segment in every direction. Though seemingly large through its inclusion of many lines, such a set can be Lebesgue-null. This has led to a vibrant line of inquiry that aims to quantify the minimal size of these sets  \cite{{Davies}, {Bourgain1}, {Bourgain2}, {Wolff}, {Wolff2}, {KatzLabaTao}, {KatzTao}, {KatzTao2}, {KatzZhal}, {Tao-blog}, {Guth1}, {Guth2}, {WangZahl1}, {WangZahl2}, {WangZahl3}, {GuthWangZahl}}.  
\vskip0.1in
\noindent Certain counter-intuitive geometric properties enjoyed by special classes of Kakeya-Besicovitch sets have played a pivotal role in many analytical problems \cite{{Tao-NotAMS},{Hickman}}, making such sets a valuable source for counter-examples in Fourier-analytic and operator-theoretic literature. One such property concerns incidences among thin rectangles \cite[Chapter X, Theorem 1]{SteinHA}: 
\vskip0.1in
\begin{center}
\begin{minipage}{0.9\textwidth}
{\em{For every $\varepsilon > 0$, there exists a collection $\mathscr{R}_N$ of $N$ rectangles of dimension $1 \times N^{-1}$ oriented in roughly $N^{-1}$-separated directions, such that the sets 
\[ K_N := \bigcup \bigl\{R : R \in \mathscr{R}_N \bigr\} \; \text{ and } \; \hat{K}_N := \bigcup \bigl\{ \hat{R} : R \in \mathscr{R}_N\bigr\} \] satisfy the relation 
\[|K_N| < \varepsilon, \; \text{ but } \; |\hat{K}_N| = 1. \]  Here $\hat{R}$ is a congruent copy of $R$ obtained by sliding $R$ along its long direction by two units. The set $K_N$ is, of course, a Kakeya-Besicovitch set.}}
\end{minipage}
\end{center} 
\vskip0.1in 
In other words, the rectangles $R$ in $\mathscr{R}_N$ overlap heavily, so that their total area is arbitrarily small, yet their translates $\hat{R}$ are entirely disjoint, creating a set of maximal size. It is this disparity between the size of $K_N$ and $\hat{K}_N$ that the following definition aims to generalize.  
\section{Admission of Kakeya-type sets} 
Given $\Omega \subseteq \mathbb R$, let us recall from \eqref{R-Omega} the definition of the collection of origin-centred rectangles $\mathscr R(\Omega)$ oriented along the directions specified by $\Omega$. Let $\mathscr{R}^{\ast}(\Omega)$ denote the collection of all possible translates of rectangles in $\mathscr R(\Omega)$:
\[ \mathscr{R}^{\ast}(\Omega) := \bigcup_{x \in \mathbb R^2} x + \mathscr{R}(\Omega). \] 
In other words, a rectangle is a member of $\mathscr{R}^{\ast}(\Omega)$ if and only if the slope of its long side lies in $\Omega$, with no restriction on its centre. 
\begin{definition}\label{Kakeya-type set}
 We say that $\Omega$ \textit{admits Kakeya-type sets} if there exist a constant $A_0 > 1$ and an infinite sequence of positive integers $\mathbb N_{\Omega}$  with the following property. 
 \vskip0.1in
 \noindent For each $N \in \mathbb N_{\Omega}$, one can find 
\begin{itemize} 
\item parameters $\delta_N \in (0, 1)$, $\delta_N \searrow 0$, and 
\vskip0.1in 
\item a collection of rectangles $\mathscr R_{N} \subseteq \mathscr R^{\ast}(\Omega)$ each of length at least 1 and width at most $\delta_N$,
\end{itemize}
 such that the families of sets $\{E_N : N \in \mathbb N_{\Omega}\}$ and $\{E_N^{\ast} : N \in \mathbb N_{\Omega}\}$ given by 
 \[ E_N := \bigcup \bigl\{R : R \in \mathscr R_N \bigr\}, \quad  E_N^*(A_0) := \bigcup \bigl\{A_0 R : R \in \mathcal R_N \bigr\} \; \text{ for } \; N \in \mathbb N_{\Omega} \]
 obey the condition
\begin{equation}\label{Kakeya-type condition}
\lim_{\begin{subarray}{c}N\rightarrow\infty\\ N \in \mathbb N_{\Omega}\end{subarray}}\frac {|E^*_N(A_0)|}{|E_N|} = \infty. 
\end{equation}
Here, 
%$|\cdot|$ denotes $(d+1)$-dimensional Lebesgue measure, and $
$A_0 R$ denotes the rectangle with the same centre, orientation and width as $R$ but $A_0$ times its length; it is referred to as the reach of $R$. 
\vskip0.1in
%The tubes that constitute $E_N$ may have variable dimensions subject to the restrictions mentioned above. 
\noindent A family of sets $\{E_N\}$ with the property \eqref{Kakeya-type condition} is said to be of Kakeya type. 
\end{definition}
\noindent To paraphrase, each $E_N$ is a union of thin rectangles oriented along $\Omega$, whose eccentricity (i.e. the ratio of short to long sides) $\delta_N$ vanishes as $N \rightarrow \infty$. 
%the $N^{\text{th}}$ member of a Kakeya-type family is a set $E_N$ consisting of thin rectangles with orientations dictated by $\Omega$ and eccentricity (i.e. the ratio of short to long sides) at most $\delta_N$. 
If $\Omega$ admits Kakeya-type sets, these rectangles can be arranged to ensure that the area of $E_N$ grows progressively smaller in comparison to that of $E_N^{\ast}$, the union of the reach rectangles, as $N$ grows without bound. A schematic diagram of this phenomenon is presented in Figure \ref{figure: Kakeya type sets}. 
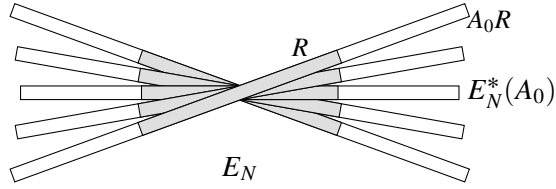
\begin{figure}[ht] 
\centering
\begin{tikzpicture}[scale=1] 

% Representative family of rectangles:
% Outer = A_0 R (outline only)
% Inner = R (light gray fill)

% Rectangle 1
\begin{scope}[rotate around={-20:(0,0)}]
    \draw[line width=0.4pt] (-3.2,-0.09) rectangle (3.2,0.09);
    \filldraw[fill=gray!25, draw=black, line width=0.4pt]
        (-1.4,-0.09) rectangle (1.4,0.09);
\end{scope}

% Rectangle 2
\begin{scope}[rotate around={-10:(0,0)}]
    \draw[line width=0.4pt] (-3.0,-0.09) rectangle (3.0,0.09);
    \filldraw[fill=gray!25, draw=black, line width=0.4pt]
        (-1.35,-0.09) rectangle (1.35,0.09);
\end{scope}

% Rectangle 3
\begin{scope}[rotate around={0:(0,0)}]
    \draw[line width=0.4pt] (-2.9,-0.09) rectangle (2.9,0.09);
    \filldraw[fill=gray!25, draw=black, line width=0.4pt]
        (-1.3,-0.09) rectangle (1.3,0.09);
\end{scope}

% Rectangle 4
\begin{scope}[rotate around={10:(0,0)}]
    \draw[line width=0.4pt] (-3.0,-0.09) rectangle (3.0,0.09);
    \filldraw[fill=gray!25, draw=black, line width=0.4pt]
        (-1.35,-0.09) rectangle (1.35,0.09);
\end{scope}

% Rectangle 5
\begin{scope}[rotate around={20:(0,0)}]
    \draw[line width=0.4pt] (-3.2,-0.09) rectangle (3.2,0.09);
    \filldraw[fill=gray!25, draw=black, line width=0.4pt]
        (-1.4,-0.09) rectangle (1.4,0.09);
\end{scope}

% Labels
\node at (0,-1) {$E_N$};
\node at (3.6,0) {$E_N^{*}(A_0)$};

% Optional labels for one representative rectangle pair
\node at (0.8,0.6) {\small $R$};
\node at (3.3,0.95) {\small $A_0R$};

\end{tikzpicture}
\caption{\small{A schematic Kakeya-type family of rectangles. The shaded rectangles represent the members of a family $\mathscr{R}_N$, whose union is $E_N$. The larger outlined rectangles represent their reaches $A_0R$, with the same center, width, and orientation but greater length; their union is $E_N^{*}(A_0)$. The rectangles in $E_N$ overlap heavily, while the corresponding reaches occupy a much larger region.}}
\label{figure: Kakeya type sets}
\end{figure}

\vskip0.1in
\noindent The main observation in this section is the following. 
\begin{lemma} \label{lemma: 1 implies 2}
Suppose that $\Omega \subseteq \mathbb R$ admits Kakeya-type sets. Then $D_{\Omega}$ and $M_{\Omega}$ are both unbounded on $L^p(\mathbb R^2)$ for every $p \in [1, \infty)$. In other words, condition \ref{condition 1: Kakeya-type sets} of Theorem \ref{thm:main}  implies condition \ref{condition 2: operators unbdd}. 	
	\end{lemma}
\begin{proof}
Condition \ref{condition 1: Kakeya-type sets} of Theorem \ref{thm:main} guarantees the existence of a Kakeya-type family of sets $\{E_N: N \in \mathbb N_{\Omega} \}$ comprising rectangles that obey \eqref{Kakeya-type condition}. We use this family of sets to construct a sequence of test functions $1_{E_N}$, whose limiting behaviour establishes the Lebesgue unboundedness of $D_{\Omega}$.  As a consequence of the pointwise inequality \eqref{DM-pointwise-ineq} and known bounds on the Hardy-Littlewood maximal function, we deduce that the $L^p$ operator norms of $D_{\Omega}$ and $M_{\Omega}$ are comparable for $p \in (1, \infty)$. Namely, these norms are bounded above and below by constant multiples of each other, therefore unboundedness of $D_{\Omega}$ implies the same for $M_{\Omega}$.   
	%Indeed, a standard argument shows that 
	\vskip0.1in
	\noindent  Let $R$ be a planar rectangle whose long side has length $\ell$ and slope $\omega \in \Omega$. We consider the line segment $L_x$ of length $2A_0 \ell$ and slope $\omega$ centred at a point $x \in A_0R$. The line
	\begin{equation} \label{defn: L_x}
	L_x = \bigl\{x + t (1, \omega) : t \in [-A_0 h, A_0 h] \bigr\}, \; h = \frac{\ell}{\sqrt{1 + \omega^2}}
	\end{equation}  passes through $R$ with essentially maximal intersection, i.e. along a maximal line segment in $R$ parallel to its long axis. Specifically, there exists an interval $T_x$ of length $h$ such that 
	\begin{equation} \label{geometry - line and rectangle}  
	L_x \cap R \supseteq \bigl\{ x + t(1, \omega) : t \in T_x \bigr\}; \;  |T_x| = h. 
	\end{equation}   
	The containment \eqref{geometry - line and rectangle} is a geometric fact best understood through a diagram (see Figure \ref{fig:line-intersection}), but we include a rigorous justification in Lemma \ref{lemma: line and rectangle} below. 
	\vskip0.1in
	\noindent Assuming \eqref{geometry - line and rectangle} for now, let us proceed to bound from below the directional average of $1_R$ over the line $L_x$. This results in the following inequality:  
\begin{align*} \frac{1}{2A_0 h} \int_{-A_0h}^{A_0 h} 1_{R}(x + t(1,\omega)) \, dt & =\frac{|T_x|}{2A_0 h} \\ &\geq \frac{h}{2A_0 h} = \frac{1}{2A_0} \;\text{ for all } x \in A_0 R. 
\end{align*} 
Choosing $R$ to be each of the constituent rectangles of the set $E_N$ defined as in \eqref{Kakeya-type condition},  the inequality above implies that for any $x \in E_N^{\ast}(A_0)$, 
\begin{equation*} 
	M_{\Omega}1_{E_N}(x) \geq D_{\Omega}1_{E_N}(x) \geq c_0 = (2A_0)^{-1} > 0.
	\end{equation*} The first inequality in the display above follows from the left inequality in \eqref{DM-pointwise-ineq}.  Although stated for smooth compactly supported functions,
the left inequality also holds almost everywhere for the
indicator functions used here, by a standard approximation argument. As a result,
	\begin{equation*}
	||M_{\Omega}||_{p \rightarrow p} \geq ||D_{\Omega}||_{p \rightarrow p} \geq \frac{c_0||1_{E_{N}^{\ast}(A_0)}||_p}{||1_{E_N}||_p} \geq c_0 \left(\frac{|E_N^{\ast}(A_0)|}{|E_N|} \right)^{\frac{1}{p}}.\end{equation*} 
If $\Omega$ admits Kakeya-type sets, Definition \ref{Kakeya-type set} ensures that the sets $E_N$ can be chosen so that the right hand side approaches infinity as $N \rightarrow \infty$ for $p \in [1, \infty)$. Hence both $M_{\Omega}$ and $D_{\Omega}$ are unbounded on $L^p(\mathbb R^2)$ for every $p \in [1, \infty)$. 
Thus, condition \ref{condition 1: Kakeya-type sets} implies condition \ref{condition 2: operators unbdd} in Theorem \ref{thm:main}. 
\end{proof} 
\noindent To complete the argument, we need to fill in the details of the step \eqref{geometry - line and rectangle}
\begin{lemma} \label{lemma: line and rectangle} 
For a rectangle $R$ and a line segment $L_x$ defined as in the proof of Lemma \ref{lemma: 1 implies 2}, the inclusion \eqref{geometry - line and rectangle} holds.  
\end{lemma} 
\begin{proof}
Without loss of generality and after a translation, we may assume that the rectangle $R$ is centred at the origin. Denoting the long and short directions of $R$ as 
\[
\vec v=\frac{(1,\omega)}{\sqrt{1+\omega^2}},
\; \text{ and } \; 
\vec v^\perp=\frac{(-\omega,1)}{\sqrt{1+\omega^2}}
\]
respectively, we may parametrize $R$ as 
\begin{equation} \label{parametrized R} 
R=\Bigl\{y: |y\cdot\vec v|\leq \frac{\ell}{2},\
          |y\cdot\vec v^\perp|\leq \frac{w}{2} \Bigr\}.
\end{equation} 
Let us fix $x\in A_0R$; then $x$ satisfies 
\begin{equation} \label{where does x lie} 
|x\cdot\vec v|\leq \frac{A_0\ell}{2},
\quad \text{ and } \quad
|x\cdot\vec v^\perp|\leq \frac{w}{2}.
\end{equation} 
Setting $s=t\sqrt{1+\omega^2}$, the line $L_x$ given by \eqref{defn: L_x} can be re-parametrized as 
\[
L_x=\{x+s\vec v: |s|\leq A_0\ell\}.
\]
In view of the defining condition \eqref{parametrized R} of $R$ and the choice \eqref{where does x lie} of $x$, the portion of the line through $x$ lying in $R$ is
\begin{align*}
\widetilde L_x &= \Bigl\{x+s\vec v : |(x+ s \vec{v}) \cdot \vec{v}| \leq \frac{\ell}{2},  \, |(x+ s \vec{v}) \cdot \vec{v}^{\perp}| \leq \frac{w}{2} \Bigr\} \\  
&=\Bigl\{x+s\vec v:
       -x\cdot\vec v-\frac{\ell}{2}
       \leq s\leq
       -x\cdot\vec v+ \frac{\ell}{2} \Bigr\} \\ 
       &= \Bigl\{x+t(1, \omega):
       -\frac{x\cdot(1, \omega)}{1 + \omega^2}-\frac{h}{2}
       \leq t \leq
       -\frac{x\cdot (1, \omega)}{1 + \omega^2}+ \frac{h}{2} \Bigr\}
\end{align*}
Then $\tilde{L}_x \subseteq R$ by construction. Moreover, \eqref{where does x lie} yields 
\[
|t|\leq \frac{|x\cdot \vec{v}|}{\sqrt{1 + \omega^2}}+ \frac{h}{2}
    \leq \frac{A_0 \ell}{2 \sqrt{1 + \omega^2}} + \frac{h}{2} \leq (A_0+1) \frac{h}{2}
     \leq A_0h,
\]
because $A_0>1$. Therefore $\widetilde L_x\subseteq L_x$.
Its Euclidean length is $\ell$, corresponding to a $t$-interval
of length $h$. This proves the claim \eqref{geometry - line and rectangle}.
\end{proof}
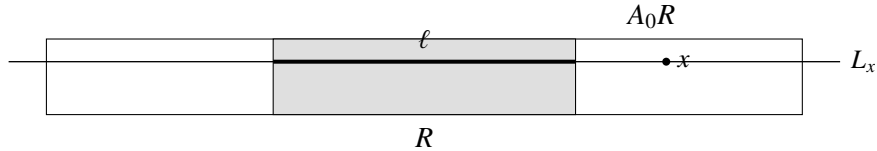
\begin{figure}[ht]
\centering
\begin{tikzpicture}[scale=1]

% Parameters
\def\length{4}      % length of R
\def\w{1}        % width of R
\def\A{2.5}        % A_0

% Outer rectangle A_0 R
\draw[line width=0.4pt]
    (-\A*\length/2, -\w/2) rectangle (\A*\length/2, \w/2);

% Inner rectangle R (filled)
\filldraw[fill=gray!25, draw=black, line width=0.4pt]
    (-\length/2, -\w/2) rectangle (\length/2, \w/2);

% Point x in A_0 R
\fill (3.2,0.2) circle (1.5pt);
\node[right] at (3.2,0.2) {\small $x$};

% Line segment L_x (long)
\draw[line width=0.5pt]
    (-5.5,0.2) -- (5.5,0.2);

\node[right] at (5.5,0.2) {\small $L_x$};

% Highlight segment inside R
\draw[line width=1.5pt]
    (-\length/2,0.2) -- (\length/2,0.2);

% Labels
\node at (0,-0.8) {$R$};
\node at (3,0.8) {$A_0R$};

\node at (0,0.5) {\small $\ell$};

\end{tikzpicture}

\caption{\small{Geometric illustration of (2.2). For any point $x$ in the expanded rectangle $A_0R$, the line segment $L_x$ centred at $x$ of length $2A_0 \ell$ parallel to the long side of $R$ intersects $R$ in a segment of length $\ell = \text{length of }R$.}}
\label{fig:line-intersection}
\end{figure}
\vskip0.1in
\noindent Reviewing the situation thus far, we have proved the implication \ref{condition 1: Kakeya-type sets} $\implies$ \ref{condition 2: operators unbdd} of Theorem \ref{thm:main}. The remaining implications require more intricate arguments. The remainder of the paper is devoted to proving these implications.

	\chapter{Admissible finite-order lacunarity} \label{section: finite-order lacunarity}
	\section{Chapter overview}
This chapter introduces the notion of admissible finite-order lacunarity, which is the core of Theorem \ref{thm:main}. The definition is recursive in nature and encodes a hierarchy of clustering within a set of slopes.
\vskip0.1in
\noindent Section \ref{basic definitions section} presents the formal definitions, beginning with lacunary sequences and building up to lacunary sets of arbitrary finite order and their admissible variants. Section \ref{section: Lp boundedness} records the main analytic consequence of this notion, namely $L^p$ boundedness of directional maximal operators for admissible finite-order lacunary slope sets.
\vskip0.1in
\noindent To build intuition, Section \ref{EXAMPLES SECTION} develops a range of examples and counterexamples, illustrating both admissible lacunary and sublacunary behaviour, and clarifying the distinctions with earlier formulations in the literature. Finally, Section \ref{SECTION: LACUNARY PROPERTIES} collects basic structural properties of admissible finite-order lacunary sets that will be used throughout the remainder of the paper.

\section{Basic definitions} \label{basic definitions section}
The notion of admissible finite-order lacunarity, which appears in condition \ref{condition 3: slopes sublacunary} of Theorem \ref{thm:main} is central to this article. Before proceeding further, we give the precise formulation of this concept. The definition is inductive on the lacunary order, with lacunary sequences serving as the foundation. Roughly speaking, higher-order lacunary sets are obtained by inserting lower-order lacunary sets into the gaps of a controlling lacunary sequence.
\subsection{Lacunary sequences and gaps}    
	\begin{definition} \label{defn: Lacunary sequence in R}
		Let $A = \{a_1,a_2,\ldots\}$ be a bounded infinite sequence of points on $\mathbb{R}$. Given a constant $0 < \lambda < 1$, we say that $A$ is a {\em{lacunary sequence}} with {\em{ lacunarity constant}} at most $\lambda$, if there exists $\alpha \in \mathbb R$ for which 
		\begin{equation}\label{lacunarity constant}
			|a_{j+1}-\alpha| \leq \lambda |a_j-\alpha| \quad \text{ for all } j \geq 1.
		\end{equation}
		 \end{definition} 
	\noindent A sequence $A$ obeying \eqref{lacunarity constant} must converge to the limit $\alpha$. Subsequences of $A$ lying on either side of $\alpha$ are themselves lacunary and additionally monotone, possibly after re-indexing. Monotone lacunary sequences will serve as the fundamental building blocks in the recursive construction of lacunary sets of higher order. The collection of all bounded, monotone, lacunary sequences on $\mathbb R$ with lacunarity constant at most $\lambda$ is denoted by ${\tt{MonLac}}(\lambda)$. 
	\vskip0.1in
	\noindent We say that an interval $I$ is a {\em{gap interval of the lacunary sequence $A$}} if the endpoints $a, b$ of $I$ are two neighbouring elements of the sequence $A$, in the Euclidean sense:
	\begin{equation} \label{gap interval}
	a, b \in A, \quad a < b, \quad (a, b) \cap A = \emptyset. 
	\end{equation}   
A gap interval may be open, closed, or half open and half closed. 	
\vskip0.1in
\noindent As we will see in the next definition, a lacunary set of order 1 is essentially a lacunary sequence of points;  gap intervals serve as the regions into which lower-order lacunary sets may be inserted to build lacunary sets of higher order. Figure \ref{fig:lacunary-gaps} contains a schematic representation of a lacunary sequence and its gap intervals. 
	\begin{figure}[ht]
\centering
\begin{tikzpicture}[scale=1]

% number line
\draw[->] (-0.5,0) -- (8.2,0);

% points of a lacunary sequence converging to alpha
\fill (0.8,0) circle (1.5pt);
\fill (2.6,0) circle (1.5pt);
\fill (4.0,0) circle (1.5pt);
\fill (5.1,0) circle (1.5pt);
\fill (6.3,0) circle (1.5pt);

% labels
\node[below] at (0.8,0) {$a_1$};
\node[below] at (2.6,0) {$a_2$};
\node[below] at (4.0,0) {$a_3$};
\node[below] at (5.1,0) {$a_4$};
\node[below] at (6.3,0) {$\alpha$};

% dashed guide to alpha
\draw[dashed] (6.3,-0.7) -- (6.3,0.7);

% gap intervals shown slightly above axis
\draw[line width=1pt] (0.8,0.45) -- (2.6,0.45);
\draw[line width=1pt] (2.6,0.75) -- (4.0,0.75);

% small endpoint ticks for gap intervals
\draw (0.8,0.38) -- (0.8,0.52);
\draw (2.6,0.38) -- (2.6,0.52);

\draw (2.6,0.68) -- (2.6,0.82);
\draw (4.0,0.68) -- (4.0,0.82);

% labels for gaps
\node at (1.7,0.8) {\small gap interval};
\node at (3.3,1.1) {\small gap interval};

% optional indication of convergence
\node[above] at (4.9,-0.05) {\small $a_j \to \alpha$};

\end{tikzpicture}
\caption{\small{A lacunary sequence $A = \{a_j\}$ converging to a limit point $\alpha$. The open intervals between consecutive terms are gap intervals; these serve as the regions into which lower-order lacunary sets are inserted in Definition 3.2.}}
\label{fig:lacunary-gaps}
\end{figure}
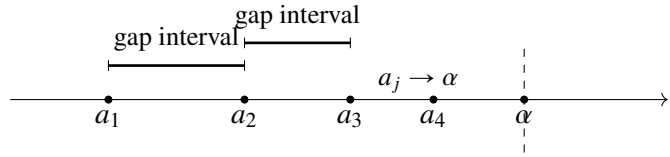

%a lacunary set of a larger order is obtained by inserting lacunary sets of lesser order into gap intervals of a lacunary sequence.  
\subsection{Lacunary sets in $\mathbb R$ of order $N$} 	\label{section: Lacunary sets of order N}
	\begin{definition} \label{defn: Lacunary sets}
		A {\em{lacunary set of order 0}} is, by definition, of cardinality at most 1, i.e., it is either empty or a singleton.  
		\vskip0.1in
		\noindent Recursively, given a constant $0<\lambda<1$ and an integer $N \geq 1$, we say that a bounded subset $U$ of $\mathbb{R}$ is a {\em{lacunary set of order at most $N$}} with lacunarity constant at most $\lambda$, and write 
		\[ U \in \Lambda(N, \lambda), \] if there exists a sequence $A \in {\tt{MonLac}}(\lambda)$ obeying the properties below:
		\begin{itemize}
			\item The set $U$ is contained in the interval spanned by the extreme values of $A$: 
			\begin{equation} \label{AU-bounds}
			U \subseteq \bigl[\inf(A), \sup(A) \bigr]. 
			%U \cap (\sup(A), \infty) = \emptyset, \quad U \cap (-\infty, \inf(A)) = \emptyset. 
			\end{equation}  
			\item The portion of $U$ lying in any gap interval of $A$ is of lower order. In other words,  for any interval $I = [a, b)$ with $a, b$ as in \eqref{gap interval}, the set $U$ satisfies   
			\begin{equation} \label{U in a gap}
			U \cap I \in \Lambda(N-1, \lambda). 
			\end{equation}    
		\end{itemize} 
		The {\em{order of lacunarity}} of $U$ is said to be exactly $N$ if $U \in \Lambda(N, \lambda) \setminus \Lambda(N-1, \lambda)$. 
		\vskip0.1in
		\noindent In other words, a lacunary set of order $N$ is built by placing sets of order $(N-1)$ inside the gaps of a controlling lacunary sequence $A$.
		\end{definition}
		\noindent A sequence $A \in {\tt{MonLac}}(\lambda)$ that obeys conditions \eqref{AU-bounds} and \eqref{U in a gap} above  for a given set $U \in \Lambda(N, \lambda)$ will be called a {\em{special sequence}} for $U$ and its limit will be termed a {\em{special point}}. The sequence $A$ serves as a structural scaffold for $U$: its extremal points bound $U$, and its gap intervals determine where lower-order clustering may occur.
		\vskip0.1in
		\noindent Thus the recursive definition enforces a strict hierarchy: clustering at level $N$ is controlled by a lacunary backbone of level 1, while any further clustering is confined to gap intervals and decreases the lacunarity order by at least one at each step. See Figure \ref{fig: finite order lacunary} for a diagram depicting a set of lacunary order 2. It is crucial that the lacunarity constants associated with the lacunary components $U \cap I$ of lower order remain under uniform control as $I$ ranges over the infinitely many gaps of the special sequence $A$.  
\begin{figure}[ht]
\centering
\begin{tikzpicture}[x=12cm,y=1cm,>=stealth]

% Real line
\draw[thick] (0,0) -- (1.02,0);
\node[above] at (1.01,0) {$\mathbb{R}$};

% Outer lacunary sequence points accumulating to the right
\foreach \x in {0.08,0.20,0.36,0.52,0.66,0.77,0.85,0.91,0.95,0.975} {
    \fill (\x,0) circle (1.2pt);
}

% Labels for first few a_j
\node[below=4pt] at (0.08,0) {$a_1$};
\node[below=4pt] at (0.20,0) {$a_2$};
\node[below=4pt] at (0.36,0) {$a_3$};
\node[below=4pt] at (0.52,0) {$a_4$};
\node[below=4pt] at (0.975,0) {$a_\infty$};

% Dots for continuation
\node[below=2pt] at (0.71,0) {$\cdots$};

% Highlighted gap intervals J1, J2, J3
%\draw[line width=1.2pt] (0.10,0.05) -- (0.18,0.05);
%\draw[line width=1.2pt] (0.38,0.05) -- (0.50,0.05);
%\draw[line width=1.2pt] (0.79,0.05) -- (0.89,0.05);

\node[above=6pt] at (0.14,-0.8) {\tiny{$I_1$}};
\node[above=6pt] at (0.44,-0.8) {\tiny{$I_2$}};
\node[above=6pt] at (0.84,-0.8) {\tiny{$I_3$}};

% First-order lacunary subsets inside the highlighted intervals
% J1
\foreach \x in {0.112,0.132,0.148,0.161,0.171} {
    \fill (\x,0) circle (1.8pt);
}

% J2
\foreach \x in {0.392,0.420,0.442,0.460,0.474,0.485} {
    \fill (\x,0) circle (1.8pt);
}

% J3
\foreach \x in {0.802,0.828,0.848,0.864,0.876,0.885} {
    \fill (\x,0) circle (1.8pt);
}

% Main labels on top
%\node at (0.50,1.35) {$U \in \Lambda(2,\lambda)$};
\node at (0.50,1.2) {\small{$U \cap I_r \in \Lambda(1,\lambda)$}};

% Arrows from secondary label to intervals
\draw[->] (0.42,0.92) -- (0.14,0.16);
\draw[->] (0.50,0.90) -- (0.44,0.16);
\draw[->] (0.58,0.92) -- (0.84,0.16);

\end{tikzpicture}
\caption{\small{A schematic picture of a set $U \in \Lambda(2,\lambda)$: elements of $U$ are interspersed within the special lacunary sequence $A = \{a_j\}$. For gap intervals $I_r$ that intersect $U$, the restricted set $U\cap I_r$ is itself in $\Lambda(1,\lambda)$.}}
\label{fig: finite order lacunary}
\end{figure}
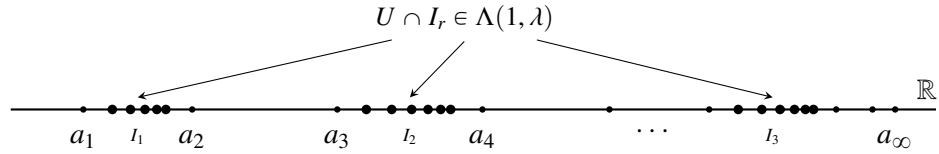
		\subsection{Admissible finite-order lacunarity, and sublacunarity} 
		Admissible finite-order lacunary sets referenced in Theorem \ref{thm:main}  are those generated by finite unions of sets of the form described in Definition \ref{defn: Lacunary sets}. The following definition of an admissible class ensures stability under finite unions.
	\begin{definition} \label{defn: Admissible finite order lacunarity} 
	Fix $\lambda \in (0,1)$ and integers $N, R \geq 1$.  
	\begin{enumerate}[1.] 
	\item \label{AdFinLac bounded sets}Suppose first that $U \subseteq \mathbb R$ is a bounded set. We write 
		\begin{equation} 
		U \in \Lambda(N, \lambda; R) \text{ and call $U$ an admissible lacunary set of order at most $N$ } \end{equation} with lacunarity constant at most $\lambda$ and covering number at most $R$ if $U$ is contained in the union of at most $R$ sets, each in $\Lambda(N, \lambda)$. The integer $R$ will be referred to as the covering constant. 
		\vskip0.1in
		\noindent The lacunarity order of $U \in \Lambda(N, \lambda; R)$ is said to be exactly $N$ if  
\[ U \notin \bigcup_{n=1}^{\infty} \Lambda(N-1, \lambda; n).  \] 
		\vskip0.1in
		 \noindent A bounded set $U \subseteq \mathbb{R}$ is termed {\em{admissible finite-order lacunary}}, and denoted 
		 \[ U \in {\tt{AdFinLac}}, \] if  there is a choice of $\lambda \in (0,1)$ and integers $N, R \geq 1$ for which $U \in \Lambda(N, \lambda; R)$.  In other words, 
		 \begin{equation}  
		 {\tt{AdFinLac}} := \bigcup \left\{ \Lambda(N, \lambda; R) : \lambda \in (0,1), N, R \geq 1\right\}. \label{defn: AdFinLac}  \end{equation}   
		%lacunary of order at most $N_2$ with lacunarity constant $\leq \lambda$.  
		If $U \not\in {\tt{AdFinLac}}$, we call it {\em{sublacunary}}, and denote it as $U \in {\tt{SubLac}}$.   
		\vskip0.1in
		\item Suppose next that $U \subseteq \mathbb R$ is unbounded. As shown in \eqref{Omega bounded justification} and
		\eqref{Omega bounded justification}, the slope set $U$ can then be written as 
		\[ U \subseteq \bigcup \bigl[U \cap [0, \infty) \bigr] \cup \bigl[U \cap (- \infty, 0] \bigr] \subseteq \bigcup_{n=0}^{7} U_n, \]
		where the slope set $U_n$ is of the form 
		\[ U_n = \left\{\tan \theta : \theta \in \frac{n \pi}{4} + \Theta_n \right\}, \quad 0 \leq n \leq 7. \] 
		We define 
		\[ \widetilde{U}_n := \left\{ \tan \theta : \theta \in \Theta_n \right\} \subseteq [0, 1]. \]  
		 For unbounded $U$, say that $U \in {\tt{AdFinLac}}$ if each normalized piece $\widetilde{U}_n$, is admissible finite-order lacunary in the sense of part \ref{AdFinLac bounded sets} of this definition; otherwise $U$ is sublacunary.
		\end{enumerate}
	\end{definition} 
	\noindent The covering constant $R$ controls how many lacunary components of maximal order are permitted; the definition of {\tt{AdFinLac}} ensures that admissible lacunarity is stable under finite unions.
	
	\section{$L^p$ bounds on $D_{\Omega}$ for $\Omega \in \Lambda(N, \lambda)$} \label{section: Lp boundedness} 
A central analytic consequence of admissible finite-order lacunarity, introduced in Section \ref{basic definitions section}, is the following $L^p$-boundedness result (Theorem \ref{THM: DB}). It says that admissible finite-order lacunarity of the slope set $\Omega$ leads to $L^p$-bounds for the directional maximal operators $D_{\Omega}$ and $M_{\Omega}$, with quantitative norm dependence only on the lacunarity order $N$, the lacunarity constant $\lambda$, and the Lebesgue exponent $p$. Related boundedness results have appeared in the work of Alfonseca \cite{Alfonseca}, Sj\"ogren and Sj\"olin \cite{SjogrenSjolin}, and Bateman \cite{Bateman}, though these do not explicitly address the notion of finite-order lacunarity introduced in the previous section. There are several notable differences in approach. 
\vskip0.1in 
\begin{itemize} 
\item First, \cite{{Alfonseca}, {SjogrenSjolin}} use a set of angles rather than a set of slopes. 
\vskip0.1in
\item Second, the proof in \cite{Bateman} relies on a tree structure of the slope set rather than the geometric property of {\tt{AdFinLac}} described above. 
\vskip0.1in
\item Most importantly, the finite-order lacunarity hypotheses of \cite{{Alfonseca}, {Bateman}, {SjogrenSjolin}} are formulated differently and their equivalence with Definition \ref{defn: Admissible finite order lacunarity} is not established in the literature. 
\vskip0.1in
\end{itemize} 
In view of this, we state and prove the result in the form needed for Theorem \ref{thm:main}.   
%Ultimately, these different approaches are equivalent. However, this equivalence is a consequence of the corrected lacunarity framework developed here, rather than an existing result in the literature.
%However, \cite{Alfonseca} does not explicitly provide a definition of finite-order lacunarity, although it introduces the notion of sets whose constituents lie between lacunary separators. On the other hand,  Definition \ref{defn: Admissible finite order lacunarity} is somewhat different from the one used in \cite{Bateman}.
 To keep the exposition self-contained, 
\begin{theorem}\label{THM: DB} 
For every Lebesgue exponent $p \in (1, \infty)$, non-negative integer $N$ and lacunarity constant $\lambda \in (0,1)$, there exists a positive constant $C(p, N, \lambda)$ depending only on these quantities such that 
\begin{equation} \label{D-Omega-bounded}
||D_{\Omega}||_{p \rightarrow p} \leq C(p, N, \lambda) \; \; \text{ for all } \Omega \in \Lambda(N, \lambda),   
\end{equation}
with $\Lambda(N, \lambda)$ as in Definition \ref{defn: Lacunary sets}.
\vskip0.1in  
\noindent Suppose that $\Omega \subseteq [0,1]$, $\Omega \in {\tt{AdFinLac}}$, in the sense of Definition \ref{defn: Admissible finite order lacunarity}, i.e., $\Omega \in \Lambda(N, \lambda; R)$ for some choice of $N, R \geq 1$ and $\lambda \in (0,1)$.  Then for all $p \in (1, \infty)$, 
\begin{equation} \label{just finiteness of op norm}  
\max \left[ ||D_{\Omega}||_{p \rightarrow p},  \; ||M_{\Omega}||_{p \rightarrow p}  \right] \leq R C(p, N, \lambda) < \infty. 
\end{equation} 
\noindent As a result, admissible finite-order lacunarity of $\Omega$ implies $L^p$-boundedness of $D_{\Omega}$ and $M_{\Omega}$; this corresponds to the implication \eqref{condition 2: operators unbdd} $\implies$ \eqref{condition 3: slopes sublacunary}  in Theorem \ref{thm:main}.   
\end{theorem} 
%	\noindent In particular, a slope set 
%	$\Omega \subseteq [0,1]$ is called admissible lacunary of finite order if it satisfies Definition \ref{defn: Admissible finite order lacunarity}. A sublacunary slope set fails Definition \ref{defn: Admissible finite order lacunarity} for every choice of parameters $N_1, N_2, \lambda$. 
%We will use this admissible notion of slopes throughout; Theorem \ref{thm:main} asserts that this condition is precisely equivalent to boundedness of the directional maximal operators $D_{\Omega}$ and $M_{\Omega}$. We will prove one implication of the theorem in the next section. 
\vskip0.1in
\noindent The proof of Theorem \ref{THM: DB} is deferred to Chapter \ref{positive direction chapter}, so as to not interrupt the flow of the main ideas. The reader may proceed there directly for the proof, which does not rely on any new material beyond what has been presented up until this point. They may also continue to Chapters \ref{trees-section} -- \ref{chapter: far from root} for the proof of the implication ``\eqref{condition 3: slopes sublacunary} $\implies$ \eqref{condition 1: Kakeya-type sets}'' of Theorem \ref{thm:main}, the other main contribution of this article. 
\vskip0.1in
\noindent The remainder of this chapter is for a more in-depth exploration of finite-order lacunarity. Section \ref{EXAMPLES SECTION} provides some examples and non-examples to anchor the concept introduced in Section \ref{basic definitions section} and points out key differences with prior work. Section \ref{SECTION: LACUNARY PROPERTIES} records some basic properties of  sets in {\tt{AdFinLac}} that will be helpful in the proof of Theorem \ref{THM: DB}. 

\section{Examples} \label{EXAMPLES SECTION}
Definitions \ref{defn: Lacunary sets} and \ref{defn: Admissible finite order lacunarity} of admissible finite-order lacunarity, presented in Section \ref{basic definitions section} may appear somewhat abstract at initial glance. 	In this section, we present several examples to illustrate these notions in practice, 
	highlighting both finite-order lacunary and sublacunary sets. 
	These examples serve to clarify the scope of the theory and emphasize the need for the more 
	general framework of finite-order lacunarity. Proofs associated with this example base are relegated to Appendix B,  Chapter \ref{1d lacunary examples section}. 
\vskip0.1in
\noindent The examples in this section fall into three categories, serving complementary roles: 
\begin{itemize} 
\item Canonical examples of admissible finite-order lacunarity illustrating the defining criteria are presented in Sections \ref{subsection: lacunary examples 1} and \ref{subsection: lacunary examples 2}. 
\vskip0.1in 
\item Examples of sublacunary sets identifying mechanisms of failure of the definition are highlighted in Sections \ref{subsection: sublacunary examples 1}, \ref{HRS-application-proof-section} and \ref{subsection: sublacunary examples 2}.  
\vskip0.1in 
\item Counterexamples emphasizing the necessity of uniformity and gap structure in Definitions \ref{defn: Lacunary sets} and \ref{defn: Admissible finite order lacunarity} and limitations of earlier definitions appear in Section \ref{section: comparison with Bateman lacunarity}.
\end{itemize} 
%	\section{Iterative constructions of finite-order lacunary sets}
	
	\subsection{Lacunary sets versus lacunary sequences} \label{subsection: lacunary examples 1} 
	\noindent The prototypical example of a lacunary set of order 1 and lacunarity constant $\lambda \in (0,1)$ is 
	\begin{equation} \label{std-lac-seq} 
	U = \{\lambda^{j} : j \geq 1\}, 
	\end{equation}  or any nontrivial subsequence thereof. Indeed $U$ is itself a lacunary sequence, and hence its own special sequence. However, it is not necessarily true that every $\Lambda(1, \lambda)$ set coincides with a lacunary sequence. The next lemma clarifies the precise relationship between these two notions. 
	%even though it has to be contained in the union of at most two such sequences, as proved in Lemma \ref{Lemma: Lacunary 1}) \eqref{Lemma: Lacunary 1 (c)}.  
	\begin{lemma} \label{Lemma: Lacunary 1}
		Lacunary sequences and sets in $\Lambda(1, \lambda)$ are not identical, but closely related; indeed, the former may be viewed as representative components of the latter, in the following sense.  
		\vskip0.1in 
		\begin{enumerate}[(a)]
		\item \label{Lemma: Lacunary 1 (a)}If $A \in {\tt{MonLac}}(\lambda)$, then $A \in \Lambda(1, \lambda)$. 
		\vskip0.1in
		\item \label{Lemma: Lacunary 1 (b)} An arbitrary lacunary sequence $A$ obeying \eqref{lacunarity constant} is contained in the union of at most two monotone lacunary subsequences, one on each side of the limit, each with lacunarity constant at most $\lambda$.
		\vskip0.1in
		\item \label{Lemma: Lacunary 1 (c)} Every bounded set $U \in \Lambda (1, \lambda)$ is contained in the union of at most two monotone, lacunary sequences $B_1$ and $B_2$, each with lacunarity constant $\leq \lambda$. Both $B_1$ and $B_2$ converge to the same limit point $\alpha$ as the special sequence $A$ of $U$, in the sense of Definition \ref{defn: Lacunary sets}. 
		\end{enumerate}
		\end{lemma} 
	\noindent The following example shows that the class $\Lambda(1, \lambda)$ 
	is strictly larger than the class of lacunary sequences. 
	This illustrates why the more general definition is needed.
	\begin{lemma} \label{Lemma: Lac set vs seq}
		The monotone decreasing sequence 
		\begin{equation}  U = \{2^{-2j} \pm 4^{-2j} : j \geq 1\} \in \Lambda \bigl(1, {1}/{2} \bigr), \label{lac set but not seq} \end{equation}  i.e. it is lacunary set of order 1 in the sense of Definition \ref{defn: Lacunary sets}, with special sequence $A = \{ 2^{-j} : j \geq 1\}$. It is not, however, a lacunary sequence. 
		\end{lemma}
		\noindent {\em{Remarks: }} 
		\begin{enumerate}[1.]
		\item  Lemmas \ref{Lemma: Lacunary 1} and \ref{Lemma: Lac set vs seq} are proved in Section \ref{sets vs sequences proof}. 
		\vskip0.1in
\item		In short, although $\Lambda(1, \lambda)$ is strictly larger than the class of lacunary sequences, its elements are still governed by at most two monotone geometric backbones converging to a common limit.
	\end{enumerate} 	
\subsection{Lacunary sets of arbitrarily large finite order} \label{subsection: lacunary examples 2}
\noindent Sets like \eqref{std-lac-seq} or \eqref{lac set but not seq}  extend naturally to  constructions exhibiting arbitrarily deep hierarchical clustering. The examples below show that lacunarity of order $N$ is stronger than that of order $(N-1)$, and cannot be reduced to a lower-order description. This underscores the hierarchy of admissible finite-order lacunarity. We begin with a canonical construction.
	\begin{lemma} \label{Lemma: iterated sums 1}
		Given any integer $N \geq 1$ and an ordered sequence of constants $\mathbf M = (M_1, \ldots, M_N)$ obeying  $N < M_1 \leq M_2 \leq \cdots \leq M_N$, the set 
		\begin{equation} \overline{U}_{N}(\mathbf M) = \Bigl\{ \sum_{r=1}^{N} M_r^{-j_r} : 1 \leq j_1 \leq j_2 \leq \cdots \leq j_N \Bigr\} \in \Lambda(N, M_1^{-1}), \label{U-order-N} \end{equation} i.e., is lacunary of order $N$ with lacunarity constant $\leq M_1^{-1}$ in the sense of Definition \ref{defn: Lacunary sets}. The special sequence of $\overline{U}_N(\mathbf M)$ can be chosen to be $A = \{ M_1^{-j} : j \geq 1\}$.
		\end{lemma}
		\noindent  This construction realizes a genuine $N$-level hierarchy, where each successive index $j_r$ introduces clustering at a smaller scale, confined within gaps determined by the previous levels. We next show that this structure persists under greater flexibility.
		\vskip0.1in
\noindent In Lemma \ref{Lemma: iterated sums 1}, bases $M_r$ and indices $j_r$ are ordered to enforce a nested structure. The next construction allows iterated sums with independent variation of the parameters, showing that such sets still lie within the framework of 
			admissible finite-order lacunarity. 
			This demonstrates the flexibility of the definition in capturing a wide class of examples.
			\begin{lemma} \label{Lemma: iterated sums 2}
				Given an integer $N \geq 1$ and constants $\mathbf M = \{M_1, \ldots, M_N\} \subseteq  (1, \infty)$,  there exist an integer $R_N$ and a constant  $\lambda_N \in (0, 1)$, depending only on $N, \mathbf M$, such that the set 
				\begin{equation} \label{Example: U Lac N} 
					U^{\ast}_{N} = U^{\ast}_N(\mathbf M) := \Bigl\{ \sum_{r=1}^{N} M_r^{-j_r} : \mathbf j = (j_1, \ldots, j_N) \in \mathbb N^{N} \Bigr\}  \end{equation}
				lies in the union of at most $R_N$ sets of type $\Lambda(N, \lambda_N)$ described in Definition \ref{defn: Lacunary sets}. Thus, $U^{\ast}_N(\mathbf M)$ is admissible lacunary of order $N$ in the sense of Definition \ref{defn: Admissible finite order lacunarity}.  
				\end{lemma}
				\vskip0.1in 
				 \noindent Lemmas \ref{Lemma: iterated sums 1} and \ref{Lemma: iterated sums 2} are proved in Section \ref{Lambda proofs}. Heuristically, they show that admissible finite-order lacunarity accommodates structured iterated clustering across multiple scales, provided that both the depth of iteration and the geometric decay remain uniformly controlled.
				 \vskip0.1in
			\noindent Finally, we illustrate that large lacunarity order can arise even in finite sets; in other words, one can construct lacunary sets of large order through finite approximations of infinite sublacunary sets. We furnish an example. 
			\begin{lemma} \label{lemma: dyadic rationals}
			\begin{enumerate}[(a)] 
			\item \label{lemma: dyadic rationals parta} Any non-empty finite set $\mathbb F \subseteq \mathbb R$ has the property that 
			\begin{equation} 
			\mathbb F \in \Lambda \left( \#(\mathbb F) -1 , \frac{1}{2}\right). 
			\end{equation} 
			 \item \label{lemma: dyadic rationals partb} For each integer $m \geq 1$, let us consider the collection of dyadic rationals of the form 
			\begin{equation} \label{dyadic Q_m} 
			\mathbb Q_m := \left\{\frac{r}{2^m} : 1 \leq r \leq 2^m  \right\}. 
			\end{equation}  
			Then $q_m \rightarrow \infty$, where $q_m$ is the smallest lacunarity order of $\mathbb Q_m$ with respect to the lacunarity constant $\frac12$, in the sense that 
			\begin{equation} 
			q_m := \min \left\{k \geq 1 : \mathbb Q_m \in \Lambda\Bigl(k, \frac{1}{2} \Bigr) \right\}. 
			\end{equation} 
			The integer $q_m$ is well-defined and finite, by part \eqref{lemma: dyadic rationals parta}. 
			\end{enumerate} 
			\end{lemma} 
			\noindent Lemma \ref{lemma: dyadic rationals} is proved in Section \ref{proof: dyadic rationals lemma}. Any finite set is, of course, a finite union of singletons, i.e., $\Lambda(0, \lambda)$ sets and therefore in {\tt{AdFinLac}}, but this example shows that even finite sets can exhibit arbitrarily high lacunarity order if the covering number $R$ is kept under control. The lacunarity order detects the number of distinct scales present in such a set.
			
\subsection{Sublacunary sets} \label{subsection: sublacunary examples 1} 
We next  examine examples of sublacunary sets. In contrast to the hierarchical clustering of finite-order lacunary sets, sublacunary sets exhibit a form of combinatorial thickness, where clustering cannot be organized into finitely many controlled scales. 
%For contrast, we next examine examples of sublacunary sets, which fall outside this framework and give rise to Kakeya-type behavior. 
The first result states that the class {\tt{AdFinLac}}, defined as in \eqref{defn: AdFinLac}, is stable under topological closure. This means that taking the closure of a finite-order lacunary set cannot generate new accumulation structure beyond that already present. This can be used to generate many sublacunary examples. A guiding principle in this section is that sublacunarity is driven by uncontrolled accumulation: any mechanism that significantly enlarges the accumulation set must destroy finite-order lacunarity.
\begin{lemma} \label{lemma: topological closure} 
For every $N \geq 1$ there is a constant $C_N \geq 1$ such that for every $\lambda \in (0, 1)$, 
\begin{equation} \label{induction: topological closure} 
\text{if $U \in \Lambda(N,\lambda)$, then $\bar{U}$ is a union of $C_N$ closed sets in $\Lambda(N, \lambda)$.} 
\end{equation} 
In particular, a set $U \in $ {\tt{AdFinLac}}  if and only if its closure $\bar{U} \in $ {\tt{AdFinLac}}. Equivalently stated, a set is sublacunary if and only if its closure is.   
\end{lemma} 
\noindent Since {\tt{AdFinLac}} sets are necessarily countable, Lemma \ref{lemma: topological closure} shows that a set whose closure introduces uncountably many accumulation points must be sublacunary.
%Sets in ${\tt{AdFinLac}}$ are countable by definition, hence Lemma \ref{lemma: topological closure} offers many examples of sublacunary sets. 		
\begin{corollary} \label{corollary: uncountable sublacunary}
The following conclusions hold for all sets $U \subseteq \mathbb R$. \label{corollary: sublacunary}
\begin{enumerate}[(a)]
\item A set $U$ that is dense in an interval is sublacunary according to Definition \ref{defn: Admissible finite order lacunarity}.  
\item A set $U$ whose topological closure $\bar{U}$ is uncountable is sublacunary. \label{corollary: sublacunary (b)}
\end{enumerate}				
\end{corollary}
\noindent {\em{Remarks: }}
\begin{enumerate} [1.]
\item Lemma \ref{lemma: topological closure} and Corollary \ref{corollary: uncountable sublacunary} are proved in Section \ref{section: sublacunary set proofs}. 
\vskip0.1in 
\item In light of the discussion in item \ref{negative item} of page \pageref{negative item}, the examples produced by Corollary \ref{corollary: uncountable sublacunary} foreshadow the role of sublacunarity in generating Kakeya-type configurations, where combinatorial richness replaces structured sparsity.
\end{enumerate} 
\subsection{Example of Hagelstein, Radillo-Murguia and Stokolos} \label{HRS-application-proof-section} 
	The set $\Omega_{\text{HRS}}$ in \eqref{HRS-example} of Section \ref{HRS-application-section} is an example of a slope set, originating in  \cite{HRS2024b}, without the Euclidean separation condition \eqref{Bateman-separation}. It consists of points in $[0,1]$ with special arrangements of binary digit blocks. The set $\Omega_{\text{HRS}}$ is itself countably infinite, being a subset of the dyadic rationals, but its closure is not; it generates an uncountable set of accumulation points through infinite concatenations of admissible digit blocks. As a consequence of Corollary \ref{corollary: uncountable sublacunary} and Theorem \ref{thm:main}, we are able to deduce  Corollary \ref{Corollary: HRS-example}		
					\begin{proof}[Proof of Corollary \ref{Corollary: HRS-example}] 
						We observe that the topological closure of $\Omega_{\text{HRS}}$ is uncountable, since 
\[\bar{\Omega}_{\text{HRS}} \supseteq  \left\{ \sum_{k=1}^{\infty} \frac{\varepsilon_k}{2^k} \; : \;  
				\begin{aligned} &{\pmb{\varepsilon}} = (\varepsilon_1, \varepsilon_2, \ldots ) = \bigl({\pmb{\kappa}}_1, {\pmb{\kappa}}_2, \ldots, {\pmb{\kappa}}_j,\ldots \bigr), \\  &{\pmb{\kappa}}_j = \text{ either } {\pmb{\eta}}_j \text{ or } {\pmb{\zeta}}_j \in \{0, 1\}^{N_j}\text{ for each } j \geq 1 \end{aligned} 
				\right\}. \] 
				The set on the right hand side allows infinitely many block choices, producing a Cantor-type structure of accumulation points.
Therefore, by Corollary \ref{corollary: sublacunary} \eqref{corollary: sublacunary (b)}, the set $\Omega_{\text{HRS}}$ is sublacunary. The conclusions about the directional maximal operators follow, of course, from Theorem \ref{thm:main}. 
						\end{proof}
\noindent Thus, the failure of Euclidean separation in $\Omega_{\text{HRS}}$ reflects a deeper combinatorial thickness detected by Definition \ref{defn: Admissible finite order lacunarity} as sublacunarity.

\subsection{Non-closure of finite-order lacunarity under algebraic sums} \label{subsection: sublacunary examples 2} 
Examples of finite-order lacunary sets of the type \eqref{U-order-N} or \eqref{Example: U Lac N} are built through controlled hierarchical sums. This gives rise to a natural question: 
\[\text{{\em{Is ${\tt{AdFinLac}}$ closed under algebraic sums?}}} \] 
We provide a counter-example in this section, based on \cite{KrocThesis}, that answers this question in the negative.  Let  $N_j \nearrow \infty$ be a fast-growing sequence, and $M_j = 2^{m_j}$ a slower growing one, so that 
\begin{equation} \label{growth condition Mj Nj}
M_j <  N_j - N_{j-1}; \quad \text{ for example, } N_j = 2^{j^2} \text{ and } M_j = 2^j \text{ will do. } 
\end{equation} 
For $j \geq 1$, and $1 \leq k \leq M_j = 2^{m_j}$ we  set 
\begin{align} 
&q_{jk} = 2^{-N_j} (1 + k2^{-m_j}), \; \text{ and define }  U := \bigcup_{j=1}^{\infty} U_j, \text{ where } \label{U+V counterexample1} \\
&U_j := \bigl\{2^{-N_j + k} + q_{jk} : 0 \leq k < M_j \bigr\}; \quad V := \left\{ - 2^{-j} : j \geq 1 \right\}.   \label{U+V counterexample2} 
\end{align} 
The construction is designed so that the dominant scales of $U$, namely  $\bigl\{ 2^{-N_j+k} : 1 \leq k \leq M_j, j \geq 1 \bigr\}$ are distinct members of a lacunary sequence $\{2^{-r} : r \geq 1\}$. Once $V$ is added, the dominant lacunary terms cancel, and the residual fine-scale perturbations, previously negligible, i.e.
\[ \{q_{jk} : 1 \leq k \leq M_j \} \subseteq U_j+V, \]  become the leading contribution. The lacunarity order of $U_j+V$ therefore increases in an uncontrolled manner with $j$,  resulting in the eventual sublacunarity for $U+V$.
\begin{lemma} \label{lemma: nonclosure under algebraic sums}
For sets $U, V$ as in \eqref{U+V counterexample1}, \eqref{U+V counterexample2}, the following conclusions hold: 
\begin{equation} 
U, V \in \Lambda\Bigl(1, \frac{1}{2}\Bigr), \; \text{ but } \; U + V \text{ is sublacunary.}  
\end{equation}  
\end{lemma} 
\noindent The proof of Lemma \ref{lemma: nonclosure under algebraic sums} is given in Section \ref{section: nonclosure under algebraic sums proof}. It shows that sublacunarity may arise not only from large accumulation in the closure, but also from algebraic interactions that disrupt hierarchical clustering. 

\subsection{Comparison with Bateman's notion of finite-order lacunarity} \label{section: comparison with Bateman lacunarity} 
We conclude this section with a discussion of our definition of admissible finite-order lacunarity with other definitions in the literature. In \cite{Bateman}, finite-order lacunarity was defined in the following way: 
\begin{equation} \label{Bateman-fol}
\begin{minipage}{4in} {\em{``Let us say that a set $\Omega$ is
		lacunary of order $N$ if it is covered by the union of a lacunary sequence $L$ of order
		$N-1$ with lacunary sequences converging to every point of $L$."}} 
		\end{minipage}
		\end{equation} 
		\vskip0.1in
		\noindent There are two subtle points in this definition that require precise quantification: 
		\begin{itemize} 
		\item {\em{Uniform control on the lacunarity constants: }} The lacunary constant associated with all the sequences involved in the definition should be uniformly bounded by a single constant $\lambda \in (0,1)$. Without uniform control on the lacunarity constant, the resulting set may be sublacunary. A counter-example appears in Lemma \ref{lemma: Bateman counterexample 1}.
		\vskip0.1in   
		\item {\em{Geometric separation of lower-order lacunary blocks: }} Even assuming the above uniformity, the positioning of the lacunary sequences relative to $L$ is equally important. Without additional assumptions separating the physical separation of these lacunary sequences, unintentional clustering may occur elsewhere, i.e. away from the lacunary limits. Such accumulation may lead to sublacunarity, as shown in Lemma \ref{lemma: Bateman counterexample 2}.     
		\end{itemize} 
		\begin{lemma} \label{lemma: Bateman counterexample 1}  
		The set defined by 
		\begin{equation}  \label{example: uncontrolled lacunary blocks}
		U  = \{ u_{jk} : j, k \geq 1\}, \quad u_{jk} = 2^{-j} + 2^{-j} \Bigl(1 - \frac{1}{j}\Bigr)^k \end{equation} 
		satisfies the definition \eqref{Bateman-fol} of finite-order lacunarity of order 2, but is sublacunary according to Definition \ref{defn: Admissible finite order lacunarity}. 
		\end{lemma} 
		\noindent Lemma \ref{lemma: Bateman counterexample 1} is proved in Section \ref{section: Bateman counterexample 1 proof}. The critical point here is that the lacunary sequence $\{u_{jk} : k \geq 1\}$ has lacunarity constant $1 - 1/j$, which approaches 1 as $j \rightarrow \infty$. Each component sequence is individually lacunary, but the degeneration of the lacunarity constant destroys uniform control, allowing accumulation at arbitrarily slow rates.
		\begin{lemma}\label{lemma: Bateman counterexample 2}  
			Let $\{q_j : j \geq 0 \}$ be an enumeration of the rationals in $[\frac{9}{10}, 1]$. Then the set 
			\[ U = \{u_{jk} : j, k \geq 0\} \quad u_{jk} = 2^{-j} + (q_j - 2^{-j})3^{-k} \]
			 satisfies the definition \eqref{Bateman-fol} of finite-order lacunarity of order 2, but is sublacunary according to Definition \ref{defn: Admissible finite order lacunarity}. 
			\end{lemma}
			\noindent The proof of Lemma \ref{lemma: Bateman counterexample 2} is given in Section \ref{section: Bateman counterexample 2 proof}. Here the lacunary components originate from a dense set of base points, allowing accumulation away from the designated lacunary limits. This example demonstrates that geometric placement of lower-order components within lacunary gaps is equally essential.

\section{Properties of finite-order lacunary sets} \label{SECTION: LACUNARY PROPERTIES}
This section records a few fundamental properties of admissible finite-order lacunarity that will be used in subsequent chapters. These properties describe how admissible finite-order lacunary sets behave under inclusion, transformations, and changes of parametrization, and will serve as useful tools in later arguments. The lengthier proofs are relegated to Appendix C, Chapter \ref{section: lacunary properties proofs}. 
	\subsection{Nesting of lacunary classes of finite order} We start with the monotonicity of the class $\Lambda(N, \lambda)$, which is a direct consequence of its definition.  
\begin{lemma} Finite-order lacunary sets are nested. Increasing the lacunarity order or the covering number, or relaxing the lacunarity constant enlarges the admissible class.\label{Lemma : lacunarity monotonicity} 
\begin{enumerate}[(a)]
\item For integers $N_i \geq 0$, lacunarity constants $\lambda_i \in (0,1)$ and covering numbers $R_i$, $i=1,2$,  \label{lacunary monotonicity part a}
\begin{equation} \label{lacunarity-monotonicity} 
\Lambda(N_1, \lambda_1; R_1) \subseteq \Lambda(N_2, \lambda_2; R_2) \; \text{ if } \; N_1 \leq N_2, \, \lambda_1 \leq \lambda_2, \, R_1 \leq R_2.
\end{equation} 
\item For a fixed $\lambda \in (0,1)$, we have a finer nesting property, as follows: \label{lacunary monotonicity part b}
\begin{equation} \label{nesting} {\tt{MonLac}}(\lambda) \subsetneq \Lambda(1, \lambda), \quad \Lambda(N, \lambda)\subseteq \bigcup_{R=1}^{\infty} \Lambda(N, \lambda; R) \; \text{ for all } N \geq 1. \end{equation} 
\end{enumerate} 
\end{lemma} 
\begin{proof} 
For part \eqref{lacunary monotonicity part a}, we ask the reader to verify, from Definitions \ref{defn: Lacunary sets} and \ref{defn: Admissible finite order lacunarity},  that 
\[ \Lambda(N_1, \lambda_1; R_1) \subseteq \Lambda(N_2, \lambda_1; R_1) \subseteq \Lambda(N_2, \lambda_2; R_1)  \subseteq \Lambda(N_2, \lambda_2; R_2) \] for $(N_1, \lambda_1, R_1)$ and $(N_2, \lambda_2, R_2)$ as in \eqref{lacunarity-monotonicity}.
\vskip0.1in
\noindent For part \eqref{lacunary monotonicity part b}, the first inclusion in \eqref{nesting} follows from Lemma \ref{Lemma: Lac set vs seq}. The second inclusion is a consequence of the identity $\Lambda(N, \lambda) = \Lambda(N, \lambda; 1)$. 
\end{proof}
\subsection{Non-closure of  {\tt{AdFinLac}} under countable unions}
An immediate corollary of Lemma \ref{Lemma : lacunarity monotonicity} is that the class {\tt{AdFinLac}} defined by \eqref{defn: AdFinLac} is closed under finite unions. Specifically, 
\begin{equation} 
\left\{
\begin{aligned} 
&\text{if } U \subseteq \bigcup_{i=1}^{K} U_i, \quad \text{ where } U_i \in \Lambda(N_i, \lambda_i; R_i) \text{ for } 1 \leq i \leq K, \text{then }\\ 
& U \in \Lambda(N, \lambda; R), \text{ with } N = \max_{1\leq i \leq K} N_i, \; \lambda = \max_{1 \leq i \leq K} \lambda_i, \; R = \sum_{i=1}^{K} R_i.  
\end{aligned}
\right\}
\end{equation} 
This property, however, is not preserved under countable unions, even for fixed $N$ and $\lambda$; the next lemma provides the necessary examples. 
\begin{lemma} \label{lemma: nonclosure under countable union} 
For each $N \geq 0, \lambda \in (0,1)$, there exists a sublacunary set $U$ that is contained in a countable union of sets of the form $\Lambda(N, \lambda; \cdot)$, i.e.,
\begin{equation} \label{nonclosure} 
U  \in {\tt{SubLac}}, \quad \text{ but } \quad U \subseteq \bigcup_{R=1}^{\infty} \tilde{U}_R \text{ with } \tilde{U}_R \in \Lambda(N, \lambda; R). 
\end{equation}
\end{lemma} 
The proof of this lemma appears in Section \ref{nonclosure proof section}. 
\subsection{Choice of the special sequence} 
As we know from Definition \ref{defn: Admissible finite order lacunarity}, special lacunary sequences are crucial to the construction of finite-order lacunary sets. They determine the point of reference for the hierarchical structure of a lacunary set. It is natural to ask whether the choice of a special sequence is canonical. While the special point may be unique under certain circumstances, our next result establishes the generic non-uniqueness of special sequences.	
\begin{lemma} \label{lemma: special sequence choice} 
For a set $U \in \Lambda(N, \lambda)$, the accompanying special sequence 
\[A = \{a_j  : j \geq 1\} \in {\tt{MonLac}}(\lambda), \]  
as prescribed by Definition \ref{defn: Lacunary sets}, is not unique in general. In fact, for any $U \in \Lambda(N, \lambda)$, the special sequence $A$ can be chosen to obey the stronger condition 
\begin{equation} \label{lacunary above and below}
\lambda^2 |a_j-a| < |a_{j+1}-a| \leq \lambda |a_{j}-a| \; \; \text{ for all } j \geq1. 
\end{equation} 
\end{lemma} 
\noindent The proof of the lemma is included in Section \ref{section: proof of special sequence choice}. Thus, one may always choose a special sequence whose decay is quantitatively controlled from both sides; this will be useful in later arguments; see for example the proof of Theorem \ref{THM: DB}, Section \ref{section: decomposition of A}.

\subsection{Invariant operations for finite-order lacunary classes} For any fixed $N$ and $\lambda$, the class $\Lambda(N, \lambda)$ as described in Definition \ref{defn: Lacunary sets} is invariant under set inclusion and affine transformations; these properties, summarized in the following lemma, are easy to verify and left to the reader. 
	\begin{lemma}\label{lacunarity under linear operations}
		Let $U \in \Lambda(N, \lambda)$. Then 
\begin{align} 		
&V \in \Lambda(N,\lambda) \text{ for any } V \subseteq U, \text{ and } \nonumber \\   
&c_1U + c_2 \in \Lambda(N,\lambda) \text{ for any $c_1 \ne 0$, $c_2 \in \mathbb R$}. \label{linear-invariance}
		\end{align} 
	\end{lemma} 
	
	\subsection{Bi-Lipschitz invariance} \label{section: bi-Lipschitz invariance} More generally, admissible finite-order lacunarity is stable under bi-Lipschitz transformations. Let $F: I \rightarrow \mathbb R$ be a bi-Lipschitz function on a bounded interval $I$; i.e., there exists a constant ${\tt{L}}> 0$ such that   
		\begin{equation} \label{def: bi-Lipschitz}
			{\tt{L}}^{-1}|x-y| \leq |F(x) - F(y)| \leq {\tt{L}} |x-y| \quad \text{ for all $x, y \in I$.}
		\end{equation}
The defining condition \eqref{def: bi-Lipschitz} forces such maps to be strictly monotone. The next proposition posits that the lacunarity order is preserved under bi-Lipschitz maps. 
\begin{proposition}
		The class of sets {\tt{AdFinLac}} (Definition \ref{defn: Admissible finite order lacunarity}) is invariant under bi-Lipschitz maps. 
		\vskip0.1in
		\noindent Precisely, choose any function $F$ as in \eqref{def: bi-Lipschitz}. Then for any $N \geq 1$ and $\lambda \in (0, 1)$, there exists a constant $C_N = C_N(\lambda, \mathtt L) \geq 1$ such that 
		\begin{equation}  U \in \Lambda(N, \lambda) \; \implies \; F(U) \in \Lambda(N, \lambda; C_N). \label{lacunary-induction} \end{equation} 
		\label{prop: AdFinLac under F}
		\end{proposition}
		\noindent Thus, while the lacunarity order is preserved under a bi-Lipschitz transformation, the covering constant may increase in a controlled manner. Proposition \ref{prop: AdFinLac under F} is proved in Section \ref{section: bi-Lipschitz maps proof}. 
		
		\subsection{Equivalence of finite-order lacunarity for slopes and angles} \label{section: slopes and angles} Since the map 
		\[ \theta \in [0, \frac{\pi}{4}] \mapsto \tan \theta \in [0,1]  {\text{ is bi-Lipschitz, }} \]
		a consequence of Proposition \ref{prop: AdFinLac under F} is that finite-order lacunarity of slopes and angles are equivalent;  Both approaches appear in the literature. For instance, \cite{Alfonseca} defines $D_{\Omega}, M_{\Omega}$ for 
		\begin{equation} \label{Theta and Omega}  
		\Omega = \{\tan \theta : \theta \in \Theta \} \subseteq [0,1], \quad \Theta \in \bigl[ 0, \frac{\pi}{4}\bigr],
		\end{equation}  and proves their $L^p$-boundedness in terms of the lacunarity order of the set of angles $\Theta$. On the other hand, \cite{Bateman} studies the unboundedness of these operators based on the sublacunarity of the slope set $\Omega$. The following corollary bridges the gap between the two approaches.  
	\begin{corollary}
	Suppose that $\Theta$ is a set of angles, $\Omega$ its corresponding set of slopes; in other words, they are related by the equation \eqref{Theta and Omega}. Then $\Theta \in {\tt{AdFinLac}}$ if and only if $\Omega \in {\tt{AdFinLac}}$.  
	\end{corollary}
\noindent This shows that the notion of admissible finite-order lacunarity is independent of whether one works with slopes or angles, justifying the interchangeability of these perspectives in the literature.
\vskip0.1in
\noindent In summary, the properties presented in this section show that admissible finite-order lacunarity defines a robust geometric class: it is stable under inclusion, flexible in its representation, invariant under natural transformations, and independent of the choice of parametrization.

\chapter{Rooted, labelled trees} \label{trees-section} 
The framework of rooted labelled trees in the context of directional maximal operators first appeared in the work of Bateman and Katz \cite{BatemanKatz}. It was subsequently developed by Bateman \cite{Bateman} and extended to a higher dimensional setting by two of the authors \cite{KrocPramanik}. The language of rooted, labelled trees continues to play a central role in our construction of Kakeya-type sets. It provides a combinatorial encoding of the set-theoretic architecture of finite-order lacunarity introduced in Chapter \ref{section: finite-order lacunarity}, translating questions about
clustering and separation in $\mathbb{R}$ into questions about branching
behaviour in a tree. 
\vskip0.1in
\noindent This perspective is central to the proof of the implication \eqref{condition 3: slopes sublacunary} $\implies$ \eqref{condition 1: Kakeya-type sets} in Theorem \ref{thm:main}. In particular, as we will see in Chapter \ref{section: split implies lacunarity}, the notion of splitting number introduced in \cite{Bateman} serves as a combinatorial proxy for lacunarity order: finite-order lacunarity corresponds to controlled branching along rays, while sublacunarity manifests as unbounded splitting behaviour.
\vskip0.1in
\noindent In Sections \ref{trees} and \ref{tree encoding section}, we recall the basic terminology of trees and describe a canonical representation of bounded subsets of $\mathbb{R}$ via $M$-adic trees. The interested reader is referred to \cite{LyonsPeres} for a comprehensive treatise on this material. This construction allows us to pass from sets of slopes to tree structures in a systematic way. A key structural invariant in this setting is the splitting number, which measures the branching complexity of a tree and will serve as a combinatorial proxy for lacunarity order. 
\vskip0.1in
\noindent Section \ref{splitting number examples} illustrates this correspondence through examples, recasting some of the sets from Section \ref{EXAMPLES SECTION} in the language of trees and computing their splitting numbers. 
\vskip0.1in
\noindent The main new contribution of this chapter appears in Sections \ref{section: compressed trees} and \ref{examples of compressed trees section}, where we introduce {\em{compressed trees}}. These provide a reduced representation of the $M$-adic tree by collapsing non-branching segments, while preserving the essential branching structure. This compression reflects the fact that, in later arguments,
only certain distinguished scales contribute to the combinatorial complexity of
the slope set. The compressed tree formalism will be crucial for the pruning
procedure in Chapter \ref{Chapter: Pruning of the slope tree} and for the construction of Kakeya-type configurations in subsequent chapters.

%Sections \ref{section: compressed trees} and \ref{examples of compressed trees section} contain the new contributions of this chapter; they introduce alternative representations of slope sets, called {\em{compressed trees}}. 
%In certain contexts, these provide a more efficient encoding of the underlying geometry and will be essential for the Kakeya-type constructions developed later in the paper.

\section{The terminology of trees}\label{trees}
\subsection{Vertex sets}
A connected, undirected graph is, by definition, a pair $(\mathcal V, \mathcal E)$, where the elements of $\mathcal V$ are called vertices, and the set $\mathcal E$ is a collection of unordered pairs of vertices. A member of $\mathcal E$ is called an edge.  A \textit{tree} is defined to be a connected undirected graph with no cycles; this means any two distinct vertices are connected by exactly one path.   
\vskip0.1in
\noindent For our purposes, trees will always be realized in a rooted, labelled form. A {\em{rooted, labelled}} tree is a connected, cycle-free graph equipped with a distinguished vertex, called the root, together with a labelling of its edges (or vertices) from a finite alphabet. Precisely, the vertex set of a rooted, labelled tree $\mathcal{T}$ is, by definition, a non-empty collection of finite sequences of non-negative integers of the form $v_n = \langle i_1, \ldots, i_n \rangle$, where  $n \geq 0$ is a running index, with the following property: 
\begin{equation} \label{defn: vertex} 
\begin{aligned}  
&\text{ if } v_n = \langle i_1,\ldots,i_n\rangle\in \mathcal{T},  \text{ then } v_k = \langle i_1,\ldots,i_k\rangle\in \mathcal{T} \text{ for all   } 0 \leq k \leq n; \\ 
&\text{ here } k = 0 \text{ corresponds to the empty sequence $\emptyset$, called the {\em{root}}.} 
\end{aligned} 
\end{equation} 
%\begin{enumerate}
%	\item[(i)] for any $k$, $0\leq k\leq n$, $\langle i_1,\ldots,i_k\rangle\in \mathcal{T}$, where $k=0$ corresponds to the empty sequence, and
%	\item[(ii)] for every $j\in \{0,1,\ldots,i_n\}$, we have $\langle i_1,\ldots,i_{n-1},j\rangle\in \mathcal{T}$.  
%\end{enumerate}
The integer $n$, which represents the length of the sequence $\langle i_1,\ldots,i_n\rangle$, is called the {\em{height}} or {\em{generation}} of the vertex $v_n = \langle i_1, \ldots, i_n \rangle$ in $\mathcal T$. We denote the height of a vertex by $h_{\mathcal T}( \cdot)$, omitting the suffix $\mathcal T$ if the ambient tree is clear from the context.
\begin{equation} \label{defn: height} 
h_{\mathcal T}(v_n) = h(v_n) = h(\langle i_1, \ldots, i_n \rangle) = n. 
\end{equation} 
The vertex $v_n$ is said to be of an $n^{\text{th}}$ generation vertex of $\mathcal T$.
The root $\emptyset$ is, by convention, a vertex of the $0^{\text{th}}$ generation. The {\em{height of the tree}} $\mathcal T$ is the supremum of the heights of its vertices. This may be finite or infinite. 
\subsection{Edge sets and parent-child relationship} \label{section: edges}
\noindent The edges of a rooted, labelled tree $\mathcal T$ are dictated by lineage, a by-product of the property \eqref{defn: vertex}. We say that 
\begin{equation} 
\begin{aligned} 
&{\text{$v_{n-1} = \langle i_1,\ldots, i_{n-1}\rangle$ is the \textit{parent} of $v_n = \langle i_1,\ldots,i_{n-1}, i_n\rangle$, and}} \\ 
&{\text{$v_n = \langle i_1,\ldots,i_{n-1}, i_n\rangle$ is a \textit{child} of $v_{n-1} = \langle i_1,\ldots,i_{n-1}\rangle$}}.  
\end{aligned} 
\end{equation} 
An {\em{edge}} $e$ of the tree $\mathcal T$ connects a parent with its child. Every non-root vertex of a tree has a unique parent, and is connected to the parent by an edge.  A {\em{directed edge}} is the {\em{ordered}} parent-child pair $(v_{n-1}, v_n)$. A vertex $v_n \in \mathcal T$ is said to be a {\em{terminal vertex of $\mathcal T$}} if it does not have a child.  
\vskip0.1in
\noindent A {\em{ray}} $\mathcal R$ of $\mathcal T$ is a sequence of connected, directed edges
\begin{equation}
\begin{aligned} 
&\mathcal R = (e_1, e_2, \ldots), \text{ where the child vertex of $e_i$ } \\
&\text{is the parent vertex of $e_{i+1}$ for all $i \geq 1$. }
\end{aligned} 
\end{equation} 
Thus every non-root vertex $v$ of $\mathcal T$ is joined to the root by a unique ray and the height $h(v)$ is the number of edges in this ray. The root has height $0$, and vertices at height $k$ form the $k$-th generation of the tree.
\vskip0.1in
\noindent A ray $\mathcal R$ of  $\mathcal T$ is of {\em{maximal length}} if it is not strictly contained in any other ray. Let $\partial \mathcal{T}$ denote the collection of all rays in $\mathcal{T}$ of maximal (possibly infinite) length. Rays in  $\partial \mathcal{T}$ correspond to infinite or maximal finite paths starting at the root. Lengths of rays in $\partial \mathcal{T}$ need not be uniform. 

\subsection{Tree ancestry} If $u$ and $v$ are two vertices in $\mathcal{T}$ that lie on a common ray, with $h(u) > h(v)$, then we say $u$ is a \textit{descendant} of $v$ (or that $v$ is an \textit{ancestor} of $u$), and we write $u\subset v$.  For the rooted labelled trees considered in this article, $u$ and $v$ will represent intervals (see Section \ref{tree encoding section}), so that the ancestry relation $u \subset v$ will coincide with actual set inclusion. The \textit{youngest common ancestor} of $u$ and $v$ in the tree $\mathcal T$, denoted by $D(u,v) = D_{\mathcal T}(u, v)$, is the vertex of maximal height common to all rays passing through both $u$ and $v$:
 \begin{equation} \label{defn: youngest common ancestor}
D_{\mathcal T}(u, v) = w \in \mathcal T \text{ where } h(w) = \max \left\{ h(w') : u \subset w', \, v \subset w', \, w' \in \mathcal T\right\}.
\end{equation} 
%The empty sequence $\emptyset$ is the designated \textit{root} of the tree $\mathcal{T}$ and all vertices of the form $\langle i_1 \rangle\in\mathcal{T}$ are children of this root.  
 % For a fixed vertex $v\in\mathcal{T}$, we define the \textit{subtree of $\mathcal{T}$ generated by the vertex} $v$ to be the maximal subtree of $\mathcal{T}$ with $v$ as the root.

%The \textit{height} of the tree is taken to be the supremum of the lengths of all the sequences in the tree.  Further, we define the height $h(\cdot)$ of a vertex to be the length of its identifying sequence.  If the height of a vertex $v$ is equal to $k$, we say that $v$ is a $k$\textit{th generation vertex} of the tree.  The height of the root is always taken to be zero. 
\subsection{Subtrees and truncated trees} 
A {\em{subtree}} $\mathcal T' \subseteq \mathcal T$ consists of a connected subset of the vertex set of $\mathcal T$, with the same ancestry relations, i.e., the edge structure induced by $\mathcal T$. For a fixed vertex $v \in\mathcal T$, the {\em{subtree of $\mathcal T$ generated by $v$}} consists of all vertices of $\mathcal T$ that descend from $v$; in other words, it is the maximal subtree of $\mathcal T$ with $v$ as its root. 
\vskip0.1in 
\noindent If $n\in\mathbb N = \{1, 2, \ldots\}$ and $\mathcal T$ is a tree of height at least $n$, the \textit{truncation} of $\mathcal{T}$ to height $n$ is the subtree of $\mathcal{T}$ consisting of all vertices of height at most $n$.  A tree is called \textit{locally finite} if its truncation to every level is finite; i.e. consists of finitely many vertices.  All our trees will have this property.  In the remainder of this article, when we speak of a \textit{tree} we will always mean a \textit{locally finite, rooted, labelled tree}, whose height may be finite or infinite depending on the context. 
\subsection{Sticky maps} \label{sticky maps section} 
\noindent  The following definition, concerning mappings between trees, will be important later.
\begin{definition}\label{D:stickiness}
	Let $\mathcal{T}$ and $\mathcal{T}'$ be two trees with equal (possibly infinite) heights.  A map $\sigma: \mathcal{T}\rightarrow \mathcal{T}'$ is called {\em sticky} if it preserves heights and lineages: 
\begin{enumerate}
	\item[$\bullet$] for all $v\in \mathcal{T}$, $h(v) = h(\sigma(v))$, and
	\item[$\bullet$] $u\subset v$ implies $\sigma(u)\subset\sigma(v)$ for all $u,v\in \mathcal{T}$.
\end{enumerate}
%We often say that $\sigma$ is sticky if it preserves heights and lineages.
\end{definition}
\noindent A one-to-one and onto sticky map between two trees, when it exists, is said to be an \textit{isomorphism} and the two trees are said to be \textit{isomorphic}.  Two isomorphic trees will be treated as essentially identical.
\subsection{The splitting number of a tree} \label{section: splitting number}
%There are many ways to quantify the ``size" or ``spread'' of a tree (see~\cite{LyonsPeres}). Of these, 
The concept of a \textit{splitting number} is used to measure the ``size" or ``spread" of a tree by quantifying its branching complexity. It was first identified in the work of Bateman \cite{Bateman} as a key tool for studying sets that are finite-order lacunary and directional maximal operators associated to them.  We recall its definition here. 
%and finite order lacunarity of slope sets was first  
%proved to be the most relevant in the planar characterization of directions that admit Kakeya-type sets~\cite{Bateman}. Not surprisingly, it will turn out to be equally important for us. One of its applications is the explicit restatement of finite order lacunarity of a set $\Omega$ in terms of the structure of the tree encoding $\Omega$.  We define the notion of splitting number below, then collect some fundamental results about this quantity that will allow us to prove Theorem~\ref{MainThm1}, which is also the first forward implication in Theorem~\ref{MainThm2}.
\vskip0.1in 
\noindent We say that a vertex $v$ \textit{splits in} $\mathcal{T}$ if it has at least two children in $\mathcal{T}$.  When the context is clear, we simply say that $v$ \textit{splits}, and call it a \textit{splitting vertex}.  Define split$_{\mathcal{T}}(\mathcal R)$, the \textit{splitting number of a ray} $\mathcal R$ \textit{in} $\mathcal{T}$ to be the number of splitting vertices in $\mathcal{T}$ along that ray.  Figure \ref{fig:splitting-ray} depicts a ray $\mathcal R \in \partial \mathcal T$ with splitting number 4.  
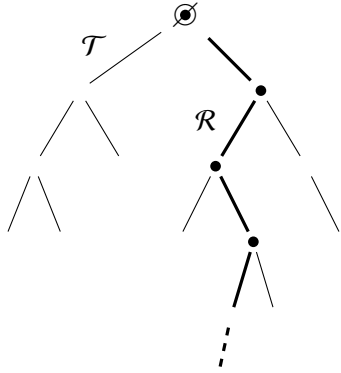
\begin{figure}[ht]
\centering
\begin{tikzpicture}[x=1cm,y=1cm,>=stealth]

% Main vertices
\node (r)   at (0,4) {$\emptyset$};

\node (a)   at (-1.4,3) {};
\node (b)   at ( 1.0,3) {};

\node (a1)  at (-2.0,2) {};
\node (a2)  at (-0.8,2) {};

\node (c)   at ( 0.4,2) {};
\node (d)   at ( 1.6,2) {};

\node (c1)  at (-0.1,1) {};
\node (c2)  at ( 0.9,1) {};

\node (e)   at ( 0.6,0) {};
\node (f)   at ( 1.2,0) {};

% Extra side branches
\node (d1)  at (2.1,1) {};
\node (a11) at (-2.4,1) {};
\node (a12) at (-1.6,1) {};

% Thin edges
\draw (r) -- (a);
\draw (a) -- (a1);
\draw (a) -- (a2);
\draw (a1) -- (a11);
\draw (a1) -- (a12);

\draw (b) -- (d);
\draw (d) -- (d1);

\draw (c) -- (c1);
\draw (c2) -- (f);

% Highlighted ray (connected path)
\draw[very thick] (r) -- (b);
\draw[very thick] (b) -- (c);
\draw[very thick] (c) -- (c2);
\draw[very thick] (c2) -- (e);

% Optional dashed continuation
\draw[very thick,dashed] (e) -- (0.45,-0.7);

% Splitting vertices along the ray
\fill (r)  circle (2pt);
\fill (b)  circle (2pt);
\fill (c)  circle (2pt);
\fill (c2) circle (2pt);

% Label for the ray
\node[right] at (0,2.6) {$\mathcal R$};
\node[right] at (-1.5,3.6) {$\mathcal T$};
\end{tikzpicture}
\caption{\small{A highlighted ray $\mathcal R$ in a rooted tree $\mathcal T$. The filled vertices are splitting vertices of $\mathcal T$ lying on $\mathcal R$, so the splitting number of the ray is $\mathrm{split}_{\mathcal T}(\mathcal R)=4$.}}
\label{fig:splitting-ray}
\end{figure}
\vskip0.1in
\noindent The \textit{splitting number of a vertex} $v$ \textit{with respect to a tree} $\mathcal{T}$ is defined to be 
\begin{equation}\label{splitting vertex}
\text{split}_{\mathcal{T}}(v) := \max_{\mathcal{S}_v\subseteq\mathcal{T}}\ \min_{\mathcal R_v\in \partial \mathcal{S}_v}\text{split}_{\mathcal{S}_v}(\mathcal R_v),
\end{equation}
where the maximum is taken over all subtrees $\mathcal{S}_v\subseteq\mathcal{T}$ rooted at $v$, and the minimum is taken over all rays $\mathcal R_v$ in $\mathcal{S}_v$ that originate at the vertex $v$. Thus the splitting number of $v \in \mathcal T$ is at least $m$ if there is a subtree $\mathcal S_v \subseteq \mathcal T$ rooted at $v$, every ray of which has at least $m$ splitting vertices in $\mathcal S_v$. See Figure \ref{fig:splitting-vertex} for a visual rendition of this property.  
\begin{figure}[ht]
\centering
\begin{tikzpicture}[x=1cm,y=1cm,>=stealth]

%------------------------------------------------
% Ambient tree T (dotted)
%------------------------------------------------

% root
\node (r) at (0,4) {};
\node[above right=-1pt] at (r) {$v$};

% level 1
\node (a) at (-1.7,3) {};
\node (b) at ( 1.3,3) {};

% level 2
\node (a1) at (-2.6,2) {};
\node (a2) at (-0.8,2) {};
\node (b1) at ( 0.6,2) {};
\node (b2) at ( 2.1,2) {};

% level 3
\node (a11) at (-3.0,1) {};
\node (a12) at (-2.2,1) {};
\node (a21) at (-1.2,1) {};
\node (a22) at (-0.4,1) {};
\node (b11) at ( 0.2,1) {};
\node (b12) at ( 1.0,1) {};
\node (b21) at ( 1.7,1) {};
\node (b22) at ( 2.5,1) {};

% level 4 (a few extra ambient branches)
\node (a121) at (-2.4,0) {};
\node (a122) at (-2.0,0) {};
\node (b111) at ( 0.0,0) {};
\node (b112) at ( 0.4,0) {};
\node (b221) at ( 2.3,0) {};
\node (b222) at ( 2.7,0) {};

% ambient tree edges
\draw[densely dotted] (r) -- (a);
\draw[densely dotted] (r) -- (b);

\draw[densely dotted] (a) -- (a1);
\draw[densely dotted] (a) -- (a2);
\draw[densely dotted] (b) -- (b1);
\draw[densely dotted] (b) -- (b2);

\draw[densely dotted] (a1) -- (a11);
\draw[densely dotted] (a1) -- (a12);
\draw[densely dotted] (a2) -- (a21);
\draw[densely dotted] (a2) -- (a22);
\draw[densely dotted] (b1) -- (b11);
\draw[densely dotted] (b1) -- (b12);
\draw[densely dotted] (b2) -- (b21);
\draw[densely dotted] (b2) -- (b22);

\draw[densely dotted] (a12) -- (a121);
\draw[densely dotted] (a12) -- (a122);
\draw[densely dotted] (b11) -- (b111);
\draw[densely dotted] (b11) -- (b112);
\draw[densely dotted] (b22) -- (b221);
\draw[densely dotted] (b22) -- (b222);

%------------------------------------------------
% Highlighted subtree S_v (thick, asymmetric)
%------------------------------------------------

% thick subtree edges
\draw[very thick] (r) -- (a);
\draw[very thick] (r) -- (b);

\draw[very thick] (a) -- (a1);
\draw[very thick] (a) -- (a2);

\draw[very thick] (a1) -- (a12);   % left branch splits again after one step
\draw[very thick] (a12) -- (a121);
\draw[very thick] (a12) -- (a122);

\draw[very thick] (a2) -- (a21);   % terminal branch

\draw[very thick] (b) -- (b1);     % right branch does not split immediately at next node
\draw[very thick] (b1) -- (b11);
\draw[very thick] (b1) -- (b12);

\draw[very thick] (b12) -- (b112); % one side continues
\draw[very thick] (b11) -- (b111); % other side continues

% splitting vertices in S_v
\fill (r)   circle (2pt);
\fill (a)   circle (2pt);
\fill (b)   circle (2pt);
\fill (a12) circle (2pt);
\fill (b1)  circle (2pt);

% subtree label
\node at (2.5,2.35) {$\mathcal S_v$};

\end{tikzpicture}
\caption{\small{A subtree $\mathcal S_v$ rooted at $v$, shown in thick lines inside an ambient tree $\mathcal T$ indicated by dotted lines. Every maximal ray in $\mathcal S_v$ passes through at least two splitting vertices. Hence $\text{split}_{\mathcal T}(v)\ge 2$.}}
\label{fig:splitting-vertex}
\end{figure}
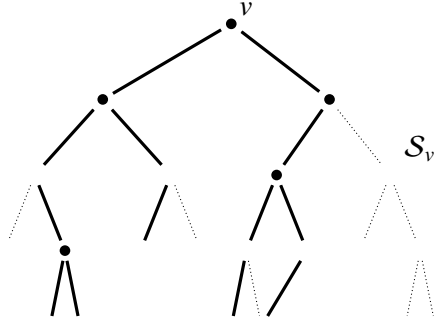
\vskip0.1in
\noindent Finally, the \textit{splitting number of the tree} $\mathcal{T}$ is defined as 
\begin{equation}\label{tree split}
\text{split}(\mathcal{T}) := \max_{v\in\mathcal{T}}\ \text{split}_{\mathcal{T}}(v).
\end{equation}
Two illustrative examples are worth noting. A tree consisting of a single ray has splitting number zero. For a full binary tree of height $N$, where each parent vertex has exactly two children, the splitting number is $N$.

\subsection{Preliminary facts about splitting numbers}
Let us take a moment to explore a few properties of splitting numbers that will be useful in the sequel. Our first result about splitting numbers (of vertices) says that they are monotone along lineages; the splitting number of an ancestor (a vertex of lesser height) is at least as large as that of one of its descendants. 
\begin{lemma}\label{monotonicity}
	Let $u,v\in\mathcal{T}$ with $u\subset v$.  Then 
	\begin{equation} \label{monotonicity: splitting number of vertices}  
	\text{split}_{\mathcal{T}}(u) \leq \text{split}_{\mathcal{T}}(v).
	\end{equation} 
	The inequality in \eqref{monotonicity: splitting number of vertices}  is strict if $v$ is a splitting vertex of $\mathcal T$, and has a descendant $u'$ obeying 
	\begin{equation}  \label{uu':hypotheses} 
	u \not\subset u', \quad u' \not\subset u, \quad \text{split}_{\mathcal T}(u') \geq \text{split}_{\mathcal T}(u). \end{equation} 
\end{lemma}
\begin{proof}
The idea is to turn a subtree rooted at $u$ to one rooted at $v$ by extending the rays upward.
Let $\mathcal{S}_u$ be a subtree of $\mathcal{T}$ rooted at $u$.  Define $\mathcal{S}_{v \rightarrow u}$ to be the union of the tree $\mathcal{S}_u$ with the path in $\mathcal{T}$ connecting $v$ to $u$.  This is a subtree of $\mathcal{T}$ rooted at $v$.  Since $v$ does not split in $\mathcal{S}_{v \rightarrow u}$ and there are no splitting vertices in $\mathcal{S}_{v \rightarrow u}$ between $v$ and $u$, we find that for any ray $\mathcal R$ in $\mathcal S_u$,  \begin{equation} \label{split equality} \text{split}_{\mathcal S_u}(\mathcal R) = \text{split}_{\mathcal S_{v \rightarrow u}}(\mathcal R_v), \end{equation} 
where $\mathcal R_v$ is the ray in $\mathcal S_{v \rightarrow u}$ rooted at $v$ obtained by extending $\mathcal R$ to $v$.  Conversely, if $\mathcal R_v$ is a ray in $\mathcal S_{v \rightarrow u}$, then \eqref{split equality} holds for $\mathcal R = \mathcal R_v \cap \mathcal S_u$.  Maximizing over all subtrees $ \mathcal{S}_u\subseteq\mathcal{T}$ rooted at $u$, we have that 
\begin{align}
	\text{split}_{\mathcal{T}}(u) &= \max_{\mathcal{S}_u\subseteq\mathcal{T}}\min_{\mathcal R\in \partial \mathcal{S}_u}\text{split}_{\mathcal S_{u}}(\mathcal R) \nonumber \\
	&= \max_{\mathcal{S}_{v \rightarrow u}\subseteq\mathcal{T}}\min_{\mathcal R_v\in \partial \mathcal{S}_{ v \rightarrow u}}\text{split}_{\mathcal S_{v \rightarrow u}}(\mathcal R_v) \nonumber \\
%	&\leq \max_{\mathcal{S}_v\subseteq\mathcal{T}} \min_{R_v\subseteq\mathcal{S}_{v}}\text{split}_{S_{u v}}(R_v) \nonumber \\
	&\leq \max_{\mathcal{S}_v\subseteq\mathcal{T}}\min_{\mathcal R\in \partial \mathcal{S}_v}\text{split}_{\mathcal S_{v}}(\mathcal R) = \text{split}_{\mathcal{T}}(v). \label{split-uv} 
\end{align}
The last inequality \eqref{split-uv} is a consequence of \eqref{splitting vertex}; the class of subtrees of the form $\mathcal S_{v \rightarrow u}$ is a sub-collection of all sub-trees rooted at $v$, as a result of which the set depicting the minimal number of splitting vertices along a ray of $\partial \mathcal S_{v \rightarrow u}$ obeys the following inclusion: 
\[ \Bigl\{ \min_{\mathcal R\in \partial \mathcal{S}_{v \rightarrow u}}\text{split}_{\mathcal S_{v \rightarrow u}}(\mathcal R) : \mathcal S_u \subseteq \mathcal T \Bigr\}  \subseteq \Bigl\{ \min_{\mathcal R\in \partial \mathcal{S}_v} \text{split}_{\mathcal S_{v}}(\mathcal R) :  \mathcal S_v \subseteq \mathcal T \Bigr\} \subseteq \mathbb R. \]  
Taking the maximum over both sets leads to the claimed inequality \eqref{split-uv}, proving \eqref{monotonicity: splitting number of vertices}.
\vskip0.1in
\noindent Let us continue to the second part of the lemma. Set $s := \text{split}_{\mathcal T}(u)$. The definition \eqref{splitting vertex} of the splitting number of a vertex and the hypotheses \eqref{uu':hypotheses} on $u, u'$ imply the existence of subtrees $\mathcal S_u, \mathcal S_{u'}$, rooted respectively at $u, u'$, with the property that every ray in $\partial \mathcal S_{u}, \partial \mathcal S_{u'}$ has at least $s$ splitting vertices within their respective subtrees. Let $\mathcal S^{\ast}$ denote the smallest subtree of $\mathcal T$ rooted at $v$ that contains both $\mathcal S_u, \mathcal S_{u'}$. Equivalently stated, $\mathcal S^{\ast}$ consists of  $\mathcal S_u, \mathcal S_{u'}$ and two rays joining $v$ to $u, u'$ respectively. 
\vskip0.1in
\noindent We claim that every ray $\mathcal R \in \partial \mathcal S^{\ast}$ has at least $(s+1)$ splitting vertices in $\mathcal S^{\ast}$. Suppose that the last vertex on the ray $\mathcal R$ lies in $\mathcal S_u$. Since $v$ is a splitting vertex in $\mathcal S^{\ast}$, the splitting vertices on $\mathcal R$ consist of those on $\mathcal R \cap \mathcal S_u$, as well as $v$. Thus $\mathcal R$ contains at least $(s+1)$ splitting vertices of $\mathcal S^{\ast}$. The same argument applies to rays $\mathcal R \in \partial \mathcal S^{\ast}$ that end in a vertex of $\mathcal S_{u'}$. Summarizing this discussion, we find that 
\[\text{split}_{\mathcal T}(v) \geq \min_{\mathcal R \in \mathcal S^{\ast}} \text{split}(\mathcal R) \geq (s+1) > s = \text{split}_{\mathcal T}(u). \] This completes the proof of the lemma.   
\end{proof}
\noindent An immediate consequence of \eqref{monotonicity: splitting number of vertices}  in Lemma~\ref{monotonicity} is the following. 
\begin{corollary} \label{corollary-monotonicity}
For any tree $\mathcal T$, 
\[ \text{split}(\mathcal{T}) = \text{split}_{\mathcal{T}}(v_0), \]
 where $v_0$ is the root of $\mathcal{T}$.  
 \end{corollary} 
 \noindent Our next result says that, like splitting numbers of vertices, splitting numbers of trees are also monotonic in an  appropriate sense.
\begin{lemma}\label{monotone trees}
	Let $\mathcal{S}\subseteq\mathcal{T}$.  Then $\text{split}(\mathcal{S})\leq \text{split}(\mathcal{T})$.
\end{lemma}
\begin{proof}
	We use Corollary \ref{corollary-monotonicity} to compare sub-trees of $\mathcal S$ and $\mathcal T$ that appear in the definitions \eqref{splitting vertex}, \eqref{tree split} of $\text{split}(\mathcal{S})$ and $\text{split}(\mathcal{S})$. By Corollary~\ref{corollary-monotonicity}, split$(\mathcal{S}) = \text{split}_{\mathcal{S}}(v_0)$, where $v_0$ is the root of $\mathcal{S}$.  Since $v_0\in\mathcal{S}\subseteq\mathcal{T}$ and any subtree of $\mathcal S$ is also a subtree of $\mathcal T$, we find that 
\begin{align*}
\text{split}_{\mathcal S}(v_0) &= \max_{\mathcal S_{v_0} \subseteq \mathcal S} \min_{R_{v_0}\in \partial \mathcal{S}_{v_0}}\text{split}_{S_{v_0}}(R_{v_0})\\ 
	&\leq \max_{\mathcal S_{v_0} \subseteq \mathcal T} \min_{R_{v_0}\in \partial \mathcal{S}_{v_0}}\text{split}_{S_{v_0}}(R_{v_0}) \\
	&\leq  \text{split}_{\mathcal T}(v_0)  \leq \text{split}(\mathcal T).
\end{align*}
The last two inequalities in the display above follow respectively from the definitions \eqref{splitting vertex} and \eqref{tree split} of the splitting number of a vertex and  a tree. This completes the proof of Lemma~\ref{monotone trees}.
\end{proof}
 
\subsection{The special ray with highest order splits}  A key structural feature of trees with finite splitting number, originally observed in \cite [Lemma 5]{Bateman}, is that the collection of vertices with maximal splitting number is either a singleton or arranged sequentially along a single ray. This specialized ray will turn out to be critical in the detection of lacunary limits. 
\begin{lemma}\label{SplitsOnARay}
Let $\mathcal T$ be a tree rooted at $v_0$ with split$(\mathcal T) = N$, with the property that the collection of vertices 
\begin{equation} \label{highest split}
\mathcal V := \left\{ v \in \mathcal T : \text{split}_{\mathcal T}(v) =N \right\}  
\end{equation}  
contains more than one element. Then there is a unique ray $\mathcal R$ of $\mathcal T$ rooted at $v_0$ containing all the vertices of $\mathcal V$. This means 
\[ v \in \mathcal V \; \text{ if and only if }  \; v \in \mathcal R. \]
This ray need not be of maximal length. 
%a vertex $v$ lies on $R$ if and only if split$_{\mathcal T}(v) = N$, provided the latter collection contains more than one element. 
\end{lemma} 
%\noindent {\em{Remarks: }} 
%\begin{enumerate}[1.]
%\item If $\mathcal V$ in \eqref{highest split} is a singleton, then this single vertex must be the root. 
%\vskip0.1in
%\item If $\mathcal V$ consists of two elements, there is only one ray joining the two; thus, Lemma \ref{SplitsOnARay} is meaningful only when $\#(\mathcal V) \geq 3$. 
%\end{enumerate} 
\begin{proof} 
Let $v \neq v_0$ be any vertex of $\mathcal V$, $h(v) > h(v_0)$. It follows from the monotonicity of splitting numbers (Lemma \ref{monotonicity} and Corollary \ref{corollary-monotonicity}) that 
\[ N = \text{split}_{\mathcal T}(v)  \leq \text{split}_{\mathcal T}(u) \leq \text{split}_{\mathcal T}(v_0) = N, \text{ for any } v \subsetneq u \subsetneq v_0.\]
In other words, if $v \in \mathcal V$, then all the vertices of the ray joining $v$ to the root lie in $\mathcal V$ as well. Thus, in order to prove Lemma \ref{SplitsOnARay}, it suffices to establish the uniqueness of the ray $\mathcal R$.  
\vskip0.1in
\noindent Towards a contradiction, let us assume the existence of a tree $\mathcal T$ of splitting number $N$ satisfying the hypothesis of the lemma, for which the ray $\mathcal R$ containing the highest-order splitting vertices is non-unique. That means there are two vertices $u,v\in\mathcal{V}$ that cannot be joined by a ray; in other words,  
\begin{equation}  \label{no ray connection}
\text{split}_{\mathcal T}(u) = \text{split}_{\mathcal T}(v) = N, \quad u \not\subset v, \; v \not\subset  u. \end{equation}  
Then their youngest common ancestor $D(u,v)$ is neither $u$ nor $v$.  By the first part of Lemma~\ref{monotonicity}, we know that $\text{split}_{\mathcal{T}}(D(u,v)) \geq N$. In view of the non-containment relation \eqref{no ray connection} between $u$ and $v$, the vertex $D(u,v)$ is actually a splitting vertex.  Therefore, the second part of Lemma~\ref{monotonicity} gives that 
\[  \text{split}_{\mathcal{T}}(D(u,v)) > \text{split}_{\mathcal T}(u), \; \text{ i.e., } \; \text{split}_{\mathcal{T}}(D(u,v)) \geq N+1. \]  This contradicts the assumption that $\text{split}(\mathcal{T}) = N$, establishing our claim.
\end{proof} 
\section{Encoding bounded subsets of $\mathbb R$ by $M$-adic trees}\label{tree encoding section}
Up to this point, trees have been treated as abstract combinatorial objects without a priori geometric meaning. We now introduce a concrete realization that encodes bounded subsets of $\mathbb{R}$ as rooted, labelled trees. This correspondence allows geometric information about a set - such as its distribution across $M$-adic scales - to be translated into combinatorial properties of the associated tree, such as branching and splitting. 
\vskip0.1in
\noindent A given set may admit many tree representations, and this flexibility will become important later. In this section, however, we focus on the most natural and canonical model: the $M$-adic tree representation. 
%So far, we have recorded a few properties of rooted, labelled trees in the abstract. At the moment, a tree is a graph-theoretic object whose vertices and edges a priori no further interpretation. However, from an application point of view, the language of rooted, labelled trees is especially convenient for representing bounded sets in Euclidean spaces. The set-to-tree representation is not unique; the same set may be encoded in a variety of trees. Indeed, a key novelty of our proof of \eqref{condition 3: slopes sublacunary} $\implies$ \eqref{condition 1: Kakeya-type sets} is the flexibility offered by different trees in capturing the geometry of a set. 
%\vskip0.1in
%\noindent In this section, we present a tree representation of a set that is the most natural and canonical.   
%This connection is well-studied in the literature.  In this section, we summarize the main points that will be relevant for the study of slope sets and Kakeya-type configurations, referring the interested reader to \cite{LyonsPeres} for more information. 
\subsection{The full $M$-adic tree} \label{section: full M-adic tree} Fix any integer $M\geq 2$.  For any non-negative integer $i$ and positive integer $k$ such that $i<M^k$, there exists a unique representation of the form 
\begin{equation}\label{M-adic representation}
	i = i_1M^{k-1} + i_2M^{k-2} + \cdots + i_{k-1}M + i_k,
\end{equation}
where the integers $i_1,\ldots,i_k$ take values in \[ \mathbb{Z}_M := \{0,1,\ldots,M-1\}. \] These are the digits of the $M$-adic expansion of $i$. As an immediate consequence of \eqref{M-adic representation},  there is a one-to-one correspondence 
\[ \frac{i}{M^k} \mapsto \langle i_1,\ldots, i_k\rangle \] between $M$-adic rationals in $[0,1]$ and finite integer sequences with entries in $\mathbb{Z}_M$. 
%%More generally, for any $\mathbf{i} = (j_1,\cdots,j_d) \in\mathbb{Z}^d$ such that $\mathbf{i}\cdot M^{-k}\in [0,1)^d$, we can apply \eqref{M-adic representation} to each component of $\mathbf{i}$ to obtain 
%%\begin{equation} \label{i expansion} \frac {\mathbf{i}}{M^k} = \frac{1}{M^k}(j_1, \cdots, j_d) = \frac {\mathbf{i}_1}{M} + \frac {\mathbf{i}_2}{M^2} + \cdots + \frac {\mathbf{i}_k}{M^k}, \end{equation}  with $\mathbf{i}_j\in\mathbb{Z}_M^d$ for all $j$.  In this way, we identify $\mathbf{i}$ with $\langle \mathbf{i}_1,\ldots,\mathbf{i}_k\rangle$. Let $\phi : \mathbb{Z}_M^d \rightarrow \{0,1,\ldots, M^d-1\}$ be an enumeration of $\mathbb{Z}_M^d$.  
\vskip0.1in
\noindent We therefore define 
\begin{equation}\label{tree encoding}
\mathcal{T}([0,1];M) = \left\{\langle {i}_1,\ldots, {i}_k \rangle : k\geq 0,\ \mathbf{i}_j \in\mathbb{Z}_M\right\}
\end{equation}
as the {\em{full $M$-adic tree}} encoding $[0,1]$. The geometric meaning of this definition will be made precise in Section \ref{section: rays as points}, where a vertex will be identified with an $M$-adic interval in ([0,1]), and a ray with a nested interval chain. For now, let us note that:
\vskip0.1in
\begin{itemize} 
\item Every vertex of the full $M$-adic tree has exactly $M$ children. 
\vskip0.1in
\item The vertex $\langle i_1, \ldots, i_k \rangle$ is called the $i_k^{\text{th}}$ child of $\langle i_1, \ldots i_{k-1}\rangle$. Thus the children of a given vertex are naturally ordered from left to right according to the value of $i_k$. 
\vskip0.1in 
\item There are exactly $M^{k}$ vertices of the $k^{\text{th}}$ generation, each an $M$-adic interval of length $M^{-k}$.  
\vskip0.1in
\item The tree is of infinite height. 
\vskip0.1in
\item The tree encodes the standard $M$-adic partition of $[0,1]$ at all scales.
\end{itemize} 
\vskip0.1in
%For our purposes, it will suffice to fix $\phi$ to be the lexicographic ordering, and so we will omit the notation for $\phi$ in \eqref{tree encoding}, writing simply, and with a slight abuse of notation, 
%%\begin{equation}\label{better tree encoding}
%%	\mathcal{T}([0,1)^d;M) = \left\{\langle \mathbf{i}_1,\ldots,\mathbf{i}_k\rangle : k\geq 0,\ \mathbf{i}_j \in\mathbb{Z}%_M^d\right\}.
%%\end{equation}
%%We will refer to the tree in \eqref{better tree encoding} by the notation $\mathcal{T}([0,1)^d)$ once the base $M$ has been fixed.
\vskip0.1in
\subsection{Rays as points} \label{section: rays as points} The full $M$-adic tree carries a natural geometric interpretation. Every vertex $u$ corresponds to a unique $M$-adic interval $Q_u$ in $[0,1]$ via 
\begin{equation} \label{cube} 
u = \langle i_1,\ldots, i_k\rangle \in \mathcal{T}([0,1];M) \longleftrightarrow Q_u := \bigl[ i M^{-k}, (i+1) M^{-k} \bigr],
%\left[\frac {j_1}{M^k},\frac {j_1+1}{M^k}\right)\times\cdots\times \left[\frac {j_d}{M^k},\frac {j_d+1}{M^k}\right).
\end{equation}   
%Here $\langle \mathbf i_1, \cdots, \mathbf i_k \rangle$ is related to $(j_1, \cdots, j_d)$ by \eqref{i expansion}. 
with $i$ is determined by the $M$-adic expansion \eqref{M-adic representation}.  Thus vertices of generation $k$ correspond precisely to the intervals in the standard $M$-adic partition of $[0,1]$ at scale $M^{-k}$.
\vskip0.1in
\noindent Under this correspondence, the ancestry relation $u \subset v$ in the tree coincides with the
set inclusion $Q_u \subset Q_v$.
%With this interpretation, tree ancestry among vertices $u, v \in \mathcal T([0,1]; M)$ translates to set containment of the intervals $Q_u, Q_v \subseteq [0,1]$. 
Specifically, suppose that 
\begin{equation} \label{tree + geometry}
\begin{aligned} 
 &u = \langle i_1, \ldots, i_k \rangle \; \text{ and } \; v = \langle j_1, \ldots, j_{\ell}\rangle \in \mathcal T([0,1]; M), \text{ and let } \\  
 &Q_{u} := \sum_{r=1}^{k}\frac{i_r}{M^r}+ [0, M^{-k}] \; \text{ and } \;  Q_v := \sum_{r=1}^{\ell} \frac{j_r}{M^r}  + \bigl[0, M^{-k+1}] 
 \end{aligned} 
\end{equation}   
denote their corresponding intervals. Then
\begin{equation} \label{tree-M-adic-tree} 
\begin{aligned} 
&u  = \langle i_1, \ldots, i_k \rangle \subset v =  \langle j_1, \ldots, j_{\ell}\rangle \text{ as vertices of the tree, i.e. } \\
&\ell \leq k  \text{ and } i_r = j_r  \text{ for all } 1 \leq r \leq \ell,  \; \text{ if and only if } Q_{u} \subset Q_v \text{ as sets.}  
\end{aligned} 
\end{equation} 
Accordingly, every ray determines a nested chain of $M$-adic intervals.
\vskip0.1in
\noindent In view of \eqref{tree-M-adic-tree}, if $\mathcal R \in \partial \mathcal T([0,1]; M)$ is a maximal ray in the sense of Section \ref{section: edges}:
\[ \mathcal R: u(0) = \emptyset \rightarrow u(1) = \langle i_1 \rangle \rightarrow u(2) = \langle i_1, i_2 \rangle  \rightarrow \cdots \rightarrow u(k) = \langle i_1, \ldots, i_k \rangle \rightarrow \cdots \] 
then the associated intervals form a nested chain
 \begin{equation}  \label{nested intervals in k} 
 Q_{u(1)} \supset Q_{u(2)} \supset \ldots \supset Q_{u(k)} \ldots, \quad |Q_{u(k)}| = M^{-k},  
 \end{equation}  
 indexed by $k$. The intersection of this infinite nested sequence consists of a single point in $[0,1]$. 
 \begin{equation} \label{limit of nested intervals} 
 \bigcap_{k=1}^{\infty} Q_{u(k)} = \{x\}. 
 \end{equation}  
In this way, a ray of the $M$-adic tree encodes a point $x \in [0,1]$. We will denote
\begin{equation} \label{ray&point} 
\alpha(\mathcal R) := \text{the point in $[0,1]$ specified by the infinite ray $\mathcal R$}. 
\end{equation} 
Conversely, every $x\in[0,1]$ can be realized in the form \eqref{limit of nested intervals}, i.e., as the intersection of an infinite, nested sequence of $M$-adic intervals \eqref{nested intervals in k}, where 
\[ x \in Q_{u(k)} \text{ for every } k \geq 1. \] If $x$ is not an $M$-adic rational in $(0,1)$, there is exactly one choice of $Q_{u(k)}$ for a given $k \geq 1$. In contrast, an $M$-adic rational $x$ lies on the boundary of some $M$-adic interval of length $M^{-k}$ for all sufficiently large $k$. In this case, there are exactly two  sequences of intervals containing $x$, depicting convergence to the boundary point $x$ from the left and from the right. In this case, we adopt the convention of choosing the ray that depicts convergence to $x$ from the right. This ray corresponds to the finitary expansion of $x$, namely the one with only finitely many nozero digits. Said differently, 
\begin{equation}  \label{canonical ray convention}
\left\{
\begin{minipage}{0.8 \textwidth}
For $b <1$, a {\em{canonical ray}} of the $M$-adic tree rooted at $v = [a, b]$ represents a point in $v^{\sharp} = [a, b)$ but not $b$ itself. Only the canonical ray representing 1 is allowed to approach 1 from the left. 
\end{minipage}   
\right\}
\end{equation} 
In summary, we view the tree $\mathcal T([0,1]; M)$ in \eqref{tree encoding} as an encoding of the set $[0,1]$ with respect to the base $M$, via the identification \eqref{cube} and the requirement that the chosen canonical rays of the tree are in one-to-one correspondence with $[0,1]$.  
%In the remainder of this paper, a vertex $u = \langle i_1, i_2, \cdots, i_k \rangle  \in \mathcal T([0,1]; M)$ is always identified with the corresponding $M$-adic interval $Q_u$ as in \eqref{tree +  geometry} lying on $[0,1]$. 
With this understanding, the notation $u \subset v$ stands both for tree ancestry and set inclusion.
\vskip0.1in
\subsection{Tree representation of sets} \label{section: tree rep of sets} A few observations emerge in light of the above discussion.  Any subset $E\subseteq[0,1]$ is represented by a subtree $\mathcal T(E;M)$ of $\mathcal{T}([0,1];M)$ of infinite height, in the following way. 
\vskip0.1in 
\noindent  {\em{The tree $\mathcal T(E;M)$ is defined to be the smallest subtree of  $\mathcal{T}([0,1];M)$ that contains, for every $x \in E$, the unique maximal ray rooted at $[0,1]$ identifying $x$.}} The identification is unique subject to the convention regarding $M$-adic rationals stated above. Thus, each vertex $u$ on the tree $\mathcal T(E;M)$ represents an $M$-adic interval of the form \eqref{cube} that has non-trivial intersection with $E$, i.e. for which \[Q_u^{\sharp} \cap E \neq \emptyset; \; \text{ here } \; Q_u^{\sharp} := \begin{cases} 
[a, b) &\text{ if } Q_u = [a, b], b < 1, \\ [a, 1] &\text{ if } Q_u = [a, 1], b = 1. 
\end{cases}  \] As a result, a maximal ray in $\mathcal T(E;M)$ identifies a point in $E$ or its closure. We will refer to $\mathcal T(E; M)$ as {\em{the $M$-adic tree representing $E$.}} 
\vskip0.1in
\noindent Conversely, any subtree $\mathcal T$ of the full $M$-adic tree $\mathcal T([0,1] ; M)$ uniquely identifies a closed subset $E \subseteq [0, 1]$, via the defining relation
\[ E := \left\{ \alpha(\mathcal R) : \mathcal R \in \partial \mathcal T \right\}. \] 

\subsection{Trees for general bounded sets} \label{section: trees for bounded sets} The correspondence between a set and its representative $M$-adic tree generalizes to all bounded subsets of $[0, \infty)$, not just subsets of $[0,1]$. The key idea is that every bounded subset of $[0, \infty)$ lies in a sufficiently large $M$-adic interval, which can be mapped to $[0,1]$ using an $M$-adic affine transformation. Such a transformation induces an isomorphism of the $M$-adic trees, leaving the tree structure of the original set unchanged.  
\begin{lemma} \label{lemma: trees under affine maps} 
Given any $M$-adic interval of the form
\[ I = [\ell, \ell+1] M^r \text { for some } \ell, r \in \mathbb Z, \]
let  
\begin{equation} T_{I}: I \rightarrow [0,1], \quad T_I(x) := xM^{-r} - \ell  \label{T: affine tf} \end{equation}  be the affine linear transformation that maps $I$ onto $[0,1]$. Let $\mathcal T(I; M)$ denote the $M$-adic tree for $I$, with $I$ as its root. Then 
\begin{enumerate}[(a)]
\item The trees $\mathcal T(I; M)$ and $\mathcal T([0,1]; M)$ are isomorphic. \label{isomorhism-1}
\vskip0.1in
\item Given any set $E \subseteq I$, the trees \label{isomorhism-2}
\[ \mathcal T(E; M) \subseteq \mathcal T(I; M) \; \text{ and } \mathcal T(T_I(E); M) \subseteq \mathcal T([0,1]; M) \text{ are isomorphic}.   \]
\end{enumerate} 
%$E \subseteq [0, M^r]$ for some large $r \geq 1$, one has to consider a subtree rooted at $v_0$ representing $[0, M^r]$, which is isomorphic to the tree of $E/M^r \subseteq [0,1]$. 
\end{lemma}
%It is worth noting that the tree representation of the set $E$ is, in general, coordinate-sensitive. Indeed, trees representing the same set in two systems of coordinates, i.e. after an arbitrary affine transformation, may possess different features - an issue that we will need to take into account shortly. 
\begin{proof} 
The affine transformation $T_I$ given by \eqref{T: affine tf} induces a tree map 
\[ T_I: \mathcal T(I; M) \rightarrow \mathcal T([0,1]; M) \]  
in the following way: if $v = \ell M^r + s M^{r-k} + [0, M^{-k+r}] \in \mathcal V_k(\mathcal T(I; M))$, then 
\begin{equation}  
T_I(v) =[s, s+1]M^{-k} \in \mathcal V_k([0,1]; M).  
\end{equation} 
In other words, $T_I$ maps every $M$-adic sub-interval of $I$ of length $M^{r-k}$ to an $M$-adic interval of $[0,1]$ of length $M^{-k}$, preserving both height and ancestry. The same is true for the inverse tree map $T_I^{-1}$.  As a result, every sequence of nested $M$-adic sub-intervals of $I$ corresponds to a similar one in $[0,1]$, leading to a bijection between the maximal rays of the two trees. This establishes the isomorphism claimed in both parts \eqref{isomorhism-1} and \eqref{isomorhism-2}.   
\end{proof} 
\vskip0.1in
\noindent We end this section with a summary of the correspondence between the tree $\mathcal T([0,1];M)$ and the interval $[0,1]$:
\begin{itemize} 
\item vertices correspond to $M$-adic subintervals,
\item descendants correspond to nested subintervals,
\item rays correspond to points in $[0,1]$,
\item subtrees correspond to subsets of $[0,1]$.
\end{itemize} 
A visual representation of this correspondence appears in Figure \ref{fig:tree-interval-correspondence}.
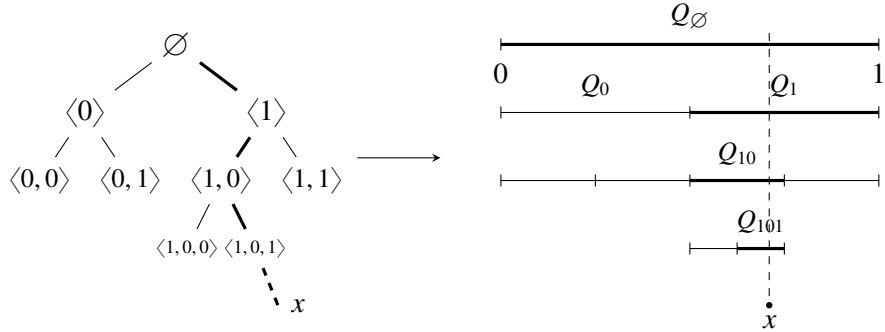
\begin{figure}[ht]
\centering
\begin{tikzpicture}[x=1cm,y=1cm,>=stealth]

%--------------------------
% LEFT: TREE
%--------------------------
% nodes
\node (r)   at (0,3.0) {$\emptyset$};

\node (n0)  at (-1.2,2.1) {$\langle0\rangle$};
\node (n1)  at ( 1.2,2.1) {$\langle1\rangle$};

\node (n00) at (-1.8,1.2) {\small{$\langle0, 0\rangle$}};
\node (n01) at (-0.6,1.2) {\small{$\langle0, 1\rangle$}};
\node (n10) at ( 0.6,1.2) {\small{$\langle1,0\rangle$}};
\node (n11) at ( 1.8,1.2) {\small{$\langle 1,1\rangle$}};

\node (n100) at (0.15,0.3) {\tiny{$\langle1,0,0\rangle$}};
\node (n101) at (1.05,0.3) {\tiny{$\langle1,0,1\rangle$}};

% edges
\draw (r) -- (n0);
\draw[very thick] (r) -- (n1);

\draw (n0) -- (n00);
\draw (n0) -- (n01);

\draw[very thick] (n1) -- (n10);
\draw (n1) -- (n11);

\draw (n10) -- (n100);
\draw[very thick] (n10) -- (n101);

% dashed continuation of highlighted ray
\draw[very thick,dashed] (n101) -- (1.35,-0.45);
\node[right] at (1.38,-0.45) {$x$};

%--------------------------
% CENTER: SEPARATION
%--------------------------
\draw[->] (2.4,1.5) -- (3.5,1.5);

%--------------------------
% RIGHT: INTERVALS
%--------------------------
% level 0
\draw[very thick] (4.3,3.0) -- (9.3,3.0);
\draw (4.3,2.92) -- (4.3,3.08);
\draw (9.3,2.92) -- (9.3,3.08);
\node[below] at (4.3,2.88) {$0$};
\node[below] at (9.3,2.88) {$1$};
\node[above] at (6.8,3.08) {\small $Q_{\emptyset}$};

% level 1
\draw (4.3,2.1) -- (9.3,2.1);
\foreach \x in {4.3,6.8,9.3}
  \draw (\x,2.02) -- (\x,2.18);
\draw[very thick] (6.8,2.1) -- (9.3,2.1);
\node[above] at (5.55,2.18) {\small $Q_0$};
\node[above] at (8.05,2.18) {\small $Q_1$};

% level 2
\draw (4.3,1.2) -- (9.3,1.2);
\foreach \x in {4.3,5.55,6.8,8.05,9.3}
  \draw (\x,1.12) -- (\x,1.28);
\draw[very thick] (6.8,1.2) -- (8.05,1.2);
\node[above] at (7.425,1.28) {\small $Q_{10}$};

% level 3 inside Q_10
\draw (6.8,0.3) -- (8.05,0.3);
\foreach \x in {6.8,7.425,8.05}
  \draw (\x,0.22) -- (\x,0.38);
\draw[very thick] (7.425,0.3) -- (8.05,0.3);
\node[above] at (7.7375,0.38) {\small $Q_{101}$};

% dashed vertical convergence line
\draw[dashed] (7.85,3.15) -- (7.85,-0.45);
\fill (7.85,-0.45) circle (1.2pt);
\node[below] at (7.85,-0.45) {$x$};

\end{tikzpicture}
\caption{\small{The dyadic encoding of $[0,1]$, where $M=2$. Vertices correspond to dyadic intervals, ancestry corresponds to inclusion, and an infinite ray determines a point $x\in[0,1]$ via the intersection of the associated nested intervals.}}
\label{fig:tree-interval-correspondence}
\end{figure}
\section{Examples: Splitting numbers of some $M$-adic trees} \label{splitting number examples}
In this section, we illustrate the correspondence between subsets of $\mathbb{R}$
and their associated $M$-adic trees through a series of examples. Our primary goal
is to compute the splitting numbers of these trees and relate them to the lacunarity
properties of the underlying sets. These examples complement those in Section \ref{EXAMPLES SECTION} by recasting familiar lacunary and sublacunary sets in the language of trees.
%The notions introduced in this chapter and the last one are intimately related. This was first anticipated in \cite{Bateman}, and we will prove the precise connection in the context of our work in Chapter \ref{section: split implies lacunarity}. The property of finite order lacunarity of a set $\Omega \subseteq [0,1]$, as given by Definition \ref{defn: Lacunary sets}, is reflected in the finiteness of the splitting number of its $M$-adic tree $\mathcal T(\Omega; M)$, as defined in \eqref{splitting vertex} and \eqref{tree split}. To motivate this connection before a formal proof, let us record a few examples of trees related to sets we have encountered in Section \ref{EXAMPLES SECTION} and compute their splitting numbers.  
\subsection{Finite sets} \label{section: split finite}
We begin with the simplest example of a finite set. The $M$-adic tree for such a set has only finitely many rays, so the splitting behaviour of such a tree  is controlled by the cardinality of the set. This reflects the absence of hierarchical clustering in the set, which in turn corresponds to lacunarity of order $0$.
\begin{lemma} \label{lemma: finite set implies finite splitting number}  
A finite set has finite splitting number; for any integer $M \geq 2$, 
\begin{equation} \label{finite set finite split}
\text{ if } \#(\Omega) < \infty, \text{ then } \text{split} \bigl( \mathcal T(\Omega; M)\bigr) \leq \log_2 \bigl(\#(\Omega) \bigr).  
\end{equation} 
\end{lemma} 
\begin{proof} 
Suppose, towards a contradiction, that there exists a finite set $\Omega$ with 
\begin{equation} \label{finite Omega}
\#(\Omega) \leq 2^m \quad \text{ but } \quad \text{split} \bigl( \mathcal T(\Omega; M)\bigr) > m.
\end{equation} 
It follows from the definition \eqref{splitting vertex}, \eqref{tree split} of a splitting number that $\mathcal T(\Omega; M)$ contains a sub-tree, say $\mathcal T_0$ such that every maximal ray $\mathcal R \in \partial \mathcal T_0$ contains at least $(m+1)$ splitting vertices. An inductive argument based on $m$ shows that there must be a height $h_0$ where 
\[ \# \bigl\{v \in \mathcal T_0 : h(v) = h_0 \bigr\} \geq 2^{m+1}. \]  
Since the collection of vertices of $\mathcal T_0$ at a given height represents a set of distinct elements of $\Omega$, this implies $\#(\Omega) \geq 2^{m+1}$, contradicting \eqref{finite Omega}. This proves \eqref{finite set finite split}.
\end{proof} 
\vskip0.1in
\noindent For $M=2$, the inequality in \eqref{finite set finite split} is sharp. For instance, 
\begin{equation} \label{splitting number sharpness}  
\left\{
\begin{aligned} 
&\text{ if } \Omega_m :=  \left\{ \frac{r}{M^m} : 0 \leq r < M^m \right\}, \text{ then } \#(\Omega_m) = M^m, \text{ and } \\  
& \text{ split}(\mathcal T\bigl(\Omega_m; M) \bigr) = \log_M\bigl(\#(\Omega_m) \bigr) = m \\
 \end{aligned}
\right\}. 
\end{equation}
In this case, branching occurs at every level of the tree upto height $m$, so each ray contains $m$ splitting vertices. Therefore, the splitting number of the tree is equal to $m = \log_M\bigl(\#(\Omega_m) \bigr)$. Figure \ref{fig:Q3-tree} depicts the tree for $\mathbb Q_3$. 
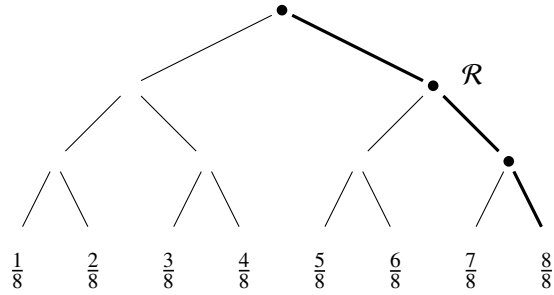
\begin{figure}[ht]
\centering
\begin{tikzpicture}[x=1cm,y=1cm,>=stealth]

% Root
\node (r) at (0,4) {};

% Level 1
\node (a) at (-2,3) {};
\node (b) at ( 2,3) {};

% Level 2
\node (a1) at (-3,2) {};
\node (a2) at (-1,2) {};
\node (b1) at ( 1,2) {};
\node (b2) at ( 3,2) {};

% Level 3 (leaves)
\node (x1) at (-3.5,1) {};
\node (x2) at (-2.5,1) {};
\node (x3) at (-1.5,1) {};
\node (x4) at (-0.5,1) {};
\node (x5) at ( 0.5,1) {};
\node (x6) at ( 1.5,1) {};
\node (x7) at ( 2.5,1) {};
\node (x8) at ( 3.5,1) {};

% Edges
\draw (r) -- (a);
\draw (r) -- (b);

\draw (a) -- (a1);
\draw (a) -- (a2);
\draw (b) -- (b1);
\draw (b) -- (b2);

\draw (a1) -- (x1);
\draw (a1) -- (x2);
\draw (a2) -- (x3);
\draw (a2) -- (x4);
\draw (b1) -- (x5);
\draw (b1) -- (x6);
\draw (b2) -- (x7);
\draw (b2) -- (x8);

% Highlight one ray
\draw[very thick] (r) -- (b);
\draw[very thick] (b) -- (b2);
\draw[very thick] (b2) -- (x8);

% Splitting vertices on highlighted ray
\fill (r)  circle (2pt);
\fill (b)  circle (2pt);
\fill (b2) circle (2pt);

% Leaf labels
\node[below=2pt] at (x1) {$\frac18$};
\node[below=2pt] at (x2) {$\frac28$};
\node[below=2pt] at (x3) {$\frac38$};
\node[below=2pt] at (x4) {$\frac48$};
\node[below=2pt] at (x5) {$\frac58$};
\node[below=2pt] at (x6) {$\frac68$};
\node[below=2pt] at (x7) {$\frac78$};
\node[below=2pt] at (x8) {$\frac88$};

% Small label
\node[right] at (2.25,3.15) {$\mathcal R$};

\end{tikzpicture}
\caption{\small{The dyadic tree encoding of $\mathbb Q_3=\{{r}/{2^3} :0\le r< 2^3\}$. This is a full binary tree of height $3$, with leaves corresponding to the points of $\mathbb Q_3$. Every non-terminal vertex splits, so each maximal ray $\mathcal R$ passes through three splitting vertices. Hence the splitting number of the tree $\mathcal T(\mathbb Q_3; 2)$ is $3$.}}
\label{fig:Q3-tree}
\end{figure}
\vskip0.1in   
\noindent Building on this example, let
\[ \Omega_{\infty} : = \bigcup_{m=0}^{\infty} \Omega_m\] 
denote the set of all $M$-adic rationals in $[0,1]$. The set is dense in $[0,1]$, and hence sublacunary according to Corollary \ref{corollary: uncountable sublacunary}. A consequence of \eqref{splitting number sharpness} is that $\mathcal T(\Omega_{\infty}; M)$ has infinite splitting number; indeed, in contrast with the previous example $\Omega_m$, where branching terminates after a few levels, here branching persists at every generation. 

 \subsection{Infinite sets whose $M$-adic trees have finite splitting number} \label{section: split finite order lacunary}
 The converse of Lemma \ref{lemma: finite set implies finite splitting number} is not true in general. In this section, we consider examples of infinite lacunary sets of finite order, and examine how their hierarchical structure is reflected in the associated trees. 
\begin{lemma} \label{lemma: splitting number of a geometric sequence} 
If $\Omega_1 = \{ M^{-j} : j \geq 1\}$, then split$(\mathcal T(\Omega_1;M)) = 1$.
\end{lemma} 
\begin{proof} 
%A tree with zero splitting number is a ray, and therefore depicts a single point. Since $\Omega_1$ is infinite, it follows that 
%\[\text{split}(\mathcal T(\Omega; M)) \geq 1. \] For the converse inequality, let us observe that 
Let us describe the structure of the tree $\mathcal T(\Omega_1;M)$ rooted at $v_0 = [0,1]$. It consists of a unique special ray $\mathcal R_0$, corresponding to the limit point $\alpha(\mathcal R_0) = 0$ of the set $\Omega_1$. The ray $\mathcal R_0$ contains all the splitting vertices of the tree. Conversely, every vertex $v$ of $\mathcal R_0$ is a splitting vertex. Suppose that 
\[ h(v) = k-1 \geq 0, \quad k \geq 1. \] Then $v$ has two children in $\mathcal T(\Omega_1;M)$, namely $[0, M^{-k}]$ and $[M^{-k}, 2M^{-k}]$. The former is a vertex of $\mathcal R_0$. The latter is the root of a single non-splitting ray $\mathcal R \neq \mathcal R_0$, leading to $M^{-k}$. 
\vskip0.1in
\noindent To summarize, every maximal ray in $\mathcal T(\Omega_1;M)$ contains at least one splitting vertex, namely $v_0$. Further, there is one maximal ray, specifically the one corresponding to $M^{-1}$, that contains exactly one splitting vertex. It follows from Corollary \ref{corollary-monotonicity} and the definition \eqref{splitting vertex} that  
\[ \text{split}(\mathcal T(\Omega_1; M)) = \text{split}_{\mathcal T}(v_0) = 1. \] 
This establishes the conclusion of the lemma. 
\end{proof} 
\noindent Figure \ref{fig:split-one} depicts the dyadic tree for the geometric sequence $\{2^{-j} : j \geq 1\}$, clarifying its branching structure that leads to its splitting number 1. 
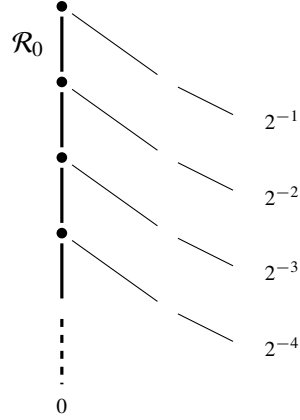
\begin{figure}[ht]
\centering
\begin{tikzpicture}[x=1cm,y=1cm,>=stealth]

% Special ray corresponding to 0
\node (r0) at (0,4.5) {};
\node (r1) at (0,3.5) {};
\node (r2) at (0,2.5) {};
\node (r3) at (0,1.5) {};
\node (r4) at (0,0.5) {};

% Side branches corresponding to 2^{-j}
\node (a1) at (1.4,3.5) {};
\node (a2) at (1.4,2.5) {};
\node (a3) at (1.4,1.5) {};
\node (a4) at (1.4,0.5) {};

% Optional continuation of nonsplitting rays
\node (b1) at (2.4,3.0) {};
\node (b2) at (2.4,2.0) {};
\node (b3) at (2.4,1.0) {};
\node (b4) at (2.4,0.0) {};

% Special ray edges
\draw[very thick] (r0) -- (r1);
\draw[very thick] (r1) -- (r2);
\draw[very thick] (r2) -- (r3);
\draw[very thick] (r3) -- (r4);
\draw[very thick,dashed] (r4) -- (0,-0.5);

% Branches splitting off from the special ray
\draw (r0) -- (a1);
\draw (r1) -- (a2);
\draw (r2) -- (a3);
\draw (r3) -- (a4);

% Continuation of the nonsplitting rays
\draw (a1) -- (b1);
\draw (a2) -- (b2);
\draw (a3) -- (b3);
\draw (a4) -- (b4);

% Splitting vertices on the special ray
\fill (r0) circle (2pt);
\fill (r1) circle (2pt);
\fill (r2) circle (2pt);
\fill (r3) circle (2pt);

% Labels
\node[left] at (-0.1,4.0) {$\mathcal R_0$};
\node[right] at (2.55,3.0) {\scriptsize $2^{-1}$};
\node[right] at (2.55,2.0) {\scriptsize $2^{-2}$};
\node[right] at (2.55,1.0) {\scriptsize $2^{-3}$};
\node[right] at (2.55,0.0) {\scriptsize $2^{-4}$};
\node[below] at (0,-0.55) {\scriptsize $0$};

\end{tikzpicture}
\caption{\small{A schematic dyadic tree associated to the infinite set $\{2^{-j}: j\geq 1\}$. The leftmost distinguished ray $\mathcal R_0$ corresponds to the limit point $0$. At each vertex of $\mathcal R_0$, a non-splitting branch separates off to identify one point $2^{-j}$. Thus every ray of the tree contains at least one splitting vertex, and the splitting number of the tree is $1$.}}
\label{fig:split-one}
\end{figure}
\noindent Lemma \ref{lemma: splitting number of a geometric sequence} generalizes to iterated sumsets of the type considered in Lemma \ref{Lemma: iterated sums 1}. 
\begin{lemma} \label{lemma: splitting number of iterated sumsets} 
The set 
\begin{equation} \label{Omega_N tree} 
\Omega_N = \Bigl\{ \sum_{k=1}^{N} M^{-j_k} : 1 \leq j_1 < j_2 < \ldots < j_N \Bigr\} \in \Lambda(N, M^{-1})  
\end{equation} 
has the property split$(\mathcal T(\Omega_N; M)) = N$.  
\end{lemma} 
\begin{proof} 
We proceed by induction on $N$, the preceding example serving as the base case $N=1$. The ray $\mathcal R_0$ with 
\[ \alpha(\mathcal R_0) = 0 \in \text{closure} (\Omega_N) \] plays a special role in  the tree  $\mathcal T_{N} = \mathcal T(\Omega_N; M)$, similar to $\mathcal R_0$ in Lemma \ref{lemma: splitting number of a geometric sequence}. Specifically, for $j_1 \geq 1$, suppose that 
\[ v(j_1-1) = [0, M^{-j_1+1}] \text{ is the vertex on $\mathcal R_0$ with $h(v(j_1)) = j_1-1$.}  \] 
Let us observe that $v(j_1-1)$ splits into two children in $\mathcal T_N$, one of which is $v(j_1)$, and the other is the $M$-adic interval 
\[ v'(j_1) := [M^{-j_1}, 2M^{-j_1}] \text{ to the immediate right of $v(j_1)$}. \] This follows from the estimate: for all $j_1 < j_2 < \ldots < j_N$, 
\begin{align*} 
M^{-j_1} &< M^{-j_1} + M^{-j_2} + \ldots + M^{-j_N}  \\ & < M^{-j_1} \Bigl[1 + \sum_{k=1}^{\infty} M^{-k} \Bigr] = M^{-j_1} \left[ \frac{M}{M-1}\right] < 2M^{-j_1}. \end{align*}  
The maximal subtree of $\mathcal T_N' \subseteq \mathcal T_N$ rooted at $v'(j_1)$ represents an affine copy of the set $\Omega_{N-1}$; specifically, 
\[ \mathcal T_N' = \mathcal T \bigl( M^{-j_1} + M^{-j_1} \Omega_{N-1}; M \bigr). \] 
 The affine transformation that maps $[M^{-j_1}, 2M^{-j_1}]$ onto $[0,1]$ is of the form \eqref{T: affine tf} mentioned in Lemma \ref{lemma: trees under affine maps}. According to this lemma, $\mathcal T_N'$ and $\mathcal T( \Omega_{N-1}; M)$ are isomorphic, and therefore share the same splitting number. The induction hypothesis then yields
\begin{equation} \label{split lower order} \text{split}_{\mathcal T_N}(v'(j_1)) = \text{split}(\mathcal T_N') = N-1 \text{ for all } j_1 \geq 1. \end{equation}  
We claim that \eqref{split lower order} implies $\text{split}(\mathcal T_N) = N$. This follows from the second part of Lemma \ref{monotonicity}. Indeed, the root $v(0)$ of $\mathcal T_N$ has two children $u = v(1)$ and $u' = v'(1)$, so that $u \not\subset u'$ and $u' \not\subset u$. On one hand, $v(1)$ is the parent of $v'(2)$, therefore applying \eqref{split lower order} with $j_1=2$ gives  
\[ \text{split}_{\mathcal T_N}(v(1)) \geq \text{split}_{\mathcal T_N}(v'(2))=  N-1. \] 
On the other hand, applying \eqref{split lower order} with $j_1=1$ leads to 
\begin{equation}  \label{split v'(1)} 
\text{split}_{\mathcal T_N}(v'(1)) =  N-1. 
\end{equation}  
Lemma \ref{monotonicity} then allows us to conclude that 
\begin{align}  
\text{split}(\mathcal T_N) &= \text{split}_{\mathcal T_N}(v(0)) \nonumber \\ 
&> \min \left[ \text{split}_{\mathcal T_N}(v(1)),  \text{split}_{\mathcal T_N}(v'(1)) \right] \geq  N-1, \text{ i.e., } \nonumber \\
\text{split}(\mathcal T_N) &= \text{split}_{\mathcal T_N}(v(0)) \geq N. \label{root-lower} 
\end{align}  
At the same time, $v(0)$ is the parent of $v'(1)$, i.e., only one generation above it; therefore \eqref{split v'(1)} ensures 
\begin{equation} \label{root-upper} 
\text{split}(\mathcal T_N) = \text{split}_{\mathcal T_N}(v(0)) \leq \text{split}_{\mathcal T_N}(v'(1)) + 1 = N.
\end{equation} 
 Combining \eqref{root-lower} and \eqref{root-upper}, we obtain $\text{split}(\mathcal T_N) = N$. This completes the induction and hence the proof of the lemma.
\end{proof}  
\noindent  Figure \ref{fig:split-two-local} depicts the binary tree representation of $\{2^{-j} + 2^{-k} : k > j\}$, showing that its splitting number is 2. 
\begin{figure}[ht]
\centering
\begin{tikzpicture}[x=1cm,y=1cm,>=stealth]

%------------------------------------------------
% Main special ray corresponding to 0
%------------------------------------------------
\node (r0) at (0,5.2) {};
\node (r1) at (0,4.2) {};
\node (r2) at (0,3.2) {};
\node (r3) at (0,2.2) {};
\node (r4) at (0,1.2) {};

\draw[very thick] (r0) -- (r1);
\draw[very thick] (r1) -- (r2);
\draw[very thick] (r2) -- (r3);
\draw[very thick] (r3) -- (r4);
\draw[very thick,dashed] (r4) -- (0,0.3);

\fill (r0) circle (2pt);
\fill (r1) circle (2pt);
\fill (r2) circle (2pt);
\fill (r3) circle (2pt);

\node[left] at (-0.1,4.8) {$\mathcal R_0$};
\node[below] at (0,0.25) {\scriptsize $0$};

%------------------------------------------------
% Main highlighted branch corresponding to j = 1
%------------------------------------------------
\node (a1) at (1.5,4.2) {};
\draw (r0) -- (a1);

% inner special ray converging to 2^{-1}
\node (a11) at (2.4,3.6) {};
\node (a12) at (2.4,3.0) {};
\node (a13) at (2.4,2.4) {};

\draw[very thick] (a1) -- (a11);
\draw[very thick] (a11) -- (a12);
\draw[very thick] (a12) -- (a13);
\draw[very thick,dashed] (a13) -- (2.4,0.8);

\fill (a1)  circle (2pt);
\fill (a11) circle (2pt);
\fill (a12) circle (2pt);

\node[right] at (2.55,0.8) {\scriptsize $2^{-1}$};
\node[right] at (1.8,1.6) {\scriptsize $2^{-2}$};
%------------------------------------------------
% Side branches encoding 2^{-1}+2^{-k}, k>1
%------------------------------------------------
\node (b1) at (3.4,3.5) {};
\node (b2) at (3.4,3.0) {};
\node (b3) at (3.4,2.4) {};

\draw (a1)  -- (b1);
\draw (a11) -- (b2);
\draw (a12) -- (b3);

\node[right] at (3.55,3.6) {\scriptsize $2^{-1}+2^{-2}$};
\node[right] at (3.55,3.0) {\scriptsize $2^{-1}+2^{-3}$};
\node[right] at (3.55,2.4) {\scriptsize $2^{-1}+2^{-4}$};

%------------------------------------------------
% Ghost subtree indicating similar structure for j = 2,3,...
%------------------------------------------------
\node (g1)  at (1.2,3.2) {};
\node (g11) at (1.9,2.5) {};
\node (g12) at (1.9,1.5) {};

\draw[dashed] (r1) -- (g1);
\draw[dashed] (g1) -- (g11);
\draw[dashed] (g11) -- (g12);

\node (h1) at (2.3,2.7) {};
\node (h2) at (2.7,2.35) {};

\draw[dashed] (g1)  -- (h1);
\draw[dashed] (g11) -- (h2);

\end{tikzpicture}
\caption{\small{A schematic subtree of the dyadic tree associated to $\{2^{-j}+2^{-k}: k>j\}$. The main special ray $\mathcal R_0$ corresponds to the limit point $0$. The highlighted branch shown here corresponds to $j=1$; along it, a lower-order special ray converges to $2^{-1}$, and its side branches encode the points $2^{-1}+2^{-k}$ for $k>1$. The light dotted subtree indicates that similar structure recurs for $2^{-2}, 2^{-3}, \dots$.}}
\label{fig:split-two-local}
\end{figure}
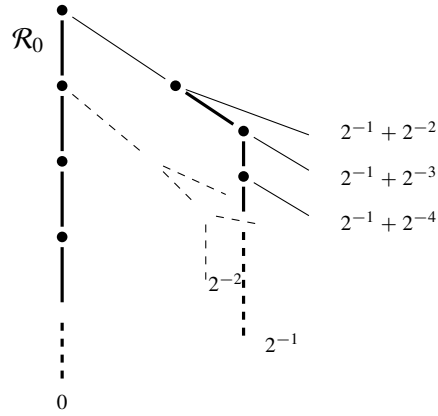
\vskip0.1in
\noindent  In summary, the splitting number provides a combinatorial measure of the depth of hierarchical clustering, and coincides in this case with the lacunarity order of the set \eqref{Omega_N tree}. More flexible constructions, such as those in Lemma \ref{Lemma: iterated sums 2}, allow independent variation of the parameters while preserving admissible finite-order lacunarity. In these cases, the associated trees still exhibit controlled splitting up to finitely many scales, and hence have finite splitting number. In other words, finiteness of lacunarity order in equivalent to finiteness of splitting number, even if the numbers may not always match. This fact, supported thus far only through examples, will be proved rigorously in Chapter \ref{section: split implies lacunarity}.   
\vskip0.1in
\noindent The ray $\mathcal R_0$ that appeared in Lemmas \ref{lemma: splitting number of a geometric sequence} and \ref{lemma: splitting number of iterated sumsets} is distinguished by the property that it is the unique ray of the respective trees containing all the vertices with the highest splitting number. The existence of such a ray is not specific to the trees considered in these examples. In Proposition \ref{TREE-LACUNARY-TRADITIONAL} of Chapter \ref{section: split implies lacunarity}, we will see that the $M$-adic tree for any set in $\Lambda(N, M^{-1})$ contains a ray on which the vertices of highest split lie.     
\subsection{Splitting number of a tree associated with $\Omega_{\text{HRS}}$} \label{section: HRS-tree} 
Let us recall from \eqref{HRS-example}, \eqref{HRS-example-unit} the building block $\Omega_{\text{HRS}}(R)$ of the counter-example $\Omega_{\text{HRS}}$ presented in \cite{{HRS2024a},{HRS2024b}}. In this section, we analyse the $M$-adic trees representing $\Omega_{\text{HRS}}(R)$ and $\Omega_{\text{HRS}}$ and compute their splitting numbers.
\vskip0.1in
\noindent The structure of $\Omega_{\text{HRS}}$ is encoded through concatenations of binary blocks. At each stage $j$, there are two possible choices for the block of length $N_j$, namely $\pmb{\eta}_j$ and $\pmb{\zeta}_j$. Each block choice introduces a bifurcation at the level corresponding to the beginning of that block. Thus, each maximal ray of $\mathcal T(\Omega_{\text{HRS}}; 2)$ corresponds to an infinite concatenation of blocks, with splitting at infinitely many prescribed heights. In contrast to the examples in Sections \ref{section: split finite} and \ref{section: split finite order lacunary}, where splitting occurs only finitely many times along any ray, the present construction exhibits splitting at arbitrarily many scales. This leads to the following lemma. 
\begin{lemma}
For $R \geq 1$, let $\Omega_{\text{HRS}}(R)$ be the set in \eqref{HRS-example-unit}. Then 
\begin{equation}  \text{split}(\mathcal T\bigl(\Omega_{\text{HRS}}(R); 2 \bigr) ) = R. \label{HRS-unit-split} \end{equation} 
As a result, 
\begin{equation}  \text{split}(\mathcal T\bigl(\Omega_{\text{HRS}}; 2 \bigr)) = \infty. \label{splitting number of HRS tree} \end{equation} 
\end{lemma} 
\begin{proof}
The definition \eqref{HRS-example-unit} implies that $\#(\Omega_{\text{HRS}}(R)) = 2^R$. Therefore by Lemma \ref{lemma: finite set implies finite splitting number}, 
\[ \text{split}(\mathcal T\bigl(\Omega_{\text{HRS}}(R); 2) \bigr) \leq R. \]
On the other hand, every maximal ray of $\mathcal  T\bigl(\Omega_{\text{HRS}}(R); 2) \bigr)$ splits exactly $R$ times, at the heights $\bar{N}_j$, $0 \leq j \leq R-1$, where  
\[ \bar{N}_0 = 0, \quad \bar{N}_j = N_1 + \ldots + N_j, \quad 1 \leq j \leq R-1. \] In view of the definition of a splitting number (see \eqref{splitting vertex}, \ref{tree split}), this yields the converse inequality, 
\[ \text{split}(\mathcal T\bigl(\Omega_{\text{HRS}}(R); 2) \bigr) \geq R. \] 
This proves \eqref{HRS-unit-split}. 
\vskip0.1in
\noindent The second claim \eqref{splitting number of HRS tree}  follows from the fact that there is no uniform bound on the number of splitting vertices along the maximal rays of $\mathcal T(\Omega_{\text{HRS}}; 2)$. Indeed, the monotonicity property in Lemma \ref{monotone trees} and the inclusion
\[ \Omega_{\text{HRS}}(R) \subseteq \Omega_{\text{HRS}} \quad \text{ imply that } \quad \mathcal T\bigl(\Omega_{\text{HRS}}(R); 2 \bigr) \subseteq \mathcal T\bigl(\Omega_{\text{HRS}}; 2 \bigr) \text{ for all } R \geq 1. \] 
This means that
\[ \text{split}\bigl(\mathcal T(\Omega_{\text{HRS}}; 2) \bigr)  \geq  \max_{R \geq 1}\text{split} \bigl( \mathcal T(\Omega_{\text{HRS}}(R); 2) \bigr) = \infty,  \] 
completing the proof.   
\end{proof}
\begin{figure}[ht]
\centering
\begin{tikzpicture}[x=1cm,y=1cm,>=stealth]

% top vertex (root)
\node (r) at (0,4.8) {};

% first generation
\node (L1) at (-1.8,3.2) {};
\node (R1) at ( 1.8,3.2) {};

% second generation
\node (L2L) at (-2.6,1.4) {};
\node (L2R) at (-1.0,1.4) {};
\node (R2L) at ( 1.0,1.4) {};
\node (R2R) at ( 2.6,1.4) {};

% dashed continuation below
\node (D1a) at (-2.9,0.2) {};
\node (D1b) at (-2.3,0.2) {};

\node (D2a) at (-1.3,0.2) {};
\node (D2b) at (-0.7,0.2) {};

\node (D3a) at ( 0.7,0.2) {};
\node (D3b) at ( 1.3,0.2) {};

\node (D4a) at ( 2.3,0.2) {};
\node (D4b) at ( 2.9,0.2) {};

% solid edges
\draw[thick] (r) -- (L1);
\draw[thick] (r) -- (R1);

\draw[thick] (L1) -- (L2L);
\draw[thick] (L1) -- (L2R);

\draw[thick] (R1) -- (R2L);
\draw[thick] (R1) -- (R2R);

% circled root (level 0)
\draw[fill=white, line width=0.6pt] (r) circle (2pt);

% circled vertices at level N1
\draw[fill=white, line width=0.6pt] (L1) circle (2pt);
\draw[fill=white, line width=0.6pt] (R1) circle (2pt);

% circled vertices at level N1 + N2
\draw[fill=white, line width=0.6pt] (L2L) circle (2pt);
\draw[fill=white, line width=0.6pt] (L2R) circle (2pt);
\draw[fill=white, line width=0.6pt] (R2L) circle (2pt);
\draw[fill=white, line width=0.6pt] (R2R) circle (2pt);

% dashed continuation
\draw[dashed] (L2L) -- (D1a);
\draw[dashed] (L2L) -- (D1b);

\draw[dashed] (L2R) -- (D2a);
\draw[dashed] (L2R) -- (D2b);

\draw[dashed] (R2L) -- (D3a);
\draw[dashed] (R2L) -- (D3b);

\draw[dashed] (R2R) -- (D4a);
\draw[dashed] (R2R) -- (D4b);

% horizontal level lines
\draw[dashed] (-2.6,4.8) -- (2.6,4.8);   % level 0
\draw[dashed] (-2.5,3.2) -- (3.2,3.2);   % N1
\draw[dashed] (-3.2,1.4) -- (3.2,1.4);   % N1+N2

% labels for levels
\node[right] at (2.8,4.8) {$0$};
\node[right] at (3.3,3.2) {$N_1$};
\node[right] at (3.3,1.4) {$N_1+N_2$};

% block labels
\node[left]  at (-0.8,4.1) {$\pmb{\eta}_1$};
\node[right] at ( 0.8,4.1) {$\pmb{\zeta}_1$};

\node[left]  at (-2.2,2.25) {$\pmb{\eta}_2$};
\node[right] at (-1.4,2.25) {$\pmb{\zeta}_2$};

\node[left]  at (1.4,2.25) {$\pmb{\eta}_2$};
\node[right] at ( 2.2,2.25) {$\pmb{\zeta}_2$};

\end{tikzpicture}
\caption{\small{A schematic block-level picture for $\mathcal T(\Omega_{HRS}; 2)$, which looks like an elongated version of the full binary tree. The root lies at level $0$. There are two vertices at height $N_1$ corresponding to the two block choices $\pmb{\eta}_1$ and $\pmb{\zeta}_1$. Each of these vertices splits into two children, corresponding to the block choices $\pmb{\eta}_2$ and $\pmb{\zeta}_2$. The circled vertices mark splitting points at successive block levels, and the dashed continuation indicates that this branching persists at deeper levels.}}
\label{fig:hrs-block-tree}
\end{figure}
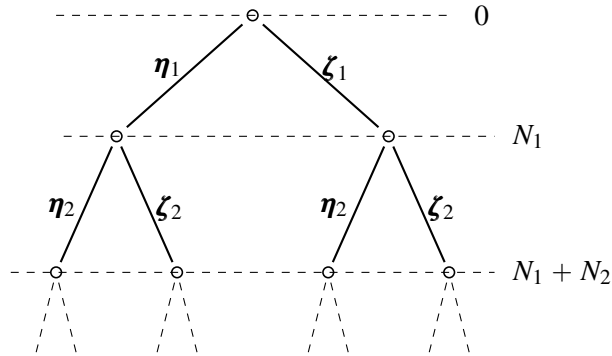
\section{Encoding certain sets in $\mathbb R$ by compressed trees}  \label{section: compressed trees}	
While the $M$-adic tree is canonical, it is sometimes too fine for our purposes. We therefore introduce more flexible compressed tree representations in which the relevant generations need not correspond to a single geometric scale. All vertices of an $M$-adic tree at generation $k$ represent intervals of the same length $M^{-k}$. In the compressed setting below, this uniformity will no longer be required.

%The $M$-adic tree $\mathcal T(E; M)$ is a natural representation of a set $E \subseteq [0,1]$, but not the only one. Alternative tree representations of $E$ also exist, and in fact play a key role in our analysis. All vertices of an $M$-adic tree of generation $k$ represent $M$-adic intervals of the same length $M^{-k}$. One can envision trees where this need not be the same. For instance, 
\subsection{A nested family of $M$-adic partitions} Let us consider a family of intervals $\{ \mathcal P_{\ell} :  \ell \geq 1  \}$ satisfying the properties below. The interval family $\mathcal P_{\ell}$ should be thought of as the collection of vertices corresponding to the $\ell^{\text{th}}$ generation of the compressed tree for $[0,1]$. 
\vskip0.1in
\begin{itemize}
\item {\em{Partitioning of $[0,1]$:}} Each $\mathcal P_{\ell}$ is a finite partition of $[0,1]$ into closed $M$-adic intervals with disjoint interiors.  
\vskip0.1in
\item {\em{Nesting of intersecting intervals: }} For any choice of indices $\ell > k \geq 1$ and any two intervals $I, I'$ with $I \in \mathcal P_{k}$, $I' \in \mathcal P_{\ell}$, one of the following two relations must hold: 
\[ \text{ either } I' \subseteq I \quad \text{ or } \quad \text{int}(I) \cap \text{int}(I') = \emptyset. \]   
In other words, an interval in $\mathcal P_{\ell}$ at a finer partition level $\ell$ lies inside the unique coarser interval of $\mathcal P_k$ whose interior it meets.
\vskip0.1in
\item {\em{Limiting fineness of the partitions:}} As $\ell$ increases, the mesh of the partition $\mathcal P_{\ell}$ becomes arbitrarily fine: 
\[ \text{diam}(\mathcal P_{\ell}) = \max \left\{ |I| : I \in \mathcal P \right\} \rightarrow 0 \text{ as } \ell \rightarrow \infty. \] 
\end{itemize}  
These three properties suffice to ensure that every $x \in [0,1]$ can be realized as the limit of a nested sequence of intervals from $\mathcal P_{\ell}$. For a number $x \in [0,1]$ that is not an $M$-adic rational, let $I_{\ell}(x)$ denote the unique interval of $\mathcal P_{\ell}$ containing $x$. The preceding assumptions on nesting and partitioning imply 
\begin{equation} \label{nesting-x-1} 
I_1(x) \supseteq I_2(x) \supseteq \cdots \supseteq I_{\ell}(x) \supseteq \cdots, \quad I_{\ell}(x) \in \mathcal P_{\ell}, \end{equation}  
whereas the vanishing mesh property implies 
\begin{equation} \bigcap_{\ell=1}^{\infty} I_{\ell}(x) = \{x\}. \label{nesting-x-2} \end{equation} 
For an $M$-adic rational number $x \in [0,1]$, there may be at most two choices of $I_{\ell}(x) \in \mathcal P_{\ell}$ containing $x$. This happens when $x$ is a boundary point of two of the intervals in $\mathcal P_{\ell}$. In such cases, we choose as $I_{\ell}(x)$ the interval for which $x$ is the left end point. This ensures uniqueness of the intervals $I_{\ell}(x)$, leading to properties \eqref{nesting-x-1} and \eqref{nesting-x-2} even for such $x$. 

\subsection{The compressed tree $\mathscr{P}$} Using the partitions $\{\mathcal P_{\ell} \}$, let us define a compressed tree representation $\mathscr{P} = \mathscr{P}[\mathcal P_1, \mathcal P_2, \ldots]$ of $[0,1]$ as follows: 
\vskip0.1in 
\begin{itemize}
\item A vertex of $\mathscr{P}$ of generation $\ell$ is an $\ell$-long nested chain $\langle I_1, \ldots, I_{\ell}\rangle$, with 
\[ I_{r+1} \subseteq I_{r}, \; 1 \leq r \leq \ell-1, \quad I_{r} \in \mathcal P_r \text{ for } 1 \leq r \leq \ell. \]
\vskip0.1in
\item Parent-child relations are defined by extension of nested chains, through set inclusion: we say that a vertex $\langle I_1, \ldots, I_{\ell}, I_{\ell+1}\rangle$ of the $(\ell+1)^{\text{th}}$ generation is a child of $\langle J_1, \ldots, J_{\ell}\rangle$ if 
\[ I_1 = J_1, \ldots, I_{\ell} = J_{\ell}, \quad  I_{\ell+1} \subseteq I_{\ell}. \] 
\end{itemize} 
\vskip0.1in
\noindent Given any set $E \subseteq [0,1]$, its tree $\mathscr{P}(E)$ is a subtree of $\mathscr{P}$ consisting of all maximal rays whose nested intervals converge to points of $E$; namely, a maximal ray of $\mathscr{P}(E)$ is of the form 
\[ \mathcal R: \langle I_1\rangle \rightarrow \langle I_1, I_2 \rangle \rightarrow \cdots \text{ such that } \bigcap_{\ell=1}^{\infty}I_{\ell} = \{x \} \text{ for some } x \in E.  \]
We will refer to $\mathscr{P}(E)$ as a {\em{compressed tree representation}} of $E$, aligned with the sequence of partitions $\{\mathcal P_{\ell} : \ell \geq 1\}$. Figure \ref{fig:binary-vs-compressed-tree} compares a compressed tree with $M = 2$ with the full binary tree. The depiction of a sample set $\Omega$ is also recorded in the two trees. 
\vskip0.1in
\noindent In the special case when $\mathcal P_{\ell}$ consists of all $M$-adic intervals of length $M^{-\ell}$, the tree $\mathscr{P}(E)$ reduces to the $M$-adic tree 
$\mathcal T(E; M)$. In general though, $\mathcal P_{\ell}$ could be a non-uniform partition of $[0,1]$. In Section \ref{section: compressions of the unit interval}, compressed tree representations of $[0,1]$ and $\Omega$ using partition scales of specific geometric significance will emerge as essential tools in the construction of Kakeya-type configurations. The selection of the scales aligned with the geometry of the slopes marks one of the main novelties of this article. 

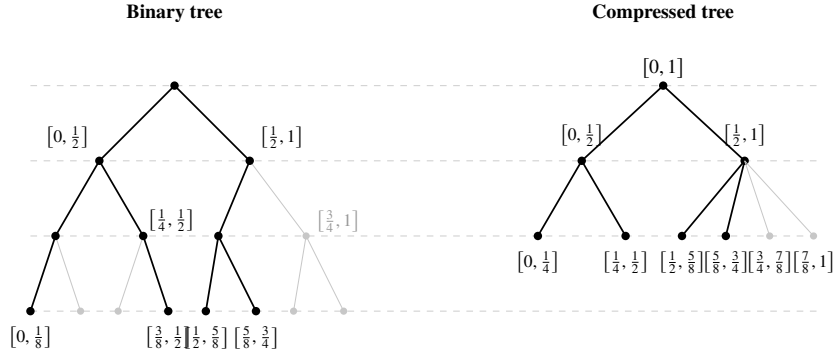
\begin{figure}[ht]
\centering
\resizebox{0.88\textwidth}{!}{%
\begin{tikzpicture}[
    x=1cm, y=1cm,
    edge/.style={line width=0.8pt},
    faded/.style={draw=black!22},
    guide/.style={draw=black!18, dashed, line width=0.5pt},
    nodepoint/.style={circle, fill=black, inner sep=1.35pt},
    fadednode/.style={circle, fill=black!22, inner sep=1.15pt},
    lab/.style={font=\footnotesize},
    intlab/.style={font=\scriptsize}
]

% =========================================================
% TITLES
% =========================================================
\node[lab] at (-4.1,1.15) {\textbf{Binary tree}};
\node[lab] at (3.7,1.15) {\textbf{Compressed tree}};

% =========================================================
% SHARED HORIZONTAL GUIDE LINES
% =========================================================
\draw[guide] (-6.4, 0.00) -- ( 6.2, 0.00);
\draw[guide] (-6.4,-1.20) -- ( 6.2,-1.20);
\draw[guide] (-6.4,-2.40) -- ( 6.2,-2.40);
\draw[guide] (-6.4,-3.60) -- ( 6.2,-3.60);

% =========================================================
% LEFT PANEL: FULL BINARY TREE
% =========================================================

% coordinates
\coordinate (r)   at (-4.1,0.00);

\coordinate (a1)  at (-5.3,-1.20);
\coordinate (a2)  at (-2.9,-1.20);

\coordinate (b1)  at (-6.0,-2.40);
\coordinate (b2)  at (-4.6,-2.40);
\coordinate (b3)  at (-3.4,-2.40);
\coordinate (b4)  at (-2.0,-2.40);

\coordinate (c1)  at (-6.4,-3.60);
\coordinate (c2)  at (-5.6,-3.60);
\coordinate (c3)  at (-5.0,-3.60);
\coordinate (c4)  at (-4.2,-3.60);
\coordinate (c5)  at (-3.6,-3.60);
\coordinate (c6)  at (-2.8,-3.60);
\coordinate (c7)  at (-2.2,-3.60);
\coordinate (c8)  at (-1.4,-3.60);

% edges
\draw[edge]  (r)--(a1);
\draw[edge]  (r)--(a2);

\draw[edge]  (a1)--(b1);
\draw[edge]  (a1)--(b2);
\draw[edge]  (a2)--(b3);
\draw[faded] (a2)--(b4);

\draw[edge]  (b1)--(c1);
\draw[faded] (b1)--(c2);
\draw[faded] (b2)--(c3);
\draw[edge]  (b2)--(c4);

\draw[edge]  (b3)--(c5);
\draw[edge]  (b3)--(c6);
\draw[faded] (b4)--(c7);
\draw[faded] (b4)--(c8);

% highlighted internal nodes and active leaves
\foreach \p in {r,a1,a2,b1,b2,b3,c1,c4,c5,c6}
    \node[nodepoint] at (\p) {};

% faded inactive nodes
\foreach \p in {b4,c2,c3,c7,c8}
    \node[fadednode] at (\p) {};

% NEW: labels at second level
\node[intlab, above left=2pt and 1pt] at (a1) {$\left[0,\tfrac12\right]$};
\node[intlab, above right=2pt and 1pt] at (a2) {$\left[\tfrac12,1\right]$};

% NEW: labels at third level
\node[intlab, above right =-2pt] at (b2) {$\left[\tfrac14,\tfrac12\right]$};
\node[intlab, above right=-2pt and 1pt, text=black!40] at (b4) {$\left[\tfrac34,1\right]$};

% selective labels only for highlighted leaves
\node[intlab, below=3pt] at (c1) {$\left[0,\tfrac18\right]$};
\node[intlab, below=3pt] at (c4) {$\left[\tfrac38,\tfrac12\right]$};
\node[intlab, below=3pt] at (c5) {$\left[\tfrac12,\tfrac58\right]$};
\node[intlab, below=3pt] at (c6) {$\left[\tfrac58,\tfrac34\right]$};

% =========================================================
% RIGHT PANEL: COMPRESSED TREE
% =========================================================

% root
\coordinate (R) at (3.7,0.00);
\node[nodepoint] at (R) {};
\node[lab] at (3.7,0.28) {$[0,1]$};

% first generation
\coordinate (L1) at (2.4,-1.20);
\coordinate (R1) at (5.0,-1.20);

\draw[edge] (R)--(L1);
\draw[edge] (R)--(R1);

\node[nodepoint] at (L1) {};
\node[nodepoint] at (R1) {};

\node[intlab, above=3pt] at (L1) {$\left[0,\tfrac12\right]$};
\node[intlab, above=3pt] at (R1) {$\left[\tfrac12,1\right]$};

% children of [0,1/2]
\coordinate (LL) at (1.7,-2.40);
\coordinate (LR) at (3.1,-2.40);

\draw[edge] (L1)--(LL);
\draw[edge] (L1)--(LR);

\node[nodepoint] at (LL) {};
\node[nodepoint] at (LR) {};

\node[intlab, below=3pt] at (LL) {$\left[0,\tfrac14\right]$};
\node[intlab, below=3pt] at (LR) {$\left[\tfrac14,\tfrac12\right]$};

% children of [1/2,1]
\coordinate (Q1) at (4.0,-2.40);   % [1/2,5/8]
\coordinate (Q2) at (4.7,-2.40);   % [5/8,6/8]
\coordinate (Q3) at (5.4,-2.40);   % [6/8,7/8]
\coordinate (Q4) at (6.1,-2.40);   % [7/8,1]

\draw[edge]  (R1)--(Q1);
\draw[edge]  (R1)--(Q2);
\draw[faded] (R1)--(Q3);
\draw[faded] (R1)--(Q4);

\node[nodepoint] at (Q1) {};
\node[nodepoint] at (Q2) {};
\node[fadednode] at (Q3) {};
\node[fadednode] at (Q4) {};

\node[intlab, below=3pt] at (Q1) {$\left[\tfrac12,\tfrac58\right]$};
\node[intlab, below=3pt] at (Q2) {$\left[\tfrac58,\tfrac34\right]$};
\node[intlab, below=3pt] at (Q3) {$\left[\tfrac34,\tfrac78\right]$};
\node[intlab, below=3pt] at (Q4) {$\left[\tfrac78,1\right]$};

\end{tikzpicture}%
}
\caption{\small{A schematic comparison between a full binary tree of height 3 on the left and its compressed counterpart of height 2  corresponding to 
$ \mathcal P_1 = \{ [0, \frac{1}{2}], [\frac{1}{2}, 1]\}$, $\mathcal P_2 = \mathcal P_1/2 \cup \left[ {1}/{2} + {\mathcal P_1}/{4} \right] 
$ on the right. On the left, only selected dyadic intervals are labelled. If $E = \{ e_i : i = 1, 2, 3, 4\}$ is a 4-point set with $e_1 \in [0, \frac{1}{8}]$, $e_2 \in [\frac{3}{8}, \frac{1}{2}]$, $e_3 \in [\frac{1}{2}, \frac{5}{8} ]$, $e_4 \in [\frac{5}{8}, \frac{3}{4} ]$, then the highlighted rays on either side constitute the truncated subtrees of $\mathcal T(E;2)$ and $\mathscr P(E)$ respectively.}}
\label{fig:binary-vs-compressed-tree}
\end{figure}

\section{Examples of compressed trees} \label{examples of compressed trees section}
We conclude this section with examples comparing two tree representations of the same set: one the ambient $M$-adic tree and the other a compressed tree adapted to the scales where branching occurs. The purpose of these examples is to demonstrate compressed trees that record only branching events rather than the full $M$-adic hierarchy. 
\subsection{The set $\Omega_{\text{HRS}}$} A natural test case for this comparison is the set $\Omega_{\text{HRS}}(R)$ defined in \eqref{HRS-example-unit}, whose dyadic structure is organized into $R$ successive binary blocks.
\vskip0.1in
\noindent The ambient dyadic tree for $\Omega_{\text{HRS}}$, namely  
\[ \mathcal T_{\text{HRS}}(R) := \mathcal T\bigl(\Omega_{\text{HRS}}(R); 2\bigr) \]
has the following characteristic features:
\vskip0.1in
\begin{itemize}  
\item Every ray of  $\mathcal T_{\text{HRS}}(R)$ becomes non-splitting after height 
\[ \bar{N}_R := N_1 + \ldots + N_R. \] 
Thus all branching occurs below height $\bar{N}_R$, and the relevant combinatorial structure is already visible in the truncation of $\mathcal T_{\text{HRS}}(R)$ at this level. 
%The splitting behaviour of $ \mathcal T_{\text{HRS}}(R)$ is therefore captured by its truncation to height $\bar{N}_R$. 
\vskip0.1in
\item Any ray of $ \mathcal T_{\text{HRS}}(R)$ from the root to height $\bar{N}_R$ contains exactly $R$ splitting vertices in total. These occur precisely at the block boundaries, i.e., at heights 
\[h \in \{ 0, \bar{N}_1, \bar{N}_2 \ldots, \bar{N}_{R-1} \}. \]
\vskip0.1in
\item Between two consecutive splitting heights $h<h'$, the tree $ \mathcal T_{\text{HRS}}(R)$ does not branch, but the two surviving rays evolve in opposite deterministic directions. If $v$ is the splitting vertex at height $h$, then its left child continues through a chain of forced right moves until height $h'$, while the right child continues through a chain of forced left moves up to the same height. Equivalently, each split at height $\bar{N}_{j-1}$ is followed by a long non-branching corridor of length $N_j$ whose composition is determined by the chosen block $\pmb{\eta}_j$ or $\pmb{\zeta}_j$.
%The non-splitting segment of the ray $\mathcal R$ between two consecutive splitting heights $h < h'$ has a special structure. If $v \in \mathcal R$ is the splitting vertex at height $h$, then the left child of vertex $v$ only generates right descendants up until height $h'$, whereas the right child generates left descendants only up to the same height. Visually, this means the left (respectively right) ray from $v$ keeps moving to the right (respectively left) until the next splitting height.     
\end{itemize}   
\vskip0.1in
\noindent We now compress this picture by retaining only the block-boundary scales $\bar{N}_1, \ldots, \bar{N}_R$ where branching occurs. 
%In contrast, let us now define a compressed tree representation for $\Omega_{\text{HRS}}(R)$, focusing only up to the dyadic scale $\bar{N}_R$.  
The resulting compressed tree $\mathscr{P}_{\text{HRS}}(R)$ has height $R$; its $\ell^{\text{th}}$ height records only the dyadic information visible at scale $2^{-\bar{N}_{\ell}}$, that is, with respect to the uniform partition 
\begin{equation}  \mathcal P_{\ell} \text{ of $[0,1]$ into dyadic intervals of length $2^{-\bar{N}_{\ell}}$. } \label{uniform partition} \end{equation}  
More explicitly, a vertex of generation $\ell$ in $\mathscr{P}_{\text{HRS}}(R)$ is the (necessarily unique) dyadic interval of length $2^{-\bar{N}_{\ell}}$ that contains an element of $\Omega_{\text{HRS}}(R)$. After all intermediate non-splitting levels are suppressed, the remaining tree is simply a full binary tree of height $R$. Thus the compressed tree records only the $R$ genuine branching events of $ \mathcal T_{\text{HRS}}(R)$, while discarding the long deterministic stretches between them.
%Written in this compressed form, $\mathscr{P}_{\text{HRS}}(R)$ is isomorphic to a full binary tree of height $R$. 
\vskip0.1in
\noindent Figure \ref{fig:HRS-3} provides a three-part depiction of the $\Omega_{\text{HRS}}(3)$: one using shrinking dyadic intervals at the block-boundary scales, the second as a dyadic tree with long non-splitting corridors, and the third as a compressed tree recording only the branching history.  
\begin{figure}[ht]
\centering
\begin{tikzpicture}[x=1cm,y=1cm,>=stealth]

% -------------------------------------------------
% Styles
% -------------------------------------------------
\tikzset{
  active/.style={draw=black, very thick},
  inactive/.style={draw=gray!35, thin},
  guide/.style={draw=gray!60, dashed, thin},
  trans/.style={draw=gray!60, thick, ->},
  pt/.style={circle,fill=black,inner sep=1.15pt},
  lab/.style={black},
  smalllab/.style={font=\small},
}

% =================================================
% PANEL A: interval picture
% =================================================
\def\XA{0}
\def\Awidth{5.6}
\pgfmathsetmacro{\midA}{\XA+\Awidth/2}

\def\yA0{0}
\def\yAone{-1.0}
\def\yAtwo{-2.0}
\def\yAthree{-3.0}

\def\a{1.15}
\def\b{0.36}
\def\c{0.11}

\node[lab] at (\midA,0.60) {\textbf{Intervals}};

% Top interval [0,1]
\draw[inactive] (\XA,\yA0) -- (\XA+\Awidth,\yA0);
\draw[inactive] (\XA,\yA0+0.07) -- (\XA,\yA0-0.07);
\draw[inactive] (\XA+\Awidth,\yA0+0.07) -- (\XA+\Awidth,\yA0-0.07);
\node[smalllab, below] at (\XA,\yA0) {$0$};
\node[smalllab, below] at (\XA+\Awidth,\yA0) {$1$};

% Midpoint guide
\draw[guide] (\midA,\yA0+0.16) -- (\midA,\yAthree-0.18);
\node[smalllab, below] at (\midA,\yA0) {$\frac12$};

% Level N1
\draw[active] (\midA-\a,\yAone) -- (\midA,\yAone);
\draw[active] (\midA,\yAone) -- (\midA+\a,\yAone);
\foreach \x in {\midA-\a,\midA,\midA+\a}{
  \draw[active] (\x,\yAone+0.06) -- (\x,\yAone-0.06);
}

% Braces for N1 interval lengths
\draw[decorate,decoration={brace,amplitude=4pt}]
(\midA-\a,\yAone+0.20) -- (\midA,\yAone+0.20)
node[midway,above=4pt,smalllab] {$2^{-\bar N_1}$};
\draw[decorate,decoration={brace,amplitude=4pt}]
(\midA,\yAone+0.20) -- (\midA+\a,\yAone+0.20);

% Level N2
\pgfmathsetmacro{\mL}{\midA-\a/2}
\pgfmathsetmacro{\mR}{\midA+\a/2}
\draw[active] (\mL-\b,\yAtwo) -- (\mL,\yAtwo);
\draw[active] (\mL,\yAtwo) -- (\mL+\b,\yAtwo);
\draw[active] (\mR-\b,\yAtwo) -- (\mR,\yAtwo);
\draw[active] (\mR,\yAtwo) -- (\mR+\b,\yAtwo);
\foreach \x in {\mL-\b,\mL,\mL+\b,\mR-\b,\mR,\mR+\b}{
  \draw[active] (\x,\yAtwo+0.05) -- (\x,\yAtwo-0.05);
}

% Brace for N2 interval length
\draw[decorate,decoration={brace,amplitude=3.5pt}]
(\mL-\b,\yAtwo+0.18) -- (\mL,\yAtwo+0.18)
node[midway,above=4pt,smalllab] {$2^{-\bar N_2}$};

% Level N3
\pgfmathsetmacro{\mLL}{\mL-\b/2}
\pgfmathsetmacro{\mLR}{\mL+\b/2}
\pgfmathsetmacro{\mRL}{\mR-\b/2}
\pgfmathsetmacro{\mRR}{\mR+\b/2}
\foreach \m in {\mLL,\mLR,\mRL,\mRR}{
  \draw[active] (\m-\c,\yAthree) -- (\m,\yAthree);
  \draw[active] (\m,\yAthree) -- (\m+\c,\yAthree);
  \foreach \x in {\m-\c,\m,\m+\c}{
    \draw[active] (\x,\yAthree+0.04) -- (\x,\yAthree-0.04);
  }
}

% Brace for N3 interval length
\draw[decorate,decoration={brace,amplitude=3pt}]
(\mLL-\c,\yAthree+0.16) -- (\mLL,\yAthree+0.16)
node[midway,above=3pt,smalllab] {$2^{-\bar N_3}$};

% Transition arrow with more space below panel A
\draw[trans] (\midA,\yAthree-0.45) -- (\midA,\yAthree-1.28);

% =================================================
% PANEL B: ambient dyadic tree
% Corrected and de-overlapped
% =================================================
\def\yBtop{-5.15}

\node[lab] at (\midA,\yBtop+0.55) {\textbf{Ambient dyadic tree}};

% Root
\coordinate (R0) at (\midA,\yBtop);

% Splitting vertices at height N̄1
\coordinate (S1L) at (\midA-1.05,\yBtop-0.95);
\coordinate (S1R) at (\midA+1.05,\yBtop-0.95);

\draw[active] (R0) -- (S1L);
\draw[active] (R0) -- (S1R);

\node[pt] at (R0) {};
\node[pt] at (S1L) {};
\node[pt] at (S1R) {};

% Immediate split from the two N̄1 vertices
% Spread these more to avoid visual overlap
\coordinate (A1) at (\midA-1.95,\yBtop-1.55);
\coordinate (A2) at (\midA-0.70,\yBtop-1.55);
\coordinate (A3) at (\midA+0.70,\yBtop-1.55);
\coordinate (A4) at (\midA+1.95,\yBtop-1.55);

\draw[active] (S1L) -- (A1);
\draw[active] (S1L) -- (A2);
\draw[active] (S1R) -- (A3);
\draw[active] (S1R) -- (A4);

% Continue without branching to N̄2
\coordinate (S2LL) at (\midA-1.95,\yBtop-2.40);
\coordinate (S2LR) at (\midA-0.70,\yBtop-2.40);
\coordinate (S2RL) at (\midA+0.70,\yBtop-2.40);
\coordinate (S2RR) at (\midA+1.95,\yBtop-2.40);

\draw[inactive] (A1) -- (S2LL);
\draw[inactive] (A2) -- (S2LR);
\draw[inactive] (A3) -- (S2RL);
\draw[inactive] (A4) -- (S2RR);

\node[pt] at (S2LL) {};
\node[pt] at (S2LR) {};
\node[pt] at (S2RL) {};
\node[pt] at (S2RR) {};

% Immediate split from the N̄2 vertices
% Also spread these more
\coordinate (B1) at (\midA-2.55,\yBtop-3.05);
\coordinate (B2) at (\midA-1.35,\yBtop-3.05);

\coordinate (B3) at (\midA-1.20,\yBtop-3.05);
\coordinate (B4) at (\midA-0.20,\yBtop-3.05);

\coordinate (B5) at (\midA+0.20,\yBtop-3.05);
\coordinate (B6) at (\midA+1.20,\yBtop-3.05);

\coordinate (B7) at (\midA+1.35,\yBtop-3.05);
\coordinate (B8) at (\midA+2.55,\yBtop-3.05);

\draw[active] (S2LL) -- (B1);
\draw[active] (S2LL) -- (B2);

\draw[active] (S2LR) -- (B3);
\draw[active] (S2LR) -- (B4);

\draw[active] (S2RL) -- (B5);
\draw[active] (S2RL) -- (B6);

\draw[active] (S2RR) -- (B7);
\draw[active] (S2RR) -- (B8);

% Continue without branching to N̄3
\coordinate (L1) at (\midA-2.55,\yBtop-3.90);
\coordinate (L2) at (\midA-1.35,\yBtop-3.90);

\coordinate (L3) at (\midA-1.20,\yBtop-3.90);
\coordinate (L4) at (\midA-0.20,\yBtop-3.90);

\coordinate (L5) at (\midA+0.20,\yBtop-3.90);
\coordinate (L6) at (\midA+1.20,\yBtop-3.90);

\coordinate (L7) at (\midA+1.35,\yBtop-3.90);
\coordinate (L8) at (\midA+2.55,\yBtop-3.90);

\draw[inactive] (B1) -- (L1);
\draw[inactive] (B2) -- (L2);
\draw[inactive] (B3) -- (L3);
\draw[inactive] (B4) -- (L4);
\draw[inactive] (B5) -- (L5);
\draw[inactive] (B6) -- (L6);
\draw[inactive] (B7) -- (L7);
\draw[inactive] (B8) -- (L8);

\node[pt] at (L1) {};
\node[pt] at (L2) {};
\node[pt] at (L3) {};
\node[pt] at (L4) {};
\node[pt] at (L5) {};
\node[pt] at (L6) {};
\node[pt] at (L7) {};
\node[pt] at (L8) {};

% Guides and labels
\draw[guide] (\XA+0.15,\yBtop) -- (\XA+\Awidth-0.15,\yBtop);
\draw[guide] (\XA+0.15,\yBtop-0.95) -- (\XA+\Awidth-0.15,\yBtop-0.95);
\draw[guide] (\XA+0.15,\yBtop-2.40) -- (\XA+\Awidth-0.15,\yBtop-2.40);
\draw[guide] (\XA+0.15,\yBtop-3.90) -- (\XA+\Awidth-0.15,\yBtop-3.90);

\node[left, smalllab] at (\XA+0.02,\yBtop) {$0$};
\node[left, smalllab] at (\XA+0.02,\yBtop-0.95) {$\bar N_1$};
\node[left, smalllab] at (\XA+0.02,\yBtop-2.40) {$\bar N_2$};
\node[left, smalllab] at (\XA+0.02,\yBtop-3.90) {$\bar N_3$};

\node[smalllab, gray!70] at (\midA,\yBtop-1.95) {corridors};
\node[smalllab, gray!70] at (\midA,\yBtop-3.48) {corridors};

% Transition arrow
\draw[trans] (\midA,\yBtop-4.35) -- (\midA,\yBtop-4.95);

% =================================================
% PANEL C: compressed tree
% =================================================
\def\yCtop{-10.90}

\node[lab] at (\midA,\yCtop+0.55) {\textbf{Compressed tree}};

\coordinate (Q0)  at (\midA,\yCtop);
\coordinate (Q1L) at (\midA-0.90,\yCtop-0.95);
\coordinate (Q1R) at (\midA+0.90,\yCtop-0.95);

\draw[active] (Q0) -- (Q1L);
\draw[active] (Q0) -- (Q1R);

\node[pt] at (Q0) {};
\node[pt] at (Q1L) {};
\node[pt] at (Q1R) {};

\coordinate (Q2LL) at (\midA-1.35,\yCtop-1.95);
\coordinate (Q2LR) at (\midA-0.45,\yCtop-1.95);
\coordinate (Q2RL) at (\midA+0.45,\yCtop-1.95);
\coordinate (Q2RR) at (\midA+1.35,\yCtop-1.95);

\draw[active] (Q1L) -- (Q2LL);
\draw[active] (Q1L) -- (Q2LR);
\draw[active] (Q1R) -- (Q2RL);
\draw[active] (Q1R) -- (Q2RR);

\node[pt] at (Q2LL) {};
\node[pt] at (Q2LR) {};
\node[pt] at (Q2RL) {};
\node[pt] at (Q2RR) {};

\coordinate (Q3a) at (\midA-1.75,\yCtop-2.95);
\coordinate (Q3b) at (\midA-1.05,\yCtop-2.95);
\coordinate (Q3c) at (\midA-0.68,\yCtop-2.95);
\coordinate (Q3d) at (\midA-0.20,\yCtop-2.95);
\coordinate (Q3e) at (\midA+0.20,\yCtop-2.95);
\coordinate (Q3f) at (\midA+0.68,\yCtop-2.95);
\coordinate (Q3g) at (\midA+1.05,\yCtop-2.95);
\coordinate (Q3h) at (\midA+1.75,\yCtop-2.95);

\draw[active] (Q2LL) -- (Q3a);
\draw[active] (Q2LL) -- (Q3b);
\draw[active] (Q2LR) -- (Q3c);
\draw[active] (Q2LR) -- (Q3d);
\draw[active] (Q2RL) -- (Q3e);
\draw[active] (Q2RL) -- (Q3f);
\draw[active] (Q2RR) -- (Q3g);
\draw[active] (Q2RR) -- (Q3h);

\node[pt] at (Q3a) {};
\node[pt] at (Q3b) {};
\node[pt] at (Q3c) {};
\node[pt] at (Q3d) {};
\node[pt] at (Q3e) {};
\node[pt] at (Q3f) {};
\node[pt] at (Q3g) {};
\node[pt] at (Q3h) {};

\node[left, smalllab] at (\XA+0.02,\yCtop) {$0$};
\node[left, smalllab] at (\XA+0.02,\yCtop-0.95) {$1$};
\node[left, smalllab] at (\XA+0.02,\yCtop-1.95) {$2$};
\node[left, smalllab] at (\XA+0.02,\yCtop-2.95) {$3$};

\end{tikzpicture}
\caption{\small{Three views of the model set $\Omega_{\mathrm{HRS}}(3)$. Top: the dyadic intervals at the block-boundary scales $2^{-\bar{N}_j}$ containing the elements of $\Omega_{\mathrm{HRS}}(3)$. Middle: the ambient dyadic tree, where branching occurs at the distinguished heights $\bar N_1,\bar N_2$, and the intervening rays continue without further branching. Bottom: the compressed tree obtained by suppressing those non-splitting corridors and retaining only the genuine branching events.}}
\label{fig:HRS-3}
\end{figure}
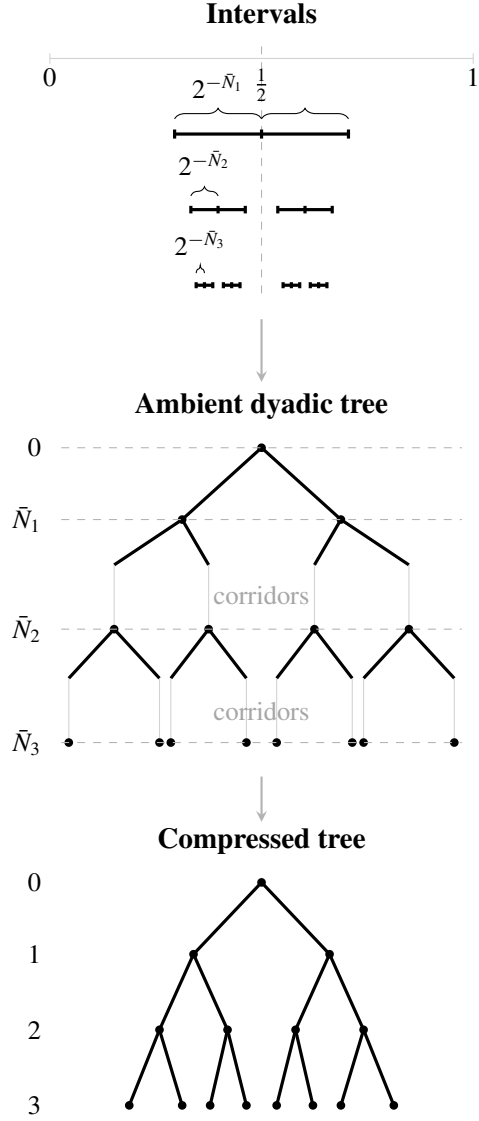

\subsection{A non-uniform variant of $\Omega_{\text{HRS}}$}
The example  $\Omega_{\text{HRS}}(R)$ in \eqref{HRS-example-unit} is unusually uniform.  Every ray in $\partial \mathcal T_{\text{HRS}}(R)$ splits at the same sequence of dyadic heights.  Accordingly, its compressed representation $\mathscr{P}_{\text{HRS}}(R)$ can be built from uniform partitions at each generation, 
%This is reflected in the uniformity of the partition $\mathcal P_{\ell}$ in the compressed tree , 
as we see from \eqref{uniform partition}.  
\vskip0.1in
\noindent More generally, however, different branches may split at different dyadic heights. It is not difficult to envision an $M$-adic tree where one branch splits at height $N_1$ for the first time, and another does so at height $N_1' \neq N_1$. In that setting, no single uniform partition can capture all branching events simultaneously. This means that a useful compressed tree can no longer be based on a single common scale at each generation; the compression must adapt to the branch. 
\vskip0.1in
\noindent To illustrate this, let us modify the example in \eqref{HRS-example-unit}, by allowing the two possible blocks $\pmb{\eta}_j$ and ${\pmb{\zeta}}_j$ at stage $j$ to have different lengths. More precisely, suppose that $\{p_j, q_j : 1 \leq j \leq R\}$ is a set of integers with the following properties: 
\begin{align} 
&p_j < q_j \text{ for } 1 \leq j \leq R,  \quad \sum_{k=1}^{j} q_k < \sum_{k=1}^{j+1} p_k \text{ for } 1 \leq j \leq R-1, \label{pq-1} \\  
&\; \#(\mathscr{H}_j) = 2^j \text{ where } \mathscr{H}_j := \Bigl\{\sum_{k=1}^{j} a_k \; \bigl| \; a_k \in \{p_k, q_k\}, 1 \leq k \leq j \Bigr\}. \label{pq-2}
\end{align}
The two inequalities in \eqref{pq-1} jointly imply that any $j$-fold sum, where the $k^{\text{th}}$ summand is either $p_k$ or $q_k$, is strictly less than any $(j+1)$-fold sum of the same type. Condition \eqref{pq-2} implies that all such $j$-fold sums, numbering $2^j$, are distinct.
\vskip0.1in
\noindent Equipped with these integers, we now define binary strings 
\begin{equation} 
\pmb{\eta}_j = (0, 1, \ldots, 1) \in \{0, 1\}^{p_j}, \qquad \pmb{\zeta}_j = (1, 0, \ldots, 0) \in \{0, 1\}^{q_j}.  
\end{equation}  
The modified set $\Omega(R)$ is given by 
\begin{equation} \label{nonuniform HRS}
\Omega(R) := \left\{ \sum_{k=1}^{\infty} \frac{\varepsilon_k}{2^k} \; \Bigl| \;\begin{aligned} &\pmb{\varepsilon} = (\varepsilon_1, \varepsilon_2, \ldots) = (\pmb{\kappa}_1, \ldots, \pmb{\kappa}_R, 0, 0, \ldots), \\ & \text{where } \pmb{\kappa}_j \text{ is either } \pmb{\eta}_j  \text{ or } \pmb{\zeta}_j \text{ for } 1 \leq j \leq R  \end{aligned} \right\}  
\end{equation} 
Let us compare and contrast this situation with the previous example. 
\vskip0.1in
\begin{itemize} 
\item Each maximal ray in the dyadic tree $\mathcal T(\Omega(R); 2)$ splits exactly $R$ times. 
\vskip0.1in
\item However, the rays now split at different heights, with the $j^{\text{th}}$ split on any ray occurring only at a height lying in $\mathscr{H}_j$.  
\vskip0.1in
\item Condition \eqref{pq-2} ensures that each of the $2^j$ splitting vertices of $j^{\text{th}}$ split lies at distinct heights. 
\vskip0.1in 
\item The assumption \eqref{pq-1} ensures that the height of any $j^{\text{th}}$ splitting vertex  of $\mathcal T(\Omega(R); 2)$ is strictly less than that of any $(j+1)^{\text{th}}$ splitting vertex.
\vskip0.1in 
\item All splits are exhausted by height $q_1 + \ldots + q_R$, so rays originating from this height are non-splitting.  
\end{itemize} 
\vskip0.1in
To summarize, the dyadic tree for $\Omega(R)$ still has binary branching at $R$ successive stages, but the locations of those stages now depend on the branch. We will now express $\Omega(R)$ as a compressed tree that encodes this branching pattern by forcing all the $\ell^{\text{th}}$ splitting vertices on each ray of the dyadic tree to lie in the $\ell^{\text{th}}$ generation. This of course means that vertices of the same generation can no longer be uniformly sized.
\vskip0.1in
\noindent To make this precise, let us first define a compression $\mathscr{P}([0,1])$ of the full binary tree adapted to the branching of $\Omega(R)$, then use it to study the  compressed tree $\mathscr{P}(R)$ that depicts the dyadic tree $\mathcal T(\Omega(R); 2)$, truncated to height $q_1 + \ldots + q_R$. 
\vskip0.1in
\begin{itemize}
\item The first generation already reflects the asymmetry of the construction. The partition $\mathcal P_1$ consists of two sub-collections of intervals:
\[ \mathcal P_1 = \mathcal I_0 \sqcup \mathcal I_1. \] 
The collection $\mathcal I_0$ contains dyadic intervals in $[0, \frac{1}{2}]$ of length $2^{-p_1}$, while $\mathcal I_1$ contains dyadic intervals in $[\frac{1}{2}, 1]$ of length $2^{-q_1}$. The tree $\mathscr{P}(R)$ has exactly two first generation vertices, one from $\mathcal I_0$ and the other from $\mathcal I_1$.  Thus the two children of the root live at different dyadic scales.
\vskip0.1in
\item This difference in lengths persists at subsequent heights. For $2 \leq j \leq R$, the partition $\mathcal P_j$ is the union of $2^j$ sub-collections, indexed by the binary string ${\pmb{\varepsilon}} = (\varepsilon_1, \ldots, \varepsilon_j) \in \{0, 1\}^j$. More concretely, 
\begin{equation}  
\mathcal P_j = \bigsqcup \Bigl\{ \mathcal I_{\pmb{\varepsilon}} : \pmb{\varepsilon} = (\varepsilon_1, \ldots, \varepsilon_j) \in \{0, 1\}^j  \Bigr\}. \label{partition j} \end{equation}  
For $\pmb{\varepsilon}' = (\varepsilon_1, \ldots, \varepsilon_{j-1})$, an interval $I' \in \mathcal I_{\pmb{\varepsilon}'}$  is partitioned  into sub-intervals $I \in \mathcal I_{\pmb{\varepsilon}}$, with 
\[ 
|I| = 2^{-\ell_j(\pmb{\varepsilon})} \; \text{ so that } \;\#(\mathcal I_{\pmb{\varepsilon}}) = 2^{\ell_j(\pmb{\varepsilon}) - \ell_j(\pmb{\varepsilon'})}. \] The length of each child $I$ of $I'$ is determined by the cumulative branch-dependent depth: 
\begin{align*} \ell_j(\pmb{\varepsilon}) &= \sum_{k=1}^j \Bigl[ (1-\varepsilon_k) p_k + \varepsilon_k q_k \Bigr] \\
&= \ell_{j-1}(\pmb{\varepsilon}') + \begin{cases} p_j &\text{ if } \varepsilon_j = 0, \\  q_j &\text{ if } \varepsilon_j = 1. \end{cases} \end{align*}
In other words, the compressed tree $\mathscr{P}([0,1])$ generated by the partitions \eqref{partition j} has the property that the left (respectively right) half of every $(j-1)^{\text{th}}$ level interval $I' \in \mathcal I_{\pmb{\varepsilon}'}$ decomposes into $2^{p_j}$ (respectively $2^{q_j}$) equal pieces at the $j^{\text{th}}$ level. 
\vskip0.1in
\item In this formulation, the tree $\mathscr{P}(R)$ contains exactly $2^j$ vertices of the $j^{\text{th}}$ generation, one element from each $\mathcal I_{\pmb{\varepsilon}}$. Every vertex of $\mathscr{P}(R)$ of generation $0 \leq j < R$ has exactly two children.    
\vskip0.1in 
\item The compressed tree $\mathscr{P}(R)$ is a full binary tree of height $R$ and splitting number $R$. Thus, despite the non-uniform scales at a given height, the compressed tree $\mathscr{P}(R)$ again retains only the essential binary branching pattern. 
\end{itemize} 
\vskip0.1in
\noindent These examples foreshadow the compressed root and slope trees of Section \ref{section: compressions of the unit interval}, where the relevant generations are chosen by geometric significance rather than by a fixed ambient dyadic scale. 

\chapter{Lacunarity order and $M$-adic splitting number} \label{section: split implies lacunarity} 
\section{Chapter overview} 
The goal of this chapter is to establish a precise connection between the notion of lacunarity order of a set, introduced in Chapter \ref{section: finite-order lacunarity}, and the combinatorial structure of its $M$-adic tree.
\vskip0.1in
\noindent The key invariant that mediates this connection is the \emph{splitting number} of the tree, as defined in Section \ref{section: splitting number}. We have seen that, informally, the splitting number measures the maximal number of branching events encountered along any ray of the tree. In doing so, it provides a quantitative description of how complex the associated set is from a combinatorial standpoint.
\vskip0.1in
\noindent A central theme of this chapter is that the splitting number captures exactly the recursive structure underlying finite-order lacunarity. {\em{A set has finite lacunarity order if and only if its associated $M$-adic tree has finite splitting number}}. This correspondence forms the backbone of the arguments developed later in the article. From a geometric perspective, each branching event corresponds to the introduction of a new scale of clustering. Controlling the number of such events ensures that the set exhibits only finitely many layers of hierarchical structure.
\vskip0.1in
\noindent The main result of this chapter, Proposition \ref{TREE-LACUNARY-TRADITIONAL}, makes one direction of this equivalence precise. This is the direction we need for proving the implication \eqref{condition 3: slopes sublacunary} $\implies$ \eqref{condition 1: Kakeya-type sets} in Theorem \ref{thm:main}. It shows that an $M$-adic tree with bounded splitting number corresponds to a set that can be decomposed into finitely many lacunary sets of finite order. The other direction of the  equivalence (finite order lacunarity $\implies$ finite splitting number) is relegated to Appendix D Chapter \ref{chapter: lac-split-appendix}.
\vskip0.1in
\noindent The proof of Proposition \ref{TREE-LACUNARY-TRADITIONAL} proceeds by induction on the splitting number, adapting the ideas developed in \cite{{Bateman},{KrocThesis}}. The argument is organized as follows. In Section \ref{section: unit split}, we analyse the base case of trees with unit splitting number, showing that they correspond to lacunary sets of order one. This section isolates a key structural lemma (Lemma \ref{BASE CASE}) that allows us to translate branching behaviour of an $M$-adic tree to lacunary clustering of its associated set. Lemma \ref{BASE CASE} is proved in Section \ref{proof-lemma-base-case}. Finally, in Section \ref{completion of proof: Tree-Lacunary-Traditional}, we complete the inductive step and prove Proposition \ref{TREE-LACUNARY-TRADITIONAL}.
%The connection between the order of lacunarity of a set and the splitting number of its representative $M$-adic tree originated in \cite{Bateman}. It was shown that one of these quantities is finite if and only the other one is. Following the ideas developed in \cite{{Bateman},{KrocThesis}}, we recast the concept of finite order lacunarity of a set of real numbers, as introduced in Section \ref{section: finite-order lacunarity}, using the structure of the splitting vertices of its tree. This is the content of Proposition \ref{TREE-LACUNARY-TRADITIONAL}, the main objective of this chapter. 
\section{Sublacunarity $\implies$ infinite splitting number} 
Let us recall the definition of finite-order lacunarity $\Lambda(N, \lambda; R)$ from Section \ref{basic definitions section}, and the $M$-adic tree representation $\mathcal T(U; M)$ of bounded subsets $U$ of the positive half line from Section \ref{section: trees for bounded sets}.
\begin{proposition}\label{TREE-LACUNARY-TRADITIONAL}
	For any $M\geq 2$, $N \geq 1$, there is a positive integer constant $R = R(N,M)$ with the following property.  
	\vskip0.1in
	\noindent Any bounded set $U \subseteq [0, \infty)$ with split$(\mathcal T(U;M)) = N$ must obey
	\[ U \in \Lambda(N, M^{-1}; R) \text{ in the sense of Definition \ref{defn: Admissible finite order lacunarity}}.  \]
In other words, any set $U$ whose $M$-adic tree has splitting number $N$ can be covered by the $R$-fold union of sets in $\Lambda(N;M^{-1})$, i.e. sets of lacunarity order at most $N$ and lacunarity constant at most $M^{-1}$, as described in Definition \ref{defn: Lacunary sets}. 
\end{proposition}
\noindent A by-product of Proposition \ref{TREE-LACUNARY-TRADITIONAL} is the following corollary.
\begin{corollary} \label{corollary: sublacunary slopes imply infinite split}
If $U \subseteq [0, \infty)$ is a bounded set, 
\begin{equation} \label{finite split implies finite order} \text{split}\bigl(\mathcal T(U; M)\bigr) < \infty \; \Longrightarrow \; U \in {\tt{AdFinLac}} \text{ as defined in Definition \ref{defn: Admissible finite order lacunarity}}.\end{equation} 
Equivalently, if $U$ is sublacunary (i.e. $U \in {\tt{SubLac}}$ according to Definition \ref{defn: Admissible finite order lacunarity}), then the splitting number of the tree $\mathcal T(U; M)$ is infinite. 
\end{corollary}
\vskip0.1in 
\noindent {\em{Remarks:}} \begin{enumerate}[1.] 
\item The converse of Proposition \ref{TREE-LACUNARY-TRADITIONAL}, 
\[ U \in {\tt{AdFinLac}} \; \Longrightarrow \; \text{split}(\mathcal T(U; M)) < \infty \] 
also holds. We do not require this implication for the proof Theorem \ref{thm:main}, but include the relevant statement in Section \ref{section: lacunarity implies split} for completeness. 
\vskip0.1in
\item The remainder of this chapter is devoted to the proof of Proposition \ref{TREE-LACUNARY-TRADITIONAL}, which proceeds by induction on $N$. Lemma \ref{BASE CASE} in the next section establishes the base case $N=1$; its proof sets the stage for the general argument. The inductive step is completed in Section \ref{completion of proof: Tree-Lacunary-Traditional}.
\end{enumerate} 
\section{Analysis of $M$-adic trees of unit splitting number} \label{section: unit split}
The goal of this section is to verify Proposition \ref{TREE-LACUNARY-TRADITIONAL} in the simplest case $N=1$. In this setting, we show that an $M$-adic tree with splitting number one corresponds to a highly rigid geometric structure: its associated set can be expressed as a finite union of lacunary sequences of order one.
%The goal of this section is to verify Proposition  in the simplest case of a non-trivial tree structure, when $N=1$. In this situation, all branching in the tree occurs along a single ray, while all other vertices lie along non-branching paths emanating from it. This rigid structure forces the associated set to exhibit strong clustering toward a single limit point identified by the special ray, with uniform geometric separation between successive elements. As we show below, this is precisely the structure of a finite union of lacunary sequences.
\begin{lemma} \label{BASE CASE}
Fix a base $M \geq 2$, and let $A \subseteq [0, \infty)$ be a bounded set with 
\begin{equation} \label{A split number 1}
\text{split}(\mathcal T(A;M)) = 1. 
\end{equation}  Then the conclusion of Proposition \ref{TREE-LACUNARY-TRADITIONAL} holds with $N=1$ and $R(1, M) = 6M$; namely, 
\[A \in \Lambda(1, M^{-1}; 6M) \text{ in the notation of Definition \ref{defn: Lacunary sets}}. \] 
More precisely, $A$ is contained in the union of at most $6M$ lacunary sequences in {\tt{MonLac}}$(M^{-1})$, as defined in Definition \ref{defn: Lacunary sequence in R}). 
\end{lemma}
\noindent The proof of Lemma \ref{BASE CASE} is presented in Section \ref{proof-lemma-base-case}, and closely follows the line of reasoning in \cite[Remark 2, page 60]{Bateman} and \cite[Lemma 3.8]{KrocThesis}. Here is a brief overview of the argument: 
\begin{itemize} 
\item The idea is to  first isolate a special ray of the tree $\mathcal T(A; M)$, denoted $\mathcal R^{\ast}$ in Section \ref{step 1 base case}, containing all the splitting vertices. This identifies a candidate for the special point of $A$, which we call $a^{\ast}$. 
\vskip0.05in 
\item The ray $\mathcal R^{\ast}$ helps in decomposing the tree $\mathcal T(A; M)$ into sub-trees $\mathcal T^{[i]}$, each representing a sequence $A^{[i]} = \{a_r : r\geq 1\}$ converging to $a^{\ast}$. This decomposition is carried out in Section \ref{Step 2: base case proof}. 
\vskip0.05in 
\item Additionally and crucially, a quantifiable Euclidean gap separates each pair of points $a_r, a_{r+3}$ in $A^{[i]}$; this feature is essential for establishing both monotonicity and lacunarity. We execute this step in Sections \ref{step 3 base case} and \ref{step 4 base case}, establishing that $A^{[i]}$ is contained in the three-fold union of sequences in ${\tt{MonLac}}(M^{-1})$. According to Lemma \ref{Lemma: Lacunary 1} \eqref{Lemma: Lacunary 1 (a)}, such sequences are members of $\Lambda(1, M^{-1})$. 
\end{itemize} 
\section{Proof of Lemma \ref{BASE CASE}} \label{proof-lemma-base-case}
Let us proceed to fill in the details. Let $A \subseteq [0, \infty)$ be a bounded set obeying \eqref{A split number 1}. Scaling $A$ by an appropriate power of $M$, we may assume without loss of generality that $A \subseteq [0, 1]$. 
\subsection{Step 1: Identification of a special ray and a lacunary limit.} \label{step 1 base case}
To reduce notational baggage, let us abbreviate $\mathcal T(A; M)$ as $\mathcal T$ in this proof. Our first claim is: 
\begin{equation} \label{what is R-star}
\exists \text{a ray $\mathcal R^{\ast} \in \partial \mathcal T$ containing all the splitting vertices of $\mathcal T$.}
\end{equation} 
Let us confirm this claim. If the set 
\[ \mathcal V = \{v \in \mathcal T \, : \,\text{split}_{\mathcal T}(v) = 1 \} \] is a singleton, its sole member has to be the root $v_0$. In this case, any maximal ray $\mathcal R^{\ast}$ originating at $v_0$ will do. If $\mathcal V$ has at least two elements, Lemma \ref{SplitsOnARay} guarantees the existence of a single ray $\mathcal R$, rooted at $v_0$ but not necessarily maximal, whose vertices are precisely all the splitting vertices of $\mathcal T$. In this case, let $\mathcal R^{\ast} \in \partial \mathcal T$ be a maximal ray containing $\mathcal R$. Unlike $\mathcal R$, which is unique, the ray $\mathcal R^{\ast}$ is non-unique in general; however, it is unique if $\mathcal R$ is of infinite length. This establishes the claim. 
\vskip0.1in
\noindent We have seen in Section \ref{tree encoding section} that any infinite ray of an $M$-adic tree identifies a unique point through a nested sequence of intervals given by its vertices. Thus the ray $\mathcal R^{\ast}$ ensured by \eqref{what is R-star} determines a distinguished point in the closure $\bar{A}$ of $A$. Let us call this point 
\[ a^{\ast} := \alpha(\mathcal R^{\ast}), \text{ in the notation of } \eqref{ray&point}. \] 
The aim is to find a cover of $\bar{A}$, and hence of $A$, consisting of $6M$ sequences in ${\tt{MonLac}}(M^{-1}) \subseteq \Lambda(1, M^{-1})$. All of these sequences will be shown to converge to $a^{\ast}$. Figure \ref{fig:split 1} shows a schematic diagram of the structure of $\mathcal T$. 
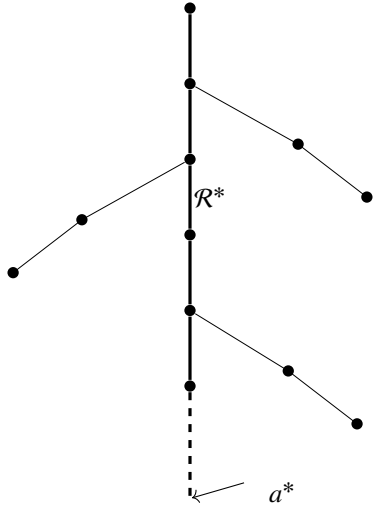
\begin{figure}[ht]
\centering
\begin{tikzpicture}[
    x=1.3cm, y=1cm,
    dot/.style={circle, fill=black, inner sep=1.5pt},
    every node/.style={font=\small}
]

% Main distinguished ray R^*
\node[dot] (v0) at (0,0) {};
\node[dot] (v1) at (0,-1) {};
\node[dot] (v2) at (0,-2) {};
\node[dot] (v3) at (0,-3) {};
\node[dot] (v4) at (0,-4) {};
\node[dot] (v5) at (0,-5) {};

\draw[very thick] (v0)--(v1)--(v2)--(v3)--(v4)--(v5);
\draw[very thick,dashed] (v5)--(0,-6.5);

% Side rays: three splitting vertices, one non-splitting vertex on R^*
\draw[very thin] (v1)--(1.1,-1.8) node[dot] {} -- (1.8,-2.5) node[dot] {};
\draw[very thin] (v2)--(-1.1,-2.8) node[dot] {} -- (-1.8,-3.5) node[dot] {};
% v3 intentionally does not split
\draw[very thin] (v4)--(1.0,-4.8) node[dot] {} -- (1.7,-5.5) node[dot] {};

% Labels
\node[left] at (0.48,-2.5) {$\mathcal{R^{*}}$};

\node[right] at (0.7,-6.4) {$a^{*}$};
\draw[->] (0.55,-6.28) -- (0.02,-6.48);

\end{tikzpicture}
\caption{\small{Step 1 of Lemma \ref{BASE CASE}: A schematic representation of a tree $\mathcal T$ with unit splitting number. The distinguished infinite ray $\mathcal R^{*}$ is indicated in bold. All splitting vertices lie on $\mathcal R^{*}$, while the rays emanating from it do not split further. The dotted continuation indicates that $\mathcal R^{*}$ is infinite and determines the limit point $a^{*}=\alpha(\mathcal R^{*})$.}}
\label{fig:split 1}
\end{figure}
\vskip0.1in 
\noindent We define two subsets $A_{\pm}$ of $\bar{A}$, containing respectively points to the right and to the left of $a^{\ast}$: 
\begin{equation} A_{+} := \bar{A} \cap [a^{\ast}, 1], \quad  A_{-} := \bar{A} \cap [0, a^{\ast}]. \label{A-plus-minus}\end{equation}  
 The conclusion of Lemma \ref{BASE CASE} is then a consequence of the following claim: 
 \begin{equation} \label{Lacunary cover}
 \left\{
 \begin{aligned} 
&\text{ Each of the sets $A_{\pm}$ is contained in the $3M$-fold union of } \\
&\text{ sequences in {\tt{MonLac}}$(M^{-1})$, each converging to $a^{\ast}$.}  
\end{aligned}
\right\} 
\end{equation} 
We will focus on proving the claim \eqref{Lacunary cover} for $A_{+}$ only. The proof for $A_{-}$ is identical. 
\subsection{Step 2: Decomposition of $\mathcal T(A_{+}; M)$} \label{Step 2: base case proof} 
Since $a^{\ast} \in A_{+}$, the tree 
\begin{equation} \label{T+}
\mathcal T_{+} := \mathcal T(A_{+}; M) \subseteq \mathcal T 
\end{equation}  contains the ray $\mathcal R^{\ast}$, and has splitting number at most 1. If the splitting number of $\mathcal T_{+}$ is zero, then it consists of the single ray $\mathcal R^{\ast}$, which means 
\[A_{+} = \{a_{+} \} \in \Lambda(0, M^{-1}) \subseteq \Lambda(1; M^{-1}). \] 
The last inclusion is a result of the monotonicity of lacunarity classes with respect to lacunarity order, as shown in Lemma \ref{Lemma : lacunarity monotonicity}. 
\vskip0.1in 
\noindent The non-trivial scenario, therefore, occurs when 
\[ \text{split}(\mathcal T_{+}) = 1.  \]
Since every splitting vertex of $\mathcal T_{+}$ is also a splitting vertex of $\mathcal T$, the ray $\mathcal R^{\ast}$ given by \eqref{what is R-star} contains all the splitting vertices of $\mathcal T_{+}$ as well. 
\vskip0.1in
\noindent A consequence of the above statement is that any ray of $\mathcal T_{+}$ that is not contained in $\mathcal R^{\ast}$ but is rooted at a vertex of $\mathcal R^{\ast}$ is non-splitting. More explicitly, suppose that the ray $\mathcal R^{\ast}$ connects the sequence of vertices
\[ \mathcal R^{\ast}:  v_0 = [0,1] \rightarrow v_1 \rightarrow \cdots \rightarrow v_j \rightarrow v_{j+1} \rightarrow \cdots,  \]
with $h(v_j) = j$ for every  $j = 0, 1, 2, \cdots$. Since a vertex of an $M$-adic tree can have at most $M$ children, each vertex $v_j$ is the root of at most $(M-1)$ rays in $\mathcal T_{+}$ that are not contained in $\mathcal R^{\ast}$. 
%Let $\mathfrak R(v_j)$ denote this collection:  
%\[ \mathfrak R(v_j) := \bigl\{ \mathcal R \text{ ray of } \mathcal T_{+} \text{ rooted at } v_j : \mathcal R \not\subseteq \mathcal R^{\ast}  %\bigr\}, \text{ so that } \#(\mathfrak R(v_j)) \leq M-1.  \]   
%The non-splitting nature of each ray in $\mathfrak R(v_j)$ 
To paraphrase, for $i \in \mathbb Z_M := \{0, 1, \ldots, M-1\}$ and $j \in \{0, 1, 2, \ldots \}$,
\begin{equation} \label{uniqueness of non-splitting ray}
\left\{ 
\begin{aligned} 
&\text{ there is at most one infinite ray $\mathcal R_{ij}$ of $\mathcal T_{+}$ such that } \\  
&\text{ $\mathcal R_{ij} \not\subseteq \mathcal R^{\ast}$, and $\mathcal R_{ij}$ is rooted at the $i^{\text{th}}$ child of $v_j$. }
\end{aligned} 
\right\}
\end{equation} 
This observation is key to decomposing the tree $\mathcal T(A; M)$ into a bounded number of lacunary components, in the following way:  
\begin{equation} \mathcal T_{+} = \bigcup_{i=0}^{M-1} \mathcal T^{[i]}.  \label{T+decomp}\end{equation}   
For each $i \in \mathbb Z_M$, the tree  $\mathcal T^{[i]}$ consists of the special ray $\mathcal R^{\ast}$ and the rays $\mathcal R_{ij}$ defined by \eqref{uniqueness of non-splitting ray}, for all indices $j \in \{0, 1, 2 \ldots \}$ where such a ray exists. Let $A^{[i]} \subseteq A_{+}$ denote the closed set represented by $\mathcal T^{[i]}$:
\begin{equation} \label{Ai+}
\mathcal T^{[i]}=: \mathcal T(A^{[i]}; M). 
\end{equation}  
Thus, in view of \eqref{T+}, \eqref{T+decomp} and \eqref{Ai+}, the tree decomposition of $\mathcal T_{+}$ leads to a set decomposition of $A_{+}$, namely 
\[ A_{+}= \bigcup_{i=0}^{M-1} A^{[i]}. \] 
Figure \ref{fig: step 2} depicts the construction of $A^{[i]}$ in the context of a trinary tree.  
\begin{figure}[ht]
\centering
\begin{tikzpicture}[x=1cm,y=1cm,every node/.style={font=\small}]

% Styles
\tikzset{
  dot/.style={circle, fill=black, inner sep=1.3pt},
  Adot/.style={circle, fill=blue, inner sep=1.3pt},
  Bdot/.style={circle, fill=red, inner sep=1.3pt}
}

% ---------------- Left panel: tree ----------------

% Special ray R^*
\node[dot] (v0) at (0,0) {};
\node[dot] (v1) at (0,-1) {};
\node[dot] (v2) at (0,-2) {};
\node[dot] (v3) at (0,-3) {};
\node[dot] (v4) at (0,-4) {};
\node[dot] (v5) at (0,-5) {};

\node[left] at (0.5,0) {$v_0$};
\node[left] at (0.5,-1) {$v_1$};
\node[left] at (0.5,-2) {$v_2$};
\node[left] at (0.5,-3) {$v_3$};

\draw[very thick] (v0)--(v1)--(v2)--(v3)--(v4)--(v5);
\draw[very thick,dashed] (v5)--(0,-6.3);

\node[left] at (0.2,-2.5) {$\mathcal R^{*}$};

% Label a^* at bottom of dotted ray
\node[right] at (0.2,-6.3) {$a^{*}$};

% Splitting vertices: v1, v2, v4; v3 is non-splitting

% From v1
\draw[blue, thick] (v1)--(-1.1,-1.7) node[Adot] {}
                  --(-1.8,-2.4) node[Adot] {};
\draw[red, thick]  (v1)--( 1.1,-1.7) node[Bdot] {}
                  --( 1.8,-2.4) node[Bdot] {};

% From v2
\draw[blue, thick] (v2)--(-1.1,-2.7) node[Adot] {}
                  --(-1.8,-3.4) node[Adot] {};
\draw[red, thick]  (v2)--( 1.1,-2.7) node[Bdot] {}
                  --( 1.8,-3.4) node[Bdot] {};

% v3: non-splitting

% From v4 (only A-ray)
\draw[blue, thick] (v4)--(-1.1,-4.7) node[Adot] {}
                  --(-1.8,-5.4) node[Adot] {};

% Arrow to right panel
\draw[->, thick] (2.7,-3.0) -- (4.6,-3.0);

% ---------------- Right panel: number line ----------------

\draw[thick] (5.1,-3.0) -- (11.2,-3.0);

% a^*
\filldraw (7.7,-3.0) circle (1.8pt);
\node[below] at (7.7,-3.0) {$a^{*}$};

% Blue points (A^{[1]}) -- uneven spacing
\filldraw[blue] (10.4,-3.0) circle (1.4pt);
\filldraw[blue] (9.3,-3.0) circle (1.4pt);
\filldraw[blue] (8.15,-3.0) circle (1.4pt);

% Red points (A^{[2]}) -- uneven spacing
\filldraw[red] (9.8,-3.0) circle (1.4pt);
\filldraw[red] (8.55,-3.0) circle (1.4pt);

% Labels for collections
\node[blue] at (9.7,-2.2) {$A^{[1]}$};
\draw[blue,->] (9.5,-2.3) to[out=-120,in=40] (9.0,-2.95);

\node[red] at (9.2,-3.8) {$A^{[2]}$};
\draw[red,->] (9.0,-3.7) to[out=120,in=-40] (8.8,-3.05);

\end{tikzpicture}
\caption{\small{Step 2: decomposition of $A_{+}$ based on the distinguished ray $\mathcal R^{*}$. The left panel shows a subtree of the full trinary tree, where $M = 3$. The set $A_{+}$ consists of the points represented by the rays of this subtree. The distinguished ray $\mathcal R^{\ast}$ is shown in bold.  
For simplicity, we have chosen every vertex on $\mathcal R^{\ast}$ as the $0^{\text{th}}$ child of its parent. A blue (respectively red) ray originating at $v_j$ is descended from the $1^{\text{st}}$ (respectively $2^{\text{nd}}$) child of $v_j$. This gives rise to corresponding collections $A^{[1]}$ and $A^{[2]}$ of points on the real line (the right panel), each converging to the distinguished limit $a^{*}=\alpha(R^{*})$.
}}
\label{fig: step 2}
\end{figure}
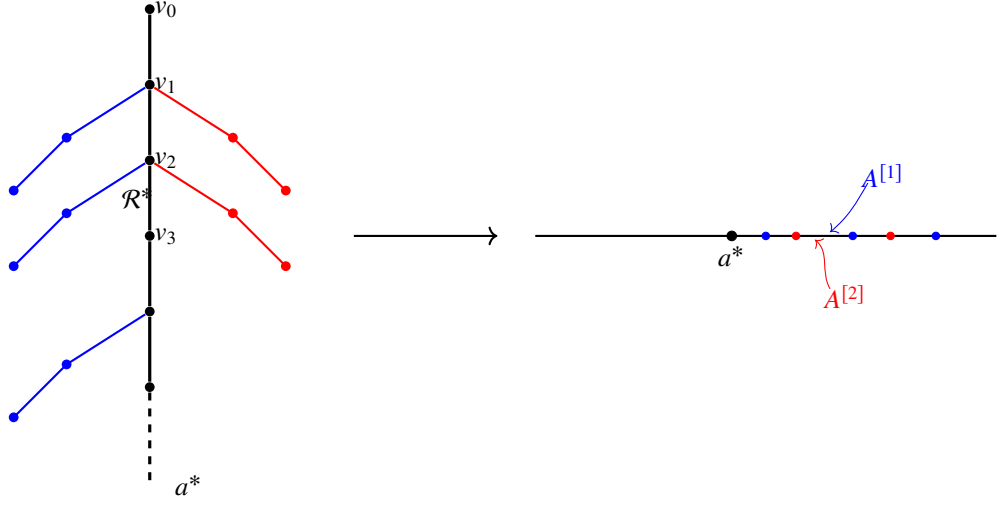
\subsection{Step 3: Lacunarity of $A^{[i]}$.} \label{step 3 base case}
We now fix $i \in \mathbb Z_M$ and proceed to cover $A^{[i]}$ by a threefold union of lacunary sequences converging to $a^{\ast}$. Let us consider the ordered subsequence
\begin{equation}   j_1 < j_2 < \cdots < j_r < \cdots  \label{strictly ordered j} \end{equation} 
of all non-negative integers $j$, possibly infinite, for which the ray $\mathcal R_{ij}$ defined in \eqref{uniqueness of non-splitting ray} exists; it is therefore unique and lies in $\mathcal T^{[i]}$.
%In other words, an index $j$ belongs to the subsequence \eqref{strictly ordered j} if and only if the $M$-adic interval represented by the $i^{\text{th}}$ child of $v_j$ has non-empty intersection with $A^{[i]}$. This intersection to be a singleton in view of \eqref{uniqueness of non-splitting ray}. 
Let $a_r := \alpha(\mathcal R_{ij_r})$ denote the point in $A^{[i]}$ specified by the ray $\mathcal R_{ij_r}$, so that \[ A^{[i]} = \{a_r: r \in \mathbb{N}_i\} \cup \{a^{\ast} \}. \] 
Here $\mathbb{N}_i$ is a (finite or infinite) subset of $\mathbb N$. 
To establish \eqref{Lacunary cover}, it suffices to show that the three subsequences of the form 
\begin{equation}  
\mathcal L_{\ell} := \{a^{\ast} \} \cup \{ a_{3r + \ell + 1} : r \geq 0, 3r + \ell + 1 \in \mathbb N_i \}, \; \ell = 0, 1,2,  
\end{equation}   which cover $A^{[i]}$ by definition,  are each contained in a member of {\tt{MonLac}}$(M^{-1})$ converging to $a^{\ast}$.  
\vskip0.1in
\noindent Let us prove this statement. The definition \eqref{uniqueness of non-splitting ray} of $\mathcal R_{ij}$ implies that 
\begin{equation} \label{ycm: a&a_r} h \bigl(D(a^{\ast}, a_r) \bigr) =  h \bigl(D(\alpha(\mathcal R^{\ast}), \alpha(\mathcal R_{ij_r}) ) \bigr) = j_r.  \end{equation}  
Here $D(a^{\ast}, a_r)$ denotes the youngest common ancestor of $a^{\ast}$ and $a_r$ in the $M$-adic tree, as defined in \eqref{defn: youngest common ancestor}. In other words, $a^{\ast}$ and $a_r$ lie in the same $M$-adic interval of length $M^{-j_r}$, as a result of which the distance between them is at most the length of the interval: 
\begin{equation}  \label{Madic distance}
0 \leq a_{r} - a^{\ast} \leq M^{-j_r} \text{ for } r \in \mathbb N_i. 
\end{equation} 
More importantly, $(a_r-a^{\ast})$ also obeys a lower bound: 
\begin{equation} \label{distance calculation} 
a_{r} - a^{\ast} \geq {M^{- j_{r+2}}} \text{ whenever } (r+2) \in \mathbb{N}_i.  
\end{equation}  
We will justify this lower bound in the next step, but a consequence of \eqref{Madic distance} and \eqref{distance calculation} is that for fixed $\ell = 0,1,2$,  and any $r \geq 0$, 
\begin{align*} 
a_{3(r+1)+\ell+1} - a^{\ast} &\leq M^{-j_{3r+4 + \ell}} = M^{-j_{3r+4 + \ell} + j_{3r+3 + \ell}} M^{-j_{3r+3+\ell}} \\  &\leq M^{-1} (a_{3r+\ell+1} - a^{\ast}), \end{align*}
where the last inequality follows from the strict ordering \eqref{strictly ordered j}. This holds whenever the displayed indices exist in $\mathbb{N}_i$. We therefore obtain both the monotonicity as well as the lacunarity of $\mathcal L_{\ell}$. Appending a geometric tail to $\mathcal L_{\ell}$ if it is finite, we find that $\mathcal L_{\ell}$ is contained in an infinite sequence in ${\tt{MonLac}}(M^{-1})$. This completes the proof of \eqref{Lacunary cover}.  
%Thus for every fixed $\ell = 0,1,2$, the sequence $\mathfrak A_{\ell} = \{ a_{n_{3k + \ell}} : k \geq 0 \}$ is covered by a lacunary sequence with constant $\leq M^{-1}$ converging to $a^{\ast}$. Since $A_{i+}$ is the union of $\{\mathfrak A_\ell : \ell = 0, 1,2\}$, the result follows.  
%Let $\{ \mathcal R_{\ell}(v_j) : \ell \geq 1 \}$ be an enumeration of $\mathfrak R(v_j)$. If $a_{\ell}(j) := \alpha(\mathcal R_{\ell}(v_j))$ is the point in $\bar{A}$ identified by $\mathcal R_{\ell}(v_j)$, we conclude that for every non-negative integer $j$, there are at most $(M-1)$ distinct points $a_{\ell}(j) \ne a^{\ast}$ in $\bar{A}$ such that 
\vskip0.1in
\subsection{Step 4: Euclidean separation between $a_r$ and $a^{\ast}$.} \label{step 4 base case} 
It remains to prove the estimate \eqref{distance calculation} from the previous step, which played an important role in establishing the lacunarity of the sequence $\mathcal L_{\ell}$.
Let 
\[ I_r = [p_r, q_r] = p_r + [0, M^{-j_r}]  \] denote the vertex of $\mathcal R^{\ast}$ at height $j_r$. This is an $M$-adic interval of length $M^{-j_r}$ containing $a^{\ast}$. We claim that the intervals $I_r$ of this nested decreasing sequence move in quantifiable steps to the left as $r$ increases, in the sense that  for all $r \geq 1$ with $r, r+1 \in \mathbb N_i$,  
\begin{equation}  q_{r+1} < q_{r}; \; \text{ in other words, } \; I_{r+1} \text{ cannot share a right endpoint with $I_{r}$}. \label{interval stacking} \end{equation}  
%An alternative way of expressing \eqref{interval stacking} is that $a^{\ast}$ cannot lie in the rightmost descendant of $I_r$ for any $r$ such that $r, r+1 \in \mathbb N_i$. 
\vskip0.1in
\noindent We prove \eqref{interval stacking} by contradiction.  Suppose, if possible, that $q_{r+1} = q_{r}$ for some index $r$. This would make $I_{r+1}$ the rightmost descendant of $I_{r}$ of the $j_{r+1}^{\text{th}}$ generation. On one hand, the definition of $I_r, I_{r+1}$ implies that 
\[ a^{\ast} \in I_r \cap I_{r+1} \; \text{ and therefore }  \; I_{r+1} \subseteq I_{r}. \] 
On the other hand, the definition \eqref{A-plus-minus} of $A_{+}$ implies that $a_{r} \geq a^{\ast}$. Jointly, these two statements force $a_{r}$ to lie in $I_{r+1}$, since the subinterval of $I_r$ to the right of $a^{\ast}$ is contained in $I_{r+1}$:
\[ a_r \in [a^{\ast}, 1] \cap I_r \subseteq I_{r+1}. \] Thus the youngest common ancestor of $a^{\ast}$ and $a_{r}$ is contained in $I_{r+1}$, and therefore has height at least $j_{r+1}$, contradicting the defining property \eqref{ycm: a&a_r} of $a_{r}$. The contradiction proves the claim \eqref{interval stacking}. 
\vskip0.1in 
\noindent The strict inequality \eqref{interval stacking} ensures quantifiable separation between  $a^{\ast}$ and $a_r$. These two points already lie in distinct children of $I_{r}$ by virtue of \eqref{ycm: a&a_r}, but apriori, this does not prevent them from lying in adjacent intervals of length $M^{-j_{r+1}}$ with arbitrarily small separation. Let $I'_{r+1}$ denote the $M$-adic interval of generation $j_{r+1}$ containing $a_r$; the ancestry relations \eqref{ycm: a&a_r} and \eqref{Madic distance}, together with the condition $a_r \geq a^{\ast}$, ensure 
\begin{equation}
\left\{
\begin{aligned}   
&I'_{r+1} \neq I_{r+1},\text{ $I'_{r+1}$ is to the right of $I_{r+1}$,} \\ 
& \text{ with } |I_{r+1}'| = |I_{r+1}| = M^{-j_{r+1}}
\end{aligned} 
\right\}.
\end{equation} 
Let us consider the $j_{r+2}^{\text{th}}$-generation descendants of $I_{r+1}$ and $I_{r+1}'$, and let $J$ denote the rightmost descendant of $I_{r+1}$ at this level. The condition \eqref{interval stacking} guarantees that $a^{\ast} \not\in J$. In other words, 
\[ \text{ $a^{\ast}$ lies to the left of $J$, whereas $a_r$ lies in $I'_{r+1}$, which is to the right of $J$. } \] 
A pictorial rendition of this appears in Figure \ref{fig: step 4}. As a result, $J$ separates $a^{\ast}$ from $a_r$, yielding the estimate 
\[a_r - a^{\ast} \geq |J| = M^{-j_{r+2}}. \] This leads to the conclusion  \eqref{distance calculation} and completes the proof of Lemma \ref{BASE CASE}. 
%in $I_{r+2}$ and $a_{n_k}$ (which is to the right of $I_{k+1}$) must lie on opposite sides of $J$, the rightmost $M$-adic subinterval of length $M^{-n_{k+2}}$ in $I_{k+1}$. This implies $a_{n_k} - a^{\ast} \geq |J|$, which is the conclusion of \eqref{distance calculation}.  
\qed
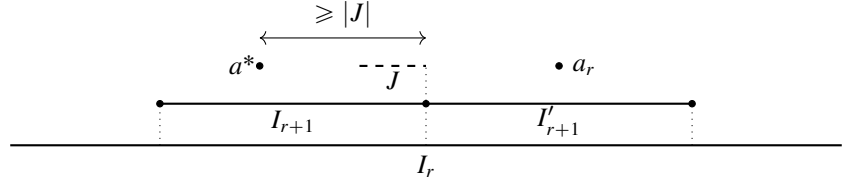
\begin{figure}[ht]
\centering
\begin{tikzpicture}[x=11cm,y=1cm,every node/.style={font=\small}]

% Main interval I_r
\draw[thick] (0,0) -- (1,0);
\node[above] at (0.5,-0.56) {$I_r$};

% Adjacent child intervals I_{r+1} and I'_{r+1}
\draw[thick] (0.18,0.55) -- (0.50,0.55);
\node[above] at (0.34,0.01) {$I_{r+1}$};

\draw[thick] (0.50,0.55) -- (0.82,0.55);
\node[above] at (0.66,-0.02) {$I'_{r+1}$};

% Mark endpoints so the two intervals are visually distinct
\filldraw (0.18,0.55) circle (1.2pt);
\filldraw (0.50,0.55) circle (1.2pt);
\filldraw (0.82,0.55) circle (1.2pt);

% Rightmost descendant J of I_{r+1}
\draw[thick,dashed] (0.42,1.05) -- (0.50,1.05);
\node[above] at (0.46,0.6) {$J$};

% Vertical guides
\draw[dotted] (0.18,0) -- (0.18,0.55);
\draw[dotted] (0.50,0) -- (0.50,1.05);
\draw[dotted] (0.82,0) -- (0.82,0.55);

% Points a^* and a_r
\filldraw (0.30,1.05) circle (1.2pt);
\node[above] at (0.28,0.80) {$a^{*}$};

\filldraw (0.66,1.05) circle (1.2pt);
\node[above] at (0.69,0.80) {$a_r$};

% Visual indicator of separation
\draw[<->] (0.30,1.42) -- (0.50,1.42);
\node[above] at (0.40,1.42) {$\geq |J|$};

\end{tikzpicture}
\caption{\small{Step 4: Euclidean separation between $a_r$ and $a^{*}$, as claimed in \eqref{distance calculation}. The intervals $I_{r+1}$ and $I'_{r+1}$ are distinct descendants of the ambient interval $I_r$, with $a^{*}\in I_{r+1}$ and $a_r\in I'_{r+1}$, so that $I'_{r+1}$ is to the right of $I_{r+1}$. If $J$ denotes the rightmost descendant of $I_{r+1}$ in generation $j_{r+2}$, then $a^{*}$ and $a_r$ lie on opposite sides of $J$. In other words, $J$ separates $a^{*}$ from $a_r$, yielding $a_r-a^{*}\ge |J|$.}}
\label{fig: step 4}
\end{figure}

% I_{r+1}

%%%%%%%%%%%%%%%%%%%%%%%%%%%%%%%%%%%%%%%%%%%%%%%%%%%%%%%%%%%%%%%%%%%%%%%
\section{Proof of Proposition \ref{TREE-LACUNARY-TRADITIONAL}: the inductive step} \label{completion of proof: Tree-Lacunary-Traditional}
The proof of Proposition  \ref{TREE-LACUNARY-TRADITIONAL} proceeds by induction on $N$. The base case $N = 1$ has been covered by Lemma \ref{BASE CASE}. To complete the induction step, we will assume that the statement of Proposition \ref{TREE-LACUNARY-TRADITIONAL} holds for trees of splitting number at most $(N-1)$. 
\vskip0.1in
\noindent Let us consider a bounded set $U \subseteq [0, \infty)$ obeying the hypothesis of the proposition. As in Lemma \ref{BASE CASE} and after a suitable scaling, it suffices to assume that $U$ is contained in $[0,1)$, and is  represented by the tree $\mathcal T := \mathcal T(U; M)$, with splitting number  
\[ \text{split}(\mathcal T) = \text{split}(\mathcal T(U; M)) = N \geq 2.  \] 
We aim to show that 
\[ U \in \Lambda(N, M^{-1}; R_N) \text{ for some controlled constant $R_N = R(N, M)$}. \] In view of Definition \ref{defn: Admissible finite order lacunarity} describing $\Lambda(N, M^{-1}; R_N)$, the objective is to identify a finite cover of $U$ using at most $R_N$ lacunary sets of order $N$: namely, there exists
\begin{equation}  \label{induction statement} 
\left\{
\begin{aligned} 
&\mathscr{U} \subseteq \Lambda(N, M^{-1}) \text{ with }  \#(\mathscr U) \leq R_N, \text{ such that } \\
& U \subseteq \bigcup \bigl\{ \hat{U} : \hat{U} \in \mathscr{U} \bigr\}, \text{ and }  \; R_N := R(N, M) = 18M R_{N-1}. 
\end{aligned} 
\right\}
\end{equation}    
Our analysis of $U$ is largely modelled on the framework of Lemma \ref{BASE CASE}. A visual depiction of the proof strategy appears in Figure \ref{fig: induction}.  
\subsection{The special ray $\mathcal R^{\ast}$}
As in Lemma \ref{BASE CASE}, we begin by selecting a ray $\mathcal R^{\ast} \in \partial \mathcal T$ containing all the vertices of $\mathcal T$ with highest splitting number; i.e. 
\begin{equation} \label{ray of highest split} 
\text{every vertex of $v^{\ast} \in \mathcal T$ with split$_{\mathcal T}(v^{\ast}) = N$  lies on $\mathcal R^{\ast}$.} 
\end{equation}  This ray is constructed exactly as in Section \ref{step 1 base case}. If there are multiple vertices $v^{\ast}$ with highest splitting number $N$, Lemma \ref{SplitsOnARay} is used to select a maximal ray $\mathcal R^{\ast}$ containing all such vertices. We omit the routine verification of
these details. Let us set $a^{\ast} := \alpha(\mathcal R^{\ast})$, which plays the role of the distinguished limit point. 
\vskip0.1in
\noindent We cover the closure $\bar{U}$ of $U$ by the union of its subsets lying to the left and to the right of $a^{\ast}$; specifically, 
\begin{align*}
&\hskip1in U \subseteq \bar{U} \subseteq U_{+} \cup U_{-}, \text{ where } \\ 
&U_{+} := \{a^{\ast} \} \cup \{x \in \bar{U}: x \geq a^{\ast} \}, \; U_{-} := \{a^{\ast} \} \cup \{x \in \bar{U}: x \leq a^{\ast}\}.
\end{align*}
Then, $\mathcal T(U_{\pm}) = \mathcal T(U_{\pm}; M)$ are subtrees of $\mathcal T$, and therefore, by the monotonicity of splitting numbers (Lemma \ref{monotone trees}),  
\[ \text{split}(\mathcal T(U_{\pm})) \leq \text{split}(\mathcal T) \leq N. \] 
Thus, in order to establish \eqref{induction statement} for $U$, it suffices to prove the same for $U_{+}$ and $U_{-}$, with $R_N$ replaced by $12MR_{N-1}$ and $6MR_{N-1}$ respectively. Henceforth, we focus only on $U_{+}$. The proof for $U_{-}$ involves minor modifications; we indicate these changes at the end of the proof, in Section \ref{lac split end of proof}.  
\vskip0.1in
\noindent If the splitting number of $\mathcal T(U_{+})$ is $\leq N-1$, then the induction hypothesis already covers $U_{+}$. Thus, without loss of generality we assume that 
\[ \text{split}(\mathcal T(U_{+})) = N, \] 
so that $\mathcal R^{\ast}$ is still a special ray of $\mathcal T(U_{+})$. Let $\hat{\mathcal V}$ denote the non-empty collection of vertices of $\mathcal T_{+} = \mathcal T(U_{+}; M)$ that are one generation removed from $\mathcal R^{\ast}$, namely 
\begin{equation} \label{first generation away from the special ray}
\hat{\mathcal V} := \left\{v \in \mathcal T_{+} : v \notin \mathcal R^{\ast}, \text{ but the parent of } v \text{ lies on } \mathcal R^{\ast} \right\}. 
\end{equation}  
The vertices of $\hat{\mathcal V}$ help in the decomposition of $U_{+}$ into lacunary components $\hat{U} \in \Lambda(N, M^{-1})$ as required by \eqref{induction statement}, and also identify their associated special sequences. This is clarified in the next step.  
\subsection{Choice of special sequences for $\hat{U}$} 
For every vertex $v \in \hat{\mathcal V}$, let $a_v$ denote the left endpoint of the $M$-adic interval represented by $v$. 
\begin{equation} v = a_v + [0, M^{-h(v)}]; \; \text{ we set } \; A := \{ a_v : v \in \hat{\mathcal V} \} \cup \{a^{\ast} \}. \label{lacunary-1-A} \end{equation}  
Each point $x\in U\setminus\{a^*\}$ has a unique ancestor
$v\in\hat{\mathcal V}$, and the choice of $a_v$ gives $a_v\leq x$.
These left-endpoint anchors will allow us to control the splitting
number of $\mathcal T(A;M)$. 
\vskip0.1in
\noindent The set $A$ is the source from which a lacunary special sequence of each $\hat{U}$ will be extracted. To make this precise, we claim that 
\begin{equation} \text{split}(\mathcal T(A;M)) =1. \label{claim-A-split-1} \end{equation}  
\vskip0.1in
\noindent The claim \eqref{claim-A-split-1} follows from two observations: 
\vskip0.1in
\begin{itemize} 
\item First, we deduce from the definition \eqref{lacunary-1-A} of $A$ that $a^{\ast} \in A$ ; therefore $\mathcal T(A; M)$ contains the ray $\mathcal R^{\ast}$. If $w$ is the parent of a vertex $v \in \hat{\mathcal V}$, then $w$ is a splitting vertex of $\mathcal T(A; M)$ lying on $\mathcal R^{\ast}$. Indeed, the definition \eqref{first generation away from the special ray} of $\hat{\mathcal V}$ implies that $w$ has at least two children in $\mathcal T(A; M)$; one lies on $\mathcal R^{\ast}$ and contains $a^{\ast}$, the other is $v \not\in \mathcal R^{\ast}$ which contains $a_v$. 
%Let us recall that each number $a_v$ is an $M$-adic rational, and our convention from Section \ref{section: rays as points} dictates that the ray encoding such a rational should reflect its finitary expansion (depicting convergence from the right). Therefore $a^{\ast}$ and $a_v$ lie in distinct $M$-adic intervals at height $h(v)$. 
This ensures that there is at least one splitting vertex in $\mathcal T(A; M)$, so that 
\begin{equation} \text{split}(\mathcal T(A;M)) \geq 1. \label{claim-A-split-more} \end{equation} 
\vskip0.1in
\item Second, the only splitting vertices of $\mathcal T(A;M)$ occur on $\mathcal R^{\ast}$. Indeed, for each $v\in\hat{\mathcal V}$, the canonical-ray convention \eqref{canonical ray convention} implies that among the points of $A$, only $a_v$ has a
representing ray through $v$.  Thus the subtree of
$T(A;M)$ rooted at $v$ is non-splitting.  This implies that every subtree $\mathcal S \subseteq \mathcal T(A;M)$ has a maximal ray that contains at most one splitting vertex, so that  
\begin{equation} \text{split}(\mathcal T(A;M)) \leq 1. \label{claim-A-split-less} \end{equation}  
Jointly, \eqref{claim-A-split-more} and \eqref{claim-A-split-less} prove the claim \eqref{claim-A-split-1}. 
\end{itemize} 
\vskip0.1in 
The relation \eqref{claim-A-split-1} enables us to invoke the induction hypothesis on $A$ with $N=1$. By Lemma \ref{BASE CASE}, $A$ can be covered by the union of at most $6M$ infinite lacunary sequences: 
\begin{equation} \label{decomp-A} 
A \subseteq \bigcup_{r=1}^{r_0} A_r, \quad 1 \leq r_0 \leq 6M, \quad A_r \in {\tt{MonLac}}(M^{-1}) \text{ converges  to $a^{\ast}$}. 
\end{equation}
The monotonicity of $A_r$ implies that it must be monotonically decreasing to $a^{\ast}$, i.e., $\inf A_r = a^{\ast}$.  
\vskip0.1in
\noindent The covering \eqref{decomp-A} of $A$ leads automatically to a decomposition of its defining index set of vertices $\hat{\mathcal V}$. Namely, if $\hat{\mathcal V}_r \subseteq \hat{\mathcal V}$ denotes the subset of vertices that correspond to $A_r$, in the sense that
\begin{equation} \hat{\mathcal V}_r := \bigl\{v \in \hat{\mathcal V}: a_v \in A_r \bigr\}, \text{ then \eqref{lacunary-1-A} and \eqref{decomp-A}} \implies  \hat{\mathcal V} = \bigcup_{r=1}^{r_0} \hat{\mathcal V}_r.   \label{V-hat-r}\end{equation}  
%Let us note in passing that if $\mathcal T_v(E; M)$ is  a subtree of $\mathcal T$ rooted at a vertex $v \in \hat{\mathcal V}_r$ and representing a set $E \subseteq U$, then the same argument that led to \eqref{relation between inf sup UA} also shows 
%\begin{equation} \label{relation between inf sup EA}
%[\inf E, \sup E] \subseteq [\inf A_r, \sup A_r] = [a^{\ast}, \sup A_r]. 
%\end{equation} 
\vskip0.1in
\noindent Using the vertex sets $\hat{\mathcal V}_r$ for $1 \leq r\leq r_0$, we now proceed to identify a class $\mathscr{U}$ of subsets of $U_{+}$ lying in $\Lambda(N, M^{-1})$, with special sequence $A_r$.
\subsection{First-order decomposition of $U_{+}$ based on sub-trees off $\mathcal R^{\ast}$}  For every $v \in \hat{\mathcal V}$, let $\mathcal T_v(U_{+})$ denote the maximal subtree of $\mathcal T_{+}$ rooted at $v$. Let $U_{v} \subseteq \bar{U}_{+}$ denote the closed set represented by the $M$-adic tree $\mathcal T_{\nu}$: 
\begin{equation}  \label{defn: T_v}
\mathcal T_v =: \mathcal T(U_v;M). 
\end{equation} 
The fact that $\mathcal T_v$ is rooted at $v$ implies that $U_v \subseteq v$. As $v$ ranges over $\hat{\mathcal V}$, the subtrees $\mathcal T_v \cup \{\mathcal R^{\ast}\}$ generate a decomposition of $\mathcal T$:
\begin{align*} 
\mathcal T &= \bigcup \Bigl\{\mathcal T_v \cup \{\mathcal R^{\ast} \} : v \in \hat{\mathcal V} \Bigr\} 
\\ &= \bigcup \Bigl\{\mathcal T_v \cup \{\mathcal R^{\ast} \} : v \in \hat{\mathcal V}_r, 1 \leq r \leq r_0 \Bigr\}, 
\end{align*}
where the last step follows from the decomposition of $\hat{\mathcal V}$ in \eqref{V-hat-r}. 
The decomposition of the tree $\mathcal T$ leads, in turn, to a decomposition of $U_{+}$:  
\begin{equation} \label{decomposition of U}
U_{+} \subseteq \bigcup \Bigl\{ U_{v} \cup \{a^{\ast} \} : v \in \hat{\mathcal V}_r, \; 1 \leq r \leq r_0 \Bigr\}.
\end{equation}  

\subsection{Second-order decomposition: the induction hypothesis} 
The definition \eqref{ray of highest split} of the ray $\mathcal R^{\ast}$ as the distinguished ray containing all vertices of highest split means that   
\begin{equation}  \label{lesser split}
\text{split}_{\mathcal T}(v) \leq N-1 \text{ for every } v \not\in \mathcal R^{\ast}.   
\end{equation}  
In particular, this is true for all vertices $v \in \hat{\mathcal V}$ as in \eqref{first generation away from the special ray}. We therefore deduce from \eqref{defn: T_v} and \eqref{lesser split} that the set $U_{v}$ satisfies
\[ \text{split}(\mathcal T(U_v;M)) = \text{split}(\mathcal T_v) = \text{split}_{\mathcal T}(v) \leq N-1 \text{ for each } v \in \hat{\mathcal V}, \] 
Applying the induction hypothesis to $U_v$, one finds a constant $R_{N-1} = R(N-1, M)$ such that $U_v$ is covered by the union of at most $R_{N-1}$ sets in $\Lambda(N-1;M^{-1})$:
 \begin{equation} \label{decomposition of U_v} 
 U_{v} \subseteq \bigcup_{i=1}^{R_{N-1}} U_v(i), \quad U_v(i) \in \Lambda(N-1, M^{-1}). 
 \end{equation} 
 Without loss of generality, intersecting with $v$ if necessary, we may assume that $U_v(i) \subseteq v$. The sets $U_v(i)$ will need to be carefully assembled to obtain a cover of $U$ with $R_N$ lacunary sets of order $N$. 
  \vskip0.1in
 \noindent For each $r$ with $\hat{\mathcal V}_r\neq \emptyset$, let $v_r$ be the
vertex in $\hat{\mathcal V}_r$ for which the left endpoint $a_{v_r}$ is largest. Such a
vertex exists because the points $\{a_v: v \in \hat{\mathcal V}_r \}$ form a subset of the decreasing sequence $A_r$. As a result, 
\[ a_{v_r} \leq \sup A_r, \text{ and any $x \in v$, $v \in \hat{\mathcal V}_r \setminus \{v_r\}$ obeys $x \leq a_{v_r}$}. \] The rightmost vertex $v_r \in \hat{\mathcal V}_r$ needs separate consideration, as we shall see presently. 
\vskip0.1in 
\noindent With these notations in place, we are now ready to define the cover for $U_{+}$. For
$1\leq i\leq R_{N-1}$, we set
\begin{equation} \label{defn: W(r,i)} 
  W(r,i):=\{a^*\}\cup
     \bigcup \bigl\{ U_v(i) \cap v^{\sharp} : v\in \hat{\mathcal V}_r\setminus\{v_r\} \bigr\},
  \quad
  E(r,i):=U_{v_r}(i).
\end{equation} 
Here $v^{\sharp}$ is the canonical half-open cell associated to the closed interval $v$: 
\[ \text{if } v = [a, b] \; \text{ then } \; v^{\sharp} = [a, b). \] 
Since the half-open cells of two distinct $M$-adic intervals are disjoint, defining $W(r, i)$ in this way ensures that its components $U_v(i) \cap v^{\sharp}$ do not overlap. This allows us to avoid endpoint issues in the covering statement, as we will see in Section \ref{lac split end of proof}.
\vskip0.1in
\noindent The first observation is that $U_{+}$ is covered by the sets $W(r, i)$ and $E(r, i)$. To see this, let $x\in U\setminus\{a_*\}$, and let $v\in\hat{\mathcal V}$ be the first vertex at which the canonical ray representing $x$ leaves $\mathcal R^*$. Then $x\in U_v\cap v^\sharp$. Choose $r$ such that $v\in\hat{\mathcal V}_r$, and $i$ such that $x\in U_v(i)$, as permitted by \eqref{decomposition of U_v}. If $v=v_r$, then $x\in E(r,i)$; otherwise $x\in W(r,i)$. The point $a_*$ belongs
to every \(W(r,i)\). Hence, we arrive at the following cover for $U_{+}$: 
\begin{equation} \label{cover U}
 U_{+}\subseteq \bigcup_{r=1}^{r_0} \bigcup_{i=1}^{R_{N-1}} \bigl\{W(r,i)\cup E(r,i) : \hat{\mathcal V}_r \neq \emptyset \bigr\}.
\end{equation} 
There are at most $2r_0R_{N-1}\leq 12MR_{N-1}$ sets in this cover.
% Combining \eqref{decomposition of U} with \eqref{decomposition of U_v}, we obtain
% \begin{align} U &\subseteq \bigcup \Bigl\{ U_{v}^{[i]} \cup \{a^{\ast} \} : v \in \hat{\mathcal V}_r, 1 \leq i \leq R_{N-1}, 1 \leq r \leq r_0 \Bigr\} \nonumber \\ 
% &\subseteq \bigcup \bigl\{ U(r, i) :  1 \leq r \leq r_0, \; 1 \leq i \leq R_{N-1} \bigr\} \; \text{ where } \nonumber  \\ &U(r, i) := \bigcup \Bigl\{ U_v^{[i]} \cup 
% \{a^{\ast} \} : v \in \hat{\mathcal V}_r \Bigr\}. \label{defn: U(r, i)}
%\end{align}  
%In other words, the set $U(r, i) \setminus \{a^{\ast} \}$ consists of a collection of disjoint subsets $U_v^{[i]}$ indexed by $v\in \hat{\mathcal V}_r$:  
%\begin{equation} \label{decomp-Uri} 
%U(r, i) \setminus \{a^{\ast} \} = \bigsqcup_{v \in \hat{\mathcal V}_r} U_v^{[i]}, \quad U_v^{[i]} \subseteq v.
%\end{equation} 
%The sets $U(r, i)$ given by \eqref{defn: U(r, i)} comprise the collection $\mathscr{U}$, whose existence is claimed in \eqref{induction statement}. We have shown that 
%\[ \#(\mathscr{U}) \leq r_0 R_{N-1} \leq 6MR_{N-1}, \] and the set $U$ can be covered by the union of the sets $\hat{U} \in \mathscr{U}$. 
 %where each $U^{[i]}$ shares a tree structure similar to $U$: it contains the point identified by $R^{\ast}$, with the additional feature that now $U_v^{[i]} \in \Lambda(N-1;M^{-1})$ for every $v \in \mathcal V^{[i]}$, where \[\mathcal V^{[i]} := \bigl\{ v \in \mathcal T(U^{[i]}; M) : v  \notin R^{\ast}\text{ but parent of $v$ is in } R^{\ast} \bigr\}. \]   
 \subsection{Lacunarity order of $W(r, i)$} \label{lac split end of proof}
By \eqref{decomposition of U_v}, we already have $E(r, i) \in \Lambda(N-1, M^{-1})$.  To complete the proof of the induction statement \eqref{induction statement}, it therefore suffices to show that 
 \begin{equation} \label{induction statement 2}
 W(r, i) \in \Lambda(N, M^{-1}) 
  \end{equation} 
  with the sequence $A_r \in {\tt{MonLac}}(M^{-1})$ from \eqref{decomp-A} as its special sequence. In view of Definition \ref{defn: Lacunary sets}, this means checking two things: first,
  \begin{equation} \label{W(r,i) extremes}
  \bigl[\inf W(r, i), \sup W(r, i) \bigr] \subseteq [a^{\ast}, a_{v_r}] \subseteq \bigl[\inf A_r, \sup A_r \bigr],
  \end{equation} 
and second, the portion of $W(r, i)$ lying within a gap interval of $A_r$ is of lower lacunary order.
\vskip0.1in
\noindent The containment \eqref{W(r,i) extremes} is immediate. Given any $x \in W(r, i)$, $x \neq a^{\ast}$, the definition \eqref{defn: W(r,i)} gives 
\[ x \in v^{\sharp} \text{ for some $v \in \hat{\mathcal V}_r$, $v \neq v_r $, and therefore } x < a_{v_r} \leq \sup A_r. \] 
\vskip0.1in
\noindent Let us confirm the second requirement. Indeed, suppose that $c$ and $d$ are two consecutive elements of $A_r$, so that $I = [c, d)$ is a gap interval as defined in \eqref{gap interval}.  We claim that 
\begin{equation} \label{W in gap}
W(r, i) \cap I \neq \emptyset \implies \exists \text{ a unique $v \in \hat{\mathcal V}_r \setminus \{v_r\}$ with } W(r, i) \cap I \subseteq U_v(i). 
\end{equation} 
The fact that $W(r, i) \cap I$ is non-trivial means some vertex $v$ obeying \eqref{W in gap} must exist. For uniqueness, suppose towards a contradiction that two distinct half-open cells $v^\sharp$, $w^\sharp$ contributing to $W(r,i)$ both meet $[c,d)$. Let $v^\sharp$ be the interval to the left of $w^{\sharp}$. Because $v^\sharp$ excludes its right endpoint, the left endpoint $a_w$ of $w$ would have to satisfy $c<a_w<d$. Since $a_w\in A_r$, this contradicts the gap condition, proving \eqref{W in gap}.
%Without loss of generality, let us assume that $A_r$ is monotone decreasing. That means 
%\[c = a_v \; \text{ and } \; d = a_{v'}, \text{ for vertices }  v, v' \in \hat{\mathcal V}_r \text{ with } h(v) > h(v'). \]  %Let $Q_v$ and $Q_{v'}$ denote the $M$-adic intervals representing $v$ and $v'$ respectively. 
%The definition of $\hat{\mathcal V}_r$ implies that $v$ and $v'$ are distinct vertices one generation removed from $\mathcal R^{\ast}$. Thus $v \not\subset v'$ and $v' \not\subset v$, which means the interiors of the intervals $v$ and $v'$ are disjoint. The monotone decreasing nature of $A_r$ means that the interval $v$ lies to the left of the interval $v'$ on the real line. The fact that no element of $A_r$ lies strictly in between $c$ and $d$ precludes the existence of a vertex $w \in \hat{\mathcal V}_r$ with $h(v') < h(w) < h(v)$, 
\vskip0.1in
\noindent The last set occurring in \eqref{W in gap}, namely $U_v(i)$, is in $\Lambda(N-1, M^{-1})$ as confirmed by the induction hypothesis \eqref{decomposition of U_v}. Therefore so does its subset $W(r, i)\cap I$.  This proves that $U_{+} \in \Lambda(N, M^{-1}; 12 M R_{N-1})$. 
\vskip0.1in
\noindent The argument for $U_{-}$ is analogous. The special sequences for the latter increase toward $a^*$, and the left-endpoint anchors already give the required span condition, so there is no need to separate an outermost vertex $v_r$. The same half-open-cell argument verifies the gap condition. In other words the sets $W(r, i)$ alone cover $U_{-}$; $E(r, i)$ is not needed. Thus the cover uses at most $6MR_{N-1}$ sets, compared with $12MR_{N-1}$ for $U^+$. Thus we may take $R_N=18MR_{N-1}$. This completes the proof of \eqref{induction statement} and hence of Proposition \ref{TREE-LACUNARY-TRADITIONAL}. 
%Let $\mathcal T^{[i]}$ denote the $M$-adic tree representing $U^{[i]}$. The definition of $U^{[i]}$ confirms that $\mathcal R^{\ast} \in \partial \mathcal T^{[i]}$. Indeed, $\mathcal T^{[i]}$ consists of the ray $\mathcal R^{\ast}$, along with a subtree rooted at every $v \in \mathcal V$ that represents a set of $\Lambda(N-1, M^{-1})$.  
%\vskip0.1in
%If $a=a_v$ and $b$ are two successive elements of this sequence with $a < b$, then $U^{[i]} \cap [a,b) = U_v^{[i]}$, which is in $\Lambda(N-1;M^{-1})$. Thus $U^{[i]}$ is in $\Lambda(N;M^{-1})$ according to Definition \ref{defn: Lacunary sets}, completing the proof. 
\qed
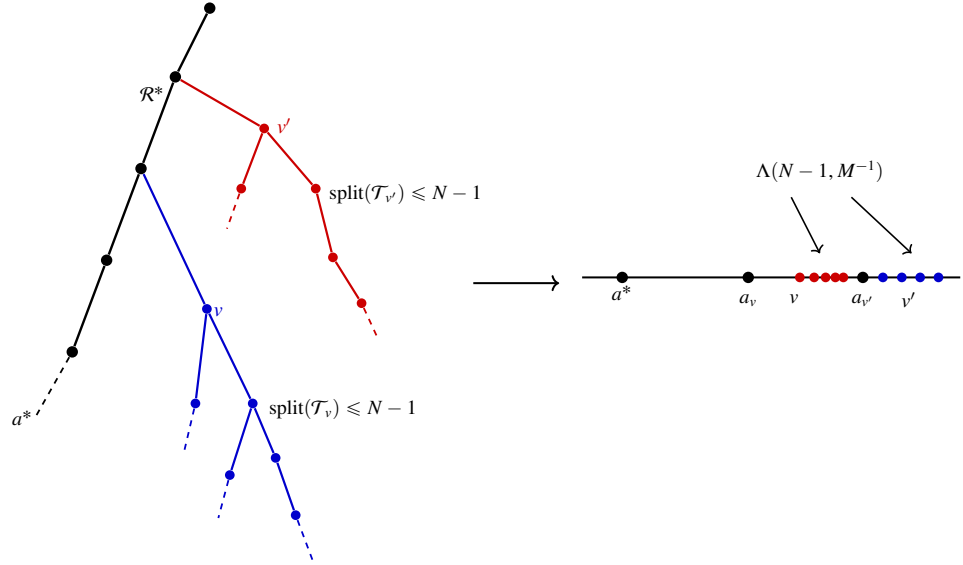
\begin{figure}[ht]
\centering
\begin{adjustbox}{max width=\textwidth}
\begin{tikzpicture}[x=1cm,y=1cm,every node/.style={font=\small}]
\tikzset{
  blackdot/.style={circle, fill=black, inner sep=2pt},
  reddot/.style={circle, fill=red!80!black, inner sep=1.7pt},
  bluedot/.style={circle, fill=blue!80!black, inner sep=1.7pt}
}

% =========================
% Left panel: main ray R^*
% =========================

\node[blackdot] (r4) at (0.8,5.8) {};
\node[blackdot] (r3) at (0.2,4.6) {};
\node[blackdot] (r2) at (-0.4,3.0) {};
\node[blackdot] (r1) at (-1.0,1.4) {};
\node[blackdot] (r0) at (-1.6,-0.2) {};

\draw[line width=1.2pt] (r4)--(r3)--(r2)--(r1)--(r0);
\draw[dashed, line width=0.9pt] (r0)--(-2.25,-1.35);

\node at (-0.2,4.3) {$\mathcal R^{*}$};
\node[left] at (-2.15,-1.35) {$a^{*}$};

% =========================
% Upper red subtree rooted at v'
% =========================

\node[reddot] (vp)  at (1.75,3.70) {};
\node[reddot] (rp1) at (1.35,2.65) {};
\node[reddot] (rp2) at (2.65,2.65) {};
\node[reddot] (rp3) at (2.95,1.45) {};
\node[reddot] (rp4) at (3.45,0.65) {};

\draw[red!80!black, line width=1.1pt] (r3)--(vp);
\draw[red!80!black, line width=1.1pt] (vp)--(rp1);
\draw[red!80!black, line width=1.1pt] (vp)--(rp2);
\draw[red!80!black, line width=1.1pt] (rp2)--(rp3);
\draw[red!80!black, line width=1.1pt] (rp3)--(rp4);

% dashed terminal rays
\draw[red!80!black, dashed, line width=0.9pt] (rp1)--(1.10,1.95);
\draw[red!80!black, dashed, line width=0.9pt] (rp4)--(3.75,0.00);

\node[red!80!black] at (2.10,3.75) {$v'$};
\node[right] at (2.75,2.55) {$\operatorname{split}(\mathcal T_{v'})\leq N-1$};

% =========================
% Lower blue subtree rooted at v
% =========================

\node[bluedot] (v)  at (0.75,0.55) {};
\node[bluedot] (b1) at (0.55,-1.10) {};
\node[bluedot] (b2) at (1.55,-1.10) {};
\node[bluedot] (b3) at (1.15,-2.35) {};
\node[bluedot] (b4) at (1.95,-2.05) {};
\node[bluedot] (b5) at (2.30,-3.05) {};

\draw[blue!80!black, line width=1.1pt] (r2)--(v);
\draw[blue!80!black, line width=1.1pt] (v)--(b1);
\draw[blue!80!black, line width=1.1pt] (v)--(b2);
\draw[blue!80!black, line width=1.1pt] (b2)--(b3);
\draw[blue!80!black, line width=1.1pt] (b2)--(b4);
\draw[blue!80!black, line width=1.1pt] (b4)--(b5);

% dashed terminal rays
\draw[blue!80!black, dashed, line width=0.9pt] (b1)--(0.35,-1.95);
\draw[blue!80!black, dashed, line width=0.9pt] (b3)--(0.95,-3.10);
\draw[blue!80!black, dashed, line width=0.9pt] (b5)--(2.60,-3.85);

\node[blue!80!black] at (0.95,0.5) {$v$};
\node[right] at (1.7,-1.20) {$\operatorname{split}(\mathcal T_v)\leq N-1$};

% Arrow between panels
\draw[->, line width=1pt] (5.40,1.00) -- (6.90,1.00);

% =========================
% Right panel: real line
% =========================

\draw[line width=1pt] (7.30,1.10) -- (13.90,1.10);

% a^*
\node[blackdot] (astar) at (8.00,1.10) {};
\node[below] at (8.00,1.10) {$a^{*}$};

% a_{v'} and a_v as black dots
\node[blackdot] (avp) at (10.20,1.10) {};
\node[below] at (10.20,0.95) {$a_{v}$};
\node[below] at (11,0.95) {$v$};

\node[blackdot] (av) at (12.20,1.10) {};
\node[below] at (12.20,0.95) {$a_{v'}$};
\node[below] at (13,0.98) {$v'$};

% red cluster: between a_{v'} and a_v, closer to a_v
\node[reddot] at (11.10,1.10) {};
\node[reddot] at (11.35,1.10) {};
\node[reddot] at (11.55,1.10) {};
\node[reddot] at (11.72,1.10) {};
\node[reddot] at (11.86,1.10) {};

% blue cluster: to the right of a_v, more spread out
\node[bluedot] at (12.55,1.10) {};
\node[bluedot] at (12.88,1.10) {};
\node[bluedot] at (13.20,1.10) {};
\node[bluedot] at (13.52,1.10) {};

% lacunary annotation
\node at (11.45,2.95) {$\Lambda(N-1,M^{-1})$};
\draw[->, line width=0.8pt] (10.95,2.50) -- (11.45,1.52);
\draw[->, line width=0.8pt] (12.00,2.50) -- (13.05,1.52);

\end{tikzpicture}
\end{adjustbox}
\caption{\small{A schematic diagram of the inductive step in Section \ref{completion of proof: Tree-Lacunary-Traditional}. The left panel shows the tree $\mathcal T = \mathcal T(U; M)$, with the distinguished ray $\mathcal R^*$ containing all vertices of maximal splitting number. For vertices in $\hat{\mathcal V}$ one generation away from $\mathcal R^*$, e.g. $v, v'$, the rooted subtrees $\mathcal T_v, \mathcal T_{v'}$, indicated in blue and red,  have splitting number at most $N-1$. The corresponding portions of $U$ in the intervals $v, v'$ on the real line therefore organize into lower-order lacunary components.}}
\label{fig: induction}
\end{figure}

\chapter{Pruning of the slope tree} \label{Chapter: Pruning of the slope tree} 
\section{Chapter overview}
This chapter develops a key structural refinement of the slope tree associated with a given direction set $\Omega$, with the goal of enforcing a quantitative form of Euclidean separation among its vertices.
\vskip0.1in
\noindent As discussed in the introduction (Chapter \ref{chapter: intro}), Euclidean separation among slopes is an essential ingredient in the construction of Kakeya-type configurations. However, an arbitrary slope set, particularly a sublacunary one, need not satisfy any uniform separation property at the outset. In particular, two vertices that separate early in the $M$-adic slope tree may still be arbitrarily close in Euclidean distance, highlighting the mismatch between $M$-adic and Euclidean notions of separation. The central objective of this chapter is to show that, despite this lack of initial structure, it is always possible to extract a sufficiently large subtree in which an appropriate separation condition holds.
\vskip0.1in
\noindent The main result of this chapter is Proposition \ref{PRUNING STAGE 1}, which asserts that any slope tree associated with a sublacunary set contains a pruned subtree enjoying controlled Euclidean separation. The pruning procedure is iterative in nature: at each stage, vertices that violate the desired separation condition are removed, while preserving enough of the tree structure to maintain a large splitting number.
\vskip0.1in
\noindent To carry out this procedure, we introduce in Section \ref{section: iterative pruning} a technical tool (Lemma \ref{PRUNING BUILDING BLOCK}) that allows us to pass from local separation properties to a global pruning mechanism. This lemma serves as the inductive engine of the construction.
\vskip0.1in
\noindent The structure of the chapter is as follows. In Section \ref{section: model trees}, we introduce model trees and sets that capture the essential combinatorial features of a desirable sublacunary slope set. Section \ref{section: sublac model subsets} shows that any sublacunary set contains such a model subset of arbitrarily large lacunarity order. In Sections \ref{section: iterative pruning}--\ref{Pruning building block proof section}, we develop and implement the pruning process, culminating in the proof of Proposition \ref{PRUNING STAGE 1}. Finally, Section \ref{section: HRS pruning} illustrates the procedure through an explicit example.

\section{Model trees and sets} \label{section: model trees}
Let us begin by identifying the structural features that our pruning procedure is designed to produce. These features will be encoded in a class of subtrees of the full $M$-adic tree, which we refer to as {\em{model trees}}.
\vskip0.1in
\noindent Suppose that $\mathcal S$ is a subtree of the full $M$-adic tree $\mathcal T([0,1];M)$ introduced in Section \ref{section: full M-adic tree}. The height of $\mathcal S$ may be finite or infinite. A priori, the rays of $\partial \mathcal S$ could be of different lengths. Of course, if $\mathcal S$ represents a subset of $[0,1]$, then the height of $\mathcal S$ is infinite, and every ray uniquely identifies a point in the set. However, we do not assume this at the moment. 
\vskip0.1in
\noindent Our goal is to isolate a subtree of $\mathcal S$ that balances two distinct, and potentially competing, requirements: a controlled and uniform branching structure, together with a quantitative form of Euclidean separation among certain descendants. These properties will ensure that the tree retains sufficient combinatorial richness while supporting the geometric separation needed for later constructions. We now formalize these requirements.
%Let us first specify a list of requirements that our pruning procedure will aim to meet. Suppose that $\mathcal S$ is a sub-tree of the full $M$-adic tree $\mathcal T([0,1]; M)$ introduced in Section \ref{section: full M-adic tree}. The height of $\mathcal S$ could be finite or infinite, and the rays of $\partial \mathcal T$ could be of different lengths. 
\vskip0.1in
\noindent Given an integer $N \geq 1$ and a constant $C_0 \geq 1$, we say that $\mathcal S$ is a {\em{model $(N, C_0)$-tree}} if it satisfies all the following three properties:  
\vskip0.1in
\begin{enumerate}[1.] 
\item \label{N split per ray} (Uniformity of splitting vertices per ray) {\em{Each maximal ray of $\mathcal S$ contains the same number of splitting vertices: }} 
\vskip0.1in 
\begin{itemize} 
\item There are exactly $2^N$ rays in $\partial \mathcal S$.
\vskip0.1in
\item Each ray in $\partial \mathcal S$ contains exactly $N$ splitting vertices, at heights that could possibly differ from one ray to another.
\end{itemize} 
\vskip0.1in 
To paraphrase,  
\begin{equation} \label{N split per ray eq}  
\#(\partial \mathcal S) = 2^N, \; \text{ and } \; \#\left\{ v \in \mathcal R : v \text{ is a splitting vertex of } \mathcal S \right\} = N 
\end{equation}
for every $R \in \partial \mathcal S$. This ensures that every maximal ray undergoes the same number of branching events.
\vskip0.1in 
\item \label{two children per split} (Controlled progeny) {\em{Every vertex of $\mathcal S$ has at most two children: }} 
\vskip0.1in
\begin{itemize} 
\item This means that every splitting vertex in $\mathcal S$ has exactly two children:
\begin{equation} \label{two children per split eq} 
\# \bigl\{ u \in \mathcal S : u \text{ is a child of } v\bigr\} = 2 \text{ for every splitting vertex } v \in \mathcal S.  
\end{equation}  
\end{itemize}
This keeps the local proliferation of vertices under tight control, regardless of the base $M$. 
\vskip0.1in
\item \label{separation condition} (Controlled Euclidean separation) {\em{Every splitting vertex $v \in \mathcal S$ admits two descendants of a certain subsequent generation, separated by a Euclidean distance that is compatible with their length. This height occurs before reaching the next splitting descendants of $v$. }} More precisely, 
\vskip0.1in
\begin{itemize} 
\item For every splitting vertex $v$ of $\mathcal S$, there is a minimal height $h^{\ast}_v > h(v)$ satisfying the following property: $v$ has exactly two descendants $w_1(v), w_2(v)$ at this height,    
\begin{align} 
&h(w_1(v)) = h(w_2(v)) = h_v^{\ast}, \text{ and the vertices $w_i(v)$ obey } \label{what is h_v-star} \\
&C_0 M^{-h^{\ast}_v} \leq \text{dist}(w_1(v), w_2(v)) \leq (C_0 + 2) M^{-h_{v}^{\ast}+1}.  \label{Euclidean distance condition}
\end{align} 
The separating height $h_v^{\ast}$ satisfies the bounds 
\begin{equation} \label{h-star-v}
h(v) < h^{\ast}_v < \min \bigl\{h(w) : w \subsetneq v, \; w \text{ is a splitting vertex of $\mathcal S$}  \bigr\},  
\end{equation} 
provided the set on the right hand side is non-empty. This property encodes the Euclidean separation that the pruning process aims to enforce.
\end{itemize}
%there exists a such that  In other words, neither of the two children of $v$ has a splitting descendant of a generation strictly smaller %than $h^{\ast}_v$. Further. In fact, $h_{v}^{\ast}$ is the smallest integer exceeding $h(v)$ with the property \eqref{Euclidean distance condition}. 
\end{enumerate}
\vskip0.1in
\noindent In item \ref{separation condition} of the above definition, and throughout this section, $h(v)$ denotes the height of a vertex $v$ in the $M$-adic tree rooted at $[0,1]$, so that $v$ is identified with an $M$-adic interval of length $M^{-h}$, as discussed in Section \ref{tree encoding section} of Chapter \ref{trees-section}. We include a schematic example of a model $(N, C_0)$-tree in Figure \ref{fig:model-3-tree} with $N=3$.
\begin{figure}[ht]
\centering
\begin{tikzpicture}[x=1.0cm,y=1.0cm, every node/.style={font=\small}]

%---------------------------------
% Column headers
%---------------------------------
\node at (-2.7,6.15) {\bf{Height}};
\node at (0.8,6.15) {\bf{$M$-adic tree}};
\node at (4.7,6.15) {\bf{Height}};

%---------------------------------
% Left/right height labels
%---------------------------------
\node[left]      at (-2.55,5.35) {\scriptsize{$h(v_0)$}};
\node[left,blue] at (-2.55,4.55) {\scriptsize{$h_{v_0}^{\ast}$}};

\node[left]      at (-2.55,3.65) {\scriptsize{$h(v_1)$}};
\node[left,blue] at (-2.55,2.85) {\scriptsize{$h_{v_1}^{\ast}$}};

\node[right]      at (4.55,4.15) {\scriptsize{$h(v_1')$}};
\node[right,blue] at (4.55,3.35) {\scriptsize{$h_{v_1'}^{\ast}$}};

%---------------------------------air
% Shift tree slightly right to center it
%---------------------------------
\begin{scope}[xshift=0.55cm]

% Root splitting vertex v0
\fill (0,5.35) circle (2.2pt);
\draw[red, thick] (0,5.35) circle (0.18);
\node[right] at (0.12,5.48) {\scriptsize{$v_0$}};

% First split from v0
\draw[thick] (0,5.35) -- (-1.7,3.65);
\draw[thick] (0,5.35) -- ( 1.7,4.15);

%---------------------------------
% h_{v0}^* : points exactly on the branches
%---------------------------------
% On left branch from (0,5.35) to (-1.7,3.65), y=4.55 gives x=-0.8
% On right branch from (0,5.35) to (1.7,4.15), y=4.55 gives x≈1.1333
\fill[blue] (-0.80,4.55) circle (1.7pt);
\fill[blue] ( 1.1333,4.55) circle (1.7pt);
\draw[blue, dotted, thick] (-3.15,4.55) -- (1.1333,4.55);

\node[blue,left]  at (-0.72,4.72) {\scriptsize{$w_1(v_0)$}};
\node[blue,right] at (1.05,4.72) {\scriptsize{$w_2(v_0)$}};

% First-order splitting vertices
\fill (-1.7,3.65) circle (2.2pt);
\draw[red, thick] (-1.7,3.65) circle (0.18);
\node[left] at (-1.8,3.78) {\scriptsize{$v_1$}};

\fill (1.7,4.15) circle (2.2pt);
\draw[red, thick] (1.7,4.15) circle (0.18);
\node[right] at (1.8,4.28) {\scriptsize{$v_1'$}};

%---------------------------------
% Second split below v1 (left branch)
%---------------------------------
\draw[thick] (-1.7,3.65) -- (-2.6,2.10);
\draw[thick] (-1.7,3.65) -- (-0.8,2.10);

% h_{v1}^* : points exactly on the branches
% Left edge from (-1.7,3.65) to (-2.6,2.10), y=2.85 gives x≈-2.1645
% Right edge from (-1.7,3.65) to (-0.8,2.10), y=2.85 gives x≈-1.2355
\fill[blue] (-2.1645,2.85) circle (1.7pt);
\fill[blue] (-1.2355,2.85) circle (1.7pt);
\draw[blue, dotted, thick] (-2.15,2.85) -- (-1.2355,2.85);

\node[blue,left]  at (-2.0,3.02) {\scriptsize{$w_1(v_1)$}};
\node[blue,right] at (-1.33,3.02) {\scriptsize{$w_2(v_1)$}};
\draw[blue, dotted, thick] (-3.15,2.85) -- (-2.15,2.85);
%---------------------------------
% Second split below v1' (right branch)
%---------------------------------
\draw[thick] (1.7,4.15) -- (0.9,2.60);
\draw[thick] (1.7,4.15) -- (2.6,2.60);

% h_{v1'}^* : points exactly on the branches
% Left edge from (1.7,4.15) to (0.9,2.60), y=3.35 gives x≈1.2871
% Right edge from (1.7,4.15) to (2.6,2.60), y=3.35 gives x≈2.1645
\fill[blue] (1.2871,3.35) circle (1.7pt);
\fill[blue] (2.1645,3.35) circle (1.7pt);
\draw[blue, dotted, thick] (1.2871,3.35) -- (2.95,3.35);

\node[blue,left]  at (1.38,3.52) {\scriptsize{$w_1(v_1')$}};
\node[blue,right] at (2.0,3.52) {\scriptsize{$w_2(v_1')$}};
\draw[blue, dotted, thick] (2.95,3.35) -- (4.0,3.35);
%---------------------------------
% Second-generation splitting vertices (unlabeled)
%---------------------------------
\fill (-2.6,2.10) circle (2.2pt);
\draw[red, thick] (-2.6,2.10) circle (0.18);

\fill (-0.8,2.10) circle (2.2pt);
\draw[red, thick] (-0.8,2.10) circle (0.18);

\fill (0.9,2.60) circle (2.2pt);
\draw[red, thick] (0.9,2.60) circle (0.18);

\fill (2.6,2.60) circle (2.2pt);
\draw[red, thick] (2.6,2.60) circle (0.18);

%---------------------------------
% Third split from second-generation splitting vertices
%---------------------------------
\draw[thick] (-2.6,2.10) -- (-3.05,0.65);
\draw[thick] (-2.6,2.10) -- (-2.15,0.65);

\draw[thick] (-0.8,2.10) -- (-1.20,0.65);
\draw[thick] (-0.8,2.10) -- (-0.35,0.65);

\draw[thick] (0.9,2.60) -- (0.45,1.15);
\draw[thick] (0.9,2.60) -- (1.35,1.15);

\draw[thick] (2.6,2.60) -- (2.15,1.15);
\draw[thick] (2.6,2.60) -- (3.05,1.15);

% terminal nodes
\fill (-3.05,0.65) circle (2pt);
\fill (-2.15,0.65) circle (2pt);
\fill (-1.20,0.65) circle (2pt);
\fill (-0.35,0.65) circle (2pt);

\fill (0.45,1.15) circle (2pt);
\fill (1.35,1.15) circle (2pt);
\fill (2.15,1.15) circle (2pt);
\fill (3.05,1.15) circle (2pt);

\end{scope}
\end{tikzpicture}
\caption{\small{A schematic model $(3,C_0)$-tree. Every maximal ray contains exactly three splitting vertices, shown in red, and each splitting vertex has two children. For a splitting vertex $v \in \{v_0, v_1, v_1' \}$, the vertices $\{w_i(v): i =1, 2\}$ denote the descendants of $v$ at height $h_v^{\ast}$. The blue dotted lines indicate the auxiliary heights $h_v^{\ast}$ at which Euclidean separation is recorded for these vertices. Each such height occurs after its associated splitting vertex and before any subsequent splitting descendant.}}
\label{fig:model-3-tree}
\end{figure}
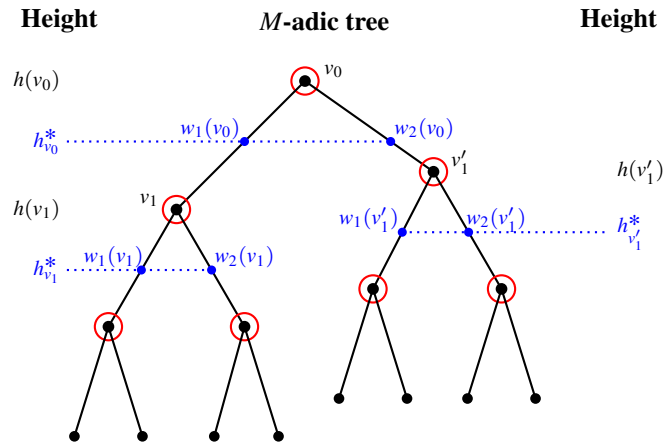
\vskip0.1in
\noindent The separation condition \eqref{Euclidean distance condition}, along with \eqref{what is h_v-star}, means that 
%\[ h(w_1(v)) = h(w_2(v)) = h_v^{\ast}, \] and 
there are at least $C_0$ and at most $(C_0+2)M$ intervals of length $M^{-h_v^{\ast}}$ separating the two intervals $w_1(v)$ and $w_2(v)$, which are also of the same length.  The condition \eqref{h-star-v} implies that the height $h_v^{\ast}$ realizing this separation appears before any further splitting has occurred; namely, any ray of $\mathcal S$ that is rooted at one of the two children of $v$ contains no splitting vertices between the heights $h(v) +1$ and $h_{v}^{\ast}$.  Figure \ref{fig:separation-property} gives a schematic illustration of the Euclidean separation phenomenon for a single splitting vertex $v \in \mathcal S$. 
\vskip0.1in
\noindent A set $U \subseteq [0,1]$ whose $M$-adic tree $\mathcal T(U; M)$ is a model $(N, C_0)$-tree is called a {\em{model $(N, C_0)$-set}}.  
\begin{figure}[ht]
\centering
\begin{tikzpicture}[x=0.9cm,y=1.0cm, every node/.style={font=\small}]

%---------------------------------
% Column headers
%---------------------------------
\node at (-2.6,4.55) {\bf{Height}};
\node at (0.0,4.55) {\bf{Vertices of $\mathcal S$}};
\node at (6.2,4.55) {\bf{Interval}};

%---------------------------------
% Left / middle panel: tree
%---------------------------------

% Height labels
\node[left] at (-3.0,3.70) {\scriptsize{$h(v)$}};
\node[left] at (-3.0,2.95) {\scriptsize{$h(v)+1$}};
\node[left] at (-3.0,1.55) {\scriptsize{$h_v^*$}};

% Guide lines
\draw[dotted] (-2.4,3.70) -- (2.5,3.70);
\draw[dotted] (-2.4,2.95) -- (2.5,2.95);
\draw[dotted] (-2.4,1.55) -- (2.5,1.55);

% Vertex v
\fill (0,3.70) circle (2.2pt);
\draw[red, thick] (0,3.70) circle (0.17);
\node[right] at (0.13,3.78) {\scriptsize{$v$}};

% Edges
\draw[thick] (0,3.70) -- (-1.10,1.55);
\draw[thick] (0,3.70) -- ( 1.10,1.55);

% Descendants
\fill (-1.10,1.55) circle (2.2pt);
\fill ( 1.10,1.55) circle (2.2pt);
\draw[blue, thick] (-1.10,1.55) circle (0.15);
\draw[blue, thick] ( 1.10,1.55) circle (0.15);

\node[blue,left]  at (-1.25,1.55) {\scriptsize{$w_1(v)$}};
\node[blue,right] at ( 1.25,1.55) {\scriptsize{$w_2(v)$}};

%---------------------------------
% Arrow
%---------------------------------
\draw[->, thick] (3.0,2.65) -- (3.9,2.65);

%---------------------------------
% Right panel: interval pictures
%---------------------------------

% Chosen x-positions for the same two points, vertically aligned
\def\xA{5.70}
\def\xB{6.35}

% Top line: v at height h(v)
\draw[thick] (4.8,3.70) -- (7.8,3.70);
\draw[thick] (4.8,3.58) -- (4.8,3.82);
\draw[thick] (7.8,3.58) -- (7.8,3.82);
\node at (6.3,4.02) {\scriptsize{$v$}};

\fill[green!60!black] (\xA,3.70) circle (2.5pt);
\fill[green!60!black] (\xB,3.70) circle (2.5pt);

% Second line: v subdivided into three equal parts at height h(v)+1
\draw[thick] (4.8,2.95) -- (7.8,2.95);
\draw[thick] (4.8,2.83) -- (4.8,3.07);
\draw[thick] (7.8,2.83) -- (7.8,3.07);

% subdivision into thirds
\draw[thick] (5.8,2.83) -- (5.8,3.07);
\draw[thick] (6.8,2.83) -- (6.8,3.07);
\node at (6.3,3.27) {\scriptsize{children of $v$}};

% Same dots, vertically aligned
\fill[green!60!black] (\xA,2.95) circle (2.5pt);
\fill[green!60!black] (\xB,2.95) circle (2.5pt);

% Bottom line: descendants at height h_v^*
\draw[thick] (4.8,1.55) -- (7.8,1.55);

% Same dots, vertically aligned
\fill[green!60!black] (\xA,1.55) circle (2.5pt);
\fill[green!60!black] (\xB,1.55) circle (2.5pt);

% Intervals around descendants, separated at the final level
\draw[green!60!black, thick] (5.45,1.40) rectangle (5.95,1.70);
\draw[green!60!black, thick] (6.10,1.40) rectangle (6.60,1.70);

\node[blue] at (5.50,1.08) {\scriptsize{$w_1(v)$}};
\node[blue] at (6.55,1.08) {\scriptsize{$w_2(v)$}};

% Separation bracket at final level
\draw[blue, thick] (5.95,2.00) -- (6.10,2.00);
\draw[blue, thick] (5.95,1.92) -- (5.95,2.08);
\draw[blue, thick] (6.10,1.92) -- (6.10,2.08);
\node[blue!70] at (6.025,2.27) {\scriptsize{scale $M^{-h_v^{\ast}}$}};

\end{tikzpicture}
\caption{\small{Schematic illustration of Property \ref{separation condition} of a model tree. The same two points (in green) are tracked across the interval pictures on the right. Their youngest common ancestor is $v$, but they lie in adjacent subintervals in the generation immediately after the splitting vertex $v$; the distance between these two sibling intervals is zero. However, a visible Euclidean separation is recorded at the later generation $h_v^*$, where the descendants $w_1(v)$ and $w_2(v)$ are separated. After separating at $v$, the two resulting rays undergo no further splitting through height $h_v^{\ast}$.}}
\label{fig:separation-property}
\end{figure}
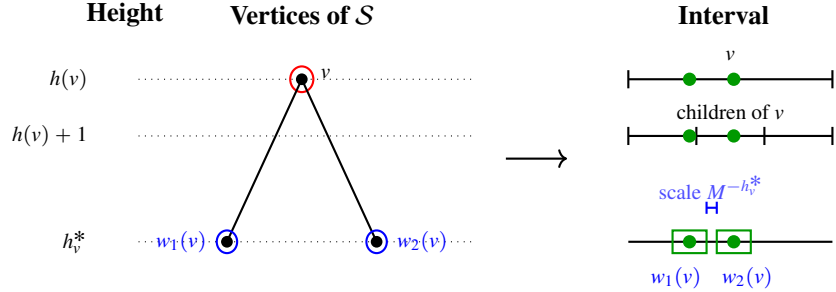

\section{Sublacunary sets contain model subsets of large lacunarity order} \label{section: sublac model subsets}
The main result of this chapter, Proposition \ref{PRUNING STAGE 1}, says that every tree of large splitting number contains a model sub-tree of a comparably large splitting number. Stated in terms of sets, every set of large lacunarity order admits a model subset of comparable lacunarity order. The idea of curating a tree to meet Euclidean constraints originates in Kroc's thesis \cite[Section 4]{KrocThesis}, which we have adapted to our setting. 
\begin{proposition}\label{PRUNING STAGE 1}
	Given a fixed base integer $M \geq 2$, an integer constant $C_0 \geq 1$, and an integer $N \geq 1$, let $\Omega \subseteq [0,1]$ be a set whose $M$-adic tree obeys the hypothesis 
	\begin{equation} \label{pre-pruning split}
	\text{split}\bigl(\mathcal T(\Omega;M) \bigr) > (N+1)(C_0+3). 
	\end{equation} Then the following conclusions hold. 
	\vskip0.1in 
	\begin{enumerate}[(a)]
	\item The set $\Omega$ contains a model $(N, C_0)$ subset $\Omega_N$. In other words, $\mathcal T(\Omega; M)$ has a subtree 
	\begin{equation} \label{S and Omega_N}
	\mathcal S_N^{\ast} := \mathcal T(\Omega_N; M) 
	\end{equation} satisfying the defining criteria \ref{N split per ray}, \ref{two children per split} and \ref{separation condition} of a model $(N, C_0)$-tree stated in Section \ref{section: model trees}. 
	\vskip0.1in
	\item \label{partJN} While $\mathcal S_N^{\ast}$ is of infinite height, one can find an integer $J > N$ such that the truncation $\mathcal S_N = \mathcal S_N[J]$ of $\mathcal S_N^{\ast}$ to height $J$ obeys the following properties:
	\vskip0.1in 
	\begin{itemize} 
	\item The truncated tree $\mathcal S_N$ continues to be a model $(N, C_0)$-tree.
	\vskip0.1in 
	\item Each terminal vertex of $\mathcal S_N$ contains exactly one element of $\Omega_N$: 
	\begin{equation} \label{one element per terminal vertex}
	\#\bigl(\omega \cap \Omega_N \bigr) = 1 \text{ for any } \omega \in \mathcal S_N, \; h(\omega) = J.   
	\end{equation}    
	\vskip0.1in 
	\item The terminal vertices of $\mathcal S_N$ obey the separation condition 
	\begin{equation} \label{Euclidean separation terminal vertices} 
	\min \bigl\{\text{dist}(\omega, \omega') : \omega, \omega' \in \mathcal S_N, \;\omega \neq \omega', \; h(\omega) = h(\omega') = J \bigr\} \geq C_0 M^{-J}. 
	\end{equation} 
	\end{itemize} 
	\end{enumerate}
	\vskip0.1in 
	%Additionally, there exists an integer $J = J(\Omega, N)$ such that 
%\begin{equation} 
%\label{Defn of J} C_0 M^{-J} \leq \min\{ |\omega - \omega'| : \omega \ne \omega', \; \omega, \omega' \in \Omega_N \}. 	
%\end{equation} 	
\end{proposition} 
\noindent The condition \eqref{Euclidean separation terminal vertices} means that the elements of $\Omega_N$ are well-separated relative to the final scale $M^{-J}$, in addition to the crucial Euclidean separation in older generations as identified by the condition \eqref{Euclidean distance condition}. 

\section{A tool for iterative pruning} \label{section: iterative pruning}
The pruning process leading to the outcome claimed in Proposition \ref{PRUNING STAGE 1} is based on an iterative algorithm that prescribes each splitting while ensuring Euclidean separation. The key strategy of the iteration is contained in Lemma \ref{PRUNING BUILDING BLOCK} below. 
\begin{lemma} \label{PRUNING BUILDING BLOCK}
Fix a vertex $v_0$ of the full $M$-adic tree $\mathcal T([0,1]; M)$, and two positive integers $N_0, C_0$ with 
\[ N_0 > C_0+3. \] Let $\mathcal T_{0}$ be a subtree rooted at $v_0$ with no terminal vertices (i.e. every maximal ray of $\mathcal T_0$ is infinite), and with the property that 
\begin{equation} \label{hypotheses on T_0}
\text{every ray in $\partial \mathcal T_{0}$ contains more than $N_0$ splitting vertices}.
\end{equation} 
\vskip0.1in
\noindent Then there exist an integer $K_0 = K_0(v_0) \geq 2$ and a subtree $\mathcal T_{1} \subseteq \mathcal T_{0}$ rooted at $v_0$ such that all the following properties hold: 
\vskip0.1in 
\begin{enumerate}[1.]
\item \label{nonadjacency - two descendants} There are exactly two terminal vertices $w_1$ and $w_2$ in $\mathcal T_1$, with
\begin{equation} \label{height of T_1}  
h(w_1) = h(w_2) = K_0 + h(v_0); 
\end{equation}   
in other words, both vertices lie $K_0$ generations below $v_0$. The unique splitting vertex of $\mathcal T_{1}$, namely $D(w_1, w_2)$, occurs at a height strictly in between that of the root and the terminal vertices: 
\begin{equation} \label{one splitting vertex} 
h(v_0) < h(D(w_1, w_2)) < K_0 + h(v_0). 
\end{equation} 
\vskip0.1in
\item \label{nonadjacency - Euclidean separation} The integer $K_0$ is the smallest positive integer $K$ for which there exist two vertices $u_1, u_2 \in \mathcal T_0$ satisfying 
\[ h(u_1​)=h(u_2​)=h(v_0​)+K, \quad \text{dist}(u_1​,u_2​) \geq C_0 ​M^{-K-h(v_0​)}. \]
The terminal vertices realize this property for $K = K_0$. In particular,  
\begin{align} 
&\text{dist}(w_1, w_2) \geq C_0 M^{-K_0 - h(v_0)}, \text{ and }  \label{Euclidean separation lower} \\
&\text{dist}(w_1, w_2) \leq (C_0 + 2) M^{-K_0 - h(v_0) + 1}.  \label{Euclidean separation upper} 
\end{align} 
\vskip0.1in
\item \label{nonadjacency - remaining splits} If $\mathcal T_{0}(w_i)$ denotes the maximal subtree of $\mathcal T_{0}$ rooted at $w_i$ then each ray in $\mathcal T_{0}(w_i)$ has more than $N_0 - (C_0+3)$ splitting vertices. 
\end{enumerate}    
\end{lemma}
\noindent Let us first see how Lemma \ref{PRUNING BUILDING BLOCK} contributes to the construction of the set $\Omega_N$ depicted in Proposition \ref{PRUNING STAGE 1}. We will complete the proof of Lemma \ref{PRUNING BUILDING BLOCK} in Section \ref{Pruning building block proof section}.   

\section{The pruning process: Proof of Proposition \ref{PRUNING STAGE 1}, given Lemma \ref{PRUNING BUILDING BLOCK}} 
\subsection{Proof overview} 
Given a set $\Omega \subseteq [0,1]$ obeying the hypotheses \eqref{pre-pruning split} of Proposition \ref{PRUNING STAGE 1}, let $\mathcal T$ denote its $M$-adic tree $\mathcal T(\Omega; M)$. This section reduces the proof of  Proposition \ref{PRUNING STAGE 1} to an inductive statement, whose justification is relegated to Sections \ref{base case for pruning} and \ref{inductive step for pruning}. Modulo this statement, the proof of  Proposition \ref{PRUNING STAGE 1} is completed here. 
\vskip0.1in
\noindent The splitting number of a tree is, by definition, the largest splitting number of its vertices. Additionally, the monotonicity of splitting numbers with respect to heights, as proved in Corollary \ref{corollary-monotonicity}, allows us to transfer the large splitting number assumption \eqref{pre-pruning split} from the tree $\mathcal T$ to its root vertex $[0,1]$; namely,  
\begin{equation} \label{large split reduced to the root}   
\text{split}_{\mathcal T}([0,1]) = \text{split}(\mathcal T) > (N+1)(C_0+3). 
\end{equation}  
According to the definition \eqref{splitting vertex}, $\text{split}_{\mathcal T}([0,1])$ is given by a min-max principle. It provides an optimal lower bound on the number of splits per ray that can occur in some subtree rooted at $[0,1]$. The condition \eqref{large split reduced to the root} then ensures the existence of a subtree $\bar{\mathcal T}[N]$ of $\mathcal T$ rooted at $[0,1]$, such that 
\begin{equation} \label{splitting property}
{\text{every ray of $\partial \bar{\mathcal T}[N]$ contains $>(N+1)(C_0 + 3)$ splitting vertices}}. 
\end{equation} 
We may additionally assume that $\bar{\mathcal T}[N]$ has no terminal vertices, i.e., every maximal ray of $\bar{\mathcal T}[N]$ is of infinite length. Indeed, every vertex of $\mathcal T=\mathcal T(\Omega;M)$ lies on an infinite ray of $\mathcal T$. We therefore extend each finite maximal ray of $\bar{\mathcal T}[N]$, from its terminal vertex onward, along a single infinite ray of $\mathcal T$. These extensions introduce no new splitting vertices, and every resulting infinite maximal ray retains the splitting vertices of the original finite maximal ray. Thus \eqref{splitting property} remains valid.
\vskip0.1in
\noindent The subtree $\bar{\mathcal T}[N]$ is pruned using a recursive  scheme that implements Lemma \ref{PRUNING BUILDING BLOCK} at every stage. Specifically, we will prove the following slightly more general statement for every $N \geq 1$:
\begin{equation} \label{inductive pruning}
\left\{
	\begin{minipage}{0.8\textwidth}
{\em{Let $\bar{\mathcal T}[N] \subseteq \mathcal T([0,1]; M)$ be a sub-tree obeying the splitting condition \eqref{splitting property}, rooted at an arbitrary vertex and having no terminal vertices. Then $\bar{\mathcal T}[N]$ has a model $(N, C_0)$-subtree $\bar{\mathcal S}_N$ with finite height and the same root. Further, the root is non-splitting in $\bar{\mathcal S}_N$.}}
\end{minipage}
%\end{aligned} 
\right\}
\end{equation}
The rays of $\partial \bar{\mathcal S}_N$ may have varying lengths. 
\vskip0.1in
\noindent Let us pause for a moment to clarify how the statement \eqref{inductive pruning} leads to the set $\Omega_N$ and the subtrees $\mathcal S_N^{\ast}, \mathcal S_N$ claimed in Proposition \ref{PRUNING STAGE 1}. By the requirement \eqref{N split per ray eq} of a model tree, the pruned tree $\bar{\mathcal S}_N$ consists of $2^N$ rays of maximal length, whose terminal vertices may lie at varying heights. Since $\bar{\mathcal S}_N$ is a subtree of $\mathcal T$, each of the $M$-adic intervals represented by the terminal vertex of such a ray has non-trivial intersection with $\Omega$. Fixing a point $x_v \in v \cap \Omega$ for every terminal vertex $v \in \bar{\mathcal S}_N$ allows us to define the pruned slope set 
\begin{equation}  \label{pruned set} 
\Omega_N := \bigl\{x_v : v \text{ a terminal vertex of } \bar{\mathcal S}_N \bigr\}. 
\end{equation}
Clearly $\#(\Omega_N) = \#(\partial \bar{\mathcal S}_N)=  2^N$. 
\vskip0.1in  
\noindent The tree $\mathcal S_N^{\ast} := \mathcal T(\Omega_N; M)$, which is of infinite height by definition, is derived from $\bar{\mathcal S}_N$ by attaching,  to every terminal vertex $v$ of $\bar{\mathcal S}_N$, a single non-splitting ray $\mathcal R_v$ representing $x_v$:
\[ x_v = \alpha(\mathcal R_v). \] 
Since this extension introduces no new splitting vertices and leaves all vertices of $\bar{\mathcal S}_N$ unchanged, the two trees have the same splitting structure. Moreover, for every splitting vertex $v$, the descendants $w_1(v),w_2(v)$ and the separating height $h_v^{\ast}$ selected in $\bar{\mathcal S}_N$ remain unchanged in $\mathcal S_N^{\ast}$. Consequently, $\mathcal S_N^{\ast}$ inherits all the defining properties \eqref{N split per ray eq}--\eqref{h-star-v} of a model $(N,C_0)$-tree. This, in turn, makes $\Omega_N$ a model $(N, C_0)$-set, which is the conclusion of Proposition \ref{PRUNING STAGE 1}.
\vskip0.1in
\noindent Finally, let us describe how the finitary tree $\mathcal S_N$ is obtained from $\mathcal S_N^{\ast}$ and $\Omega_N$. Since $x_v \neq x_{v'}$ for any two terminal vertices $v, v'$ of $\bar{\mathcal S}_N$, there must exist an integer $J \geq 1$ such that 
\begin{equation} \label{points in Omega_N are distinct}  
\min \bigl\{ |x_v - x_{v'}| : v \neq v' \bigr\} \geq (C_0+2) M^{-J}.  
\end{equation} 
Without loss of generality, we may choose $J$ to be at least as large as the height of $\bar{\mathcal S}_N$. We define $\mathcal S_N$ as the truncation of ${\mathcal S}^{\ast}_N$ to height $J$. Thus $\mathcal S_N$ is a finitary extension of $\bar{\mathcal S}_N$ to height $J$, with no new splits. Since $J$ is at least the height of $\bar{\mathcal S}_N$, this truncation retains every splitting vertex, every separating height $h_v^*$, and the corresponding descendants $w_1(v),w_2(v)$. Hence $\mathcal S_N$ remains a model $(N,C_0)$-tree.
\vskip0.1in
\noindent The choice \eqref{pruned set} of $\Omega_N$ enforces \eqref{one element per terminal vertex}, namely that each terminal vertex of $\mathcal S_N$ contains exactly one element of $\Omega_N$.      
Finally, if $\omega, \omega'$ denote the $M$-adic intervals of length $M^{-J}$ containing $x_v, x_{v'}$ respectively, then the point separation condition \eqref{points in Omega_N are distinct} implies 
\[ \text{dist}(\omega, \omega') \geq |x_v - x_{v'}| - 2M^{-J} \geq C_0 M^{-J}, \]
which is the claimed interval separation condition \eqref{Euclidean separation terminal vertices}.
\vskip0.1in 
\noindent It therefore remains to establish the statement \eqref{inductive pruning}. As is common in this line of reasoning, the proof proceeds by induction on $N$; the base case and the inductive step are carried out in Sections \ref{base case for pruning} and \ref{inductive step for pruning} respectively. 
\vskip0.1in
\noindent \subsection{Proof of \eqref{inductive pruning}: The base case $N=1$.} \label{base case for pruning} 
Let us consider a tree $\bar{\mathcal T}[1]$ where every ray in $\partial \bar{\mathcal T}[1]$ has more than $2(C_0+3)$ splitting vertices. Applying Lemma \ref{PRUNING BUILDING BLOCK} with the choice of parameters 
\[ \mathcal T_{0} = \bar{\mathcal T}[1], \quad v_0 = \text{ root of } \bar{\mathcal T}[1], \; \text{ and } \; N_0 = 2(C_0 + 3) \] yields a subtree of $\mathcal T_0$ rooted at $v_0$ and with $K_0 = K_0(v_0)$ generations. The resulting subtree has exactly two terminal vertices $w_1$ and $w_2$ whose absolute height within the full $M$-adic tree is $i_0 : =K_0(v_0) + h(v_0)$ . We denote this sub-tree of $\bar{\mathcal T}[1]$ by $\bar{\mathcal S}_1$. 
\vskip0.1in
\noindent 
Clearly, $\#(\partial \bar{\mathcal S}_1) = 2$. By part \ref{nonadjacency - two descendants} of Lemma \ref{PRUNING BUILDING BLOCK}, there is exactly one splitting vertex \[ v = D(w_1, w_2) \in \bar{\mathcal S}_{1}, \]  hence each ray of $\bar{\mathcal S}_1$ has exactly one splitting vertex. This confirms that $\bar{\mathcal S}_1$ obeys conditions \ref{N split per ray} and \ref{two children per split} of a model $(1, C_0)$-tree. Part \ref{nonadjacency - Euclidean separation} of Lemma \ref{PRUNING BUILDING BLOCK} states that the two descendants $w_1, w_2$ of $v$ at $M$-adic height $i_0$ obey 
\begin{equation} C_0 M^{-i_0} \leq \text{dist}(w_1, w_2) \leq (C_0 + 2) M^{-i_0+1}, \; \text{ by } \eqref{Euclidean separation lower} \text{ and } \eqref{Euclidean separation upper}. \label{base case upper and lower} \end{equation} 
The minimality of $K_0$ in Lemma \ref{PRUNING BUILDING BLOCK} shows that $i_0$ is the smallest height greater than $h(v)$ at which two descendants of $v$ satisfy \eqref{Euclidean distance condition}. Thus $h_v^*=i_0$. 
\vskip0.1in
\noindent Next, let us verify \eqref{h-star-v} for the only splitting vertex $v = D(w_1, w_2)$ in $\bar{\mathcal S}_1$. The right inequality of the pruning criterion \eqref{one splitting vertex} gives 
\[ h(v) = h \bigl( D(w_1, w_2)\bigr)< K_0 + h(v_0) = i_0 = h_v^{\ast}.\] 
This establishes the left inequality in \eqref{h-star-v}. Since there are no splitting vertices $\subsetneq v$, the right hand inequality of \eqref{h-star-v} imposes no condition, since the relevant set is empty. To summarize,  the tree $\bar{\mathcal S}_1$ satisfies the last defining requirement \ref{separation condition} of a model $(1, C_0)$-tree, with $h_v^{\ast} = i_0$.
\vskip0.1in
\noindent The left inequality in \eqref{one splitting vertex} shows that the root $v_0$ is non-splitting in $\bar{\mathcal S}_1$. This completes the proof of \eqref{inductive pruning} in the base case. 
\vskip0.1in
\noindent \subsection{Proof of \eqref{inductive pruning}: The inductive step} \label{inductive step for pruning} We assume now that the induction statement \eqref{inductive pruning} holds for all $N \leq n-1$. Let us consider a tree $\bar{\mathcal T}[n]$ rooted at a vertex $v_0$ in which every maximal ray contains more than $(n+1)(C_0+3)$ splitting vertices. We apply Lemma \ref{PRUNING BUILDING BLOCK} on this tree with
the following choice of parameters 
\[ \mathcal T_{0} = \bar{\mathcal T}[n], \quad v_0 = \text{ root of } \bar{\mathcal T}[n], \; \text{ and } \quad N_0 = (n+1)(C_0 + 3). \]
As in Section \ref{base case for pruning}, this yields a subtree $\mathcal T_1[n] \subseteq \bar{\mathcal T}[n]$ rooted at $v_0$ with two vertices $w_1$ and $w_2$ located $K_0 = K_0(v_0)$ generations away from $v_0$. As in the base case, 
\[h(w_1) = h(w_2) = i_0 := h(v_0) + K_0, \text{ so that } \;  |w_1| = |w_2| = M^{-i_0}. \] 
These vertices obey the Euclidean distance conditions \eqref{Euclidean separation lower} and \eqref{Euclidean separation upper}. As before, the pruning criterion \eqref{one splitting vertex} confirms that $v_0$ is a non-splitting root of $\mathcal T_1[n]$. 
\vskip0.1in
\noindent Let $\bar{\mathcal T}(w_i)$ denote the maximal subtree of $\bar{\mathcal T} = \bar{\mathcal T}[n]$ rooted at $w_i$. By part \ref{nonadjacency - remaining splits} of Lemma \ref{PRUNING BUILDING BLOCK}, we deduce that each ray in $\partial \bar{\mathcal T}(w_i)$ has more than $n(C_0+3)$ splitting vertices. Moreover, since $\bar{\mathcal T}[n]$ has no terminal vertices, neither does the maximal subtree $\bar{\mathcal T}(w_i)$, for $i=1,2$.
\vskip0.1in
\noindent We may therefore apply the induction hypothesis \eqref{inductive pruning} on each $\bar{\mathcal T}(w_i)$, $i=1,2$, obtaining in the process a model $(n-1, C_0)$-subtree rooted at $w_i$ of finite height
\[\bar{\mathcal S}_{n-1}(w_i) \subseteq \bar{\mathcal T}(w_i). \]
Appending the subtree $\bar{\mathcal S}_{n-1}(w_i)$ to $w_i$ for $i =1, 2$ in $\mathcal T_1[n]$ generates a new tree, which we call $\bar{\mathcal S}_n$. The induction hypothesis ensures that $w_i$ is a non-splitting root of $\bar{\mathcal S}_{n-1}(w_i)$. The two trees $\bar{\mathcal S}_{n-1}(w_i)$ may have different heights, which could possibly result in rays of $\partial \bar{\mathcal S}_n$ with varying lengths. We claim that the tree $\bar{\mathcal S}_n$, which is of finite height by construction, satisfies all the requirements for a model $(n, C_0)$-subtree. Let us verify the defining conditions \eqref{N split per ray}, \eqref{two children per split} and \eqref{separation condition}. 
\vskip0.1in
\noindent Each ray $\mathcal R \in \partial \bar{\mathcal S}_n$ is the concatenation of a ray of $\partial \mathcal T_1[n]$ and a ray of $\partial \bar{\mathcal S}_{n-1}(w_i)$; the former splits exactly once, the latter exactly $(n-1)$ times, and each split results in exactly two children, from the model nature of the respective trees. It follows that $\mathcal R$ splits exactly $n$ times in $\bar{\mathcal S}_n$, and each splitting vertex has exactly two children. This confirms two of the conditions for $\bar{\mathcal S}_n$ to be a model $(n, C_0)$-tree, namely \eqref{N split per ray eq} and \eqref{two children per split eq}. 
\vskip0.1in
\noindent It remains to verify the Euclidean separation criterion \eqref{Euclidean distance condition} and the height condition \eqref{h-star-v} for every splitting vertex $v$ of $\bar{\mathcal S}_n$. We verify this in two parts -- first for the oldest splitting vertex $v \in \bar{\mathcal S}_n$, and then the younger ones. 
\vskip0.1in
\noindent The oldest splitting vertex $v \in \bar{\mathcal S}_n$ is the only splitting vertex of $\mathcal T_1[n]$. For this $v$,  the desired estimate \eqref{Euclidean distance condition} has already been  established through the application of Lemma \ref{PRUNING BUILDING BLOCK} with $h_v^{\ast}$ replaced by $i_0$, via \eqref{base case upper and lower}. 
In fact, the minimality of $K_0$ in Lemma \ref{PRUNING BUILDING BLOCK} shows that $i_0$ is the smallest height greater than $h(v)$ at which two descendants of $v$ satisfy the separation condition \eqref{Euclidean distance condition}. This confirms $h_v^*=i_0$.
\vskip0.1in
\noindent  Lemma \ref{PRUNING BUILDING BLOCK} also gives
\[ h(v) < K_0 + h(v_0) = i_0 = h_v^{\ast}, \text{ which is the left inequality in \eqref{h-star-v}}.\] 
Since $\mathcal T_1[n]$ has no splitting vertex other than $v$, every splitting vertex $w \subsetneq v$ in $\bar{\mathcal S}_n$ therefore satisfies 
\[ h(w) > i_0 = h_v^{\ast}, \text{ which gives the right inequality in \eqref{h-star-v}.} \]  
\vskip0.1in
\noindent If $v$ is not the oldest splitting vertex of $\bar{\mathcal S}_n$, it must belong to one of the trees $\bar{\mathcal S}_{n-1}(w_i)$. Since $\bar{\mathcal S}_{n-1}(w_i)$ is a model $(n-1, C_0)$ subtree on the full $M$-adic tree, the induction hypothesis ensures both \eqref{Euclidean distance condition} and \eqref{h-star-v} for a certain choice of $h_v^{\ast}$. 
\vskip0.1in
\noindent Finally, $v_0$ is the non-splitting root of $\mathcal T_1[n]$, and adjoining the subtree $\bar{\mathcal S}_{n-1}(w_i)$ below $w_i$ does not alter this fact. Thus $\bar{\mathcal S}_n$ has the same root as $\mathcal T_1[n]$, and this root is non-splitting in $\bar{\mathcal S}_n$. This completes the inductive step, and hence the proof of Proposition \ref{PRUNING STAGE 1}. 
\qed

\section{Proof of Lemma \ref{PRUNING BUILDING BLOCK}} \label{Pruning building block proof section}
As we saw in the proof of Proposition \ref{PRUNING STAGE 1}, Lemma \ref{PRUNING BUILDING BLOCK} was the key tool in the pruning process. We prove the lemma in this section.
\begin{proof}
Let $w_0 \in \mathcal T_0$ be a child of $v_0$, and let $\mathcal T_0(w_0)$ be the maximal subtree of $\mathcal T_0$ rooted at $w_0$.  
\subsection{Identifying a generation with more than $C_0$ vertices} According to the hypothesis \eqref{hypotheses on T_0}, the number of splitting vertices on each ray of $\mathcal T_{0}$ is larger than $N_0$. 
Every maximal ray $\mathcal R \in \partial \mathcal T_0(w_0)$, when extended upward to $v_0$, gives rise to a maximal ray $\mathcal R' \in \partial \mathcal T_0$. According to the hypothesis \eqref{hypotheses on T_0}, this extended ray $\mathcal R'$ has more than $N_0$ splitting vertices. Since $h(w_0) = h(v_0)+1$, passing to the portion of $\mathcal R'$ rooted at $w_0$ removes at most the vertex $v_0$, and hence at most one splitting vertex. Therefore, by \eqref{hypotheses on T_0}, every ray $\mathcal R \in \partial \mathcal T_0(w_0)$ contains at least $N_0$ splitting vertices. 
\vskip0.1in
\noindent The last observation means there exists some generation in the tree
$\mathcal T_0(w_0)$ consisting of at least $2^{N_0}$ vertices. 
Indeed, by iterating through the first $N_0$ splitting vertices along each branch, one obtains $2^{N_0}$ pairwise incomparable vertices of $\mathcal T_0(w_0)$. Since $\mathcal T_0$, and hence $\mathcal T_0(w_0)$, has no terminal vertices, each of these incomparable vertices has descendants at every later generation. Passing, if necessary, to a common later generation therefore produces at least $2^{N_0}$ vertices at that generation.
\vskip0.1in
\noindent Our choice of parameters $N_0, C_0$ dictates that  
\[ 2^{N_0} > N_0 > C_0 + 3 > C_0+1, \] so there exists some level in $\mathcal T_0(w_0)$ where the number of vertices in $\mathcal T_0(w_0)$ exceeds $(C_0+1)$ for the first time. Let $h_0$ denote the separation in height from $v_0$ for this occurrence; in other words, $h_0$ is the smallest integer such that $\mathcal T_0(w_0)$ has more than $(C_0+1)$ vertices at absolute height $h(v_0)+h_0$ in the full $M$-adic tree. 
\begin{equation} \label{minimality of h_0}
h_0 := \min \Bigl\{ h \geq 2: \# \{v \in \mathcal T_0(w_0) : h(v) = h + h(v_0)  \} > C_0 + 1 \Bigr\}.  
\end{equation} 
Clearly $h_0 \geq 2$, since the tree $\mathcal T_0(w_0)$ has only one vertex at height $h(v_0)+1$, namely $w_0$. 
\subsection{Existence of many vertices forces Euclidean separation among some} We choose $K_0$ to be the smallest positive integer $K$ with the property that there exist at least two vertices $w_1$ and $w_2$ of $\mathcal T_{0}(w_0)$ obeying the relation 
\begin{equation} \label{defining condition K}  
 h(w_1) = h(w_2) = K + h(v_0), \quad \text{dist}(w_1, w_2) \geq C_0 M^{-K - h(v_0)}. 
\end{equation}  
It is important to justify the existence of $K$ and $K_0$. We will show momentarily in Lemma \ref{nonadjacency} below that the property \eqref{defining condition K} holds for $K = h_0$, hence the set of integers $K$ obeying \eqref{defining condition K} is not vacuous. In particular, $K_0$ exists and  
\begin{equation} \label{bound on K_0}
2 \leq K_0 \leq h_0. 
\end{equation} 
Assuming this fact for now, let us continue with the rest of the proof. 
\subsection{Definition of $\mathcal T_1$} Let us choose any two vertices $w_1$ and $w_2$ obeying \eqref{defining condition K} with $K=K_0$. Let $\mathcal T_{1}$ be defined as the union of the two ancestral paths joining $v_0$ to $w_1$ and $w_2$. Then $\mathcal T_{1} \subseteq \mathcal T_0$ is a subtree rooted at $v_0$, with $K_0$ generations and exactly two terminal vertices $w_1, w_2$. To clarify, the root vertex $v_0$ is non-splitting in $\mathcal T_1$, whose maximal rays are $K_0$-long and of the form 
\[ v_0 \rightarrow w_0 \rightarrow \cdots \rightarrow w_i, \; h(w_i) = h(v_0) + K_0, \quad i=1,2.\] 
This confirms \eqref{height of T_1}.  We will verify the claims in Lemma \ref{PRUNING BUILDING BLOCK} for this tree $\mathcal T_1$. 

\subsection{Height and number of terminal vertices} The tree $\mathcal T_1$ obeys condition (\ref{nonadjacency - two descendants}) of Lemma \ref{PRUNING BUILDING BLOCK}. In fact, the construction gives 
\[ D(w_1, w_2) \subseteq w_0, \text{ so that } h \bigl(  D(w_1, w_2) \bigr) \geq h(w_0) = h(v_0) + 1 > h(v_0). \] 
On the other hand, $w_1 \neq w_2$, so by the definition of a common ancestor,
\[ h(D(w_1, w_2)) < h(w_i) = K_0 + h(v_0).  \] 
This establishes \eqref{one splitting vertex}. 
\subsection{Euclidean separation between $w_1$ and $w_2$} The lower bound \eqref{Euclidean separation lower} on dist$(w_1, w_2)$, required in part (\ref{nonadjacency - Euclidean separation}), is built into the construction, through \eqref{defining condition K}. 
\vskip0.1in
\noindent To obtain the upper bound \eqref{Euclidean separation upper}, let us consider for $i=1,2$ the parent of $w_i$, which we denote by $w_i'$. Then 
\[ h(w_i') = h(v_0) + K_0-1, \quad i = 1,2. \] 
%$w_i'$ is a vertex of $\mathcal T_0$ at height $h(v_0) + K_0-1$. 
It follows from the minimality of $K_0$ that \eqref{defining condition K} fails to hold at this height; therefore, 
\[\text{dist}(w_1', w_2') < C_0 M^{-K_0 - h(v_0) + 1}. \] 
The upper bound on the distance between $w_1'$ and $w_2'$ allows us to control the distance between $w_1$ and $w_2$: 
\begin{align*}  
\text{dist}(w_1, w_2) &\leq \text{diam}(w_1') + \text{diam}(w_2') + \text{dist}(w_1', w_2')  \\
&\leq 2 M^{-K_0 - h(v_0)+1} + C_0 M^{-K_0 - h(v_0) + 1} = (C_0+2) M^{-K_0 - h(v_0)+1},
\end{align*} 
yielding \eqref{Euclidean separation upper}.

\subsection{Splits on rays rooted at $w_1, w_2$}  It remains to verify part (\ref{nonadjacency - remaining splits}) of Lemma \ref{PRUNING BUILDING BLOCK}; specifically that each ray in the maximal subtree $\mathcal T_{0}(w_i) \subseteq \mathcal T_0$ rooted at $w_i$ contains more than $N_0 - (C_0+3)$ splitting vertices. The minimality \eqref{minimality of h_0} of $h_0$ and the bound \eqref{bound on K_0} on $K_0$ imply 
\begin{equation}  \#\bigl\{v  \in \mathcal T_0(w_0) : h(v) = h(v_0) + K_0 - 1\bigr\}  \leq C_0+1. \label{number of vertices at a given height} \end{equation}  
We argue that a bound like \eqref{number of vertices at a given height}, which places a cap on the number of vertices at a given height, also controls the number of splitting vertices on rays of $\mathcal T_0$ up to that height. Namely, for any ray $\mathcal R \in \partial \mathcal T_0(w_0)$, we claim that 
\begin{equation} \label{control on splitting vertices per ray} 
\#\left\{ v \in \mathcal R \; \Biggl| \; \begin{aligned} &v \text{ is a splitting vertex of $\mathcal T_0$,} \\ &h(v) \leq h(v_0) + K_0-2
\end{aligned} \right\}  \leq C_0. 
\end{equation} 
The collection on the left is taken to be empty if $K_0 = 2$, since the range $h(w_0) \leq h(v) \leq h(v_0)+K+0-2$ is vacuous in this case. 
\vskip0.1in
\noindent Let us prove the claim. Indeed, suppose that a ray $\mathcal R \in \partial \mathcal T_0(w_0)$ contains $s$ splitting vertices of height at most $h(v_0)+K_0-2$. Each such splitting vertex has a child not belonging to $\mathcal R$; choosing a descendant of that child at height $h(v_0)+K_0-1$, and adjoining the vertex of $\mathcal R$ at that height, produces at least $s+1$ distinct vertices of $\mathcal T_0(w_0)$ at that height. If $s+1>C_0+1$, this contradicts \eqref{number of vertices at a given height}. Hence $s\le C_0$, proving \eqref{control on splitting vertices per ray}.
%Indeed, each splitting vertex contributes at least two descendants to each subsequent generation; therefore contrary to \eqref{control on splitting vertices per ray}, if there is a ray $\mathcal R$ of $\mathcal T_0$ rooted at $v_0$ that contained more than $C_0+1$ splitting vertices of height at most $h(v_0) + K_0 -2$, then that ray alone would generate at least $C_0+2$ descendants at height $h(v_0) + K_0-1$, contradicting \eqref{number of  vertices at a given height}
\vskip0.1in
\noindent 
%Each splitting vertex of height $\leq K_0-2$ gives rise to at least one new element (different among themselves and distinct from the terminating vertex of the ray) at height $K_0-1$. 
We are now in a position to estimate the number of splitting vertices on a ray of $\mathcal T_0(w_i)$. Let us choose any $\mathcal R_i\in\partial \mathcal T_0(w_i)$, and extend $\mathcal R_i$ upward through $w_0$ and $v_0$ to a ray $\mathcal R\in\partial \mathcal T_0$. By \eqref{hypotheses on T_0}, $\mathcal R$ contains more than $N_0$ splitting vertices. According to \eqref{control on splitting vertices per ray}, at most $C_0$ of these occur at vertices $v \in \mathcal R$ satisfying 
\[ h(w_0) \leq h(v) \leq h(v_0) + K_0-2. \] In addition, $v_0$ contributes at most one splitting vertex, while the two remaining heights $h(v_0)+K_0-1$ and $h(v_0)+K_0$ contribute at most two. Thus $\mathcal R$ contains at most $C_0+3$ splitting vertices through the height of $w_i$. 
%In other words, 
%\[ \max_{\mathcal R \in \partial \mathcal T_0} \#\{v \in \mathcal R:  v \text{ is a splitting vertex of $\mathcal T_0$,} \; h(v) \leq h(v_0) + K_0 \}  \leq (C_0 + 3). \] 
\vskip0.1in
\noindent The vertices $w_1, w_2$ live at height  $h(v_0) + K_0$. Every ray rooted at $w_i$ for $i=1,2$ is thus left with more than $N_0 - (C_0+3)$ splitting vertices of $\mathcal T_0$, which is the conclusion claimed in part \ref{nonadjacency - remaining splits}. This completes the proof of Lemma \ref{PRUNING BUILDING BLOCK}. 
\end{proof} 
\vskip0.1in
\noindent In the proof of Lemma \ref{PRUNING BUILDING BLOCK}, we needed to establish \eqref{defining condition K}, namely the existence of two well-separated $M$-adic intervals of the same length among a large class of such intervals. It remains to prove this fact, which we record as Lemma \ref{nonadjacency}. In simple terms, it is a generalization of the easy geometric property that among three $M$-adic intervals of the same length, at least two must be non-adjacent, i.e. have Euclidean separation. 
\begin{lemma} \label{nonadjacency}
Fix any integer $r$ and a positive integer $C_0$. Let $\mathscr{I}$ be a collection of $M$-adic intervals of length $M^{-r}$, with the property
\begin{equation} \label{interval collection cardinality}  
\#(\mathscr{I}) \geq C_0+2. 
\end{equation} 
Then there exist at least two intervals $I, I' \in \mathscr{I}$ such that 
\[ \text{dist}(I, I') \geq C_0 M^{-r}. \] 
\end{lemma} 
\begin{proof}
Towards a contradiction, suppose that there exists a collection $\mathscr{I}$ of intervals obeying \eqref{interval collection cardinality}, any two members of which are less than $C_0 M^{-r}$-separated. Let $I_1, I_2, \ldots, I_{C_0+2}$ denote $(C_0+2)$ intervals of $\mathscr{I}$ arranged from left to right. In other words, there exist integers $a_k$ such that 
\[ I_k = [a_k, a_{k} + 1] M^{-r}, \; \text{ with } \; a_k + 1 \leq a_{k+1} \text{ for all } 1 \leq k \leq C_0+2. \] 
The distance assumption dist$(I_k, I_{k'}) < C_0 M^{-r}$ means that 
\begin{align}
&0 \leq \bigl[ a_{k'} - (a_k +1) \bigr] M^{-r} < C_0 M^{-r} \nonumber, \text{ or } \\  
\label{difference-akal} 
&1 \leq a_{k'} - a_k \leq C_0 \; \text{ for all } 1 \leq k < k' \leq C_0+2. 
\end{align} 
Let us apply the inequality \eqref{difference-akal} $C_0+1$ times, with $k' = k+1$ and $k =1, \ldots, C_0+1$. Adding them up, we obtain from the left side of the resulting inequality that 
\[ a_{C_0+2} - a_{1} \geq C_0+1.  \]
This contradicts the rightmost inequality in \eqref{difference-akal}, completing the proof.  
%We first treat the case $r=0$, when all the $M$-adic intervals have unit length. The extreme configuration occurs when $(2C_0+2)$-many $M$-adic intervals are packed side by side, all contained inside a larger interval of the form $Q_0 = n + [0, 2C_0+2]$ for some $n \in \mathbb Z$. The central interval $Q = [C_0+1, C_0+2]$ maintains a minimum distance of $C_0$ from the first interval $n + [0,1]$. Rephrasing this, any $M$-adic interval $Q$ with vertices in $\mathbb Z$ and of length one admits at most $(2C_0+1)$ intervals of the same type whose distance from itself is $\leq C_0$. The case of a general $r \geq 0$ follows by scaling $Q_0$ by a factor of $M^{-r}$.    
\end{proof}

\section{The pruning process: an example} \label{section: HRS pruning}
We end this chapter by demonstrating the pruning process, etablished in Proposition \ref{PRUNING STAGE 1}, in the context of a specific example. The sublacunary set to be studied is $\Omega_{\text{HRS}}$ as in \eqref{HRS-example}, crafted out of the building blocks $\Omega_{\text{HRS}}(R)$ given by \eqref{HRS-example-unit}. In Section \ref{section: HRS-tree}, we analysed the tree structure of 
\[ \mathcal T_{\text{HRS}}(R) := \mathcal T \bigl( \Omega_{\text{HRS}}(R); M\bigr), \] establishing in particular that it has splitting number $R$.     
\vskip0.1in
\noindent Given a sequence of positive integers $\{N_j : j \geq 1\}$ satisfying the summability criterion \eqref{N_j summability}, let us set $n_0 := 0$, 
\begin{equation} 
n_j = N_{2j-1} + N_{2j}, \quad \bar{n}_0 := 0, \quad \bar{n}_j := n_1 + \ldots + n_j, \quad 1 \leq j \leq R. 
\end{equation} 
Let us recall from \eqref{eta-zeta} the $N_j$-long sequences $\pmb{\eta}_j = (0, 1, 1, \ldots 1)$ and $\pmb{\zeta}_j = (1, 0, \ldots, 0)$ that were used to define $\Omega_{\text{HRS}}$.  Concatenating two of these strings of the same type, we introduce two new $n_j$-long binary strings,
\begin{align*} 
&{\pmb{\alpha}}_j := \bigl(\pmb{\eta}_{2j-1}, \pmb{\eta}_{2j} \bigr) \; \text{ and } \; {\pmb{\beta}}_j = \bigl( \pmb{\zeta}_{2j-1}, \pmb{\zeta}_{2j} \bigr) \in \{0, 1\}^{n_j}, \text{ so that } \\ 
&{\pmb{\alpha}}_j= (\underbrace{0, 1, \ldots, 1}_{N_{2j-1} \text{entries}}, \underbrace{0, 1, \ldots 1}_{N_{2j} \text{entries}}), \quad {\pmb{\beta}}_j= (\underbrace{1, 0, \ldots, 0}_{N_{2j-1} \text{entries}}, \underbrace{1, 0, \ldots 0}_{N_{2j}. \text{entries}}). 
\end{align*} 
Using ${\pmb{\alpha}}_j$ and ${\pmb{\beta}}_j$ as building blocks, we define the set
\begin{align} 
&\Omega^{\ast}_{\text{HRS}}(R) := \left\{v (\pmb{\varepsilon}) \; \Bigl| \; \begin{aligned} &(\varepsilon_1, \varepsilon_2, \ldots) = (\pmb{\kappa}_1, \pmb{\kappa}_2, \ldots, \pmb{\kappa}_R, 0, 0, \ldots), \text{ where } \\ &{\pmb{\kappa}}_j \in \{0, 1\}^{n_j}, \; {\pmb{\kappa}}_j = \text{ either } {\pmb{\alpha}}_j \text{ or } {\pmb{\beta}}_j, \; 1 \leq j \leq R  \end{aligned} \right\}, \label{HRS-star} \\  
&\text{ where }  v(\pmb{\varepsilon}) :=  \sum_{k=1}^{\infty} \frac{\varepsilon_k}{2^k} \; \text{ for any infinite binary string } \pmb{\varepsilon} = (\varepsilon_1, \varepsilon_2, \ldots). 
\end{align}  
It follows from the definition that 
\begin{align*} 
\#\bigl(\Omega_{\text{HRS}}^{\ast}(R)\bigr) = 2^R, \quad  \Omega_{\text{HRS}}^{\ast}(R) \subsetneq \Omega_{\text{HRS}}(2R) \subseteq \Omega_{\text{HRS}}, \quad \text{ for all } R \geq 1. 
\end{align*}  
Further, every vertex of height $\bar{n}_j$ on the dyadic tree 
\[ \mathcal T^{\ast}_{\text{HRS}}(R) := \mathcal T \bigl( \Omega_{\text{HRS}}^{\ast}(R); 2\bigr) \]  
is a splitting vertex, for $0 \leq j \leq R-1$. These are the only splitting vertices of $ \mathcal T^{\ast}_{\text{HRS}}(R)$.
Each (non-splitting) ray on the tree only moves left after height $\bar{n}_R$. 
\vskip0.1in
\noindent Even though the elements of $\Omega_{\text{HRS}}$ are not well-separated in terms of Euclidean distance as verified in \cite{{HRS2024a},{HRS2024b}}, the next lemma shows that its subset $\Omega_{\text{HRS}}^{\ast}(R)$ inherits the weak Euclidean separation property \eqref{Euclidean distance condition}. 
\begin{lemma} \label{lemma: HRS pruning}
Suppose that 
\begin{equation} \label{HRS splitting vertex}   
v = v(\pmb{\varepsilon}) + [0, 2^{-\bar{n}_j}], \quad \pmb{\varepsilon} = (\pmb{\kappa}_1, \ldots, \pmb{\kappa}_j, 0, \ldots, 0), \; 0 \leq j \leq R-1
\end{equation} 
is a (splitting) vertex of $\mathcal T^{\ast}_{\text{HRS}}(R)$ at height $\bar{n}_j$.  the separating height $h_v^{\ast}$ prescribed by \eqref{Euclidean distance condition} obeys
\begin{equation} \label{lambda for HRS} 
h_{v}^{\ast} = \bar{n}_j + N_{2j+1} + 1, \quad \text{ so that } \quad h(v) = \bar{n}_j < h_v^{\ast} < \bar{n}_{j+1}.
\end{equation}
As a result, for each $R \geq 1$, the set $\Omega_{\text{HRS}}^{\ast}(R)$ defined in \eqref{HRS-star} is a model $(R, 1)$-subset.   
\end{lemma} 
\begin{proof} 
Let us start with the case $j=0$, where $v_0 = [0,1]$. We will verify that  
\begin{equation} h_{v_0}^{\ast} = \bar{n}_0 + N_1 + 1 = N_1+1. \label{HRS-root-separation} \end{equation} 
Let us denote by $w_1(v_0)$ and $w_2(v_0)$ the two descendants of $v_0$ at the height $h = (N_1+1)$ of $\mathcal T_{\text{HRS}}^{\ast}(R)$. These are $M$-adic intervals of length $2^{-N_1-1}$ given by
\begin{align} 
w_1(v_0) &= \sum_{k=2}^{N_1} \frac{1}{2^k} + \bigl[0, 2^{-N_1-1} \bigr] = \frac{1}{2} - 2^{-N_1} + \bigl[0, 2^{-N_1-1} \bigr]  \nonumber \\ &
= \Bigl[ \frac{1}{2} - 2^{-N_1}, \frac{1}{2} - 2^{-N_1-1}\Bigr], \label{w1-root} \\ 
w_2(v_0) &= \frac{1}{2} + 2^{-N_1-1} +  \bigl[0, 2^{-N_1-1} \bigr] = \Bigl[ \frac{1}{2} + 2^{-N_1-1}, \frac{1}{2} + 2^{-N_1}\Bigr].
\label{w2-root} \end{align}
The descriptions \eqref{w1-root}, \eqref{w2-root} of $w_1(v_0), w_2(v_0)$ jointly imply  
\[ \text{dist}\bigl(w_1(v_0), w_2(v_0) \bigr)  = 2^{-N_1}, \] which lies between  $C_0 2^{-h} = 2^{-N_1-1} \text{ and } (C_0+2) 2^{-h+1} = 3 \times 2^{-N_1+1}$. This confirms the Euclidean separation condition \eqref{Euclidean distance condition} with \[ v = v_0, \; M=2, \; C_0=1 \text{ and } h_v^{\ast} = N_1+1, \text{ establishing \eqref{HRS-root-separation}}. \]
\vskip0.1in
\noindent The proof of \eqref{lambda for HRS} for a general splitting vertex is very similar. For $1 \leq j \leq R-1$, let $v$ be a splitting vertex of $\mathcal T^{\ast}_{\text{HRS}}(R)$:
\[ v = v(\pmb{\varepsilon}) + [0, 2^{-\bar{n}_j}] \; \text{ given by } \; \eqref{HRS splitting vertex}, \quad h(v) = \bar{n}_j. \] If $w_1(v), w_2(v)$ are the descendants of this vertex at height $h = \bar{n}_j + N_{2j+1}+1$, then 
\begin{align*}    
w_1(v) -  v(\pmb{\varepsilon}) &=  2^{-\bar{n}_j}\sum_{k=2}^{N_{2j+1}} 2^{-k} + [0, 2^{-h}] \\
&=  2^{-\bar{n}_j-1} - 2^{-h+1} + [0, 2^{-h}] \\ &=  2^{-\bar{n}_j-1} +\bigl[-2^{-h+1}, - 2^{-h} \bigr], \text{ whereas }  \\ 
w_2(v) -  v(\pmb{\varepsilon}) &= 2^{-\bar{n}_j-1} + 2^{-\bar{n}_j - N_{2j+1}-1} + [0, 2^{-h}] \\ 
&= 2^{-\bar{n}_j-1} + [2^{-h}, 2^{-h+1}].
\end{align*} 
Once again we observe that 
\[ \text{dist}\bigl(w_1(v), w_2(v) \bigr)  = 2^{-h+1}, \] 
which verifies \eqref{Euclidean distance condition} with $h_{v}^{\ast}$ as in \eqref{lambda for HRS}, completing the proof.  
\end{proof}
\chapter{Fundamental heights of the pruned tree} \label{Chapter: Fundamental heights} 
 Let us briefly review our progress thus far. Given a slope set $\Omega \subseteq [0,1]$, whose $M$-adic tree has a splitting number exceeding a large multiple of $N$, we now have a pruning mechanism that isolates a model subset $\Omega_N \subseteq \Omega$. The truncated $M$-adic tree $\mathcal S_N$ of $\Omega_N$ is a model tree; it has a special splitting structure that retains a large splitting number at least $N$, and hence comparable to the original tree, while ensuring Euclidean separation among certain vertices. This pruning algorithm was presented in Proposition \ref{PRUNING STAGE 1} of the last chapter.
\vskip0.1in
\noindent 
%The pruned slope set $\Omega_N$ is the basis for the Kakeya-type configuration that we are aiming to construct for the proof of the implication \eqref{condition 3: slopes sublacunary} $\implies$ \eqref{condition 1: Kakeya-type sets} in Theorem \ref{thm:main}. The constituent rectangles of this configuration will be oriented in the directions specified by $\Omega_N$, in a way that promotes intersection in certain regions of the plane. 
This chapter is devoted to a deeper exploration of the pruned and truncated tree $\mathcal S_N$. Its primary objective is to develop the necessary terminology and notation related to the splitting behaviour of $\mathcal S_N$, such as {\em{splitting index}}, {\em{fundamental height}} and {\em{basic slope intervals}}. These concepts are instrumental in the construction of Kakeya-like configurations to follow in Chapter \ref{chapter: sticky rectangles chapter}. 
%For instance, each maximal ray in $\mathcal S_N$ has $N$ splitting vertices, and a {\em{splitting index}} of such a vertex $v$ identifies its position among the splitting vertices on its ray. The {\em{fundamental height}} $\lambda(v)$ of $v$ isolates the $M$-adic height $h_v^{\ast}$, ensured by Proposition \ref{PRUNING STAGE 1}, where the Euclidean distance between the two descendants of $v$ can be precisely quantified. The intervals of $\mathcal S_N$ at a fundamental height are termed {\em{basic slope intervals}}.  

\section{Splitting index} \label{section: splitting index}
The pruned slope tree $\mathcal S_N$ produced by Proposition \ref{PRUNING STAGE 1} looks like the full binary tree of height $N$, though possibly elongated and asymmetric. Maximal rays in this tree, which consist of $J$ consecutive edges, may have long segments with no splits. However, only the splitting vertices of $\mathcal S_N$ and certain other vertices related to these are of central importance to the subsequent analysis. With this in mind and to aid in quantification later on, we introduce the class of splitting vertices of $\mathcal S_N$: 
\begin{align} 
{\tt{SplitV}} = {\tt{SplitV}}(\mathcal S_N) &:= \bigcup_{j=1}^N \tt{SplitV}_j, \text{ where for every $1 \leq j \leq N$, } \label{splitting vertex collection} \noindent \\ 
{\tt{SplitV}}_j = {\tt{SplitV}}_j(\mathcal S_N) &:= \bigcup_{\mathcal R \in \partial \mathcal S_N} \left\{v: v {\text{ is the $j^{\text{th}}$ splitting vertex on $\mathcal R$}} 
%\begin{aligned} 
%&\text{ there exists $v \in \Omega_N$ such that $\gamma$ is the $j$th splitting} \\  &\text{ vertex on the ray identifying $v$ in $\mathcal T_J(\Omega_N;M)$} 
%\end{aligned} 
\right\}. \label{generation j splitting vertices}
\end{align} 
In other words, $\tt{SplitV}_j$ is the collection of all splitting vertices in the pruned tree $\mathcal S_N$ that have exactly $(j-1)$ splitting ancestors. 
The elements in ${\tt{SplitV}}_j$ will be termed the {\em{$j^{\text{th}}$ splitting vertices}} of $\mathcal S_N$. 
\vskip0.1in
\begin{itemize} 
\item Each vertex $v \in {\tt{SplitV}}$ has exactly two children in $\mathcal S_N$, one to the left of the other on the real line. We call these {\em{the $0^{\text{th}}$ and the $1^{\text{st}}$ child}} of $v$ respectively.  
\vskip0.1in
\item Similarly, if $v \in {\tt{SplitV}}_j$, $1 \leq j \leq N$, then there are exactly two vertices $v_1, v_2 \in   {\tt{SplitV}}_{j+1}$ with $v_1, v_2 \subsetneq v$. Said differently, 
\begin{equation} \label{2 splitting children per vertex} 
\# \bigl\{ w \in {\tt{SplitV}}(\mathcal S_N): w \subsetneq v, \; \iota(w) = \iota(v) + 1 \bigr\} = 2. 
\end{equation} This in turn means that 
\begin{align} 
&\# \bigl\{ w \in \mathcal S_N : w \subseteq v, \; h(w) = h \bigr\} = 2  \nonumber \\ &\hskip0.3in \text{ for all heights $h$ with $h(v) < h \leq \min \bigl\{h(v_1), h(v_2) \bigr\}$}. \label{nonsplitting heights} \end{align} 
In other words, there are exactly two descendants of $v$ at all the heights $h$ given by \eqref{nonsplitting heights}. As in the earlier item, we call these the  {\em{the $0^{\text{th}}$ and the $1^{\text{st}}$ descendant}} of $v$ at this height, the former being to the left of the latter. 
\vskip0.1in
\item For a fixed $v \in  {\tt{SplitV}}_j$, the two descendants of $v$ at any of the heights $h$ in \eqref{nonsplitting heights} represent $M$-adic intervals of the same length. It is important to note that the next splitting descendants $v_1, v_2$ of $v$ need not be of the same length. 
\vskip0.1in
\item On the other hand, the vertices $v \in {\tt{SplitV}}_j$ may occur at different heights of $\mathcal S_N$ even for a fixed $1 \leq j \leq N$, as dictated by the pruning mechanism. As a result, ${\tt{SplitV}}_j$ may contain $M$-adic intervals of varying sizes. Thus the index $j$, which encodes the number of splitting vertices on the ray leading up to and including $v \in {\tt{Split}}_j$, should not be confused with the height $h(v)$ of $v$ in $\mathcal S_N$. 
\end{itemize} 
\vskip0.1in
\noindent  Given $v \in {\tt{SplitV}}$, we write 
\begin{equation} \label{defn iota} 
\iota(v) = j \quad \text{ if } v \in {\tt{SplitV}}_j,
\end{equation} 
and refer to $\iota(v)$ as the {\em{splitting index}} of $v$. The maximal subtree of $\mathcal S_N$ rooted at $v \in {\tt{SplitV}}_j$ contains $2^{N-j}$ maximal
rays, with $(N-j)$ splitting vertices per ray. One may therefore alternatively view $\iota(v)$ as 
\[ N-\iota(v) = \text{split}_{\mathcal S_N}(v), \text{ the latter defined as in \eqref{splitting vertex}}.\] 
%is the splitting number of $\gamma$ with respect to $\mathcal T_J(\Omega_N;M)$,  
Note that ${\tt{SplitV}}_1$ consists of a single element, namely the unique splitting vertex of $\mathcal S_N$ of minimal height.  In general, 
\begin{equation} \label{number of generation j splitting vertices} 
\#({\tt{SplitV}}_j) = 2^{j-1}, \text{ i.e., there are $2^{j-1}$ splitting vertices of index $j$}. 
\end{equation}  
%We declare $\mathcal G_{N+1}(\Omega_N) \equiv \Omega_N$.  

\section{Fundamental heights} \label{section: fundamental heights} 

\subsection{The fundamental height of a splitting vertex} 
The {\em{fundamental height}} of a  vertex $v \in {\tt{SplitV}}(\mathcal S_N)$ is defined as
\begin{equation} \label{defn: fundamental height}
\lambda(v)  := h_v^{\ast}, 
\end{equation} 
where $h_v^{\ast}$ is the separating height of $v$ in the model $(N, C_0)$-tree $\mathcal S_N$, whose existence and properties are ensured by the conditions \eqref{Euclidean distance condition} and \eqref{h-star-v}. It is, by definition, the smallest height in $\mathcal S_N$ strictly larger than $h(v)$ where the two descendants of $v$ are separated, in the Euclidean sense,  by a distance that is comparable with their lengths. 
\vskip0.1in
\noindent The collection of all fundamental heights of $\mathcal S_N$ will be denoted by {\tt{FundHt}}:
\begin{equation} 
{\tt{FundHt}} = {\tt{FundHt}}(\mathcal S_N) := \bigl\{\lambda(v) : v \in {\tt{SplitV}}\bigr\}. 
\end{equation} 
\subsection{Euclidean distance between distinct pruned slopes} \label{section: Euclidean distance + splitting descendants}
The fundamental height is the key determinant of Euclidean distance between any two slopes in $\Omega_N$, as the next lemma shows. Suppose that 
\begin{equation}  \label{yca} \omega, \omega' \in \Omega_N, \; \omega \neq \omega';  \quad \text{ then }  v = D(\omega, \omega') \in {\tt{SplitV}}(\mathcal S_N), \end{equation} 
where $D(\omega, \omega')$ denotes the youngest common ancestor of $\omega, \omega'$ in $\mathcal S_N$. This is therefore also the youngest common ancestor of $\omega$ and $\omega'$ in the full $M$-adic tree.    
\begin{lemma}[Pruned slopes are separated at a scale given by the fundamental height of their youngest common ancestor] 
\label{lemma: distance between pruned slopes} 
There are two positive constants $C_1, C_2$ depending only on $C_0$ and $M$ such that  
\begin{equation} \label{distance between slopes}
C_1 M^{- \lambda(v)} \leq |\omega - \omega'| \leq C_2 M^{-\lambda(v)},
\end{equation} 
for any choice of $\omega, \omega' \in \Omega_N$, and $v$ as in \eqref{yca}. The constants $C_1, C_2$ can be chosen as $C_1=C_0$ and $C_2 = M(C_0+2)+2$. 
\end{lemma} 
\begin{proof} 
Since $\omega \neq \omega'$ and $v = D(\omega, \omega')$ by \eqref{yca}, the slopes $\omega, \omega'$ must be contained in distinct children of $v$, and therefore in distinct descendants of $v$ at every scale $h > h(v)$.  Without loss of generality, let us assume  
\[ \omega \in w_1(v), \; \omega' \in w_2(v), \]
where $w_i(v)$ are the two descendants of $v$ at the fundamental height $\lambda(v) := h_v^{\ast} > h(v)$. It then follows from \eqref{Euclidean distance condition} and \eqref{h-star-v} that 
\begin{equation}  |\omega - \omega'| \geq \text{dist}(w_1(v), w_2(v))  \geq C_0 M^{- \lambda(v)}. \label{omega-omega'-lower-bound} \end{equation}  
Conversely, 
\begin{align}
|\omega - \omega'| &\leq \text{dist}\bigl(w_1(v), w_2(v) \bigr) + \text{diam}(w_1(v)) + \text{diam}(w_2(v)) \nonumber \\ 
&\leq (C_0+2) M^{- \lambda(v) + 1} + 2M^{-\lambda(v)} = \Bigl[ M(C_0+2) + 2 \Bigr] M^{- \lambda(v)}.   \label{omega-omega'-upper-bound}
\end{align}  
Combining \eqref{omega-omega'-lower-bound} and \eqref{omega-omega'-upper-bound} establishes the inequality \eqref{distance between slopes} with 
\[ C_1 = C_0, \quad C_2 = M(C_0+2) + 2, \]
completing the proof. 
\end{proof}

\section{Basic slope intervals} \label{section: basic slopes} 
\noindent The notions of splitting index and fundamental heights lead naturally to the concept of basic slopes. Loosely speaking, these are the vertices of the pruned slope tree $\mathcal S_N$ that correspond to a fundamental height. 
\vskip0.1in
\noindent As a consequence of \eqref{h-star-v} and \eqref{nonsplitting heights}, any splitting vertex $v \in {\tt{SplitV}}_j(\mathcal S_N)$ has exactly two descendants in the tree $\mathcal S_N$ at height $\lambda(v) = h_v^{\ast}$; these are called the {\em{$j^{\text{th}}$ basic slope intervals}} descended from $v$. 
\vskip0.1in
\noindent The collection of all $j^{\text{th}}$ basic slope intervals as $v$ ranges over the $2^{j-1}$ splitting vertices of ${\tt{SplitV}}_j$ is termed ${\tt{BasicSl}}_j$. More precisely,  
\begin{equation} \label{jth basic slope cubes} 
{\tt{BasicSl}}_j := \left\{ \theta \in \mathcal S_N \, \bigl| \exists v \in {\tt{SplitV}}_j  \text{ such that } \theta \subsetneq v, \; h(\theta) = \lambda(v)  \right\}. 
\end{equation} 
The definition of a model tree implies that every $v \in {\tt{SplitV}}_j$ contributes exactly two intervals to ${\tt{BasicSl}}_j$, both of the same length. Hence \eqref{number of generation j splitting vertices} implies 
\begin{equation}  \#\bigl( {\tt{BasicSl}}_j \bigr) = 2^{j}.  \label{basic slope count j} \end{equation} 
%We adopt the convention 
%\begin{equation} \label{basic slope convention}  
%{\tt{BasicSl}}_0 := {\tt{SplitV}}_1, \quad {\tt{BasicSl}}_{N+1} := \bigl\{ v \in \mathcal S_N : h(v) = J \bigr\}. 
%\end{equation} 
Let us note that for $1 \leq j \leq N-1$, every $j^{\text{th}}$ basic slope interval $\theta$ with 
\begin{equation}  \label{split-basic-split}
\theta \subsetneq v \in {\tt{SplitV}}_j \text{ uniquely identifies $\bar{v} \in {\tt{SplitV}}_{j+1}$, with $\bar{v} \subsetneq \theta \subsetneq v$.} \end{equation}  In other words, the ray in $\mathcal S_N$ joining $v$ and $\theta$ permits a unique extension terminating at a splitting vertex $\bar{v}$ of splitting index $\iota(\bar{v}) = \iota(v)+1$. 
\vskip0.1in
\noindent If $1 \leq j \leq N-1$, the next splitting descendant $\bar{v}$ of $v$ goes on to generate two intervals in ${\tt{BasicSl}}_{j+1}$. Thus, 
\[ \# \Bigl\{ \theta' \in {\tt{BasicSl}}_{j+1} : \theta' \subsetneq \theta \Bigr\} = 2 \text{ for all } \theta \in {\tt{BasicSl}}_j, \; 1 \leq j \leq N-1.  \] 
In other words, each vertex in ${\tt{BasicSl}}_j$, which is an $M$-adic interval, contributes exactly two vertices, i.e. sub-intervals, of the same length to  ${\tt{BasicSl}}_{j+1}$. This dyadic splitting of basic slopes will be exploited heavily in subsequent analysis; see Chapter \ref{chapter: compressions} for specific applications. 
%Note that every $j^{\text{th}}$ basic slope interval $\theta$ strictly contained within $v \in {\tt{SplitV}}_j$ with $h(\theta) = \lambda(v)$ has to fall in one of two categories: 
%\vskip0.1in
%\begin{itemize} 
%\item Either $\theta \in {\tt{SplitV}}_{j+1}$, i.e. $\theta$ is itself a splitting vertex of index $(j+1)$, or 
%\vskip0.1in
%\item It uniquely identifies such a vertex, in the sense that 
%\begin{equation} \label{BasicSl-FundHt+1}  \text{there is a single $\gamma \in {\tt{SplitV}}_{j+1}$, with $\gamma \subsetneq \theta \subsetneq v$, $h(\theta) = \lambda(v)$}. \end{equation}  In other words, the ray in $\mathcal S_N$ joining $v$ and $\theta$ permits a unique extension terminating at a splitting vertex of index $(j+1)$.
%\end{itemize}
%\vskip0.1in  
%In either event, we say that {\em{$\gamma \in {\tt{SplitV}}_{j+1}$ is identified by $\theta \in {\tt{BasicSl}}_{j}$}}. 

\section{Chains of fundamental heights}  \label{section: fundamental height chains} 
In this section, we will use the notion of fundamental heights developed in Section \ref{section: fundamental heights} to define chains of such heights that are realized in the pruned slope tree $\mathcal S_N$. Each such chain will be used in the next section to define a compressed tree representation of $[0,1]$. 
\vskip0.1in
\noindent Given an $(N-1)$-long binary sequence 
\[ \pmb{\varepsilon}_{N-1} := (\varepsilon_1, \ldots, \varepsilon_{N-1}) \in \{0, 1\}^{N-1}, \] 
we define a chain of fundamental heights 
\begin{equation}  \label{chain of fundamental heights} 
 \pmb{\lambda}(\pmb{\varepsilon}_{N-1}) := \bigl(\lambda_2(\pmb{\varepsilon}_1), \lambda_3(\pmb{\varepsilon}_2), \cdots, \lambda_N(\pmb{\varepsilon}_{N-1}) \bigr) \in {\tt{FundHt}}^{N-1}
\end{equation}  
 using the following inductive procedure: in view of properties \eqref{generation j splitting vertices} and \eqref{number of generation j splitting vertices},  ${\tt{SplitV}}_1$ consists of a single vertex $v_0$; we set 
\begin{equation} \label{defn: lambda_1} 
\lambda_1 := \lambda(v_0) \text{ for } v_0 \in {\tt{SplitV}}_1. 
\end{equation}  
%In view of the properties \eqref{generation j splitting vertices} and \eqref{number of generation j splitting vertices}, ${\tt{SplitV}}_1$ consists of a single vertex $v_0$, therefore $\lambda_1$ is uniquely well-defined. 
%\vskip0.1in
 Let us recall from \eqref{2 splitting children per vertex} that 
\[ \text{$v_0 \in {\tt{SplitV}}_1$ has exactly two descendants in ${\tt{SplitV}}_2$}, \] 
descended from the $0^{\text{th}}$ and the $1^{\text{st}}$ child of $v_0$ respectively. For $\varepsilon_1 \in \{0, 1\}$, we set  
\begin{align} 
&\lambda_2(\varepsilon_1) := \lambda(v(\varepsilon_1)), \quad v(\varepsilon_1) \in {\tt{SplitV}}_2, \\ 
& v(\varepsilon) \subseteq \text{$\varepsilon^{\text{th}}$ child of $v_0 \in {\tt{SplitV}}_1$, for $\varepsilon = 0, 1$.} \nonumber
\end{align} 
In general, for $1 \leq j \leq N-1$, and $\pmb{\varepsilon}_j \in \{0, 1\}^j$, suppose that $v(\pmb{\varepsilon}_j)$ is the splitting vertex of $\mathcal S_N$ determined by the nested sequence: 
\begin{equation}  \label{split-nesting-1}
v(\pmb{\varepsilon}_j) \subsetneq  v(\pmb{\varepsilon}_{j-1}) \subsetneq \ldots \subsetneq v(\pmb{\varepsilon}_1) \subsetneq v_0, 
\end{equation}
where for $1 \leq k \leq j$, and $\varepsilon_k \in \{0, 1\}$,
\begin{equation} \label{split-nesting-2} 
v(\pmb{\varepsilon}_k) \in {\tt{SplitV}}_{k+1} \text{ is descended from the $\varepsilon_k^{\text{th}}$ child of  $v(\pmb{\varepsilon}_{k-1}) \in {\tt{SplitV}}_{k}$.} 
\end{equation}  
The quantity $\lambda_{j+1}$ is a positive integer defined on all $j$-long binary sequences in the following way: 
\begin{equation} \label{defn: lambda_j}
\lambda_{j+1}: \{0, 1\}^j \rightarrow \mathbb N, \quad \lambda_{j+1}(\pmb{\varepsilon}_j) = \lambda \bigl(v(\pmb{\varepsilon}_j) \bigr).
\end{equation} 
A basic slope interval $\theta \in {\tt{BasicSl}}_{j+1}$ with 
\begin{equation} \label{basic slope - fundamental ht}
v(\pmb{\varepsilon}_{j}) \subsetneq \theta \subsetneq v(\pmb{\varepsilon}_{j-1}) \text{ must be of length } |\theta| = M^{-\lambda_{j}(\pmb{\varepsilon}_{j-1})}. 
\end{equation}  
\vskip0.1in
\noindent A chain of fundamental heights is naturally associated with a nested sequence of basic slope intervals. These chains of basic slopes ${\tt{BasicSl}}_j$ permit an alternative tree representation of $\Omega_N$ as a full binary tree of height $N$. We describe this in the next chapter.

\section{The model set $\Omega_{\text{HRS}}^{\ast}(R)$ revisited}
 Let us recall from \eqref{HRS-example} the important example $\Omega_{\text{HRS}}$, and from Section \ref{section: HRS pruning} and Lemma \ref{lemma: HRS pruning} the pruned tree $\mathcal T^{\ast}_{\text{HRS}}(R)$ corresponding to its model subset $\Omega^{\ast}_{\text{HRS}}(R)$ given by \eqref{HRS-star}.  We conclude this chapter with an explicit computation of the fundamental heights and basic slopes of $\mathcal T^{\ast}_{\text{HRS}}(R)$.
\vskip0.1in
\noindent It follows from Lemma \ref{lemma: HRS pruning} that $\mathcal T^{\ast}_{\text{HRS}}(R)$ has exactly $R$ fundamental heights: 
\[ \lambda(v) = \bar{n}_j + N_{2j+1} + 1 \text{ for } v \in {\tt{SplitV}}_j \text{ given by } \eqref{HRS splitting vertex}, \; 1 \leq j \leq R. \] 
It is worth noting that, in general, the number of fundamental heights in a pruned tree can be as large as $2^{R-1}$. The small number of fundamental heights in $\mathcal T^{\ast}_{\text{HRS}}(R)$ reflects the symmetry of the tree. The two basic intervals in ${\tt{BasicSl}}_j$ descended from $v = v(\pmb{\varepsilon})$ are:  
\begin{align*} 
w_1(v) &= v(\pmb{\varepsilon}) + 2^{-\bar{n}_j-1} + 2^{- \lambda(v)}[-2, -1], \\  
w_2(v) &= v(\pmb{\varepsilon}) + 2^{-\bar{n}_j-1} + 2^{- \lambda(v)}[1, 2].
\end{align*} 
%{\bf{\color{magenta} Should we include the binary compression of $\mathcal T^{\ast}_{\text{HRS}}$ here, as a precursor to the next chapter?}}

\chapter{Compressed root and slope trees} \label{chapter: compressions}
\section{Chapter overview}
Let us summarize the situation so far, and present a roadmap for this chapter. Starting from Chapter \ref{trees-section}, we have developed a tree-based framework for analysing slope sets, culminating in the pruned slope tree $\mathcal S_N$ in Chapter \ref{Chapter: Pruning of the slope tree} and a study of its finer properties in Chapter \ref{Chapter: Fundamental heights}. Our goal remains to establish the implication 
\[ \text{\eqref{condition 3: slopes sublacunary} $\implies$ \eqref{condition 1: Kakeya-type sets} in Theorem \ref{thm:main}}, \]
namely that sublacunarity of a slope set $\Omega$ gives rise to Kakeya-type sets composed of rectangles with slopes from $\Omega$. 
\vskip0.1in
\noindent This chapter marks a conceptual pivot in that argument. It bridges the structural analysis of slope sets carried out in Chapters \ref{trees-section}--\ref{Chapter: Fundamental heights} with the geometric construction of Kakeya-type sets that will be undertaken in Chapter \ref{chapter: sticky rectangles chapter}.
\vskip0.1in
\noindent Up to this point, the slope set $\Omega$, and in particular its pruned model $\mathcal S_N$, have been encoded using $M$-adic trees. While this representation is natural and precise, it retains more information than is needed for the purposes of the Kakeya construction. What ultimately drives the geometry of the argument is not the full $M$-adic resolution, but the pattern of splitting exhibited by the pruned tree.
%Thus far, the tree representations related to the sublacunary slope set $\Omega$, such as the slope tree $\mathcal T(\Omega; M)$ or its pruned version $\mathcal S_N$, have been $M$-adic trees. These are sub-trees of $\mathcal T([0,1]; M)$ as defined in Section \ref{tree encoding section}. It encodes a set at all scales - namely, any $k^{\text{th}}$ generation vertex of such a tree is an $M$-adic interval of length $M^{-k}$. For the purposes of constructing Kakeya-type configurations, this level of resolution is unnecessarily fine. What governs the geometry of the construction is not the full $M$-adic structure, rather the pattern of splitting in the pruned slope tree.
\vskip0.1in
\noindent The aim of this chapter is to isolate this essential structure and recast it in a more economical form. To this end, we introduce a compressed representation of the slope tree, in which the hierarchy is governed by splitting events rather than fixed scales. This shift retains the combinatorial features relevant to the construction while discarding inessential detail. The strategy of this chapter proceeds in three steps:
%For the first time, in this chapter, we depart from the $M$-adic framework, and offer a different description of trees subordinate to the splitting structure of the pruned tree $\mathcal S_N$. 
\vskip0.1in
\begin{itemize} 
\item First, we encode the pruned $M$-adic slope tree $\mathcal S_N$ as a compressed tree $\mathscr{B}_N$ that captures its splitting structure. 
\vskip0.1in
\item Second, we construct compatible compressed versions of the root tree. 
\vskip0.1in
\item Finally, we define a family of sticky maps linking the compressed root trees with the compressed slope tree $\mathscr{B}_N$. which will later be used to assign slopes to spatial locations.
\end{itemize} 
\vskip0.1in
\noindent These constructions provide a new framework for relating spatial locations (“roots”) to slope directions in a manner adapted to the pruned geometry. This framework will serve as the foundation for the Kakeya-type construction in Chapter \ref{chapter: sticky rectangles chapter}.
%This compressed slope tree, in turn, is instrumental in developing a family of compressed and truncated versions of the root tree $\mathcal T([0,1]; M)$, and in crafting a new family of sticky maps, each mapping a compressed root tree into the compressed slope tree. This new formulation of roots and slopes will be the foundation for the Kakeya-type construction to follow.  
\vskip0.1in
\noindent The primary objective of this chapter is therefore twofold, with compression being the common theme. 
%A generation in the compressed tree will correspond to a specific splitting index and the scales of the $M$-adic intervals occurring in that generation will reflect the fundamental heights of that splitting index.     
%\vskip0.1in
%\noindent In Chapter \ref{section: split implies lacunarity}, we established that the finiteness of lacunarity order of a slope set is equivalent to the finiteness of the splitting number of its $M$-adic tree. Thus the sublacunary slope set given by the hypothesis \eqref{condition 3: slopes sublacunary} of Theorem \ref{thm:main} is represented by an $M$-adic tree with infinite splitting number. Equivalently stated, the $M$-adic slope tree contains sub-trees whose splitting numbers are arbitrarily large. Let $N$ be a running index quantifying increasingly large splitting numbers. 

 \subsection{Re-framing the model set $\Omega_N$ as a full binary tree} 
Let us describe how the pruned slope tree $\mathcal S_N$, introduced in Chapter \ref{Chapter: Pruning of the slope tree} and analysed further in Chapter \ref{Chapter: Fundamental heights}, can be recast in a compressed form that reflects only its essential combinatorial structure.
 \vskip0.1in
\noindent The key observation is that the relevant hierarchy in $\mathcal S_N$ is not determined by all ambient $M$-adic scales, but by the specific locations and scales at which splitting occurs. The terminology developed in Sections \ref{section: splitting index}, \ref{section: fundamental heights}, \ref{section: basic slopes} and \ref{section: fundamental height chains} -- particularly the notions of splitting index, fundamental heights, and basic slope intervals -- allows us to isolate these critical levels.
\vskip0.1in
\noindent Using this structure, we associate to $\mathcal S_N$ a new tree $\mathscr{B}_N$, whose vertices encode the basic slope intervals determined by the fundamental heights. The tree $\mathscr{B}_N$ is organized so that each generation corresponds to a successive splitting stage in $\mathcal{S}_N$ along with the critical Euclidean separation condition, rather than to a fixed $M$-adic scale. In this sense, $\mathscr{B}_N$ is a compressed tree in the terminology of Section \ref{section: compressed trees}.
\vskip0.1in
\noindent This representation retains the hierarchical organization of the slope set while suppressing intermediate scales. In particular, vertices of $\mathscr{B}_N$ correspond to slope intervals of varying lengths, and the ancestry relation reflects inclusion among these intervals, as inherited from the original $M$-adic structure. The resulting non-$M$-adic tree $\mathscr{B}_N$ may be viewed as a combinatorial skeleton of the pruned slope tree $\mathcal S_N$. Not only does it precisely capture the pattern of branching that governs the subsequent construction, it focuses only on the basic slope intervals of $\mathcal S_N$ which have been created to embody key Euclidean separation.  The precise properties of $\mathscr{B}_N$, along with its structural properties, is given in Section \ref{section: binary tree recast}.

% The terminology formulated in Sections \ref{section: splitting index}, \ref{section: fundamental heights}, \ref{section: basic slopes} and \ref{section: fundamental height chains}  plays a crucial role in recasting the model set $\Omega_N$ and its $M$-adic tree $\mathcal S_N$ as a full binary tree $\mathscr{B}_N$ of height $N$. The tree $\mathscr{B}_N$ is not a traditional $M$-adic tree, although it encodes similar geometric information. Rather, {\em{$\mathscr{B}_N$ is a compressed tree in the sense of Section \ref{section: compressed trees}, where each generation consists of basic slope intervals associated with vertices of the same splitting index in $\mathcal S_N$}}. Intervals that are of the same generation in $\mathscr{B}_N$ may have different lengths. This analysis is carried out in Section \ref{section: binary tree recast}. Any perceived disadvantage resulting from variable dimensions at the same height is offset by the key Euclidean separation property enjoyed by its vertices. 

\subsection{Constructing compressed root trees and inductively sticky maps} 
Having recast the slope set $\Omega_N$ in terms of the compressed tree $\mathscr{B}_N$, we now turn to the corresponding structure on the spatial side and the mechanism that links the two.
\vskip0.1in
\noindent In the geometric construction to follow, slope directions must be assigned to a collection of spatial locations (or ``roots'') in a manner that reflects the hierarchical structure of $\mathscr{B}_N$. To organize these allocations, we introduce a family of root trees, whose vertices encode spatial positions at varying scales. As with the slope tree, it is convenient to work with a compressed representation of roots in which the hierarchy reflects only the relevant structural features inherited from $\mathscr{B}_N$.
\vskip0.1in
\noindent The key object connecting the slope and root trees is a class of maps that assign  a slope to each root in a manner consistent with the respective tree structures. Informally, such a map associates to each vertex of a root tree a vertex of the slope tree of the same length, preserving ancestry: descendants are mapped to descendants. As described in Section \ref{sticky maps section}, we refer to such maps as sticky maps, emphasizing that once two vertices are associated at a given level, their descendants remain consistently linked.
\vskip0.1in
\noindent This compatibility condition ensures that slope assignments vary in a controlled way across space, so that finer tubes with thinner roots remain clustered within their thicker predecessors; this occurs for a fixed sequence of scales identified by the fundamental heights of $\mathcal S_N$. It is precisely this coherence or ``stickiness'' that will allow us, in the next chapter, to construct families of tubes or rectangles whose orientations are governed by the compressed slope tree. In summary, the compressed slope tree $\mathscr{B}_N$, together with the root trees and the sticky maps linking them, provides a unified framework for assigning directions to spatial configurations in a manner adapted to the pruned structure of $\Omega_N$. The slope tree encodes available directions, while the root trees and sticky maps determine how these directions are distributed to roots across space. 
\vskip0.1in
\noindent The notion of tree-map pairs adapted to $\mathscr{B}_N$ appears in Section \ref{section: tree-map pair defn}. The construction of such root trees and sticky maps will be given in Section \ref{section: compressions of the unit interval}.

\section{Main results on compressed trees} 
The main results of this chapter are Propositions \ref{PROP: COMPRESSED SLOPES} and \ref{PROP: COMPRESSED ROOTS}. 
\subsection{Representation of the pruned slope set $\Omega_N$ as a full binary tree $\mathscr{B}_N$}  \label{section: binary tree recast} 
\begin{proposition} \label{PROP: COMPRESSED SLOPES} 
Let $\mathcal S_N$ be the model $(N, C_0)$-tree representing the pruned slope set $\Omega_N$ given by Proposition \ref{PRUNING STAGE 1}. Then there exists a full binary tree $\mathscr{B}_N$ of height $N$ that represents $\mathcal S_N$ in the sense that the vertices of $\mathscr{B}_N$ correspond to chains of basic slope intervals in $\mathcal S_N$.  
\vskip0.1in
\noindent More concretely, $\mathscr{B}_N$ satisfies the following properties:  
\vskip0.1in
\begin{enumerate}[(a)] 
\item \label{root-part-B_N} The {\em{root of the tree $\mathscr{B}_N$}}, i.e. the only vertex of the $0^{\text{th}}$ generation, is the unit interval $[0,1]$:
\begin{equation} \label{root of B_N} 
\mathcal V_0(\mathscr{B}_N) := \{ [0,1]\}, \; \text{ so that } \;  \# \bigl[ \mathcal V_0(\mathscr{B}_N) \bigr] = 1. 
\end{equation} 
%$v_0$, the unique element of ${\tt{BasicSl}}_0 = {\tt{SplitV}}_1$. It is an $M$-adic interval of length $M^{-\lamda_1}$.
\vskip0.1in 
\item \label{generation-part-B_N} For $1 \leq j \leq N$, 
\begin{equation} \mathcal V_j(\mathscr{B}_N) \cong {\tt{BasicSl}}_j.  \label{generations in B_N} \end{equation}  
The relation $\cong$ means that each {\em{vertex of the $j^{\text{th}}$ generation of $\mathscr{B}_N$}} is a nested chain of basic intervals of $\mathcal S_N$, of the form $\langle \theta_1, \ldots, \theta_j \rangle$, where 
\begin{equation} \label{chain of basic slopes} \theta_k \in {\tt{BasicSl}}_k(\mathcal S_N), \quad 1 \leq k \leq j \leq N, \quad \theta_1 \supsetneq \theta_{2} \supsetneq \cdots \supsetneq \theta_{j-1} \supsetneq \theta_j. \end{equation}  
The final element $\theta_j \in {\tt{BasicSl}}_j$ uniquely identifies the chain $\langle \theta_1, \ldots, \theta_j \rangle$. There are exactly $2^j$ vertices of $\mathscr{B}_N$ of the $j^{\text{th}}$ generation. 
\vskip0.1in
\item\label{ancestry-part-B_N}  {\em{Ancestry in the tree $\mathscr{B}_N$}} is determined by set containment: 
\begin{equation} \label{ancestry in B_N} 
\langle \theta_1, \ldots, \theta_j, \theta_{j+1} \rangle \text{ is a child of } \langle \theta_1, \ldots, \theta_j \rangle \text{ if } \theta_{j+1} \subsetneq \theta_j. 
\end{equation} 
Each vertex of $\mathscr{B}_N$ has exactly two children. 
\vskip0.1in 
\item {\em{Identification with $\mathcal S_N$ and $\Omega_N$:}} Every terminal vertex in $\mathscr{B}_N$, which is a member of ${\tt{BasicSl}}_N$, contains a unique $M$-adic interval of length $M^{-J}$ from $\mathcal S_N$. This, in turn, contains a unique element of $\Omega_N$. Thus the collection $\mathcal V_N(\mathscr{B}_N)$ is in one-to-one correspondence with $\mathcal V_J(\mathcal S_N)$, and therefore also with $\Omega_N$.   \label{B_N identification} 
\vskip0.1in 
\item {\em{Geometric interpretation of $\mathscr{B}_N$:}}  The $(N+1)$ vertices on a {\em{maximal ray}} in $\mathscr{B}_N$ 
\[ [0,1] \mapsto \langle \theta_1 \rangle \mapsto \langle \theta_1, \theta_2 \rangle \mapsto \cdots \mapsto \langle \theta_1, \ldots, \theta_N \rangle\] \label{geometry of B_N}
are elements of a nested sequence of $M$-adic intervals whose lengths correspond to a chain $\pmb{\lambda}(\cdot)$ of fundamental heights of $\mathcal S_N$, given by \eqref{chain of fundamental heights}:
\[ 1 \mapsto M^{-\lambda_1} \mapsto M^{-\lambda_2(\pmb{\varepsilon}_1)} \mapsto \cdots \mapsto M^{-\lambda_N(\pmb{\varepsilon}_{N-1})}. \] 
Here the binary sequence ${\pmb{\varepsilon}}_{N-1} \in \{0, 1\}^{N-1}$ is specified by the nesting property 
\begin{equation} \label{interspersed} 
v_0 \supsetneq \theta_1 \supsetneq v(\pmb{\varepsilon}_1) \supsetneq \theta_2 \supsetneq \cdots \supsetneq v(\pmb{\varepsilon}_{N-1}) \supsetneq \theta_N,\end{equation} 
where the splitting vertices $v(\pmb{\varepsilon}_j) \in {\tt{SplitV}}_j$ obey \eqref{split-nesting-1} and \eqref{split-nesting-2}. 
\vskip0.1in 
\item {\em{Scales of intervals in $\mathscr{B}_N$:}}  The two children of each non-terminal vertex of $\mathscr{B}_N$ have the same length, although this length may vary across vertices at the same generation. 
\vskip0.1in
\noindent More precisely, for $1 \leq j \leq N-1$, the two children of a single vertex $\langle \theta_1, \ldots, \theta_j \rangle$ in $\mathscr{B}_N$ represent $M$-adic intervals of length  $M^{-k}$  where
\begin{equation} 
k = \lambda_{j+1}(\pmb{\varepsilon}_{j}) = \lambda(v(\pmb{\varepsilon}_j)) \text{ is defined as in \eqref{defn: lambda_j}}. 
\end{equation}  
However, the quantity $\lambda_{j+1}(\pmb{\varepsilon}_{j})$ may vary as $\pmb{\varepsilon}_j$ ranges over $\{0, 1\}^j$, resulting in possibly varying scales over a single generation of intervals in $\mathscr{B}_N$.  \label{equidimensional children} 
\end{enumerate} 
\end{proposition} 
{\em{Remarks: }} 
\vskip0.1in 
\begin{enumerate}[1.] 
\item Proposition \ref{PROP: COMPRESSED SLOPES} is proved later in this chapter, in Section \ref{section: compressed slopes proof}.
\vskip0.1in 
\item It is important to note that the tree $\mathscr{B}_N$ is {\em{not}} an $M$-adic tree in general, since different vertices of a single generation $j$ in $\mathscr{B}_N$ may correspond to $M$-adic intervals of varying length. However, it does represent $\Omega_N$, in the sense of Proposition \ref{PROP: COMPRESSED SLOPES} part \eqref{B_N identification}. Furthermore, $\mathscr{B}_N$ has a cleaner structure than $\mathcal S_N$, which is convenient for certain purposes. 
\vskip0.1in
\item  In the pruned $M$-adic slope tree $\mathcal S_N$, the chain of basic slopes is interspersed with vertices in ${\tt{SplitV}}(\mathcal S_N)$, as in \eqref{interspersed}. However, these splitting vertices do not appear in $\mathscr{B}_N$. 
\end{enumerate} 
\subsection{Tree-map pairs adapted to $\mathscr{B}_N$} \label{section: tree-map pair defn}
The compressed tree $\mathscr{B}_N$ given by Proposition \ref{PROP: COMPRESSED SLOPES} provides a basis for building other trees and maps that respect its structure. Informally, a tree-map pair consists of a compressed root tree $\mathscr{U}$, whose vertices represent spatial intervals, together with a map $\sigma$ that assigns to each such interval a slope in $\mathscr{B}_N$ at the same level, in a way that preserves ancestry and scale.
\vskip0.1in
\noindent To make this precise, let us denote by $\mathcal Q^{\ast}$ the collection of all $M$-adic subintervals of $[0,1]$, and by $\mathcal Q(k)$ the collection of such intervals of length $M^{-k}$:
\begin{align} 
\mathcal Q(k) &:= \left\{ Q : h(Q) = k  \right\} = \left\{\left[ r, r+1\right] M^{-k}  : 0 \leq r < M^k \right\}, \label{what is Q(k)} \\  
\mathcal Q^{\ast} &:= \bigsqcup_{k=0}^{J_0} \mathcal Q(k) = \left\{ \left[ \frac{r}{M^k}, \frac{r+1}{M^k}\right] : 0 \leq r < M^k, \; 0 \leq k \leq J_0  \right\}, \label{what is Q-star} \\
J_0 &:= \max \bigl\{ \lambda(v) : v \in {\tt{SplitV}}_N (\mathcal S_N) \bigr\}. 
\end{align} 
\vskip0.1in
\noindent Given a full binary tree $\mathscr{B}_N$ representing $\Omega_N$ as in Proposition \ref{PROP: COMPRESSED SLOPES}, we say that a tree $\mathscr{U}$ and a map 
\begin{equation} \label{inductively sticky sigma}  
\sigma: \mathscr{U} \rightarrow \mathscr{B}_N
\end{equation}  is a {\em{tree-map pair $(\mathscr{U}, \sigma)$ adapted to $\mathscr{B}_N$}} if the following conditions hold:  
\vskip0.1in 
\begin{enumerate}[1.] 
\item (Height of $\mathscr{U}$) The tree $\mathscr{U}$, called a {\em{compressed root tree}}, is of height $N$. 
\vskip0.1in
\item (Vertices of $\mathscr{U}$) Each vertex of $\mathscr{U}$ is an $M$-adic interval, specifically an element of $\mathcal Q^{\ast}$. Ancestry in $\mathscr{U}$ is determined by set containment: if $Q, Q'$ are vertices of $\mathscr{U}$, then
\begin{equation}  Q \text{ is said to be a descendant of $Q'$ if } Q \subsetneq Q'. \end{equation} 
\vskip0.1in
\item (Scales at a fixed height) Vertex intervals of $\mathscr{U}$ at each level cover $[0,1]$, but unlike an $M$-adic tree, may have different lengths. However, the lengths are not entirely arbitrary, and must coincide with a fundamental height in $\mathcal S_N$. 
\vskip0.1in
\noindent Specifically, for $1 \leq j \leq N$, let  us denote 
\[ \mathcal V_j(\mathscr{U}) := \text{ the collection of $j^{\text{th}}$ generation vertices of $\mathscr{U}$}. \] Then $\mathcal V_j(\mathscr{U})$ forms a partition of $[0,1]$, with the $j^{\text{th}}$ fundamental heights of $\mathcal S_N$ as its contributing scales:
%\vskip0.1in
%\noindent To paraphrase, the $M$-adic intervals of $\mathcal V_j(\mathscr{U})$ have disjoint interiors, cover $[0,1]$, and unlike an $M$-adic tree, could have possibly varying lengths:
\begin{equation} \label{cover-Q}
\bigcup\bigl\{ Q : Q \in \mathcal V_j(\mathscr{U}) \bigr\} = [0,1]; \; \mathcal V_j(\mathscr{U}) \subseteq \bigsqcup \Bigl\{ \mathcal Q(\lambda(v)) : v \in {\tt{SplitV}_j} \Bigr\}.
\end{equation} 
This means the length of a vertex in $\mathcal V_j(\mathscr{U})$  must coincide with the length of some $j^{\text{th}}$ basic slope interval of $\mathcal S_N$, and is therefore restricted to lie in 
\[
 \bigl\{ M^{-\lambda(v)} : v \in {\tt{SplitV}}_j(\mathcal S_N )\bigr\} = \bigl\{ |w| : w \in {\tt{BasicSl}}_j(\mathcal S_N) \bigr\}.
\]
The covering property \eqref{cover-Q} translates to the relation 
\begin{equation} \label{cover-Q sum}
\sum_{Q} \bigl\{|Q|:  Q \in \mathcal V_j(\mathscr{U} \bigr\}= 1, \quad 1 \leq j \leq N.   
\end{equation}
\vskip0.1in 
\item (The map $\sigma$ preserves height, ancestry and length) The map $\sigma$ as in \eqref{inductively sticky sigma}, called the {\em{compressed root-to-slope}} map, is sticky in the sense defined in Section \ref{sticky maps section}. This means it preserves heights and lineages. Moreover, $\sigma$ preserves dimension, i.e., the image of any interval $Q$ in $\mathscr{U}$ under $\sigma$ has the same length as $Q$. 
%This is equivalent to saying that for $1 \leq j \leq N$, the map $\sigma_{\mathbb X}$ allocates to each $Q \in \mathcal V_j(\mathscr{U}_{\mathbb X})$ a $j^{\text{th}}$ basic slope interval of the same length, in other words a $j^{\text{th}}$ generation vertex of $\mathscr{B}_N$, in a way that preserves lineage:
\begin{align} 
&\mathcal V_j(\mathscr{U}) \ni Q \mapsto \sigma_{\mathbb X}(Q) \in {\tt{BasicSl}}_j(\mathcal S_N), \quad |Q| = |\sigma_{\mathbb X}(Q)|, \label{sigma preserves length} \\ 
&Q_1, Q_2 \in \mathscr{U}_{\mathbb X}, \; Q_1 \subseteq Q_2 \; \implies \; \sigma_{\mathbb X}(Q_1) \subseteq \sigma_{\mathbb X}(Q_2). \label{sigma preserves lineage} 
\end{align} 
\end{enumerate}  
\subsection{Extension of a tree-map pair to the finest scale} \label{section: finest scale extension} 
In light of the discussion in the preceding section, the length of intervals in ${\tt{BasicSl}}_N$ is, a priori, the smallest scale where a tree-map pair $(\mathscr{U}, \sigma)$ adapted to $\mathscr{B}_N$ is defined. A posteriori, however, a tree-map pair $(\mathscr{U}, \sigma)$ extends uniquely to a map on all of $\mathcal Q(J)$
\begin{equation} \label{penultimate sigma defn}
\bar{\sigma} : \mathcal Q(J) \rightarrow \Omega_N 
\end{equation}  
in the following way: for every $Q_0 \in \mathcal Q(J)$, we set 
\begin{equation} 
\bar{\sigma}(Q_0) := \sigma(Q) \cap \Omega_N, \label{final sigma defn}
\end{equation} 
where $Q$ is the terminal vertex of $\mathscr{U}$ containing $Q_0$.  Since $\sigma_{\mathbb X}(Q) \in {\tt{BasicSl}}_N$ contains exactly one element of $\Omega_N$ by the construction of the pruned tree $\mathcal S_N$, the map $\bar{\sigma}$ is well-defined. It preserves ancestry, though not length.  
\vskip0.1in
\noindent To reduce notational baggage and with a slight abuse of notation, we will rename $\bar{\sigma}$ on $\mathcal Q(J)$ as $\sigma$. 
\subsection{Adapted pairs: compressed root trees and inductively sticky maps} \label{section: compressions of the unit interval}
Having defined in the previous section a tree-map pair adapted to a pruned and compressed slope tree $\mathscr{B}_N$, we now show that such tree–map pairs exist in abundance, and can be constructed inductively from binary data.  Unlike classical sticky maps, where stickiness is enforced uniformly across all scales, the trees and maps constructed here exhibit {\em{inductive stickiness}}: the assignment of slopes is determined progressively across generations, past slope assignments are used to specify new generations of roots, and compatibility is enforced only at scales dictated by the compressed structure of $\mathscr{B}_N$. 
\vskip0.1in
\noindent Let us fix an ordered binary sequence of digits 0 and 1, indexed by the elements of $\mathcal Q^{\ast}$. The sequence $\mathbb X$ prescribes, for each vertex $Q$ that is a child of $Q'$, which of the two children of $\sigma_{\mathbb X}(Q')$ is selected as its image. Thus $\mathbb X$ encodes a global choice of branching across the tree.
\begin{equation} \label{binary X}
\mathbb X := \bigl\{ X(Q) : Q \in \mathcal Q^{\ast} \bigr\}, \quad X(Q) \in \{0, 1\}. 
\end{equation} 
For each $\mathbb X$,  we will inductively define a tree-map pair $(\mathscr{U}_{\mathbb X}, \sigma_{\mathbb X})$ in $N$ steps. The construction proceeds as follows: starting from the root, we define the first generation of $\mathscr{U}_{\mathbb X}$ using the first fundamental height. At each subsequent step, given a vertex $Q'$ and its assigned slope $\sigma_{\mathbb X}(Q')$. we determine its children and their slope assignments using the next fundamental height and the value of $X(Q)$. The admissible choices of scales and branching are governed by the fundamental heights and splitting structure of the pruned tree $\mathcal{S}_N$.
\vskip0.1in
\noindent The important features of the construction are summarized in the proposition below --  the second main result of this chapter.
\begin{proposition} \label{PROP: COMPRESSED ROOTS} 
For every binary sequence $\mathbb X$ as in \eqref{binary X}, there exists a tree-map pair $(\mathscr{U}_{\mathbb X}, \sigma_{\mathbb X})$ adapted to $\mathscr{B}_N$, in the sense of Section \ref{section: tree-map pair defn}, with the following properties: 
\vskip0.1in
\begin{enumerate}[(a)]
\item \label{compressed root} The root of $\mathscr{U}_{\mathbb X}$ is $[0,1]$, same as the root of $\mathscr{B}_N$. At the $0^{\text{th}}$ level, the map $\sigma_{\mathbb X}$ maps the former to the latter. 
\vskip0.1in 
\item \label{compressed roots - first generation}  The first generation of vertices $\mathcal V_1(\mathscr{U}_{\mathbb X})$ consists of all $M$-adic intervals in $\mathcal Q(\lambda_1)$, where $\lambda_1$ is the first fundamental height of $\mathcal S_N$. 
\vskip0.1in 
\item \label{compressed roots - children} At each step, the slope assigned to a child in $\mathscr{U}_{\mathbb X}$  is chosen as one of the two children of the parent slope, with the choice determined by the binary variable $X(Q)$. 
\vskip0.1in
\noindent Specifically, if $Q$ is a child of $Q'$ in $\mathscr{U}_{\mathbb X}$, then the slope allocation rule $\sigma_{\mathbb X}$ is defined as follows, 
\begin{equation}  |Q| = | \sigma_{\mathbb X}(Q)| \; \text{ and } \;  \sigma_{\mathbb X}(Q) := X(Q)^{\text{th}} \text{ child of } \sigma_{\mathbb X}(Q') \text{ in } \mathscr{B}_N. \label{defn: mapX} \end{equation}  
\vskip0.1in
\item \label{compressed roots - map} The children of any vertex If $Q' \in \mathscr{U}_{\mathbb X}$ are equi-dimensional, and constitute a uniform partition of $Q'$.  However, the children may have varying lengths as the parent $Q'$ ranges over all vertices of a fixed generation in $\mathscr{U}_{\mathbb X}$. 
\vskip0.1in
\noindent Specifically, suppose that 
\[ \sigma_{\mathbb X}(Q') = \theta, \] 
where $\theta$ is a basic slope interval descended from $v \in {\tt{SplitV}}(\mathcal S_N)$. Then 
\begin{align} 
&\# \left[ \text{children of } Q' \text{ in } \mathscr{U}_{\mathbb X}\right] = M^{\lambda(\bar{v}) - \lambda(v)}, \label{compressed children 1} \\ 
&\bigl\{|Q| : Q \text{ is a child of } Q' \text{ in } \mathscr{U}_{\mathbb X} \bigr\} = M^{-\lambda(\bar{v})}, \label{compressed children 2} 
\end{align} 
where $\bar{v}$ is the unique splitting vertex of $\mathcal S_N$ determined by \eqref{split-basic-split}, i.e. 
\[ \bar{v} \subsetneq \theta \subsetneq v, \quad \iota(\bar{v}) = \iota (v) + 1, \quad i=1, 2.\]    
\end{enumerate} 
\end{proposition} 
\vskip0.1in
\noindent {\em{Remarks: }} 
\vskip0.1in 
\begin{enumerate}[1.] 
\item Proposition \ref{PROP: COMPRESSED ROOTS} is proved in Section \ref{section: compressed roots proof} of this chapter.
\vskip0.1in
\item The root-to-slope map \eqref{final sigma defn} is induced by $\sigma_{\mathbb X}$ at the finest scale $\mathcal Q(J)$ will be the foundation of the Kakeya-type construction in the next chapter. 
\end{enumerate} 
\section{The compressed slope tree: Proof of Proposition \ref{PROP: COMPRESSED SLOPES}} \label{section: compressed slopes proof} 
\begin{proof} 
The compressed tree $\mathscr{B}_N$ is constructed by organizing chains of basic slope intervals into a tree structure indexed by splitting events, using the properties listed in parts \eqref{root-part-B_N}, \eqref{generation-part-B_N} and \eqref{ancestry-part-B_N}.  For instance, the root of $\mathscr{B}_N$ is defined by \eqref{root of B_N}, each vertex in $\mathcal V_j(\mathscr{B}_N)$ is defined as a chain of nested basic slope intervals given by \eqref{chain of basic slopes} and ancestry in $\mathscr{B}_N$ is described via \eqref{ancestry in B_N}. It then follows from the definition that $\mathscr{B}_N$ is a rooted labelled tree of height $N$. The properties listed in parts \eqref{root-part-B_N}, \eqref{generation-part-B_N} and \eqref{ancestry-part-B_N} are embedded in its construction. The proof then proceeds by verifying that this construction yields a {\em{full binary tree}} of height $N$, and that it satisfies properties \ref{B_N identification}, \eqref{geometry of B_N} and \eqref{equidimensional children} of Proposition \ref{PROP: COMPRESSED SLOPES}. 
\vskip0.1in
\noindent  Let us start by addressing the three latter parts. Part \ref{B_N identification} follows from the property of the pruned tree $\mathcal S_N$. Namely, Proposition \ref{PRUNING STAGE 1} and the definition of basic slopes imply that a ray originating at a vertex of ${\tt{BasicSl}}_N$ to height $J$ is non-splitting. In view of \eqref{one element per terminal vertex}, such a ray uniquely identifies a vertex of $\mathcal S_N$ at height $J$, which in turn contains a unique element of $\Omega_N$. This establishes the bijection between $\mathcal V_N(\mathscr{B}_N)$ and $\mathcal V_J(\mathcal S_N)$, and also a bijection between $\mathcal V_N(\mathscr{B}_N)$ and $\Omega_N$. 
\vskip0.1in
\noindent The discussion surrounding basic slope intervals and fundamental height chains in Sections \ref{section: basic slopes} and \ref{section: fundamental height chains}  ensures that $\mathscr{B}_N$ obeys the properties in parts \eqref{geometry of B_N} and \eqref{equidimensional children} of Proposition \ref{PROP: COMPRESSED SLOPES}. For instance, let us interpret the geometry of rays in $\mathscr{B}_N$, which establishes \eqref{geometry of B_N}. The property \eqref{split-basic-split} says that each basic slope interval of $\mathcal S_N$ is sandwiched between a youngest splitting ancestor and a uniquely defined oldest splitting descendant. This gives rise to the nested chain in \eqref{interspersed} in which splitting vertices and basic slope intervals alternate. The chain is identified through a binary sequence $\pmb{\varepsilon}_{N-1} = (\varepsilon_1, \ldots, \varepsilon_{N-1})$ determined by the location of the splitting vertices in $\mathcal S_N$. The $j^{\text{th}}$ element $\varepsilon_j$ of the sequence is 0 or 1 depending on whether $\theta_j$ is descended from the left or the right child of its last splitting ancestor $v(\pmb{\varepsilon}_{j-1})$. The length of a basic slope interval $\theta_j \in {\tt{BasicSl}}_j$ in this chain is 
\[ M^{- \lambda_j(\pmb{\varepsilon}_{j-1})},  \text{ as defined in \eqref{defn: lambda_j} and \eqref{basic slope - fundamental ht}; this proves \eqref{geometry of B_N}. } \] 
Next, we analyse how children are generated at each vertex. We have noted following the definition \eqref{jth basic slope cubes} that each basic slope in ${\tt{BasicSl}}_j$ has exactly two descendants in ${\tt{BasicSl}}_j$, and they are of the same length. This verifies \eqref{equidimensional children}. 
\vskip0.1in
\noindent In order to complete the proof of Proposition \ref{PROP: COMPRESSED SLOPES}, it remains to establish $\mathscr{B}_N$ as a full binary tree of height $N$. This may be clear to the discerning reader from the structural discussion of $\mathscr{B}_N$ thus far, but we include a proof for completeness. Let us recall the definition of the standard binary tree $\mathcal B_N$ of height $N$ whose vertex set at height $j$ consists of $j$-long binary sequences of $0$ and $1$. Every vertex of $\mathcal B_N$ is a splitting vertex and has exactly two children. In the remainder of the proof, we will identify a tree isomorphism 
\begin{equation} \Psi : \mathcal B_N \rightarrow \mathscr{B}_N \text{ in the sense of Section \ref{sticky maps section}}, \label{tree isomorphism Psi} \end{equation} 
which will establish $\mathscr{B}_N$ as a full binary tree of height $N$.
\vskip0.1in 
\noindent Let us recall from Proposition \ref{PRUNING STAGE 1} and Section \ref{section: splitting index} that  every splitting vertex of $\mathcal S_N$ has exactly two children, distinguished by the indices 0 and 1.  For $1 \leq j \leq N$, we define a bijective map 
\begin{equation}  \Psi_j: \{0,1\}^j \rightarrow {\tt{BasicSl}}_j \label{Psi bijection}  \end{equation}  inductively as follows. For $j=1$, 
\begin{equation} \label{defn Psi_1} \Psi_1(\varepsilon_1) := \left\{\begin{aligned}
&\text{ the unique element of ${\tt{BasicSl}}_1$} \\ &\text{ descended from $v(\varepsilon_1)$, the $\varepsilon_1^{\text{th}}$ child of $v_0$ }
\end{aligned} \right\},  \quad \varepsilon_1 =0,1,   \end{equation} 
where $v_0$ is the single element in ${\tt{BasicSl}}_0 = {\tt{SplitV}}_1$. The two intervals corresponding to $\Psi_1(\varepsilon_1)$ for $\varepsilon_1=0,1$ have the same length $M^{-\lambda_1}$, with $\lambda_1$ as in \eqref{defn: lambda_1}. Clearly, $\Psi_1$ is a bijection from one two-point set to another. 
\vskip0.1in
\noindent In general, suppose that $\Psi_{j}$ has been defined as a bijection of $\{0, 1\}^j$ onto ${\tt{BasicSl}}_j$. Then for $\pmb{\varepsilon}_j \in \{0,1\}^{j}$ and $\varepsilon =0,1$, we set 
\begin{equation} 
\Psi_{j+1}(\pmb{\varepsilon}_j, \varepsilon) := \left\{\begin{aligned}
&\text{ the unique element of ${\tt{BasicSl}}_{j+1}$ } \\ &\text{ descended from the $\varepsilon^{\text{th}}$ child of $v(\pmb{\varepsilon}_j)$}  \end{aligned} \right\}, \label{defn Psi_j} \end{equation} 
where $v(\pmb{\varepsilon}_j)$ is the splitting vertex of index $(j+1)$ defined as in \eqref{split-nesting-1} and \eqref{split-nesting-2}. 
%$\gamma_{j+1}$ is the unique element of ${\tt{SplitV}}_{j+1}$ identified by $\Psi_{j}(\bar{\epsilon}) \in {\tt{BasicSl}}_j$, as described in Section \ref{section:  fundamental heights and basic slopes}. 
Since $\Psi_j$ is a bijection, and $v(\pmb{\varepsilon}_j)$ has exactly two children, we conclude that $\Psi_{j+1}$ is a bijection as well. Indeed, $\Psi_{j+1}$ preserves ancestry, in the sense that 
\begin{equation} \label{Psi preserves ancestry}   
\Psi_{j+1}(\pmb{\varepsilon}_j, \varepsilon)  \subsetneq \Psi_{j}(\pmb{\varepsilon}_j); \;    \Psi_{j+1}(\pmb{\varepsilon}_j, \varepsilon) \text{ has length } M^{-\lambda_{j+1}(\pmb{\varepsilon}_j)} \text{ for } \varepsilon = 0, 1. 
\end{equation}   
Proceeding in this way, we arrive at a sequence of maps $(\Psi_1, \ldots, \Psi_N)$, where each $\Psi_j$ is a bijection as in \eqref{Psi bijection}, with a nesting property as in \eqref{Psi preserves ancestry}. The final map in the sequence, namely $\Psi_N$, provides a bijection of $\{0,1\}^N$ onto ${\tt{BasicSl}}_N$.  For $\Psi_j$ as in \eqref{defn Psi_j}, the tree isomorphism $\Psi$ claimed in \eqref{tree isomorphism Psi} is therefore given by  
\begin{equation} \label{splitting vertex binary tree isomorphism}
\begin{aligned} 
\Psi(\emptyset) &= \text{ $v_0 \in {\tt{BasicSl}}_0$}, \\ \Psi(\bar{\pmb{\varepsilon}}) &= \; \Psi_{j}(\bar{\pmb{\varepsilon}}) \in {\tt{BasicSl}}_j \; \text{ if } \; \bar{\pmb{\varepsilon}} \in \{0,1\}^j, \; 1 \leq j \leq N.
\end{aligned}
\end{equation}  
This completes the proof of Proposition \ref{PROP: COMPRESSED SLOPES}.
%This, in turn, is in one-to-one correspondence with the pruned set of slopes $\Omega_N$, since each vertex in ${\tt{BasicSl}}_N$ contains exactly one element of $\Omega_N$. 
\end{proof} 

\section{Compressed root trees: Proof of Proposition \ref{PROP: COMPRESSED ROOTS}} \label{section: compressed roots proof} 
\begin{proof}
We construct $\mathscr{U}_{\mathbb X}$ and $\sigma_{\mathbb X}$ inductively on generations, and verify that the defining criteria of an adapted tree-map pair and the properties \eqref{compressed root}-\eqref{compressed roots - map} of Proposition \ref{PROP: COMPRESSED ROOTS} are satisfied up to each generation. It will be helpful to remember the big picture: $X(Q)$ determines the slope assigned to $Q \in \mathscr{U}_{\mathbb X}$, and this slope $\sigma_{\mathbb X}(Q) \in \mathscr{B}_N$ in turn prescribes the subsequent branching of $Q$ in the next generation of $\mathscr{U}_{\mathbb X}$. A schematic depiction of the process is given in Figure \ref{fig:adapted-tree-map}.

\begin{figure}[ht]
\centering

\begin{tikzpicture}[
    x=1cm,
    y=1cm,
    font=\small,
    blueedge/.style={
        blue!75!black,
        line width=1.25pt,
        line cap=round
    },
    rededge/.style={
        red!75!black,
        line width=1.25pt,
        line cap=round
    },
    guide/.style={
        gray!50,
        densely dashed,
        line width=0.45pt
    },
    vertex/.style={
        circle,
        fill=black,
        inner sep=1.7pt
    },
    maparrow/.style={
        ->,
        line width=0.9pt
    },
    verticalarrow/.style={
        <->,
        line width=0.8pt
    }
]

% =================================================
% Adjustable coordinates
% =================================================

% Upper diagram levels
\def\toprootheight{8.10}
\def\topfirstheight{6.55}
\def\topsecondheight{5.20}
\def\toplowestheight{3.95}

% Lower diagram levels
\def\bottomrootheight{2.55}
\def\bottomfirstheight{1.55}
\def\bottomleafheight{0.50}

% Main horizontal centres
\def\leftcenter{2.55}
\def\rightcenter{8.75}

% =================================================
% GUIDE LINES
% Drawn first so all trees and labels appear above them
% =================================================

% Upper diagram
\draw[guide]
    (2.80,\toprootheight) -- (8.55,\toprootheight);

\draw[guide]
    (1.40,\topfirstheight) -- (10.00,\topfirstheight);

\draw[guide]
    (0.75,\topsecondheight) -- (8.70,\topsecondheight);

\draw[guide]
    (1.95,\toplowestheight) -- (10.70,\toplowestheight);

% Lower diagram
\draw[guide]
    (2.95,\bottomrootheight) -- (8.60,\bottomrootheight);

\draw[guide]
    (1.70,\bottomfirstheight) -- (9.90,\bottomfirstheight);

\draw[guide]
    (1.35,\bottomleafheight) -- (10.25,\bottomleafheight);

% =================================================
% UPPER LEFT TREE: J_3([0,1];2)
% =================================================

\node[font=\large] at (2.35,9.05)
    {$\mathcal{T}_3([0,1];2)$};
\node[font=\large] at (9.35,9.05)
    {$\mathcal S_2$};

% Root and first generation
\node[vertex] (ULroot) at (\leftcenter,\toprootheight) {};
\node[vertex] (ULleft)  at (1.70,\topfirstheight) {};
\node[vertex] (ULright) at (3.35,\topfirstheight) {};

% Root branches: red left, blue right
\draw[rededge]  (ULroot) -- (ULleft);
\draw[blueedge] (ULroot) -- (ULright);

% Left subtree: two blue branches
\node[vertex] (ULL1) at (1.05,\topsecondheight) {};
\node[vertex] (ULL2) at (2.05,\topsecondheight) {};

\draw[blueedge] (ULleft) -- (ULL1);
\draw[blueedge] (ULleft) -- (ULL2);

% Right subtree: four branches, red-blue-red-blue
\node[vertex] (ULR1) at (2.25,\toplowestheight) {};
\node[vertex] (ULR2) at (3.00,\toplowestheight) {};
\node[vertex] (ULR3) at (3.72,\toplowestheight) {};
\node[vertex] (ULR4) at (4.48,\toplowestheight) {};

\draw[rededge]  (ULright) -- (ULR1);
\draw[blueedge] (ULright) -- (ULR2);
\draw[rededge]  (ULright) -- (ULR3);
\draw[blueedge] (ULright) -- (ULR4);

% Interval labels
\node[
    anchor=west,
    xshift=5pt,
    yshift=2pt,
    fill=white,
    inner sep=1pt
] at (ULroot)
    {$[0,1]$};

\node[
    anchor=east,
    xshift=-5pt,
    yshift=3pt,
    fill=white,
    inner sep=1pt
] at (ULleft)
    {$[0,\frac12]$};

\node[
    anchor=west,
    xshift=5pt,
    yshift=3pt,
    fill=white,
    inner sep=1pt
] at (ULright)
    {$[\frac12,1]$};

\node[
    anchor=east,
    xshift=-4pt,
    yshift=-1pt,
    fill=white,
    inner sep=1pt
] at (ULL1)
    {$[0,\frac14]$};

\node[
    anchor=west,
    xshift=4pt,
    yshift=-1pt,
    fill=white,
    inner sep=1pt
] at (ULL2)
    {$[\frac14,\frac12]$};

% =================================================
% UPPER RIGHT TREE
% =================================================

\node[vertex] (URroot) at (\rightcenter,\toprootheight) {};

\node[
    anchor=west,
    xshift=5pt,
    yshift=2pt,
    fill=white,
    inner sep=1pt
] at (URroot)
    {$[0,1]$};

% First binary generation
\node[vertex] (UR0) at (7.90,\topfirstheight) {};
\node[vertex] (UR1) at (9.60,\topfirstheight) {};

% Root colours: blue for 0, red for 1
\draw[blueedge] (URroot) -- (UR0);
\draw[rededge]  (URroot) -- (UR1);

\node[
    left=4pt,
    fill=white,
    inner sep=1pt
] at ($(URroot)!0.52!(UR0)$)
    {$0$};

\node[
    right=4pt,
    fill=white,
    inner sep=1pt
] at ($(URroot)!0.52!(UR1)$)
    {$1$};

% Descendants of 0 lie on the upper terminal level
\node[vertex] (UR00) at (7.20,\topsecondheight) {};
\node[vertex] (UR01) at (8.35,\topsecondheight) {};

\draw[blueedge] (UR0) -- (UR00);
\draw[rededge]  (UR0) -- (UR01);

\node[
    left=3pt,
    fill=white,
    inner sep=1pt
] at ($(UR0)!0.54!(UR00)$)
    {$0$};

\node[
    right=3pt,
    fill=white,
    inner sep=1pt
] at ($(UR0)!0.54!(UR01)$)
    {$1$};

% Descendants of 1 lie on the lower terminal level
\node[vertex] (UR10) at (8.75,\toplowestheight) {};
\node[vertex] (UR11) at (10.35,\toplowestheight) {};

\draw[blueedge] (UR1) -- (UR10);
\draw[rededge]  (UR1) -- (UR11);

\node[
    left=4pt,
    fill=white,
    inner sep=1pt
] at ($(UR1)!0.58!(UR10)$)
    {$0$};

\node[
    right=4pt,
    fill=white,
    inner sep=1pt
] at ($(UR1)!0.58!(UR11)$)
    {$1$};

% =================================================
% VERTICAL RELATION BETWEEN UPPER AND LOWER PAIRS
% =================================================

\draw[verticalarrow]
    (6.05,3.52) -- (6.05,2.95);

% =================================================
% LOWER LEFT TREE: U_x
% =================================================

\node[font=\large] at (2.80,3.15)
    {$\mathscr{U}_{\mathbb X}$};

\node[vertex] (BLroot) at (2.75,\bottomrootheight) {};
\node[vertex] (BLleft)  at (2.00,\bottomfirstheight) {};
\node[vertex] (BLright) at (3.50,\bottomfirstheight) {};

% Root branches: red left, blue right
\draw[rededge]  (BLroot) -- (BLleft);
\draw[blueedge] (BLroot) -- (BLright);

% Left child: two blue branches
\node[vertex] (BLL1) at (1.58,\bottomleafheight) {};
\node[vertex] (BLL2) at (2.32,\bottomleafheight) {};

\draw[blueedge] (BLleft) -- (BLL1);
\draw[blueedge] (BLleft) -- (BLL2);

% Right child: four branches, red-blue-red-blue
\node[vertex] (BLR1) at (2.88,\bottomleafheight) {};
\node[vertex] (BLR2) at (3.30,\bottomleafheight) {};
\node[vertex] (BLR3) at (3.72,\bottomleafheight) {};
\node[vertex] (BLR4) at (4.25,\bottomleafheight) {};

\draw[rededge]  (BLright) -- (BLR1);
\draw[blueedge] (BLright) -- (BLR2);
\draw[rededge]  (BLright) -- (BLR3);
\draw[blueedge] (BLright) -- (BLR4);

% =================================================
% LOWER RIGHT TREE: B_2
% =================================================

\node[font=\large] at (8.80,3.15)
    {$\mathscr{B}_2$};

\node[vertex] (BRroot) at (8.80,\bottomrootheight) {};

\node[vertex] (BR0) at (8.05,\bottomfirstheight) {};
\node[vertex] (BR1) at (9.55,\bottomfirstheight) {};

% Root branches: blue left, red right
\draw[blueedge] (BRroot) -- (BR0);
\draw[rededge]  (BRroot) -- (BR1);

% Each binary child has blue-left and red-right branches
\node[vertex] (BR00) at (7.62,\bottomleafheight) {};
\node[vertex] (BR01) at (8.43,\bottomleafheight) {};
\node[vertex] (BR10) at (9.15,\bottomleafheight) {};
\node[vertex] (BR11) at (9.98,\bottomleafheight) {};

\draw[blueedge] (BR0) -- (BR00);
\draw[rededge]  (BR0) -- (BR01);

\draw[blueedge] (BR1) -- (BR10);
\draw[rededge]  (BR1) -- (BR11);

% =================================================
% LOWER MAP ARROW
% =================================================

\node[font=\large] at (6.00,1.72)
    {$\sigma_{\mathbb X}$};

\draw[maparrow]
    (5.42,1.30) -- (6.68,1.30);

\end{tikzpicture}

\caption{\small{
A schematic adapted tree-map pair $(\mathscr{U}_{\mathbb X}, \sigma_{\mathbb X})$ associated with
$\mathcal{T}_3([0,1];2)$. The top two trees are $M$-adic sub-trees; the left tree corresponds to roots in $[0,1]$; the right one represents the pruned slope tree $\mathcal S_N$, $N=2$. The two trees at the lower level are their compressed versions. The lower tree $\mathscr{U}_x$ retains the branching and colour pattern inherited from the upper-left tree, while $\sigma_{\mathbb X}$ maps it to the compressed binary tree $\mathcal B_2$ in a manner that preserves length and ancestry. Edges on the root side map to edges of the same colour on the slope side. Thus $\mathbb X$ in this figure obeys $X([0, \frac{1}{2}]) = 1$, $X([\frac12, 1]) = 0$, $X([0, \frac14]) = X([\frac14, \frac 12]) = 0$, etc. The dashed horizontal lines indicate successive generations in both sets of trees.
}}
\label{fig:adapted-tree-map}
\end{figure}
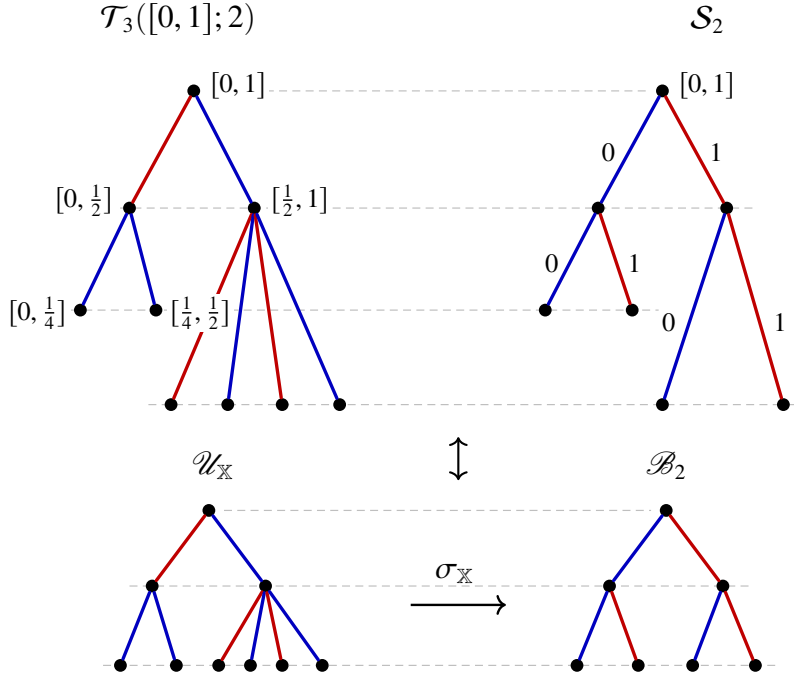

\subsection{The root of $\mathscr{U}_{\mathbb X}$} The root vertex of $\mathscr{U}_{\mathbb X}$, which is the unique vertex of the $0^{\text{th}}$ generation, corresponds to the interval $[0,1]$. This coincides with the root of the pruned and compressed slope tree $\mathscr{B}_N$; see \eqref{root of B_N}. We define 
\[ \sigma_{\mathbb X}([0,1]) := [0,1]; \quad \text{ thus \eqref{compressed root} is satisfied.}   \]
\subsection{The first generation vertices of $\mathscr{U}_{\mathbb X}$ and the first level sticky map} The first generation of vertices in $\mathscr{U} = \mathscr{U}_{\mathbb X}$, denoted $\mathcal V_1(\mathscr{U})$, is defined as 
	\begin{equation} \label{V1U} \mathcal{V}_1(\mathscr{U}) := \mathcal Q(\lambda_1), \; \text{ so that } \; \#\bigl(\mathcal V_1(\mathscr{U})\bigr) = M^{\lambda_1}. \end{equation} 
	Here $\mathcal Q(\cdot)$ is the collection of $M$-adic intervals given in \eqref{what is Q(k)}. In other words, the first generation of $\mathscr{U}$ consists of all $M$-adic intervals in $[0,1]$ of length $M^{-\lambda_1}$, with $\lambda_1$ being the first fundamental height of $\mathcal S_N$, given by \eqref{defn: lambda_1}. This verifies the requirement \eqref{compressed roots - first generation}.  
	\vskip0.1in
	\noindent We define  a map
\begin{align*} 
&\sigma_{\mathbb X}: \mathcal V_1(\mathscr{U}) = \mathcal Q(\lambda_1) \rightarrow {\tt{BasicSl}_1}, \text{ by the rule } \\  
&\sigma_{\mathbb X}(Q) := \Psi_1(X(Q)) \text{ for $Q \in \mathcal V_1(\mathscr{U})$},   
\end{align*}
with $\Psi_1$ defined as in \eqref{defn Psi_1}. Said differently, $\sigma_{\mathbb X}$ maps $Q$ onto a first-order basic slope interval. There are two such intervals, both of length $M^{-\lambda_1}$, and $\sigma_{\mathbb X}$ picks the one descended from the left or right child of $v_0$, according as $X(Q) =0$ or 1. Since \eqref{generations in B_N} identifies ${\tt{BasicSl}}_1$ with the first generation of $\mathscr{B}_N$, we conclude that $\sigma_{\mathbb X}$ preserves length, height and lineage between $\mathscr{U}_{\mathbb X}$ and $\mathscr{B}_N$ up to level 1.  Thus property \eqref{compressed roots - children} is confirmed. 
\vskip0.1in
\noindent Finally, all the children of the root $Q' = [0,1]$ are equi-dimensional, and they form a uniform cover of $[0,1]$, by \eqref{V1U}. Therefore \eqref{compressed roots - map} holds. Thus all the stated conditions for an adapted tree-map pair have been verified for $(\mathscr{U}_{\mathbb X}, \sigma_{\mathbb X})$ up until the first generation. 
\subsection{The second generation of roots and maps} Let us continue to the construction of the second generation of the tree-map pair. The class $\mathcal V_2(\mathscr{U}_{\mathbb X})$ of second generation vertices of $\mathscr{U}_{\mathbb X}$ is defined as the collection of intervals $Q_2$ satisfying 
\begin{align} 
&Q_2 \subsetneq Q_1 \in \mathcal V_1(\mathscr{U}), \quad Q_2 \in \mathcal Q(\lambda_{2}(\pmb{\varepsilon}_1)) \text{ if } X(Q_1) = \varepsilon_1 \in \{0, 1\};  \label{Q1Q2} \\  
& \text{ this means } \# \bigl\{ Q_2 \in \mathcal V_2(\mathscr{U}_{\mathbb X}) : Q_2 \subseteq Q_1\bigr\} = M^{\lambda_2(\pmb{\varepsilon}_1) - \lambda_1} \text{ if } X(Q_1) = \varepsilon_1. \nonumber 
\end{align}  
To be more explicit, the second generation of $\mathscr{U}_{\mathbb X}$ is dictated by the mapping behaviour of $\sigma_{\mathbb X}$ in the first generation. An interval $Q_1 \in \mathcal V_1(\mathscr{U}_{\mathbb X})$ that is assigned by $\sigma_{\mathbb X}$ to the left (respectively right) basic slope of ${\tt{BasicSl}}_1$ splits into sub-intervals of length $M^{-\lambda_2(0)}$ (respectively $M^{-\lambda_2(1)}$) at the second level. This leads to the following conclusions: 
\vskip0.1in  
\begin{itemize} 
\item The number of second generation children of a first generation parent $Q_1$ is determined by the slope $\sigma_{\mathbb X}(Q_1)$ that the parent receives.
\vskip0.1in 
\item The second generation children originating from the same parent are of the same length. Children originating from different parents need not be of the same length.
\vskip0.1in  
\item  A maximum of two distinct interval lengths $M^{-\lambda_2(\varepsilon_1)}$, $\varepsilon_1 = 0, 1$, may be represented among the vertices of the second generation. 
\vskip0.1in
\item The second generation vertices form a (possibly non-uniform) partition of $[0,1]$, with disjoint interiors.  
\end{itemize} 
These observations verify the covering and dimension properties \eqref{cover-Q}, \eqref{cover-Q sum} of a tree-map pair up to height 2. They also verify part \eqref{compressed roots - map} of Proposition \ref{PROP: COMPRESSED ROOTS}. Properties \eqref{compressed root} and \eqref{compressed roots - first generation} of course carry over from the previous steps. 
\vskip0.1in
\noindent We now define $\sigma_{\mathbb X}(Q_2)$ as the second basic interval descended from the $X(Q_2)^{\text{th}}$ child of $\sigma_{\mathbb X}(Q_1)$:
\[ \sigma_{\mathbb X}(Q_2) := \Psi_2 \bigl(X(Q_1), X(Q_2) \bigr) \in {\tt{BasicSl}}_2, \text{ with $Q_1$ as in \eqref{Q1Q2}. } \] 
This means that $\sigma_{\mathbb X}(Q_2)$ is the second-order basic interval descended from the $X(Q_2)^{\text{th}}$ child of $\sigma_{\mathbb X}(Q_1) = \Psi_1(X(Q_1))$. In view of \eqref{Psi preserves ancestry} and \eqref{Q1Q2}, we deduce 
\[ |\sigma_{\mathbb X}(Q_2)| = |Q_2| = M^{-\lambda_2(\pmb{\varepsilon}_1)} \text{ if } \varepsilon_1 = X(Q_1). \] 
Thus $\sigma_{\mathbb X}$ preserves length, height and lineage among $\mathscr{U}_{\mathbb X}$ and $\mathscr{B}_N$ up to height 2. This confirms \eqref{defn: mapX}, \eqref{compressed children 1} and \eqref{compressed children 2}, verifying Proposition \ref{PROP: COMPRESSED ROOTS} \eqref{compressed roots - map}. 
\subsection{The inductive step} Proceeding inductively, suppose that $\mathcal V_k(\mathscr{U}_{\mathbb X})$ and 
\[ \sigma_{\mathbb X}: \mathcal V_k(\mathscr{U}_{\mathbb X}) \rightarrow {\tt{BasicSl}}_k  \] 
have been defined for $1 \leq k \leq j$, obeying the requirements of parts \eqref{compressed roots - children} and \eqref{compressed roots - map} of Proposition \ref{PROP: COMPRESSED ROOTS} at all levels up to at including $j$. Precisely, this means that 
\begin{align} 
&|\sigma_{\mathbb X}(Q_k)| = |Q_k| \quad \text{ for each }  Q_k \in \mathcal V_k(\mathscr{U}_k), \; 1\leq k \leq j; \text{ further, }\nonumber \\ 
&\text{if } Q_j \subsetneq Q_{j-1} \subsetneq \cdots \subsetneq Q_1,  \;\text{ then } \; \sigma_{\mathbb X}(Q_j) \subsetneq \sigma_{\mathbb X}(Q_{j-1}) \subsetneq \cdots \subsetneq \sigma_{\mathbb X}(Q_1). \label{nesting at j-1}
\end{align} 
\vskip0.1in
\noindent Based on this,  the $(j+1)^{\text{th}}$ generation  of vertices $\mathcal V_{j+1}(\mathscr{U}_{\mathbb X})$ of $\mathscr{U}_{\mathbb X}$ is defined in the following way.  Suppose that \eqref{nesting at j-1} holds for a nested sequence $(Q_1, \ldots, Q_{j})$, with 
\begin{align} 
&\pmb{\varepsilon}_{j} := (X(Q_1), \ldots, X(Q_{j})), \; \text{ so that } \label{eps-X} \\ 
&\sigma_{\mathbb X}(Q_{j}) = \Psi_{j}(\pmb{\varepsilon}_{j}) \; \text{ and } \; |Q_{j}| = M^{- \lambda_{j}(\pmb{\varepsilon}_{j-1})}. \label{Q_j length} 
\end{align}  
Then we set
\begin{equation} \label{vertex set - induction}  
\mathcal V_{j+1}(\mathscr{U}_{\mathbb X}) := \bigcup_{Q_j \in \mathcal V_j(\mathscr{U}_{\mathbb X})} \left\{ Q \in \mathcal Q^{\ast} \; \Biggl| \; \begin{aligned} 
&Q \subsetneq Q_j, \; Q \in \mathcal Q(\lambda_{j+1}(\pmb{\varepsilon}_{j})), \\ &\pmb{\varepsilon}_j \text{ as in } \eqref{nesting at j-1} \text{ and } \eqref{eps-X} \end{aligned} \right\}. 
\end{equation}  
In other words, the children of $Q_{j}$ in $\mathscr{U}_{\mathbb X}$ are all of length $M^{-\lambda_{j+1}(\pmb{\varepsilon}_{j})}$, and together they cover $Q_{j}$. Here $\lambda_{j+1}(\cdot)$ is the $(j+1)^{\text{th}}$ fundamental height defined in \eqref{defn: lambda_j}. Once $\mathcal V_{j+1}(\mathscr{U}_{\mathbb X})$ is specified, the map 
\[\sigma_{\mathbb X} : \mathcal V_{j+1}(\mathscr{U}_{\mathbb X}) \rightarrow {\tt{BasicSl}}_{j+1} \] is defined to be
\begin{equation} \label{image slope}   
\begin{aligned}
\sigma_{\mathbb X}(Q_{j+1}) &:= \Psi_{j+1}(X(Q_1), \ldots, X(Q_{j}), X(Q_{j+1}))  \\
&\,= \Psi_{j+1}(\pmb{\varepsilon}_{j+1}) = \varepsilon_{j+1}^{\text{th}} \text{ child of } \sigma_{\mathbb X}(Q_{j}),  
\end{aligned} 
\end{equation} 
with $\Psi_{j+1}$ as in \eqref{defn Psi_j}, $Q_{j+1}$ a child of $Q_j$,  $\pmb{\varepsilon}_{j+1} = (\pmb{\varepsilon}_j, \varepsilon_{j+1})$ and $\varepsilon_{j+1} = X(Q_{j+1})$. 
\vskip0.1in
\noindent We now verify that the above construction of $(\mathscr{U}_{\mathbb X}, \sigma_{\mathbb X})$ satisfies \eqref{compressed root}--\eqref{compressed roots - map} up until height $(j+1)$. Properties \eqref{compressed root} and \eqref{compressed roots - first generation} are inherited from the first two steps of the construction. Part \eqref{compressed roots - children} is a consequence of \eqref{image slope}, since 
\[ |Q_{j+1}| = |\sigma_{\mathbb X}(Q_{j+1})| = M^{-\lambda_{j+1}(\pmb{\varepsilon}_j)}, \text{ by \eqref{Psi preserves ancestry}}.  \] 
Equi-dimensionality and number of children of a fixed vertex, as stated in \eqref{compressed roots - map} follows from \eqref{Q_j length} and the definition \eqref{vertex set - induction}, since 
\[ \# \left\{ Q_{j+1} \in \mathcal V_{j+1}(\mathscr{U}_{\mathbb X}) : Q_{j+1} \subsetneq Q_j \right\} = M^{\lambda_{j+1}(\pmb{\varepsilon}_j) - \lambda_j(\pmb{\varepsilon}_{j-1})}. \] 
This closes the induction, completing the proof of Proposition \ref{PROP: COMPRESSED ROOTS}.  
\end{proof}

\chapter{Kakeya-type sets: sticky tube families with pruned slopes} \label{chapter: sticky rectangles chapter} 
Let us consider the present juncture of the argument, based on the findings of the preceding three chapters.  
\vskip0.1in
\begin{itemize} 
\item In Chapter \ref{Chapter: Pruning of the slope tree}, we have identified a {\em{model subset}} $\Omega_N$ of the slope set $\Omega$ whose finitary $M$-adic tree $\mathcal S_N$ can be viewed as a possibly elongated and asymmetric version of a full binary tree of splitting number $N$. 
\vskip0.1in
\item In Chapter \ref{Chapter: Fundamental heights}, we have isolated further special features of $\mathcal S_N$, such as {\em{fundamental heights}} and {\em{basic slope intervals}}, which helped compress the slope tree as a full binary tree $\mathscr{B}_N$, with possibly variable scales at each height.  
\vskip0.1in
\item In Chapter \ref{chapter: compressions}, the compressed slope tree $\mathscr{B}_N$ turned out to be instrumental in determining an inductively defined family of compressed root trees $\mathscr{U}_{\mathbb X}$ together with sticky maps
\[ \sigma_{\mathbb X} : \mathscr{U}_{\mathbb X} \longrightarrow \mathscr{B}_N. \]  Both $\mathscr{U}_{\mathbb X}$ and $\sigma_{\mathbb X}$ are indexed by binary sequences $\mathbb X$. As described in Section \ref{section: finest scale extension}, each such $\sigma_{\mathbb X}$ permits a natural extension to the finest scale $M^{-J}$, yielding an extended map
\[\sigma_{\mathbb X} : \mathcal Q(J) \longrightarrow \Omega_N.  \] 
%Each of these trees encodes $\mathcal Q(J)$, the $M$-adic intervals of the finest scale $M^{-J}$ given by \eqref{what is Q(k)}, \eqref{what is Q-star}.
%The binary sequence $\mathbb X$ also helps in defining a sticky map $\sigma_{\mathbb X}: \mathscr{U}_{\mathbb X} \rightarrow \mathscr{B}_N$, which assigns every $Q \in \mathcal Q(J)$ a slope from $\Omega_N$.    
\end{itemize}
\vskip0.1in
 \noindent Equipped with these, we are finally in a position to describe families of rectangles that will ultimately lead to the desired Kakeya-like configuration.  The sets of interest that will eventually verify condition \ref{condition 1: Kakeya-type sets} of Theorem \ref{thm:main} are unions of families of thin rectangles, whose slopes and locations in the plane are determined by $\sigma_{\mathbb X}$ in a way that preserves ``geometric stickiness''.

\section{Families of sticky tubes} \label{general facts about tube families}
Let us fix two absolute constants $c_0, A_0 > 0$, with  
\begin{equation} \label{c0A0} 
0 < c_0 < 1/2 < 2 < A_0.
\end{equation} 
For instance, choosing $A_0 = c_0^{-1} = 10^3$ will suffice.  
\begin{definition}[A tube rooted at an interval] 
Given an $M$-adic interval $Q$ contained in $[0,1]$ and a slope $\omega \in [0,1]$, we define a tube rooted at $Q$ with orientation $\omega$ as a  parallelogram of the form 
 \begin{equation} \label{defn: tube}
{\tt{Tube}} =  {\tt{Tube}}(Q, \omega, I) := \bigl\{ (0, y) + \rho (1, \omega) : y \in \tilde{Q}, \; \rho \in I \bigr\} 
  \end{equation} 
Here $\tilde{Q}$ denotes the $c_0$-dilate of $Q$, i.e., the interval with the same centre as $Q$ but $c_0$ times its length; $I$ is an interval contained in $[0, 10 A_0]$ of length larger than $|Q|$.   
\end{definition} 
\vskip0.1in
\noindent A few observations emerge from the definition: 
\vskip0.1in
\begin{itemize} 
\item Each of the three entries in the argument of ${\tt{Tube}}(Q, \omega, I)$ records identifying geometric information about it; $\tilde{Q}$ determines both its short side and its position in the plane, $\omega$ specifies the orientation of its long side and $I$ records its length. 
\vskip0.1in
\item The interval $\{0\} \times Q$ is called the {\em{root interval}} associated with ${\tt{Tube}}$, even though the actual intersection of ${\tt{Tube}}$ with the $y$-axis is $\{0\} \times \tilde{Q}$. By a slight abuse of notation, we will often refer to $Q$ as the root of ${\tt{Tube}}$.  
\vskip0.1in 
\item The segment $\{0\} \times [0,1]$, which is the union of all possible root intervals, will be referred to as the {\em{root line}}. 
\vskip0.1in
\item Each {\tt{Tube}} is contained in a larger parallelogram of comparable dimensions
\[ {\tt{Tube}}(Q, \omega, I) \subseteq  {\tt{Tube}}\bigl(Q, \omega, [0, 10A_0] \bigr) \] 
whose long side is a large multiple of $A_0$ oriented in the direction $(1, \omega)$, and whose short side is positioned on the $y$-axis along the interval $\tilde{Q}$. 
\vskip0.1in
\item All tubes are located within a large square of the positive quadrant with a corner at the origin, e.g. $A_0^2 [0,1]^2$.    
\end{itemize} 
%However, the assumption $\omega \in [0,1]$ ensures that {\tt{Tube}} contains and is contained in actual rectangles of dimensions comparable to ${\tt{Tube}}$, oriented in the same direction, with absolute implicit constants. With this minor adjustment, we will henceforth work with tubes rather than rectangles.  

\begin{definition}
Given a binary sequence $\mathbb X$ as in \eqref{binary X}, we define a {\em{sticky tube}} rooted at $t$ associated with $\mathbb X$ as follows:  
\begin{equation} \label{defn: sticky tube} 
{\tt{Tube}}_{\mathbb X}[t] := {\tt{Tube}} \bigl(t, \sigma_{\mathbb X}(t), [0, 10A_0] \bigr),  \text{ where } \\
t \in \mathcal Q(J) \text{ and } \sigma_{\mathbb X}(t) \in \Omega_N. 
\end{equation} 
Here ${\tt{Tube}}$ is defined as in \eqref{defn: tube}; $\mathscr{U}_{\mathbb X}$ is the tree of height $N$ and $\sigma_{\mathbb X}: \mathscr{U}_{\mathbb X} \rightarrow \mathscr{B}_N$ is the sticky map defined in Section \ref{section: compressions of the unit interval}.
\end{definition} 
 \vskip0.1in
 \noindent Let us recall from \eqref{what is Q(k)} that $\mathcal Q(J)$ is the collection of $M$-adic subintervals of $[0,1]$ of length $M^{-J}$. The terminal vertices of $\mathscr{B}_N$ are the basic slope intervals in ${\tt{BasicSl}}_N(\mathcal S_N)$; by Proposition \ref{PROP: COMPRESSED SLOPES} \eqref{B_N identification}, each such vertex corresponds to a unique interval in $\mathcal V_J(\mathcal S_N)$, and hence to a unique slope in $\Omega_N$. 
 Accordingly, as described in Section \ref{section: finest scale extension} and \eqref{penultimate sigma defn}-\eqref{final sigma defn}, the tree map $\sigma_{\mathbb X}$ extends to a map
\[ \sigma_{\mathbb X}: \mathcal Q(J) \longrightarrow \Omega_N. \]
 Thus ${\tt{Tube}}_{\mathbb X}[t]$ is rooted at the $M$-adic interval $t$ of length $M^{-J}$ on the root line, with slope $\sigma_{\mathbb X}(t)$ assigned by the extension of the inductively sticky map $\sigma_{\mathbb X}$. We will refer to the collection 
 \begin{equation} \label{tube family X}   
 \mathbb T_{\mathbb X} := \{{\tt{Tube}}_{\mathbb X}[t] : t \in \mathcal Q(J) \} 
 \end{equation}  
 as the {\em{sticky tube family associated with $\mathbb X$}}. The union of this sticky tube family, denoted $\mathtt K_{\mathbb X}$
\begin{equation} \label{defn: extended Kakeya set K} 
\mathtt K_{\mathbb X} := \bigcup_{T \in \mathbb T_{\mathbb X}} T = \bigcup_{t \in \mathcal Q(J) } {\tt{Tube}}_{\mathbb X}[t] \subseteq A_0^2 [0,1]^2
\end{equation} 
will provide the source of the Kakeya-like configuration that we seek. 
\subsection{Two key properties of $\mathtt K_{\mathbb X}$} 
The proof of the implication \eqref{condition 3: slopes sublacunary} $\implies$ \eqref{condition 1: Kakeya-type sets} in Theorem \ref{thm:main} rests on two main pillars, summarized in the following propositions.
\subsection{The portion of $\tt{K}_{\mathbb X}$ near the root line is always large} 
 The geometric intuition behind the next result is easy to understand. The roots of the tube family $\mathbb T_{\mathbb X}$ are disjoint intervals by definition; this restricts their ability to overlap heavily in a suitable neighbourhood of the root line, even if their directions are well-separated. Thus their collective size is large, quantified by the following lower bound. 
\begin{proposition} \label{prop: near the roots}
There is a constant $c_1 > 0$ depending only on $C_0$ and $M \geq 2$ that satisfies the following property for all sufficiently large positive integers $N$. 
\vskip0.1in
\noindent Let $\Omega_{N} \subseteq [0,1]$ be a model $(N, C_0)$ slope set  with the properties listed in Proposition \ref{PRUNING STAGE 1}. Then for any choice of the binary sequence $\mathbb X$ given by \eqref{binary X}, the set ${\tt{K}}_{\mathbb X}$, defined as in \eqref{defn: extended Kakeya set K} obeys the estimate
\begin{equation} 
\bigl| {\tt{K}}_{\mathbb X} \cap \Bigl[[0,1] \times \mathbb R \Bigr] \bigr| \geq  c_1 \frac{\log N}{N}. 
\end{equation}  
\end{proposition}
\noindent Proposition \ref{prop: near the roots} is proved in Chapter \ref{chapter: near the roots}, Section \ref{section: proof of prop near roots}. 
\subsection{The portion of $\tt{K}_{\mathbb X}$ far from the root line can be small}
In contrast, as one moves away from the root line, tubes are much more likely to intersect, especially if many of their slopes are well-separated. This creates the opportunity for the tubes in $\mathbb T_{\mathbb X}$ to generate a much smaller combined area far from their roots. Such a phenomenon need not occur for every choice of the binary sequence $\mathbb X$, since different choices can produce substantially different patterns of overlap. The next proposition asserts that there is at least one admissible choice of $\mathbb X$ for which the required smallness occurs.
%Clearly, such a phenomenon cannot occur for every choice of the binary sequence $\mathbb X$; one can easily envision scenarios where all tubes are parallel, so there is no intersection anywhere in the plane. The next proposition confirms that while such extreme examples exist, they are not the norm.     
\begin{proposition} \label{prop: away from the roots} 
There is a constant $C_2 > 0$, depending only on $M \geq 2$ and $A_0$, that satisfies the following property.  
\vskip0.1in
\noindent For any positive integer $N$, let $\Omega_{N} \subseteq [0,1]$ be a model $(N, C_0)$ slope set  with the properties listed in Proposition \ref{PRUNING STAGE 1}. Then there exists a choice of a binary sequence $\mathbb X$ as in \eqref{binary X}, for which the set ${\tt{K}}_{\mathbb X}$ defined by \eqref{defn: extended Kakeya set K} obeys the estimate
\begin{equation} \label{estimate: away from the roots} 
\bigl| {\tt{K}}_{\mathbb X} \cap \Bigl[[A_0, A_0+1] \times \mathbb R \Bigr] \bigr| \leq \frac{C_2}{N}. 
\end{equation}
\end{proposition} 
\noindent Proposition \ref{prop: away from the roots} is proved in Chapter \ref{random construction section}. 
\vskip0.1in
\noindent Assuming these two propositions for now, we are in a position to complete the proof of one of the key statements in Theorem \ref{thm:main}. This is accomplished in Section \ref{section: Sublac-Kakeya}.  

\subsection{Tubes versus rectangles} \label{section: tubes & rectangles}
\noindent Another point is worth clarifying at this juncture. Each ${\tt{Tube}}(Q, \omega, I)$ is strictly speaking a parallelogram rather than a rectangle, because its short sides are vertical rather than perpendicular to $(1, \omega)$. However, the passage from a tube to a rectangle is benign, as we explain below. 
\vskip0.1in
\noindent By modifying each end of the tube within a region of diameter at most a constant multiple of $c_0|Q|$, we obtain two  actual rectangles, one contained in ${\tt{Tube}}$ and the other containing it: 
\[ {\tt{InnerRec}}(Q, \omega, I) \subseteq {\tt{Tube}}(Q, \omega, I)\subseteq {\tt{OuterRec}}(Q, \omega, I). \] 
Both rectangles have the same centre and orientation as ${\tt{Tube}}(Q, \omega, I)$, and their long-side lengths differ from that of the tube by at most a constant multiple of $c_0|Q|$. Their width is comparable to $c_0|Q|$; their short sides are perpendicular to $(1, \omega)$, as shown in Figure \ref{fig:tube-inner-outer-rectangles}. The area between the outer and inner rectangles is no more than a constant multiple of $(c_0 |Q|)^2$:
 \[ |{\tt{OuterRec}} \setminus {\tt{InnerRec}}| \leq (c_0 |Q|)^2.  \] 
\vskip0.1in
\noindent At the finest scale $|Q|=M^{-J}$, summing this error over all $M^J$ tubes contributes an error that is no more than a constant multiple of $M^{-J}$. Specifically, if ${\tt{InnerK}}_{\mathbb X}$ and ${\tt{OuterK}}_{\mathbb X}$ are defined exactly as in \eqref{defn: extended Kakeya set K}, but with tubes replaced by the corresponding rectangles {\tt{InnerRec}} and {\tt{OuterRec}},  then   
\[ {\tt{InnerK}}_{\mathbb X} \subseteq \mathtt{K}_{\mathbb X} \subseteq {\tt{OuterK}}_{\mathbb X}, \text{ with } \bigl|{\tt{OuterK}}_{\mathbb X} \setminus {\tt{InnerK}}_{\mathbb X} \bigr| \leq c_0^2 M^{-J}. \]
In particular, this means that the symmetric difference of the tube family $\mathtt{K}_{\mathbb X}$ with either one of the rectangle families ${\tt{InnerK}}_{\mathbb X}$ or ${\tt{OuterK}}_{\mathbb X}$ has area no more than a constant multiple of $M^{-J}$. Since $J > N$ by Proposition \ref{PRUNING STAGE 1}\eqref{partJN}, this error is negligible compared with both $N^{-1}$ and $N^{-1}\log N$, the two main bounds on portions of $\mathtt K_{\mathbb X}$ derived in Propositions \ref{prop: near the roots} and \ref{prop: away from the roots}. In view of this, we will continue to use the tube families for the estimates below; when verifying Definition \ref{Kakeya-type set} in Section \ref{section: Sublac-Kakeya}, we will pass to the corresponding family of outer rectangles.
\begin{figure}[ht]
\centering

\begin{tikzpicture}[
    scale=1.35,
    tube/.style={
        draw=orange!85!black,
        fill=orange!25,
        very thick
    },
    outer/.style={
        draw=blue!70!black,
        thick
    },
    inner/.style={
        draw=black,
        thick,
        dashed
    },
    construction/.style={
        draw=gray!70,
        thin,
        dashed
    },
    pointlabel/.style={
        font=\small
    }
]

% Original parallelogram ABCD.
\coordinate (B) at (0,1.15);
\coordinate (A) at (0,0);
\coordinate (C) at (6,2.45);
\coordinate (D) at (6,1.30);

% Feet of the perpendiculars.
% E is the projection of A onto the line BC.
% G is the projection of B onto the line AD.
% F and H are the corresponding points at the right endpoint.
\coordinate (E) at (-0.238,1.098);
\coordinate (G) at ( 0.238,0.052);
\coordinate (F) at ( 6.238,1.352);
\coordinate (H) at ( 5.762,2.398);

% Extensions of the two long sides of the tube.
\draw[construction] (E) -- (C);
\draw[construction] (A) -- (F);

% Outer rectangle E-C-F-A.
\draw[outer]
    (E) -- (C) -- (F) -- (A) -- cycle;

% Original tube ABCD.
\draw[tube]
    (B) -- (C) -- (D) -- (A) -- cycle;

% Inner rectangle B-H-D-G.
\draw[inner]
    (B) -- (H) -- (D) -- (G) -- cycle;

% Perpendiculars determining the rectangles.
\draw[outer] (A) -- (E);
\draw[inner] (B) -- (G);

% Right-angle markers.
\pic[
    draw=blue!70!black,
    angle radius=4mm
] {right angle=A--E--B};

\pic[
    draw=black,
    angle radius=4mm
] {right angle=B--G--A};

% Labels for the vertices of the tube.
\node[pointlabel,left]  at (A) {$A$};
\node[pointlabel,left]  at (B) {$B$};
\node[pointlabel,right] at (C) {$C$};
\node[pointlabel,right] at (D) {$D$};

\end{tikzpicture}

\caption{\small{
The orange parallelogram is the tube
$\operatorname{Tube}(Q,\omega,I)$.
The solid blue boundary gives the outer rectangle, while the dashed
black boundary gives the inner rectangle. The perpendicular from $A$
to the extension of $BC$ determines the outer rectangle; the
perpendicular from $B$ to $AD$ determines the inner rectangle.
}}
\label{fig:tube-inner-outer-rectangles}
\end{figure}
\section{Sublacunary slope sets generate Kakeya-type configurations} \label{section: Sublac-Kakeya} 
\subsection{Proof of \eqref{condition 3: slopes sublacunary} $\implies$ \eqref{condition 1: Kakeya-type sets} in Theorem \ref{thm:main}, given Propositions \ref{prop: near the roots}, \ref{prop: away from the roots}}
\begin{proof}
Condition \eqref{condition 3: slopes sublacunary} of Theorem \ref{thm:main} allows us to assume that the slope set $\Omega \subseteq [0,1]$ is sublacunary, in the sense of Definition \ref{defn: Admissible finite order lacunarity}.  By Corollary \ref{corollary: sublacunary slopes imply infinite split}, 
\[ \text{split}\bigl( \mathcal T(\Omega; M)\bigr) = \infty, \text{ where $ \mathcal T(\Omega; M)$ is the $M$-adic tree representing $\Omega$.}\]
Therefore, it obeys the hypothesis \eqref{pre-pruning split} of Proposition \ref{PRUNING STAGE 1} for any choice of $N$. The pruning mechanism specified by Proposition \ref{PRUNING STAGE 1} then allows us to extract a model $(N, C_0)$ slope set $\Omega_N$, which in turn leads to sticky tube families $\mathbb T_{\mathbb X}$ and their unions $\mathtt K_{\mathbb X}$ indexed by $\mathbb X$, defined respectively by \eqref{tube family X} and \eqref{defn: extended Kakeya set K}.  
\vskip0.1in 
\noindent We first construct a tube model for the desired Kakeya-type set $E_N$. For each $N \geq 1$, let $\mathbb X^{\ast}$ be the binary sequence identified by Proposition \ref{prop: away from the roots}.  We define 
\begin{align} 
&{\tt{KTube}}_{\mathbb X^{\ast}}[t] = {\tt{Tube}}_{\mathbb X^{\ast}}[t] \cap [A_0, A_0+1] \times \mathbb R, \text{ and set } \\   
&E_N := \bigcup \Bigl\{{\tt{KTube}}_{\mathbb X^{\ast}}[t] : t \in \mathcal Q(J) \Bigr\} = {\tt{K}}_{\mathbb X^{\ast}} \cap [A_0, A_0+1] \times \mathbb R. \label{our E_N}
\end{align} 
In other words, $E_N$ is a portion of ${\tt{K}}_{\mathbb X^{\ast}}$ away from the root line. It is built of sub-tubes of the form ${\tt{KTube}}_{\mathbb X^{\ast}}[t]$, which is the restriction of the longer ${\tt{Tube}}_{\mathbb X^{\ast}}[t]$ to the vertical strip above $[A_0, A_0+1]$. Each ${\tt{KTube}}_{\mathbb X^*}[t]$ is a parallelogram tube whose long side has slope $\sigma_{\mathbb X^*}(t)\in\Omega_N\subseteq\Omega$, whose
length is at least one, and whose width (i.e. perpendicular distance between its long sides) is at most 
\[
\delta_N:=c_0M^{-J}.
\]
Since the height $J$ of $\mathcal S_N$ satisfies $J > N$ by Proposition \ref{PRUNING STAGE 1} \eqref{partJN}, we have
$\delta_N\to 0$ as $N\to\infty$.
%Each ${\tt{KTube}}_{\mathbb X^{\ast}}[t]$ is essentially a rectangle of unit length oriented with slope $\sigma_{\mathbb X^{\ast}}(t) \in \Omega_N \subset \Omega$. 
\vskip0.1in
\noindent It remains to establish for the tube model $E_N$ given by \eqref{our E_N} the two measure estimates required in  \eqref{Kakeya-type condition}. On one hand, Proposition \ref{prop: away from the roots} gives that 
\begin{equation} \label{E_N: upper}
|E_N| = \bigl|{\tt{K}}_{\mathbb X^{\ast}} \cap \Bigl[[A_0, A_0+1] \times \mathbb R \Bigr] \bigr| \leq \frac{C_2}{N}. 
\end{equation}  
On the other hand, if $R =  {\tt{KTube}}_{\mathbb X^{\ast}}[t]$ is one of the constituent tubes of $E_N$, then we set
\[ A_1 :=2A_0+1.\]
Thus the constant $A_1$ is independent of $N$. As in the case for a rectangle, let $A_1R$ denotes the parallelogram obtained from $R$ by preserving its centre, orientation and width, while multiplying its long-side length by $A_1$. Then 
\[ A_1 R  \supseteq {\tt{Tube}}_{\mathbb X^{\ast}}[t] \cap \Bigl[[0, 1] \times \mathbb R \Bigr]. \]
For the tube model, let us write 
\[E_N^\ast(A_1):= \bigcup_{t\in\mathcal Q(J)} A_1 {\tt{KTube}}_{\mathbb X^\ast}[t]. \]
Therefore 
\[ E_N^{\ast}(A_1) \supseteq \bigcup_{t \in \mathcal Q(J)} {\tt{Tube}}_{\mathbb X^{\ast}}[t] \cap \Bigl[[0, 1] \times \mathbb R \Bigr] = {\tt{K}}_{\mathbb X^{\ast}} \cap \Bigl[[0,1] \times \mathbb R \Bigr]. \]
It follows from Proposition \ref{prop: near the roots} that
\begin{equation}  \label{E_n-star: lower} 
|E_N^{\ast}(A_1)| \geq c_1 \frac{\log N}{N}. 
\end{equation}
Combining \eqref{E_N: upper} and \eqref{E_n-star: lower}, we arrive at 
\[ \frac{|E_N^{\ast}(A_1)|}{|E_N|} \geq \frac{c_1}{C_2} \log N, \text{ which approaches } \infty \text{ as } N \rightarrow \infty.  \]    
\vskip0.1in
\noindent Strictly speaking, the set $E_N$ in \eqref{our E_N} is a union of parallelogram tubes. Replace each constituent tube by its corresponding outer rectangle from Section \ref{section: tubes & rectangles}, and continue to denote the resulting union by $E_N$. This increases its measure by at most a constant multiple of $c_0^2M^{-J}$. Moreover, since each outer rectangle contains the corresponding tube and has the same centre and long-side direction, its $A_1$-reach contains the reach used above. Thus, in view of the discussion in Section \ref{section: tubes & rectangles}, the bounds in \eqref{E_N: upper} and \eqref{E_n-star: lower} remain valid for the rectangle family $E_N$, with possibly different absolute constants. The resulting rectangles have slopes in $\Omega$ and lengths at least one. Their widths are at most a constant multiple of $\delta_N = c_0M^{-J}$, which converges to 0 as $N \rightarrow \infty$. They therefore satisfy all the requirements of a Kakeya-type set, as given by Definition \ref{Kakeya-type set}.
\vskip0.1in
\noindent This establishes that $\Omega$ admits Kakeya-type sets, proving \eqref{condition 1: Kakeya-type sets}. 
\end{proof}

\chapter{Tubes near the root line} \label{chapter: near the roots} 
The goal of this chapter is to prove Proposition \ref{prop: near the roots}, which provides a lower bound on the size of a sticky tube family near the root line. 
\section{Tubes in a fixed $M$-adic strip} 
\noindent Let us denote by 
\begin{equation} \label{defn: stripr}   
{\tt{Strip}}_r := \bigl[M^{-r}, M^{-r+1} \bigr] \times \mathbb R 
\end{equation} 
the vertical strip at distance roughly $M^{-r}$ from the root line.  For ${\tt{Tube}}_{\mathbb X}[t]$ defined as in \eqref{defn: sticky tube}, we set 
\begin{equation} \label{tube in strip} 
{\tt{Tube}}_{\mathbb X}[t; r] := {\tt{Tube}}_{\mathbb X}[t] \cap {\tt{Strip}}_r, 
\end{equation} 
so that ${\tt{Tube}}_{\mathbb X}[t; r]$ denotes the portion of ${\tt{Tube}}_{\mathbb X}[t]$ within the $r^{\text{th}}$ vertical strip. Thus ${\tt{Tube}}_{\mathbb X}[t; r]$ is itself a tube with long and short sides comparable to $M^{-r}$ and $M^{-J}$ respectively. The main observation concerning the collective size of these truncated tubes is a uniform bound in $r$, which says that for a large number of consecutive strips near the root line (a number that increases with $N$), the total area of the truncated tubes remain bounded below by a constant multiple of $1/N$. 
\begin{proposition} \label{prop: tubes on a strip}
Let us assume the hypothesis of Proposition \ref{prop: near the roots}. Then for $c_1^{-1} \leq r \leq \log_{M}(N)$, the following estimate holds for all sufficiently large integers $N$: 
\begin{equation} \label{K in a strip} 
\bigl| \mathtt K_{\mathbb X} \cap {\tt{Strip}}_r\bigr| = \Bigl| \bigcup_{t \in \mathcal Q(J)} {\tt{Tube}}_{\mathbb X}[t; r] \Bigr| \geq \frac{c_1}{N} 
\end{equation}  
\end{proposition}   
\noindent We prove Proposition \ref{prop: tubes on a strip} in the next section. Assuming this, Proposition \ref{prop: near the roots} is proved as follows. 
\subsection{Proof of Proposition \ref{prop: near the roots}, assuming Proposition \ref{prop: tubes on a strip}}  \label{section: proof of prop near roots} 
\begin{proof} 
Since the sets $\{{\tt{Strip}}_r: r \geq 0 \}$ have disjoint interiors, we may write 
\begin{align*} \Bigl| \mathtt K_{\mathbb X} \cap \Bigl[ [0,1] \times \mathbb R \Bigr] \Bigr| 
&\geq \sum_{r=1/c_1}^{\log_M N} \bigl| \mathtt K_{\mathbb X} \cap {\tt{Strip}}_r \bigr|  \\
& \geq \sum_{r} \Bigl\{ \bigl| \mathtt K_{\mathbb X} \cap {\tt{Strip}}_r \bigr| : \frac{1}{2} \log_M N \leq r \leq \log_M N \Bigr\} \\ 
&\geq c_1 \frac{\log_M N}{N}, \text{ whenever } \log_{M}N > \frac{2}{c_1}.
\end{align*}
The last step uses the conclusion \eqref{K in a strip} of Proposition \ref{prop: tubes on a strip}.  
\end{proof}

\section{Pairwise intersections of truncated tubes near the root line} 
The main result in this section is an upper bound on the total size of pairwise intersections taking place in ${\tt{Strip}}_r$ among the tubes in $\mathtt K_{\mathbb X}$. 
\begin{proposition} \label{PAIRWISE PROP} 
There exists a constant $C_1 = C_1(M, C_0) > 0$ such that, in the set-up of Proposition \ref{prop: tubes on a strip}, the following estimate holds for all $r \geq C_1$ and all sufficiently large $N$:
\begin{equation} \label{pairwise upper bound}
\sum_{\begin{subarray}{c}t_1, t_2 \in \mathcal Q(J) \\ t_1 \neq t_2 \end{subarray}}\Bigl| {\tt{Tube}}_{\mathbb X}[t_1; r] \cap {\tt{Tube}}_{\mathbb X}[t_2; r] \Bigr| \leq C_1 \frac{N}{M^{2r}}.
\end{equation} 
\end{proposition}
\noindent We will prove Proposition \ref{PAIRWISE PROP} momentarily. Assuming it for now, let us complete the proof of Proposition \ref{prop: tubes on a strip} in Section \ref{proof: tubes on a strip}. The intermediate section supplies a necessary ingredient. 
\subsection{Controlled pairwise intersection implies a large union}
A family of sets whose pairwise intersections are controlled cannot overlap too much, hence their combined size should be large.  This general measure-theoretic observation has been used in \cite{{BatemanKatz},{Bateman}} to estimate the size of tube families near the root line. We will use it for the same purpose.    
\begin{lemma}[Lemma 8, \cite{Bateman}] \label{measure theory lemma} 
Let $\{\mathfrak A_r : 1 \leq r \leq \mathfrak R\}$ be measurable sets in a general measure space $(X, \mathfrak X, \mu)$, with 
\begin{equation}  \mu(\mathfrak A_r) \geq \alpha \text{ for all } 1 \leq r \leq \mathfrak R. \end{equation} 
If 
\begin{equation} 
\sum_{r, s = 1}^{\mathfrak R} \mu \bigl( \mathfrak A_r \cap \mathfrak A_s \bigr) \leq \mathfrak M, 
\end{equation} 
then 
\begin{equation} \label{lower bound on union}
\mu \Bigl( \bigcup_{r=1}^{\mathfrak R} \mathfrak A_r \Bigr) \geq \frac{\alpha^2 \mathfrak R^2}{16 \mathfrak M}. 
\end{equation} 
\end{lemma} 
\subsection{Proof of Proposition \ref{prop: tubes on a strip}, assuming Proposition \ref{PAIRWISE PROP}} \label{proof: tubes on a strip} 
\begin{proof} 
The desired estimate \eqref{K in a strip} is a direct consequence of the measure-theoretic Lemma \ref{measure theory lemma}. Let us fix an index $r$ with $C_1 \leq r \leq \log_M N$. Since 
\[ K_{\mathbb X} \cap {\tt{Strip}}_r = \bigcup_{t \in \mathcal Q(J)} {\tt{Tube}}_{\mathbb X}[t ; r], \]
we are in the set-up of Lemma \ref{measure theory lemma}, with the following choice of parameters: 
\begin{align} 
\mathfrak A_{t} &= {\tt{Tube}}_{\mathbb X}[t ; j], \quad \alpha = c_0 M^{-J-r}, \quad \mathfrak R = M^J, \text{ and } \label{parameter set 1} \\
\mathfrak M &= \sum_{t_1, t_2 \in \mathcal Q(J)} \Bigl| {\tt{Tube}}_{\mathbb X}[t_1; r] \cap {\tt{Tube}}_{\mathbb X}[t_2; r] \Bigr| \nonumber \\
 &\leq \Bigl[ \sum_{t_1 = t_2} + \sum_{t_1 \neq t_2} \Bigr]   \Bigl| {\tt{Tube}}_{\mathbb X}[t_1; r] \cap {\tt{Tube}}_{\mathbb X}[t_2; r] \Bigr| \nonumber \\
 & \leq \#(\mathcal Q(J)) c_0 M^{-J-r} + C_1 \frac{N}{M^{2r}} = c_0 M^{-r} + C_1 \frac{N}{M^{2r}} \nonumber 
\end{align}   
where the last estimate uses \eqref{pairwise upper bound} from Proposition \ref{PAIRWISE PROP}. Since 
\[ M^{-r} \leq \frac{N}{M^{2r}} \; \text{ for } \; 0 \leq r \leq \log_{M}N, \text{ and in particular for } C_1 \leq r \leq \log_{M}N, \]
we deduce that for $r$ in the latter range, 
\begin{equation} \label{parameter set 2}   
\mathfrak M \leq 2C_1 \frac{N}{M^{2r}}. 
\end{equation} 
Substituting \eqref{parameter set 1} and \eqref{parameter set 2} into \eqref{lower bound on union}, we obtain
\begin{align*} 
\bigl| K_{\mathbb X} \cap {\tt{Strip}}_r \bigr| &= \Bigl| \bigcup_{t \in \mathcal Q(J)} {\tt{Tube}}_{\mathbb X}[t ; r] \Bigr| \\
&\geq \frac{\alpha^2 \mathfrak R^2}{16 \mathfrak M^2} \geq c_1 \frac{M^{-2J-2r} M^{2J}}{NM^{-2r}} \geq \frac{c_1}{N}.  
\end{align*}  
This is the claimed estimate \eqref{K in a strip}. 
\end{proof}
\vskip0.1in
\noindent To recap, we have reduced the proof of Proposition \ref{prop: near the roots} from Proposition \ref{prop: tubes on a strip} to Proposition \ref{PAIRWISE PROP}. The remainder of this chapter is devoted to the proof of the latter. The main geometric content concerning overlapping tubes appears here. We gather the necessary facts concerning tube intersections in Sections \ref{section: tube geometry} and \ref{section: counting tube intersections} before embarking on the actual proof in Section \ref{section: pairwise prop proof}. 
\section{Intersection of two tubes} \label{section: tube geometry}
\subsection{Conditions for intersection} \label{section: intersection of tubes} 
The summand appearing on the left hand side of \eqref{pairwise upper bound} is non-zero only if the two tubes ${\tt{Tube}}_{\mathbb X}[t_1]$ and ${\tt{Tube}}_{\mathbb X}[t_2]$ intersect in ${\tt{Strip}}_r$. The following lemma derives a necessary criterion for this. 
\begin{lemma} \label{intersection lemma 2} 
Given any binary sequence $\mathbb X$, and tubes ${\tt{Tube}}_{\mathbb X}[t_i; r]$, $i= 1, 2$ defined as in \eqref{tube in strip}, suppose that 
\begin{equation}  
{\tt{Tube}}_{\mathbb X}[t_1; r] \cap {\tt{Tube}}_{\mathbb X}[t_2; r] \neq \emptyset, \; t_1 \neq t_2. \label{intersection in strip}  \end{equation} 
Then there is a constant $C_1 > 1$ depending only on $M$ and $C_0$ such that 
\begin{equation} \label{Euclidean-distance t1t2}
C_1^{-1} \leq M^{r + \lambda(v)} \bigl|\text{cen}(t_1) - \text{cen}(t_2) \bigr| \leq C_1.  
\end{equation}  
Here $\text{cen}(t)$ denotes the centre of the interval $t$, and $\lambda(v)$ denotes the fundamental height corresponding to the splitting vertex
\begin{equation} \label{yca in B_N} 
v = D_{\mathcal S_N}(\omega_1, \omega_2), \quad \omega_i := \sigma_{\mathbb X}(t_i), \; i= 1,2. 
\end{equation}
In other words, $v$ is the youngest common ancestor of $\omega_1$ and $\omega_2$ in the pruned slope tree $\mathcal S_N$. 
\end{lemma}
\vskip0.1in
\noindent {\em{Remark: }} Distinct roots $t_1,t_2\in\mathcal Q(J)$ have disjoint $c_0$-dilates. If their assigned slopes coincide, the corresponding tubes are parallel and disjoint. Thus \eqref{intersection in strip} implies $\omega_1\ne\omega_2$, so the splitting vertex $v$ in \eqref{yca in B_N} is well-defined.
\begin{proof} 
According to the definition \eqref{defn: tube} of a tube, each tube ${\tt{Tube}}_{\mathbb X}[t_i]$ is a union of parallel line segments of slope $\omega_i$ originating from the interval $\{0 \} \times \tilde{t}_i$ on the root line. The intersection condition \eqref{intersection in strip} therefore implies the existence of three points 
\begin{equation} \label{PA1A2}  \left\{
\begin{aligned} 
&P \in {\tt{Strip}}_r \text{ and $A_1 = (0, y_1)$, $A_2 = (0, y_2)$ on the root line, with} \\
&\text{$y_i \in \tilde{t}_i$ such that  the line $A_i P$ has slope $\omega_i = \tan \theta_i$ for $i=1,2$. }
\end{aligned} 
\right\}
\end{equation}
Applying the law of sines to the triangle $\Delta PA_1A_2$ leads to 
\begin{equation}  \label{law of sines} 
\frac{|A_1-A_2|}{|\sin(\theta_1-\theta_2)|} = \frac{|PA_1|}{|\sin(\frac{\pi}{2} \pm \theta_2)|}. 
\end{equation}  
Since $\omega_i = \tan \theta_i \in [0, 1]$, we conclude that $\theta_i \in [0, \frac{\pi}{4}]$. This restriction on $\theta_i$ permits the following estimates 
\begin{align} \label{cosine bounds}
&C_1^{-1} \leq |\sin(\frac{\pi}{2} \pm \theta_i| = |\cos \theta_i|  \leq 1, \; \text{ and } \\ 
&C_1^{-1} |\omega_1-\omega_2| \leq |\sin(\theta_1-\theta_2)| \leq C_1 |\omega_1 - \omega_2|. \label{cosine bounds 1}
\end{align} 
The location of $P$ given by \eqref{PA1A2} and the definition \eqref{defn: stripr} of ${\tt{Strip}}_r$ yield 
 \begin{equation} \label{PA1}
 |PA_1| \cos\theta_1 = \text{dist($P$, root line)} \in [M^{-r}, M^{-r+1}]. 
 \end{equation} 
Incorporating the relation \eqref{cosine bounds} into \eqref{PA1} leads to a size bound on $|PA_1|$:
\begin{equation} \label{cosine bounds 2} 
C_1^{-1} M^{-r} \leq |PA_1| \leq C_1 M^{-r}.    
\end{equation}   
 Substituting \eqref{cosine bounds}, \eqref{cosine bounds 1} and \eqref{cosine bounds 2} into \eqref{law of sines}, we arrive at the two-sided inequality 
 \begin{equation} 
 \begin{aligned} 
  |y_1-y_2| = |A_1 - A_2| &= \frac{|PA_1| \times |\sin(\theta_1 - \theta_2)|}{|\cos \theta_2|}  \\
 &\in M^{-r} |\omega_1- \omega_2|\left[ C_1^{-1}, C_1\right]. 
\end{aligned} \label{pre-conclusion}
\end{equation} 
We claim that the desired conclusion \eqref{Euclidean-distance t1t2} follows from \eqref{pre-conclusion}. On one hand, we have the two inequalities  
\begin{align*} 
&\Bigl| |\text{cen}(t_1) - \text{cen}(t_2)| -  |y_1-y_2| \Bigr| \leq \bigl|\text{cen}(t_1) - y_1 \bigr| + \bigl|\text{cen}(t_2) - y_2 \bigr| \\
&\hskip2in \leq \text{diam}(t_1) + \text{diam}(t_2) \leq 2c_0 M^{-J}, \text{ and }  \\
&|\text{cen}(t_1) - \text{cen}(t_2)| \geq M^{-J}, \text{ since }  t_1, t_2 \in \mathcal Q(J), \; t_1 \neq t_2.
\end{align*} 
Jointly, they imply 
\begin{equation} \label{y-centres} 
C_1^{-1} |y_1-y_2| \leq  |\text{cen}(t_1) - \text{cen}(t_2)| \leq C_1 |y_1-y_2|.
\end{equation} 
On the other hand, $\omega_i \in \Omega_N$, therefore by the conclusion \eqref{distance between slopes} of Lemma \ref{lemma: distance between pruned slopes},
\begin{equation} \label{difference-omega}
C_1 \leq |\omega_1 - \omega_2| M^{\lambda(v)} \leq C_1, \text{ with $v$ as in \eqref{yca in B_N}}. 
\end{equation}
Combining \eqref{pre-conclusion} with \eqref{y-centres} and \eqref{difference-omega} leads to \eqref{Euclidean-distance t1t2}, completing the proof.   
\end{proof} 
\subsection{Size of intersection} 
Lemma \ref{intersection lemma 2} provides a criterion for intersection of two tubes, in the form of an algebro-geometric inequality. We will also need to know the size of this intersection. This estimate is by now standard in the literature, dating back to the work of C\'ordoba \cite{Cordoba}. The result below is easily verifiable, but the reader may consult \cite[Lemma 10.3.6, p.~374]{Grafakos} as a reference.   
\begin{lemma} \label{intersection size lemma} 
If ${\tt{Tube}}$ and ${\tt{Tube}}'$ are any two intersecting tubes of the form \eqref{defn: tube}, then 
\begin{equation}  \label{intersection size 1}
| {\tt{Tube}} \cap {\tt{Tube}}'| \leq \frac{C_0 M^{-2J}}{M^{-J} + |\omega - \omega'|},
\end{equation} 
where $C_0$ is an absolute constant.   
\end{lemma}
\vskip0.1in
\noindent Combining Lemma \ref{intersection size lemma} with the Euclidean separation properties embedded in the pruned slope tree $\mathcal S_N$, we arrive at the following consequence. 
\begin{corollary} \label{corollary: intersection size} 
There is a constant $C_1 > 0$ depending only on $M$ and $C_0$ with the following property. If $\omega, \omega' \in \Omega_N$, $\omega \neq \omega'$, then
\begin{equation} \label{intersection size 2} 
 | {\tt{Tube}} \cap {\tt{Tube}}'| \leq C_1 M^{-2J + \lambda(v)},
\end{equation} 
where $v$ is given by \eqref{yca in B_N}. 
\end{corollary} 
\begin{proof} 
The estimate \eqref{intersection size 2} is a direct consequence of \eqref{intersection size 1}. Since $\omega \neq \omega'$, the choice of the height $J$ and the property \eqref{Euclidean separation terminal vertices} of $\mathcal S_N$ ensure that 
\[ |\omega - \omega'| \geq C_0 M^{-J}, \; \text{ so that } \; M^{-J} + |\omega-\omega'| \geq c_0 |\omega - \omega'|  \]  
for some small constant $c_0 > 0$.  This means that 
\[ | {\tt{Tube}} \cap {\tt{Tube}}'| \leq C_1 \frac{M^{-2J}}{|\omega-\omega'|}. \] 
Lemma \ref{lemma: distance between pruned slopes} then specifies that $|\omega - \omega'|$ is comparable to $M^{-\lambda(v)}$, completing the proof. 
\end{proof} 

\section{Counting pairs of intersecting tubes} \label{section: counting tube intersections} 
Section \ref{section: intersection of tubes} provided criteria for two tubes to intersect in ${\tt{Strip}}_r$. This section is focused on estimating the number of tube pairs, constructed via a single root-to-slope map $\sigma_{\mathbb X}$, that obey this intersection condition. Specifically, we are interested in counting the number of tube pairs identified by $t_1, t_2 \in \mathcal Q(J)$ that satisfy 
\begin{equation} \label{Tube1 and Tube2} 
{\tt{Tube}}_{\mathbb X}[t_1; r] \cap  {\tt{Tube}}_{\mathbb X}[t_2; r] \cap \bigl[M^{-r}, M^{-r+1} \bigr] \neq \emptyset, 
\end{equation}
where ${\tt{Tube}}_{\mathbb X}[t; r]$ denotes the tube defined by \eqref{defn: sticky tube}. Accordingly, we define the collection ${\tt{TP}}_{\mathbb X}[r]$, consisting of distinct tube pairs with this intersection property:  
\begin{equation}
{\tt{TP}}_{\mathbb X}[r] := \\
\left\{ (t_1, t_2) \in \mathcal Q(J)^2 : t_1 \neq t_2, \; \eqref{Tube1 and Tube2} \text{ holds} \right\}.  
\end{equation} 
We will also define sums associated with families of tube pairs.  For instance, the left hand side of \eqref{pairwise upper bound} is denoted by  
\begin{equation}  \label{defn: SumTP}
{\tt{SumTP}}_{\mathbb X}[r] := \sum \Bigl\{ \left| {\tt{Tube}}_{\mathbb X}[t_1; r] \cap {\tt{Tube}}_{\mathbb X}[t_2; r] \right|  \, : \, (t_1, t_2) \in {\tt{TP}}_{\mathbb X}[r] \Bigr\}.
\end{equation}
The aim of Proposition \ref{PAIRWISE PROP} is to provide an  upper bound for ${\tt{SumTP}}_{\mathbb X}[r]$.  

\subsection{Decomposition of ${\tt{TP}}_{\mathbb X}[r]$} With this goal in mind, let us partition 
\begin{equation} \label{TP-partition}
{\tt{TP}}_{\mathbb X}[r] = \bigsqcup \left\{ {\tt{TP}}_{\mathbb X}[k, u, t; r] \; \Biggl| \; \begin{aligned} &1 \leq k \leq N, \; t \in \mathcal V_k(\mathscr{U}_{\mathbb X}) \\ &u \in {\tt{SplitV}}_k(\mathcal S_N) \end{aligned} \right\},
\end{equation} 
based on the youngest common ancestor $t$ of $t_1$ and $t_2$ in the compressed root tree $\mathscr{U}_{\mathbb X}$. The decomposition \eqref{TP-partition} uses the height $k$ of $t$ in $\mathscr{U}_{\mathbb X}$ and its image under $\sigma_{\mathbb X}$. The latter, by definition of a tree-map pair, is a basic slope interval at the $k^{\text{th}}$ fundamental height of $\mathcal S_N$, of the same length as $t$. Specifically, we define 
\begin{equation} \label{TP-piece} 
{\tt{TP}}_{\mathbb X}[k, u, t; r] := \left\{(t_1, t_2) \in {\tt{TP}}_{\mathbb X}[r] \, \Biggl| \,
\begin{aligned} 
&t = D_{\mathscr{U}_{\mathbb X}}(t_1, t_2) \text{ obeys } \;  \\ 
& t \in \mathcal V_k(\mathscr{U}_{\mathbb X}), \; \sigma_{\mathbb X}(t) = \theta(u)  
\end{aligned} 
\right\}. 
\end{equation}
Here $\theta(u)$ is a basic slope interval descended from the $k^{\text{th}}$ splitting vertex $u \in \mathcal S_N$; the length of $\theta(u)$ therefore corresponds to the fundamental height $\lambda(u)$.  
To clarify,  $u$ signifies the youngest splitting ancestor of $\sigma_{\mathbb X}(t)$ in $\mathcal S_N$, so that
\[|t| =  |\sigma_{\mathbb X}(t)| = |\theta(u)| = M^{- \lambda(u)}. \]   
The decomposition \eqref{TP-partition} results in a decomposition of the sum \eqref{defn: SumTP}, indexed by $k, u$ and $t$: 
\begin{align} 
&{\tt{SumTP}}_{\mathbb X}[r] = \sum_{k=1}^N \sum_{u \in {\tt{SplitV}}_k} \sum_{\begin{subarray}{c} t \in \mathcal V_k(\mathscr{U}_{\mathbb X}) \\ \sigma_{\mathbb X}(t) \subsetneq u \end{subarray}} {\tt{SumTP}}_{\mathbb X}[k, u, t; r], \text{ where } \label{subsum-1}\\ 
&{\tt{SumTP}}_{\mathbb X}[k,u, t; r] := \sum_{(t_1, t_2)}' \left| {\tt{Tube}}_{\mathbb X}[t_1; r] \cap {\tt{Tube}}_{\mathbb X}[t_2; r] \right| \label{subsum-2}
\end{align}
%is the sum over all pairs $(t_1, t_2) \in {\tt{TP}}_{\mathbb X}[t,u,k; r]$.   
%\vskip0.1in
%\noindent In order to reduce \eqref{subsum-2} to blocks where each summand is of uniform size, we make one last round of decomposition on ${\tt{TP}}_{\mathbb X}[t, u, k; r]$. The estimate \eqref{intersection size 2} specifies the size of the summand in terms of the youngest common ancestor of the slopes $\sigma_{\mathbb X}(t_1)$ and $\sigma_{\mathbb X}(t_2)$ in the pruned slope tree $\mathcal S_N$. Motivated by this, we set  
%\begin{align} 
%{\tt{TP}}_{\mathbb X}[t, u, k; r] &= \bigsqcup_{v \in {\tt{SplitV}}} {\tt{TP}}_{\mathbb X}[v, t, u, k; r], \text{ where } \label{TP-1} \\  
%{\tt{TP}}_{\mathbb X}[v, t, u, k; r] &= \left\{(t_1, t_2) \in {\tt{TP}}_{\mathbb X}[t, u, k; r]  : v = D_{\mathcal S_N} \bigl(\sigma_{\mathbb X}(t_1), \sigma_{\mathbb X}(t_2)\bigr) \right\} \label{TP-2}
%\end{align}  
%The corresponding sum ${\tt{SumTP}}_{\mathbb X}[v,k,u, t; r]$ is defined analogous to \eqref{subsum-2}. 
%\vskip0.1in
%\noindent The sticky nature of $\sigma_{\mathbb X}$ provides a mechanism for comparing the fundamental heights of $u$ and $v$. 
The following lemma captures the relation between the fundamental heights associated with the youngest common ancestors of the roots $(t_1, t_2)$ and their images. The strictness of the inequality \eqref{lambda relation} is important. 
\begin{lemma} \label{lemma: lambda relation} 
Suppose that $(t_1, t_2)$ lies in the collection ${\tt{TP}}_{\mathbb X}[k, u, t; r]$ given by \eqref{TP-piece}. If
\begin{equation}  \label{lambda relation} 
 v = D_{\mathcal S_N}\bigl(\sigma_{\mathbb X}(t_1), \sigma_{\mathbb X}(t_2) \bigr) \; \text{ then } \; \lambda(u) < \lambda(v). \end{equation}  
\end{lemma} 
\begin{proof} 
Proposition \ref{PROP: COMPRESSED ROOTS} posits that $(\mathscr{U}_{\mathbb X}, \sigma_{\mathbb X})$ is a tree-map pair adapted to $\mathscr{B}_N$. According to the property \eqref{sigma preserves lineage} of such a pair, 
\[ t_1, t_2 \subset t \; \text{ implies } \; \sigma_{\mathbb X}(t_1), \sigma_{\mathbb X}(t_2) \subset \sigma_{\mathbb X}(t).  \]
Since $v$ is the youngest common ancestor of $\sigma_{\mathbb X}(t_1)$, $\sigma_{\mathbb X}(t_2)$ in the $M$-adic tree, we obtain 
\[ v \subseteq   \sigma_{\mathbb X}(t), \; \text{ which in turn implies } \; M^{-h(v)} =  |v| \leq  |\sigma_{\mathbb X}(t)| = M^{- \lambda(u)}. \]
The last statement leads to 
\[ \lambda(u) \leq h(v) < \lambda(v), \]
where the final strict inequality is a consequence of \eqref{h-star-v} and the definition of a fundamental height, namely $\lambda(v) = h_v^{\ast}$. This completes the proof.  
\end{proof} 
\vskip0.1in
\noindent The conclusion \eqref{lambda relation} of Lemma \ref{lemma: lambda relation} allows us to define a new integer parameter $\ell \geq 1$, which will be a key ingredient in the estimation:
\begin{equation} \label{lambda difference} 
\lambda(v) = \lambda(u) + \ell. 
\end{equation}   
Let us point out that $\ell$ does not necessarily run over all positive integers. Since $v \subsetneq \theta(u)$, it must be contained in the first splitting descendant of $u$, since the ray of $\mathcal S_N$ joining the latter to $\theta(u)$ is non-splitting. Specifically, 
\begin{align} \label{what is l} 
&\ell \geq \lambda(u') - \lambda(u), \; \text{ where } u' \in {\tt{SplitV}}(\mathcal S_N),\\ 
& u' \subsetneq \theta(u) \subsetneq u, \; \iota(u') = \iota(u) + 1 = k+1.  \label{what is u'}
\end{align} 
We make a last round of decomposition based on $\ell$: 
\begin{align} 
{\tt{TP}}_{\mathbb X}[k, u, t; r] &= \bigsqcup_{v \in {\tt{SplitV}}} {\tt{TP}}_{\mathbb X}[\ell, k, u, t; r], \text{ where } \label{TP-1} \\  
{\tt{TP}}_{\mathbb X}[\ell, k, u, t; r] &= \left\{(t_1, t_2) \in {\tt{TP}}_{\mathbb X}[k, u, t; r]  : \eqref{lambda difference} \text{ holds} \right\} \label{TP-2}
\end{align}  
The corresponding quantity ${\tt{SumTP}}_{\mathbb X}[\ell,k,u, t; r]$ is defined analogously to \eqref{subsum-2}. We have arrived at the identity:  
\begin{equation} \label{subsum-3} 
{\tt{SumTP}}_{\mathbb X}[k,u, t; r] = \sum_{\ell=1}^{\infty} {\tt{SumTP}}_{\mathbb X}[\ell, k,u, t; r].
\end{equation}  
\subsection{A counting lemma} Our next result estimates the cardinality of ${\tt{TP}}_{\mathbb X}[v, k, u, t; r]$ and the resulting bound on the sub-sums  \eqref{subsum-3} and \eqref{subsum-2}. 
\begin{lemma} \label{subsum counting lemma} 
For the collection of root pairs given by \eqref{TP-2}, the following size estimate holds:
\begin{equation} \label{TP-counting-estimate} 
\# \Bigl[ {\tt{TP}}_{\mathbb X}[\ell, k, u, t; r]\Bigr] \leq C_1 M^{2J-2r} \times M^{- 3\lambda(u) + \lambda(u') - 2 \ell},
\end{equation} 
where $u'$ as in \eqref{what is u'} is the unique splitting vertex of index $(k+1)$ in $\mathcal S_N$ descended from $\theta(u)$.  
\vskip0.1in
\noindent As a result, we obtain the following estimate of the quantity \eqref{subsum-2}: 
\begin{equation} \label{subsum-estimate}
{\tt{SumTP}}_{\mathbb X}[k,u, t; r]  \leq C_1 M^{-2r - \lambda(u)}.  
\end{equation} 
\end{lemma} 
\begin{proof} 
Let us first establish \eqref{subsum-estimate} from \eqref{TP-counting-estimate}, combining the latter with \eqref{intersection size 2}. Indeed, Corollary \ref{corollary: intersection size} implies that for any $(t_1, t_2) \in  {\tt{TP}}_{\mathbb X}[v, k, u, t; r]$, 
\[ \left| {\tt{Tube}}_{\mathbb X}[t_1; r] \cap {\tt{Tube}}_{\mathbb X}[t_2; r] \right| \leq C_1 M^{-2J + \lambda(v)} \leq C_1 M^{-2J + \lambda(u) + \ell}. \] 
This means that 
\begin{align*} 
{\tt{SumTP}}_{\mathbb X}[\ell, k,u, t; r]  &\leq C_1 M^{-2J + \lambda(v)} \times \# \Bigl[ {\tt{TP}}_{\mathbb X}[v, k, u, t; r]\Bigr]  \\
&\leq C_1 M^{-2J + \lambda(u) + \ell} \times M^{2J-2r} \times M^{- 3 \lambda(u) + \lambda(u') - 2 \ell} \\
&\leq C_1 M^{-2r - 2\lambda(u) + \lambda(u') - \ell}.
\end{align*}  
Inserting the last estimate on the right hand side of \eqref{subsum-3}, and summing in $\ell$ subject to the restriction \eqref{what is l} , we find that 
\begin{align*} 
{\tt{SumTP}}_{\mathbb X}[k,u, t; r]  &\leq C_1 M^{-2r - 2 \lambda(u) + \lambda(u') }\sum_{\ell = \lambda(u') - \lambda(u)}^{\infty} M^{-\ell} \\ 
&\leq C_1 M^{-2r - 2 \lambda(u) + \lambda(u') } \times M^{-\lambda(u') + \lambda(u)} \leq C_1 M^{-2r-\lambda(u)}.
\end{align*} 
This completes the proof of \eqref{subsum-estimate}, given \eqref{TP-counting-estimate}. 
\vskip0.1in
\noindent We now turn to the proof of \eqref{TP-counting-estimate}, the cardinality estimate for ${\tt{TP}}_{\mathbb X}[\ell, k, u, t; r]$. This requires us to count the number of pairs $(t_1, t_2)$ obeying three properties: 
\vskip0.1in
\begin{itemize} 
\item the non-trivial intersection condition \eqref{Tube1 and Tube2} for the corresponding tubes,
\vskip0.1in 
\item the size and ancestry relations
\begin{equation} \label{t-t1t2} 
t = D_{\mathscr{U}_{\mathbb X}}(t_1, t_2) \in \mathcal V_k({\mathscr U_{\mathbb X}}),  \quad |t| = M^{-\lambda(u)}.  
\end{equation} 
\vskip0.1in 
\item the separation of fundamental heights given by \eqref{lambda difference}, where $u, v$ are as in \eqref{TP-piece} and \eqref{lambda relation}.  
\end{itemize}
The strategy is as follows: we will first count the number of pairs $(Q_1, Q_2)$, where each $Q_i$ is a child of $t$ with the potential to generate an intersecting pair $(t_1, t_2) \in {\tt{TP}}_{\mathbb X}[k, u, t; r]$; this means  
\begin{equation} \label{tQ}
t_i \subseteq Q_i, \quad i =1,2.
\end{equation}  We will then count the number of descendants $t_i$ from each $Q_i$ to complete the estimate.   
\vskip0.1in
\noindent To make this precise, let us define for fixed $\ell, k, u, t$ the collection of all possible ancestral tube-pairs at height $(k+1)$:
\begin{equation} 
\widetilde{\tt{TP}} := \left\{ (Q_1, Q_2) \, \Biggl| \, \begin{aligned}   & Q_i \in \mathcal V_{k+1}(\mathscr{U}_{\mathbb X}) \text{ for } i = 1, 2; \, \exists (t_1, t_2)  \\ & \text{in } {\tt{TP}}_{\mathbb X}[\ell, k, u, t ; r] \ni \eqref{tQ} {\text{ holds}} \end{aligned} \right\}. 
\end{equation} 
By the properties \eqref{compressed children 1} and \eqref{compressed children 2} in Proposition \ref{PROP: COMPRESSED ROOTS}, the vertex $t$ produces $M^{\lambda(u') - \lambda(u)}$ children at generation $(k+1)$, each an interval of length $M^{-\lambda(u')}$. On one hand, in view of \eqref{t-t1t2}, the roots $t_1, t_2$ must lie in separate children at this height, and must obey 
\begin{equation}  \label{dist-t1t2} 
\text{dist}(t_1, t_2) \leq C_1 M^{-r - \lambda(v)} \leq C_1 M^{-r - \lambda(u) - \ell}. 
\end{equation}    
This last condition uses \eqref{lambda difference} and follows from \eqref{Euclidean-distance t1t2} in Lemma \ref{intersection lemma 2}, as a consequence of 
the non-trivial intersection \eqref{Tube1 and Tube2} between the two tubes ${\tt{Tube}}_{\mathbb X}[t_i; r]$, $i=1, 2$. 
\vskip0.1in 
\noindent On the other hand, the right hand side of \eqref{dist-t1t2} is at most $C_1 M^{-r - \lambda(u')}$, which is strictly smaller than $M^{-\lambda(u')}/10$ if 
\[ C_1 M^{-r} < \frac{1}{10} \; \text{ i.e. if } \; r \geq \log_M (10C_1). \]  
To meet the separation restriction \eqref{dist-t1t2}, the vertices $t_1, t_2 \in \mathscr{U}_{\mathbb X}$ are thus forced to lie in adjacent children $Q_1, Q_2$ of $t$, separated by a distance of at most $C_1 M^{-r - \lambda(u) - \ell}$ across the boundary of $Q_1$ and $Q_2$. Specifically, this means 
\begin{equation} \label{TP-tilde count} 
\# \bigl[ \widetilde{\tt{TP}}\bigr]
%&\quad \leq \# \left\{ (Q_1, Q_2) \in \in \Bigl[\mathcal V_{k+1}(\mathscr{U}_{\mathbb X}) \Bigr]^2 \; \Biggl| \; \begin{aligned}  &\exists (t_1, t_2) \text{obeying \eqref{Tube1 and Tube2} and \eqref{t-t1t2}} \\  &\text{such that } t_i \subseteq Q_i \text{ for } i = 1, 2 \end{aligned} \right \} \\
\leq \#\{Q \in \mathcal V_{k+1}(\mathscr{U}_{\mathbb X}) : Q \subsetneq t \} = M^{\lambda(u') - \lambda(u)}.
\end{equation}  
Given a pair $(Q_1, Q_2) \in \widetilde{\tt{TP}}$ and in view of \eqref{t-t1t2} and \eqref{dist-t1t2}, an interval tuple $(t_1, t_2) \in {\tt{TP}}_{\mathbb X}[\ell, k, u, t; r]$ descended from this pair has to satisfy the criteria 
\begin{align*}  
&t_i \subseteq I_i \subseteq Q_i, \text{ where $I_i$ is the unique interval with } \\ 
&|I_i| = C_1 M^{-r - \lambda(u) - \ell}, \; \text{dist}(I_i, \partial Q_i) \leq C_1 M^{-r - \lambda(u) - \ell}. 
\end{align*}  
Let us recall that $|t_i| = M^{-J}$, and therefore  
\begin{equation} \label{final descendant count} 
\# \left\{ q \in \mathcal Q(J) : q \subseteq I_i \right\}  \leq C_1 M^{J - r - \lambda(u) - \ell}.
\end{equation}  
Putting \eqref{TP-tilde count} and \eqref{final descendant count} together, we arrive at  
\begin{align*} 
\# \left[ {\tt{TP}}_{\mathbb X}[\ell, k, u, t; r]\right] &\leq \# \bigl[  \widetilde{\tt{TP}} \bigr] \times \# \left\{ q \in \mathcal Q(J) : q \subseteq I_i \right\}  \\ 
&\leq M^{\lambda(u') - \lambda(u)} \times \left( M^{-r - \lambda(u) - \ell + J}\right)^2 = M^{2J - 2r - 3 \lambda(u) + \lambda(u')}. 
\end{align*}  
This gives the desired estimate \eqref{TP-counting-estimate}, completing the proof. 
\end{proof}

\section{Proof of Proposition \ref{PAIRWISE PROP}} \label{section: pairwise prop proof} 
\begin{proof} 
To prove \eqref{pairwise upper bound}, we appeal to the splitting \eqref{subsum-1} of ${\tt{SumTP}}_{\mathbb X}[r]$ into sub-sums ${\tt{SumTP}}_{\mathbb X}[k, u, t; r]$, then estimate each summand using  \eqref{subsum-estimate} from Lemma \ref{subsum counting lemma}. This leads to the inequality
\begin{align*}
&{\tt{SumTP}}_{\mathbb X}[r] \\
&\quad = \sum \Bigl\{ {\tt{SumTP}}_{\mathbb X}[k, u, t; r]\, \Bigl| \, t \in \mathcal V_k(\mathscr{U}_{\mathbb X}), \, u \in {\tt{SplitV}}_k, \, \sigma_{\mathbb X}(t) \subsetneq u, \, 1 \leq k \leq N   \Bigr\} \\ &\quad \leq C_1 M^{-2r}\sum_{k=1}^{N} \sum_{u \in {\tt{SplitV}}_k} \sum_{t} \Bigl\{ M^{-\lambda(u)} \, \Bigl| \, t \in \mathcal V_k(\mathscr{U}_{\mathbb X}), \, \sigma_{\mathbb X}(t) \subsetneq u \Bigr\}  \\ 
& \quad \leq C_1 M^{-2r} \sum_{k=1}^{N} \sum_{t \in \mathcal V_k(\mathscr{U}_{\mathbb X})} |t| 
%& \leq C_1 M^{-2r} \sum_{k=1}^{N} \sum_{u \in {\tt{SplitV}}_k} M^{- \lambda (u)} \# \left\{t \in {\tt{BasicSl}}_k : |t| = M^{- \lambda(u)} \right\} \\
\leq  C_1 M^{-2r} \sum_{k=1}^{N} 1 \leq C_1N M^{-2r}. 
\end{align*}
Each $\sigma_{\mathbb X}(t)$ has a unique youngest splitting ancestor $u\in {\tt{SplitV}}_k$, and $|t|=M^{-\lambda(u)}$.
The sum over $t$ equals 1 by \eqref{cover-Q sum}. This completes the proof of the desired inequality \eqref{pairwise upper bound}.
\end{proof}

\chapter{Random construction of Kakeya-type sets} \label{random construction section}	
\section{Chapter overview} 
The goal of this chapter is to complete the proof of Proposition \ref{prop: away from the roots} by constructing a random Kakeya-type set $\mathtt K_{\mathbb X}$ whose measure is small away from the root line. This is achieved by assigning a probability measure to the collection of binary sequences of the form $\mathbb X$  described in \eqref{binary X}. Each realization $\mathbb X$ determines a slope assignment $t \mapsto \sigma_{\mathbb X}(t)$ to the family of root intervals $t$, resulting in a random family of tubes $\mathbb T_{\mathbb X}$ given by \eqref{tube family X}. This random assignment is guided by the combinatorial structure of the pruned slope tree. Thus the randomness in $\mathbb X$ induces randomness in the slope assignment; this in turn influences the geometry of the random set $\mathtt K_{\mathbb X}$, given by \eqref{defn: extended Kakeya set K} as the union of the tube configuration $\mathbb T_{\mathbb X}$.
\vskip0.1in
\noindent The central idea is to randomize the selection of slopes in a manner that preserves the inductive stickiness properties introduced in Chapter \ref{chapter: compressions}, while simultaneously promoting transverse intersections among tubes away from the root line. We show that such intersections create significant overlap, resulting in a small expected measure of $\mathtt K_{\mathbb X}$ in those regions. 
\vskip0.1in
\noindent The proof of Proposition \ref{prop: away from the roots} proceeds through a three-step reduction. First, we translate the problem into a succession of two estimates; one combinatorial, the other probabilistic. Assuming these two estimates, the proof is completed in this chapter. The probabilistic estimate is an application of known results; this is included here as well. The combinatorial cum geometric estimate that bridges the probabilistic estimate to the conclusion of Proposition \ref{prop: away from the roots} will be treated in the next chapter.   
\vskip0.1in
\begin{itemize} 
\item In Section \ref{section: av small away from root}, we state Proposition \ref{prop: probability x away from root}, the main result of this chapter. Specifically, it says that for a random choice of slope assignments, the likelihood of a Kakeya-type set to include a fixed point away from the root line is small. A consequence of this is a small expected measure of the set $\mathtt K_{\mathbb X}$ in regions far from the root line. This reduces the problem of measuring the size of $\mathtt K_{\mathbb X}$ to a pointwise probabilistic estimate, namely bounding $\mathbb P(x \in \mathtt K_{\mathbb X})$. Section \ref{proof section: prop away from the roots} contains the proof of Proposition \ref{prop: away from the roots}, assuming Proposition \ref{prop: probability x away from root}.
\vskip0.1in
\noindent The remainder of the chapter charts a proof strategy for Proposition \ref{prop: probability x away from root}. The proof rests on two statements, Propositions \ref{PROP: RELATED TO SURVIVAL} and \ref{prop: computing survival probability}. Assuming these two statements, the proof of Proposition \ref{prop: probability x away from root} is completed in Section \ref{section: completion of probability proof}.  
\vskip0.1in
\item In Proposition \ref{PROP: RELATED TO SURVIVAL} of Section \ref{section: percolation intro}, we frame the event of inclusion of a point $x$ in $\mathtt K_{\mathbb X}$ as a Bernoulli percolation process on a new tree called ${\tt{TrPR}}_x$, the tree of possible roots associated with $x$. This  is the key structural step: it converts the geometric inclusion event $x \in \mathtt K_{\mathbb X}$ into a combinatorial event on a tree dictated by $x$. The precise description of ${\tt{TrPR}}_x$ relies on the compressed structure of the pruned slope tree. We defer this description and the proof of Proposition \ref{PROP: RELATED TO SURVIVAL}  to Chapter \ref{chapter: far from root}.  
\vskip0.1in
\item 
%{\color{red} I think some parts of this remark overlaps with the previous remark.}
The combinatorial reformulation mentioned above allows us to interpret the probability $\mathbb P(x \in \mathtt K_{\mathbb X})$ in terms of survival of paths in the new tree ${\tt{TrPR}}_x$ after percolation. There is an extensive literature devoted to the computation of survival probabilities under percolation, and many of the techniques developed there are now available as a result of the reduction. Following this, Proposition \ref{prop: computing survival probability} provides a quantitative estimate for the necessary probability. 
%and estimate it using the extensive literature on the latter subject. 
\vskip0.1in 
\noindent In Section \ref{section: electrical networks}, we estimate the survival probability using an electrical network analogy. This follows the work of Lyons \cite{{Lyons1}, {Lyons2}} relating the survival probability of a percolation process to the effective resistance of an associated electrical network on the tree. This yields quantitative bounds sufficient to establish Proposition \ref{prop: computing survival probability}. 
\end{itemize} 
\vskip0.1in
\noindent To summarize, Proposition \ref{prop: probability x away from root} directly implies Proposition \ref{prop: away from the roots}; see Section \ref{proof section: prop away from the roots}. The proof of Proposition \ref{prop: probability x away from root} is itself reduced to two intermediate results: Proposition \ref{PROP: RELATED TO SURVIVAL}, established in Chapter \ref{chapter: far from root}, and Proposition \ref{prop: computing survival probability}, proved in Section \ref{section: electrical networks}.

\section{Average smallness of $\mathtt K_{\mathbb X}$ far from the root line} \label{section: av small away from root} 
\subsection{Likelihood of a distant point to lie in $\mathtt K_{\mathbb X}$} 
\label{section: probabilistic setup} 
%Having completed the proof of Proposition \ref{prop: near the roots}, which captures the behaviour of $\mathtt K_{\mathbb X}$ near the root line, we now turn our attention to proving Proposition \ref{prop: away from the roots}. This is the statement that ensures the existence of {\em{some}} root-to-slope map $\sigma_{\mathbb X}$ for which $\mathtt K_{\mathbb X}$ is small away from the root line, in the sense of \eqref{estimate: away from the roots}. 
%\vskip0.1in 
%\noindent We identify the desired map $\sigma_{\mathbb X}$ through a probabilistic mechanism. Precisely, 
Let us describe the probabilistic framework underlying the construction. We consider the probability measure space given by the Cartesian product 
 \[ \prod_{Q \in \mathcal Q^{\ast}} \mathfrak X_Q. \quad \text{ with } \quad \mathfrak X_{Q} = \{0, 1\} \text{ for each } Q \in \mathcal Q^{\ast}.  \] 
The probability measure on each binary space $\mathfrak X_Q$ assigns equal value to $0$ and $1$, i.e., each $\mathfrak X_Q$ is equipped with a random variable $X(Q)$ such that 
\[ \mathbb P(X(Q) = 0) = \mathbb P(X(Q)=1) = \frac{1}{2}. \]
In this measure space, the binary sequence $\mathbb X = \{X(Q) : Q \in \mathcal Q^{\ast}\}$ given by \eqref{binary X} is an independent and identically distributed (i.i.d) sequence of unbiased Bernoulli random variables. Every realization of $\mathbb X$ in this measure space then determines a slope assignment $t \mapsto \sigma_{\mathbb X}(t)$ as explained in Proposition \ref{PROP: COMPRESSED ROOTS}. This, in turn, generates a set $\mathtt K_{\mathbb X}$ according to the prescription \eqref{defn: extended Kakeya set K}. Thus $\mathtt K_{\mathbb X}$ becomes a random subset of the plane. For each point $x \in \mathbb R^2$, the membership event $\{x \in \mathtt K_{\mathbb X} \}$ is a measurable set in the underlying probability space.  
\vskip0.1in
\noindent To quantify the size of $\mathtt K_{\mathbb X}$, it is therefore natural to study the probability that a given point $x$ lies in $\mathtt K_{\mathbb X}$. This leads to the consideration of the quantities 
\[ \mathbb P\bigl( x \in \mathbb K_{\mathbb X}\bigr) \]
for points $x$ lying far from the root line. The main result of this chapter provides a uniform (in $x$) bound for these probabilities. 
%The main result of this section, Proposition \ref{prop: probability x away from root} below, asserts that for any point $x$ far away from the root line, the probability of $x$ lying in $\mathtt K_{\mathbb X}$ is small. This is a key step in establishing that, on average away from the root line, $\mathtt K_{\mathbb X}$ is small, i.e., displays Kakeya-like behaviour.   
\begin{proposition} \label{prop: probability x away from root}
Let $A_0 > 1$ be the absolute constant used to define ${\tt{Tube}}$ as in \eqref{defn: tube}. Then there is a constant  $A_1 > 0$ depending only on $A_0$ with the following property. 
\vskip0.1in
\noindent If $\mathbb X$ given by \eqref{binary X} is an i.i.d. sequence of Bernoulli$(\frac{1}{2})$ random variables, then  the random set $\mathtt K_{\mathbb X}$ given by \eqref{defn: extended Kakeya set K}, whose tubes have slopes in $\Omega_N$, obeys the estimate 
\begin{equation} \label{probability x away from root} 
\mathbb P \left( x \in \mathtt K_{\mathbb X}\right) \leq \frac{A_1}{N} \quad \text{ for all } x \in [A_0, A_0 +1] \times \mathbb R. 
\end{equation} 
\end{proposition} 
\noindent Proposition \ref{prop: probability x away from root} is proved in Section \ref{section: completion of probability proof}, assuming two intermediate results that will be addressed later in this chapter and the next. Proposition \ref{prop: away from the roots} is a direct consequence of Proposition \ref{prop: probability x away from root}. We complete this in the next section.    

\subsection{Proof of Proposition \ref{prop: away from the roots}, assuming Proposition \ref{prop: probability x away from root}} \label{proof section: prop away from the roots} 
\begin{proof} 
The size of $\mathtt K_{\mathbb X}$ far from the root is given by 
\begin{align}
|\mathtt K_{\mathbb X} \cap [A_0, A_0+1] \times \mathbb R| &= \int_{x_2 \in \mathbb R} \int_{x_1=A_0}^{A_0+1} \mathbf 1_{\mathtt K_{\mathbb X}}(x) \, dx \nonumber \\ 
&= \int_{x_2=0}^{A_0+2} \int_{x_1=A_0}^{A_0+1} \mathbf 1_{\mathtt K_{\mathbb X}}(x) \, dx. \label{pre-expectation} 
\end{align} 
In the display above, $\mathbf 1_E(x)$ denotes the indicator function of a set $E$:
\[ \mathbf 1_E(x) := \begin{cases} 1 &\text{ if } x \in E, \\ 0 &\text { otherwise.}\end{cases} \]  The last step \eqref{pre-expectation} uses the fact that $\Omega_{N} \subseteq [0,1]$, in the following way: a tube rooted on the vertical line segment $\{0\} \times [0,1]$ with slope in $[0,1]$ cannot extend outside the region $0 \leq x_2 \leq A_0+2$ on the vertical strip $A_0 \leq x_1 \leq A_0+1$. The shaded region in Figure \ref{fig: trapezoid} shows portion of the plane where the tubes away from the root line are geometrically confined.  
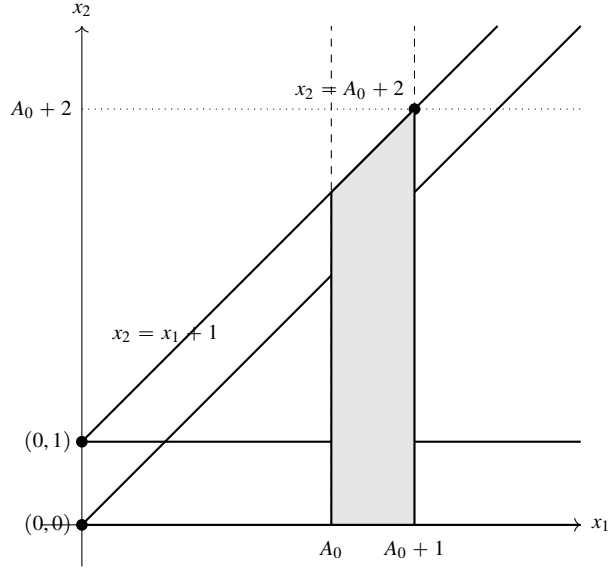
\begin{figure}[ht]
\centering
\begin{tikzpicture}[scale=1.1]

% Parameters
\def\A{3}

% Axes
\draw[->] (-0.5,0) -- (6,0) node[right] {\scriptsize{$x_1$}};
\draw[->] (0,-0.5) -- (0,6) node[above] {\scriptsize{$x_2$}};

% Points on y-axis
\fill (0,0) circle (2pt) node[left] {\scriptsize{$(0,0)$}};
\fill (0,1) circle (2pt) node[left] {\scriptsize{$(0,1)$}};

% Horizontal rays from (0,0) and (0,1)
\draw[thick] (0,0) -- (6,0);
\draw[thick] (0,1) -- (6,1);

% Slope 1 rays
\draw[thick] (0,0) -- (6,6);
\draw[thick] (0,1) -- (5,6);

% Vertical lines x = A0 and x = A0+1
\draw[dashed] (\A,0) -- (\A,6);
\draw[dashed] (\A+1,0) -- (\A+1,6);

\node at (\A,-0.3) {\scriptsize{$A_0$}};
\node at (\A+1,-0.3) {\scriptsize{$A_0+1$}};

% Key point (A0+1, A0+2)
\fill (\A+1,\A+2) circle (2pt);
\node[above left] at (\A+1,\A+2) {\scriptsize{$x_2 = A_0+2$}};

% Dotted horizontal line y = A0+2
\draw[dotted] (0,\A+2) -- (6,\A+2);
\node[left] at (0,\A+2) {\scriptsize{$A_0+2$}};

% Shaded trapezoid region
\fill[gray!20]
(\A,0) --
(\A,\A+1) --
(\A+1,\A+2) --
(\A+1,0) -- cycle;

% Boundary of trapezoid
\draw[thick]
(\A,0) --
(\A,\A+1) --
(\A+1,\A+2) --
(\A+1,0) -- cycle;

% Label for slope-1 boundary
\node at (1,2.3) {\scriptsize{$x_2=x_1+1$}};

\end{tikzpicture}

\caption{\small{The set $\mathtt K_{\mathbb X}$ is restricted to lie in the shaded region.}}
\label{fig: trapezoid} 
\end{figure}
In other words, 
\[ \mathtt K_{\mathbb X} \cap \Bigl[[A_0, A_0+1] \times \mathbb R \Bigr] \subseteq [A_0, A_0+1] \times [0, A_0 + 2]. \] 
We now take expected value of both sides of the equation \eqref{pre-expectation} with respect to $\mathbb X$. By linearity of expectation, 
 \begin{align} 
\mathbb E_{\mathbb X} \bigl[ |\mathtt K_{\mathbb X} \cap [A_0, A_0+1] \times \mathbb R| \bigr] &= \int_{x_2=0}^{A_0+2} \int_{x_1=A_0}^{A_0+1} \mathbb E_{\mathbb X} \left[ \mathbf 1_{\mathtt K_{\mathbb X}}(x) \right] \, dx \nonumber \\ 
&= \int_{x_2=0}^{A_0+2} \int_{x_1=A_0}^{A_0+1} \mathbb P \left( x \in \mathtt K_{\mathbb X}\right) \, dx \nonumber \\ 
&\leq \frac{A_1(A_0+2)}{N}.  \label{bound on average}
\end{align} 
The last step in the display above uses the probability estimate \eqref{probability x away from root} from Proposition \ref{prop: probability x away from root}. The upper bound \eqref{bound on average} on the expected value ensures the existence of a realization $\mathbb X$ for which the same estimate holds. The desired claim \eqref{estimate: away from the roots} holds for this choice of $\mathbb X$, with $C_2 = A_1(A_0+2)$. 
\end{proof} 
\vskip0.1in
\noindent The remainder of this chapter and the next are devoted to the proof of Proposition \ref{prop: probability x away from root}. 
%\noindent To sum up the situation thus far, we have reduced the proof of Proposition \ref{prop: away from the roots} to that of Proposition \ref{prop: probability x away from root}. In this chapter and the next, we focus on justifying the probabilistic estimate \eqref{probability x away from root}. The argument proceeds in two steps:
%\vskip0.1in 
%\begin{itemize} 
%\item We first relate the probability $\mathbb P(x \in \mathtt K_{\mathbb X})$ to the survival probability of a Bernoulli percolation on a certain deterministic tree called ${\tt{TrPR}}_x$, the tree of possible roots that can reach $x$. This step is summarized in Proposition \ref{PROP: RELATED TO SURVIVAL} below. 
%\vskip0.1in
%\item Next, we appeal to the structure of ${\tt{TrPR}}_x$ and known results from probability theory to compute the survival probability. The relevant statement appears in Proposition \ref{prop: computing survival probability}. 
%\end{itemize} 
\section{Bernoulli percolation and probability of survival} \label{section: percolation intro} 
Let us recall the definition of a standard Bernoulli$(p)$ percolation on a rooted, labelled tree $\mathscr{T}$ of finite height. 
\vskip0.1in
\begin{itemize}
\item Each edge of the tree is declared ``retained'' with probability $p$ and ``removed'' with probability $1-p$. 
\vskip0.1in
\item These choices are made independently for every edge. 
\end{itemize} 
\vskip0.1in
\noindent A {\em{path}} in $\mathscr{T}$ under the percolation model is a connected sequence of retained edges. The {\em{survival of the tree $\mathscr{T}$}} under percolation refers to the event that a path of maximal length exists, i.e. some ray of $\partial \mathscr{T}$ remains in the random graph after the percolation process. The {\em{survival probability}} $\varrho(\mathscr{T}; p)$ of $\mathscr{T}$ under a standard Bernoulli$(p)$ percolation is defined to be 
\begin{equation} 
\varrho(\mathscr{T}; p) := \left\{ \begin{aligned} &{\text{the probability of existence of a path starting from}} \\ &{\text{the root and ending at a terminal vertex of $\mathscr{T}$.}}
\end{aligned} 
\right\}
\end{equation}  
Figure \ref{fig: percolation} gives a schematic illustration of a percolation on an arbitrary tree $\mathscr{T}$. 
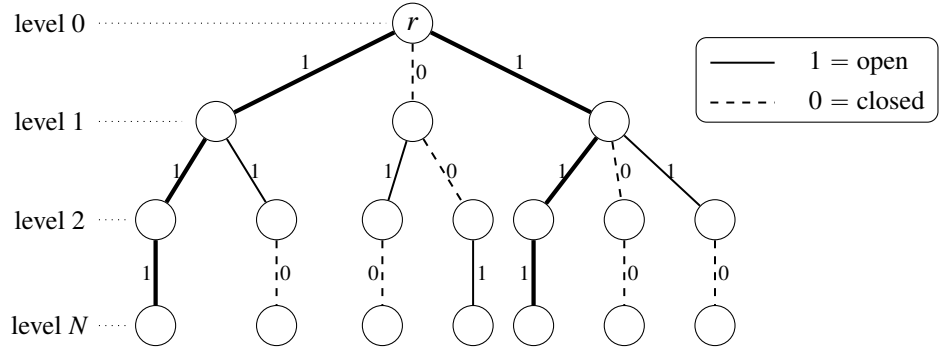
\begin{figure}[ht]
\centering
\begin{tikzpicture}[
    vertex/.style={circle, draw, minimum size=5.3mm, inner sep=0pt},
    open/.style={black, line width=0.75pt},
    closed/.style={black, dashed, line width=0.75pt},
    openpath/.style={black, line width=1.6pt},
    lab/.style={font=\scriptsize, inner sep=1pt},
    level/.style={font=\small},
    legend/.style={draw, rectangle, rounded corners, inner sep=5pt, font=\small}
]

% Vertices
\node[vertex] (r) at (0,0) {$r$};

\node[vertex] (a) at (-2.6,-1.3) {};
\node[vertex] (b) at (0,-1.3) {};
\node[vertex] (c) at (2.6,-1.3) {};

\node[vertex] (a1) at (-3.4,-2.6) {};
\node[vertex] (a2) at (-1.8,-2.6) {};

\node[vertex] (b1) at (-0.4,-2.6) {};
\node[vertex] (b2) at (0.8,-2.6) {};

\node[vertex] (c1) at (1.6,-2.6) {};
\node[vertex] (c2) at (2.8,-2.6) {};
\node[vertex] (c3) at (4.0,-2.6) {};

\node[vertex] (a11) at (-3.4,-4.0) {};
\node[vertex] (a21) at (-1.8,-4.0) {};

\node[vertex] (b11) at (-0.4,-4.0) {};
\node[vertex] (b21) at (0.8,-4.0) {};

\node[vertex] (c11) at (1.6,-4.0) {};
\node[vertex] (c21) at (2.8,-4.0) {};
\node[vertex] (c31) at (4.0,-4.0) {};

% Root to level 1
\draw[open]   (r) -- (a) node[midway, above left, lab] {$1$};
\draw[closed] (r) -- (b) node[midway, right, lab] {$0$};
\draw[open]   (r) -- (c) node[midway, above right, lab] {$1$};

% Level 1 to level 2
\draw[open]   (a) -- (a1) node[midway, left, lab] {$1$};
\draw[open]   (a) -- (a2) node[midway, right, lab] {$1$};

\draw[open]   (b) -- (b1) node[midway, left, lab] {$1$};
\draw[closed] (b) -- (b2) node[midway, right, lab] {$0$};

\draw[open]   (c) -- (c1) node[midway, left, lab] {$1$};
\draw[closed] (c) -- (c2) node[midway, right, lab] {$0$};
\draw[open]   (c) -- (c3) node[midway, right, lab] {$1$};

% Level 2 to level N
\draw[open]   (a1) -- (a11) node[midway, left, lab] {$1$};
\draw[closed] (a2) -- (a21) node[midway, right, lab] {$0$};

\draw[closed] (b1) -- (b11) node[midway, left, lab] {$0$};
\draw[open]   (b2) -- (b21) node[midway, right, lab] {$1$};

\draw[open]   (c1) -- (c11) node[midway, left, lab] {$1$};
\draw[closed] (c2) -- (c21) node[midway, right, lab] {$0$};
\draw[closed] (c3) -- (c31) node[midway, right, lab] {$0$};

% Emphasize surviving open paths only
\draw[openpath] (r) -- (a) -- (a1) -- (a11);
\draw[openpath] (r) -- (c) -- (c1) -- (c11);

% Level labels
\node[level] at (-4.8,0) {level $0$};
\node[level] at (-4.8,-1.3) {level $1$};
\node[level] at (-4.8,-2.6) {level $2$};
\node[level] at (-4.8,-4.0) {level $N$};

% Dotted guides
\draw[dotted] (-4.15,0) -- (-0.35,0);
\draw[dotted] (-4.15,-1.3) -- (-3.05,-1.3);
\draw[dotted] (-4.15,-2.6) -- (-3.75,-2.6);
\draw[dotted] (-4.15,-4.0) -- (-3.75,-4.0);

% Legend
\node[legend, align=left] at (5.35,-0.75) {
\begin{tikzpicture}[baseline=-0.5ex]
\draw[open] (0,0) -- (0.85,0);
\end{tikzpicture}
\quad $1=$ open\\[3pt]
\begin{tikzpicture}[baseline=-0.5ex]
\draw[closed] (0,0) -- (0.85,0);
\end{tikzpicture}
\quad $0=$ closed
};

\end{tikzpicture}

\caption{\small{Bernoulli percolation on an arbitrary tree $\mathscr{T}$, where different vertices may have different numbers of children. Edges are independently retained (marked 1) with probability $p$ and removed (marked 0) with probability $(1-p)$. The specific realization of percolation in the diagram has two surviving paths of maximal length.}
\label{fig: percolation}}
\end{figure}

\subsection{The tree of possible roots ${\tt{TrPR}}_x$}
The following proposition is the main pillar of the proof of Proposition \ref{prop: probability x away from root}. It identifies an auxiliary, deterministic, compressed tree ${\tt{TrPR}}_x$ of roots (the ``tree of possible roots'') whose survival under an unbiased Bernoulli percolation is equivalent to the inclusion of $x$ in $\mathbb K_{\mathbb X}$.
\begin{proposition} \label{PROP: RELATED TO SURVIVAL}
Assume the set-up of Proposition \ref{prop: probability x away from root}. Let $x \in [A_0, A_0+1] \times \mathbb R$ be a point that is reachable by a tube with slope from $\Omega_N$. In other words, 
\begin{equation} \label{foreshadowing Possx}  
\exists t \in \mathcal Q(J), \omega \in \Omega_N \text{ such that } x \in {\tt{Tube}}\bigl(t, \omega, [A_0, A_0+1] \bigr).  
\end{equation}   
Then there exists a deterministic tree ${\tt{TrPR}}_x$ with the following properties:
\vskip0.1in 
\begin{enumerate}[(a)]
\item \label{survival-part1} The tree ${\tt{TrPR}}_x$ has height $N$; in fact, every maximal ray of ${\tt{TrPR}}_x$ consists of exactly $N$ edges. 
\vskip0.1in
\item \label{survival-part2} Every non-terminal vertex of ${\tt{TrPR}}_x$ has at most $C_1$ children, where $C_1 \geq 2$ is an integer constant depending only on $A_0$ and uniform in $x$. Moreover, there are at most $C_1 2^j$ vertices of ${\tt{TrPR}}_x$ at height $j$, $1 \leq j \leq N$, i.e., 
\begin{equation} \label{vertex count in TrPR} 
\# \bigl[ \mathcal V_j({\tt{TrPR}}_x) \bigr] \leq C_1 2^j. 
\end{equation}  
\vskip0.1in
\item \label{survival-part3} The probability occurring on the left hand side of \eqref{probability x away from root}, i.e., the probability of the random set $\mathtt K_{\mathbb X}$ containing $x$ equals the survival probability of a Bernoulli$(\frac{1}{2})$ percolation on ${\tt{TrPR}}_x$:
\begin{equation} \label{prob identity}
\mathbb P(x \in \mathtt K_{\mathbb X}) = \varrho \bigl({\tt{TrPR}}_x; \frac{1}{2}\bigr). 
\end{equation} 
\end{enumerate}
\end{proposition} 
\vskip0.1in
\noindent {\em{Remarks: }} 
\begin{enumerate}[1.] 
\item If \eqref{foreshadowing Possx} does not hold, then $x \in \mathtt K_{\mathbb X}$ is an impossible event, and 
\[ \mathbb P(x \in \mathtt K_{\mathbb X}) = 0. \] 
\vskip0.1in
\item In Chapter \ref{chapter: far from root}, we give an explicit description of the tree ${\tt{TrPR}}_x$ and list its key properties. See Section \ref{section: reference trees}. 
\vskip0.1in 
\item Based on this information, Proposition \ref{PROP: RELATED TO SURVIVAL} is proved in Section \ref{survival prop proof}. 
\vskip0.1in 
\item The features of  ${\tt{TrPR}}_x$ important for the proof are established in Chapter \ref{chapter: reference tree geometry}.
\end{enumerate} 
\subsection{An estimate for the survival probability} 
Having translated $\mathbb P(x \in \mathtt K_{\mathbb X})$ to the survival probability of the tree ${\tt{TrPR}}_x$ under a Bernoulli percolation, we are tasked with estimating the latter quantity, namely the right hand side of \eqref{prob identity}. The following proposition provides the required estimate. 
\begin{proposition} \label{prop: computing survival probability} 
Let ${\tt{TrPR}}_x$ be the tree described in Proposition \ref{PROP: RELATED TO SURVIVAL}. Then the survival probability for a Bernoulli$(\frac{1}{2})$ percolation on this tree obeys the estimate:
\begin{equation} \label{upper bound on survival prob} 
\varrho \bigl({\tt{TrPR}}_x; \frac{1}{2}\bigr) \leq \frac{A_1}{N}. 
\end{equation} 
\end{proposition} 
\noindent Let us take a moment to verify that Proposition \ref{prop: probability x away from root} follow directly from Propositions \ref{PROP: RELATED TO SURVIVAL} and \ref{prop: computing survival probability}.
\subsection{Proof of Proposition \ref{prop: probability x away from root}, assuming Propositions \ref{PROP: RELATED TO SURVIVAL} and \ref{prop: computing survival probability}} \label{section: completion of probability proof}
%{\color{red} I think we don't need this.}
The auxiliary tree ${\tt{TrPR}}_x$ guaranteed by Proposition \ref{PROP: RELATED TO SURVIVAL} is the natural bridge between the left and right sides of the desired estimate \eqref{probability x away from root}. Indeed, combining \eqref{prob identity} and \eqref{upper bound on survival prob} leads to the desired estimate: 
\[ \mathbb P(x \in \mathtt K_{\mathbb X}) = \varrho \bigl( {\tt{TrPR}}_x; \frac{1}{2}\bigr) \leq \frac{A_1}{N},\]
completing the proof. \qed 

\section{Estimation of survival probability via electrical networks} \label{section: electrical networks} 
The goal of this section is to prove Proposition \ref{prop: computing survival probability}. For this, we will invoke a result of Lyons \cite{Lyons1}. Let us describe the set-up that leads to this result.  The interested reader may consult~\cite{Grimmett} for a discussion of percolation processes in much greater generality.
\vskip0.1in
\noindent Given a tree $\mathscr{T}$ with an edge set $\mathcal{E}$, we define an \textit{edge-dependent Bernoulli (bond) percolation process} to be a collection of independent random variables $\{Y_e : e\in\mathcal{E}\}$, where $Y_e$ is Bernoulli$(p_e)$ with $p_e<1$.  For a given edge $e\in\mathcal{E}$, the event $\{Y_e=0\}$ indicates \textit{removal} of the edge $e$ from the edge set $\mathcal{E}$ during percolation, and the event $\{Y_e=1\}$ signifies \textit{retention} of this edge.  Thus, 
\[ p_e := \text{Pr}(Y_e=1)  \text{ is the retention probability for $e$}. \] 
We define the family of edge retention probabilities by 
\[ \mathscr{P}(\mathscr{T}) := \bigl\{p_e : e \in \mathcal E \bigr\}. \]  As mentioned in Section \ref{section: percolation intro}, survival of the tree $\mathscr{T}$ is defined to be the event that at least one ray remains from the root of the tree to its bottommost level. The probability of this event is referred to as the \textit{survival probability} of the corresponding percolation process, denoted $\varrho(\mathscr{T}; \mathscr{P}(\mathscr{T}))$.  If the random variables $\{Y_e : e\in\mathcal{E}\}$ are mutually independent and identically distributed Bernoulli$(p)$ random variables, with a constant $p < 1$ independent of the edge $e$, then the process is called a \textit{standard Bernoulli$(p)$ percolation}, and we denote the survival probability simply by $\varrho(\mathscr{T}; p)$.
\vskip0.1in
\noindent In \cite{Lyons1}, Lyons visualized percolation on a tree as a certain electrical network, in which the survival probability under percolation can be directly related to the effective resistance of the network.  The natural electrical network is defined as follows. Suppose that the tree $\mathscr{T}$ is of finite height. We place the positive node of a battery at the root of $\mathscr{T}$, and connect each of the terminating vertices to the negative node of the battery. A resistor is placed on every edge $e$ of $\mathscr{T}$ with resistance $\mathtt R_e$ defined by
\begin{equation}\label{resistance}
	\frac 1{\mathtt R_e} = \frac 1{1-p_e}\prod_{\begin{subarray}{c} e' \in \mathcal E \\ v(e)\subseteq v(e')\end{subarray}} p_{e'},
\end{equation}
where $v(e)$ is the vertex in $\mathscr{T}$ at which $e$ terminates.  Notice that the conductance $\mathtt C_e$ for the edge $e$ (which is by definition the reciprocal of the resistance $\mathtt R_e$) is essentially the probability that a path remains from the root of $\mathscr{T}$ to the vertex $v(e)$ after percolation. For a standard Bernoulli$(p)$ percolation, these edge resistances take on a particularly simple form:
\begin{equation} \label{resistances} 
\begin{aligned} 
&\frac{1}{\mathtt R_{e}} = \frac{1}{1 - p_e} p^{h(v(e))}, \; \text{ where } h(v(e)) = \text{ height of $v(e)$ in $\mathscr{T}$.}\\
&\text{For $p = \frac{1}{2}$, this reduces to } \mathtt R_e = 2^{h(v(e))-1}.  
\end{aligned} 
\end{equation}   
\vskip0.1in
\noindent The following result quantifies the relation between the survival probability under percolation and the effective resistance of a tree. 
\begin{letteredtheorem}{\cite[Theorem 2.1]{Lyons1}} \label{Lyons survival prob}
Let $\mathscr{T}$ be a tree equipped with a Bernoulli percolation process, and its associated electrical network, as described above. Then  
\begin{equation} \label{probability and resistance} 
\varrho\bigl(\mathscr{T}; \mathscr{P}(\mathscr{T}) \bigr) \leq \frac{2 \mathtt{C}(\mathscr{T})}{1 + \mathtt{C}(\mathscr{T})} = \frac{2}{1 + \mathtt R(\mathscr{T})},
\end{equation}
where $\mathtt{C}(\mathscr{T})$ and $\mathtt R(\mathscr{T})$ denote respectively the effective conductance and effective resistance of the electrical network on the tree $\mathscr{T}$.  
\end{letteredtheorem}
\subsection{Proof of Proposition \ref{prop: computing survival probability}} 
\begin{proof} 
In light of Theorem \ref{Lyons survival prob}, we see that
\[ \varrho \left( {\tt{TrPR}}_x; \frac{1}{2}\right) \leq \frac{2}{1 + \mathtt R({\tt{TrPR}}_x)}.\] 
Thus it is sufficient to bound the resistance of the electrical network ${\tt{TrPR}}_x$ from below. In its original form, the network ${\tt{TrPR}}_x$ is a complex system of resistors, some in series, others in parallel; however, the following convenient fact allows us to replace it by a simpler network, whose overall resistance is lower than that of ${\tt{TrPR}}_x$ and is more amenable to computation. This is achieved by connecting any two vertices at a given height by an ideal conductor. 
%is useful in estimating the effective resistance of this network by replacing it by one we need the useful fact that connecting any two vertices at a given height by an ideal conductor (i.e. one with zero resistance) only decreases the overall resistance of the circuit.  
\begin{lemma}[{\cite[Proposition 5.1]{KrocPramanik}}]\label{resistance prop}
	Let $\mathscr{T}$ be a tree of height $N$, each of whose maximal rays consists of $N$ edges. Consider the corresponding electrical network generated by a standard Bernoulli$(\frac 12)$ percolation process on $\mathscr{T}$.  For a fixed $k<N$, we connect two vertices at height $k$ by a conductor with zero resistance.  Then the total resistance of the resulting electrical network is no greater than that of the original.
\end{lemma}
%For a proof of this fact, see~\cite[Proposition 5.1]{KrocPramanik}.  
\noindent The main consequence of the lemma 
%that we draw upon  in Lemma~\ref{computing survival probability} 
is the following corollary.
\begin{corollary} \label{survival probability reduced}
	Given a tree $\mathscr{T}$ of height $N$ as in Lemma \ref{resistance prop}, let $\mathtt R(\mathscr{T})$ denote the total resistance of the electrical network that corresponds to the standard Bernoulli$(\frac{1}{2})$ percolation on this tree, in the sense of Theorem \ref{Lyons survival prob}.   Then
\begin{equation}\label{ResistBound}
	\mathtt R(\mathscr{T}) \geq \sum_{k=1}^N\frac {2^{k-1}}{n_k} \text{ where } n_k = \#\mathcal V_k(\mathscr{T})
\end{equation}
denotes the number of $k^{\text{th}}$ generation vertices of $\mathscr{T}$.
\end{corollary} 
\begin{proof}
To show \eqref{ResistBound}, we construct an auxiliary electrical network from the one naturally associated to our tree $\mathscr{T}$, as follows.  For every $k\geq 1$, we connect all vertices at height $k$ by an ideal conductor to make one node $\mathtt V_k$.  We call this new circuit $\mathscr{E}$.  The resistance of $\mathscr{E}$ cannot be greater than the resistance of the original circuit, by Lemma~\ref{resistance prop}. 
\vskip0.1in
\noindent Fix $k$, $1\leq k\leq N$, and let $\mathtt R_k$ denote the resistance in $\mathscr{E}$ between $\mathtt V_{k-1}$ and $\mathtt V_k$.  The number of edges between $\mathtt V_{k-1}$ and $\mathtt V_k$ is equal to the number $n_k$ of $k$th generation vertices in $\mathscr{T}$, and each edge is endowed with resistance $2^{k-1}$ by \eqref{resistances}.  Since these resistors are in parallel, we obtain 
\begin{equation}
	\frac 1{\mathtt R_k} = \sum_{1}^{n_k} \frac{1}{2^{k-1}} = \frac{n_k}{2^{k-1}}.\nonumber
\end{equation}
This holds for every $1\leq k\leq N$.  Since the resistors $\{\mathtt R_k : 1 \leq k \leq N \}$ are in series, 
\[ \mathtt R(\mathcal{T}_N) \geq \mathtt R(\mathscr{E}) = \sum_{k=1}^{N} \mathtt R_k = \sum_{k=1}^{N} \frac{2^{k-1}}{n_k}, \] establishing inequality \eqref{ResistBound}.
\end{proof}
\vskip0.1in
\noindent The proof of Proposition \ref{prop: computing survival probability} follows from Proposition \ref{PROP: RELATED TO SURVIVAL} and Corollary \ref{survival probability reduced}. According to part \eqref{survival-part2} of Proposition \ref{PROP: RELATED TO SURVIVAL}, 
\[ n_k = \# \left[ \mathcal V_k({\tt{TrPR}}_x) \right] \leq C_1 2^k.  \] 
Substituting this into \eqref{ResistBound}, we obtain a lower bound on the total resistance of ${\tt{TrPR}}_x$:
\[ \mathtt R({\tt{TrPR}}_x) \geq \sum_{k=1}^{N} \frac{2^{k-1}}{C_1 2^k} = \frac{N}{2C_1}.\] 
Finally, inserting this lower bound into \eqref{probability and resistance} yields
\[ \varrho \left({\tt{TrPR}}_x; \frac{1}{2}\right) \leq \frac{2}{1 + \frac{N}{2C_1}} \leq \frac{A_1}{N} \text{ for a suitably chosen $A_1 > 0$}, \] 
proving \eqref{upper bound on survival prob}. 
\end{proof} 

\chapter{Tubes far from the root line} \label{chapter: far from root} 
\section{Chapter overview} 
Chapter \ref{random construction section} reduced the construction of Kakeya-type sets to a single outstanding ingredient, namely Proposition \ref{PROP: RELATED TO SURVIVAL}. That proposition controls the behaviour of tubes far from the root line. The purpose of the present chapter, together with Chapter \ref{chapter: reference tree geometry}, is to establish this remaining step.

%In Chapter \ref{random construction section}, we presented a probabilistic construction of Kakeya-type configurations, modulo a key step, namely Proposition \ref{PROP: RELATED TO SURVIVAL}. This emerged as the final ingredient in this construction, controlling the behaviour of tubes far from the root line. In this chapter and the next, we complete the proof of Proposition \ref{PROP: RELATED TO SURVIVAL}. 
%This chapter isolates the key contributors, based on which Proposition \ref{PROP: RELATED TO SURVIVAL} is proved in Section \ref{survival prop proof}. The proofs of the key contributors are relegated to Chapter \ref{chapter: reference tree geometry} 
\vskip0.1in
\noindent For each point $x$ away from the root line, Proposition \ref{PROP: RELATED TO SURVIVAL} posits the existence of an auxiliary tree ${\tt{TrPR}}_x$. This tree is designed so that its survival under percolation signals the inclusion of $x$ in $\mathtt K_{\mathbb X}$. The central idea is to encode, for each point $x$, the set of slopes that can produce a tube passing through $x$, and to compare it with the random slope assignment. This is  framed through a percolation process on trees.
\vskip0.1in
\noindent After developing the necessary geometric preliminaries in Section \ref{section: Poss_x}, we construct this tree in Section \ref{section: reference trees}. Its essential structural properties are recorded in this section in a sequence of lemmas. Assuming these, we establish Proposition \ref{PROP: RELATED TO SURVIVAL} in Section \ref{survival prop proof}. The proofs of these lemmas are given later in Chapter \ref{chapter: reference tree geometry}.
\vskip0.1in
\noindent The first task is to determine which roots could possibly generate a tube passing through a fixed point $x$. This leads to the deterministic reference set ${\tt{Poss}}(x)$, consisting of all ``geometrically possible'' roots for $x$. The tree ${\tt{TrPR}}_x$ may be viewed as a compressed encoding of this set. Roughly speaking, for a given point $x$, the set ${\tt{Poss}}(x)$ consists of all possible roots such that a rectangle with this root and an orientation drawn from the pruned slope set $\Omega_N$ can pass through $x$.
\vskip0.1in
\noindent As we will see in Section \ref{section: Poss_x}, the construction of ${\tt{Poss}}(x)$ is naturally tied to a deterministic map $\tau_x$ that assigns to each root $t \in {\tt{Poss}}(x)$ a uniquely defined ``correct'' or reference slope $\tau_x(t) \in \Omega_N$ such that the corresponding tube contains $x$. For $x$ to lie in $\mathtt K_{\mathbb X}$, at least one of the constituent tubes ${\tt{Tube}}_{\mathbb X}[t]$ must contain $x$; in other words, the slope $\sigma_{\mathbb X}(t)$ of the tube ${\tt{Tube}}_{\mathbb X}[t]$ must agree with this reference slope $\tau_x(t)$. This leads to a matching condition between the deterministic slope $\tau_x(t)$ and the randomly assigned slope $\sigma_{\mathbb X}(t)$. In other words, there must exist
\begin{equation} \label{sigma tau match}
t \in {\tt{Poss}}(x) \text{ and } \tau_x(t) = \sigma_{\mathbb X}(t).
\end{equation}
In this way, $\tau_x(t)$ provides a reference against which the random assignments $\sigma_{\mathbb X}(t)$ are compared as $\mathbb X$ varies. Each assignment is deemed ``correct'' or ``incorrect'' for the inclusion event $x \in \mathtt K_{\mathbb X}$ depending on whether it matches the reference slope $\tau_x(t)$.
\vskip0.1in
\noindent Recall from Section \ref{section: compressions of the unit interval} that the slope allocation map $t \mapsto \sigma_{\mathbb X}(t)$ arises as the outcome of an iterative procedure based on an adapted tree-map pair $(\mathscr{U}_{\mathbb X}, \sigma_{\mathbb X})$. The binary sequence $\mathbb X$ governs the sequence of left and right turns along a path in the compressed tree $\mathscr{B}_N$, ultimately determining the value of $\sigma_{\mathbb X}(t)$. It is therefore natural and essential to reinterpret the matching criterion \eqref{sigma tau match} in the language of trees. Section \ref{section: reference trees} is devoted to making the connection precise.
\vskip0.1in
\noindent A priori, as seen in Sections \ref{section: ref set and ref map def} and \ref{section: ref map injectivity}, the set ${\tt{Poss}}(x) \subseteq \mathcal Q(J)$ and its image $\tau_x({\tt{Poss}}(x)) \subseteq \Omega_N$ are defined only at the level of roots and slopes. It is not clear, in general, that the reference map $\tau_x$ can be lifted to a sticky tree-map between trees representing ${\tt{Poss}}(x)$ and $\tau_x({\tt{Poss}}(x))$, in a manner that permits an edge-by-edge comparison with $\sigma_{\mathbb X}(t)$. Indeed, for arbitrary tree representations of these sets, such a lifting need not exist. The main objective of Lemma \ref{LEMMA: MAP LIFT} is to show that, for compressed trees whose generations are determined by the fundamental heights, this lifting of $\tau_x$ can in fact be achieved.
\vskip0.1in
\noindent To be specific, suppose that ${\tt{TrPS}}_x \subseteq \mathscr{B}_N$ is a compressed tree representation of $\tau_x({\tt{Poss}}(x))$. Lemmas \ref{CONSISTENCY LEMMA} and \ref{TREE SET LEMMA} construct compressed trees ${\tt{TrPR}}_x$ and ${\tt{TrPS}}_x$ encoding
${\tt{Poss}}(x)$ and its image under $\tau_x$. Lemma \ref{LEMMA: MAP LIFT} then lifts
\[ \tau_x : {\tt{TrPR}}_x \rightarrow {\tt{TrPS}}_x \text{ to a length-preserving sticky map. } \]  
This makes it possible to carry out an edge-by-edge comparison between the deterministic reference map $\tau_x$ with the random slope allocation map $\sigma_{\mathbb X}$. In turn, this leads naturally to a percolation process on ${\tt{TrPR}}_x$, in which an edge is retained precisely when $\tau_x$ and $\sigma_{\mathbb X}$ agree at its terminating vertex. Survival of a full ray in ${\tt{TrPR}}_x$ under this percolation is then equivalent to the inclusion $x \in \mathtt K_{\mathbb X}$. The remainder of the chapter develops this deterministic reference structure and shows how it interfaces with the random slope assignment.

\section{The reference set ${\tt{Poss}}(x)$ of possible slopes} \label{section: Poss_x}	
\noindent We now begin the deterministic part of the analysis. Throughout this section, we fix a point $x \in I_0 \times \mathbb R$, and ask a simple geometric question: which root intervals can generate a tube, with slope drawn from $\Omega_N$, that passes through $x$? 
\subsection{A criterion for tube intersections} 
The first ingredient in addressing this question is a geometric estimate describing when two tubes can intersect inside the strip $I_0 \times \mathbb R$. This estimate ultimately allows us to recover the slope of a tube from its root and a point lying inside it.
\begin{lemma} \label{intersection criterion lemma}
For slopes $\omega, \omega' \in \Omega_N$, roots $t, t' \in \mathcal Q(J)$,  and any interval $I \subseteq [0, 10A_0]$, let
\begin{equation} \label{tube and tube'}
{\tt{Tube}} = {\tt{Tube}}(t, \omega, I) \; \text{ and } \; {\tt{Tube}}' = {\tt{Tube}}(t', \omega', I), 
\end{equation} be any two parallelograms defined as in \eqref{defn: tube}. If there exists 
\begin{equation} 
x = (x_1, x_2) \in {\tt{Tube}} \cap {\tt{Tube}}', \label{x in intersection} 
\end{equation} then the inequality
 \begin{equation} \label{intersection criterion inequality} \bigl| \text{cen}(t') - \text{cen}(t) + x_1(\omega'-\omega) \bigr| \leq c_0 M^{-J} \end{equation} 
holds, where $\text{cen}(t)$ denotes the centre of the interval $t$. Here $c_0 \in (0,1)$ is the absolute constant used in \eqref{defn: tube} to define {\tt{Tube}}.   
\end{lemma}
\begin{proof}
Let $\tilde{Q}_t, \tilde{Q}_{t'}$ denote respectively the $c_0$-dilates of the intervals 
\[ t = Q_t \in \mathcal Q(J), \quad t' = Q_{t'} \in \mathcal Q_J. \] The assumption $x \in {\tt{Tube}} \cap {\tt{Tube}}'$ implies that there exist points $y \in \tilde{Q}_t$, $y' \in \tilde{Q}_{t'}$, and $\rho, \rho' \in [0, 10A_0]$ such that
\[ x = (x_1, x_2) = (0, y) + \rho (1, \omega) = (0, y') + \rho' (1, \omega'), \] Comparing coordinates in the vector-valued equation above, we obtain $\rho = \rho' = x_1$, leading to the relation 
\begin{equation} 
x_2 = y + x_1 \omega = y' + x_1 \omega', \; \text{ and therefore } \; x_1(\omega'- \omega) = y - y'. \label{intersection criterion inequality-1} 
\end{equation}  
The location of $y, y'$ within $t, t'$ allow us to write
\begin{align*} 
&y = \text{cen}(t) + z, \;  y' = \text{cen}(t') + z', \text{ with } \\
&2|z|, 2|z'| \leq \text{diam}(\tilde{Q}_t) = \text{diam}(\tilde{Q}_{t'}) = c_0 M^{-J}.
\end{align*}  
Substituting these into \eqref{intersection criterion inequality-1}, we obtain  
\[ \bigl| \text{cen}(t') - \text{cen}(t) + x_1(\omega'-\omega) \bigr| = |z - z'| \leq |z| + |z'| \leq c_0 M^{-J}, \] 
which is the claimed inequality (\ref{intersection criterion inequality}). Figure \ref{fig: two intersecting tubes} pictorially demonstrates the argument. This completes the proof.    
\end{proof} 
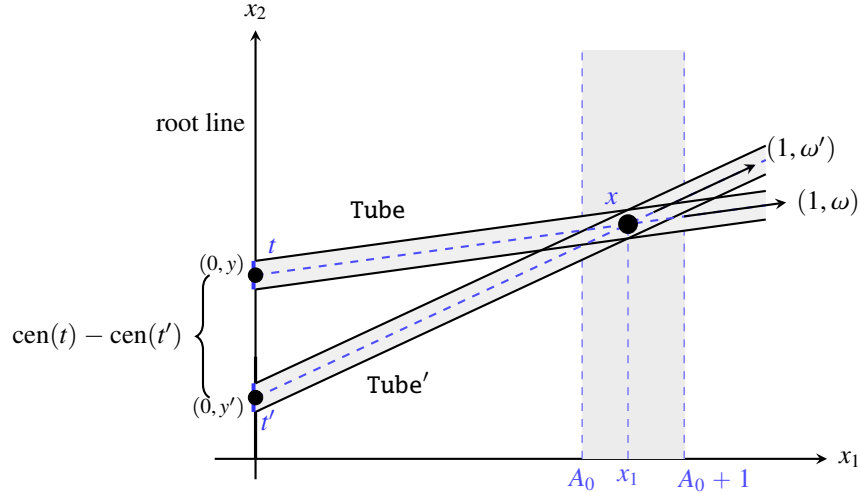
\begin{figure}[h]
\centering
\begin{tikzpicture}[scale=1.35, >=stealth, font=\small]

\def\A{3.2}
\def\B{4.2}
\def\xone{3.65}
\def\xtwo{2.3}

% Spread out y and y'
\def\y{1.8}
\def\yp{0.6}

% slopes so centre lines pass through x
\pgfmathsetmacro{\w}{(\xtwo-\y)/\xone}
\pgfmathsetmacro{\wp}{(\xtwo-\yp)/\xone}

% thicker tubes
\def\eps{0.14}
\def\L{5.0}

% Axes
\draw[->, thick] (-0.4,0) -- (5.6,0) node[right] {$x_1$};
\draw[->, thick] (0,-0.2) -- (0,4.2) node[above] {$x_2$};

% Root line
\draw[very thick] (0,0) -- (0,1);
\node[left] at (0,3.3) {root line};

% Shaded strip
\fill[gray!15] (\A,0) rectangle (\B,4.0);
\draw[dashed, blue!70] (\A,0) -- (\A,4.0);
\draw[dashed, blue!70] (\B,0) -- (\B,4.0);
\node[below, blue!70] at (\A,0) {$A_0$};
\node[below, blue!70] at (\B+0.3,0) {$A_0+1$};

% Root intervals (larger)
\draw[line width=3pt, blue!80] (0,\y-\eps) -- (0,\y+\eps);
\node[left, blue!80] at (0.3,\y+\eps+0.15) {$t$};

\draw[line width=3pt, blue!80] (0,\yp-\eps) -- (0,\yp+\eps);
\node[left, blue!80] at (0.3,\yp-\eps-0.08) {$t'$};

% Brace
\draw[decorate, decoration={brace, amplitude=5pt}, thick]
(-0.45,\yp) -- (-0.45,\y)
node[midway,left=6pt] {$\operatorname{cen}(t)-\operatorname{cen}(t')$};

% Tubes
\fill[gray!12]
  (0,\y-\eps) --
  (\L,{\y-\eps+\w*\L}) --
  (\L,{\y+\eps+\w*\L}) --
  (0,\y+\eps) -- cycle;

\fill[gray!12]
  (0,\yp-\eps) --
  (\L,{\yp-\eps+\wp*\L}) --
  (\L,{\yp+\eps+\wp*\L}) --
  (0,\yp+\eps) -- cycle;

\draw[thick] (0,\y-\eps) -- (\L,{\y-\eps+\w*\L});
\draw[thick] (0,\y+\eps) -- (\L,{\y+\eps+\w*\L});

\draw[thick] (0,\yp-\eps) -- (\L,{\yp-\eps+\wp*\L});
\draw[thick] (0,\yp+\eps) -- (\L,{\yp+\eps+\wp*\L});

% Centre lines
\draw[dashed, blue!70, thick] (0,\y) -- (\L,{\y+\w*\L});
\draw[dashed, blue!70, thick] (0,\yp) -- (\L,{\yp+\wp*\L});

% Points y and y'
\filldraw[black] (0,\y) circle (2pt);
\node[left=4pt] at (0.1,\y+0.1) {\scriptsize{$(0,y)$}};

\filldraw[black] (0,\yp) circle (2pt);
\node[left=4pt] at (0.1,\yp-0.1) {\scriptsize{$(0,y')$}};

% Labels (more breathing room)
\node[above, yshift=6pt] at (1.2,{\y+\w*1.2+\eps}) {${\tt Tube}$};
\node[below, yshift=-6pt] at (1.4,{\yp+\wp*1.4-\eps}) {${\tt Tube}'$};

% Direction arrows
\draw[->, thick]
  (4.2,{\y+\w*4.2}) -- (5.2,{\y+\w*5.2});
\node[right] at (5.2,{\y+\w*5.2}) {$(1,\omega)$};

\draw[->, thick]
  (3.9,{\yp+\wp*3.9}) -- (4.9,{\yp+\wp*4.9});
\node[right] at (4.9,{\yp+\wp*4.9+0.15}) {$(1,\omega')$};

% x and projection (last!)
\draw[dashed, blue!70] (\xone,0) -- (\xone,\xtwo);
\node[below, blue!70] at (\xone,0) {$x_1$};

\filldraw[black] (\xone,\xtwo) circle (2.6pt);
\node[above left, blue!80] at (\xone,\xtwo+0.1) {$x$};

\end{tikzpicture}
\caption{\small{Two tubes rooted at intervals $t,t'$ on the root line and intersecting at a point $x$ in the strip $[A_0,A_0+1]\times \mathbb R$.}}
\label{fig: two intersecting tubes}
\end{figure}
%\begin{corollary} \label{which is bigger corollary}
%If the positive constant $c_0$ in the definition \eqref{defn: tube} of a tube is chosen sufficiently small, then under the hypotheses \eqref{x in intersection} of distinct intersecting tubes in Lemma \ref{intersection criterion lemma}, 
%\begin{equation} \label{which is bigger} 
%|x_1||\omega-\omega'| \geq \frac{1}{2} M^{-J}. 
%\end{equation} 
%\end{corollary}
%\begin{proof}
%Since $t \ne t'$, we know that $|\text{cen}(t') - \text{cen}(t)| \geq M^{-J}$. The inequality in \eqref{intersection criterion inequality} therefore implies that 
%\[ |x_1| |\omega -\omega'| \geq |\text{cen}(t) - \text{cen}(t')| - c_0  M^{-J} \geq (1 - c_0) M^{-J} \geq \frac{1}{2} M^{-J},\]
%provided $c_0$ is chosen to satisfy $c_0 < \frac{1}{2}$.  
%\end{proof}  
\noindent The intersection estimate \eqref{intersection criterion inequality} has an important consequence. Once the point $x$ and the root interval $t$ are fixed, there is essentially no freedom left in the choice of slope of the tube.
%The algebraic inequality  resulting from the intersection geometry \eqref{x in intersection} is critical in identifying at most a single correct slope per root that can yield a tube containing $x$. This is the content of the next corollary. 
\begin{corollary} \label{cor: uniqueness of possible slopes} 
Fix constants $C_0, A_0 > 1$, with $C_0$ as in Section \ref{section: model trees} and $A_0$ as in Section \ref{general facts about tube families}. Let $x = (x_1, x_2) \in \mathbb R^2$, and suppose that 
\begin{equation}  
x_1 \in I_0 := [A_0, A_0+1]. \label{distance from the root line} 
\end{equation}  
Then for every $t \in \mathcal Q(J)$, there exists at most one slope $\omega \in \Omega_N$ such that 
\begin{equation} x \in {\tt{Tube}}(t, \omega, I_0), \label{x containment}  \end{equation}
with {\tt{Tube}} defined as in \eqref{defn: tube}.  
\end{corollary} 
\begin{proof} 
Towards a contradiction, suppose that $\omega, \omega' \in \Omega_N$ are two distinct slopes with the property that 
\begin{equation*} 
 x \in {\tt{Tube}}(t, \omega, I_0) \cap {\tt{Tube}}(t, \omega', I_0). 
 \end{equation*}  
Invoking \eqref{intersection criterion inequality} from Lemma \ref{intersection criterion lemma} with $t = t'$, we obtain 
\begin{equation}  \label{two slopes upper bound} 
|x_1(\omega'-\omega)| \leq c_0 M^{-J}, \; \text{ and therefore } \quad A_0 |\omega'-\omega| \leq c_0 M^{-J}. 
\end{equation} 
The last inequality uses the hypothesis \eqref{distance from the root line}, namely $x_1 \geq A_0$.  
On the other hand, the pruning mechanism in Proposition \ref{PRUNING STAGE 1} leading to the construction of the slope set $\Omega_N$ ensures a certain separation between $\omega$ and $\omega'$. Not only do they lie in distinct $M$-adic intervals of length $M^{-J}$, the finest scale $J$ was chosen to obey the stronger separation conditions \eqref{one element per terminal vertex} and \eqref{Euclidean separation terminal vertices}:
\begin{equation}  \label{two slopes lower bound} 
|\omega - \omega'| \geq C_0 M^{-J}. 
\end{equation}   
Combining \eqref{two slopes upper bound} and \eqref{two slopes lower bound} we obtain
\[ C_0 M^{-J} \leq |\omega-\omega'| \leq c_0 A_0^{-1}M^{-J}, \]
leading to the contradiction $1 < A_0C_0 \leq c_0 < 1$. This completes the proof of the lemma.
\end{proof} 

\subsection{Definition of the reference set and the reference map} \label{section: ref set and ref map def} 
Corollary \ref{cor: uniqueness of possible slopes} suggests reversing the perspective. Rather than fixing a root and asking which slopes are possible, we now fix the point $x$ and collect all roots for which such a slope exists.
\vskip0.1in
\noindent This leads to a natural definition: for $I_0$ as \eqref{distance from the root line} and  any $x \in I_0 \times \mathbb R$ 
\begin{equation} \label{defn: Poss(x)}
{\tt{Poss}}(x) := \bigl\{ t \in \mathcal Q(J) : \exists \omega \in \Omega_N \text{ such that } x \in {\tt{Tube}}(t, \omega, I_0) \bigr\}.
\end{equation} 
Thus ${\tt{Poss}}(x)$ records precisely the roots from which a tube of admissible slope can reach $x$. Figure \ref{fig:poss-x} illustrates both membership and non-membership in this set.
\vskip0.1in
\noindent Membership in ${\tt{Poss}}(x)$ identifies the admissible roots for $x$. The uniqueness statement of Corollary \ref{cor: uniqueness of possible slopes} now allows us to associate a unique reference slope to each such root.
%${\tt{Poss}}(x)$ is the collection of all roots $t$ such that a tube rooted at $t$ and oriented along some slope $\omega \in \Omega_N$ contains $x$. Such a slope $\omega$, if it exists, is necessarily unique by Corollary \ref{cor: uniqueness of possible slopes}. 
One can therefore define a map 
\begin{equation} \label{what is tau_x}
\tau_x: {\tt{Poss}}(x) \rightarrow \Omega_N \; \text{ by } \; \tau_x(t) := \omega \text{ obeying } \eqref{x containment}.
\end{equation}  
We call $\tau_x$ the {\em{reference map}} for $x$. It maps a root $t$ to its correct slope, namely the only possible slope for a tube rooted at $t$ to pass through $x$. 

\begin{figure}[ht]
\centering

\begin{tikzpicture}[
    x=2.0cm,
    y=1.35cm,
    font=\small,
    axis/.style={
        line width=0.7pt
    },
    boundary/.style={
        densely dashed,
        line width=0.6pt
    },
    possline/.style={
        blue!75!black,
        line width=0.9pt
    },
    nonpossline/.style={
        red!70!black,
        line width=0.8pt,
        densely dashed
    },
    possinterval/.style={
        blue!75!black,
        line width=3pt,
        line cap=round
    },
    nonpossinterval/.style={
        red!70!black,
        line width=3pt,
        line cap=round
    }
]

% -------------------------------------------------
% Parameters
% -------------------------------------------------
\def\A{2}
\def\xcoord{2.5}
\def\ycoord{1.5}

% Backward root heights corresponding to slopes 0, 1/2, 1
\def\tone{1.5}
\def\ttwo{0.25}
\def\tthree{-1.0}

% A root inside the allowed root segment but not in Poss(x)
\def\tfour{2.65}

% -------------------------------------------------
% Coordinate axes and root line
% -------------------------------------------------
\draw[axis,->]
    (0,-1.45) -- (0,4.45)
    node[above] {$x_2$};

\draw[axis,->]
    (-0.15,0) -- (3.55,0)
    node[right] {$x_1$};

% Allowed root segment
\draw[line width=1.2pt]
    (0,0) -- (0,4);

%\node[above left=2pt] at (-0.2,4.35)
 %   {root line};

\node[left=3pt] at (0,0) {$0$};
\node[left=3pt] at (0,4) {$A_0+2$};

% Brace indicating the allowed root segment
\draw[
    decorate,
    decoration={brace,amplitude=5pt},
    line width=0.6pt
]
(-0.12,0) -- (-0.12,4)
node[midway,left=7pt]
{$\begin{array}{c}\text{root line} \\ \{0 \} \times [0,A_0+2] \end{array}$};

% -------------------------------------------------
% Strip [A_0,A_0+1] x [0,A_0+2]
% -------------------------------------------------
\fill[gray!8]
    (\A,0) rectangle ({\A+1},4);

\draw[boundary]
    (\A,0) -- (\A,4);

\draw[boundary]
    ({\A+1},0) -- ({\A+1},4);

\node[below] at (\A,0) {$A_0$};
\node[below] at ({\A+1},0) {$A_0+1$};

\node[above=8pt,align=center]
    at ({\A+0.5},4)
    {$[A_0,A_0+1]\times[0,A_0+2]$};

% -------------------------------------------------
% Roots t_1 and t_2 in Poss(x), and their tubes
% -------------------------------------------------

% Tube from t_1, corresponding to slope 0
\draw[
    blue!18,
    line width=6pt,
    line cap=round
]
    (0,\tone) -- (\xcoord,\ycoord);

\draw[possline]
    (0,\tone) -- (\xcoord,\ycoord);

% Tube from t_2, corresponding to slope 1/2
\draw[
    blue!18,
    line width=6pt,
    line cap=round
]
    (0,\ttwo) -- (\xcoord,\ycoord);

\draw[possline]
    (0,\ttwo) -- (\xcoord,\ycoord);

% Root intervals t_1 and t_2
\draw[possinterval]
    (0,{\tone-0.09}) -- (0,{\tone+0.09});

\draw[possinterval]
    (0,{\ttwo-0.09}) -- (0,{\ttwo+0.09});

\node[left=5pt,blue!75!black]
    at (0,\tone) {$t_1$};

\node[left=5pt,blue!75!black]
    at (0,\ttwo) {$t_2$};

% -------------------------------------------------
% Root t_3 not in Poss(x):
% its backward trace lies below the allowed root segment
% -------------------------------------------------
\draw[
    red!15,
    line width=6pt,
    line cap=round
]
    (0,\tthree) -- (\xcoord,\ycoord);

\draw[
    red!70!black,
    line width=0.9pt
]
    (0,\tthree) -- (\xcoord,\ycoord);

\draw[nonpossinterval]
    (0,{\tthree-0.09}) -- (0,{\tthree+0.09});

\node[left=5pt,red!70!black]
    at (0,\tthree) {$t_3$};

\node[
    red!70!black,
    anchor=west
]
at (0.18,-1.10)
{$t_3\notin {\tt{Poss}}(x)$};

% -------------------------------------------------
% Root t_4 not in Poss(x):
% none of the three permissible slopes reaches x
% -------------------------------------------------
\draw[nonpossinterval]
    (0,{\tfour-0.09}) -- (0,{\tfour+0.09});

\node[left=5pt,red!70!black]
    at (0,\tfour) {$t_4$};

% Slope 0
\draw[nonpossline]
    (0,\tfour) -- (3.1,\tfour);

% Slope 1/2
\draw[nonpossline]
    (0,\tfour) -- (2.65,{\tfour+0.5*2.65});

% Slope 1
\draw[nonpossline]
    (0,\tfour) -- (1.55,{\tfour+1.55});

% Placed below the horizontal dashed tube
\node[
    red!70!black,
    anchor=north west,
    fill=white,
    inner sep=1pt
]
at (0.15,{\tfour-0.10})
{$t_4\notin {\tt{Poss}}(x)$};

% -------------------------------------------------
% Remaining label
% -------------------------------------------------
\node[anchor=west]
    at (3.18,3.55)
    {$\Omega=\left\{0,\frac12,1\right\}$};

% -------------------------------------------------
% Fixed point x
% Drawn last to remain visible above all tubes
% -------------------------------------------------
\filldraw[
    fill=black,
    draw=white,
    line width=0.9pt
]
(\xcoord,\ycoord) circle (3.2pt);

\fill[black]
    (\xcoord,\ycoord) circle (2.3pt);

\node[
    right=5pt,
    fill=white,
    inner sep=1pt
]
at (\xcoord,\ycoord)
{$x$};

\end{tikzpicture}

\caption{\small{
A schematic depiction of the root collection ${\tt{Poss}}(x)$ for $\Omega = \Omega_N =\{0,\frac12,1\}$. The blue tubes of permissible slopes passing through $x$ trace backward to the root intervals $t_1, t_2, t_3$. Of these, the first two lie on the allowed root segment $\{0\} \times [0, A_0 + 2]$, and constitute ${\tt{Poss}}(x)$, with $\tau_x(t_1) = 0$, $\tau_x(t_2) = \frac12$. The tube through $x$ associated with $t_3$ traces backward to an interval
below the allowed root segment. The interval $t_4$ lies on the allowed
root segment, but none of the three permissible tubes rooted at $t_4$
passes through $x$. Neither $t_3$ nor $t_4$ lie in ${\tt{Poss}}(x)$. 
}}
\label{fig:poss-x}
\end{figure}
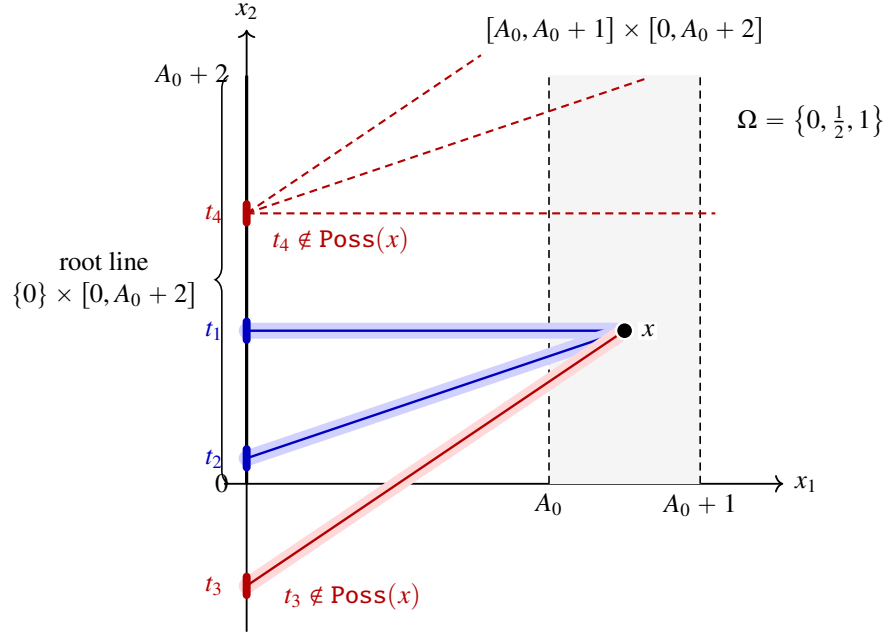

\subsection{Injectivity of the reference map} \label{section: ref map injectivity}
The next observation shows that distinct roots correspond to distinct reference slopes. This injectivity will later allow us to compare the geometry of the reference root tree with that of the reference slope tree.
\begin{lemma} \label{lemma: injectivity} 
The reference map $\tau_x$ given by \eqref{what is tau_x} is injective:
\begin{equation} \label{reference map injectivity} 
\text{ If } t \neq t', \text{ then } \tau_x(t) \neq \tau_x(t').  
\end{equation} 
\end{lemma} 
\begin{proof}
Aiming for a contradiction, let us assume if possible that there exist $t \neq t'$ such that 
\[ \omega = \tau_x(t) = \omega' = \tau_x(t'). \] Since $x$ lies in the intersection of ${\tt{Tube}}(t, \omega, I_0)$ and ${\tt{Tube}}(t', \omega', I_0)$,  the intersection condition \eqref{intersection criterion inequality} must hold. Setting $\omega = \omega'$ in this inequality results in the estimate
\begin{equation}  \label{centre-diff-upper} 
|\text{cen}(t) - \text{cen}(t')| \leq c_0 M^{-J}. 
\end{equation}   
However, $t \neq t'$, $t, t' \in \mathcal Q(J)$ implies that  
\begin{equation} \label{centre-diff-lower}
|\text{cen}(t) - \text{cen}(t') | \geq M^{-J}
\end{equation}   
Combining the two inequalities \eqref{centre-diff-upper} and \eqref{centre-diff-lower} leads to a contradiction of the assumption $c_0 < 1$. This completes the proof of the lemma. 
\end{proof} 
\noindent The objects introduced in this section currently remain sets and maps between sets. In the next section they will be reorganized into compressed trees and sticky tree maps, providing the deterministic framework needed for the probabilistic argument.
 
\section{The reference trees of slopes and roots} \label{section: reference trees}
Having identified the relevant roots and slopes, we now turn to their organization. Section \ref{section: Poss_x} identified two deterministic objects associated with a fixed point $x$, namely the reference set ${\tt{Poss}}(x)$ and the reference map $\tau_x$. The next step is to organize these objects into the tree structure anticipated in Proposition \ref{PROP: RELATED TO SURVIVAL}. Specifically, we shall show that ${\tt{Poss}}(x)$ naturally gives rise to a compressed root tree, while the reference map extends to a sticky map between this tree and a corresponding compressed slope tree.
%So far, ${\tt{Poss}}(x)$ has been defined as a sub-collection of root intervals $\mathcal Q(J)$ and $\tau_x$ as a map defined on ${\tt{Poss}}(x)$ with range in $\Omega_N$. The objective of this section is to identify the former as a compressed tree, and the latter as a sticky map defined on this tree. 
Most quantities in this section depend on $x$, which is held fixed throughout. We indicate this dependence once when each object is introduced, then subsequently suppress $x$ from the notation.  
\vskip0.1in   
\noindent To construct compressed trees, we must identify the natural ancestors of a possible slope and of its corresponding root at each fundamental height of the pruned slope tree. With this in mind, for $t \in {\tt{Poss}}(x)$ and any index $1 \leq j \leq N$, let $v_j(t) = v_j(t; x)$ and $\theta_j(t) = \theta_j(t; x)$ denote respectively the $j^{\text{th}}$ splitting vertex and the $j^{\text{th}}$ basic slope interval in the pruned slope tree $\mathcal S_N$ containing $\tau_x(t)$: 
\begin{equation} \label{theta-jtx}   
\left\{ 
\begin{aligned} 
&\tau_x(t) \in \theta_j(t) \subseteq v_j(t), \text{ where } v_j(t) = v_j(t; x) \in {\tt{SplitV}}_j,   \\ &\text{ and } \theta_j(t) = \theta_j(t; x) \in {\tt{BasicSl}}_j.
\end{aligned} 
\right\}
\end{equation} 
As we have seen in Sections \ref{section: fundamental heights} and \ref{section: basic slopes}, the length of $\theta_j(t)$ corresponds to a fundamental height, namely 
\[ |\theta_j(t)| = M^{-\lambda(v_j(t))}.   \] 
On the root side, let $Q_j^{\ast}(t) = Q_j^{\ast}(t; x)$ denote the unique $M$-adic interval, equi-dimensional with $\theta_j(t)$, satisfying
\begin{equation}  \label{Q-jtx} 
t \subseteq Q_j^{\ast}(t), \quad |Q_j^{\ast}(t)| = |\theta_j(t)|. 
\end{equation}  
%Here $\lambda_j(\cdot)$ is a $j^{\text{th}}$ fundamental height in the pruned slope tree $\mathcal S_N$; see Section \ref{section: fundamental heights}.
\vskip0.1in
\noindent The intervals $Q_j^{\ast}(t)$ and $\theta_j(t)$ naturally assemble into decreasing nested chains, one on the root side and one on the slope side. We express them as follows:
\begin{align} 	
&{\tt{PossRoot}}_j(t) = {\tt{PossRoot}}_j(t; x) := \langle Q_1^{\ast}(t), \ldots, Q_j^{\ast}(t)  \rangle,  \label{defn: PossRoot} \\ 	
&{\tt{PossSlope}}_j(t) = {\tt{PossSlope}}_j(t; x) := \langle \theta_1(t), \ldots, \theta_j(t)  \rangle, \label{defn: PossSlope} \\ 
&\tau_x(t) \in  \theta_j(t) \subsetneq \cdots \subsetneq \theta_1(t), \quad
t \subsetneq  Q_{j}^{\ast}(t) \subsetneq \cdots \subsetneq Q_1^{\ast}(t). \label{PossRoot chains}
\end{align} 
The $j^{\text{th}}$ fundamental height of $\mathcal S_N$ is non-unique in general, consequently the intervals in each collection
\[ \bigl\{ \theta_j(t) : t \in {\tt{Poss}}(x)\bigr\} \; \text{ and } \; \bigl\{ Q_j^{\ast}(t) : t \in {\tt{Poss}}(x) \bigr\}\] need not be of the same length.  Our goal is to define compressed trees ${\tt{TrPR}}_x$ and ${\tt{TRPS}}_x$ of roots and slopes, in which the $j^{\text{th}}$-level vertices are given respectively by ${\tt{PossRoot}}_j(t)$ and ${\tt{PossSlope}}_j(t)$, as $t$ ranges over ${\tt{Poss}}(x)$. The remainder of this section establishes that these nested interval chains are compatible across different choices of $t$, allowing them to be assembled into well-defined compressed trees of roots and slopes that encode the deterministic geometry associated with $x$.
\subsection{Description of ${\tt{TrPR}}_x$ and ${\tt{TrPS}}_x$ as trees} 
The collections of nested interval sequences introduced in \eqref{defn: PossRoot} and \eqref{defn: PossSlope} suggest the structure of rooted trees. However, this is not automatic. Since the intervals at the $j^{\text{th}}$ level may have different lengths, an arbitrary family of such chains need not satisfy the consistency conditions required of a rooted labelled tree in the sense of Chapter~\ref{trees-section}.
%A priori, an arbitrary collection of nested interval sequences of the form \eqref{defn: PossRoot} or \eqref{defn: PossSlope}, where the $j^{\text{th}}$ entry intervals may have variable lengths, need not generate a rooted, labelled tree in the sense of Chapter \ref{trees-section}. 
As we saw in that chapter, a rooted labelled tree structure is determined by consistent heights and lineages of vertices. For a collection of chains of the form \eqref{defn: PossRoot} or \eqref{defn: PossSlope}, inconsistencies in heights and lineages may occur, contradicting the tree structure in one of the following ways: for $t \neq t'$, one could have
\begin{align} 
& Q_j^{\ast}(t) = Q_j^{\ast}(t') \; \text{ but }  Q_{k}^{\ast}(t) \neq Q_{k}^{\ast}(t') \text{ for some } k < j, \; \text{ or } \label{scenario 1} \\
& Q_j(t) = Q_k(t') \;\text{ for some } j \neq k. \label{scenario 2} 
\end{align} 
Either scenario would contradict the defining requirements of a tree: 
\vskip0.1in 
\begin{itemize} 
\item The first condition \ref{scenario 1} fails the tree criterion of unique ancestry: here the same interval occurs as the final entry of two $j$-long distinct chains ${\tt{PossRoot}}_j(t)$ and ${\tt{PossRoot}}_j(t')$. It corresponds to the situation of a non-unique ray, where a single vertex at height $j$ has different ancestors at some height $k$ of the tree, $k < j$. 
\vskip0.1in
\item The second condition \eqref{scenario 2} fails the tree criterion of a well-defined height for every vertex. This scenario occurs when the same interval appears in different positions in different chains, giving rise to ambiguity in its height. 
%(situation ) or has two different parents (situation \eqref{scenario 1}). 
%The latter scenario occurs when two $M$-adic ancestors of different lengths appear in the same generation. 
 \end{itemize} 
 \vskip0.1in
\noindent The next lemma shows that no such inconsistencies occur.  
 \begin{lemma} \label{CONSISTENCY LEMMA}
 Let ${\tt{Poss}}(x)$, $\tau_x$, ${\tt{PossRoot}}_j(t)$ and ${\tt{PossSlope}}_j(t)$ be defined as above in \eqref{defn: Poss(x)}, \eqref{what is tau_x}, \eqref{defn: PossRoot}--\eqref{PossRoot chains}. Then the following conclusions hold: 
\vskip0.1in
 \begin{enumerate}[(a)] 
\item {\em{(Consistency of rays with common entries)}} If two chains of the form \eqref{defn: PossRoot} (or \eqref{defn: PossSlope}) have a common entry, then that entry must occur at the same position within both chains, and the chains must coincide on all entries up to that point. \label{consistency part1} 
\vskip0.1in
\noindent Specifically, for $t, t' \in {\tt{Poss}}(x)$ and $1 \leq j, j' \leq N$, suppose that 
 \begin{equation} \label{Qjj'-equal}
 Q_j^{\ast}(t) = Q_{j'}^{\ast}(t'). 
 \end{equation}
 Then, 
 \begin{equation} \label{consistency check}
 j = j' \; \; \text{ and } \; \; Q_k^{\ast}(t) = Q_k^{\ast}(t'), \; \theta_k(t) = \theta_k(t') \; \text{ for } 1 \leq k \leq j. 
 \end{equation}  
 In particular, neither \eqref{scenario 1} nor \eqref{scenario 2} can hold. 
 \vskip0.1in
\item {\em{(Definition of trees ${\tt{TrPR}}_x$ and ${\tt{TrPS}}_x$)}} \label{consistency part2}  There exist deterministic rooted, labelled trees ${\tt{TrPR}}_x$ and ${\tt{TrPS}}_x$ of height $N$, with vertex sets $\mathcal V_j(\cdot)$ defined by 
\begin{align} 
&\mathcal V_0({\tt{TrPR}}_x) := \mathcal V_0({\tt{TrPS}}_x) = \{[0,1] \}, \label{level-0} \\ 
&\mathcal V_j({\tt{TrPR}}_x) := \left\{ {\tt{PossRoot}}_j(t) : t \in {\tt{Poss}}(x) \right\},  \; 1 \leq j \leq N, \label{level-j} \\  
&\mathcal V_j({\tt{TrPS}}_x) := \left\{ {\tt{PossSlope}}_j(t) : t \in {\tt{Poss}}(x) \right\}, \; 1 \leq j \leq N. \label{possible slopes j} 
\end{align}  
Ancestry within each tree is given by the usual containment relation. 
 \end{enumerate} 
 \end{lemma} 
\vskip0.1in
\noindent {\em{Remarks: }} 
\begin{enumerate}[1.]
\item The proof of Lemma \ref{CONSISTENCY LEMMA} is given in Section \ref{section: proof of consistency lemma}. 
\vskip0.1in
\item The lemma describes the tree ${\tt{TrPR}}_x$ whose existence was asserted in Proposition \ref{PROP: RELATED TO SURVIVAL}\eqref{survival-part1}. 
\vskip0.1in
\item Lemma \ref{CONSISTENCY LEMMA} marks the point at which the deterministic objects ${\tt{Poss}}(x)$ and its image under $\tau_x$ acquire the compressed tree structure needed for the remainder of the argument. In particular, it permits the application of the tree-theoretic machinery developed in Chapter~\ref{trees-section}. 
\vskip0.1in
%\item Every vertex in a tree uniquely determines the ray leading up to it. Thus in light of Lemma \ref{CONSISTENCY LEMMA}, and to simplify notation, we will henceforth identify the sequences ${\tt{PossRoot}}_j(t)$ and ${\tt{PossSlope}}_j(t)$ with their last entries $Q_j^{\ast}(t)$ and $\theta_j(t)$ respectively; so,  
%\begin{align*} 
%{\tt{PossRoot}}_j(t) &= \langle Q_1^{\ast}(t), \ldots, Q_j^{\ast}(t) \rangle \longleftrightarrow Q_j^{\ast}(t), \\ 
%{\tt{PossSlope}}_j(t) &= \langle \theta_1(t), \ldots, \theta_j(t) \rangle \longleftrightarrow \theta_j(t). 
%\end{align*} 
%These final entries uniquely determine the sequences \eqref{defn: PossRoot} and \eqref{defn: PossSlope}, through the nesting property of $M$-adic intervals, and the defining conditions \eqref{theta-jtx} and \eqref{Q-jtx}. Proposition \ref{PROP: REFERENCE TREES} is proved in Section \ref{proof: reference trees}.  
\end{enumerate} 
 
\subsection{Identification of ${\tt{Poss}}(x)$ and $\tau_x({\tt{Poss}}(x))$ through ${\tt{TrPR}}_x$ and ${\tt{TrPS}}_x$}  
Now that ${\tt{TrPR}}_x$ and ${\tt{TrPS}}_x$ have been confirmed as compressed trees, it is natural to ask whether they continue to faithfully encode the root and the slope sets they originated from. The following lemma answers this question in the affirmative, clarifying the scope of this correspondence. 
%which geometric sets they correspond to, in the sense of Section \ref{section: compressed trees}. Each maximal ray of ${\tt{TrPR}}_x$ (respectively ${\tt{TrPS}}_x$) ends in an interval of the $N^{\text{th}}$ fundamental height. Such an interval, by itself, need not be an element of ${\tt{Poss}}(x)$ (respectively $\Omega_N$). Nonetheless, each such interval contains a unique element of ${\tt{Poss}}(x)$ (respectively a unique slope in $\Omega_N$), as the following lemma shows.  
\begin{lemma} \label{TREE SET LEMMA} 
Let ${\tt{TrPR}}_x$ and ${\tt{TrPS}}_x$ be the two compressed trees of height $N$ established in Lemma \ref{CONSISTENCY LEMMA} \eqref{consistency part2}. 
\vskip0.1in
\begin{enumerate}[(a)] 
\item The tree ${\tt{TrPR}}_x$ represents ${\tt{Poss}}(x)$ in the following sense:  \label{consistency part2a}
\vskip0.1in 
\begin{itemize} 
\item Every terminal vertex of ${\tt{TrPR}}_x$ identifies a single root interval $t \in {\tt{Poss}}(x)$; in other words, for every vertex $\langle Q_1, \ldots , Q_N\rangle \in \mathcal V_N({\tt{TrPR}}_x)$ there exists a unique  $t \in {\tt{Poss}}(x)$ such that 
\begin{equation}  \label{tree to root interval} 
 t \subsetneq Q_N = Q_N^{\ast}(t), \quad Q_j = Q_j^{\ast}(t) \text{ for } 1 \leq j \leq N.  
\end{equation} 
\item Conversely, every $t \in {\tt{Poss}}(x)$ can be associated with a unique vertex  
\[{\tt{PossRoot}}_N(t) = \langle Q_1^{\ast}(t), \ldots, Q_N^{\ast}(t) \rangle  \text{ of ${\tt{TrPR}}_x$ obeying } \eqref{tree to root interval}.  \] 
\end{itemize} 
\vskip0.1in 
\item Similarly, the tree ${\tt{TrPS}}_x$ represents $\tau_x \bigl({\tt{Poss}}(x) \bigr)$ in the following sense: \label{consistency part2b}
\vskip0.1in
\begin{itemize} 
\item Every terminal vertex of ${\tt{TrPS}}_x$ contains a unique slope of $\tau_x \bigl({\tt{Poss}}(x)\bigr)$. In other words, for every ${\tt{PossSlope}}_N(t) = \langle \theta_1(t) , \ldots, \theta_N(t) \rangle \in \mathcal V_N({\tt{TrPS}}_x)$, there is a unique $\omega \in \tau_x({\tt{Poss}}(x))$ such that 
\begin{equation}  \label{tree to slope} 
\omega = \tau_x(t) \in \theta_N(t).  \end{equation}   
\vskip0.1in
\item Conversely, for every $\omega \in \tau_x \bigl({\tt{Poss}}(x) \bigr)$, there exists $t \in {\tt{Poss}}(x)$ such that
\[  {\tt{PossSlope}}_N(t) = \langle \theta_1(t) ,  \ldots, \theta_N(t) \rangle  \text{ obeys  \eqref{tree to slope}}.  \]   
\end{itemize} 
\end{enumerate} 
\end{lemma} 
 \vskip0.1in
\noindent {\em{Remarks: }} 
\begin{enumerate}[1.]
\item Lemma \ref{TREE SET LEMMA} is proved in Section \ref{section: tree set lemma proof}.
\vskip0.1in
\item It is important to note that while Lemma \ref{TREE SET LEMMA} isolates $Q_N^{\ast}(t)$ as a unique identifier of the root $t \in {\tt{Poss}}(x)$, the two intervals are not equal. The interval $Q_N^{\ast}(t)$ may contain many root intervals from $\mathcal Q_J$, of which only one lies in ${\tt{Poss}}(x)$. 
\vskip0.1in
\item The analogous interpretation for ${\tt{TrPS}}_x$ is easier. Since each terminal vertex $\theta_N(t) \in {\tt{BasicSl}}_N$, it contains exactly one slope of $\Omega_N$ by virtue of our pruning process.  This slope is $\tau_x(t)$. 
%\vskip0.1in
%\item The correspondence claimed in Lemma \ref{TREE SET LEMMA} allows an unambiguous representation of ${\tt{Poss}}(x)$ and $\tau_x({\tt{Poss}}(x))$ using the respective trees ${\tt{TrPR}}_x$ and ${\tt{TrPS}}_x$. 
\end{enumerate}

%\item {\em{(Consistency of rays with common entries)}} If two chains of the form \eqref{defn: PossRoot} (or \eqref{defn: PossSlope}) have a common entry, then that entry must occur at the same position within both chains, and the chains must coincide on all entries up to that point. \label{consistency part1} 
%\vskip0.1in
%\noindent Specifically, for $t, t' \in {\tt{Poss}}(x)$ and $1 \leq j, j' \leq N$, suppose that 
% \begin{equation} \label{Qjj'-equal}
% Q_j^{\ast}(t) = Q_{j'}^{\ast}(t'). 
% \end{equation}
% Then, 
% \begin{equation} \label{consistency check}
% j = j' \; \; \text{ and } \; \; Q_k^{\ast}(t) = Q_k^{\ast}(t'), \; \theta_k(t) = \theta_k(t') \; \text{ for } 1 \leq k \leq j. 
% \end{equation}  
% \vskip0.1in
%\item {\em{(Definition of trees ${\tt{TrPR}}_x$ and ${\tt{TrPS}}_x$)}} \label{consistency part2}  There exist well-defined deterministic trees ${\tt{TrPR}}_x$ and ${\tt{TrPS}}_x$ of height $N$, whose vertex sets $\mathcal V_j(\cdot)$ defined by 
%\begin{align} 
%&\mathcal V_0({\tt{TrPR}}_x) := \mathcal V_0({\tt{TrPS}}_x) = \{[0,1] \}, \label{level-0} \\ 
%&\mathcal V_j({\tt{TrPR}}_x) := \left\{ {\tt{PossRoot}}_j(t) : t \in {\tt{Poss}}(x) \right\},  \; 1 \leq j \leq N, \label{level-j} \\  
%&\mathcal V_j({\tt{TrPS}}_x) := \left\{ {\tt{PossSlope}}_j(t) : t \in {\tt{Poss}}(x) \right\}, \; 1 \leq j \leq N,  \label{possible slopes j} 
%\end{align}  
%follow the defining criteria for height and ancestry required of a rooted, labelled tree, equipped with the usual containment relations. 
%\vskip0.1in
\subsection{Definition of $\tau_x$ as a sticky map} 
Having established ${\tt{TrPR}}_x$ and ${\tt{TrPS}}_x$ as compressed trees representing ${\tt{Poss}}(x)$ and $\tau_x({\tt{Poss}}(x))$ respectively, our next step is to study the reference map $\tau_x$ from Section \ref{section: ref set and ref map def} in relation to these trees. The following lemma asserts that $\tau_x$ extends to a height and ancestry-preserving map between the two trees in a natural way. 
\begin{lemma} \label{LEMMA: MAP LIFT} The map $\tau_x$, initially defined on ${\tt{Poss}}(x)$ via \eqref{what is tau_x}, lifts to a well-defined, length-preserving sticky map between compressed trees, according to the following definition and in the sense of Section \ref{sticky maps section}:  
\begin{equation}  \label{tau sticky tree map}
\tau_x : {\tt{TrPR}}_x \longrightarrow {\tt{TrPS}}_x, \quad \tau_x \bigl({\tt{PossRoot}}_j(t) \bigr) :=  {\tt{PossSlope}}_{j}(t),
\end{equation}  
for $t \in {\tt{Poss}}(x)$. Accordingly, we write $\tau_x(Q_j^{\ast}(t)) = \theta_j(t)$. 
\vskip0.1in
\noindent In particular, there exists a deterministic binary sequence $\mathbb Z_x$ with entries either $0$ or $1$, indexed by the vertices of ${\tt{TrPR}}_x$:
\begin{align}  
&\mathbb Z_x  = \left\{ Z_x(Q) : Q = Q^{\ast}_j(t), \; t \in {\tt{Poss}}(x) \right\} \text{ such that } \label{Z_x1} \\  
&\tau_x \bigl(Q_j^{\ast}(t)\bigr) = \theta_j(t) = \text{ the } Z_x\bigl(Q_j^{\ast}(t) \bigr)^{\text{th}} \text{ child  of } \theta_{j-1}(t) \label{Z_x2} \\
&\hskip1.1in = \text{ the } Z_x\bigl(Q_j^{\ast}(t) \bigr)^{\text{th}} \text{ child  of } \tau_x\bigl( Q_{j-1}^{\ast}(t)\bigr) \noindent
 \end{align}
 for all $t \in {\tt{Poss}}(x)$.
\end{lemma} 
\vskip0.1in
\noindent {\em{Remarks:}}
\vskip0.1in
\begin{enumerate}[1.]
\item Lemma \ref{LEMMA: MAP LIFT} is proved in Section \ref{proof section: lift map}.
\vskip0.1in
\item Apart from its intrinsic geometric significance, Lemma \ref{LEMMA: MAP LIFT} identifies $({\tt{TrPR}}_{x}, \tau_x)$ as a restriction of an adapted tree-map pair $(\mathscr{U}, \tau)$, in the sense of Section \ref{section: tree-map pair defn}. This interpretation will turn out to be important in Lemma \ref{PROP: REFERENCE TREES} below.  
\end{enumerate} 
\subsection{Uniformity of progeny in ${\tt{TrPR}}_x$ and ${\tt{TrPS}}_x$} 
In Lemma \ref{lemma: injectivity}, we established that the map $\tau_x: {\tt{Poss}}(x) \rightarrow \Omega_N$ is injective. Our next result establishes that the lifted map $\tau_x: {\tt{TrPR}}_x \rightarrow {\tt{TrPS}}_x$ is essentially injective, in the sense that the pre-image of every vertex is of bounded cardinality. 
\begin{lemma} \label{LEMMA: UNIFORMITY OF PROGENY} There is a constant $C_2 \geq 1$, depending only on $A_0$, such that \begin{equation} 
\# \Bigl[ \tau_x^{-1}\bigl(\theta \bigr) \Bigr] \leq C_2 \; \text{ for all } \theta \in {\tt{BasicSl}}_j, \; 1 \leq j \leq N. \label{uniformity of progeny}
\end{equation}  
In other words, at any given height $j$ of ${\tt{TrPR}}_x$ and ${\tt{TrPS}}_x$, there are at most $C_2$ vertices of ${\tt{TrPR}}_x$ that map under $\tau_x$ to a single vertex of ${\tt{TrPS}}_x$ at that height. 
\vskip0.1in
\noindent As a result, each non-terminal vertex of ${\tt{TrPR}}_x$ has at most $2C_2$ children, and 
\begin{equation} \label{cardinality of V_j} 
\# \bigl[ \mathcal V_j \bigl( {\tt{TrPR}}_x \bigr)\bigr] \leq C_2 2^{j}. 
\end{equation}  
\end{lemma} 
\vskip0.1in
\noindent {\em{Remarks: }} 
\begin{enumerate}[1.]
\item Lemma \ref{LEMMA: UNIFORMITY OF PROGENY} is proved in Section \ref{proof section: uniformity of progeny}.
\vskip0.1in
\item The reader will recognize the assertion of Lemma \ref{LEMMA: UNIFORMITY OF PROGENY} as identical to that of Proposition \ref{PROP: RELATED TO SURVIVAL} \eqref{survival-part2}. 
\end{enumerate} 
\subsection{Characterization of inclusion in $\mathtt K_{\mathbb X}$} 
Having identified the geometry and structure of the tree ${\tt{TrPR}}_x$, we are now in a position to specify its role in the inclusion event $x \in \mathtt K_{\mathbb X}$. 
\vskip0.1in
\noindent Let us recall from Proposition \ref{PROP: COMPRESSED ROOTS} the compressed (random) root tree $\mathscr{U}_{\mathbb X}$ associated with the binary sequence $\mathbb X$. As \eqref{defn: sticky tube} and \eqref{defn: extended Kakeya set K} attest, the tree $\mathscr{U}_{\mathbb X}$ and the slope map $\sigma_{\mathbb X}$ are instrumental in the definition of the random Kakeya set $\mathtt K_{\mathbb X}$. The next result identifies the relation between the reference pair $({\tt{TrPR}}_x, \tau_x)$ and the random pair $(\mathscr{U}_{\mathbb X}, \sigma_{\mathbb X})$ that ensures $x$ is contained in $\mathtt K_{\mathbb X}$. 
\begin{lemma}
\label{PROP: REFERENCE TREES}
 %Let ${\tt{Poss}}(x)$, $\tau_x$, ${\tt{PossRoot}}_j(t)$ and ${\tt{PossSlope}}_j(t)$ be defined as in \eqref{defn: Poss(x)}, \eqref{what is tau_x}, \eqref{defn: PossRoot}--\eqref{PossRoot chains}. 
For $x \in I_0 \times \mathbb R$, the following conclusions hold: 
\begin{enumerate}[(a)] 
\item \label{consistency part5} The point $x$ lies in $\mathtt K_{\mathbb X}$ if and only if  
\begin{equation} \label{sigma = tau} 
\left\{
\begin{aligned} 
&\text{the two trees ${\tt{TrPR}}_x$ and $\mathscr{U}_{\mathbb X}$ share at least one maximal ray} \\ 
&\mathcal R \in \partial {\tt{TrPR}}_x \cap \partial \mathscr{U}_{\mathbb X}  \; \text{ such that }  \; \sigma_{\mathbb X} \equiv \tau_x \text{ on every vertex of } \mathcal R. 
\end{aligned} 
\right\} 
\end{equation}  
\vskip0.1in 
\item \label{consistency part6} The inclusion characterization above is also equivalent to $\mathbb X$ satisfying the following property: 
\begin{equation} \label{XZ} 
\left\{
\begin{aligned} 
&\exists t \in {\tt{Poss}}(x), \text{ i.e. }  {\tt{PossRoot}}_N(t) = \langle Q_1^{\ast}(t), \ldots, Q_N^{\ast}(t)\rangle \text{ such that } \\
& X(Q) = Z_x(Q) \quad \text{ for all } Q = Q_j^{\ast}(t), \text{ and all } 1 \leq j \leq N.    
\end{aligned} 
\right\}
\end{equation} 
Here $Z_x(\cdot)$ is the deterministic binary sequence given by \eqref{Z_x1}, \eqref{Z_x2}. 
\end{enumerate}  
\end{lemma}
\vskip0.1in
\noindent {\em{Remarks: }}
\begin{enumerate}[1.]
\item Lemma \ref{PROP: REFERENCE TREES} is proved in Section \ref{proof: reference trees}. 
\vskip0.1in 
\item Lemma \ref{PROP: REFERENCE TREES} plays a key role in the proof of Proposition \ref{PROP: RELATED TO SURVIVAL} \eqref{survival-part3}, by identifying the inclusion event $x \in \mathtt K_{\mathbb X}$ with the matching of two $N$-long ordered binary sequences
\[ \bigl\{X(Q_j^{\ast}(t)) : 1 \leq j \leq N \bigr\} = \bigl\{Z_x(Q_j^{\ast}(t)) : 1 \leq j \leq N \bigr\}  \]  
for some $t \in {\tt{Poss}}(x)$. Let us now head into the proof of this proposition.
\end{enumerate} 
\section{Proof of Proposition \ref{PROP: RELATED TO SURVIVAL}, assuming Lemmas \ref{CONSISTENCY LEMMA}--\ref{PROP: REFERENCE TREES}} \label{survival prop proof}
\noindent The properties of the tree ${\tt{TrPR}}_x$ listed in Lemmas \ref{CONSISTENCY LEMMA}--\ref{PROP: REFERENCE TREES} are essential ingredients in the proof of Proposition \ref{PROP: RELATED TO SURVIVAL}. We complete it here. 
\begin{proof} 
Once Lemma \ref{CONSISTENCY LEMMA} \eqref{consistency part2} identifies ${\tt{TrPR}}_x$ and ${\tt{TrPS}}_x$ as rooted, labelled trees, it becomes clear from \eqref{level-0} and \eqref{level-j} that both trees have height $N$, and all maximal rays are of the same length $N$. Proposition \ref{PROP: RELATED TO SURVIVAL}\eqref{survival-part1}  follows from this.
\vskip0.1in 
\noindent Proposition \ref{PROP: RELATED TO SURVIVAL}\eqref{survival-part2}, which bounds the number of children of every non-terminal vertex of ${\tt{TrPR}}_x$ and the cardinality of $\mathcal V_j \bigl({\tt{TrPR}}_x\bigr)$, is identical to Lemma \ref{LEMMA: UNIFORMITY OF PROGENY}, with $C_1 = 2C_2$. 
\vskip0.1in
\noindent Proposition \ref{PROP: RELATED TO SURVIVAL} \eqref{survival-part3} contains the heart of the proof. For $x \in I_0 \times \mathbb R$, Lemma \ref{PROP: REFERENCE TREES} \eqref{consistency part6} implies 
\begin{equation} \label{prob1}
\mathbb P \bigl(x \in \mathtt K_{\mathbb X}\bigr) = \mathbb P(\,\eqref{XZ} \text{ holds} \,). 
\end{equation}  
We claim that
\begin{equation} \label{prob2}
\mathbb P(\,\eqref{XZ} \text{ holds} \,) = \varrho \bigl({\tt{TrPR}}_x; \frac{1}{2} \bigr).  
\end{equation}   
The desired conclusion \eqref{prob identity} follows by combining \eqref{prob1} and \eqref{prob2}. 
\vskip0.1in
\noindent To justify the claim \eqref{prob2}, we construct a Bernoulli percolation on the tree ${\tt{TrPR}}_x$, as follows: if $e$ is an edge of ${\tt{TrPR}}_x$ terminating at the vertex $Q$, then we set 
 \begin{equation} 
 Y_{e} := \begin{cases} 1 &\text{ if } X(Q) = Z_x(Q), \\ 0 &\text{ otherwise. }\end{cases} 
 \end{equation} 
Said differently, the edge $e$ is retained if and only if the slope allocated  to $Q$ by $\sigma_{\mathbb X}$ matches the reference slope $\tau_x(Q)$. Let us confirm that this is a standard Bernoulli$(\frac{1}{2})$ percolation on ${\tt{TrPR}}_x$:
\vskip0.1in
\begin{itemize}
\item According to the probabilistic set-up described in Section \ref{section: probabilistic setup}, the collection 
\begin{equation} \label{XQ-Bernoulli}  
\{X(Q) : Q \in \mathcal V({\tt{TrPR}}_x) \} 
\end{equation}  is an i.i.d. Bernoulli$(\frac{1}{2})$ sequence. Indeed, the tree structure of ${\tt{TrPR}}_x$ ensures that an $M$-adic interval $Q \in \mathcal Q^{\ast}$ appears in ${\tt{TrPR}}_x$ at most once. This ensures that the random variables $X(Q)$ appearing in the collection \eqref{XQ-Bernoulli} are independent. 
\vskip0.1in
\item On the other hand, Lemma \ref{LEMMA: MAP LIFT} identifies $\{Z_x(Q) : Q \in \mathcal V({\tt{TrPR}}_x) \}$ as a deterministic sequence derived from $\tau_x$. 
\end{itemize} 
\vskip0.1in
\noindent Thus the random variables $\{Y_e : e \in \mathcal E({\tt{TrPR}}_x) \}$ are independent, with 
 \begin{align*} 
 &\mathbb P\bigl(Y_e = 1 \bigr) = \mathbb P\bigl(X(Q) = Z_x(Q)\bigr) = \frac{1}{2}, \text{ and hence } \\ 
& \mathbb P \bigl( Y_e = 0 \bigr) = \mathbb P\bigl(X(Q) \neq Z_x(Q)\bigr) = \frac12. 
 \end{align*}  
In this set-up, the event \eqref{XZ} encoding a ``full match'', where every vertex of a maximal ray of ${\tt{TrPR}}_x$ is mapped to its correct reference slope at the same height under $\sigma_{\mathbb X}$, is equivalent to the retention of a maximal ray of ${\tt{TrPR}}_x$ after the retention-removal process, i.e., to the survival of ${\tt{TrPR}}_x$ under a standard Bernoulli$(\frac{1}{2})$ percolation. This concludes the proof of Proposition \ref{PROP: RELATED TO SURVIVAL} \eqref{survival-part3}. 
\end{proof} 

\chapter{Geometry of reference trees} \label{chapter: reference tree geometry} 
In Chapter~\ref{chapter: far from root}, the construction of Kakeya-type sets was reduced to a collection of structural properties of the reference trees ${\tt{TrPR}}_x$ and ${\tt{TrPS}}_x$, summarized in Lemmas \ref{CONSISTENCY LEMMA}--\ref{PROP: REFERENCE TREES}. The purpose of the present chapter is to establish these properties. Together with the probabilistic construction developed in Chapters~\ref{chapter: sticky rectangles chapter}--\ref{chapter: far from root}, this supplies the final ingredient in the implication \eqref{condition 3: slopes sublacunary} $\implies$ \eqref{condition 1: Kakeya-type sets} of Theorem~\ref{thm:main}.

%In this chapter, we fill in the pending proofs from Chapter \ref{chapter: far from root} concerning the structure of the reference trees ${\tt{TrPR}}_x$ and ${\tt{TrPS}}_x$. This supplies the final ingredients of the implication \eqref{condition 3: slopes sublacunary} $\implies$ \eqref{condition 1: Kakeya-type sets} in Theorem \ref{thm:main}.  
\section{Weak stickiness of $\tau_x$} 
Before embarking on the proofs, let us first record a fact about the $M$-adic scales occurring in ${\tt{Poss}}(x)$ and its image under $\tau_x$. This observation will play an important role in establishing the tree structure of ${\tt{TrPR}}_x$ and ${\tt{TrPS}}_x$ (Lemma \ref{CONSISTENCY LEMMA}), and may be of independent interest. 
\vskip0.1in
\noindent For root intervals $t, t' \in {\tt{Poss}}(x)$, $t \neq t'$, let us denote 
\begin{equation} \label{ycm of roots} 
u := D_{\mathcal T}(t, t') = \text{ the youngest common ancestor of } (t, t')  
\end{equation} 
in the full $M$-adic tree $\mathcal T = \mathcal T([0,1]; M)$ defined in \eqref{tree encoding}. This is, by definition, the smallest $M$-adic interval that contains both $t$ and $t'$. 
\vskip0.1in
\noindent At the same time, the membership of $t, t'$ in ${\tt{Poss}}(x)$ means that Lemma \ref{lemma: injectivity} may be applied, ensuring that the two corresponding slopes are distinct, i.e., $\tau_x(t) \neq \tau_x(t')$. These slopes are elements of the pruned slope set $\Omega_N$, whose $M$-adic representation is $\mathcal S_N$. The pruning criterion \eqref{one element per terminal vertex} therefore implies that $\tau_x(t)$ and $\tau_x(t')$ lie in distinct terminal vertices of $\mathcal S_N$.  Therefore we may define
\begin{equation} \label{ycm of slopes}  
\begin{aligned} 
 v &:= D_{\mathcal T}\bigl(\tau_x(t), \tau_x(t') \bigr) = D_{\mathcal S_N}\bigl(\tau_x(t), \tau_x(t') \bigr) \\
 &\;= \text{ the youngest common ancestor of } (\tau_x(t), \tau_x(t')) 
\end{aligned} 
\end{equation}  
in the pruned $M$-adic slope tree $\mathcal S_N \subseteq \mathcal T$. 
\vskip0.1in
\noindent It turns out that the $M$-adic scale of the root vertex $u$ is always strictly smaller than the fundamental height of the slope vertex $v$. This relation, proved in the following lemma, is crucial for the proof of  Lemma \ref{CONSISTENCY LEMMA}. 
 \begin{lemma} \label{uv-lemma}
For $t, t' \in  {\tt{Poss}}(x)$, $t \neq t'$, let $u, v$ be as in  \eqref{ycm of roots} and \eqref{ycm of slopes}. Then 
\begin{equation} \label{uv-relation}
h(u) < \lambda(v), \; \text{ where } \; |u| = M^{-h(u)}.
\end{equation} 
The quantity $\lambda(v)$, defined in \eqref{defn: fundamental height}, refers to the fundamental height of the splitting vertex $v \in \mathcal S_N$. 
\end{lemma} 
\vskip0.1in 
\noindent {\em{Remark: }} The conclusion \eqref{uv-relation} of Lemma \ref{uv-lemma} may be viewed as a weak form of stickiness in the reference map $\tau_x$. It is not difficult to see that the reference map 
\[ \tau_x: {\tt{Poss}}(x) \rightarrow \Omega_N \] need not extend as a sticky map between the corresponding $M$-adic trees 
\[ \tau_x: \mathcal T({\tt{Poss}}(x); M) \rightarrow \mathcal T\bigl(\tau_x({\tt{Poss}}(x)); M \bigr) \text{ in the sense of Section \ref{sticky maps section}}.  \] 
Indeed, if $\tau_x$ were sticky, it would force the relations
\begin{equation*} 
|\tau_x(t)| = |t|, \; |\tau_x(t')| = |t'|, \; v = D_{\mathcal T}\bigl(\tau_x(t), \tau_x(t') \bigr) \subseteq  \tau_x(u) = \tau_x \bigl(D_{\mathcal T}(t, t') \bigr), \end{equation*}
leading to the inequality 
\begin{equation} \label{what-is-not}  
 h(u) = h\bigl(D_{\mathcal T}(t, t') \bigr) \leq h(v) = h \bigl( D_{\mathcal T}\bigl(\tau_x(t), \tau_x(t') \bigr)\bigr), 
 \end{equation} 
which would then imply $h(u) \leq h(v) < \lambda(v)$ in view of \eqref{h-star-v}. This is precisely the relation claimed in \eqref{uv-relation}. In this sense,  \eqref{uv-relation} is a weak manifestation of stickiness. 
\vskip0.1in
\noindent However, the inequality in \eqref{what-is-not} is not true in general, as can be seen by choosing 
\[ u = \Bigl[0, \frac{1}{2} \Bigr], \; v = [0, 1], \; t \subsetneq \Bigl[0, \frac{1}{4} \Bigr], \; t' \subsetneq \Bigl[\frac{1}{4}, \frac{1}{2} \Bigr], \; \tau_x(t) = 1, \; \tau_x(t') = 0. \]
Thus $h(v) = 0 < h(u) = 1$, but tubes rooted at $t$ and $t'$ with slopes $\tau_x(t)$ and $\tau_x(t')$ intersect have non-trivial intersection away from the root line. 
\vskip0.1in    
\begin{proof} 
We will prove \eqref{uv-relation} by estimating the slope difference $|\tau_x(t) - \tau_x(t')|$ from above and below, and then comparing the two bounds. 
\vskip0.1in
\noindent On one hand, it follows from the definition \eqref{defn: Poss(x)} of ${\tt{Poss}}(x)$ and the assumption $t, t' \in {\tt{Poss}}(x)$ that $x$ lies in the two tubes ${\tt{Tube}}$ and ${\tt{Tube}}'$  given by \eqref{tube and tube'}, with $\omega = \tau_x(t)$ and $\omega' = \tau_x(t')$.  
The location of $x \in I_0$ means that the intersection condition \eqref{intersection criterion inequality} from Lemma \ref{intersection criterion lemma} holds, yielding the estimate:
\begin{align} 
A_0 |\tau_x(t) - \tau_x(t')| &\leq |x_1\bigl(\tau_x(t) - \tau_x(t') \bigr)| \nonumber \\ 
&\leq \bigl| \text{cen}(t') - \text{cen}(t) + x_1 \bigl(\tau_x(t')-\tau_x(t) \bigr) \bigr| + \bigl|  \text{cen}(t') - \text{cen}(t) \bigr| \nonumber \\
&\leq c_0 M^{-J} + M^{-h(u)} < 2M^{-h(u)}. \label{slope-diff-upper}
\end{align}  
The penultimate inequality uses the fact that $\text{cen}(t) \in t \subseteq u$, $\text{cen}(t') \in t' \subseteq u$, which leads to 
\[| \text{cen}(t') - \text{cen}(t)| \leq \text{diam}(u) = M^{-h(u)}. \] The last inequality follows from the estimate $c_0 M^{-J} < M^{-h(u)}$, which is a consequence of $c_0 < 1$ and $J > h(u)$. 
\vskip0.1in
\noindent On the other hand, it follows from Lemma \ref{lemma: distance between pruned slopes} that the distance $|\tau_x(t) - \tau_x(t')|$ is determined by the fundamental height of their youngest common ancestor in $\mathcal S_N$; indeed the conclusion \eqref{distance between slopes} of Lemma \ref{lemma: distance between pruned slopes} gives
\begin{equation} 
|\tau_x(t) - \tau_x(t')| \geq C_0 M^{-\lambda(v)}, \text{ with $v$ as in \eqref{ycm of slopes}}. \label{slope-diff-lower} 
\end{equation}  
Comparing \eqref{slope-diff-upper} and \eqref{slope-diff-lower}, we obtain 
\begin{align*}  
&C_0 M^{-\lambda(v)} \leq |\tau_x(t) - \tau_x(t')| < \frac{2}{A_0} M^{-h(u)}, \; \text{ i.e. } \\ 
&M^{- \lambda(v)} < \frac{2}{A_0C_0} M^{-h(u)} < M^{-h(u)} \; \text{ or } \; h(u) < \lambda(v).
\end{align*}
The last step uses the fact that $A_0C_0 > 2$, which is consistent with the choice of the constant $A_0$ in \eqref{c0A0}. This establishes the desired conclusion \eqref{uv-relation}.  
\end{proof}

\begin{figure}[ht]
\centering

\begin{tikzpicture}[
    x=1cm,
    y=0.82cm,
    font=\small,
    bluebranch/.style={
        blue!75!black,
        line width=1.1pt
    },
    redbranch/.style={
        red!75!black,
        line width=1.1pt
    },
    guide/.style={
        gray!45,
        densely dashed,
        line width=0.45pt
    },
    ancestry/.style={
        gray!65,
        densely dotted,
        line width=0.65pt
    },
    vertex/.style={
        circle,
        fill=black,
        inner sep=1.65pt
    },
    bluevertex/.style={
        circle,
        fill=blue!75!black,
        inner sep=1.75pt
    },
    redvertex/.style={
        circle,
        fill=red!75!black,
        inner sep=1.75pt
    }
]

% -------------------------------------------------
% Adjustable coordinates
% -------------------------------------------------

% Horizontal positions
\def\heightx{0.0}
\def\rootx{4.0}
\def\slopex{8.8}

% Vertical levels
\def\titleheight{5.55}
\def\topheight{5.05}

% Unequally spaced M-adic heights
\def\hvheight{4.45}
\def\huheight{3.30}
\def\lambdaheight{2.35}

% Common base level for both trees
\def\leafheight{0.55}

% Horizontal spread of branches
\def\rootspread{1.15}
\def\slopespread{1.25}

% -------------------------------------------------
% Headings
% -------------------------------------------------

\node[font=\normalsize] at (\heightx,\titleheight)
    {$M$-adic height};

\node[font=\normalsize] at (\rootx,\titleheight)
    {Root};

\node[font=\normalsize] at (\slopex,\titleheight)
    {Slope};

% -------------------------------------------------
% Horizontal height guides
% -------------------------------------------------

\draw[guide]
    (0.95,\hvheight) -- (10.2,\hvheight);

\draw[guide]
    (0.95,\huheight) -- (10.2,\huheight);

\draw[guide]
    (0.95,\lambdaheight) -- (10.2,\lambdaheight);

% Height labels
\node[anchor=east] at (0.55,\hvheight)
    {$h(v)$};

\node[anchor=east] at (0.55,\huheight)
    {$h(u)$};

\node[anchor=east] at (0.55,\lambdaheight)
    {$\lambda(v)$};

% -------------------------------------------------
% Root tree
% -------------------------------------------------

\node[vertex] (u) at (\rootx,\huheight) {};

% Dotted continuation above u, ending at the common top height
\draw[ancestry]
    (u) -- (\rootx,\topheight);

% Descendants, both at the common base level
\node[bluevertex] (t)
    at ({\rootx-\rootspread},\leafheight) {};

\node[redvertex] (tp)
    at ({\rootx+\rootspread},\leafheight) {};

% Straight branches
\draw[bluebranch]
    (u) -- (t);

\draw[redbranch]
    (u) -- (tp);

% Labels
\node[above left=2pt] at (u)
    {$u$};

\node[
    below left=2pt,
    blue!75!black
] at (t)
    {$t$};

\node[
    below right=2pt,
    red!75!black
] at (tp)
    {$t'$};

% -------------------------------------------------
% Slope tree
% -------------------------------------------------

\node[vertex] (v) at (\slopex,\hvheight) {};

% Dotted continuation above v, ending at the same top height
\draw[ancestry]
    (v) -- (\slopex,\topheight);

% Descendants, also at the common base level
\node[bluevertex] (taut)
    at ({\slopex-\slopespread},\leafheight) {};

\node[redvertex] (tautp)
    at ({\slopex+\slopespread},\leafheight) {};

% Straight branches
\draw[bluebranch]
    (v) -- (taut);

\draw[redbranch]
    (v) -- (tautp);

% Labels
\node[above right=2pt] at (v)
    {$v$};

\node[
    below=4pt,
    blue!75!black
] at (taut)
    {$\tau_x(t)$};

\node[
    below=4pt,
    red!75!black
] at (tautp)
    {$\tau_x(t')$};

\end{tikzpicture}

\caption{\small{Schematic comparison of the root and slope trees in Lemma \ref{uv-lemma}.
The roots $t$ and $t'$ have youngest common ancestor $u$, while their
reference slopes $\tau_x(t)$ and $\tau_x(t')$ have youngest common
ancestor $v$, both in the $M$-adic tree. The dashed lines indicate the relevant $M$-adic heights $h(v)$, $h(u)$, and $\lambda(v)$. While $h(u)$ may be strictly larger than $h(v)$, it must be less than $\lambda(v)$.  
}}
\label{fig:h-and-lambda}
\end{figure}
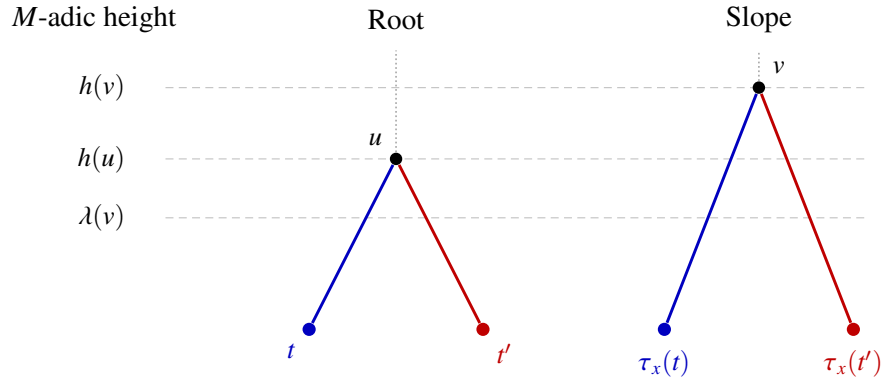

\section{From sets to trees: Proof of Lemma \ref{CONSISTENCY LEMMA}} \label{section: proof of consistency lemma} 
Suppose that $t, t'$ are two elements of ${\tt{Poss}}(x)$, whose $j$-long chains in ${\tt{TrPR}}_x$ and ${\tt{TrPS}}_x$ are given by  
\begin{align*} 
&{\tt{PossRoot}}_j(s) = \langle Q_1^{\ast}(s), \ldots, Q_{j}^{\ast}(s) \rangle,  \\
&{\tt{PosSlope}}_j(s) = \langle \theta_1(s), \ldots, \theta_{j}(s) \rangle,  \; \text{ for } s \in \{t, t'\}, 1 \leq j \leq N. 
\end{align*}  
Lemma \ref{CONSISTENCY LEMMA} asserts that a root interval of the form $Q_j^{\ast}(\cdot)$ can occur once in ${\tt{TrPR}}_x$, simultaneously fixing $j$, ${\tt{PossRoot}}_j( \cdot)$ and ${\tt{PosSlope}}_j( \cdot)$.   
%\noindent We are ready to embark on the proof of Lemma \ref{CONSISTENCY LEMMA}. The proof of each part will be handled in a separate subsection. 
\subsection{Proof of part \eqref{consistency part1}:}  We will prove this part through two lemmas. Under the assumption \eqref{Qjj'-equal}, Lemma \ref{slope equality lemma} below confirms the equality of the indices $j, j'$ and the equality of the slopes $\theta_k(t)$ and $\theta_k(t')$ for $1 \leq k \leq j$, while Lemma \ref{root equality lemma} does the same for the roots $Q_k^{\ast}(t)$ and $Q_k^{\ast}(t')$. Together, these two lemmas establish the claim \eqref{consistency check}. 
\vskip0.1in
\noindent 
\begin{lemma} \label{slope equality lemma} In the set-up described above, suppose that \eqref{Qjj'-equal} holds for some choice of indices $1 \leq j, j' \leq N$. Then the following conclusions hold:
\begin{equation} \label{j=j'}
j = j' \; \text{ and } \; \theta_j(t) = \theta_j(t'). 
\end{equation} 
In particular, this implies
\begin{equation} \label{k=k'}
 \theta_k(t) = \theta_k(t') \text{ for } 1 \leq k \leq j. 
\end{equation}  
\end{lemma} 
\begin{proof} 
Let us first check that \eqref{j=j'} implies \eqref{k=k'}. By Proposition \ref{PROP: COMPRESSED SLOPES} \eqref{generation-part-B_N}, any basic slope in ${\tt{BasicSl}}_j(\mathcal S_N)$ uniquely specifies a nested chain of basic slope ancestors in ${\tt{BasicSl}}_k$ for all $k \leq j$. Thus, equality of the slopes $\theta_j(t), \theta_j(t') \in {\tt{BasicSl}}_j$ implies equality of their ancestors at every preceding generation, proving \eqref{k=k'}.
\vskip0.1in
\noindent We turn our attention to \eqref{j=j'} next. The definition \eqref{Q-jtx} dictates that for every $s \in {\tt{Poss}}(x)$, corresponding elements of the two sequences ${\tt{PossRoot}}_N(s)$ and ${\tt{PosSlope}}_N(s)$ have the same length as $M$-adic intervals. Combined with the hypothesis \eqref{Qjj'-equal}, this implies 
\begin{equation} \label{theta=Q}  
|\theta_j(t)| = |Q_j^{\ast}(t)| = |Q_{j'}^{\ast}(t')| = |\theta_{j'}(t')|. 
\end{equation}  
Since $\theta_j(t)$ and $\theta_{j'}(t')$ are $M$-adic intervals of the same length, either they are identical, or their interiors are disjoint. The first case, when
\begin{align*} 
&\theta_j(t) = \theta_{j'}(t'),  \text{ occurs when } j=j', \text{since } \\ 
&\theta_j(t) \in {\tt{BasicSl}}_j, \;  \theta_{j'}(t') \in {\tt{BasicSl}}_{j'}, \text{ and } \\
&{\tt{BasicSl}}_j \cap {\tt{BasicSl}}_{j'} = \emptyset \text{ for } j \neq j' \text{ from the definition \eqref{jth basic slope cubes}}. 
\end{align*} 
This is precisely the desired conclusion \eqref{j=j'}. 
\vskip0.1in
\noindent In order to complete the proof of \eqref{j=j'}, one therefore has to eliminate the second possibility. For this, it suffices to show that $\theta_j(t) \cap \theta_{j'}(t')$ always contains an interval. Specifically, we will show that    
\begin{equation} \label{v in theta_j} 
v \subseteq \theta_j(t) \cap \theta_{j'}(t') \text{ where $v$ is the interval given by \eqref{ycm of slopes}}.
\end{equation} 
To prove \eqref{v in theta_j}, we show that $v \subseteq \theta_j(t)$; the proof of $v \subseteq \theta_{j'}(t')$ is identical. Let us note, in view of \eqref{ycm of slopes} and \eqref{theta-jtx}, that the $M$-adic intervals  $v$ and $\theta_j(t)$ both contain the slope $\tau_x(t)$. Therefore, one of the two (mutually exclusive) containment relations must hold: 
\begin{equation} 
\text{ either } v \subseteq \theta_j(t) \quad \text{ or } \quad \theta_j(t) \subsetneq v. 
\end{equation} Let us confirm that the last inclusion is impossible. Indeed, if $\theta_j(t) \subsetneq v$, then $\theta_j(t)$ must be a basic slope interval contained in one of the two basic slopes that are descended immediately from the splitting vertex $v \in {\tt{SplitV}}(\mathcal S_N)$. Since these latter basic slopes are of length $M^{-\lambda(v)}$, we deduce that  
\begin{equation}  \label{theta_j upper}
|\theta_j(t)| \leq M^{-\lambda(v)}. 
\end{equation}    
At the same time, $t$ and $t'$ are contained respectively in the intervals $Q_j^{\ast}(t)$ and $Q_{j'}^{\ast}(t')$. Therefore the assumption \eqref{Qjj'-equal} means  
\[ Q_j^{\ast}(t) = Q_{j'}^{\ast}(t')  \text{ contains both $t$ and $t'$}. \] Therefore it must contain their youngest common ancestor $u$ given by \eqref{ycm of roots}, i.e. $u \subseteq Q_j^{\ast}(t)$. In view of \eqref{theta=Q}, this means
\begin{equation} \label{theta_j lower}
 |\theta_j(t)| =  |Q_j^{\ast}(t)|  \geq |u| = M^{-h(u)}. 
\end{equation}  
Combining \eqref{theta_j upper} and \eqref{theta_j lower}, we obtain $\lambda(v) \leq h(u)$, which contradicts the conclusion \eqref{uv-relation} of Lemma \ref{uv-lemma}. 
\end{proof} 
\begin{lemma} \label{root equality lemma} 
Under the hypothesis \eqref{Qjj'-equal} of Lemma \ref{CONSISTENCY LEMMA}, one has the following additional set identities on the level of roots:
\begin{equation} \label{Qk=Qk'} 
Q_{k}^{\ast}(t) = Q_{k}^{\ast}(t') \; \text{ for } 1 \leq k \leq j. 
\end{equation}    
\end{lemma} 
\begin{proof} 
Let us recall from the definition \eqref{Q-jtx} that for $s \in \{t, t'\}$,
\[ s \subsetneq Q_{j}^{\ast}(s) \subsetneq \cdots \subsetneq Q_{1}^{\ast}(s), \text{ with } |Q_{k}^{\ast}(s)| = |\theta_{k}^{\ast}(s)| \text{ for all } 1 \leq k \leq j. \]
Further, in view of the conclusions \eqref{j=j'}, \eqref{k=k'} of Lemma \ref{slope equality lemma}, 
\[ |Q_{k}^{\ast}(t)| =  |\theta_{k}^{\ast}(t)| =  |\theta_{k}^{\ast}(t)| =  |Q_{k}^{\ast}(t')|. \]   
In other words, for $1 \leq k \leq j-1$, the $M$-adic intervals $Q_{k}^{\ast}(t)$ and $Q_k^{\ast}(t')$ are  of the same length, and they both contain the same interval $Q_j^{\ast}(t) = Q_j^{\ast}(t')$. This means that they must be equal, which is the intended conclusion \eqref{Qk=Qk'}. 
\end{proof} 
\subsection{Proof of part \eqref{consistency part2}:}
\begin{proof}
Part \eqref{consistency part1} of Lemma \ref{CONSISTENCY LEMMA}, which we have just proved, is key to the verification of the tree structure of ${\tt{TrPR}}_x$ and ${\tt{TrPS}}_x$. We use it to check that heights (i.e., generation levels) and ancestry relations are well-defined on each tree. Specifically, part \eqref{consistency part1} ensures that a vertex of the form 
\begin{align*} 
&{\tt{PossRoot}}_j(\cdot) = \langle Q_1, \ldots, Q_j \rangle \; \\  &(\text{respectively } {\tt{PossSlope}}_j(\cdot) = \langle \theta_1, \ldots, \theta_j \rangle)
\end{align*}  
can only occur in a single generation $j$ of the tree ${\tt{TrPR}}_x$ (respectively ${\tt{TrPS}}_x$), regardless of the argument $\cdot$ in ${\tt{Poss}}(x)$ that it originates from. Namely, suppose that a root interval $Q_j$ and an equi-dimensional slope interval $\theta_j \in {\tt{BasicSl}}_j$ occurs as the terminal vertex of two finite sequences 
\begin{align} 
&{\tt{PossRoot}}_{\ell}(t) \; \text{ and } \; {\tt{PossRoot}}_{\ell'}(t') \text{ in } {\tt{TrPR}}_x, \label{poss root pair} \\
&(\text{respectively } {\tt{PossSlope}}_{\ell}(t) \; \text{ and } \; {\tt{PossSlope}}_{\ell'}(t') \text{ in } {\tt{TrPS}}_x). \label{poss slope pair}
\end{align}
Then part \eqref{consistency part1} allows us to conclude that $\ell = \ell' = j$ and that the two sequences in \eqref{poss root pair} (respectively \eqref{poss slope pair}) are identical. This verifies that the notion of height of $Q_j$ in ${\tt{TrPR}}_x$ and the height of $\theta_j$ in ${\tt{TrPS}}_x$ are both well-defined.  
\vskip0.1in
\noindent Ancestry in either tree is determined by set containment. Specifically,  
\begin{align*}
&{\tt{PossRoot}}_{j+1}(t) = \langle Q_1, \ldots, Q_{j+1} \rangle = \langle Q_1^{\ast}(t), \ldots, Q_{j+1}^{\ast}(t) \rangle \text{ is a child of } \\
&{\tt{PossRoot}}_{j}(t') =  \langle Q_1', \ldots, Q_j'\rangle	= \langle Q_1^{\ast}(t'), \ldots, Q_j^{\ast}(t') \rangle \\ 
& \quad \text{ if } Q_k = Q_k' \text{ for } k \leq j\, \text{ and } Q_{j+1} \subsetneq Q_j.
\end{align*} 
Similarly, 
\begin{align*} 
&{\tt{PossSlope}}_{j+1}(t) = \langle \theta_1, \ldots, \theta_{j+1} \rangle = \langle \theta_1(t), \ldots, \theta_{j+1}(t) \rangle \text{ is a child of } \\ 
&{\tt{PossSlope}}_j(t') = \langle \theta_1', \ldots, \theta_j' \rangle = \langle \theta_1(t'), \ldots, \theta_j(t') \rangle \\ 
& \quad \text{ if } \theta_k = \theta_k' \text{ for } k \leq j \text{ and } \theta_{j+1} \subsetneq \theta_j. 
\end{align*}  
Once again, part \eqref{consistency part1} shows that even though a non-terminal vertex of ${\tt{TrPR}}_x$ (respectively ${\tt{TrPS}}_x$) may originate from distinct $t, t' \in {\tt{Poss}}(x)$, it has a unique parent in the respective trees regardless of $t, t'$. This verifies the consistency of the ancestry relation in both trees and the uniqueness of a ray rooted at the origin ending at a given vertex. This completes the verification of the tree structure of ${\tt{TrPR}}_x$ and ${\tt{TrPS}}_x$, and therefore the proof of part \eqref{consistency part2}.
\end{proof}

\section{From trees to sets: Proof of Lemma \ref{TREE SET LEMMA}} \label{section: tree set lemma proof}
We now turn to the identification of the geometric sets that the trees ${\tt{TrPR}}_x$ and ${\tt{TrPS}}_x$ represent. The goal is to establish a bijection between maximal rays of ${\tt{TrPR}}_x$ (respectively ${\tt{TrPS}}_x$) and elements of ${\tt{Poss}}(x)$ (respectively $\tau_x({\tt{Poss}}_x)$). 

\subsection{Proof of part \eqref{consistency part2b}}  We start with ${\tt{TrPS}}_x$ first. 
\vskip0.1in
\noindent It follows from Lemma \ref{CONSISTENCY LEMMA} \eqref{consistency part2} that ${\tt{TrPS}}_x$ is a subtree of $\mathscr{B}_N$, the compressed version of the pruned slope tree $\mathcal S_N$ given by Proposition \ref{PROP: COMPRESSED SLOPES}. According to part \eqref{B_N identification} of this proposition, each terminal vertex $\theta_N(t)$ of ${\tt{TrPS}}_x$, which also corresponds to a terminal vertex of a maximal ray in $\mathscr{B}_N$, contains exactly one slope in $\Omega_N$. This slope has to be $\tau_x(t)$ by virtue of \eqref{PossRoot chains}. 
\vskip0.1in
\noindent Conversely, every element $\omega \in \tau_x({\tt{Poss}}(x))$ is of the form $\omega = \tau_x(t)$ for some $t \in {\tt{Poss}}(x)$, by Lemma \ref{CONSISTENCY LEMMA} \eqref{consistency part2}. Again by Proposition \ref{PROP: COMPRESSED SLOPES} there is a unique basic slope interval $\theta_N = \theta_N(t) \in {\tt{BasicSl}}_N$ that contains $\tau_x(t) \in \Omega_N$. The maximal ray of $\mathscr{B}_N$ that identifies $\theta_N$, and therefore identifies $\omega \in \tau_x({\tt{Poss}}(x))$, is also a maximal ray of ${\tt{TrPS}}_x$ given by ${\tt{PossSlope}}_N(t)$. 
%In other words, ${\tt{TrPS}}_x$ determines $\tau_x({\tt{Poss}}(x))$.

\subsection{Proof of Lemma \ref{TREE SET LEMMA}\eqref{consistency part2a}} We consider ${\tt{TrPR}}_x$ next. 
\vskip0.1in
\noindent Every terminal vertex $Q_N$ of the tree ${\tt{TrPR}}_x$ is of the form $Q_N^{\ast}(t)$ for some $t \in {\tt{Poss}}(x)$, according to Lemma \ref{CONSISTENCY LEMMA} \eqref{consistency part2}. We argue that the choice of $t$ is unique. If $t, t' \in {\tt{Poss}}(x)$ are such that $Q_N^{\ast}(t) = Q_N^{\ast}(t')$, then Lemma \ref{CONSISTENCY LEMMA} \eqref{consistency part1} dictates  
\begin{equation}  \label{thetaNtt'}
\theta_N(t) = \theta_N(t'), \text{ with } \tau_x(s) \in \theta_N(s), \; s \in \{t, t'\}.
\end{equation} 
Since each basic slope interval contains a unique element of $\Omega_N$ by our pruning process, \eqref{thetaNtt'} implies 
\[  \tau_x(t) = \tau_x(t'). \] The injectivity Lemma \ref{lemma: injectivity} then allows us to conclude $t = t'$. Thus the ray determined by $Q_N$ identifies a unique $t \in {\tt{Poss}}(x)$.
\vskip0.1in
\noindent Conversely, every $t \in {\tt{Poss}}(x)$ is contained in a unique nested chain ${\tt{PossRoot}}_N(t)$, which represents a terminal vertex and hence a maximal ray of ${\tt{TrPR}}_x$. That means ${\tt{Poss}}(x)$ is a subset of the set represented by ${\tt{TrPR}}_x$. 
 
\section{Lifting the reference map to trees: Proof of Lemma \ref{LEMMA: MAP LIFT}} \label{proof section: lift map} 
Lemma \ref{LEMMA: MAP LIFT} asserts that the reference map $\tau_x: {\tt{Poss}}(x) \rightarrow \tau_x\bigl({\tt{Poss}}(x)\bigr)$ can be lifted to a length-preserving sticky map 
\[ \tilde{\tau}_x : {\tt{TrPR}}_x \rightarrow {\tt{TrPS}}_x \] among trees, according to the definition
\begin{align} 
&\tilde{\tau}_x(Q_j^{\ast}(t)) := \theta_j(t), \; \text{ for } 1 \leq j \leq N,  \text{ or said differently, } \label{tau_x tree map} \\ 
&\tilde{\tau}_x \bigl(\langle Q_1^{\ast}(t), \ldots, Q_j^{\ast}(t)\rangle \bigr) = \langle \theta_1(t), \ldots, \theta_j(t)\rangle. \nonumber 
\end{align}
Here $Q_j^{\ast}(t)$ and $\theta_j(t)$ are the root and the slope intervals containing $t$ and $\tau_x(t)$ respectively, as given by \eqref{Q-jtx} and \eqref{theta-jtx}
\vskip0.1in
\noindent We distinguish between the set map $\tau_x$ and the tree map $\tilde{\tau}_x$ for the purpose of the proof. Once the lemma is established, we will denote the tree map $\tilde{\tau}_x$ by $\tau_x$ going forward. 
\begin{proof}
First of all, the map $\tilde{\tau}_x$ given by \eqref{tau_x tree map} is well-defined. Namely, if $Q_{j}^{\ast}(t) = Q_{j'}^{\ast}(t')$, then Lemma \ref{CONSISTENCY LEMMA} \eqref{consistency part1} implies $j = j'$, and therefore 
 \[ \tilde{\tau}_x(Q_{j}^{\ast}(t)) = \theta_j(t) = \theta_j(t') = \tilde{\tau}_x(Q_{j'}^{\ast}(t')).  
\]
\vskip0.1in
\noindent Second, $\tilde{\tau}_x$ preserves height, since both intervals $Q_j^{\ast}(t)$ and $\theta_j(t)$ are of height $j$ in their respective trees ${\tt{TrPR}}_x$ and ${\tt{TrPS}}_x$. This is one of the defining criteria for a sticky map, as seen in Section \ref{sticky maps section}. 
\vskip0.1in
\noindent Third, $\tilde{\tau}_x$ preserves lineage. If $Q_{j+1}$ is a child of $Q_j$ in ${\tt{TrPR}}_x$, the ancestry relation in the tree implies that there exists a root interval $t \in {\tt{Poss}}(x)$ such that 
\[ Q_{j+1} = Q_{j+1}^{\ast}(t) \subsetneq Q_j = Q_j^{\ast}(t). \] 
The tree structure of ${\tt{TrPS}}_x$, proved in Lemma \ref{CONSISTENCY LEMMA} \eqref{consistency part2}, then confirms that 
\[  \tilde{\tau}_x \bigl(Q_{j+1}^{\ast}(t) \bigr) = \theta_{j+1}(t) 
\text{ is a child of } \theta_j(t) = \tilde{\tau}_x \bigl(Q_{j}^{\ast}(t) \bigr). \] This shows that $\tilde{\tau}_x$ is sticky. 
\vskip0.1in
\noindent The defining condition \eqref{Q-jtx} states that $Q_j^{\ast}(t)$ and $\theta_j(t)$ are of the same length; thus, $\tilde{\tau}_x$ preserves length.
\vskip0.1in
\noindent The map $\tilde{\tau}_x$ is a tree extension of $\tau_x$ in the following sense. Consider an $N$-long chain $\langle Q_1, \ldots, Q_N \rangle$ that occurs as a terminal vertex of ${\tt{TrPR}}_x$. It follows from Lemma \ref{TREE SET LEMMA} \eqref{consistency part2a} that there is a unique root interval $t \in {\tt{Poss}}(x)$ such that $Q_N = Q_N^{\ast}(t)$. The defining relation \eqref{tau_x tree map} dictates that  
\[ \tilde{\tau}_x(Q_N^{\ast}(t)) = \theta_N(t), \text{ which is a terminal vertex of } {\tt{TrPS}}_x. \] According to Lemma \ref{TREE SET LEMMA} \eqref{consistency part2b}, the vertex $\theta_N(t)$ uniquely identifies $\tau_x(t)$. In this way, the final $N^{\text{th}}$ level of tree-based map $\tilde{\tau}_x$ as defined in \eqref{tau_x tree map} leads to the definition \eqref{what is tau_x} of the set-based map $\tau_x$. 
\vskip0.1in
\noindent Finally, for every $t \in {\tt{Poss}}(x)$, let us define the binary quantity $Z_x(Q_j^{\ast}(t))$ as follows, 
\[ Z_x \bigl( Q_j^{\ast}(t)\bigr) = \begin{cases} 0 &\text{ if } \tilde{\tau}_x(Q_j^{\ast}(t)) \text{ is the $0^{\text{th}}$ child of } \tilde{\tau}_x(Q_{j-1}^{\ast}(t)), \\ 
1 &\text{ otherwise.} \end{cases} \] 
Let us verify that $Z_x(Q)$ is unambiguously defined for each $Q \in \mathcal V({\tt{TrPR}}_x)$. Once again, Lemma \ref{CONSISTENCY LEMMA} \eqref{consistency part1} provides the necessary ingredient. Indeed 
\[ \text{ if } Q = Q_{j}^{\ast}(t) = Q_{j'}^{\ast}(t') \text{ for some } t, t' \in {\tt{Poss}}(x), \] then the lemma implies \[j = j' \text{ and } \theta_k(t) = \theta_k(t') = \theta_k \; \text{ for all $1 \leq k \leq j$}. \] This yields the same value of $Z_x(Q)$ regardless of the choice $t, t'$. This completes the proof.   
\end{proof}
\section{Uniformity of progeny in ${\tt{TrPR}}_x$: Proof of Lemma \ref{LEMMA: UNIFORMITY OF PROGENY}}  \label{proof section: uniformity of progeny} 
Lemma \ref{lemma: injectivity} claims injectivity of the map $\tau_x$ on the set ${\tt{Poss}}(x)$. Since ${\tt{Poss}}(x)$ is a subset of $\mathcal Q(J)$, every root in it is of the finest scale, i.e. $|t| =M^{-J}$ for all $t \in {\tt{Poss}}(x)$. Having lifted $\tau_x$ to a sticky map from ${\tt{TrPR}}_x$ to ${\tt{TrPS}}_x$, it is natural to ask to what extent this injectivity persists at the coarser scales appearing in the heights of these trees. Quantifying this is the goal of Lemma \ref{LEMMA: UNIFORMITY OF PROGENY}. We prove Lemma \ref{LEMMA: UNIFORMITY OF PROGENY} in this section. 
\vskip0.1in
\noindent In preparation, we use the geometric intersection criterion \eqref{intersection criterion inequality} to study the structure of
\[ \tau_x^{-1}(\theta) = \text{ the pre-image of $\theta$ under $\tau_x$,} \quad \text{ for $\theta \in {\tt{BasicSl}}_j$, $1 \leq j \leq N$. } \] 
In view of Lemma \ref{LEMMA: MAP LIFT},  every root verted in $\tau_x^{-1}(\theta)$ is an $M$-adic interval of length $|\theta|$. While $\tau_x$ may not be injective at a non-terminal height, in other words $\tau_x^{-1}(\theta)$ need not be a singleton for $j < N$, we conclude in Lemma \ref{lemma: tree injectivity} below that the intervals in $\tau_x^{-1}(\theta)$ are constrained to lie within a single interval of length comparable to $|\theta|$; thus, the cardinality of $\tau_x^{-1}(\theta)$  is uniformly controlled by an absolute constant. This fact will feature prominently in the proof of Lemma \ref{LEMMA: UNIFORMITY OF PROGENY}, which appears in the next subsection.   
\begin{lemma} \label{lemma: tree injectivity}
For $1 \leq j \leq N$, let $\theta$ be a $j^{\text{th}}$ level basic slope interval in $\mathcal S_N$. If $t, t' \in {\tt{Poss}}(x)$ have the property that $\tau_x(t), \tau_x(t') \in \theta$, then 
\begin{equation} \label{distance between centres}
|\text{cen}(t) - \text{cen}(t')| \leq (A_0+2) |\theta|. 
\end{equation} 
In other words, there is a subinterval $J(\theta) \subseteq [0,1]$ of length at most $(A_0+2)|\theta|$ such that 
\begin{equation} \label{J-theta}
\left\{ \text{cen}(t) : \tau_x(t) \in \theta, t \in {\tt{Poss}}(x)  \right\} \subseteq J(\theta). 
\end{equation} 
\end{lemma}
\begin{proof} 
The definition of the reference map $\tau_x$ implies that 
\[ x \in {\tt{Tube}}(t, \tau_x(t), I_0) \cap {\tt{Tube}}(t', \tau_x(t'), I_0). \]
The intersection condition  \eqref{intersection criterion inequality} for tubes then yields 
\[ 
\bigl| \text{cen}(t') - \text{cen}(t) + x_1 \bigl(\tau_x(t') - \tau_x(t) \bigr) \bigr|  \leq c_0 M^{-J}. 
\]
According to the hypothesis, $\tau_x(t)$ and $\tau_x(t')$ both lie in a single interval $\theta$, therefore the Euclidean distance between is at most $\text{diam}(\theta) = |\theta|$. By the triangle inequality, this leads to 
\begin{align*} 
\bigl| \text{cen}(t') - \text{cen}(t) \bigr| &\leq c_0 M^{-J} + |x_1 \bigl(\tau_x(t') - \tau_x(t) | \\ 
&\leq c_0 M^{-J} + (A_0 + 1) |\theta| \leq (A_0+2)|\theta|.  
\end{align*} 
At the last step, we have used the fact that $c_0 M^{-J} < M^{-J} < |\theta| = M^{-\lambda(\cdot)}$, since the pruning process of the slope set ensures that the finest scale $J$ exceeds all the fundamental heights. The inclusion \eqref{J-theta} is a direct consequence of \eqref{distance between centres}. 
\end{proof} 
\subsection{Proof of Lemma \ref{LEMMA: UNIFORMITY OF PROGENY}}
The first conclusion \eqref{uniformity of progeny} follows directly from \eqref{J-theta}. Suppose that $\theta \in {\tt{BasicSl}}_j$ for some index $j$ with $1 \leq j \leq N$. Then 
\begin{align*} 
\tau_x^{-1}(\theta) &= \bigl\{Q \in \mathcal V_j \bigl( {\tt{TrPR}}_x\bigr) : \tau_x(Q) = \theta, \; |Q| = |\theta|  \bigr\} \\
&= \left\{Q = Q_j^{\ast}(t) \; \Bigl| \; \begin{aligned} &t \in {\tt{Poss}}(x), \; t \subseteq Q_j^{\ast}(t), \; \tau_x(t) \in \theta \\ &\tau_x \bigl(Q_j^{\ast}(t)\bigr) = \theta, \; |Q_j^{\ast}(t)| = |\theta| \end{aligned} \right\}  \\ 
&\subseteq \left\{Q = Q_j^{\ast}(t) \; \Bigl| \; \begin{aligned} &t \subseteq Q_j^{\ast}(t), \; \text{cen}(t) \subseteq J(\theta) \\  &t \in {\tt{Poss}}(x), \; |Q_j^{\ast}(t)| = |\theta| \end{aligned} \right\},
\end{align*} 
where $J(\theta)$ is the interval of length at most $(A_0+2)|\theta|$ specified by \eqref{J-theta}. Since each $Q = Q_j^{\ast}(t) \in \tau_x^{-1}(\theta)$ contains ${\text{cen}}(t)\in J(\theta)$,
the number of possible $M$-adic intervals $Q$ of length $|\theta|$ is at most $C_2 = \lceil A_0+2 \rceil + 1$. This proves \eqref{uniformity of progeny}.
\vskip0.1in
\noindent Let us now use \eqref{uniformity of progeny} to show that each non-terminal vertex of ${\tt{TrPR}}_x$ has at most $2C_2$ children. Suppose that $v = \langle Q_1, \ldots Q_j \rangle$ is a vertex of ${\tt{TrPR}}_x$ of generation $j$, $1 \leq j < N$. It follows from \eqref{level-j} that there exists $t \in {\tt{Poss}}(x)$ such that  
\[ v = {\tt{PossRoot}}_j(t), \; \text{ i.e., } \;  Q_k = Q_k^{\ast}(t) \quad \text{ for } 1 \leq k \leq j.  \] 
Let $\theta_j(t) \in {\tt{BasicSl}}_j$ be the $j^{\text{th}}$ basic interval containing $\tau_x(t)$, so that 
\[ \tau_x(Q_j) = \tau_x(Q_j^{\ast}(t)) = \theta_j(t), \text{ in view of Lemma \ref{LEMMA: MAP LIFT}}. \] It follows from the property of basic slopes (Section \ref{section: basic slopes}) that $\theta_j(t)$ has exactly two descendants in $\mathcal S_N$, say $\vartheta_1$ and $\vartheta_2$, that belong to ${\tt{BasicSl}}_{j+1}$. The stickiness of the map $\tau_x$ implies that the children of $Q_j^{\ast}(t)$ must map into either $\vartheta_1$ or $\vartheta_2$ under $\tau_x$, i.e., 
\[ 
\text{the children of } {\tt{PossRoot}}_j(t) \text{ are contained in } \tau_x^{-1}(\vartheta_1) \cup \tau_x^{-1}(\vartheta_2). \]
As a result, 
\begin{equation} \# \left\{ \text{children of } {\tt{PossRoot}}_j(t) \text{ in} {\tt{TrPR}}_x  \right\} \leq \sum_{\ell=1}^2 \tau_x^{-1}(\vartheta_\ell) \leq 2C_2,  \label{uniformity of progeny2}
\end{equation} 
where the last step follows from \eqref{uniformity of progeny}. 
\vskip0.1in 
\noindent Finally, we use \eqref{uniformity of progeny2} to prove \eqref{cardinality of V_j}. According to Lemma \ref{LEMMA: MAP LIFT}, every vertex of ${\tt{TrPR}}_x$ maps into a vertex of ${\tt{TrPS}}_x$ of the same generation under the tree map $\tau_x$. This leads to the set inclusion
\[ \mathcal V_j ({\tt{TrPR}}_x) = \bigcup_{\theta} \left\{ \tau_x^{-1}(\theta) : \theta \in \mathcal V_j ({\tt{TrPS}}_x)  \right\} \subseteq \bigcup_{\theta} \left\{ \tau_x^{-1}(\theta) : \theta \in {\tt{BasicSl}}_j \right\}. \]
From this we deduce that   
\begin{align*} 
\# \left[ \mathcal V_j ({\tt{TrPR}}_x) \right] &\leq \sum_{\theta} \left\{ \# \left[ \tau_x^{-1}(\theta) \right] : \theta \in {\tt{BasicSl}}_j \right\} \\
& \leq \sup_{\theta} \# \left[ \tau_x^{-1}(\theta) \right] \times \left[ \#(\tt{BasicSl}_j)\right] \leq C_22^j,
\end{align*} 
where we have used \eqref{uniformity of progeny} and \eqref{basic slope count j} at the final step.  This establishes \eqref{cardinality of V_j}, and completes the proof of Lemma \ref{LEMMA: UNIFORMITY OF PROGENY}. 

\section{Proof of Lemma \ref{PROP: REFERENCE TREES}} 
\label{proof: reference trees}
This lemma provides two characterizations for the inclusion event $x \in \mathtt K_{\mathbb X}$. The first characterization \eqref{sigma = tau} asserts a common maximal ray between ${\tt{TrPR}}_x$ and $\mathscr{U}_{\mathbb X}$ on which $\tau_x, \sigma_{\mathbb X}$ coincide. The second one \eqref{XZ} translates \eqref{sigma = tau}  into the equality of two binary sequences.
\subsection{Proof of part \eqref{consistency part5}} \label{inclusion equivalence proof section} 
\begin{proof} 
%In this part of the proposition, we are tasked with proving the following equivalence: for $x \in I_0 \times \mathbb R$, 
%\begin{equation}  \label{inclusion-equivalence-1}
%x \in \mathtt K_{\mathbb X} \Longleftrightarrow \eqref{sigma = tau}  \Longleftrightarrow \eqref{XZ}. 
%\end{equation}  
Let us recall  from \eqref{defn: extended Kakeya set K} and \eqref{tube family X} that $\mathtt K_{\mathbb X}$ is the union of tubes of the form ${\tt{Tube}}_{\mathbb X}[t]$, as $t$ ranges over $\mathcal Q(J)$. This means that 
\[ \mathtt K_{\mathbb X} \cap \bigl[ I_0 \times \mathbb R \bigr] = \bigcup \bigl\{ {\tt{Tube}}(t, \sigma_{\mathbb X}(t), I_0) : t \in \mathcal Q(J) \bigr\} \; \text{ where } I_0 = [A_0, A_0+1]. \] 
Therefore $x \in \mathtt K_{\mathbb X} \cap \bigl[ I_0 \times \mathbb R \bigr]$ if and only if 
\begin{equation} \label{special t} 
\text{there exists } t \in \mathcal Q(J) \text{ such that } x \in {\tt{Tube}}(t, \sigma_{\mathbb X}(t), I_0). 
\end{equation}  
Thus, proving Lemma \ref{PROP: REFERENCE TREES} \eqref{consistency part5} reduces to showing that 
\begin{equation}  \label{inclusion-equivalence-2}
\eqref{special t} \Longleftrightarrow \eqref{sigma = tau}. 
\end{equation} 
%We will henceforth focus on establishing the implications in \eqref{inclusion-equivalence-2}. 
\vskip0.1in
\noindent Suppose \eqref{special t} holds. Then it follows from the definition \eqref{defn: Poss(x)} that $t \in {\tt{Poss}}(x)$. On the other hand, Corollary \ref{cor: uniqueness of possible slopes} confirms that there is a unique slope $\omega = \tau_x(t)$ such that $x \in {\tt{Tube}}(t, \omega, I_0)$. This means
\begin{equation} \label{special t equiv}  
\omega = \tau_x(t)  = \sigma_{\mathbb X}(t) \in \Omega_N. 
%\text{ let us call this slope } \omega. 
\end{equation}  
Let $\theta_j$ denote the unique interval  in ${\tt{BasicSl}}_j$ that contains $\omega$, $1 \leq j \leq N$. 
\vskip0.1in
\noindent By Lemma \ref{TREE SET LEMMA} \eqref{consistency part2a}, $t$ is determined by a maximal ray $\mathcal R \in \partial {\tt{TrPR}}_x$; specifically, 
\[ \text{ if } t \in {\tt{Poss}}(x) \text{ is associated with } {\tt{PossRoot}}_N(t) = \langle Q_1^{\ast}(t), \ldots, Q_N^{\ast}(t) \rangle, \] 
then $\mathcal R$ is given by
\begin{align} 
&\mathcal R: [0,1] \rightarrow Q_1^{\ast}(t) \rightarrow Q_2^{\ast}(t) \rightarrow \cdots \rightarrow Q_N^{\ast}(t), \text{ where } \nonumber \\ 
&t \subsetneq Q_N^{\ast}(t) \subsetneq \cdots \subsetneq Q_1^{\ast}(t), \quad \tau_x(Q_j^{\ast}(t)) = \theta_j, \quad |Q_j^{\ast}(t)| = |\theta_j|.   
\label{tin-1}
\end{align} 
\vskip0.1in
\noindent At the same time, Proposition \ref{PROP: COMPRESSED ROOTS} asserts that  $\mathscr{U}_{\mathbb X}$ is a compressed tree representation of $[0,1]$, and $\sigma_{\mathbb X}: \mathscr{U}_{\mathbb X} \rightarrow \mathscr{B}_N$ is a sticky, length-preserving tree map. This means that there is a unique maximal ray $\mathcal R' \in \partial \mathscr{U}_{\mathbb X}$ identifying $t$,   
\begin{align}  
&\mathcal R': [0,1] \rightarrow Q_1 \rightarrow Q_2 \rightarrow \cdots \rightarrow Q_N, \quad t \subsetneq Q_N \subsetneq Q_{N-1} \subsetneq \cdots \subsetneq Q_1, \label{tin-2}  \\
 \label{tin-3} 
&{\text{which also satisfies }} 
\sigma_{\mathbb X}(Q_j) = \theta_j, \; \text{ with } \; |Q_j| = |\theta_j|, \text{ by \eqref{special t equiv}}.  
\end{align}  
Comparing the inclusion and size  conditions \eqref{tin-1}, \eqref{tin-2} and \eqref{tin-3}, we find that $Q_j$ and $Q_j^{\ast}(t)$ are equi-dimensional $M$-adic intervals, both containing $t$. This means
\[ Q_j^{\ast}(t) = Q_j \text{ for all } 1 \leq j \leq N, \text{ and hence } \mathcal R = \mathcal R', \text{ with } \sigma_{\mathbb X} \equiv \tau_x \text{ on } \mathcal R.  \] 
This implies the desired conclusion \eqref{sigma = tau}. 
\vskip0.1in
\noindent The converse implication is also true; namely \eqref{sigma = tau} implies \eqref{special t}. Suppose that ${\tt{TrPR}}_x$ and $\mathscr{U}_{\mathbb X}$ share a common maximal ray $\mathcal R$. This common ray identifies a unique point $t \in {\tt{Poss}}(x)$ by virtue of Lemma \ref{TREE SET LEMMA} \eqref{consistency part2a}. Moreover, the assumption \eqref{sigma = tau} says that $\sigma_{\mathbb X}$ and $\tau_x$ coincide on every vertex of this ray. If $\hat{\mathcal R}$ denotes the image of $\mathcal R$ under these maps, this means 
\[ \hat{\mathcal R} := \tau_x(\mathcal R) =\sigma_{\mathbb X}(\mathcal R), \text{ so that $\hat{\mathcal R}$ is a maximal ray in } {\tt{TrPS}}_x \subseteq \mathscr{B}_N. \] 
Comparing the terminal vertices on the rays in the above display, we obtain 
\[ \tau_x \bigl(Q_N^{\ast}(t) \bigr) = \sigma_{\mathbb X}(Q_N^{\ast}(t)) = \theta_N(t) \in {\tt{BasicSl}}_N, \; t \subseteq Q_N^{\ast}(t).   \] 
As a basic slope interval of the finest scale, $\theta_N(t)$ contains a unique slope $\omega \in \Omega_N$. On one hand, Lemmas \ref{TREE SET LEMMA} and \ref{LEMMA: MAP LIFT} imply that this slope $\omega$ must equal $\tau_x(t)$. On the other, $t$ is a subset of $Q_N^{\ast}(t)$, therefore by the definition of $\sigma_{\mathbb X}$ prescribed in Section \ref{section: finest scale extension}, $\omega = \sigma_{\mathbb X}(t)$. Combining these two observations with the definition \eqref{what is tau_x} of $\tau_x$, we arrive at the conclusion
\[ x \in {\tt{Tube}}(t, \omega, I_0) = x \in {\tt{Tube}}(t, \tau_x(t), I_0)  = {\tt{Tube}}(t, \sigma_{\mathbb X}(t), I_0). \] 
In other words, \eqref{special t} holds. This completes the proof of \eqref{consistency part5}.
\end{proof} 

%On one hand, $Q_N^{\ast}(t)$ is the terminal vertex of ${\tt{TrPR}}_x$, and therefore $\theta_N(t) = \tau_x \bigl(Q_N^{\ast}(t) \bigr)$ is a terminal vertex of ${\tt{TrPS}}_x$. By Lemma \ref{TREE SET LEMMA} \eqref{consistency part2b}, $\theta_N(t)$ contains a unique slope $\omega \in \Omega_N$. On the other hand, $\sigma_{\mathbb X}: \mathscr{U}_{\mathbb X} \rightarrow \mathcal S_N$ is also a sticky tree map,  therefore as a terminal vertex of $\mathcal S_N$, the slope interval $\theta_N(t) = \sigma_{\mathbb X}(Q_N^{\ast}(t))$ also contains   
%As we observed in Section \ref{section: finest scale extension} and in Lemma \ref{TREE SET LEMMA} \eqref{consistency part2b}, the finest slope interval $\theta_N(t)$ in ${\tt{TrPS}}_x \subseteq \mathcal S_N$ uniquely identifies a slope $\omega \in \Omega_N$ such that $\omega \in \theta_N(t)$. This means $\sigma_{\mathbb X}(t) = \tau_x(t) = \omega$. 
%%Thus the root $t$ and the slope $\omega$ obey \eqref{special t equiv}. 
%It now follows from the definition of $\tau_x$ that 

 \subsection{Proof of part \eqref{consistency part6}} 
 \begin{proof} 
 Here we need to show that 
 \begin{equation} 
 \eqref{sigma = tau} \Longleftrightarrow \eqref{XZ}.
 \end{equation} 
 \vskip0.1in
\noindent For the forward implication, where we assume \eqref{sigma = tau}, let $\mathcal R \in \partial {\tt{TrPR}}_x \cap \partial \mathscr{U}_{\mathbb X}$ be the ray identifying $t \in {\tt{Poss}}(x)$ where $\sigma_{\mathbb X}$ and $\tau_x$ coincide. If $Q_j = Q_j^{\ast}(t)$ and $\theta_j$ denote respectively the $j^{\text{th}}$ vertex on $\mathcal R$ and $\hat{\mathcal R} = \tau_x(\mathcal R) = \sigma_{\mathbb X}(\mathcal R)$, then for each $Q = Q_j^{\ast}(t)$, $1 \leq j \leq N$, we have
\begin{align*} 
\sigma_{\mathbb X}(Q_j) &= X(Q_j)^{\text{th}} \text{ child of } \theta_{j-1} 
\text{ from \eqref{defn: mapX} of Proposition \ref{PROP: COMPRESSED ROOTS}, whereas} \\ 
\tau_x(Q_j) &=  Z_x(Q_j)^{\text{th}} \text{ child of } \theta_{j-1} \text{ from \eqref{Z_x1} and \eqref{Z_x2} in Lemma \ref{LEMMA: MAP LIFT}.}
\end{align*}  
The equality $\sigma_{\mathbb X}(Q_j) = \tau_x(Q_j)$ implies $X(Q_j) = Z_x(Q_j)$, proving the claim \eqref{XZ}. 
\vskip0.1in 
\noindent The reverse implication \eqref{XZ} $\implies$ \eqref{sigma = tau} is essentially a backward retracing of the above steps. Suppose that the binary collection $\mathbb X$ given by \eqref{binary X} obeys the condition \eqref{XZ}. Let $t$ be the root in ${\tt{Poss}}(x)$ given by this hypothesis, and let $\mathcal R \in \partial {\tt{TrPR}}_x$ be the maximal ray representing $t$, in the sense of Lemma \ref{TREE SET LEMMA}.
 \vskip0.1in 
\noindent According to Proposition \ref{PROP: COMPRESSED ROOTS}, $(\mathscr{U}_{\mathbb X}, \sigma_{\mathbb X})$ is a tree-map pair adapted to $\mathscr{B}_N$. The defining property of $\mathscr{U}_{\mathbb X}$ asserts that a ray $\mathcal R'$ given by 
\[ \mathcal R': Q_0 \rightarrow Q_1 \rightarrow Q_2 \cdots \rightarrow Q_N, \quad Q_{N} \subsetneq Q_{N-1} \subsetneq \cdots \subsetneq Q_1 \subseteq Q_0 = [0,1], \]
is a maximal ray of $\mathscr{U}_{\mathbb X}$ if and only if for each $1 \leq j \leq N$, 
\begin{equation}  |Q_j| = |\sigma_{\mathbb X}(Q_j)| \; \text{ and } \; \sigma_{\mathbb X}(Q_{j}) = X(Q_{j})^{\text{th}} \text{ child of } \sigma_{\mathbb X}(Q_{j-1}).  
\end{equation}   
The assumption 
\begin{equation} \label{X=Z}  
X(Q) = Z_x(Q) \; \text{ for all } \; Q = Q_j^{\ast}(t), \; 1 \leq j \leq N, 
\end{equation}
and the definition \eqref{Z_x2} of $Z_x(Q)$ in relation to the reference map $\tau_x$ jointly imply that the entries of ${\tt{PossRoot}}_N(t) = \langle Q_1^{\ast}(t), \ldots, Q_N^{\ast}(t)\rangle$ obey this defining property with $Q_j = Q_j^{\ast}(t)$. In other words, $\mathcal R \in \partial \mathscr{U}_{\mathbb X}$, which is part of the desired conclusion \eqref{sigma = tau}.  
\vskip0.1in
\noindent Once we have established that $\mathcal R = \mathcal R' \in \partial {\tt{TrPR}}_x \cap \partial \mathscr{U}_{\mathbb X}$, the agreement of $\sigma_{\mathbb X}$ and $\tau_x$ on $\mathcal R$ follows from \eqref{X=Z}, since the sequence
\begin{align*} 
&\bigl\{X(Q) : Q = Q_j^{\ast}(t), 1 \leq j \leq N \bigr\} \text{ determines } \sigma_{\mathbb X} \text{ on $\mathcal R$, whereas } \\ 
&\bigl\{Z_x(Q) : Q = Q_j^{\ast}(t), 1 \leq j \leq N \bigr\}  \text{ determines }  \tau_x \text{ on $\mathcal R$}. 
\end{align*} This completes the proof of \eqref{sigma = tau}. 
%To establish \eqref{sigma = tau}, we need to confirm that the maximal ray $\mathcal R \in \partial {\tt{TrPR}}_x$ given by ${\tt{PossRoot}}_N(t)$ also occurs in the compressed tree $\partial \mathscr{U}_{\mathbb X}$. Once we establish this, the equality of $\sigma_{\mathbb X}$ and $\tau_x$ on $\mathcal R$ follows directly from the equality $X(Q) = Z_x(Q)$ for all vertices $Q$ on $\mathcal R$.  
%\vskip0.1in
%\noindent To show that $\mathcal R \in \mathscr{U}_{\mathbb X}$, we appeal to Proposition \ref{PROP: COMPRESSED ROOTS}, which specifies the iterative construction of $\mathscr{U}_{\mathbb X}$ and $\sigma_{\mathbb X}$. Following the prescription
\end{proof} 

\chapter{Directional maximal operators with {\tt{AdFinLac}} slopes} \label{positive direction chapter}
One of the main results of this paper, which justifies Definition \ref{defn: Admissible finite order lacunarity} as the correct formulation of finite-order lacunarity, is Theorem \ref{THM: DB}. It states that a direction set $\Omega \in {\tt{AdFinLac}}$ gives rise to a maximal directional operator $D_{\Omega}$ that is bounded on $L^p(\mathbb R^2)$ for each $p \in (1, \infty)$. This chapter is devoted to a proof of this theorem. 
 \vskip0.1in
\noindent Results similar in spirit to Theorem \ref{THM: DB} have appeared in \cite{{Alfonseca}, {AlfonsecaSoriaVargas1}, {AlfonsecaSoriaVargas2}, {Katz1}, {Katz2}}. We use a different proof strategy. The analytical backbone of the proof is an almost-orthogonality principle due to Christ \cite{Carbery}, which is dimension-free and effectively extracts the operator-theoretic ingredients of the problem from the geometric ones. We state this principle in the next section. 
\section{Christ's almost orthogonality principle} 	\label{section: Christ almost orthogonality principle}
In \cite[Section 2]{Carbery}, Carbery presented an operator-theoretic result, based on private communication with Christ, that provides an analytical framework for bootstrapping $L^p$ bounds of multi-parameter maximal operators. This result, which is central to our proof strategy, controls a $(k+1)$-parameter maximal operator using uniform $L^p$ bounds over certain $k$-parameter sub-families. We recall it here.
\vskip0.1in
\noindent Let us fix a dimension $d \geq 1$, a Lebesgue exponent $p \in (1, \infty)$, a constant $A_0 > 0$ and a collection of index sets 
\[\{\mathcal S_j : j \in \mathbb Z_{\geq 0} = \{0, 1, 2, \ldots\}\}. \] Two families of operators on $L^p(\mathbb R^d)$ are given:
\vskip0.1in 
\begin{itemize}
\item  A doubly indexed family of sub-additive operators \footnote{An operator $\mathcal T$ is {\em{subadditive}} if $|\mathcal T(f+g)| \leq |\mathcal T(f)| + |\mathcal T(g)|$ for all functions $f$ and $g$.}   
\begin{equation} \label{operator family T} 
{\pmb{\mathscr{T}}} := \{\mathscr{T}_{j \nu} : (j, \nu) \in \mathcal S^{\ast} \}, \text{ where } \mathcal S^{\ast} := \{ (j, \nu) : j \in \mathbb Z_{\geq 0}, \nu \in \mathcal S_j\}. 
\end{equation}   
\vskip0.1in
\item A singly indexed family of operators $\{\mathscr{R}_j : j \in \mathbb Z_{\geq 0}\}$.
\end{itemize}
\vskip0.1in  
We assume that these operator families obey following properties:
\begin{enumerate}[1.] 
 \item {\em{Uniform $L^p$ control on $\sup_{\nu} |\mathscr{T}_{j \nu}|$:}}
 \begin{equation} 
 \sup_{j \geq 0} \, \bigl|  \bigl| \sup_{\nu \in \mathcal S_j} |\mathscr{T}_{j \nu} f| \bigr| \bigr|_{p} \leq A_0 ||f||_p. \label{uniform op norm in j}
 \end{equation} 
 This ensures that the family of partial maximal operators  $\sup_{\nu} |\mathscr{T}_{j \nu}|$ indexed by $j$ is uniformly bounded in $L^p$ norm. 
 \vskip0.1in 
\item {\em{A square function estimate on $\{\mathscr R_j\}$:}} 
\begin{equation} 
\Bigl|\Bigl| \Bigl( \sum_{j \geq 0} \bigl| \mathscr{R}_j f\bigr|^2\Bigr)^{\frac{1}{2}}\Bigr|\Bigr|_p \leq A_0 ||f||_p.  \label{R-ORTHOGONAL} 
\end{equation} 
In other words, the operators $\mathscr{R}_j$ exhibit a form of weak orthogonality. 
\vskip0.1in
\item {\em{Strong control on the error term:}}
\begin{equation} 
\Bigl|\Bigl| \sup_{j, \nu} \bigl|\mathscr{T}_{j \nu} (I - \mathscr{R}_j)  f \bigr| \Bigr|\Bigr|_p \leq A_0 ||f||_p. \label{remainder op norm}
\end{equation}
This assumption implies that the worst behaviour of $\sup_{\nu} |\mathscr{T}_{j \nu}|$ occurs in the range of $\mathscr{R}_j$. Thus for the purpose of maximal estimates $\mathscr{T}_{j \nu}$ is essentially the same as  $\mathscr{T}_{j \nu} \mathscr{R}_j$, in the sense that the maximal error term $\sup_{j, \nu} \bigl|\mathscr{T}_{j \nu} (I - \mathscr{R}_j) \bigr|$ is $L^p$-bounded. 
\end{enumerate}  
\vskip0.1in 
We need an additional definition in order to state the almost orthogonality principle. A family ${\pmb{\mathscr{T}}}$ of the form \eqref{operator family T} is called {\em{essentially positive}} if there are operators $\mathscr{A}_{j \nu}, \mathscr{B}_{j \nu}, \mathscr{C}_{j \nu}$ with
\begin{align}
&\mathscr{T}_{j \nu} = \mathscr{A}_{j \nu} - \mathscr{B}_{j \nu}, \quad 
\bigl|\mathscr{A}_{j \nu} f \bigr| \leq \mathscr{A}_{j \nu} g \text{ whenever } |f| \leq g;  \label{ess-pos-1} \\
& |\mathscr{B}_{j \nu} f| \leq \mathscr{C}_{j \nu} |f|, \quad 0 \leq \mathscr{C}_{j \nu} f \leq \mathscr{C}_{j \nu} g \; \text{ for } 0 \leq f \leq g;  \label{ess-pos-2} \\ 
&\quad \Bigl| \Bigl| \sup_{j, \nu} \mathscr{C}_{j \nu} |f|  \Bigr| \Bigr|_p \leq A_0 ||f||_p. \label{ess-pos-3}
 \end{align} 
Intuitively, essential positivity means that each $\mathscr{T}_{j \nu}$ can be expressed as a difference of two auxiliary operators, one of which, namely $\mathscr{A}_{j \nu}$, is order-preserving. The other, namely $\mathscr{B}_{j \nu}$, is pointwise bounded by an order-preserving operator $\mathscr{C}_{j \nu}$ which satisfies a strong maximal estimate over the doubly indexed family $\{j \in \mathbb Z, \nu \in \mathcal S\}$. 
\vskip0.1in
\noindent We can now quote the almost orthogonality principle to be used in the proof. 
\begin{letteredtheorem}[M.~Christ, stated in {\cite[Theorem B]{Carbery}}] \label{Christ-theorem}
Suppose that for some $p \in [2, \infty)$, the conditions \eqref{uniform op norm in j}, \eqref{R-ORTHOGONAL} and \eqref{remainder op norm} hold. Then 
\begin{equation} \label{doubly indexed maximal inequality}  
\bigl|\bigl| \sup_{j, \nu} | \mathscr{T}_{j \nu} f| \bigr|\bigr|_p \leq C_0 ||f||_p
\end{equation}  
for some constant $C_0 >0$ depending only on $d, p$ and $A_0$. 
 \vskip0.1in
\noindent For $p \in (1, 2)$, the norm inequality \eqref{doubly indexed maximal inequality} holds under the extra hypothesis of essential positivity of ${\pmb{\mathscr{T}}}$, i.e. if one assumes that conditions \eqref{ess-pos-1}--\eqref{ess-pos-3} hold in addition to \eqref{uniform op norm in j}--\eqref{remainder op norm}.
\end{letteredtheorem}
\section{Proof  overview of Theorem \ref{THM: DB}} \label{Proof DB Overview Section} 
 We will aim to only prove the conclusion \eqref{D-Omega-bounded} of Theorem \ref{THM: DB}, since \eqref{just finiteness of op norm} is an immediate consequence of it. As Figure \ref{fig: proof flowchart DB} explains, the proof proceeds by replacing the original maximal operator by a smoother one, decomposing the latter into low- and high-frequency pieces, estimating the latter 
 through Christ's almost orthogonality principle, and finally combining the resulting bounds. 
\begin{center}
\begin{figure} 
\begin{tikzpicture}[
    box/.style={
        rectangle, draw, rounded corners,
        align=center, text badly centered,
        text width=4.6cm,
        inner sep=4.5pt,
        minimum height=0.85cm
    },
    smallbox/.style={
        rectangle, draw, rounded corners,
        align=center, text badly centered,
        text width=3.8cm,
        inner sep=4.5pt,
        minimum height=0.85cm
    },
    arrow/.style={
        -{Stealth[length=2mm]},
        very thick,
        shorten >=4pt,
        shorten <=4pt
    },
    line/.style={very thick}
]

\node[box] (goal)
{\textbf{Goal:} Prove Theorem \ref{THM: DB},\\ Show \(D_\Omega\) is \(L^p\)-bounded};

\node[box, below=7mm of goal] (roadmap)
{\textbf{\S \ref{Proof DB Overview Section}}: Proof roadmap};

\node[box, right=8mm of roadmap] (christ)
{\textbf{\S \ref{section: Christ almost orthogonality principle}}: The abstract set-up,\\ Theorem \ref{Christ-theorem}};

\node[box, below=7mm of roadmap] (regular)
{\textbf{\S \ref{section: smoothing}} Regularization:\\
Replace \(\mathcal A_{\omega,h} \mapsto \tilde{\mathcal A}_{\omega,h}\) , \(D_\Omega \mapsto \tilde D_\Omega\)};

\node[box, below=7mm of regular] (decomp)
{\textbf{\S \ref{section: decomposition of A}}: Decompositions:\\
\(\tilde{\mathcal A}_{\omega,h}
=\mathfrak B_{\omega,h}+\mathfrak C_{\omega,h}\)\\
\(\tilde D_\Omega\leq
\mathfrak B_\Omega^\ast+\mathfrak C_\Omega^\ast\)};

\node[smallbox, below left=10mm and 5mm of decomp] (low)
{\textbf{\S \ref{section: B-star}}:
Low frequency \\ bound \(\mathfrak B_\Omega^\ast\)};

\node[smallbox, below right=10mm and 5mm of decomp] (high)
{\textbf{\S \ref{section: Induction C-star}}:
High frequency \\ bound \(\mathfrak C_\Omega^\ast\)};

\node[smallbox, below=7mm of high] (square)
{\textbf{\S \ref{square function estimate section}}:\\
Square function estimate};

\node[smallbox, below=7mm of square] (loc)
{\textbf{\S \ref{section: cone geometry}}:\\
Localization identity};

\node[box, below=13mm of loc, xshift=-2.9cm] (concl)
{\textbf{\S \ref{section: DB conclusion}}:\\
Combine estimates from \S \ref{section: B-star} and \S \ref{section: Induction C-star} to derive the conclusion \eqref{D-Omega-bounded}};

% Main proof route
\draw[arrow] (goal) -- (roadmap);
\draw[arrow] (roadmap) -- (regular);
\draw[arrow] (regular) -- (decomp);

% Goal to Christ input
\draw[arrow] (goal.east) -| (christ.north);

% Christ input into high-frequency section
\coordinate (christdrop) at ($(high.north)+(0,7mm)$);
\draw[line] (christ.south) -- ++(0,-6mm) -| (christdrop);
\draw[arrow] (christdrop) -- (high.north);

% Split after decomposition
\coordinate (split) at ($(decomp.south)+(0,-5mm)$);
\coordinate (lowtop) at ($(low.north)+(0,5mm)$);
\coordinate (hightop) at ($(high.north)+(0,5mm)$);

\draw[line] (decomp.south) -- (split);
\draw[line] (split) -| (lowtop);
\draw[line] (split) -| (hightop);
\draw[arrow] (lowtop) -- (low.north);
\draw[arrow] (hightop) -- (high.north);

% High-frequency chain
\draw[arrow] (high) -- (square);
\draw[arrow] (square) -- (loc);

% Merge before conclusion
\coordinate (merge) at ($(concl.north)+(0,6mm)$);
\coordinate (lowbottom) at ($(low.south)+(0,-6mm)$);
\coordinate (locbottom) at ($(loc.south)+(0,-6mm)$);

\draw[line] (low.south) -- (lowbottom);
\draw[line] (loc.south) -- (locbottom);
\draw[line] (lowbottom) -| (merge);
\draw[line] (locbottom) -| (merge);
\draw[arrow] (merge) -- (concl.north);

\end{tikzpicture}
\caption{Proof structure of Theorem \ref{THM: DB}} \label{fig: proof flowchart DB}
\end{figure}
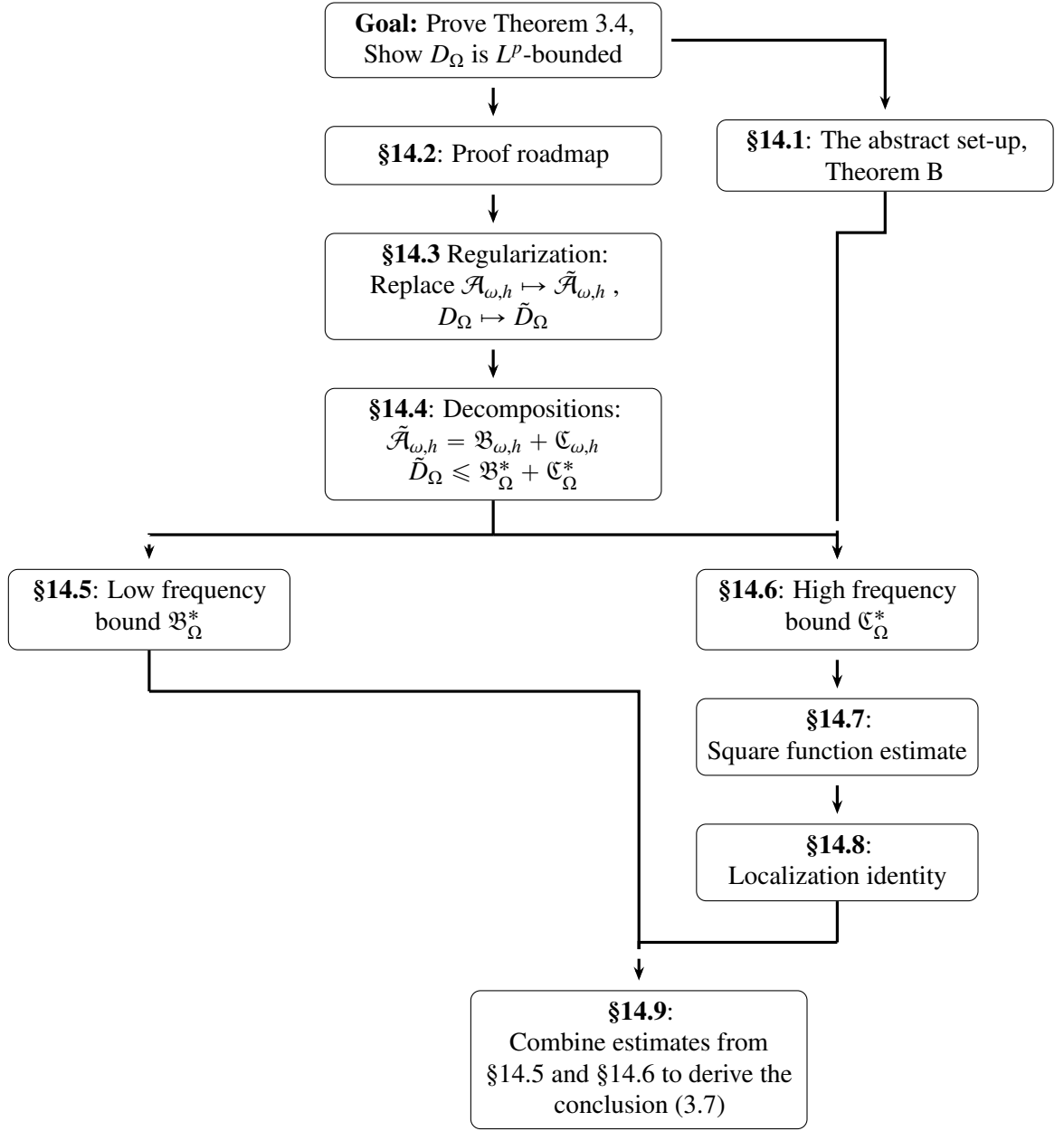 
\end{center}

\subsection{Regularization} The maximal operator $D_{\Omega}$, given by \eqref{dir-max-op-def}, is defined via averages over finite line segments:
\begin{equation} \label{old average} 
\mathcal A_{\omega, h} f(x) = \frac{1}{2h} \int_{-h}^{h} |f\bigl(x+ t(1, \omega) \bigr)| \, dt, 
\end{equation}  
which a priori involve sharp cut-offs on lines. Our first step is to bound $D_{\Omega}$ pointwise by another maximal operator $\tilde{D}_{\Omega}$ composed of smoother averages $\tilde{\mathcal A}_{\omega, h}$. This step is completed in Section \ref{section: smoothing}, with $\tilde{\mathcal A}_{\omega, h}$ and $\tilde{D}_{\Omega}$ defined in \eqref{new average} and \eqref{D<Dtilde} respectively. The smoothness inherent in $\tilde{\mathcal A}_{\omega, h}$ provides additional technical convenience in subsequent steps. 
\vskip0.1in
\noindent To establish the desired conclusion \eqref{D-Omega-bounded}, it therefore suffices to prove the existence of a constant $C(p, N, \lambda) > 0$ such that the following holds: 
\begin{equation} \label{Dtilde-Omega-bounded}
||\tilde{D}_{\Omega}||_{p \rightarrow p} \leq C(p, N, \lambda) \quad \text{ for all } \Omega \in \Lambda(N, \lambda). 
\end{equation} 
Henceforth, we focus on proving \eqref{Dtilde-Omega-bounded}. 
\subsection{Decomposition of the averaging operator $
\tilde{\mathcal A}_{\omega, h}$} 
The maximal operator $\tilde{D}_{\Omega}$ is decomposed further, to reduce to yet another maximal operator that captures the non-trivial portion of it. In Section \ref{section: decomposition of A}, we split the averaging operator $\tilde{\mathcal A}_{\omega, h}$ into two parts: an innocuous term $\mathfrak B_{\omega, h}$ that is Fourier supported within a large origin-centred ball, and a main part $\mathfrak C_{\omega, h}$ that is Fourier-localized in high frequencies. This results in a pointwise bound for $\tilde{D}_{\Omega}$ by the sum of two maximal operators $\mathfrak B^{\ast}_{\Omega}$ and $\mathfrak C^{\ast}_{\Omega}$, originating from $\mathfrak B_{\omega, h}$ and $\mathfrak C_{\omega, h}$ respectively.  
\subsection{The operator $\mathfrak B^{\ast}_{\Omega}$} In Section \ref{section: B-star}, we show that $\mathfrak B^{\ast} = \mathfrak B^{\ast}_{\Omega}$ is bounded pointwise by the classical strong maximal function (independently of $\Omega$), and is therefore benign in terms of $L^p$ bounds, for all $p \in (1, \infty)$. The onus of proving \eqref{Dtilde-Omega-bounded} therefore shifts to proving a similar statement with $\tilde{D}_{\Omega}$ replaced by $\mathfrak C^{\ast}_{\Omega}$. The reduced statement appears in Proposition \ref{prop: C*}, replacing \eqref{Dtilde-Omega-bounded} by \eqref{Maximal C inequality} as the main statement to be proved. 
\vskip0.1in
\noindent 
\subsection{The operator $\mathfrak C_{\Omega}^{\ast}$} 
In Sections \ref{section: Induction C-star}--\ref{section: cone geometry}, we complete the proof of Theorem~\ref{THM: DB}
by induction on the lacunarity order $N$. The base case $N=0$ follows from
the $L^p$ boundedness of the one-dimensional Hardy-Littlewood maximal
operator. For the inductive step, we assume the conclusion of
Theorem~\ref{THM: DB} for sets of lacunarity order at most
$N-1$. Proposition~\ref{prop: C*} then uses this induction
hypothesis, together with Theorem \ref{Christ-theorem} from Section \ref{section: Christ almost orthogonality principle}, to establish the
high-frequency estimate
\[
        \|\mathfrak C_\Omega^\ast\|_{p\to p}
        \leq C(p,N,\lambda)
\]
for $\Omega\in\Lambda(N,\lambda)$. Combining this estimate with the
low-frequency bound from Section \ref{section: B-star} proves the theorem at order $N$.
%In Sections \ref{section: Induction C-star} -- \ref{section: cone geometry}, we prove the revised estimate \eqref{Maximal C inequality}. The proof involves induction on the lacunarity order $N$, using Theorem \ref{Christ-theorem} from Section \ref{section: Christ almost orthogonality principle} at each inductive step. 
\vskip0.1in
\noindent Specifically, the operator $\mathscr{T}_{j \nu}$ in Theorem \ref{Christ-theorem} will correspond to $\mathfrak C_{\omega, h}$, the high frequency portion of $\mathcal A_{\omega, h}$ localized outside a large ball centred at the origin. We need to set up the induction so that Theorem \ref{Christ-theorem} is applicable, choosing suitable localization operators $\mathscr{R}_j$ that verify conditions \eqref{uniform op norm in j}--\eqref{ess-pos-3}. These verifications, carried out in Sections \ref{section: Induction C-star}--\ref{section: cone geometry}, occupy the remainder of the chapter. Section \ref{section: DB conclusion}  assembles the different ingredients to complete the proof of Theorem \ref{THM: DB}. 
%The induction hypothesis provides the criterion \eqref{uniform op norm in j}, but the operator $\mathscr{T}_{j \nu}$ is not exactly $\mathcal A_{h, \omega}$, rather a portion of it, termed $\mathfrak C_{\omega, h}$, supported outside a large ball in frequency space (\S \ref{section: decomposition of A}). Framing the induction on $N$ in a way so that Theorem A applies (\S \ref{section: B-star}), identifying the localization operators $\mathscr R_j$, and verifying the hypotheses required for the application take up the remainder of the proof (\S \ref{section: Induction C-star}). 

\section{The operator $\tilde{D}_{\Omega}$ related to the smoothed averages $\tilde{\mathcal A}_{h, \omega}$} \label{section: smoothing}
\subsection{Constants $c_0, R_0$ and an auxiliary function $\psi$} Let us fix a small constant $c_0 \in (0, 1)$ depending on $\lambda$:
\begin{equation} \label{little c0}
100c_0 < \min \left[ \lambda^2, 1 - \lambda^2 \right]. 
\end{equation} 
We will also need a Schwartz function $\psi$ with the following properties: 
\begin{equation} \label{psi requirements}
\left\{
\begin{aligned}  
&\psi: \mathbb R \rightarrow [0, \infty), \; \widehat{\psi}(\xi) := \int e^{-i x \xi} \psi(x) \, dx \geq 0 \text{ for all $\xi \in \mathbb R$}, \\ 
&\psi(x) \geq 1 \text{ for } x \in [-1,1], \; \; {\text{supp}}(\widehat{\psi}) \subseteq [-R_0, R_0], \; \; R_0 > 0. 
\end{aligned}
\right\}
\end{equation}  
The notation $\widehat{\psi}$ denotes the Fourier transform of $\psi$. 
Here $R_0$ is a fixed, large absolute constant; for instance $R_0 =100$ will suffice.  
\vskip0.1in
\noindent To confirm the existence of a function $\psi$ obeying all the requirements in \eqref{psi requirements}, let us choose a non-negative, smooth, even function 
\begin{equation}  \chi_0 : \mathbb R \rightarrow [0,1], \quad \chi_0 \equiv \left\{ \begin{aligned} 1 &\text{ on } [-1,1], \\ 0 &\text{ outside } [-2, 2], \end{aligned} \right\} \; \text{ and define } \psi_0 := |\widecheck{\chi}_0|^2, \label{what is chi_0} \end{equation}     
where $\widecheck{\chi}_0$ represents the inverse Fourier transform of $\chi_0$. It follows from the definition that 
\begin{align*} 
&\psi_0 \geq 0, \; \widehat{\psi}_0 = \chi_0 \ast \chi_0 \geq 0,  \;  \text{supp}(\widehat{\psi}_0) \subseteq [-4,4], \text{ and } \\
&\psi_0(0) = \bigl| \widecheck{\chi}_0(0)\bigr|^2 = \Bigl[ \frac{1}{2 \pi} \int {\chi}_0(x) \, dx \Bigr]^2 > 0.
\end{align*} 
Therefore one can find constants $c_1, c_2 > 0$ depending on $R_0$ such that the function $\psi(x) := c_2 \psi_0(c_1 x)$ obeys all the conditions in \eqref{psi requirements}.  

\subsection{The operators $\tilde{\mathcal A}_{\omega, h}$ and $\tilde{D}_{\Omega}$} For $\psi$ as in \eqref{psi requirements}, $h > 0$, $\omega \in \Omega$ and non-negative, smooth function $f: \mathbb R^2 \rightarrow \mathbb R$, let us define the smoother average: 
\begin{equation} \label{new average}
\tilde{\mathcal A}_{\omega, h} f(x) := \frac{1}{2h} \int_{\mathbb R}  f \bigl(x + t(1, \omega) \bigr) \,  \psi \left(\frac{t}{h} \right) \, dt,   
\end{equation} 
and the corresponding maximal operator:
\begin{equation} 
\tilde{D}_{\Omega}f(x) := \sup_{\omega \in \Omega} \sup_{h > 0} \tilde{\mathcal A}_{\omega, h} f(x). \label{D<Dtilde}
\end{equation} 
The operators $\tilde{\mathcal A}_{\omega, h}$ and $\tilde{D}_{\Omega}$ are essentially equivalent to $\mathcal A_{\omega, h}$ and $D_{\Omega}$ respectively, as the following lemma shows. 
\subsection{The pointwise equivalence of $D_{\Omega}$ and $\tilde{D}_{\Omega}$} 
\begin{lemma}
There exists a constant $C_0 > 0$ depending only on $\psi$ with the following property: 
\begin{enumerate}[(a)]
\item For $\omega \in \Omega$ and $h > 0$, the operators $\tilde{\mathcal A}_{\omega, h}$ and $\tilde{D}_{\Omega}$ obey the pointwise estimate 
\begin{align} 
&\mathcal A_{\omega, h} f(x)  \leq \tilde{\mathcal A}_{\omega, h}  f(x) \leq C_0 D_{\Omega}f(x),  \label{A-Atilde-pointwise}\\ 
&\text{therefore, } \; D_{\Omega} f(x) \leq \tilde{D}_{\Omega}f(x) \leq  C_0 D_{\Omega}f(x) \label{D-Dtilde-pointwise}
\end{align} 
for all $x \in \mathbb R^2$ and for all non-negative smooth functions $f$. Here $\mathcal A_{\omega, h}$ is the averaging operator in \eqref{old average} and the directional maximal operator $D_{\Omega}$ defined in \eqref{dir-max-op-def} is the supremum of $\mathcal A_{\omega, h}$ over $h >0$, $\omega \in \Omega$. 
\vskip0.1in
\item As a result, 
\begin{equation}  
||D_{\Omega}||_{p \rightarrow p} \leq ||\tilde{D}_{\Omega}||_{p \rightarrow p} \leq C_0 ||D_{\Omega}||_{p \rightarrow p} \; \text{ for all } p \in [1, \infty]. \label{D-Dtilde-norm}
\end{equation}  
In other words, $D_{\Omega}$ is bounded on $L^{p}(\mathbb R^2)$ if and only if $\tilde{D}_{\Omega}$ is. 
\end{enumerate} 
\end{lemma} 
\begin{proof} 
The inequalities in \eqref{D-Dtilde-pointwise} and \eqref{D-Dtilde-norm} follow directly from \eqref{A-Atilde-pointwise}; therefore we focus only on proving the latter. 
\vskip0.1in 
\noindent The lower bound $\psi \geq 1$ on $[-1, 1]$ gives rise to the pointwise inequality 
\[ \psi\left( \frac{t}{h}\right) \geq 1 \text{ for } t \in [-h, h]. \] 
Substituting this into the integral \eqref{old average} representing $\mathcal A_{\omega, h}$ gives rise to the left inequality in \eqref{A-Atilde-pointwise}. 
\vskip0.1in
\noindent To establish the inequality on the right, we split the integral representing $\tilde{\mathcal A}_{\omega, h}$ into dyadic annuli outside $[-h,h]$, and use the rapid decay of $\psi(t/h)$ outside this interval to dominate $\tilde{\mathcal A}_{\omega, h}$ by a geometric series whose main contribution comes from $[-h,h]$. Specifically, let us note that for every $M \geq 1$, there is a constant $C_M > 0$ depending only on $R_0$ and $\widehat{\psi}^{(\alpha)}$, $0 \leq \alpha \leq M$, such that  
\begin{align}
 \psi\left(\frac{t}{h} \right) &\leq C_M 2^{-\ell M} \; \text{ for } \; |t| \in h \bigl[2^{\ell}, 2^{\ell+1} \bigr] \text{ and all $\ell \geq 1$; hence } \nonumber    \\
\tilde{\mathcal A}_{\omega, h} f(x) &\leq \frac{1}{2h} \left[ \int_{|t| \leq h}  ||\psi||_{\infty} + \sum_{\ell \geq 1} \int_{|t| \leq 2^{\ell} h} C_{M} 2^{-\ell M} \right] f\bigl(x + t(1, \omega) \bigr)  \, dt \nonumber \\ 
&\leq C_0 \left[ 1 + \sum_{\ell \geq 1} 2^{-\ell(M-1)}\right] D_{\Omega}f(x) \leq C_0 D_{\Omega} f(x),  \nonumber
%\tilde{D}_{\Omega} f(x) &\leq C_0 D_{\Omega} f(x) \text{ where $C_0$ depends only on $\psi$.} \label{Dtilde<D}
\end{align} 
yielding the right inequality in \eqref{A-Atilde-pointwise}, for any $M \geq 2$. This completes the proof. 
\end{proof} 
\subsection{The multiplier of $\tilde{\mathcal A}_{\omega, h}$} The fact that $\tilde{\mathcal A}_{\omega, h}$ is a convolution operator ensures that it is given by a smooth multiplier $\mathfrak a_{\omega, h}$ in Fourier space: for $\xi = (\xi_1, \xi_2) \in \mathbb R^2$, 
\begin{equation} \label{multiplier a}  \bigl[ \tilde{\mathcal A}_{\omega, h} f \bigr]^{\wedge}(\xi) = \mathfrak a_{\omega, h}(\xi) \widehat{f}(\xi)  \; \text{ where } \mathfrak a_{\omega, h} (\xi) := \widehat{\psi}\bigl( h(1, \omega) \cdot \xi \bigr) . \end{equation}  
It follows from the properties of $\widehat{\psi}$ in \eqref{psi requirements} that 
\begin{equation} \label{infinite strip}
\text{supp}\bigl(\mathfrak a_{\omega, h} \bigr) \subseteq \left\{ \xi \in \mathbb R^2 :  |\xi_1+ \omega \xi_2| \leq \frac{R_0}{h}\right\}. 
\end{equation}   
Geometrically, this is an infinite planar strip of width $2R_0 /(h \sqrt{1+ \omega^2})$, whose central axis passes through the origin and is perpendicular to the direction $(1, \omega)$. A visual depiction of this is in Figure \ref{fig: multiplier support}. 

% Preamble:
% \usepackage{tikz}
% \usetikzlibrary{calc,arrows.meta,decorations.pathreplacing}

\begin{center}
\begin{figure}
\begin{tikzpicture}[
    axis/.style={-{Stealth[length=2.5mm]}, very thick},
    bluevec/.style={blue!70!black, very thick},
    support/.style={purple!70!black, very thick},
    hatch/.style={cyan!65, thin},
    brace/.style={decorate, decoration={brace, amplitude=4pt}}
]

%%%%%%%%%%%%%%%%%%%%%%%%%%%%%%%%%%%%%%%%%%%%%%%%%%%%%%%%%%%%%%
%% Parameters
%%%%%%%%%%%%%%%%%%%%%%%%%%%%%%%%%%%%%%%%%%%%%%%%%%%%%%%%%%%%%%

% Numerical value of the slope omega
\pgfmathsetmacro{\omegaslope}{0.22}

% Half-width of the strip (= R_0/h)
\pgfmathsetmacro{\halfwidth}{1.55}

% Half-length of the strip
\pgfmathsetmacro{\halflength}{2.75}

% Normalization factor
\pgfmathsetmacro{\norm}{sqrt(1+\omegaslope*\omegaslope)}

%%%%%%%%%%%%%%%%%%%%%%%%%%%%%%%%%%%%%%%%%%%%%%%%%%%%%%%%%%%%%%
%% Unit vectors
%%%%%%%%%%%%%%%%%%%%%%%%%%%%%%%%%%%%%%%%%%%%%%%%%%%%%%%%%%%%%%

% u = (1,omega)/|(1,omega)|
\coordinate (u) at ({1/\norm},{\omegaslope/\norm});

% v = (-omega,1)/|(1,omega)|
\coordinate (v) at ({-\omegaslope/\norm},{1/\norm});

%%%%%%%%%%%%%%%%%%%%%%%%%%%%%%%%%%%%%%%%%%%%%%%%%%%%%%%%%%%%%%
%% Coordinate axes
%%%%%%%%%%%%%%%%%%%%%%%%%%%%%%%%%%%%%%%%%%%%%%%%%%%%%%%%%%%%%%

\draw[axis] (-4.3,0) -- (4.8,0)
    node[right] {$\xi_1$};

\draw[axis] (0,-2.6) -- (0,3.0)
    node[above] {$\xi_2$};

%%%%%%%%%%%%%%%%%%%%%%%%%%%%%%%%%%%%%%%%%%%%%%%%%%%%%%%%%%%%%%
%% Strip vertices
%%%%%%%%%%%%%%%%%%%%%%%%%%%%%%%%%%%%%%%%%%%%%%%%%%%%%%%%%%%%%%

\coordinate (A) at ($(0,0)-\halfwidth*(u)-\halflength*(v)$);
\coordinate (B) at ($(0,0)-\halfwidth*(u)+\halflength*(v)$);
\coordinate (C) at ($(0,0)+\halfwidth*(u)+\halflength*(v)$);
\coordinate (D) at ($(0,0)+\halfwidth*(u)-\halflength*(v)$);

%%%%%%%%%%%%%%%%%%%%%%%%%%%%%%%%%%%%%%%%%%%%%%%%%%%%%%%%%%%%%%
%% Shading
%%%%%%%%%%%%%%%%%%%%%%%%%%%%%%%%%%%%%%%%%%%%%%%%%%%%%%%%%%%%%%

\begin{scope}

\clip (A) -- (B) -- (C) -- (D) -- cycle;

\foreach \t in {-2.8,-2.45,...,2.8}
{
    \draw[hatch]
        ($(0,0)-\halfwidth*(u)+\t*(v)$)
        --
        ($(0,0)+\halfwidth*(u)+\t*(v)$);
}

\end{scope}

%%%%%%%%%%%%%%%%%%%%%%%%%%%%%%%%%%%%%%%%%%%%%%%%%%%%%%%%%%%%%%
%% Strip boundaries
%%%%%%%%%%%%%%%%%%%%%%%%%%%%%%%%%%%%%%%%%%%%%%%%%%%%%%%%%%%%%%

\draw[support] (A) -- (B);
\draw[support] (D) -- (C);

%%%%%%%%%%%%%%%%%%%%%%%%%%%%%%%%%%%%%%%%%%%%%%%%%%%%%%%%%%%%%%
%% Direction (1,\omega)
%%%%%%%%%%%%%%%%%%%%%%%%%%%%%%%%%%%%%%%%%%%%%%%%%%%%%%%%%%%%%%

\draw[bluevec]
    ($(0,0)-4.2*(u)$)
    --
    ($(0,0)+4.2*(u)$)
    node[right] {$\xi_2 = \omega \xi_1$};

%%%%%%%%%%%%%%%%%%%%%%%%%%%%%%%%%%%%%%%%%%%%%%%%%%%%%%%%%%%%%%
%% Width annotation
%%%%%%%%%%%%%%%%%%%%%%%%%%%%%%%%%%%%%%%%%%%%%%%%%%%%%%%%%%%%%%

\draw[
    brace,
    blue!70!black,
    thick
]
    ($(0,0)+0.05*(u)$)
    --
    ($(0,0)+\halfwidth*(u)$)
node[midway, above=6pt, blue!70!black]
{\scriptsize{$\dfrac{R_0}{h \sqrt{1 + \omega^2}}$}};

%%%%%%%%%%%%%%%%%%%%%%%%%%%%%%%%%%%%%%%%%%%%%%%%%%%%%%%%%%%%%%
%% Label
%%%%%%%%%%%%%%%%%%%%%%%%%%%%%%%%%%%%%%%%%%%%%%%%%%%%%%%%%%%%%%

\node[purple!70!black]
    at (3.0,-1.6)
    {$\operatorname{supp}(\mathfrak a_{\omega,h})$};

\end{tikzpicture}
\caption{\small{The support of the multiplier for $\tilde{\mathcal A}_{\omega, h}$ is contained in an infinite planar strip centred at the origin, perpendicular to $(1, \omega)$.}} \label{fig: multiplier support} 
\end{figure}
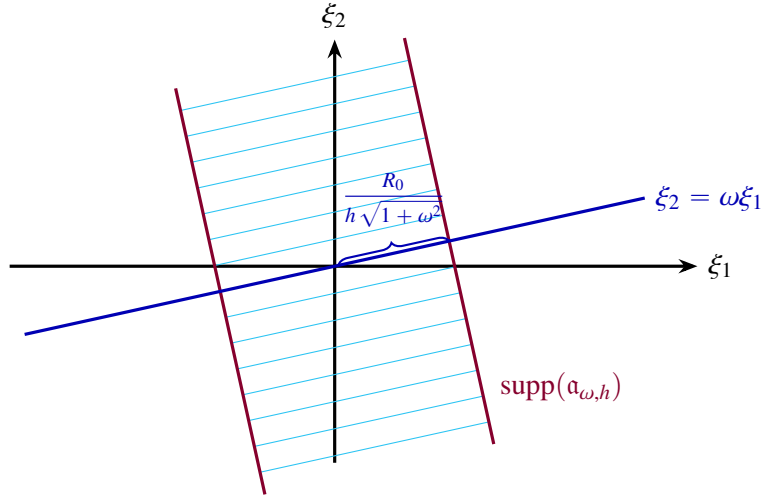
\end{center}

\section{Decomposition of $\widetilde{\mathcal A}_{\omega, h}$} \label{section: decomposition of A}
In this section, we use the geometry of the slope set $\Omega$ to decompose the averaging operator $\tilde{\mathcal A}_{\omega, h}$ into two parts. 
\vskip0.1in
\noindent Let us recall that $\Omega \subseteq [0,1]$ lies in $\Lambda(N, \lambda)$ for some $N \geq 0$ and $\lambda \in (0,1)$. According to Definition \ref{defn: Lacunary sets} of $\Lambda(N, \lambda)$ for $N \geq 1$, the set $\Omega$ has a special sequence; Lemma \ref{lemma: special sequence choice} ensures that the special sequence can be chosen to obey a two-sided inequality of the form \eqref{lacunary above and below}. After an invertible affine linear transformation on $\Omega$, which leaves $||D_{\Omega}||_{p \rightarrow p}$ unchanged, we may without loss of generality reduce to the case where $\Omega \in \Lambda(N, \lambda)$, with the special sequence of $\Omega$ of the form
\begin{align}
&A = \{a_j : j \geq 0 \} \in {\tt{MonLac}}(\lambda), \; a_j \searrow 0, \; a_0=1, \label{generic A} \\
&\; \lambda^2 a_{j} < a_{j+1} \leq \lambda a_j \text{ for all } j \geq 0. \label{special A} 
\end{align} 
\vskip0.1in
\noindent Let us denote by $\Omega_j$ the portion of $\Omega$ in   the $j^{\text{th}}$ gap interval of $A$: 
\begin{equation} \label{Omegaj}
\Omega_j := \Omega \cap [a_{j+1}, a_j), \text{ so that } \Omega_j \in \Lambda(N-1, \lambda). 
\end{equation} 
%We will require the special sequence $A$ to obey an additional property that places a stronger restriction on the lacunary gap:
%Such a choice of special sequence is permitted for any $\Omega \in \Lambda(N, \lambda)$, as we have shown in Lemma \ref{lemma: special sequence choice} in the Appendix \S \ref{section: lacunarity properties}.  
\vskip0.1in
\noindent The gap intervals $[a_{j+1},a_j)$ cover $(0,1)$. The possible exceptional slopes in $\Omega\cap\{0,1\}$
form a set of cardinality at most two. Since the maximal operators
associated with individual slopes are bounded on $L^p$, their
contribution may be absorbed into the final constant. It therefore
suffices to consider
\[
        \Omega^\circ:=\Omega\cap(0,1).
\]
Replacing $\Omega$ by $\Omega^\circ$, if necessary, we may henceforth assume that each slope $\omega \in \Omega$ lies in $\Omega_j$ for a unique index $j = j(\omega)$. The lacunary separators $a_j, a_{j+1}$ enable a decomposition of the multiplier $\mathfrak a_{\omega, h}$ depending on $j$: 
\begin{align} 
&\text{ Setting } \delta_j := a_j - a_{j+1}, \text{ we write } \mathfrak a_{\omega, h} = \mathfrak b_{\omega, h} + \mathfrak c_{\omega, h}, \text{ with }  \label{decomposition of a}\\ 
&\mathfrak b_{\omega, h}(\xi) = \mathfrak a_{\omega, h}(\xi) \Phi(h \delta_j |\xi|), \quad\mathfrak c_{\omega, h}(\xi) = \mathfrak a_{\omega, h}(\xi) \bigl[ 1 - \Phi(h \delta_j |\xi|) \bigr]. \label{multipliers b and c} \end{align}
Here $\Phi$ is a bivariate, radial, compactly supported, smooth function 
\begin{equation}  \label{defn: Phi}
\Phi(\xi) := \chi_0\left( \frac{c_0 |\xi|}{10 R_0} \right) \; \text{ with $\chi_0: \mathbb R \rightarrow [0, 1]$ and $c_0$ as in \eqref{what is chi_0}, \eqref{little c0}}. 
\end{equation}  
In other words, the multiplier $\mathfrak b_{\omega, h}$ is essentially the same as $\mathfrak a_{\omega, h}$, but localized to a ball of radius $20 R_0(c_0h \delta_j)^{-1}$ centred at the origin, capturing its low-frequency component. The high-frequency part, namely the portion of $\mathfrak a_{\omega, h}$ outside this ball, is captured by $\mathfrak c_{\omega, h}$. 
\vskip0.1in
\noindent The decomposition \eqref{decomposition of a} of $\mathfrak a_{\omega, h}$ generates a decomposition of the operator $\tilde{\mathcal A}_{\omega, h}$ in \eqref{new average}. 
\begin{align} 
&\tilde{\mathcal A}_{\omega, h} = \mathfrak B_{\omega, h} + \mathfrak C_{\omega, h} \text{ for all $\omega \in \Omega$ and $h > 0$, where } \label{A=B+C} \\   
&   \bigl[ \mathfrak B_{\omega, h} f \bigr]^{\wedge} (\xi) =  \mathfrak b_{\omega, h}(\xi) \widehat{f}(\xi), \quad  \bigl[ \mathfrak C_{\omega, h} f \bigr]^{\wedge} (\xi) =   \mathfrak c_{\omega, h}(\xi) \widehat{f}(\xi). \label{b-and-c-multipliers} 
\end{align} 
 This, in turn, enables us to control the maximal operator $\tilde{D}_{\Omega}$ in two parts: 
\begin{align}
&\tilde{D}_{\Omega} f(x) = \sup_{\omega \in \Omega} \sup_{h > 0} \tilde{\mathcal A}_{\omega, h} f(x)  \leq \mathfrak B^{\ast}_{\Omega}f(x) + \mathfrak C^{\ast}_{\Omega} f(x),  \text{ where } \label{B*C*}  \\ 
&\mathfrak B^{\ast}_{\Omega}f(x) := \sup_{\omega \in \Omega} \sup_{h > 0} \bigl|\mathfrak B_{\omega, h} f(x) \bigr|,  \quad \mathfrak C^{\ast}_{\Omega} f(x) := \sup_{\omega \in \Omega} \sup_{h > 0} \bigl|\mathfrak C_{\omega, h} f(x) \bigr|, \label{BC-maximal}
\end{align}
Thus, $\mathfrak B^{\ast} = \mathfrak B^{\ast}_{\Omega}$ and $\mathfrak C^{\ast} = \mathfrak C^{\ast}_{\Omega}$ collect the maximal contributions from the low-frequency and high-frequency components of $\tilde{D}_{\Omega}$.
\section{The low frequency maximal operator $\mathfrak B^{\ast}$} \label{section: B-star}
As in previous boundedness proofs \cite{{Alfonseca}, {Bateman}}, we show that the maximal operator $\mathfrak B^{\ast}$ corresponding to the low frequency part is relatively  easy to control. In fact, we will prove a slightly stronger statement, specifically \eqref{pointwise estimate: strong maximal function} below, that will be used in the sequel (see Section \ref{section: ess pos}).  
\vskip0.1in
\noindent Let us recall from \eqref{multiplier a} and \eqref{multipliers b and c} that $\mathfrak B_{\omega, h}$ is a convolution operator given by the Fourier multiplier $\mathfrak b_{\omega, h}$. In other words, 
\begin{equation} \label{kernel check b}
\mathfrak B_{\omega, h} f(x)  = f \ast \check{\mathfrak b}_{\omega, h}(x), \text{ with } \check{\mathfrak b}_{\omega, h}(x) :=  \int_{\mathbb R^2} e(x \cdot \xi) \mathfrak b_{\omega, h}(\xi) \, d\xi.   
\end{equation} 
Here $e(t) := e^{it}$. Let us denote by $\widetilde{\mathfrak B}_{\omega, h}$ the convolution operator with respect to the kernel $|\check{\mathfrak b}_{\omega, h}|$. In other words, 
\begin{equation} \label{newB} 
\widetilde{\mathfrak B}_{\omega, h}f(x) := f \ast |\check{\mathfrak b}_{\omega, h}|(x).
\end{equation} 
Of course, this implies the pointwise bound 
\begin{equation} \label{B and new B}
|\mathfrak B_{\omega, h} f(x)| \leq  \widetilde{\mathfrak B}_{\omega, h}|f|(x) \quad \text{ for all } \omega \in \Omega, \; h > 0. 
\end{equation}
\vskip0.1in
\noindent We will use $ \mathcal M_{\text{\tiny{str}}}$ to denote the classical strong maximal function
\[ \mathcal M_{\text{\tiny{str}}}f(x) := \sup_{R \in \mathcal R_{\text{\tiny{str}}}} \frac{1}{|R|} \int_{R} |f(x + y)| \, dy, \] 
defined as the supremal average of $|f(x + \cdot)|$ over the family $\mathcal R_{\text{\tiny{str}}}$ of all axes-parallel rectangles centred at the origin. It is well-known \cite[Chapter II, Section E]{SteinHA} that $\mathcal M_{\text{\tiny{str}}}$ is bounded on $L^p(\mathbb R^2)$ for $p \in (1, \infty]$. The next result shows that this strong maximal function dominates $\mathfrak B^{\ast}$ pointwise.  
\begin{lemma} \label{lemma : strong maximal function}
For every $\lambda \in (0, 1)$ and $R_0 > 0$ as in \eqref{psi requirements}, there is a constant $C_1 = C_1(\lambda, R_0) >0$ depending only on these quantities, with the following property.
\vskip0.1in 
\noindent Given any slope set $\Omega$ and $\Omega_j$ as in \eqref{Omegaj},  any $\omega \in \Omega_j$ and $h > 0$, the pointwise inequality  
\begin{equation}  \bigl|\mathfrak B_{\omega, h}f(x) \bigr| \leq \widetilde{\mathfrak B}_{\omega, h}|f|(x) \leq C_1 \mathcal M_{\text{\tiny{str}}}f(x)  \label{pointwise estimate: strong maximal function}  \end{equation} 
holds for all bivariate smooth functions $f$ and all $x \in \mathbb R^2$. 
\vskip0.1in
\noindent Consequently, $\mathfrak B^{\ast}f(x) \leq C_1 \mathcal M_{\text{str}}f(x)$, which implies 
\[ ||\mathfrak B^{\ast}||_{p \rightarrow p} \leq C_1 ||\mathcal M_{\text{\tiny{str}}}||_{p \rightarrow p} < \infty, \; \text{ for all } p \in (1, \infty). \]  
\end{lemma}  
\begin{proof}
The first inequality in \eqref{pointwise estimate: strong maximal function} has already been recorded in \eqref{B and new B}. It therefore suffices to establish the second inequality. 
\vskip0.1in
\noindent The key estimate for $\check{\mathfrak b}_{\omega, h}$, which we establish in Lemma \ref{IBP-lemma} below, states that for every integer $M_1, M_2 \geq 1$,
\begin{equation} \bigl| \check{\mathfrak b}_{\omega, h}(x) \bigr| \leq \frac{C_{M_1, M_2, \lambda}}{h^2 \delta_j} \left[ 1 + \frac{|x_1|}{h}\right]^{-2M_1} \times \left[1 + \frac{|x_2|}{ h \delta_j} \right]^{-2M_2}. \label{post-IBP-final}  \end{equation} 
The constant $C_{M_1, M_2, \lambda}$ depends only  on $R_0, \lambda$ and the $C^{M_1 + M_2}$-norms on $\widehat{\psi}$ and $\Phi$; it blows up to $\infty$ as $\lambda \nearrow 1$.
\vskip0.1in
\noindent Assuming \eqref{post-IBP-final} for the moment, the proof of  \eqref{pointwise estimate: strong maximal function} is completed as follows. Let $\mathcal R$ denote the origin-centred, axes-parallel rectangle of dimensions $h$ and $h \delta_j$ in the $x_1$ and $x_2$ directions respectively. Then for any large enough choice of $M_1, M_2$ (e.g. $M_1 = M_2 =3$ suffices), the inequality \eqref{post-IBP-final} yields 
\begin{align*}  \widetilde{\mathfrak B}_{\omega, h}|f|(x) &= \left[ \int_{\mathcal R} + \sum_{r=0}^{\infty} \int_{2^{r+1}\mathcal R \setminus 2^r \mathcal R}
%\int_{\begin{subarray}{c}|x_1| \leq h \\ |x_2| \leq h \delta_j \end{subarray}} + \sum_{r, s=0}^{\infty} 2^{-rM_1 - sM_2)}\int_{\begin{subarray}{c}2^r h \leq |x_1| \leq 2^{r+1}h \\ 2^s h \delta_j \leq |x_2| \leq 2^{s+1}h \delta_j \end{subarray}} 
\right] \bigl| \check{\mathfrak b}_{\omega, h}(y) \bigr| \bigl| f(x+y)\bigr| \, dy \\ 
&\leq \frac{C_{1}}{h^2 \delta_j}  \left[ \int_{\mathcal R} + \sum_{r=0}^{\infty} \int_{2^{r+1}\mathcal R \setminus 2^r \mathcal R} 2^{-2r \min\{M_1, M_2\}}\right] \bigl| f(x+y)\bigr| \, dy 
\\ &\leq  C_1 \left[ 1 + \sum_{r=0}^{\infty} 2^{-2r(\min\{M_1, M_2\}-1)}\right] \mathcal M_{\text{\tiny{str}}} f(x) \leq C_1 \mathcal M_{\text{\tiny{str}}} f(x).   \end{align*} 
Here $C_1$ is a (running) constant depending only on $C_{M_1, M_2, \lambda}$. Taking the supremum of the left hand side with respect to $\omega \in \Omega$ and $h > 0$ leads to the desired conclusion \eqref{pointwise estimate: strong maximal function}.  
\end{proof} 
\noindent It therefore remains to prove the following lemma. 
\begin{lemma} \label{IBP-lemma} 
In the set-up of Lemma \ref{lemma : strong maximal function}, the estimate \eqref{post-IBP-final} holds. 
\end{lemma} 
\begin{proof} Let us write 
\begin{align}
\check{\mathfrak b}_{\omega, h}(x) &:=  \int_{\mathbb R^2} e(x \cdot \xi) \mathfrak b_{\omega, h}(\xi) \, d\xi  =  \int_{\mathbb R^2} e(x \cdot \xi) \widehat{\psi}(h(1, \omega) \cdot \xi)\, \Phi(h \delta_j \xi) \, d\xi \nonumber \\ &= \int_{\mathbb R^2} e(x \cdot \xi)  \widehat{\psi}\bigl( h(\xi_1 + \omega \xi_2)\bigr) \Phi(h \delta_j \xi) d \xi \nonumber \\ &= \frac{1}{h^2 \delta_j}\int_{\mathbb R^2} e \Bigl[ \frac{x_1 \eta_1}{h} + \frac{x_2 \eta_2}{h \delta_j}\Bigr] \widehat{\psi}\Bigl(\eta_1 + \frac{\omega \eta_2}{\delta_j} \Bigr) \Phi(\delta_j \eta_1, \eta_2) \, d\eta.
 \label{pre-IBP}
\end{align} 
The last step is a consequence of the change of variables $\eta_1 = h \xi_1$, $\eta_2 = h \delta_j  \xi_2$. From here on, the proof relies on a fairly standard argument involving integration by parts in the variables $\eta_1$ and $\eta_2$. As we will see, each such application yields a factor of $x_1/h$ or $x_2/(h \delta_j)$ respectively, but also yields additional factors that need to be controlled through $C_{M_1, M_2, \lambda}$.
\vskip0.1in
\noindent Let us fill in the details. It follows from \eqref{psi requirements} and \eqref{defn: Phi} that 
\[ \left\{
\begin{aligned}
&\widehat{\psi}\Bigl(\eta_1 + \frac{\omega \eta_2}{\delta_j} \Bigr)  \text{ and }  \Phi(\delta_j \eta_1, \eta_2) \text{ are smooth functions}, \\ &\text{ with respective supports in the sets} \\   
&\left\{\eta : \Bigl|\eta_1 + \frac{\omega \eta_2}{\delta_j} \Bigr| \leq R_0 \right\} \text{ and } \left\{ \eta \in \mathbb R^2 : |\delta_j \eta_1|^2 + |\eta_2|^2 \leq \bigl(\frac{20 R_0}{c_0} \bigr)^{2}\right\}.  
\end{aligned} 
\right\}
\]
As a result of this, the integrand in \eqref{pre-IBP} is smooth, with support in 
\begin{equation} \label{integrand support} \Bigl\{ \eta: \bigl|\eta_1 + \frac{\omega \eta_2}{\delta_j} \bigr| \leq R_0, \; |\eta_2| \leq \frac{20 R_0}{c_0} \Bigr\}. 
\end{equation}
The domain of integration in \eqref{pre-IBP} is therefore of size at most $20 R_0^2/c_0$.    
%Inserting this into \eqref{pre-IBP}, we arrive at the estimate 
%\begin{equation} 
%\bigl|\check{\mathfrak b}_{\omega, h}(x) \bigr| \leq \frac{C_1}{h^2 \delta_j}, \end{equation}  with $C_1$ depending only on $R_0$, $||\widehat{\psi}||_{\infty}$ and $||\Phi||_{\infty}$. 
\vskip0.1in
\noindent The smoothness and compact support of $\hat{\psi}$ and $\Phi$ allow repeated applications of integration by parts  in \eqref{pre-IBP}. 
%leading to a stronger conclusion (proved in \eqref{post-IBP-final}), namely that $\check{\mathfrak b}_{\omega, h}$ decays rapidly away from the rectangle $\{ |x_1| \leq h, |x_2| \leq h \delta_j\}$. As we will see below, this decay is key to the proof of \eqref{pointwise estimate: strong maximal function}.
%\vskip0.1in
%\noindent Let us focus on the proof of \eqref{post-IBP-final} then. 
Integrating by parts using the differential operator $\mathcal D_{\eta}= (1 + \partial_{\eta_1}^{2M_1}) (1 + \partial_{\eta_2}^{2M_2})$  results in the expression
\begin{align} 
&\left[ 1 + \Bigl(\frac{|x_1|}{h} \Bigr)^{2M_1}\right] \times \left[1 + \Bigl(\frac{|x_2|}{ h \delta_j} \Bigr)^{2M_2}\right] \check{\mathfrak b}_{\omega, h}(x)  
 \label{post-IBP-1} \\
&\quad =  \frac{1}{h^2 \delta_j}\int \mathcal D_{\eta}  \left[ e \left( \frac{x_1 \eta_1}{h} + \frac{x_2 \eta_2}{h \delta_j} \right) \right] \widehat{\psi}\Bigl(\eta_1 + \frac{\omega \eta_2}{\delta_j} \Bigr) \Phi(\delta_j \eta_1, \eta_2) \, d\eta \nonumber \\
&\quad = \frac{1}{h^2 \delta_j}\int e \Bigl[ \frac{x_1 \eta_1}{h} + \frac{x_2 \eta_2}{h \delta_j}\Bigr]  \mathcal D_{\eta}^t \left[\widehat{\psi}\Bigl(\eta_1 + \frac{\omega \eta_2}{\delta_j} \Bigr) \Phi(\delta_j \eta_1, \eta_2) \right] \, d\eta. \label{post-IBP} 
\end{align}
The notation $\mathcal D_{\eta}^t$ represents the adjoint of $\mathcal D_{\eta}$, which in this case is $\mathcal D_{\eta}$ itself. 
\vskip0.1in
\noindent Each derivative in $\eta_1$ occurring in $\mathcal D_{\eta}^{t}[ \cdots]$ yields a factor of either 1 or $\delta_j \leq 1$, depending on the term of the product it acts on. Similarly, each derivative in $\eta_2$ gives rise to the factor 1 or $\frac{\omega}{\delta_j}$. Thus a moment's reflection reveals the following bound on the integrand in \eqref{post-IBP}:
\begin{equation} 
\Bigl|\mathcal D_{\eta}^t \left[\widehat{\psi}\Bigl(\eta_1 + \frac{\omega \eta_2}{\delta_j} \Bigr) \Phi(\delta_j \eta_1, \eta_2) \right] \Bigr| \leq \sum_{\gamma=1}^{C_{M_1, M_2}} \mathfrak f_{\gamma}(\eta), \label{chain rule sum}
\end{equation} 
where each summand $\mathfrak f_{\gamma}(\eta) = \mathfrak f_{\gamma}(\eta; \omega, \delta_j)$ is of the form
\begin{equation}
\label{chain rule}
\mathfrak f_{\gamma}(\eta) = \max(1, \delta_j)^{\alpha} \max \Bigl(1, \frac{\omega}{\delta_j} \Bigr)^{\beta} \mathfrak g_{\gamma}\Bigl( \eta_1 + \frac{\omega \eta_2}{\delta_j}\Bigr) \mathfrak h_{\gamma}(\delta_j \eta_1, \eta_2)
\end{equation}
for some choice of non-negative integers $\alpha = \alpha(\gamma)$ and $\beta = \beta(\gamma)$ obeying $0 \leq \alpha \leq 2M_1$, $0 \leq \beta \leq 2M_2$. The functions $\mathfrak g_{\gamma}$ and $\mathfrak h_{\gamma}$ are absolute values of certain derivatives of $\widehat{\psi}$ and $\Phi$ respectively, of order at most $2M_1+2M_2$. 
\vskip0.1in
\noindent Let us estimate the term in \eqref{chain rule}. In light of the above observation, the functions $\mathfrak g_{\gamma}$ and $\mathfrak h_{\gamma}$ share the same support as $\widehat{\psi}$ and $\Phi$ respectively, and are also uniformly bounded by their respective $C^{2M_1+2M_2}$ norms. Therefore, each function $\mathfrak f_{\gamma}$ given by \eqref{chain rule} is supported in the set given by \eqref{integrand support}. The assumption $\omega \in \Omega_j$ implies
\begin{equation} \frac{\omega}{\delta_j} \leq \frac{a_j}{a_j - a_{j+1}} \leq \frac{1}{1 - \lambda}, \text{ since } a_{j+1} \leq \lambda a_j \text{ by  \eqref{special A}}. \label{importance of lambda}  \end{equation} 
%Using all of this for each integer $M_1, M_2 \geq 0$, and combining with the integrand's support property \eqref{integrand support},  
Combining \eqref{chain rule sum}, \eqref{chain rule} and \eqref{importance of lambda}, the expression in \eqref{post-IBP} can be bounded in absolute value by a constant multiple of  
\begin{align} 
&\frac{\left[ \max(1, \delta_j)\right]^{2M_1}}{h^2 \delta_j}  \Bigl[ \max \Bigl(1, \frac{1}{1-\lambda} \Bigr) \Bigr]^{2M_2} \int_{|\eta_2| \leq R_0} \int_{|\eta_1 + \frac{\omega \eta_2}{\delta_j}| \leq \frac{20R_0}{c_0}} d\eta_1 d\eta_2, \nonumber \\ 
&\text{ which in turn is} \leq \frac{20 R_0^2}{c_0h^2 \delta_j}  (1-\lambda)^{-2M_2} \leq \frac{C_{M_1, M_2, \lambda}}{h^2 \delta_j}. \label{post-IBP+1}
\end{align} 
 %Combining \eqref{pre-IBP}, \eqref{post-IBP-1}, \eqref{post-IBP} and \eqref{post-IBP+1}, we arrive at the estimate 
This is the desired estimate \eqref{post-IBP-final}.
\end{proof} 
\noindent Lemma \ref{IBP-lemma} concludes our analysis of $\mathfrak B^{\ast}$. We tackle $\mathfrak C^{\ast}$ next.  
\section{The high frequency maximal operator $\mathfrak C^{\ast}$} \label{section: Induction C-star}
It is worth noting that the structure of the slope set $\Omega$, as manifested by its finite-order lacunarity, was somewhat implicit in the analysis of $\mathfrak B^{\ast}$. The latter relied exclusively on the special lacunary sequence $A$, and not on any special property of $\Omega_j$. In contrast, the lower order lacunarity of $\Omega_j$ turns out to be key to the analysis of $\mathfrak C^{\ast}$.  
\vskip0.1in
\noindent Let us recall from \eqref{Omegaj}  that $\Omega$ is the union of disjoint subsets $\Omega_j$ contained respectively in the intervals $[a_{j+1}, a_j)$ determined by the special sequence $A = \{a_j : j \geq 0\}$. It follows from \eqref{BC-maximal} and \eqref{b-and-c-multipliers} that
\begin{align} 
\mathfrak C^{\ast}f(x) &= \mathfrak C^{\ast}_{\Omega}f(x) := \sup_{\omega \in \Omega} \sup_{h >0} \bigl|\mathfrak C_{\omega, h} f(x) \bigr| \nonumber \\ &= \sup_{j\geq 0} \sup_{\omega \in \Omega_j} \sup_{h >0} \bigl| \mathfrak C_{\omega, h} f(x) \bigr| = \sup_{j\geq 0} \mathfrak C^{\ast}_{\Omega_j}f(x) 
\end{align} 
The main result concerning $\mathfrak C^{\ast}_{\Omega}$ is the following.  
\begin{proposition} \label{prop: C*} 
Fix an integer $N \geq 1$, a lacunarity constant $\lambda \in (0,1)$ and a Lebesgue exponent $p \in (1, \infty)$. Suppose that there is a constant $C(p, N-1, \lambda) > 0$ such that the norm estimate 
\begin{equation}   \label{Inductive D inequality}
|| D_{\Omega'}||_{p \rightarrow p} 
%= ||\sup_{j \geq 1} \mathfrak C^{\ast}_{\Omega_j}||_{p \rightarrow p} 
\leq C(p, N-1, \lambda)  \; \text{ holds for every } \Omega'\in\Lambda(N-1,\lambda).
\end{equation} 
Then there is a constant $C(p,N,\lambda)>0$ such that
the uniform norm bound  
\begin{equation}  ||\mathfrak C^{\ast}_{\Omega}||_{p \rightarrow p} < C_(p, N, \lambda) \label{Maximal C inequality}
\end{equation} 
%\text{ for any } \Omega \in {\tt{AdFinLac}}.  \] 
holds for any slope set $\Omega \in \Lambda(N, \lambda)$ whose special sequence $A$ obeys \eqref{generic A} and \eqref{special A}.  
\end{proposition}
\vskip0.1in
\noindent {\em{Remarks: }} 
\begin{enumerate}[1.]
\item Proposition \ref{prop: C*} and Lemma \ref{lemma : strong maximal function} combine to prove the inductive step of Theorem \ref{THM: DB}. The reader may head to Section \ref{section: DB conclusion} for the conclusion of the proof of the theorem, given the proposition.  
\vskip0.1in
\item The remainder of this section is given over to a proof sketch of  Proposition \ref{prop: C*}. The details are carried out in the next two sections. 
\end{enumerate} 
\begin{proof} 
For fixed $\lambda \in (0, 1)$, $p \in (1, \infty)$ and $N \geq 1$, suppose that \eqref{Inductive D inequality} holds. 
Let $\Omega \in \Lambda(N, \lambda)$ be a slope set whose special sequence $A$ obeys \eqref{generic A} and \eqref{special A}, and whose portions lying within the gaps of $A$ are denoted by $\Omega_j$, as in \eqref{Omegaj}. Our aim is to apply Theorem \ref{Christ-theorem} from Section \ref{section: Christ almost orthogonality principle} to arrive at \eqref{Maximal C inequality}. To do so, we need to align the present set-up with the framework of that theorem and verify its hypotheses. 
\subsection{Identification of the operators} 
In the notation of Section \ref{section: Christ almost orthogonality principle}, we set for each $j \geq 0$, 
\begin{align*} 
&\mathscr{T}_{j \nu} := \mathfrak C_{\omega, h} \text{ with } \nu = (\omega, h) \in  \mathcal S_j := \Omega_j \times (0, \infty); \text{  this means } \\
&\mathfrak C^{\ast}_{\Omega} = \sup_{j \geq 0} \sup_{\nu \in \mathcal S_j} \bigl|\mathscr{T}_{j \nu}\bigr| = \sup_{(j, \nu) \in \mathcal S^{\ast}} \bigl|\mathscr{T}_{j \nu}\bigr|,  \text{ with }  
\mathcal S^{\ast} = \bigl\{(j, \nu) : j \in \mathbb Z_{\geq 0}, \nu \in \mathcal S_j \bigr\}. 
\end{align*}  
For each $j\geq 0$, we have $\Omega_j\in\Lambda(N-1,\lambda)$ from \eqref{Omegaj}. 
The hypothesis \eqref{Inductive D inequality} therefore gives
\[
        \|D_{\Omega_j}\|_{p\to p}
        \leq C(p,N-1,\lambda)
        \quad\text{uniformly in }j.
\]
Moreover, by \eqref{A=B+C}, \eqref{D-Dtilde-pointwise} and Lemma \ref{lemma : strong maximal function},
\[
\begin{split}
        \mathfrak C_{\Omega_j}^{\ast}f(x) &= \sup_{\omega \in \Omega_j} \sup_{h > 0} \bigl| \mathfrak C_{\omega, h} f(x)\bigr| \\ 
        &\leq \sup_{\omega \in \Omega_j} \sup_{h > 0} \bigl| \tilde{\mathcal A}_{\omega, h} f(x)\bigr| + \sup_{\omega \in \Omega_j} \sup_{h > 0} \bigl| \mathfrak B_{\omega, h} f(x)\bigr| \\
        &\leq \tilde{D}_{\Omega_j}|f|
              + \mathfrak B_{\Omega_j}^{\ast}f 
        \leq C_0D_{\Omega_j}|f|
              +C_1 \mathcal M_{\text{str}}f .
\end{split}
\]
Here $\mathfrak C_{\Omega_j}^{\ast}$ and $\mathfrak B_{\Omega_j}^{\ast}$ are defined as in \eqref{BC-maximal}, with $\Omega$ replaced by $\Omega_j$; in other words, they are the 
restrictions of the parent averaging operators $\mathfrak C_{\omega, h}$ and $\mathfrak B_{\omega, h}$ to slopes $\omega\in\Omega_j$.
In particular, the operators $\mathfrak B_{\omega,h}$ and $\mathfrak C_{\omega,h}$ continue to be defined using the parent gap scale $\delta_j$ from \eqref{decomposition of a}.
\vskip0.1in
\noindent Consequently, the hypothesis \eqref{Inductive D inequality} gives 
\[
        \sup_{j\geq 0}\|\mathfrak C_{\Omega_j}^{\ast}\|_{p\to p}
        \leq C_0C(p,N-1,\lambda)
              +C_1\|\mathcal M_{\text{str}}\|_{p\to p}<\infty.
\]
This verifies hypothesis \eqref{uniform op norm in j} of Theorem \ref{Christ-theorem}.
%With this correspondence, the induction hypothesis \eqref{COmegaj-induction-hypothesis} is the analogue of the  hypothesis \eqref{uniform op norm in j}, ensuring $L^p$ control over the maximal operator $\sup_{\nu \in \mathcal S_j} |\mathscr{T}_{j \nu}|$, uniformly in $j$. 
\subsection{Verifying essential positivity} \label{section: ess pos} In order to apply Theorem \ref{Christ-theorem} in the range $p \in (1, 2)$, we need to confirm that the operator family \begin{equation} \label{proving ess pos}
\pmb{\mathscr{T}} = \{ \mathscr{T}_{j \nu}: (j, \nu) \in \mathcal S^{\ast}\}  \text{ with } \mathscr{T}_{j \nu} = \mathfrak C_{\omega, h} \text{ is essentially positive}, 
\end{equation}  
in the sense described on page \pageref{ess-pos-1}. This requires us to write $\mathscr{T}_{j \nu} = \mathfrak C_{\omega, h}$ as a difference of two operators, for which the three conditions \eqref{ess-pos-1}--\eqref{ess-pos-3} can be verified. 
\vskip0.1in
\noindent Towards this goal,  let us express $\mathscr{T}_{j \nu}$ as 
\[ \mathscr{T}_{j \nu} = \mathfrak C_{\omega, h} = \tilde{\mathcal A}_{\omega, h} - \mathfrak B_{\omega, h} = \mathscr{A}_{j \nu} - \mathscr{B}_{j \nu}, \text{ with } \mathscr{A}_{j \nu} := \tilde{\mathcal A}_{\omega, h} \text{ and  } \mathscr{B}_{j \nu} = \mathfrak B_{\omega, h}. \]  The non-negativity of $\psi$, mentioned in \eqref{psi requirements}, implies that for $0 \leq |f| \leq g$, 
\[
\tilde{\mathcal A}_{\omega, h}f(x) = \left|\int f(x + t(1, \omega)) \psi\bigl(\frac{t}{h} \bigr) \, dt \right| \leq \tilde{\mathcal A}_{\omega, h}g(x), \text{ confirming \eqref{ess-pos-1}}. 
\]
\vskip0.1in 
\noindent We turn now to the verification of \eqref{ess-pos-2}, which we claim holds with
\[ \mathscr{B}_{j \nu} = \mathfrak B_{\omega, h} \text{ given by \eqref{kernel check b}} \; \text{ and } \; \mathscr{C}_{j \nu} = \widetilde{\mathfrak B}_{\omega, h} \text{ given by } \eqref{newB}. \] Indeed the first inequality in \eqref{ess-pos-2} is verified in \eqref{B and new B}. 
Since $\widetilde{\mathfrak B}_{\omega, h}$ is a convolution operator with a non-negative kernel $|\check{\mathfrak b}_{\omega, h}|$, it obeys the monotonicity property:
\[ 0 \leq \widetilde{\mathfrak B}_{\omega, h}f(x) \leq \widetilde{\mathfrak B}_{\omega, h}g(x) \text{ for } 0 \leq f \leq g. \]
This confirms the second condition in \eqref{ess-pos-2}. 
\vskip0.1in
\noindent It remains to verify \eqref{ess-pos-3}. According to the conclusion \eqref{pointwise estimate: strong maximal function} of Lemma \ref{lemma : strong maximal function}, 
\[ \sup_{j \nu} \mathscr{C}_{j, \nu}|f| (x) = \sup_{\begin{subarray}{c}\omega \in \Omega \\ h > 0 \end{subarray}}  \widetilde{\mathfrak B}_{\omega, h} |f| (x) \leq C_1 \mathcal M_{\text{str}} f(x).   \] 
The desired conclusion \eqref{ess-pos-3} now follows from the $L^p$-boundedness of $\mathcal M_{\text{str}}$. This verifies that the operator family in \eqref{proving ess pos} is essentially positive. 
\subsection{Verifying the square function estimate \eqref{R-ORTHOGONAL}}  \label{section:sq function} Our next task is to identify auxiliary operators $\mathscr{R}_j$ satisfying the square function estimate \eqref{R-ORTHOGONAL}, for which $\mathscr{T}_{j\nu}\mathscr{R}_j$ approximates $\mathscr{T}_{j \nu}$ in the sense of \eqref{remainder op norm}. We define these operators using Fourier localization in conical sectors, as follows.  Let $\Delta_j$ and $\tilde{\Delta}_j$ denote the bivariate cones
\begin{align}
\Delta_j &:= \left\{ (\xi_1, \xi_2) \in \mathbb R^2 \setminus \{0\} : a_{j+1} - c_0 \delta_j \leq \frac{\xi_1}{-\xi_2} < a_j + c_0  \delta_j \right\}, \label{defn: cone Deltaj} \\
\tilde{\Delta}_j &:= \left\{ (\xi_1, \xi_2) \in \mathbb R^2 \setminus \{0\} : a_{j+1} - 2c_0 \delta_j \leq \frac{\xi_1}{-\xi_2} < a_j + 2c_0  \delta_j \right\}. \label{defn: cone Deltaj tilde} 
\end{align}  
Here $c_0$ is the constant chosen in \eqref{little c0}, and $\delta_j := a_{j} - a_{j+1}$ denotes the gap length between two consecutive elements of the lacunary sequence $A$. As we will see in Sections \ref{square function estimate section} and \ref{section: cone geometry}, this choice of $c_0$ ensures two key geometric properties for $\Delta_j$ and $\tilde{\Delta}_j$: 
\vskip0.1in 
 \begin{itemize} 
\item  First, for all $h > 0$ and $\omega \in \Omega_j$, the Fourier support for $\mathfrak C_{\omega, h}$ is contained in the conical sector $\Delta_j$.
\vskip0.1in
\item Second, the cones $\tilde{\Delta}_j$ are essentially disjoint.
\vskip0.1in
\end{itemize}
 Figure \ref{fig: conical sectors} gives a visual depiction of $\Delta_j$ and $\tilde{\Delta}_j$ for a fixed $j$.
 \vskip0.1in
\noindent For $j \geq 0$, we choose a function $\zeta_j$ that is homogeneous of degree 0, 
\begin{equation} \label{def: zeta_j}
\zeta_j \in C^{\infty}(\mathbb R^2 \setminus  \{0\}), \quad 
\zeta_j \equiv 1 \text{ on } \Delta_j, \quad  \zeta_j \equiv 0 \text{ outside } \tilde{\Delta}_j.
\end{equation} 
To ensure uniformly controlled choice of derivative estimates required to apply the Marcinkiewicz multiplier theorem later on, we choose the function $\zeta_j$ in the following way. Fix $\eta\in C_c^\infty(\mathbb R)$ such that
\begin{align} 
&\eta \equiv 1 \text{ on } [-c_0,1+c_0] \; \text{ and } \; \operatorname{supp}\eta\subseteq [-2c_0,1+2c_0], \label{ref map eta for cones} \\
&\zeta_j(\xi_1,\xi_2)
=\eta\!\left(\frac{-\xi_1/\xi_2-a_{j+1}}{\delta_j}\right)
\; \text{ for } \; \xi_2\ne0. \label{zetaj for cones}
\end{align} 
For $\xi_2 = 0$ and $\xi_1 \neq 0$, we set $\zeta_j(\xi) =0$. This choice has the properties in \eqref{def: zeta_j}; its uniform multiplier bounds will be used in Section \ref{square function estimate section}.
\vskip0.1in
\noindent The operator $\mathscr{R}_j$ corresponds to the Fourier multiplier $\zeta_j$:
\begin{equation} \label{Rj conical multipliers}
\left[ \mathscr{R}_j f\right]^{\wedge}(\xi) := \zeta_j(\xi) \widehat{f}(\xi). 
\end{equation}
 In Section \ref{square function estimate section} below, we verify the required square function estimate \eqref{R-ORTHOGONAL} for the operators $\mathscr{R}_j$ as in \eqref{Rj conical multipliers}. 
 \begin{center}
\begin{figure}
\begin{tikzpicture}[x=1cm,y=1cm]

  %------------------------------------------------------------
  % Numeric plotting parameter for c_0.
  % For BOTH inner and outer cones to stay in quadrants II and IV,
  % keep 0 < c_0 < 1/2.
  %------------------------------------------------------------
  \pgfmathsetmacro{\czeroval}{0.20}

  % Fixed data: a_{j+1}=1/4, a_j=1/2, delta_j=1/4
  \pgfmathsetmacro{\ajplusone}{0.25}
  \pgfmathsetmacro{\aj}{0.50}
  \pgfmathsetmacro{\deltaj}{\aj-\ajplusone}

  % Inner boundary parameters
  \pgfmathsetmacro{\mlow}{\ajplusone-\czeroval*\deltaj}
  \pgfmathsetmacro{\mhigh}{\aj+\czeroval*\deltaj}

  % Outer boundary parameters
  \pgfmathsetmacro{\mOuterLow}{\ajplusone-3.5*\czeroval*\deltaj}
  \pgfmathsetmacro{\mOuterHigh}{\aj+3.5*\czeroval*\deltaj}

  % Visible radial reach
  \pgfmathsetmacro{\Ymax}{3.20}

  \coordinate (O) at (0,0);

  % Inner cone points
  \coordinate (QIVlow)   at ({\mlow*\Ymax},{-\Ymax});
  \coordinate (QIVhigh)  at ({\mhigh*\Ymax},{-\Ymax});
  \coordinate (QIIlow)   at ({-\mlow*\Ymax},{\Ymax});
  \coordinate (QIIhigh)  at ({-\mhigh*\Ymax},{\Ymax});

  % Outer cone points
  \coordinate (QIVlowT)  at ({\mOuterLow*\Ymax},{-\Ymax});
  \coordinate (QIVhighT) at ({\mOuterHigh*\Ymax},{-\Ymax});
  \coordinate (QIIlowT)  at ({-\mOuterLow*\Ymax},{\Ymax});
  \coordinate (QIIhighT) at ({-\mOuterHigh*\Ymax},{\Ymax});

  % Axes
  \draw[-{Stealth[length=2.6mm]}, very thick] (-4.5,0) -- (4.8,0)
    node[right] {$\xi_1$};
  \draw[-{Stealth[length=2.6mm]}, very thick] (0,-3.8) -- (0,3.8)
    node[above] {$\xi_2$};

  % Outer cone first (lighter shading)
  \fill[blue!55!black, opacity=.12] (O) -- (QIVlowT) -- (QIVhighT) -- cycle;
  \fill[blue!55!black, opacity=.12] (O) -- (QIIlowT) -- (QIIhighT) -- cycle;

  % Inner cone on top (darker shading)
  \fill[blue!55!black, opacity=.28] (O) -- (QIVlow) -- (QIVhigh) -- cycle;
  \fill[blue!55!black, opacity=.28] (O) -- (QIIlow) -- (QIIhigh) -- cycle;

  % Outer boundary rays
  \draw[blue!70!black, thin]  (O) -- ($(O)!1.05!(QIVlowT)$);
  \draw[blue!70!black, thin]  (O) -- ($(O)!1.05!(QIVhighT)$);
  \draw[blue!70!black, thin]  (O) -- ($(O)!1.05!(QIIlowT)$);
  \draw[blue!70!black, thin]  (O) -- ($(O)!1.05!(QIIhighT)$);

  % Inner boundary rays
  \draw[blue!70!black, thick] (O) -- ($(O)!1.05!(QIVlow)$);
  \draw[blue!70!black, thick] (O) -- ($(O)!1.05!(QIVhigh)$);
  \draw[blue!70!black, thick] (O) -- ($(O)!1.05!(QIIlow)$);
  \draw[blue!70!black, thick] (O) -- ($(O)!1.05!(QIIhigh)$);

  % Labels for the regions
  \node[blue!70!black] at (0.7,-1.95) {$\Delta_j$};
  \node[blue!70!black] at (-0.7,1.95) {$\Delta_j$};

  % Place \widetilde{\Delta}_j in the outer-only band near the upper boundary
  \pgfmathsetmacro{\mBand}{0.5*(\mhigh+\mOuterHigh)}
  \node[blue!70!black] at ({\mBand*2.85},{-2.7}) {$\widetilde{\Delta}_j$};
  \node[blue!70!black] at ({-\mBand*2.85},{ 2.7}) {$\widetilde{\Delta}_j$};

  % Minimal inner-boundary labels; the table in the text carries the rest.
  %\node[fill=white, inner sep=1.5pt]
    %at ($(O)!0.76!(QIVlow)+(0.77,-0.18)$)
    %{\scriptsize$\frac14-\frac{c_0}{4}$};

  %\node[fill=white, inner sep=1.5pt]
    %at ($(O)!0.68!(QIVhigh)+(0.86,0.17)$)
    %{\scriptsize$\frac12+\frac{c_0}{4}$};

\end{tikzpicture}
\caption{Sectors $\Delta_j$ and $\tilde{\Delta}_j$} \label{fig: conical sectors} 
\end{figure}
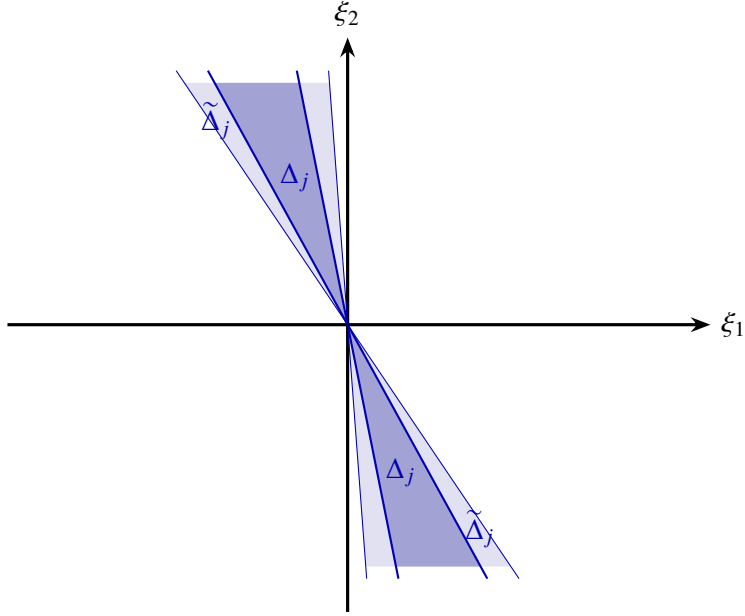
\end{center} 

%The square function estimate \eqref{R-ORTHOGONAL} for operators $\mathscr{R}_j$ is well-known; a proof can be found in 
%{\cite[equation 4]{NagelSteinWainger}}, following the arguments in \cite{Cordoba-Fefferman}. 
\subsection{Checking \eqref{remainder op norm}} \label{section: error} The geometric significance of $\mathscr{R}_j$ lies in the following claim: for $\omega \in \Omega_j$ and $h > 0$, we claim the relation
\begin{equation} \label{GEOMETRIC STATEMENT FOR CONES}
\mathfrak C_{\omega, h} \mathscr{R}_j \equiv \mathfrak C_{\omega, h};  \text{ in other words, } \; \mathscr{T}_{j \nu} \equiv \mathscr{T}_{j \nu} \mathscr{R}_j. 
\end{equation} 
In particular, the estimate \eqref{remainder op norm} holds trivially, since the left hand side of the inequality is identically zero. We will prove \eqref{GEOMETRIC STATEMENT FOR CONES} in Section \ref{section: cone geometry} below. 
\vskip0.1in
\noindent Assuming for now the two statements from Sections \ref{section:sq function} and \ref{section: error},
we may apply the conclusion \eqref{doubly indexed maximal inequality} of Theorem \ref{Christ-theorem}, which yields the
claimed norm bound \eqref{Maximal C inequality}. The resulting constant depends only
on $p$, $\lambda$, $C(p,N-1,\lambda)$ and the fixed auxiliary
constants appearing in \eqref{uniform op norm in j}--\eqref{ess-pos-3}; in particular, it
is uniform over all $\Omega\in\Lambda(N,\lambda)$. This proves
Proposition \ref{prop: C*}, subject only to the square-function estimate
\eqref{R-ORTHOGONAL} and the localization identity \eqref{GEOMETRIC STATEMENT FOR CONES}, which are
established in Sections~\ref{square function estimate section} and~\ref{section: cone geometry}, respectively.
\end{proof}

\section{Proof of the square function estimate \eqref{R-ORTHOGONAL}}
% with $\mathscr R_j$ as in \eqref{Rj conical multipliers}} 
\label{square function estimate section}
\begin{proof}  
In \cite[p.424-425]{Cordoba-Fefferman1} Cordoba and Fefferman introduced an argument involving Rademacher functions and the Marcinkiewicz multiplier theorem that has become a standard tool for deriving square function estimates for conical multipliers, such as \eqref{R-ORTHOGONAL}; in fact, Nagel, Stein and Wainger {\cite[equation 4]{NagelSteinWainger}} used this argument to prove that $D_{\Omega}$ is $L^p$ bounded for $\Omega \in \Lambda(1, \lambda)$. The key observation in \cite{{NagelSteinWainger}, {Cordoba-Fefferman1}} is  that a square function estimate of the form \eqref{R-ORTHOGONAL} for $\mathscr{R}_j$ as in \eqref{Rj conical multipliers} follows from a finitely overlapping property of the Fourier supports of $\mathscr R_j$. 
\vskip0.1in
\noindent In the current context, to use the standard Rademacher-function and Marcinkiewicz-multiplier argument cited above, we need to prove the following claim:  there exists a large integer constant $K_0 \geq 1$ obeying 
\begin{equation} \label{big C0}
\frac{1 + \lambda^{K_0}}{1 - \lambda^{K_0}} < 4
\end{equation} 
such that for any lacunary sequence $A$ obeying  \eqref{special A}, the conical sectors $\tilde{\Delta}_j$ given by \eqref{defn: cone Deltaj tilde} obey the following essential disjointness property: 
\begin{equation} \label{finitely overlapping-1}
\sup_{j \geq 0}\# \left\{ j'  \geq 0 : \tilde{\Delta}_j \cap \tilde{\Delta}_{j'} \neq \emptyset \right\} \leq 2K_0+1. 
\end{equation} 
In other words, for every $j \geq 0$, the sector $\tilde{\Delta}_j$ can intersect with at most $(2K_0+1)$ conical sectors $\tilde{\Delta}_{j'}$, $j' \geq 0$. The relation \eqref{finitely overlapping-1} in turn follows from a stronger statement: if $(\tilde{\Delta}_j, \tilde{\Delta}_{j'})$ is an intersecting pair, then the indices $j,j'$ must be proximal. Specifically,   
\begin{equation} \label{finitely overlapping} 
\left\{k \in \mathbb Z : j + k \geq 0, \;\tilde{\Delta}_j \cap \tilde{\Delta}_{j+k} \neq \emptyset  \right\} \subseteq \{ -K_0, \ldots, K_0\} \; \text{ for all } j \geq 0.  
\end{equation} 
To establish \eqref{finitely overlapping}, it suffices to show that 
\begin{equation}  \tilde{\Delta}_j \cap \tilde{\Delta}_{j+k} = \emptyset  \text{ for } k > K_0. \label{k>K0}
\end{equation}  Indeed, the case $k<-K_0$ follows by interchanging the indices
$j$ and $j+k$. Since the extremal slopes of the conical sectors $\tilde{\Delta}_j$ form a monotone sequence in $j$, proving \eqref{k>K0} is tantamount to proving the following condition: 
\begin{equation} \label{disjoint intervals} 
a_{j+k} + 2c_0 \delta_{j+k} < a_{j+1} - 2c_0 \delta_j \quad \text{ if } k > K_0.  
 \end{equation} 
See Figure \ref{fig: disjointness of sectors}.
 \vskip0.1in
 \noindent We set about verifying the inequality \eqref{disjoint intervals}. Upon rearrangement, we find that it is equivalent to 
\begin{equation}
2c_0  (\delta_j + \delta_{j+k}) < a_{j+1} - a_{j+k}. \label{disjointness-2} 
\end{equation} 
The left-hand side of \eqref{disjointness-2} equals $2c_0 (a_j - a_{j+1} + a_{j+k} - a_{j+k+1})$,  which is  
\begin{equation} \leq 2c_0 (a_j + a_{j+k}) \leq 2c_0 a_j \Bigl( 1 + \frac{a_{j+k}}{a_j} \Bigr) \leq 2c_0a_j (1 + \lambda^{K_0}), \label{lhs-upper} \end{equation} 
using the lacunarity condition \eqref{special A} of $A$ and the assumption $k > K_0$. On the other hand, the right-hand side of \eqref{disjointness-2} equals
\begin{equation}
 a_{j+1} \left( 1 - \frac{a_{j+k}}{a_{j+1}}\right), \text{ which is }  \geq a_{j+1} (1 - \lambda^{K_0}), \text{ also from \eqref{special A}.}  \label{rhs-lower}
\end{equation} 
Thus, in order to prove  \eqref{disjointness-2}, and therefore \eqref{disjoint intervals}, it suffices to verify that 
\[
2c_0 a_j (1 + \lambda^{K_0}) < a_{j+1} (1 - \lambda^{K_0}). \] 
Since $a_{j+1} > \lambda^2 a_j$ from \eqref{special A}, this would follow from  
\begin{equation} 
2c_0  \frac{1 + \lambda^{K_0}}{1- \lambda^{K_0}} < \lambda^2.   \label{disjointness-final}
\end{equation}  
Our choices of $c_0$ and $K_0$ from \eqref{little c0} and \eqref{big C0} respectively show that the left-hand side of the inequality \eqref{disjointness-final} can be estimated as: 
\[ 2c_0  \frac{1 + \lambda^{K_0}}{1- \lambda^{K_0}}  <8c_0 < \lambda^2. \] 
Thus \eqref{disjointness-final} holds, proving \eqref{finitely overlapping} and hence the bounded overlap \eqref{finitely overlapping-1}. 
\vskip0.1in
\noindent To obtain \eqref{R-ORTHOGONAL}, we also require the uniform (in $j$) derivative bounds for the cut-off functions $\zeta_j$ chosen in \eqref{zetaj for cones}. Indeed,
\begin{align}
&\delta_j\geq(1-\lambda)a_j  \text{ by \eqref{decomposition of a} and \eqref{special A}; this gives } \nonumber \\
&\bigl| \partial_r^{k} \zeta_j| \leq \|\eta^{(k)} \|_{\infty} \delta_j^{-k} \leq C_{\lambda, k} a_j^{-k} \text{ for } k \geq 1. \\ 
&\text{Additionally, } \bigl|-\frac{\xi_1}{\xi_2} \bigr|\leq C_{\lambda} a_j \text{ on
$\operatorname{supp}\zeta_j$}. 
\end{align}  Combining these estimates via the chain rule we arrive at the following bound uniformly in $j$:
\[
\sup_j \sup_{\xi \neq 0}|\xi_1|^{\ell_1} |\xi_2|^{\ell_2}
\bigl|\partial_{\xi_1}^{\ell_1}\partial_{\xi_2}^{\ell_2}
\zeta_j(\xi)\bigr| \leq C_{\ell_1, \ell_2, \lambda}
\qquad \ell_1, \ell_2 \in \mathbb Z_{\geq 0}.
\]
By \eqref{finitely overlapping-1}, only a bounded number of these cut-off functions are non-zero at
any point. Thus every signed sum $\sum_j\varepsilon_j\zeta_j$,
$\varepsilon_j\in\{-1,1\}$, obeys the same Marcinkiewicz bounds,
uniformly in the signs and in the number of terms. The Marcinkiewicz
multiplier theorem and the standard Rademacher-function argument
now yield \eqref{R-ORTHOGONAL}.
\end{proof} 
% Preamble:
% \usepackage{tikz}
% \usetikzlibrary{arrows.meta,calc}

\begin{center}
\begin{figure}
\begin{tikzpicture}[scale=1.05, >=Stealth]

% Axes
\draw[->, very thick] (-4.5,0) -- (4.8,0) node[right] {$\xi_1$};
\draw[->, very thick] (0,-3.15) -- (0,3.25) node[above] {$\xi_2$};

% Origin circle
\draw[very thick, fill=white] (0,0) circle (0.075);

% ------------------------------------------------------------
% Blue cone: \widetilde{\Delta}_j
% ------------------------------------------------------------

\coordinate (BupA) at (-3.35,0.85);
\coordinate (BupB) at (-2.65,2.05);

\coordinate (BdnA) at (3.35,-0.85);
\coordinate (BdnB) at (2.65,-2.05);

\fill[blue!12] (0,0) -- (BupA) -- (BupB) -- cycle;
\fill[blue!12] (0,0) -- (BdnA) -- (BdnB) -- cycle;

\draw[blue!75!black, very thick] (0,0) -- (BupA);
\draw[blue!75!black, very thick] (0,0) -- (BupB);
\draw[blue!75!black, very thick] (0,0) -- (BdnA);
\draw[blue!75!black, very thick] (0,0) -- (BdnB);

\node[blue!75!black] at (-2.25,1.00) {$\widetilde{\Delta}_j$};
\node[blue!75!black] at (1.95,-1.05) {$\widetilde{\Delta}_j$};

% Dashed extension and formula
\draw[blue!75!black, dashed, thick] (BdnA) -- (3.95,-1.00);

\node[blue!75!black, scale=0.72, anchor=west] at (4.02,-1.03)
{$\displaystyle \frac{\xi_1}{-\xi_2}=a_j-2c_0\delta_j$};

% ------------------------------------------------------------
% Purple cone: \widetilde{\Delta}_{j+k}
% ------------------------------------------------------------

\coordinate (PupA) at (-1.35,2.75);
\coordinate (PupB) at (-0.42,3.00);

\coordinate (PdnA) at (1.35,-2.75);
\coordinate (PdnB) at (0.42,-3.00);

\fill[purple!13] (0,0) -- (PupA) -- (PupB) -- cycle;
\fill[purple!13] (0,0) -- (PdnA) -- (PdnB) -- cycle;

\draw[purple!75!black, very thick] (0,0) -- (PupA);
\draw[purple!75!black, very thick] (0,0) -- (PupB);
\draw[purple!75!black, very thick] (0,0) -- (PdnA);
\draw[purple!75!black, very thick] (0,0) -- (PdnB);

\node[purple!75!black] at (-0.62,2.15) {$\widetilde{\Delta}_{j+k}$};
\node[purple!75!black] at (0.72,-2.10) {$\widetilde{\Delta}_{j+k}$};

% Dashed extension and horizontal formula
\draw[purple!75!black, dashed, thick] (PdnA) -- (1.62,-3.30);

\node[purple!75!black, scale=0.72, anchor=west] at (1.05,-3.48)
{$\displaystyle \frac{\xi_1}{-\xi_2}=a_{j+k}+2c_0\delta_{j+k}$};

\end{tikzpicture}
\caption{\small{Disjointness of conical sectors $\tilde{\Delta}_j$ and $\tilde{\Delta}_{j+k}$ for $|k| > K_0$}} \label{fig: disjointness of sectors}  
\end{figure}
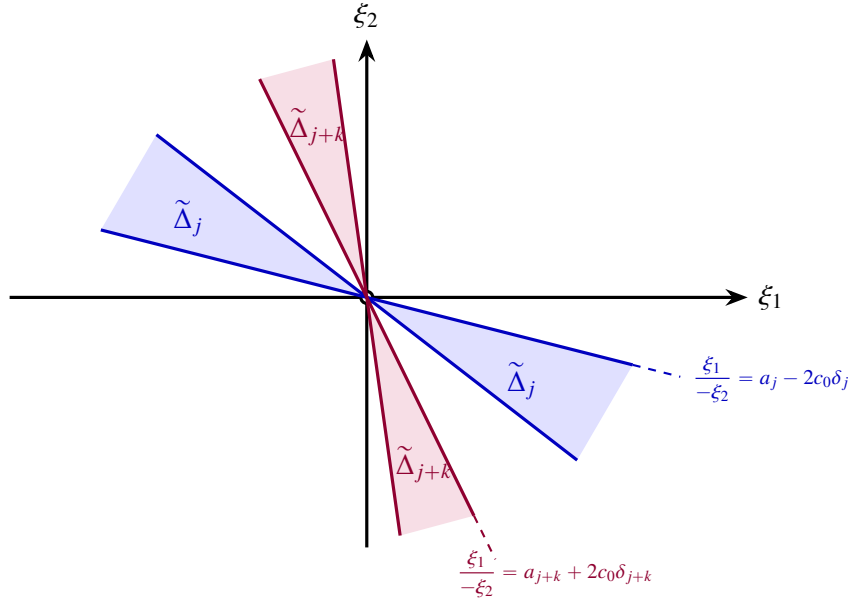
\end{center}
\section{Proof of the localization identity \eqref{GEOMETRIC STATEMENT FOR CONES}} \label{section: cone geometry} 
\begin{proof} 
Let us recall from \eqref{b-and-c-multipliers} and \eqref{Rj conical multipliers} that $\mathfrak C_{\omega, h}$ and $\mathscr{R}_j$ are convolution operators, with 
\[ \Bigl[ \mathfrak C_{\omega, h} \mathscr{R}_jf \Bigr]^{\wedge} (\xi) = \mathfrak c_{\omega, h}(\xi) \zeta_j (\xi) \widehat{f}(\xi) \quad \text{ and } \quad \bigl[ \mathfrak C_{\omega, h}f  \bigr]^{\wedge}(\xi) = \mathfrak c_{\omega, h}(\xi) \widehat{f}(\xi). \] 
To prove the claim \eqref{GEOMETRIC STATEMENT FOR CONES}, which asserts equality of the two operators in the display above, we need to show that their multipliers are identical; in other words, we aim to show that 
\[ \zeta_j \equiv 1 \text{ on the support of } \mathfrak c_{\omega, h}, \text{ for } \omega \in \Omega_j. \] In view of the defining properties \eqref{def: zeta_j} of $\zeta_j$, the last statement follows from the set inclusions 
\begin{align} 
\text{supp}\bigl(\mathfrak c_{\omega, h} \bigr) &\subseteq \left\{\xi \in \mathbb R^2 :  |\xi| \geq \frac{10 R_0}{c_0 h \delta_j}, \; |\xi_1 + \omega \xi_2| \leq \frac{R_0}{h} \right\} \label{containment-1} \\
&\subseteq \Delta_j \subseteq \left\{ \xi \in \mathbb R^2 : \zeta_j(\xi) = 1 \right\}. \label{containment-2} 
\end{align} 
See Figure \ref{fig: localization identity} for a visual depiction of this inclusion.
\vskip0.1in 
\noindent According to \eqref{multipliers b and c}, $\mathfrak c_{\omega, h}$ is a product of two factors $\mathfrak a_{\omega, h}$ and $(1 - \Phi(h \delta_j \xi))$.  Therefore $\text{supp}\bigl(\mathfrak c_{\omega, h} \bigr) $ is contained in the intersection of the supports of these two factors. The first factor $\mathfrak a_{\omega, h}$ is localized on the infinite strip \eqref{infinite strip}, whereas the second factor $1 - \Phi(h \delta_j \xi)$ vanishes identically on $|\xi| \leq 10 R_0/(c_0 h \delta_j)$, by virtue of \eqref{defn: Phi}. This establishes the inclusion \eqref{containment-1}.  
\vskip0.1in
\noindent Let us address \eqref{containment-2}. We first show that the lower bound on $|\xi|$ implies a similar one for $|\xi_2|$. Namely, for $\xi$ in the right side set of \eqref{containment-1}, we have the relation 
\begin{align*}
\left(\frac{10 R_0}{c_0 h \delta_j} \right)^2 &\leq |\xi|^2 =  \xi_1^2 + \xi_2^2 
= (\xi_1 + \omega \xi_2 - \omega \xi_2)^2 + \xi_2^2 \\ 
&\leq 2 (\xi_1 + \omega \xi_2)^2 + (2 \omega^2 + 1) \xi_2^2 \leq \frac{2R_0^2}{h^2} + 3 |\xi_2|^2, 
\end{align*} the last step using $\omega \in [0,1]$. The displayed inequality above implies 
\begin{equation} |\xi_2|^2 \geq \frac{R_0^2}{3h^2} \Bigl( \frac{100}{c_0^2 \delta_j^2}-2\Bigr), \text{ which in turn is } \geq  \frac{R_0^2}{c_0^2 \delta_j^{2} h^{2}}. \label{xi2 is large}\end{equation} 
%for a constant $c_0 > 0$ specified as in \eqref{little c0}.  
The last inequality in \eqref{xi2 is large} follows from the relation 
\[\frac{100}{c_0^2 \delta_j^2} - 2 \geq \frac{3}{c_0^2 \delta_j^2}, \text{ a consequence of } c_0 \delta_j < 1. \]
%can be justified using the lacunarity of the special sequence $A$:
%\[
%\delta_j^{2} = (a_j - a_{j+1})^2 \leq a_j^2 \leq \lambda^{2j} a_0 \leq \lambda^2 \leq 1 - c_0^2, \text{ which implies } 
%\delta_j^{-2} - 1  \geq c_0^2 \delta_j^{-2}.
%\]
Combining the lower bound on $|\xi_2|$ from \eqref{xi2 is large} with the assumption $|\xi_1 + \omega \xi_2| \leq R_0/h$ from \eqref{containment-1}, we arrive at  
\[ \left| \frac{\xi_1}{\xi_2} + \omega \right| \leq \frac{R_0}{h |\xi_2|} \leq c_0 \delta_j \quad \text{ or } \quad \omega - c_0 \delta_j \leq \frac{\xi_1}{- \xi_2} \leq \omega + c_0 \delta_j. \] 
Since $\omega \in \Omega_j \subseteq [a_{j+1}, a_j)$, the last string of inequalities leads to the defining relation \eqref{defn: cone Deltaj} of $\Delta_j$, completing the proof of \eqref{containment-2}. 
\end{proof} 
\begin{center}
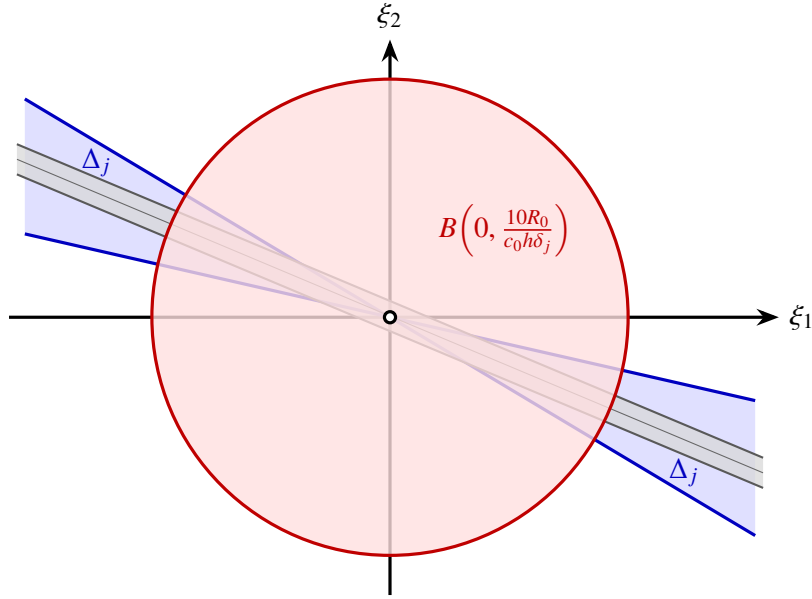
\begin{figure}
\begin{tikzpicture}[scale=1.05, >=Stealth]

%------------------------------------------------------------
% Axes
%------------------------------------------------------------

\draw[->, very thick] (-4.8,0) -- (4.9,0) node[right] {$\xi_1$};
\draw[->, very thick] (0,-3.5) -- (0,3.5) node[above] {$\xi_2$};

%------------------------------------------------------------
% Blue cone \Delta_j
%------------------------------------------------------------

\coordinate (C1) at (-4.6,1.05);
\coordinate (C2) at (-4.6,2.75);

\coordinate (C3) at (4.6,-1.05);
\coordinate (C4) at (4.6,-2.75);

\fill[blue!12] (0,0) -- (C1) -- (C2) -- cycle;
\fill[blue!12] (0,0) -- (C3) -- (C4) -- cycle;

\draw[blue!75!black, very thick] (0,0) -- (C1);
\draw[blue!75!black, very thick] (0,0) -- (C2);
\draw[blue!75!black, very thick] (0,0) -- (C3);
\draw[blue!75!black, very thick] (0,0) -- (C4);

\node[blue!75!black] at (-3.7,1.95) {$\Delta_j$};
\node[blue!75!black] at (3.7,-1.95) {$\Delta_j$};

%------------------------------------------------------------
% Grey strip
%------------------------------------------------------------

\coordinate (S1) at (-4.7,2.18);
\coordinate (S2) at (4.7,-1.77);
\coordinate (S3) at (4.7,-2.15);
\coordinate (S4) at (-4.7,1.80);

\fill[gray!30,opacity=.75] (S1)--(S2)--(S3)--(S4)--cycle;

\draw[gray!70!black,thick] (S1)--(S2);
\draw[gray!70!black,thick] (S4)--(S3);

% Optional centre line of strip
\draw[black!60,thin] (-4.7,1.99)--(4.7,-1.96);

%------------------------------------------------------------
% Red ball
%------------------------------------------------------------

\fill[red!12, opacity=.8] (0,0) circle (3.0);

\draw[red!75!black, very thick] (0,0) circle (3.0);

% Label for the ball
\node[red!75!black] at (1.45,1.10)
{$B\!\left(0,\frac{10R_0}{c_0h \delta_j}\right)$};

%------------------------------------------------------------
% Origin (drawn last)
%------------------------------------------------------------

\draw[very thick,fill=white] (0,0) circle (0.07);

\end{tikzpicture}
\caption{\small{The geometry behind the localization identity $\mathfrak C_{\omega, h} \mathscr{R}_j \equiv \mathfrak C_{\omega, h}$. The function $1 - \Phi(h \delta_j \xi)$ is supported outside the red ball of radius $10R_0/(c_0 h \delta_j)$. The grey strip outside the ball supports the multiplier $\mathfrak c_{\omega, h}$. This is contained in the cone $\Delta_j$, where the multiplier $\zeta_j$ of $\mathscr{R}_j$ is identically 1.}} \label{fig: localization identity}
\end{figure} 
\end{center}

\section{Conclusion of the proof of Theorem \ref{THM: DB}} \label{section: DB conclusion} 
\begin{proof} 
We prove \eqref{D-Omega-bounded} by induction on the lacunarity order $N$.
\vskip0.1in
\noindent Suppose first that $N=0$. By Definition \ref{defn: Lacunary sets}, a set
$\Omega\in\Lambda(0,\lambda)$ is either empty or a singleton. The empty case is trivial, so let $\Omega=\{\omega\}$. Consider the linear transformation
\[
        T_\omega(x_1,x_2):=(x_1,x_2-\omega x_1),
\]
which has determinant one and maps the direction $(1,\omega)$ to $(1,0)$. 
\begin{align*} 
&\text{If } f_\omega:=f\circ T_\omega^{-1}, \; \text{ then } \mathcal A_{\omega,h}f(x)
        =\mathcal A_{0,h}f_\omega(T_\omega x), \\ 
        &\text{and hence } \;        
        \|D_{\{\omega\}}\|_{p\to p}
        =\|D_{\{0\}}\|_{p\to p}.
\end{align*} 
The operator $D_{\{0\}}$ is the one-dimensional centred
Hardy-Littlewood maximal operator in the first coordinate, tensored
with the identity in the second coordinate. It is therefore bounded
on $L^p(\mathbb R^2)$ for every $p\in(1,\infty)$. This proves
\eqref{D-Omega-bounded} when $N=0$.
\vskip0.1in
\noindent We proceed to the inductive step. Fix $N\geq 1$ and suppose that
\eqref{D-Omega-bounded} holds at lacunarity order $N-1$; thus
\[
        \|D_{\Omega'}\|_{p\to p}
        \leq C(p,N-1,\lambda)
        \quad
        \text{for every }\Omega'\in\Lambda(N-1,\lambda).
\]
Let $\Omega\in\Lambda(N,\lambda)$. By the affine normalization
described at the beginning of Section \ref{section: decomposition of A}, together with the endpoint
reduction following \eqref{Omegaj}, it suffices to consider the
normalized setting in which the special sequence of $\Omega$ satisfies \eqref{generic A} and \eqref{special A}. The induction hypothesis is precisely
the hypothesis \eqref{Inductive D inequality} of Proposition \ref{prop: C*}, so that proposition
applies to $\Omega$. Combining the estimates \eqref{D-Dtilde-pointwise} and \eqref{B*C*} with the proposition gives
\begin{align*}
&D_{\Omega}f(x) \leq \tilde{D}_{\Omega}f(x) \leq \mathfrak B^{\ast}_{\Omega}f(x) + \mathfrak C_{\Omega}^{\ast}f(x), \\
&\text{and therefore, }       
\|D_\Omega\|_{p\to p}
        \leq C_1(\lambda)\|\mathcal M_{\text{str}}\|_{p\to p}
              +C(p,N,\lambda).
\end{align*} 
After enlarging $C(p,N,\lambda)$ if necessary, this proves
\eqref{D-Omega-bounded} at order $N$ and completes the induction.
\vskip0.1in
\noindent Finally, suppose that
$\Omega\in\Lambda(N,\lambda;R)$. Then there exist sets
$\Omega^{[r]}\in\Lambda(N,\lambda)$, $1\leq r\leq R$, such that
\[
        \Omega\subseteq\bigcup_{r=1}^{R}\Omega^{[r]}; \; \text{ consequently,} \; 
        D_\Omega f
        \leq \max_{1\leq r\leq R}D_{\Omega^{[r]}}f
        \leq \sum_{r=1}^{R}D_{\Omega^{[r]}}f.
\]
By the triangle inequality for operator norms, $\|D_{\Omega}\|_{p \to p}$ is bounded above by the $R$-fold sum of the summands $\|D_{\Omega^{[r]}}\|_{p \to p}$. Substituting \eqref{D-Omega-bounded} into each summand then yields 
\[
        \|D_\Omega\|_{p\to p}
        \leq R\,C(p,N,\lambda),
\]
proving \eqref{just finiteness of op norm} for $D_{\Omega}$, $\Omega \in \Lambda(N, \lambda; R)$. 
\vskip0.1in
\noindent The corresponding bound for $M_\Omega$ follows from \eqref{DM-pointwise-ineq},
after another harmless enlargement of $C(p,N,\lambda)$. This proves
\eqref{just finiteness of op norm}, completing the proof of Theorem \ref{THM: DB}.
\end{proof}

% The relations \eqref{D<Dtilde} and \eqref{B*C*} give 
%\begin{align*} 
%&D_{\Omega}f(x) \leq \tilde{D}_{\Omega}f(x) \leq \mathfrak B^{\ast}_{\Omega}f(x) + \mathfrak C_{\Omega}^{\ast}f(x). \\
%&\text{Lemma \ref{lemma : strong maximal function} says that } ||\mathfrak B_{\Omega}^{\ast}||_{p \rightarrow p} \leq C_1(\lambda) ||\mathcal M_{\text{str}}||_{p \rightarrow p}, \text{ whereas } \\ 
%&\text{Proposition \ref{prop: C*} provides }  ||\mathfrak C_{\Omega}^{\ast}||_{p \rightarrow p} \leq C(p, N,  \lambda). 
%\end{align*} 
%Combining these statements completes the proof of Theorem \ref{THM: DB} for $D_{\Omega}$. The bound for $M_{\Omega}$ follows immediately from \eqref{DM-pointwise-ineq}. 
% 

	\chapter{Appendix A: Directional versus rectangular averages}
	On page \pageref{DM-pointwise-ineq}, we claimed an inequality \eqref{DM-pointwise-ineq} that connected the $L^p$ behaviour of $D_{\Omega}$ and $M_{\Omega}$. We prove this inequality here. 
	\section{Relation between $D_{\Omega}$ and $M_{\Omega}$} \label{DM-ineq-proof-section}
	\begin{lemma} 
		There is an absolute constant $C_0 > 0$ with the following property: for any $\Omega \subseteq [0,1]$, the  directional operators $D_{\Omega}$ and $M_{\Omega}$ given by \eqref{dir-max-op-def} and \eqref{Nikodym-like-op-def} satisfy  the pointwise inequality \eqref{DM-pointwise-ineq} for all smooth functions $f$ with compact support. 
		\end{lemma} 
	\begin{proof}
	Let us begin by proving the left inequality in \eqref{DM-pointwise-ineq}. The underlying intuition is that the line averages defining $D_{\Omega}$ can be thought of as limiting values of certain rectangular averages that specify $M_{\Omega}$. To make this precise, we define a unit vector $V_{\omega}$ with slope $\omega \in [0, 1]$:
		\[ V_{\omega} := \frac{(1, \omega)}{\sqrt{1 + \omega^2}}, \; \text{ and denote its perpendicular unit vector by }   \; V_{\omega}^{\perp} := \frac{(-\omega, 1)}{\sqrt{1 + \omega^2}}. \]
		Given $h > 0$ and any compactly supported smooth function $f$, the averaging operator $\mathcal A_{\omega, h}$ defined by 
		\[ \mathcal A_{\omega, h}f(x) := \frac{1}{2h} \int_{-h}^{h} \bigl| f \left(x + t V_{\omega}\right) \bigr| \, dt \] 
		 is continuous in $x$. This means that for all $x \in \mathbb R^2$, 
		 \begin{align} \mathcal A_{\omega, h}f(x) &= \lim_{r \rightarrow 0}\frac{1}{2r} \int_{-r}^{r} \mathcal A_{\omega, h} f(x + u V_{\omega}^{\perp}) \, du \nonumber \\ &= \lim_{r \rightarrow 0}\frac{1}{2r} \int_{-r}^{r} \frac{1}{2h} \int_{-h}^{h} \bigl|f \left(x + u V_{\omega}^{\perp} + t V_{\omega}\right) \bigr| \, dt du \nonumber \\
		 &= \lim_{r \rightarrow 0} \frac{1}{|R|}\int_R |f(x + y)| \, dy. 
		 \label{D<M}
		 \end{align}  
		At the last step, $R$ denotes a $(2h) \times (2r)$ rectangle centred at the origin, with length $2h$ and width $2r$, long side pointing along $V_{\omega}$, and short side along $V_{\omega}^{\perp}$. Since $h$ is fixed and $r \rightarrow 0$, one can assume $r < h$. If $\omega \in \Omega$, the  expression in \eqref{D<M} is bounded above by $M_{\Omega}f(x)$, uniformly in $h$ and $\omega$. The operator $D_{\Omega} f(x)$, on the other hand, equals the supremum of $\mathcal A_{\omega, h} f(x)$ over all $h > 0$ and $\omega \in \Omega$. Taking the supremum of both sides of \eqref{D<M} over $h >0$ and $\omega \in \Omega$ leads to the estimate $D_{\Omega} f(x) \leq M_{\Omega} f(x)$.     
		\vskip0.1in
		\noindent We now focus on proving the right side inequality in \eqref{DM-pointwise-ineq}. The idea is to compare the rectangle averages in $M_\Omega$ with the line averages in $D_\Omega$, and then bound the former by a Hardy-Littlewood 
		maximal function applied to the latter. Given any rectangle $R$ centred at the origin with length $2h$, width $2r$ and long direction oriented along $V_{\omega}$, let us write $R$ as 
		\begin{align} 
			R &= \bigl\{t V_{\omega} + u V_{\omega}^{\perp} : |t| \leq h, |u| \leq r \bigr\}, \; h > r, \; \text{ so that } \nonumber \\ 
			\frac{1}{|R|} &\int_R |f(x+y)| \, dy = \frac{1}{4hr} \int_{-h}^{h} \int_{ -r}^{r} |f(x + t V_{\omega} + u V_{\omega}^{\perp})| \, du \, dt.  \label{D>M}
		\end{align} 
		%Let us analyse the rightmost integral in two complementary ranges of $h$ and $r$. 
	%	\vskip0.1in
	  %\noindent {\em{Case 1: }} If $r \leq h \leq 2r$, then the domain of integration is contained in the square $[-h, h]^2$, whereas $hr \geq h^2/4$. In this case, the integral in \eqref{D>M} is estimated as follows,  
	  %\[ \frac{1}{|R|} \int_R |f(x+y)| \, dy \leq C_0 M_{\text{HL}} f(x) \leq C_0 M_{\text{HL}} D_{\Omega}f(x).  \]
	  %\vskip0.1in 
	  %\noindent {\em{Case 2:}} If $h > 2r$, 
	  For every $v \in [-r, r]$, a change of variables $(t, u) \mapsto (s, u)$, $s = t - v$ allows us to bound the rightmost integral in \eqref{D>M} by a family of integrals, parametrized by $v$. In fact, the relation $|t| = |s+v| \leq h$ leads to $|s| \leq h + |v| \leq h + r \leq 2h$; this implies  
	  \begin{equation}
	  \begin{aligned} 
	  	\frac{1}{|R|} \int_R |f(x+y)| \, dy &= \frac{1}{4hr} \int_{|s+v| \leq h} \int_{|u| \leq r} |f(x + sV_{\omega} + v V_{\omega} + u V_{\omega}^{\perp})| \, du \; ds \\ 
	  	&\leq \frac{1}{4hr}\int_{|s| \leq 2h} \int_{|u| \leq 2r} |f(x + s V_{\omega}+ v V_{\omega} + u V_{\omega}^{\perp})| \, du \; ds.
	  	\end{aligned} \label{pre-average}
	  	\end{equation}
	  	Averaging both sides of \eqref{pre-average} with respect to $v \in [-r, r]$, we obtain 
	  	\begin{align*}  
	  		\frac{1}{|R|} \int_R |f(x+y)| \, dy &\leq \frac{1}{8hr^2}\int_{|v| \leq r} \iint_{\begin{subarray}{c} |s| \leq 2h \\ |u| \leq 2r \end{subarray}} |f(x + s V_{\omega} + v V_{\omega} + u V_{\omega}^{\perp})| \, du \, ds \, dv \\ 
	  		&\leq \frac{1}{2r^2} \iint_{\begin{subarray}{c} |u| \leq 2r \\ |v| \leq 2r \end{subarray}} \left[ \frac{1}{4h} \int_{|s| \leq 2h} |f(x + v V_{\omega} + u V_{\omega}^{\perp} + s V_{\omega})| \, ds \right] \, du \, dv \\ 
	  		&\leq \frac{1}{2r^2} \iint_{\begin{subarray}{c} |u| \leq 2r \\ |v| \leq 2r \end{subarray}} D_{\Omega} f(x + v V_{\omega} + u V_{\omega}^{\perp}) \, du dv \\ 
	  		 &\leq 8 M_{\text{HL}} D_{\Omega}f(x). 
	  		\end{align*} 
	  		Taking the supremum of both sides with respect to $R \in \mathcal R_{\omega}$ and all $\omega \in \Omega$, we establish the right hand side inequality claimed in \eqref{DM-pointwise-ineq}, with $C_0 = 8$. This completes the proof of \eqref{DM-pointwise-ineq}, showing in particular that $D_\Omega$ and $M_\Omega$ share the same $L^p$ 
	  		boundedness for $1<p<\infty$.
		\end{proof}
		
		\chapter{Appendix B: Proofs from Section \ref{EXAMPLES SECTION}}
		\label{1d lacunary examples section}
	\section{First order lacunarity: sets versus sequences} \label{sets vs sequences proof}  
	
	\begin{proof}[Proof of Lemma \ref{Lemma: Lacunary 1}]
		Part \eqref{Lemma: Lacunary 1 (a)} is a straightforward verification of $A \in \Lambda(1, \lambda)$ using Definition \ref{defn: Lacunary sets} and taking $A$ as its own special sequence. 
		\vskip0.1in
		\noindent An arbitrary lacunary sequence $A$ obeying \eqref{lacunarity constant} is composed of two subsequences $A_{+}$ and $A_{-}$, approaching the limit $\alpha$ from the right and the left respectively. Each subsequence obeys the condition \eqref{lacunarity constant}, which ensures that $A_{\pm}$ are lacunary and monotone, with $A_{+}$ monotone decreasing and $A_{-}$ monotone increasing towards $\alpha$. This proves part \eqref{Lemma: Lacunary 1 (b)}. 
%Therefore $A_{\pm} \in \Lambda(1, \lambda)$ by part \eqref{Lemma: Lacunary 1 (a)}. 
		\vskip0.1in
		 \noindent It remains to prove part \eqref{Lemma: Lacunary 1 (c)}. It follows from Definition \ref{defn: Lacunary sets} of $\Lambda (1, \lambda)$ that $U$ admits a special sequence $A$ and a special point $\alpha$; namely, there is a monotone, lacunary sequence $A = \{a_j : j \geq 1\}$ converging to $\alpha$, such that any interval in $\mathbb R$ between two neighbouring points of $A$ contains at most one point of $U$. Let us assume, without loss of generality, that $A$ is monotone decreasing; the argument for a monotone increasing $A$ is identical. 
		 %Let us decompose $A$ into two parts: $A = A_{+} \cup A_{-}$, where $A_{+}$ (respectively $A_{-}$) denotes the subsequence of $A$ to the right (resp.ly left) of $a$. 
		 We will use $A$ to show that $U$ is itself 
		 %each one of the sets $\Omega \cap A_{\pm}$ is 
		 the union of at most two lacunary sequences $B_1, B_2$ of lacunarity constant $\lambda$, as required by the lemma.  
		 \vskip0.1in
		 \noindent 
		 %After re-indexing, let us write $A_{+} = \{ \bar{a}_j : j \geq 1 \}$, with $\bar{a}_j \searrow a$.
		  Let $J$ denote the collection of indices $j \geq 1$ such that $U \cap [{a}_{j+1}, {a}_j)$ is non-empty. It follows from the definition of $\Lambda(1, \lambda)$ that this intersection is a singleton; we denote its unique element by $b_j$. Thus $U$ consists of the elements $b_j$ for $j \in J$. The sequence $U = \{ b_j : j \in J \}$ is decreasing, because $A$ is. The condition $a_{j+1} \leq b_j < a_j$ for $j \in J$ implies that $U$ converges to $\alpha$ if $J$ is infinite. The relative positions of the elements of $U$ and $A$ imply that for $j \geq k+2$, $j, k \in J$, 
		 	\begin{equation} |b_j - \alpha| \leq |{a}_{j} - \alpha| \leq \lambda |{a}_{j-1} - \alpha| \leq \lambda |{a}_{k+1}-\alpha| \leq \lambda |b_{k} - \alpha|, \label{jump-lacunary}  \end{equation}
		  so the even and odd subsequences of $U$, namely $B_1^{\ast}= \{b_{j} : j \in J, \; j \text{ odd} \}$ and $B_2^{\ast} = \{ b_{j} : j \in J,  j \text{ even}   \}$, are both decreasing and lacunary with constant at most $\lambda$. If $B_i^{\ast}$ is an infinite sequence, for some $i=1,2$, we set $B_i := B_i^{\ast}$. If $B_i^{\ast}$ is a finite set, the property \eqref{jump-lacunary} ensures that $B_i^{\ast}$ can be extended to an infinite decreasing lacunary sequence $B_i$ with constant $\leq \lambda$. In either case, $U \subseteq B_1 \cup B_2$. This completes the proof. 
		\end{proof}
	
		\begin{proof}[Proof of Lemma \ref{Lemma: Lac set vs seq}]
			Let us observe that the interval between two consecutive elements of $A$ contains exactly one element of $U$: 
			\[ 
			U \cap [2^{-j}, 2^{-j+1}) = 
			\begin{cases}
\{2^{-j} + 4^{-j} \} &\text{ if $j$ is even,} \\   
\{2^{-(j-1)} + 4^{-(j-1)} \} &\text{ if $j$ is odd.}
\end{cases}			\]
This verifies that $U \in \Lambda(1, \frac{1}{2})$. 
\vskip0.1in
\noindent To confirm that $U$ cannot be a lacunary sequence according to Definition \ref{defn: Lacunary sequence in R}, let us note that the sequence $U$ converges to $\alpha = 0$, while the ratio of two successive elements  
\[ \frac{(2^{-2j} - 4^{-2j})}{(2^{-2j} + 4^{-2j})} \rightarrow 1 \text{ as } j \rightarrow \infty. \] 
Thus a lacunarity condition of the form \eqref{lacunarity constant} cannot hold for any $\lambda < 1$. In sum, $U$ belongs to $\Lambda(1;\tfrac12)$ but is not a lacunary sequence, 
since the ratio of successive gaps tends to $1$ rather than remaining 
uniformly bounded above by $\lambda<1$. 
This demonstrates that finite-order lacunarity genuinely extends beyond 
classical lacunary sequences.
			\end{proof} 
%		\noindent Despite this distinction, monotone lacunary sequences {\em{are}} representative of the class $\Lambda(1;\lambda)$, as seen in Lemma \ref{Lemma: Lacunary 1}. For any $(a,b) \in \mathbb R^2 \setminus \{0\}\}$, Lemma \ref{lacunarity under linear operations} gives 
%		\[ U = \{a \lambda^j + b : j \geq 1 \} \in \Lambda(1, \lambda), \text{ its special sequence being itself}.  \]
%		\vskip0.1in
		
		\section{Examples of $\Lambda(N, \lambda; R)$} \label{Lambda proofs} 
		\begin{proof}[Proof of Lemma \ref{Lemma: iterated sums 1}] 
			The proof is by induction on $N$. For $N = 1$, the set $\overline{U}_1(M_1)$ is the lacunary sequence $\{M_1^{-j}\}$, which lies in $\Lambda(1, M_1^{-1})$ by Lemma \ref{Lemma: Lacunary 1}. Proceeding to the induction step, the monotonicity of the indices $M_r$ and $j_r$ imply 
			\begin{align*} &\sum_{r=1}^{N} M_r^{-j_r} \leq N M_1^{-j_1}, \text{ which is } < M_1^{-j_1+1} \text{ since } M_1 > N. \\ &\text{As a result, } \sum_{r=1}^{N} M_r^{-j_r} \in \left[ M_1^{-j_1}, M_1^{-j_1+1} \right) \text{ whenever } j_1 \leq j_2 \leq \ldots \leq j_N.   \end{align*} 
			Thus, each element of $\overline{U}_N(\mathbf M)$ lies within an interval of the form $[M^{-j_1}, M^{-j_1+1}]$; in fact, the portion of $\overline{U}_N(\mathbf M)$ in such an interval is a translate of $\overline{U}_{N-1}(\mathbf M')$, where $\mathbf M' = (M_2, \ldots, M_N)$.  
			The induction hypothesis, combined with Lemma \ref{lacunarity under linear operations} \eqref{linear-invariance}, allows us to deduce that for every $j_1 \geq 1$, 
			 \begin{align*}  \overline{U}_N(\mathbf M) \cap \left[M_1^{-j_1}, M_1^{-j_1+1} \right] &\subseteq M_1^{-j_1} + \left\{ \sum_{r=2}^{N} M_r^{-j_r} : 1 \leq j_2 \leq \cdots \leq j_N \right\} \\ &\subseteq M_1^{-j_1} + \overline{U}_{N-1}(\mathbf M') \in \Lambda(N-1, M_2^{-1}). 
			 	%\subseteq \Lambda(N-1, M_1^{-1}).
			 	\end{align*}
To clarify the last statement, the induction hypothesis implies $\overline{U}_{N-1}(\mathbf M') \in \Lambda(N, M_2^{-1})$; therefore, the same inclusion is true of its translate, by Lemma \ref{lacunarity under linear operations} \eqref{linear-invariance}. 
			 	The assumption $M_2^{-1} \leq M_1^{-1}$ implies that $ \Lambda(N-1, M_2^{-1}) \subseteq  \Lambda(N-1, M_1^{-1})$, by the monotonicity of the classes $\Lambda(\cdot, \lambda)$ in $\lambda$, as shown in Lemma \ref{Lemma : lacunarity monotonicity}. Hence $\overline{U}_N(\mathbf M) \in \Lambda(N, M_1^{-1})$, completing the induction. This establishes \eqref{U-order-N},  concluding the proof.  
			\end{proof} 
			
				\begin{proof}[Proof of Lemma \ref{Lemma: iterated sums 2}]
					As in Lemma \ref{Lemma: iterated sums 1}, the proof is by induction on $N$. 
\vskip0.1in
\noindent {\bf{Base step. }} For $N = 1$, the set $U^{\ast}_1(M_1) = \{M_1^{-j} : j \geq 1\}$ is a monotone, lacunary sequence, and therefore in $\Lambda(1, M_1^{-1})$ by Lemma  \ref{Lemma: Lacunary 1} \eqref{Lemma: Lacunary 1 (a)}. This leads to the desired conclusion with $R_1 =1$ and $\lambda_1 = M_1^{-1}$. 
\vskip0.1in
\noindent {\bf{Inductive step: Overview. }}Proceeding to the inductive step, the induction hypothesis posits that 
\[
\left\{ \begin{aligned} 
&{\text{ for any $\ell < N$, a set of the form $U_{\ell}^{\ast}(\mathbf M')$ is contained in the union of $R_{\ell}$}} \\ &{\text{sets in the class $\Lambda(\ell, \lambda_{\ell})$, where $R_{\ell}, \lambda_{\ell}$ depend only on $\ell$ and $\mathbf M' \in (1, \infty)^{\ell}$. }}
\end{aligned} \right\} \]  
					%for any $I \subset \{1, 2, \ldots, N\}$, $\#(I) < N$, there exist at most $R_{N-1}$ sets $V_i(I) \in \Lambda(N-1, M_1^{-1})$
					%\begin{equation}
					%	U(I) := \Bigl\{ \sum_{r \in I} M_r^{-j_r}: \, j_r \in \{1, 2, 3, \ldots\} \text{ for all } r \in I\Bigr\} \subseteq \bigcup_{i=1}^{R_{N-1}} V_i(I).  
					%	\end{equation}  
						To prove the statement for $\ell = N$, we decompose the set $U^{\ast}_N(\mathbf M)$ in \eqref{Example: U Lac N}  as follows. Each element of $U^{\ast}_N(\mathbf M)$ is a sum of $N$ positive summands, hence the size of the element is comparable to the largest summand. The index of the largest summand allows us to categorize $U_N^{\ast}(\mathbf M)$ into $N$ subsets, where for $1 \leq \sigma \leq N$, the $\sigma^{\text{th}}$ set contains the members of $U$ whose $\sigma^{\text{th}}$ summand is representative of its size. 
						%This allows us to classify the set of multi-indices $\mathbf j = (j_1, \ldots, j_N) \in \mathbb N^{N}$ indexing the elements of $U_N^{\ast}$, using the subscript $r$ of the representative  
						%Precisely, for $\sigma \in \mathbb S_N$, i.e. a permutation $\sigma$ of $\{1, 2, \ldots, N\}$, we define  
						%\begin{equation}
						%	\label{perm-sigma} \mathbb J_{\sigma} := \left\{ \mathbf j = (j_1, \ldots, j_N) \in \mathbb N^N: M_{\sigma(1)}^{-j_{\sigma(1)}}\leq M_{\sigma(2)}^{-j_{\sigma(2)}}
						%		\leq \cdots \leq M_{\sigma(N)}^{-j_{\sigma(N)}} \right\}. \end{equation} 
							Precisely, let
	\begin{align*} &\mathbb J_{\sigma} := \left\{ \mathbf j = (j_1, \ldots, j_N) \in \mathbb N^{N} : M_{\sigma}^{-j_{\sigma}} = \max_{1 \leq \ell \leq N} M_{\ell}^{-j_{\ell}}\right\}, \; 1 \leq \sigma \leq N, \text{ so that } \\  
	&U^{\ast}_N(\mathbf M) = \bigcup_{\sigma=1}^{N} U_{\sigma}(\mathbf M), \quad \text{ where } \quad U_{\sigma} = U_{\sigma}(\mathbf M) := \left\{ \sum_{r=1}^{N} M_r^{-j_r} \, : \, \mathbf j \in \mathbb J_{\sigma}\right\}. 
	\end{align*}
The index sets $\mathbb J_{\sigma}$ may overlap when the maximum is attained for multiple indices; this is harmless since we only seek an upper bound on the number of covering sets.  
	Let us pick $C_N \in \mathbb N$ and $\lambda_N \in (0,1)$, depending only on $N$ and $\mathbf M$ such that
	\begin{equation}  
\left\{
\begin{aligned} 		
&N < \min_{1 \leq \sigma \leq N}M_{\sigma}^{C_N}, \; \text{ e.g. } \; C_N := 2 \max_{\sigma} \left \lceil \frac{\log N}{\log M_{\sigma}}  \right \rceil, \\
		&\lambda_N := \lambda_N(\mathbf M) :=  \max_{1 \leq \sigma \leq N} \left\{\lambda_{N-1}(\mathbf M_{\sigma}'), M_{\sigma}^{-1} \right\} < 1.
\end{aligned} \right\} \label{lambdaN} \end{equation} 
		Here $\mathbf M_{\sigma}'$ denotes the vector $\mathbf M$, with its $\sigma^{\text{th}}$ entry $M_{\sigma}$ omitted.
	In the remainder of the proof, we aim to show that 
\begin{equation} \label{U-sigma-covering-claim}
\text{$U_{\sigma}$ is covered by at most $R_{N-1} C_N$ sets of $\Lambda(N, \lambda_N)$, }
\end{equation} 
with $C_N, \lambda_N$ as in \eqref{lambdaN}.
	%To complete the induction, we need certain structural properties of $U_{\sigma}$. As we show below, a part of $U_{\sigma}$, termed $U_{\sigma 2}$, is a controlled, finite union of lacunary sets of order at most $(N-1)$, using the induction hypothesis. The remainder, called $U_{\sigma 1}$, can be realized as a finite union of lacunary sets of order $N$, possibly with a larger lacunarity constant, also appealing to the induction hypothesis. In quantitative terms, we will prove the existence of positive constants $C_N \geq 1$ and $\lambda_N \in (0,1)$, depending only on $N, \mathbf M$, such that each set $U_{\sigma }$ is contained in the union of at most $2 R_{N-1}  C_N$ sets of type $\Lambda(N, \lambda_N)$. 
	This will complete the induction with $R_N = N R_{N-1} C_N$.
	\vskip0.1in
	\noindent To prove \eqref{U-sigma-covering-claim}, fix $\sigma \in \{1, 2, \ldots, N\}$. Every $\mathbf j \in \mathbb J_{\sigma}$ obeys the property 
	\begin{equation}
	 M_{\sigma}^{-j_{\sigma}} \leq \sum_{r=1}^{N} M_r^{-j_r} \leq N M_{\sigma}^{-j_{\sigma}} \quad \text{ i.e., } \quad   \sum_{r=1}^{N} M_r^{-j_r} \in \Bigl[  M_{\sigma}^{-j_{\sigma}}, N M_{\sigma}^{-j_{\sigma}}  \Bigr]. \label{U-sigma-property}
	\end{equation}
					%The set $U^{\ast}_N$ is invariant under permutations of $(M_1, \ldots, M_N)$; therefore without loss of generality, by re-indexing the constants $M_r$ if necessary, we may assume that $\sigma$ in \eqref{perm-sigma} is the identity permutation, i.e.,  
					%\begin{align} &M_1^{-j_1} \leq M_2^{-j_2} \leq \cdots \leq M_N^{-j_N} \; \text{ for all } \mathbf j \in \mathbb J_{\sigma}; \nonumber 
					%	\\ &\text{ this implies }  
					 %M_N^{-j_N} \leq \sum_{r=1}^{N} M_r^{-j_r} \leq N M_N^{-j_N}. \label{U-sigma-property}
					 %\end{align}
To paraphrase, a sum $\sum_{r} M_r^{-j_r} \in U_{\sigma}$ lies in an interval whose length is at most $(N-1)$ times its left endpoint, namely $M_{\sigma}^{-j_{\sigma}}$, which is also the largest summand.  
Our goal is to identify the interval on the right of \eqref{U-sigma-property} as the space between two (or at most a controlled number of) consecutive elements of a monotone lacunary sequence, then use the induction hypothesis to deduce that the portion of $U_{\sigma}$ in this interval is admissible lacunary of lower order.  Through the iterative definition of $\Lambda(N, \cdot)$, this will lead us to \eqref{U-sigma-covering-claim}.  We focus now on the details of the proof. 
\vskip0.1in 
\noindent {\bf{Completion of the inductive step: Proof of \eqref{U-sigma-covering-claim}}} Two cases arise, depending on the size of $M_{\sigma}$.
\vskip0.1in  
\noindent {\em{Case 1:}} Suppose $M_{\sigma} > N$. Then  it follows from \eqref{U-sigma-property} that every multi-index $\mathbf j \in \mathbb J_{\sigma}$ and its corresponding element in $U_{\sigma}$ obeys  the property 
\begin{align} 
&\sum_{r=1}^{N} M_r^{-j_r} \in \Bigl[  M_{\sigma}^{-j_{\sigma}},  M_{\sigma}^{-j_{\sigma}+1} \Bigr) =: \mathcal A_{j_{\sigma}};   \;  \text{ this in turn means that }  \nonumber \\ 
& U_{\sigma} \cap  \mathcal A_{j_{\sigma}} \subseteq \mathcal A_{j_{\sigma}} \cap \left[M_{\sigma}^{-j_{\sigma}} +   U^{\ast}_{N-1} (\mathbf M_{\sigma}') \right] \;\text{ for every fixed $j_{\sigma} \geq 1$. } \label{U-sigma-portion} 
\end{align}  
% We claim that $U_{\sigma} \in \Lambda(N, \lambda_{N})$ with special sequence $A = \{ M_{\sigma}^{-j} : j \geq 1\}$. The lacunarity constant of $A$ is $M_{\sigma}^{-1}$, which is $\leq \lambda_N$.  
To prove the claim \eqref{U-sigma-covering-claim}, we note by the induction hypothesis that the set $U^{\ast}_{N-1} (\mathbf M_{\sigma}')$ is contained in the union of $R_{N-1}$ sets of type $\Lambda(N-1, \lambda_{N-1}(\mathbf M_{\sigma}'))$. By Lemma \ref{lacunarity under linear operations}, the same conclusion is true for all subsets and translates of $U^{\ast}_{N-1} (\mathbf M_{\sigma}')$, in particular for the rightmost set in \eqref{U-sigma-portion}. In other words, for each $j \geq 1$, there is a collection of sets 
 \begin{align}
&\left\{
\begin{aligned}
 &\{ V_{i j} :  1 \leq i \leq R_{N-1}\}, \text{ with } V_{ij} \subseteq  \mathcal A_j, \; V_{ij} \in \Lambda(N-1, \lambda_{N-1}(\mathbf M'_{\sigma})) \\ 
 &U_{\sigma} \cap  \mathcal A_j \subseteq \mathcal A_j \cap \bigl[ M_{\sigma}^{-j} + U_{N-1}^{\ast}(\mathbf M_{\sigma}') \bigr] \subseteq \bigcup_{i=1}^{R_{N-1}} V_{ij} \subseteq \mathcal A_j 
\end{aligned} 
\right\}.  
\label{Vij-def}\\ 
&\text{ This means } 
U_{\sigma} \subseteq  \bigcup_{j\geq 1} \left\{ U_{\sigma} \cap  \mathcal A_j \right\}  = \bigcup_{i=1}^{R_{N-1}} V_{i},  \, \text{ where } \, V_{i} := 
\bigcup_{j \geq 1} V_{ij}. \label{Vi-def}
\end{align}
To conclude the proof of \eqref{U-sigma-covering-claim}, it remains to show that for each $1 \leq i \leq R_{N-1}$, the set $V_i \in \Lambda(N, \lambda_N)$, with special sequence $\mathcal A = \{ M_{\sigma}^{-j} : j \geq 0 \}$. Let us note that $\mathcal A$ has lacunarity constant $M_{\sigma}^{-1}$, which is $\leq \lambda_N$, by virtue of \eqref{lambdaN}. According to \eqref{Vij-def}, the portion of $V_i$ that lies in the interval $\mathcal A_j$ between two consecutive elements of the special sequence $\mathcal A$ is $V_{ij}$. By the monotonicity \eqref{lacunarity-monotonicity} of the classes $\Lambda(\cdot, \lambda)$ in $\lambda$, as proved in Lemma \ref{Lemma : lacunarity monotonicity},  the set $V_{ij}$ is of type $\Lambda(N-1, \lambda_N)$, since $\lambda_{N-1}(\mathbf M_{\sigma}') \leq \lambda_N$ by \eqref{lambdaN}. Thus, in this case, $U_{\sigma}$ is contained in the union of at most $R_{N-1}$ sets in $\Lambda(N, \lambda_N)$. This proves the statement \eqref{U-sigma-covering-claim} in Case 1, in fact with a smaller cover than claimed ($R_{N-1}$ instead of $R_{N-1} C_N$).  
\vskip0.1in 
\noindent {\em{Case 2:}} Suppose now that $M_{\sigma} \leq N$. The main distinction in this case is that  the endpoints of the interval $[M_{\sigma}^{-j_{\sigma}}, NM_{\sigma}^{-j_{\sigma}}]$ occuring in \eqref{U-sigma-property} cannot be consecutive elements of a single lacunary sequence. Nonetheless, the interval contains a controlled finite number of elements of the lacunary sequence $\mathcal A = \{ M_{\sigma}^{-j} : j \geq - C_N \}$, as we show below. Let us recall the definition of $C_N$ from \eqref{lambdaN}. The relation \eqref{U-sigma-property} then leads to
\[ \sum_{r=1}^{N} M_r^{-j_r} \in \left[ M_{\sigma}^{-j_{\sigma}}, M_{\sigma}^{-j_{\sigma} + C_{N}}\right)=  \bigcup_{\ell=0}^{C_N-1} \mathcal A_{j_{\sigma} - \ell}. \]
In other words, the sum $\sum_{r} M_r^{-j_r} \in U_{\sigma}$ can lie in at most $C_N$ adjacent intervals $\mathcal A_j$ formed by consecutive elements of the lacunary sequence $\mathcal A$. Here $j_{\sigma}$ indexes the largest summand, whereas $j$ indexes (gaps within) the lacunary sequence $\mathcal A$; the indices $j$ and $j_{\sigma}$ are related by $0 \leq j_{\sigma} - j \leq C_N-1$. At the same time, $\sum_r M_r^{-j_r} \in M_{\sigma}^{-j_{\sigma}} + U^{\ast}_{N-1}(\mathbf M_{\sigma}')$. Combining these two statements, we obtain 
\begin{align} 
 &U_{\sigma} \cap  \mathcal A_{j}  \subseteq \bigcup_{j_{\sigma} = j}^{j + C_N-1}  W(j, j_{\sigma}), \text{ for all $j > -C_N$, where }  \label{U-sigma-portion-case-2} \\ &W(j, j_{\sigma} ) := \Bigl[ M_{\sigma}^{-j_\sigma} + U^{\ast}_{N-1}(\mathbf M_{\sigma}') \Bigr] \cap  \mathcal A_{j}. \nonumber 
\end{align}
The proof now proceeds as in Case 1. For each $j \geq 1$ and $j \leq j_{\sigma} \leq  j+C_N-1$, the induction hypothesis yields subsets $\bigl\{\widetilde{W}_i(j, j_{\sigma}) : 1 \leq i \leq R_{N-1}  \bigr\}$ of $\mathcal A_j$, such that
\begin{equation}  W(j, j_{\sigma}) \subseteq \bigcup_{i=1}^{R_{N-1}} \widetilde{W}_i(j, j_{\sigma}), \quad \widetilde{W}_i(j, j_{\sigma}) \in \Lambda(N-1, \lambda_{N-1}(\mathbf M_{\sigma}')). \label{W-def}
\end{equation} 
Inserting \eqref{W-def} into \eqref{U-sigma-portion-case-2}, we arrive at 
\begin{align*} 
U_{\sigma} \subseteq \bigcup_{j=-C_N+1}^{\infty} \bigl[ U_{\sigma} \cap \mathcal A_j \bigr] &\subseteq \bigcup_{i=1}^{R_{N-1}} \left[ \bigcup_{j = -C_N+1}^{\infty} \bigcup_{\ell=0}^{C_{N}-1} \widetilde{W}_i(j, j + \ell) \right] \\ &=  \bigcup_{i=1}^{R_{N-1}}  \bigcup_{\ell=0}^{C_N-1} W_{i \ell}, \text{ with } W_{i \ell} := \bigcup_{j= - C_N + 1}^{\infty}  \widetilde{W}_i(j, j + \ell).  
\end{align*}  
By the monotonicity of $\Lambda(\cdot, \lambda)$ as in Case 1, we deduce that $W_{i \ell} \in \Lambda(N, \lambda_N)$ for each $1 \leq i \leq R_{N-1}$, $0 \leq \ell \leq C_N-1$, all of them with the special sequence $\mathcal A$. Thus $U_{\sigma}$ is contained in the $R_{N-1}C_N$-fold union of sets of type $\Lambda(N, \lambda_N)$, proving \eqref{U-sigma-covering-claim} and completing the induction.
					\end{proof}
					\subsection{Lacunarity order in arithmetic progressions} \label{proof: dyadic rationals lemma}  
Lemma \ref{lemma: dyadic rationals} has two parts. Part \eqref{lemma: dyadic rationals parta} asserts a uniform lacunarity order for all non-empty finite sets of bounded cardinality.  Part \eqref{lemma: dyadic rationals partb} gives an example of an increasing sequence of equally spaced finite sets $\mathbb Q_m$ that demonstrates that, in general, the lacunarity order grows with the cardinality. Each set $\mathbb Q_m$ is finite, and hence in {\tt{AdFinLac}}, but the smallest lacunarity order $N$ such that $\mathbb Q_m \in \Lambda(N; \frac{1}{2})$ grows without bound. We prove this lemma in this section. 
\vskip0.1in
\subsubsection{Proof of Lemma \ref{lemma: dyadic rationals} \eqref{lemma: dyadic rationals parta}} 
\begin{proof} 
Using strong induction on $n$, we will show that 
\begin{equation} 
  {\text{any finite set }} \mathbb F_n \subseteq\mathbb R  \text{ with }                   \#(\mathbb F_n )=n \text{ belongs to } \Lambda\bigl(n-1,\frac12 \bigr).
\end{equation} 
The assertion is immediate when $n=1$, since a singleton belongs to $\Lambda(0,\frac12)$. Suppose now that $n\geq 2$, and write
\[
                    \mathbb F_n=\{x_1<x_2<\cdots<x_n\}.
\]
Set \(\alpha=x_1\), and consider the monotone decreasing sequence
\[
 a_j:=\alpha+2^{-(j-1)}(x_n-\alpha),
 \qquad j\geq 1.
\]
Then $a_1=x_n$, $a_j\to\alpha$, and
\[
                    a_{j+1}-\alpha
                       =\frac12(a_j-\alpha).
\]
Thus $A=\{a_j:j\geq 1\} \in {\tt{MonLac}}(\frac12)$, and
\[
                    \mathbb F_n \subseteq[\alpha,a_1].
\]
Neither $x_1=\alpha$ nor $x_n=a_1$ belongs to any gap interval
$\bigl[a_{j+1},a_j\bigr)$. Consequently,
\[
 \#\bigl(\mathbb F_n \cap \bigl[a_{j+1},a_j \bigr)\bigr)\leq n-2
 \quad\text{for every }j\geq 1.
\]
By the strong induction hypothesis and the nesting property of the lacunary classes,
\[
 \mathbb F_n \cap[a_{j+1},a_j)\in\Lambda(n-2, \frac12)
 \quad\text{for every }j\geq 1.
\]
Here the empty intersection belongs to $\Lambda(0, \frac12)$.
It follows from Definition \ref{defn: Lacunary sets} that $\mathbb F_n\in\Lambda(n-1,\frac12)$ completing the induction.
\end{proof} 
\subsubsection{Proof of Lemma \ref{lemma: dyadic rationals} \eqref{lemma: dyadic rationals partb}} 
\begin{proof} 
Applying Lemma \ref{lemma: dyadic rationals} \eqref{lemma: dyadic rationals parta} to $\mathbb Q_m$, which has cardinality $2^m$, gives
\[
                    \mathbb Q_m\in\Lambda \Bigl(2^m-1, \frac12 \Bigr).
\]
Therefore the collection appearing in the definition of $q_m$ is
nonempty, and $q_m$ is well-defined.
\vskip0.1in 
\noindent  Lemma \ref{lemma: dyadic rationals} \eqref{lemma: dyadic rationals partb} is a consequence of the following stronger statement concerning arbitrary arithmetic progressions. 
\begin{lemma} \label{lemma: dyadic rationals refined}
For $m \geq 1$, every arithmetic progression $\mathscr{P}_m$ of length 
\begin{equation} \label{not in Lac induction} 
\# \mathscr{P}_m = 5\cdot 3^{m-1}  \; \text{ satisfies } \; \mathscr{P}_m \notin  \Lambda \Bigl(m,\frac12 \Bigr). 
\end{equation}  
			\end{lemma}
\noindent Let us first see why Lemma \ref{lemma: dyadic rationals refined} implies Lemma \ref{lemma: dyadic rationals}.
Fix an integer $k\geq 1$. If
\[
                    2^m \geq 5\cdot 3^{k-1},
\]
then $\mathbb Q_m$ contains an arithmetic progression of length
$5\cdot 3^{k-1}$. By Lemma \ref{lemma: dyadic rationals refined}, this sub-progression does not belong to $\Lambda(k,\frac12)$. If $\mathbb Q_m$ belonged to $\Lambda(k,\frac12)$, then
the subset property in Lemma \ref{lacunarity under linear operations} would imply that the sub-progression also belonged to $\Lambda(k,\frac12)$, a contradiction. This proves that 
%\[
%                    Q_m\notin\Lambda(k,1/2).
%\]
%
%By Lemma \ref{lemma: dyadic rationals refined}, this sub-progression does not
%belong to $\Lambda(k,\frac12)$. Since $\Lambda(k,\frac12)$ is preserved
%under passage to subsets, it follows that
\[
                    \mathbb Q_m\notin\Lambda \bigl(k,\frac12 \bigr).
\]
Thus $q_m>k$ whenever $2^m\geq 5\cdot 3^{k-1}$. Since $k$ is
arbitrary, we conclude that $q_m\to\infty$ as $m\to\infty$.                  
\end{proof} 
\subsubsection{Proof of Lemma \ref{lemma: dyadic rationals refined}} 
		\begin{proof}
		 In view of Lemma \ref{lacunarity under linear operations}, and after an affine transformation if necessary, we may assume 
		\begin{equation} \label{AP}
		 \mathscr{P}_m = \{n \theta : 1 \leq n \leq d_m \} \text{ for some } \theta > 0, \; d_m := \# \bigl(\mathscr{P}_m \bigr) = 5 \times 3^{m-1}. \end{equation} 
		 We record two harmless normalizations that will be used below. First, after applying the reflection
\[
                    x\longmapsto (d_m+1)\theta-x
\]
if necessary, we may assume that any special sequence under
consideration is monotone decreasing. Notice that this reflection
preserves the progression $\mathscr{P}_m$.
\vskip0.1in
\noindent Second, after discarding finitely many initial terms of the special sequence and, if necessary, adjoining one additional initial term, we may arrange that the largest element $d_m\theta$ of $\mathscr{P}_m$ lies in its first gap interval. Indeed, adjoining an initial term $a_0$ satisfying
\[
                    a_0-\alpha=2(a_1-\alpha)
\]
preserves the lacunarity constant $\frac12$. Thus, in the notation used
below, we may assume that
\begin{equation} \label{first gap}
                    d_m\theta\in[a_2,a_1).
\end{equation} 
\vskip0.1in		 
\noindent		With these normalizations in place, the proof of Lemma \ref{lemma: dyadic rationals refined} is by induction on $m$.  
		%without loss of generality, we assume
%        \[ Q_1=\{ \alpha,2\alpha,3\alpha,4\alpha,5\alpha \} \]
%        for some $\alpha>0$. 
        Let us start with the base case $m=1$. Aiming for a contradiction, suppose if possible that 
\begin{equation*} 
\left\{
\begin{aligned}
&\mathscr{P}_1 = \bigl\{ \theta, 2 \theta, \ldots, 5 \theta \bigr\} \in \Lambda \bigl(1,\frac12 \bigr), \text{with a decreasing special}  \\
&\text{sequence } A = \{a_n : n \geq 1 \} \in {\tt{MonLac}}\bigl(\frac{1}{2} \bigr) \text{ obeying \eqref{first gap}}. 
\end{aligned} 
\right\}
\end{equation*}  
%Without loss of generality, let us assume that $A$ is monotone decreasing, and converges to $\alpha$. 
The definition of $\Lambda(1, \frac{1}{2})$ then implies that 
\begin{equation}  \label{lac1 condition} 
\mathscr{P}_1 \subset [\alpha, a_1], \quad \# \Bigl[ \mathscr{P}_1 \cap [a_{k+1}, a_k) \Bigr] \leq 1. 
\end{equation}  
%Then, there exists a monotonically decreasing lacunary sequence $\{ a_n \}_{n\geq 1}$ converging to $a$ such that $Q_1\subset [a,a_1]$, and each interval $[a_{k+1},a_k)$ contains at most one element of $Q_1$ for all $k$. 
The condition \eqref{first gap} with $d_1 = 5$ gives  
\begin{equation}  \label{5 & 4}
5\theta\in [a_2,a_1), \quad \text{ which then implies } \quad 4\theta \notin [a_2,a_1) \text{ by \eqref{lac1 condition}}. 
\end{equation}  
On the other hand, the lacunarity of the sequence $A$ implies \[ a_3- \alpha \leq \frac{1}{2}(a_2-\alpha ), \; \text{ which leads to } \; a_3\leq \frac{1}{2}(a_2+\alpha). \] Since $\alpha \leq \theta$ and $a_2 \leq 5 \theta$, we deduce that 
\[ a_3 \leq \frac{5\theta}{2}+\frac{\theta}{2}=3\theta. \] Combined with \eqref{5 & 4}, this means that 
\[ \{ 3\theta,4\theta \}\subset [a_3,a_2) \cap \mathscr{P}_1, \quad \text{ contradicting \eqref{lac1 condition}}. \]
This completes the proof of the base case. 
\vskip 0.1in 
\noindent Proceeding to the inductive step, suppose that the statement \eqref{not in Lac induction}, i.e., the conclusion of Lemma \ref{lemma: dyadic rationals refined}, holds
for all $m\leq M-1$, for some $M\geq 2$. Set
\[
                    L:=5\cdot 3^{M-2}, \; \text{ so that } \;
                    d_M=5\cdot 3^{M-1}=3L. \]
Aiming for a contradiction, suppose that $\mathscr{P}_M\in\Lambda (M,\frac12)$, with a special sequence $B=\{b_n:n\geq 1\}\in {\tt{MonLac}}(\frac12)$.
By the normalizations above, we may assume that \[ b_n \searrow \beta \; \text{ and } \; d_M\theta\in[b_2,b_1). \]
The defining properties of $\Lambda(M,\frac12)$ then give $ \mathscr{P}_M\subseteq[\beta,b_1]$, so that $\beta\leq\theta$; also, 
\begin{equation} \label{AP in a gap} 
 [b_{k+1},b_k)\cap \mathscr{P}_M\in\Lambda\Bigl(M-1,\frac12\Bigr)
 \quad\text{for every }k\geq 1.
\end{equation} 
\vskip0.1in
\noindent We now claim that
\begin{equation} \label{subap}
 \#\bigl([b_{k+1},b_k)\cap \mathscr{P}_M\bigr)<L
 \; \text{for every }k\geq 1.                 
\end{equation} 
Indeed, the intersection of any interval with $\mathscr{P}_M$ is a block of
consecutive terms of the progression, which constitutes a progression in its own right. If the intersection above had
cardinality at least $L$, it would therefore contain an arithmetic progression
of length $L=5\cdot 3^{M-2}$. By the induction hypothesis, that
sub-progression would not belong to $\Lambda(M-1,\frac12)$. This would contradict both \eqref{AP in a gap}
%\[
% \bigl[b_{k+1},b_k \bigr)\cap \mathscr{P}_M\in\Lambda\Bigl(M-1,\frac12 \Bigr)
%\]
and the fact that $\Lambda(M-1,\frac12)$ is preserved under passage to subsets. This proves \eqref{subap}.
\vskip0.1in
\noindent According to \eqref{first gap}, $d_M\theta\in[b_2,b_1)$, therefore there is an integer $r\geq 1$ such that
\[
 [b_2,b_1)\cap \mathscr{P}_M
   =\{r\theta,(r+1)\theta,\ldots,d_M\theta\}.
\]
It follows from \eqref{subap} that
\begin{equation} \label{what is r}
                    r\geq 2L+2
                      =10\cdot 3^{M-2}+2.
\end{equation} 
Indeed, if $r<2L+2$, then
\[
 \#\bigl([b_2,b_1)\cap \mathscr{P}_M\bigr)
   =d_M-r+1
   =3L-r+1
   \geq L,
\]
contradicting \eqref{subap}. We consequently have $b_2\leq r\theta$, $(r-1)\theta<b_2$. The lacunarity of $B$ now gives
\begin{align}
& b_3-\beta\leq\frac12(b_2-\beta), \text{ and hence } \nonumber \\ 
& b_3
 \leq \frac12(b_2+\beta)
 \leq \frac12(r+1)\theta
 <(r-1)\theta
 <b_2. \label{between b2 and b3}
\end{align} 
Here we used $\beta\leq\theta$ and $r>3$. The chain of inequalities \eqref{between b2 and b3} says the gap
interval $[b_3,b_2)$ contains every integer multiple $n\theta$
such that
\[
                    \frac{r+1}{2}\leq n\leq r-1.
\]
It follows from the above inequality and \eqref{what is r} that
\[
 \#\bigl([b_3,b_2)\cap \mathscr{P}_M\bigr)
   \geq \frac{r-2}{2}
   \geq L
   =5\cdot 3^{M-2},
\]
again contradicting \eqref{subap}. This completes the induction and hence
the proof of the lemma.
		\end{proof}

		\section{Examples of sublacunary sets} \label{section: sublacunary set proofs} 
	%	\begin{lemma} 
	%		The set \[ U = \left\{ \frac{k}{2^m} : 0 \leq k \leq 2^m, \, k, m \in \mathbb N \right\}\] of dyadic rationals in $[0,1]$ is sublacunary. 
	%	\label{dyadic rationals example} 
	%	\end{lemma} 
	%	\begin{proof}
	%		{\color{red}{\bf{Question:}} Is there a direct proof that does not use splitting numbers?}
	%		\end{proof}
Lemma \ref{lemma: topological closure} asserts a quantitative version of the property that admissible finite order lacunarity is preserved under topological closure. We prove this lemma in this section. 		
\begin{proof}[Proof of Lemma \ref{lemma: topological closure}]
%By Lemma \ref{lacunarity under linear operations}, the property of finite-order lacunarity is preserved under set containment. Therefore, $\bar{U} \in $ {\tt{AdFinLac}} implies that $U \in$ {\tt{AdFinLac}}. For the converse, we will show by induction that 
%and therefore in {\tt{AdFinLac}}. 
The main statement to be proved in \eqref{induction: topological closure}; all other claims in the lemma follow from it. 
\noindent We prove \eqref{induction: topological closure} by induction on $N$. The base case of  for $N = 0$ is clear, since $U = \bar{U}$ for a set $U$ that is empty or a singleton. 
\vskip0.1in
\noindent For $N=1$, Lemma \ref{Lemma: Lacunary 1} says that $U \in \Lambda(1, \lambda)$ is contained in the union of two monotone lacunary sequences $B_1$ and $B_2$, with the same limit, say $\alpha$. Each $B_i$ is a member of $\Lambda(1, \lambda)$, by the same lemma. Since 
\[ \bar{B}_i = B_i \cup \{\alpha\}, \] we deduce that for any $U \in \Lambda(1, \lambda)$, its closure $\bar{U}$ is contained in the union of the two closed sets $\bar{B}_1$ and $\bar{B}_2$. 
%and the singleton $\{\alpha\}$. A singleton is in $\Lambda(0, \lambda)$, and therefore in $\Lambda(1, \lambda)$, by the monotonicity property Lemma \ref{Lemma : lacunarity monotonicity}.  
In other words, \eqref{induction: topological closure} holds for $N=1$ with $C_1=2$.  
\vskip0.1in
\noindent Let us move on to the induction step. Suppose that $N \geq 2$ and  $U$ is in $\Lambda(N, \lambda)$, with a special sequence $A = \{a_j: j \geq 1\}$. Assuming without loss of generality that $U \subseteq [0,1]$ and $A$ is monotone decreasing to a limit $\alpha$, we know from Definition \ref{defn: Lacunary sets} that 
\begin{equation*} 
U_j := U \cap [a_{j+1}, a_j) \in \Lambda(N-1, \lambda) \text{ for all $j \geq 1$}. 
\end{equation*}  
Then on one hand, 
\begin{equation} \label{barUj} 
\bar{U}_j \subseteq \bar{U} \cap [a_{j+1}, a_j], \text{ so that } \bar{U} \subseteq \bar{A} \cup \bigcup_{j=1}^{\infty} \bar{U}_j, \; \text{ where } \bar{A} = A \cup \{\alpha\}.  
\end{equation} 
On the other hand, for every $j \geq 1$ and by the induction hypothesis, there exist  closed sets $\{{V}_{ij} : 1 \leq i \leq C_{N-1}\}$ such that 
\begin{equation} \label{Vij lower order} 
\bar{U}_j \subseteq \bigcup_{i=1}^{C_{N-1}} {V}_{ij}, \;
\text{ where ${V}_{ij} \in \Lambda(N-1, \lambda)$ for all $i, j$;}
\end{equation} 
Without loss of generality, after intersecting with $[a_{j+1}, a_j]$ if necessary, we may assume that $V_{ij} \subseteq [a_{j+1}, a_j]$. 
Combining \eqref{Vij lower order} with \eqref{barUj} yields the cover 
\begin{align} \label{U closure cover} 
&\bar{U} \subseteq \bar{A} \cup \bigcup_{j=1}^{\infty} \bar{U}_j  \subseteq \bar{A} \cup \bigcup_{j=1}^{\infty}  \bigcup_{i=1}^{C_{N-1}}  V_{ij} \subseteq \bar{A} \cup \bigcup_{i=1}^{C_{N-1}} \Bigl[V_i^{\text{\tiny{odd}}} \sqcup V_i^{\text{\tiny{even}}} \Bigr], \text{where } \\
&V_i^{\text{\tiny{odd}}} := \{\alpha\} \cup \bigcup_{j=1}^{\infty} \bigl\{ V_{ij} : j \text{ odd} \bigr\}, \quad  V_i^{\text{\tiny{even}}} := \{\alpha\} \cup \bigcup_{j=1}^{\infty} \bigl\{ V_{ij} : j \text{ even } \bigr\}. \label{odd even}
\end{align} 
Clearly $\bar{A} \in \Lambda(1, \lambda) \subseteq \Lambda(N, \lambda)$. To complete the induction, it suffices to show that 
\begin{equation} \label{completing the induction}  
V_i^{\text{\tiny{odd}}}, V_i^{\text{\tiny{even}}} \in \Lambda(N, \lambda); \text{ further, each set is closed.}  
\end{equation} 
In fact, once \eqref{completing the induction} is established, one may intersect $\bar{A}$, $V_i^{\text{\tiny{odd}}}$, and $V_i^{\text{\tiny{even}}}$ with $\bar{U}$. These intersections are closed, remain in $\Lambda(N,\lambda)$ by Lemma \ref{lacunarity under linear operations}, and their union equals $\bar{U}$.
\vskip0.1in
\noindent Let us proceed to prove \eqref{completing the induction}. The placement of $V_{ij}$ in the alternate closed gaps of $A$ implies that 
\[ V_i^{\text{\tiny{odd}}} \cap \bigl[ a_{j+1}, a_j\bigr) \subseteq \begin{cases} V_{ij} \in \Lambda(N-1, \lambda) &\text{ if $j$ is odd }, \\ \{ a_{j+1}\} \in \Lambda(0, \lambda) \subseteq \Lambda(N-1, \lambda) &\text{ otherwise.}
\end{cases}  \] 
From Definition \ref{defn: Lacunary sets}, we deduce that  
\[ V_i^{\text{\tiny{odd}}} \in \Lambda(N, \lambda); \text{ similarly, } V_i^{\text{\tiny{even}}} \in \Lambda(N, \lambda). \] Thus \eqref{U closure cover} says that $\bar{U}$ is covered by at most $C_N = 2C_{N-1}+1$ sets in $\Lambda(N, \lambda)$. 
\vskip0.1in
\noindent Finally, let us confirm that each of the sets $V_i^{\text{\tiny{odd}}}$ and $V_i^{\text{\tiny{even}}}$ is closed. It suffices to only prove this when the parity is odd. Let us consider a convergent sequence $\{v_n \} \subseteq V_i^{\text{\tiny{odd}}}$,  $v_n \rightarrow v$. 
\vskip0.1in
\begin{itemize}
\item If $v = \alpha$, then $v \in V_i^{\text{\tiny{odd}}}$ by the definition \eqref{odd even} of the latter. 
\vskip0.1in
\item If $v \neq \alpha$, then it must lie is some closed gap of $A$: 
\[ v \in I_{\ell} := \bigl[a_{\ell+1}, a_\ell \bigr] \text{ for some $\ell$}. \] A small enough neighbourhood of $v$ must then meet at most two adjacent  gaps of $A$, one of which must be $I_{\ell}$. Thus the convergent sequence $\{v_n\} \subseteq V_i^{\text{\tiny{odd}}}$ must eventually lie in one of these gaps. In other words,  there exists an odd index $j$,
\[ |j - \ell| \leq 1 \text{ such that } v_n \in {I}_j \text{ for all but finitely many } n. \] 
Since $V_{ij}$ is closed, we deduce that $v \in V_{ij} \subseteq V_i^{\text{\tiny{odd}}}$.
\end{itemize}
Thus, in either case, the limit point $v$ lies in $V_i^{\text{\tiny{odd}}}$, proving that the  set is closed. This establishes \eqref{completing the induction}, closing the induction and completing the proof. 
\vskip0.1in
\noindent 
\end{proof} 
\section{Nonclosure of {\tt{AdFinLac}} under algebraic sums} \label{section: nonclosure under algebraic sums proof}
\begin{proof}[Proof of Lemma \ref{lemma: nonclosure under algebraic sums}]
The set $V$ is a monotone lacunary sequence, with lacunarity constant $\frac{1}{2}$; therefore $V \in \Lambda(1, \frac{1}{2})$ by Lemma \ref{Lemma: Lacunary 1} \eqref{Lemma: Lacunary 1 (a)}. To establish the same conclusion for $U$, it suffices to show that for every index $\ell \geq 1$, 
\begin{equation} \label{empty/singleton}
U \cap [2^{-\ell}, 2^{-\ell+1}) \text{ is either empty or a singleton. } 
\end{equation} In fact, for $0 \leq k < M_j$, we observe that 
\begin{align*} 
&2^{-N_j+k} \leq 2^{-N_j +k} + q_{jk} <  2^{-N_j + k+1}, \text{ which means that } \\
&\# \left[U_j \cap \bigl[2^{-\ell}, 2^{-\ell +1} \bigr) \right] = \begin{cases} 1 &\text{ for every $\ell$ with } 0\leq N_j - \ell < M_j, \\ 0 &\text{ otherwise. } \end{cases} 
\end{align*}
The growth condition \eqref{growth condition Mj Nj} implies that  the intervals representing the ranges of $\ell$, namely  
\[  (N_j - M_j, N_j] \text{ are disjoint in $j$, }\]
hence an interval $[2^{-\ell}, 2^{-\ell+1})$ can have non-trivial intersection with $U_j$ for at most one index $j$. This completes the proof of \eqref{empty/singleton}, establishing that $U \in \Lambda(1, \frac{1}{2})$.    
\vskip0.1in
\noindent On the other hand, 
\[ U + V \supseteq \bigcup_{j=1}^{\infty} (U_j + V) \supseteq \bigcup_{j=1}^{\infty} W_j \text{ where } W_j :=\left\{q_{jk} : 0 \leq k < M_j \right\}. \] 
In other words, for every $j \geq 1$, the set $U+V$ contains $W_j$, which is an affine copy of the dyadic rationals of the form $\{k2^{-m_j}: 0 \leq k < 2^{m_j} \}$. By Lemma \ref{lemma: topological closure}, $U+V$ is sublacunary. 
\end{proof} 
\section{Sublacunarity involving components with uncontrolled lacunarity} \label{section: Bateman counterexample 1 proof}  	
In Lemma \ref{lemma: Bateman counterexample 1}, we gave an example of a set that superficially resembles a lacunary set of order 2; it is comprised of lacunary sequences lying within the gaps of the sequence $\{2^{-j}: j \geq 0\}$. However, the lacunarity constants associated with this family of sequences grow to 1. In this section, we prove  Lemma \ref{lemma: Bateman counterexample 1} by establishing the sublacunarity of such a set.
\vskip0.1in
\noindent In preparation for the proof, let $j\geq2$, and set
\[
                    q_j:=1-\frac1j, \quad x_k = x_k(j) :=2^{-j}+2^{-j}q_j^k. \] 
                    Given a non-empty set of integers 
\begin{equation}   \label{subset A of integers}                  \mathcal A\subseteq\{n,n+1,\ldots,n+N\}, \text{ we define } 
                    \Omega_{\mathcal A} = \Omega_{\mathcal A}(j) 
                       :=\{x_k(j):k\in\mathcal A\}.
\end{equation} 
Let $B=\{b_i:i\geq1\} \in {\tt{MonLac}}(\lambda)$ be a decreasing sequence converging to $\beta$, with
\begin{equation} \label{B-criterion-1}
                    \Omega_{\mathcal A}\subseteq[\beta,b_1].
\end{equation} 
After deleting finitely many initial terms of $B$, if necessary, we may assume that
\begin{equation} \label{B-criterion-2}
                    \max\Omega_{\mathcal A}\in[b_2,b_1).
\end{equation} 
The following key lemma provides an upper bound on the number of gap intervals of $B$ that $\Omega_{\mathcal A}$ can meet. 
\begin{lemma}\label{lem:gap-count-slow-lacunary}
For every $\lambda \in (0, 1)$, there is a constant $C_{\lambda} > 0$ with the following property. Suppose that $N \geq 1$, and choose $j$ sufficiently large relative to $N$ so that 
\begin{equation} \label{q_j^N} 
        q_j^N=\Bigl(1-\frac1j\Bigr)^N\geq \frac12. 
\end{equation} 
%If \(j\) is sufficiently large depending on \(N\), then
If $B$ is any decreasing ${\tt{MonLac}}(\lambda)$ sequence normalized as above, then at most $C_{\lambda} \log(N+1)$ gap intervals of $B$ can intersect $\Omega_{\mathcal A}$:
\begin{equation} \label{how many i}
 \#\left\{i\geq1:
   \Omega_{\mathcal A}\cap[b_{i+1},b_i)\neq\varnothing
 \right\}
 \leq C_{\lambda} \log(N+1).
\end{equation} 
%with an absolute implied constant.
\end{lemma}
\noindent We will prove the lemma first, then use it to prove Lemma \ref{lemma: Bateman counterexample 1}. 
\vskip0.1in
\subsubsection{Proof of Lemma \ref{lem:gap-count-slow-lacunary}}
\begin{proof}
Given a non-empty subset $\mathcal A$ of integers as in \eqref{subset A of integers}, let us set  
\[
        k_-:=\min\mathcal A,
        \qquad
        k_+:=\max\mathcal A,
        \qquad
        D:=k_+-k_-.
\]
Thus $0\leq D\leq N$, and, since $x_k$ is decreasing in \(k\),
\[
        \max\Omega_{\mathcal A}=x_{k_-},
        \qquad
        \min\Omega_{\mathcal A}=x_{k_+}.
\]
If $D=0$, then $\Omega_{\mathcal A}$ is a singleton and the
intended conclusion \eqref{how many i} is immediate. We may therefore suppose that $D\geq 1$.
\vskip0.1in
\noindent The main claim concerning $\mathcal A$ is the bound 
\begin{equation} \label{kk-}
 \frac{x_{k_-}-\beta}{x_k-\beta}
 \leq 2D\leq 2N
 \qquad
 \text{for every \(k\in\mathcal A\setminus\{k_+\}\)}.
\end{equation} 
We will prove this estimate momentarily. Assuming it for now, the proof is completed as follows. Suppose that $i \geq 2$ is an index that contributes to the set occurring on the left side of \eqref{how many i}; in other words, 
\begin{equation} \label{where is xk} 
        x_k\in[b_{i+1},b_i) \text{ for some } k\in\mathcal A\setminus\{k_+\}. 
\end{equation} 
Since $B \in {\tt{MonLac}}(\lambda)$ and decreases to $\beta$, we have 
\begin{equation} \label{lacunary B} 
        b_i-\beta
        \leq \lambda^{i-2}(b_2-\beta).
\end{equation} 
Moreover, our normalization
\begin{equation} \label{normalization outcome} 
        x_{k_-}=\max\Omega_{\mathcal A}\in[b_2,b_1) \; \text{ gives } b_2-\beta\leq x_{k_-}-\beta. 
        \end{equation}  The relations \eqref{where is xk}, \eqref{lacunary B} and \eqref{normalization outcome} jointly imply 
\begin{equation*} 
 x_k-\beta
 <
 b_i-\beta
 \leq
\lambda^{i-2}(x_{k_-}-\beta).
\end{equation*} 
Together with \eqref{kk-}, the above estimate yields
\[
        \lambda^{-(i-2)}
        <
        \frac{x_{k_-}-\beta}{x_k-\beta}
        \leq 2N.
\]
Hence every gap containing a point of
\(\Omega_{\mathcal A}\setminus\{x_{k_+}\}\) has index
\[
        i<2+\log_{\lambda^{-1}}(2N).
\]
The remaining point \(x_{k_+}\) can meet at most one additional gap.
It follows that
\[
 \#\bigl\{i\geq1:
      \Omega_{\mathcal A}\cap[b_{i+1},b_i)\neq\varnothing
   \bigr\}
 \leq 3+\left\lceil\log_{\lambda^{-1}}(2N)\right\rceil
 \leq C_{\lambda}\log(N+1)
\]
for a suitable constant $C_{\lambda}>0$ depending only on $\lambda$. This is the desired conclusion \eqref{how many i}. 
\vskip0.1in
\noindent In order to complete the proof of the lemma, it remains to check \eqref{kk-}. Since
\(\beta\leq x_{k_+}<x_k\leq x_{k_-}\) for \(k<k_+\), we have
\begin{equation} \label{3 fractions}
 \frac{x_{k_-}-\beta}{x_k-\beta}
 \leq
 \frac{x_{k_-}-x_{k_+}}{x_k-x_{k_+}}
 \leq
 \frac{x_{k_-}-x_{k_+}}
      {x_{k_+-1}-x_{k_+}}.
\end{equation} 
Recalling that
\[
        x_k=2^{-j}+2^{-j}q_j^k,
\]
the rightmost expression in \eqref{3 fractions} equals 
\[
 \frac{q_j^{k_-}-q_j^{k_+}}
      {q_j^{k_+-1}-q_j^{k_+}}
 =
 q_j^{1-D}\frac{1-q_j^D}{1-q_j}
 =
 q_j^{1-D}\sum_{r=0}^{D-1}q_j^r.
\]
Because \(D\leq N\) and \(q_j^N\geq \frac12\), it follows that
\[ q_j^{1-D}\sum_{r=0}^{D-1}q_j^r \leq D q_j^{-N} \leq 2D.\]
This establishes \eqref{kk-}. 
\end{proof}
\vskip0.1in
\noindent {\em{Remark: }} \label{B increasing remark} We stated and proved Lemma \ref{lem:gap-count-slow-lacunary} for $B \in {\tt{MonLac}}(\lambda)$ decreasing. By reflecting about the origin, the same conclusion holds when $B \in {\tt{MonLac}}(\lambda)$ is
increasing. The gap intervals are now written as $[b_i,b_{i+1})$,
and the assumptions \eqref{B-criterion-1}, \eqref{B-criterion-2} are replaced by the following natural counterparts:
\begin{equation} \label{increasing B changes} 
\Omega_{\mathcal A} \subseteq [b_1, \beta], \quad \min \Omega_{\mathcal A} \in [b_1, b_2). 
\end{equation} 
Under these assumptions, the verification of Lemma \ref{lem:gap-count-slow-lacunary} in the $B$ increasing case is left to the reader. 
\vskip0.1in
\noindent The following consequence of Lemma \ref{lem:gap-count-slow-lacunary} will be used repeatedly in proving Lemma \ref{lemma: Bateman counterexample 1}. Let us record it here in the form that we need. 
\begin{corollary}
Fix $N\geq1$, and choose $j$ sufficiently large that \eqref{q_j^N} holds. Then for every $\lambda \in (0, 1)$, every integer $m\geq 0$, every set
$E\in\Lambda(m,\lambda)$, and every non-empty set $\mathcal A$ as in \eqref{subset A of integers}, the following statement holds: 
\begin{equation} \label{induction: E & Omega_A} 
        \#(E\cap\Omega_{\mathcal A})\leq H_N^m, \; \text{ with } \; H_N:=C_\lambda\log(N+1)+2. 
\end{equation} 
\end{corollary}
\begin{proof} 
We prove \eqref{induction: E & Omega_A} by induction on $m$. When $m=0$, the set $E$
has cardinality at most one, and hence
\[
        \#(E\cap\Omega_{\mathcal A})\leq1=H_N^0.
\]
Now suppose that $m\geq1$ and that the claim holds at level $m-1$. Let 
\[ E \in \Lambda(m, \lambda), \text{ with } B=\{b_i:i\geq1\}\in {\tt{MonLac}}(\lambda) \text{ as a special sequence.} \] 
For every gap interval $I$ of $B$, Definition \ref{defn: Lacunary sets} gives $E\cap I\in\Lambda(m-1,\lambda)$. The induction hypothesis, applied to $E \cap I$, therefore implies 
\begin{equation} \label{E in a gap}
        \#(E\cap\Omega_{\mathcal A}\cap I)
        \leq H_N^{m-1}.
\end{equation} 
Define \[ \mathcal A_{E} := \bigl\{k \in \mathcal A : x_k \in E \bigr\}, \text{ so that } \mathcal A_E \subseteq \mathcal A \subseteq \{n, n+1, \ldots, n+N\}.  \]
Then 
\[ E \cap \Omega_{\mathcal A} = \Omega_{\mathcal A_E}. \] 
If $\mathcal A_E$ is empty, there is nothing to prove, so we assume $\mathcal A_E \neq \emptyset$. Since $B$ is
a special sequence for $E$, the set $\Omega_{\mathcal A_E}$ is
contained in the interval spanned by $B$. If $B$ is decreasing, then after discarding any initial terms of $B$ that precede
$\Omega_{\mathcal A_E}$, we may ensure that the conditions \eqref{B-criterion-1} and \eqref{B-criterion-2} hold. If $B$ is increasing, a similar discarding of initial terms ensures \eqref{increasing B changes}. In either case, this deletion does not remove any gap interval meeting $\Omega_{\mathcal A_E}$. With all the hypotheses verified, we may now apply Lemma \ref{lem:gap-count-slow-lacunary} (regardless of $B$ increasing or decreasing, as per the remark on page \pageref{B increasing remark}) with $\mathcal A$ replaced by $\mathcal A_E$.  This yields   
\begin{equation}  \label{number of relevant gaps} 
\begin{aligned}
&\# \bigl\{I \text{ gap of } B : E \cap \Omega_{\mathcal A} \cap I  \neq \emptyset \bigr\} \\ &\hskip1in = \# \bigl\{I \text{ gap of } B :  \Omega_{\mathcal A_E} \cap I  \neq \emptyset \bigr\}\leq C_\lambda\log(N+1). 
\end{aligned} 
\end{equation}  
The gaps $I$ cover the interval spanned by $B$, with the possible exception of the two extreme points. Summing over the occupied gaps with the possible inclusion of the extreme
points, and invoking \eqref{E in a gap} and \eqref{number of relevant gaps} we obtain
\[
\begin{aligned}
        \#(E\cap\Omega_{\mathcal A}) &\leq 2 + \sum_{I} \left\{\# (E\cap\Omega_{\mathcal A} \cap I) : I \text{ gap of } B, E \cap \Omega_{\mathcal A} \cap I  \neq \emptyset \right\} \\
        &\leq  2 + C_\lambda\log(N+1)H_N^{m-1}  \leq H_N^m.
\end{aligned}
\]
This completes the induction, proving \eqref{induction: E & Omega_A}.
\end{proof}

\vskip0.1in
\subsubsection{Proof of Lemma \ref{lemma: Bateman counterexample 1}} 
\begin{proof}
We first verify that $U$ satisfies the definition of finite-order
lacunarity in \eqref{Bateman-fol}. Let $L:=\{2^{-j}:j\geq1\}$; 
this is a lacunary sequence with lacunarity constant $\frac12$. For each fixed $j\geq 2$, the sequence $L_j:=\{u_{jk}:k\geq1\}$ as given by \eqref{example: uncontrolled lacunary blocks} converges to $2^{-j}$ and has lacunarity constant $q_j=1-\frac1j$. For $j=1$, the set $L_1 = \{1/2\}$ is a singleton, and causes no difficulty. Thus $U$, which is a countable union of the sets $L_j$, satisfies the requirement of \eqref{Bateman-fol} with order $2$.
\vskip0.1in
\noindent It remains to prove that $U \in {\tt{SubLac}}$ in the sense of
Definition~\ref{defn: Admissible finite order lacunarity}. Suppose, to the contrary,
that $U\in{\tt{AdFinLac}}$. This means one can find $\lambda\in(0,1)$ and integers
$N_0,R_0\geq1$ such that $U\in\Lambda(N_0,\lambda;R_0)$, i.e., 
\begin{equation}   \label{U covered by Ei}
\exists \text{ sets } E_1,\ldots,E_{R_0} \in\Lambda(N_0,\lambda) \text{ with } 
        U\subseteq\bigcup_{r=1}^{R_0}  E_r.
\end{equation} 
We now choose $N$ sufficiently large that
\begin{equation} \label{choice of NN0}
        N+1>R_{0}\bigl(C_\lambda\log(N+1)+2\bigr)^{N_0}.
       \end{equation} 
Having fixed $N$, we next choose $j$ sufficiently large so as to satisfy  \eqref{q_j^N}; in other words, $q_j^N\geq \frac12$. For this choice of $j$, consider the finite block
\begin{align*}
&        \Omega_{N,j}
        :=\{u_{jk}:1\leq k\leq N+1\}\subseteq U, \; \text{ so that } \\ &\Omega_{N, j} = \Omega_{\mathcal A}(j) \text{ with } \mathcal A = \{1, \ldots, N+1\}
\end{align*} 
in the notation of \eqref{subset A of integers}. On one hand, 
this set has cardinality $N+1$. On the other, \eqref{induction: E & Omega_A} and
\eqref{U covered by Ei} give
\[N+1
        =\#\Omega_{N,j} \leq
        \sum_{r=1}^{R_0}\#(E_r\cap\Omega_{N,j})  \leq
        R_0\bigl(C_\lambda\log(N+1)+2\bigr)^{N_0},
\]
contradicting \eqref{choice of NN0}. Therefore $U\notin {\tt{AdFinLac}}$, i.e. $U$ is sub-lacunary in the sense of Definition \ref{defn: Admissible finite order lacunarity}.
\end{proof}

\section{A sublacunary set consisting of sequences with auxiliary accumulation} \label{section: Bateman counterexample 2 proof}  	
In Lemma \ref{lemma: Bateman counterexample 2}, we gave an example of a set given by a union of ${\tt{MonLac}}(\frac{1}{3})$ sequences, whose limit points form a ${\tt{MonLac}}(\frac{1}{2})$ sequence. However, the family of lacunary sequences are not restricted to lie within lacunary gaps, causing sub-lacunary build-up elsewhere. We prove this lemma in this section.   
\begin{proof}[Proof of Lemma \ref{lemma: Bateman counterexample 2}]
				The central lacunary sequence in the construction of $U$ is 
				\[ \mathscr{L} := \{ 2^{-j} \} \in {\tt{MonLac}}(1/2). \] For every fixed $j \geq 0$, the lacunary sequence 
				\[ \mathscr{L}_j :=\{u_{jk} : k \geq 0\} \in {\tt{MonLac}}(1/3) \text{ converges to $2^{-j}$, a point of $\mathscr{L}$}. \] 
				The set $U$ is the union (in $j$) of the lacunary sequences $\mathscr{L}_j$, and therefore obeys the definition of $N^{\text{th}}$ order lacunarity as given in \cite{Bateman}, with $N=2$.  
				 \vskip0.1in
				 \noindent However, $U \not\in {\tt{AdFinLac}}$, according to Definition \ref{defn: Admissible finite order lacunarity}. Indeed, 
				 \[ U \supseteq \{ u_{j0} : j \geq 0\} = \{ q_j : j \geq 0 \}. \]
				 The last set is sublacunary, as proved in Lemma \ref{lemma: topological closure} and Corollary \ref{corollary: uncountable sublacunary}, and therefore so is $U$.  Indeed, the source of sublacunarity is the starting point $u_{j0}$ of the sequence $\mathscr{L}_j$. For all large $j$, this initial point, which lies in $[\frac{9}{10}, 1]$, is far away from $[2^{-j}, 2^{-j+1}])$, the gap interval of $\mathscr{L}$ where $\mathscr{L}_j$ eventually converges. Without the enforced Euclidean separation of the gap intervals, the points are able to create sublacunary accumulation elsewhere. 
\end{proof} 
	
	\chapter{Appendix C: Proofs from Section \ref{SECTION: LACUNARY PROPERTIES}} \label{section: lacunary properties proofs}
In Section \ref{SECTION: LACUNARY PROPERTIES}, we recorded a few properties related to admissible finite-order lacunarity. We prove these properties here. 
\section{Proof of Lemma \ref{lemma: nonclosure under countable union}} \label{nonclosure proof section}
This lemma refers to the non-closure of admissible finite-order lacunarity under countable unions. For every $N \geq 0$ and $\lambda \in (0,1)$, we will establish the existence of a sublacunary set that is contained in a countably infinite union of sets in $\Lambda(N, \lambda; R)$, with varying $R$.  
\begin{proof}
For $N = 0$, let us consider the finite set of dyadic rationals $\mathbb Q_m$ given in \eqref{dyadic Q_m}. Each set $\mathbb Q_m$ has cardinality $2^m$, i.e., it is a $2^m$-fold union of singletons. This means 
\begin{equation} \label{dyadic rationals sublacunary} 
\mathbb Q_m \in \Lambda(0, \lambda; 2^m) \text{ for any } \lambda \in (0, 1). 
\end{equation} 
However, the countable union of $\mathbb Q_m$ over non-negative integers $m$ is the collection of all dyadic rationals in $[0,1]$, which is dense. It follows from Corollary \ref{corollary: uncountable sublacunary} that 
%{\color{red} This does not hold. Whenever $R$ is finite, it is adimissible lacunary.}
\[ \bigcup_{m=1}^{\infty} \mathbb Q_m \in {\tt{SubLac}}.  
%\text{ however, } ,\text{ so } \bigcup_{m=1}^{\infty} \mathbb Q_m \in \bigcup_{R=1}^{\infty} \Lambda(0, \lambda; R).
\] 
This establishes \eqref{nonclosure} for $N = 0$. 
\vskip0.1in 
\noindent For $N \geq 1$ and $\lambda \in (0,1)$, let us fix a set 
\begin{equation} \label{V} V \subseteq [0,1), \quad V \in \Lambda(N, \lambda), \quad 1 \in V. \end{equation}   
For instance, $V$ could be chosen as  
 \[ V = \Bigl\{ \lambda^{j_1} - \lambda^{j_1 + j_2}(1 - \lambda) - \ldots - \lambda^{j_1+\ldots +j_N} (1 - \lambda)^{N-1} : j_1, \ldots, j_N \geq 0\Bigr\}. \] 
 For $V$ as in \eqref{V}, and $1 \leq r \leq 2^{m}$, we define
 \[ V_r := \frac{r-1}{2^m} + 2^{-m} V \subseteq \Bigl[\frac{r-1}{2^m}, \frac{r}{2^m} \Bigr]. \]
 Since $V_r$ is an affine copy of $V$, we deduce that $V_r \in \Lambda(N, \lambda)$, by \eqref{V} and the invariance of the lacunarity order under affine maps (conclusion \eqref{linear-invariance} of Lemma \ref{lacunarity under linear operations}). Let us now set
 \[ \tilde{U}_m := \bigcup_{r=1}^{2^m} V_r, \quad \text{ so that } \quad V_r \in \Lambda(N, \lambda; 2^m).  \]  
 However, the countable union $U$ of the sets $\tilde{U}_m$ over all non-negative integers $m$ is sublacunary. This follows from \eqref{V}: the inclusion $1 \in V$ means  
 \[ \frac{r}{2^m} \in V_r, \text{ which translates to } \mathbb Q_m \subseteq \tilde{U}_m, \text{ and therefore } \bigcup_{m=0}^{\infty} \mathbb Q_m \subseteq  U. \]
In view of \eqref{dyadic rationals sublacunary}, $U$ contains a sublacunary set and is therefore itself sublacunary. 
\end{proof}

\section{Proof of Lemma \ref{lemma: special sequence choice}}	 
Lemma  \ref{lemma: special sequence choice} posits that a set $U \in \Lambda(N, \lambda)$ may possess  multiple special sequences; in particular, $U$ has a special sequence $A$ that satisfies a two-sided inequality of the form \eqref{lacunary above and below}; the ratio of the distances of the limit $a = \lim A$ from two consecutive elements $A$ is bounded above and below by a power of the lacunarity constant. We prove this lemma here.
\label{section: proof of special sequence choice} 
\begin{proof}
Given $U \in \Lambda(N, \lambda)$, let $B = \{b_{\ell} : \ell \geq 1\}$ be a special sequence of $U \in \Lambda(N, \lambda)$, converging to the limit $b$. After a linear transformation if necessary, we can assume that the sequence $B$ decreases strictly monotonically to $b = 0$. It follows from Definition \ref{defn: Lacunary sets} of $\Lambda(N, \lambda)$ that 
\begin{equation} \label{properties of B} 
 U \subseteq \bigl[0, b_1\bigr), \; \; 0 \leq b_{\ell+1} \leq \lambda b_{\ell}, \;  \; U \cap [b_{\ell+1}, b_\ell) \in \Lambda(N-1, \lambda), \; \ell \geq 1.  
\end{equation} 
A priori, $B$ only obeys the defining condition \eqref{lacunarity constant} of a lacunary sequence, and therefore could decay arbitrarily fast. To quantify this decay we define, for every $\ell \geq 1$, the unique integer $n_{\ell}$ obeying 
\begin{equation}   \lambda^{n_{\ell}+1}  b_{\ell} < b_{\ell+1} \leq \lambda^{n_{\ell}} b_{\ell} \leq \lambda b_{\ell}. \label{what is n_l} \end{equation} 
In other words, $n_{\ell}$ is the largest integer $n$ such that $\lambda^n b_{\ell}$ is at least as large as $b_{\ell+1}$. The lacunarity property \eqref{what is n_l} of $B$ ensures that $n_{\ell} \geq 1$. 
\vskip0.1in 
\noindent We now insert $(n_{\ell}-1)$ new elements between $[b_{\ell+1}, b_{\ell})$ to control the lacunarity constant from below.  Specifically, let us consider the finite sequence of points $C_{\ell} := \{c_{\ell  k} :  0 \leq k \leq n_{\ell}\}$ in between $b_{\ell+1}$ and $b_{\ell}$, including both endpoints:
\begin{equation}  
%c_{\ell 0} := b_{\ell+1} \text{ if } n_{\ell} = 1;  \; \text{ for } n_{\ell} \geq 2, \text{ we set } 
c_{\ell k} := \begin{cases} \lambda^{k} b_{\ell} &\text{ if } 0 \leq k < n_{\ell}, \\ b_{\ell+1} &\text{ if } k = n_{\ell}. \end{cases} \label{defn: c}  \end{equation} 
Let $A = \{a_j: j \geq 1\}$ denote an ordered enumeration of the combined elements of $\{C_{\ell}: \ell \geq 1\}$, with $a_{j+1} < a_j$ for all $j \geq 1$. Thus each gap interval of $A$ lies within a gap of $B$. Figure \ref{fig:special A} depicts the construction of these finer gaps. We claim that $A$ serves as a special sequence of $U$ obeying \eqref{lacunary above and below}.  
\vskip0.1in
\noindent Let us verify the condition \eqref{lacunary above and below} first. Given two consecutive elements $a_{j}, a_{j+1} \in A$, let $\ell$ be the unique index such that $a_j \in (b_{\ell+1}, b_{\ell}]$.  If $n_{\ell}=1$, then 
\[ C_{\ell} \cap (b_{\ell+1}, b_{\ell}) = \emptyset,  \text{ which means } a_{j} = b_{\ell} \text{ and } a_{j+1} = c_{\ell 0} = b_{\ell+1}, \] and the condition \eqref{lacunary above and below} follows  immediately from \eqref{what is n_l}.  However, if the index $\ell$ is such that $n_{\ell} \geq 2$, we are led to consider the following cases. 
\begin{itemize} 
\item First suppose that $a_j = b_{\ell}$ for some $\ell$. Then it follows from \eqref{defn: c} that 
\[ a_j = c_{\ell 0}, \; \text{ and } \; a_{j+1} = c_{\ell 1} = \lambda b_{\ell} = \lambda a_j. \] Thus \eqref{lacunary above and below} holds with equality on the right hand side. 
\vskip0.1in
\item This leaves the case where $b_{\ell+1} < a_j < b_{\ell}$. In this situation 
\[ a_j = c_{\ell k} \text{ for some } 1 \leq k < n_{\ell}. \] 
If $k < n_{\ell}-1$, then $a_{j+1} = c_{\ell, k+1} = \lambda^{k+1} b_{\ell} = \lambda c_{\ell k} =  a_j$ as before. 
\vskip0.1in
\noindent If $k = n_{\ell}-1$, then $a_j = \lambda^{n_{\ell}-1} b_{\ell}$, $a_{j+1} = b_{\ell+1}$. In this case \eqref{what is n_l} implies
\[ 
\frac{a_{j+1}}{a_j} = \frac{b_{\ell+1}}{\lambda^{n_{\ell}-1} b_{\ell}} > \frac{\lambda^{n_{\ell} + 1} b_{\ell}}{\lambda^{n_{\ell}-1}b_{\ell}} = \lambda^2, \] completing the verification of \eqref{lacunary above and below}. 
\end{itemize}
\vskip0.1in
\noindent Finally, let us confirm that $A$ is a special sequence for $U$. For $j \geq 1$, the interval $[a_{j+1}, a_j)$ is contained in $[b_{\ell+1}, b_{\ell})$ for some $\ell$. It then follows from \eqref{properties of B} that 
\begin{equation*}
U \cap [a_{j+1}, a_j) \subseteq U \cap [b_{\ell+1}, b_{\ell}), 
\end{equation*} 
and the latter set lies in $\Lambda(N-1, \lambda)$.  Thus $A$ serves as a special sequence for $U$.  This concludes the proof of the lemma. 
\end{proof} 
\begin{center}
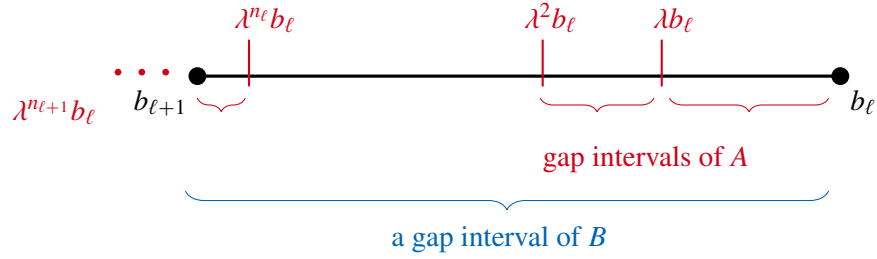
\begin{figure}
\begin{tikzpicture}[scale=1.05, line cap=round, line join=round]

% Colors
\definecolor{myred}{RGB}{210,0,25}
\definecolor{myblue}{RGB}{0,100,180}

% Main black interval
\draw[very thick] (1.7,0) -- (9.8,0);
\fill (1.7,0) circle (3.2pt);
\fill (9.8,0) circle (3.2pt);

% Left dotted red points
\foreach \x in {0.7,1.0,1.3}
  \fill[myred] (\x,0.05) circle (1.2pt);

% Labels for endpoints
\node[below left] at (1.7,-0.05) {$b_{\ell+1}$};
\node[below right] at (9.8,-0.05) {$b_\ell$};

% Red labels and ticks
\draw[myred, thick] (2.35,-0.12) -- (2.35,0.45);
\node[myred, above] at (2.55,0.42) {$\lambda^{n_\ell}b_\ell$};

\draw[myred, thick] (6.05,-0.12) -- (6.05,0.45);
\node[myred, above] at (6.15,0.42) {$\lambda^2 b_\ell$};

\draw[myred, thick] (7.55,-0.12) -- (7.55,0.45);
\node[myred, above] at (7.7,0.42) {$\lambda b_\ell$};

\node[myred, left] at (0.55,-0.45) {$\lambda^{n_{\ell+1}}b_\ell$};

% Red braces for gap intervals of A
\draw[myred, decorate, decoration={brace, amplitude=5pt, mirror}]
  (6.02,-0.28) -- (7.45,-0.28);
\draw[myred, decorate, decoration={brace, amplitude=5pt, mirror}]
  (7.65,-0.28) -- (9.65,-0.28);
\draw[myred, decorate, decoration={brace, amplitude=5pt, mirror}]
  (1.7,-0.28) -- (2.3,-0.28);
\node[myred] at (7.35,-1.05) {gap intervals of $A$};

% Blue brace for gap interval of B
\draw[myblue, decorate, decoration={brace, amplitude=7pt, mirror}]
  (1.55,-1.45) -- (9.65,-1.45);

\node[myblue] at (5.5,-2.05) {a gap interval of $B$};

\end{tikzpicture}
\caption{\small{The diagram shows how a special sequence $A$ obeying controlled lacunary decay of the form \eqref{lacunary above and below} is constructed from an arbitrary special sequence $B$, as claimed in Lemma \ref{lemma: special sequence choice}.}} \label{fig:special A}
\end{figure} 
\end{center} 
\section{Finite-order lacunarity under bi-Lipschitz maps} \label{section: bi-Lipschitz maps proof} 
\noindent A bi-Lipschitz map is a function $F: I \rightarrow \mathbb R$ obeying \eqref{def: bi-Lipschitz}, which is a two-sided smoothness condition. Proposition \ref{prop: AdFinLac under F} states that lacunarity order remains unchanged under bi-Lipschitz maps. We prove this proposition in this section. The next lemma verifies this statement for $N=1$. 
%Lemma \ref{lacunarity under linear operations} established that the property of  finite order lacunarity remains unchanged under bi-Lipschitz maps. In fact, a stronger version of Lemma \ref{lacunarity under linear operations} \eqref{linear-invariance} is true for $U \in {\tt{AdFinLac}}$, as we shall see in Proposition \ref{prop: AdFinLac under F} below. The intermediate lemmas provide the necessary steps.  

		\begin{lemma} \label{lemma: Lac1 under bi-Lipschitz}
			The collection 
			\[ \bigcup_{R=1}^{\infty} \Lambda(1, \lambda; R) \] which consists of sets that are admissible lacunary of order at most 1 according to Definition \ref{defn: Admissible finite order lacunarity}, is invariant under bi-Lipschitz maps. 
			\vskip0.1in
			\noindent Specifically, given a lacunarity constant $\lambda \in (0,1)$ and a function $F$ obeying \eqref{def: bi-Lipschitz},
			there exists an integer $C_1 = C_1(\lambda, F) \geq 1$ with the following property: if $U \subseteq I$ is such that $U \in \Lambda(1, \lambda)$, then $F(U)$ is contained in the union of at most $C_1$ monotone lacunary sequences.  In view of Lemma \ref{Lemma: Lacunary 1} \eqref{Lemma: Lacunary 1 (a)}, this means that
			\begin{equation} F(U) \in \Lambda(1, \lambda; C_1). \end{equation} 
			\end{lemma}
			\begin{proof}
				To begin, let us observe that any function $F$ obeying the bi-Lipschitz condition \eqref{def: bi-Lipschitz} must be continuous, and strictly monotone. Indeed, if $F$ were not monotone, there would be three distinct points $x < y< z$ in $I$ such that 
				\[ \text{ either } F(x), F(z) < F(y) \qquad \text{ or } \qquad F(x), F(z) > F(y).\]
				In either case one can find, by the intermediate value theorem, two points $u, v \in I$, $u \in [x, y)$, $v \in (y, z]$ such that $F(u) = F(v)$. This contradicts the injectivity of $F$, as implied by the left inequality of \eqref{def: bi-Lipschitz}.  
				\vskip0.1in
				\noindent Now suppose $U \in \Lambda(1, \lambda)$. By Lemma \ref{Lemma: Lacunary 1} \eqref{Lemma: Lacunary 1 (c)}, $U$ is contained in the union of at most two monotone, lacunary sequences of the form $A = \{a_j : j \geq 0\}$ with the same limit $\alpha$, satisfying \eqref{lacunarity constant}. By the conclusion of the preceding paragraph, $F(A)$ is also strictly monotone. To verify the statement of the lemma, it therefore suffices to prove that $F(A)$ is contained in the union of at most $\frac{C_1}{2}$ sets in $\Lambda(1, \lambda)$, for some even integer $C_1 \geq 2$ to be determined.
				\vskip0.1in
				\noindent Iterating the lacunarity condition on $A$, namely \eqref{lacunarity constant}, $n$ times and combining with the bi-Lipschitz condition \eqref{def: bi-Lipschitz}, we obtain   
				\[ \bigl| F(a_{j+n}) - F(\alpha)\bigr| \leq \mathtt L |a_{j+n} - \alpha| \leq \mathtt L \lambda^n |a_j-\alpha| \leq \mathtt L^2 \lambda^n |F(a_j) - F(\alpha)|. \]  
				Choosing $n_0 = n_0(\mathtt L, \lambda)\geq 1$ large enough so that 
				$\mathtt L^2 \lambda^{n_0-1} < 1$, we note that any sub-sequence of $F(A) = \{F(a_j)\}$ whose consecutive indices $j$ are separated by $n_0$ is lacunary. In other words, for any $k \geq 1$,  
				\[ \bigl\{ F(a_{k + jn_0}) : j \geq 0 \bigr\} \in {\tt{MonLac}}(\lambda). \]  Thus $F(A)$ is the union of $n_0$ lacunary sequences $F(A_k) \in {\tt{MonLac}}(\lambda)$:
				\begin{equation} \label{F(lacunary sequence)} F(A) = \bigcup_{k=0}^{n_0-1} F(A_k), \quad A_k := \{ a_{k + jn_0} : j \geq 0 \}.  \end{equation}
	%			Each lacunary sequence $F(A_k)$ is monotone, since the function $F$ and the sequence $A$ both are. 
	According to Lemma \ref{Lemma: Lacunary 1} \eqref{Lemma: Lacunary 1 (a)}, $F(A_k) \in \Lambda(1, \lambda)$. Thus, \eqref{F(lacunary sequence)} says that $F(A)$ is contained in an $n_0$-fold union of monotone lacunary sequences in $\Lambda(1, \lambda)$, establishing desired conclusion with $C_1 = 2n_0$.
				\end{proof} 
		\subsection{Proof of Proposition \ref{prop: AdFinLac under F}} \label{prop proof: AdFinLac under F}  
	
	\begin{proof} 
	By the affine invariance in Lemma \ref{lacunarity under linear operations}, it suffices to consider the function $F$ to be increasing. Indeed, if $F$ is decreasing, then $-F$ is increasing with the same bi-Lipschitz constant, and reflection of its image preserves the lacunarity class and covering constant.
	%	In the proof of Lemma \ref{lemma: Lacunary 1 under bi-Lipschitz}, we have already noted that $F$ must be strictly monotone. Without loss of generality and replacing $F$ by $-F$ if necessary, the function $F$ can be chosen to be monotone increasing. 
		\vskip0.1in
		\noindent We prove the desired conclusion \eqref{lacunary-induction} by induction on $N$. Lemma \ref{lemma: Lac1 under bi-Lipschitz} covers the base case $N=1$. 
 \vskip0.1in
 \noindent Let us proceed to the induction step. The idea is that the image of a special sequence under $F$ remains lacunary, possibly after a finite decomposition, by virtue of Lemma \ref{lemma: Lac1 under bi-Lipschitz}. We show that each lower-order component of $U$ lying in a gap interval of the special sequence maps into a controlled finite union of lower-order lacunary sets, each contained within a gap of the image sequence. 
 \vskip0.1in 
 \noindent To make this precise, suppose that $U \in \Lambda(N, \lambda)$. Let 
 \[ A = \{ a_j : j \geq 0\} \text{ denote a special sequence of $U$, and set $\alpha = \lim A$.} \] If $A$ obeys the lacunarity condition \eqref{lacunarity constant}, then the relation \eqref{F(lacunary sequence)} gives a finite cover of $F(A)$ by monotone lacunary sequences of the form 
 \[ B_k = F(A_k) = \{ F(a_j) : j \equiv k \text{ mod } n_0\} \in {\tt{MonLac}}(\lambda) \text{ for } 0 \leq k < n_0. \] 
 Let us observe that all the sequences $B_k$ share a common limit, namely $F(\alpha)$. 
 \vskip0.1in
 \noindent Definition \ref{defn: Lacunary sets} of $\Lambda(N, \lambda)$ asserts that $U$ is contained in $(\inf (A), \sup(A))$, i.e. it lies in the interval with endpoints $\alpha$ and $a_0$. Since $F$ is continuous and strictly monotone, its image $F(U)$ is contained in the interval with endpoints $F(\alpha)$ and $F(a_0)$, which are the extremal points of $B_0$. To complete the induction, and hence the proof of Proposition \ref{prop: AdFinLac under F}, we aim to show that 
 \begin{equation} \label{F(U) end of induction} 
 F(U) \text{ is contained in the union of at most $n_0 C_{N-1}$ sets in $\Lambda(N, \lambda)$}, 
 \end{equation} with $B_0 = F(A_0)$ serving as the special sequence for all of them. This will complete the induction step \eqref{lacunary-induction} with 
 \[C_N = n_0 C_{N-1} = \frac{C_1 C_{N-1}}{2}. \]  
\vskip0.1in
\noindent With this goal in  mind, let $J = [a, b)$ be a gap interval of $B_0$. Both $a$ and $b$ are elements of $B_0 = F(A_0)$, so they are of the form $F(a_j)$ for two consecutive indices $j$ that are 0 mod $n_0$. This means there are indices $\ell, \ell' = 0$ mod $n_0$ such that 
\[ a = F(a_{\ell}), \; b = F(a_{\ell'}) \in B_0, \quad \ell' = \begin{cases} \ell - n_0 &\text{ if  $A$ is decreasing}, \\ \ell+n_0 &\text{ if $A$ is increasing.}  \end{cases} \] 
 In either case, and in light of the strict monotonicity of $F$, the interval $J_{\ell} := [a, b)$ decomposes into $n_0$ subintervals with disjoint interiors:  
   \begin{align}  
   	J = J_{\ell} &= [a, b) = \bigcup_{r=0}^{n_0-1} \widetilde{I}_\ell(r), \; \text{ where } \; \widetilde{I}_\ell(r) = F(I_\ell(r)),  \label{interval decomposition}  \\ 
   	{I}_\ell(r) &:= \begin{cases} \bigl[a_{\ell'+r+1}, a_{\ell'+r} \bigr) &\text{ for $A$  decreasing,} \\  
   	\bigl[ a_{\ell+r}, a_{\ell+r+1} \bigr) &\text{ for $A$ increasing.}
   \end{cases} \nonumber
   %\bigl[ F(a_{\ell+n_0}), F(a_{\ell})\bigr] = \bigcup_{r=0}^{n_0-1} \bigl[ F(a_{\ell + r+1}), F(a_{\ell+r}) \bigr].  
   \end{align} 
We will show momentarily, using the induction hypothesis, that 
   \begin{equation} F \bigl(U \cap I_{\ell}(r) \bigr)\subseteq \bigcup_{i=1}^{C_{N-1}} \widetilde{U}_i(\ell, r), \; \text{ for sets } \widetilde{U}_i(\ell, r) \in \Lambda(N-1, \lambda). \label{F(portion of U)} \end{equation} 
 Assuming this for the moment, let us complete the proof of \eqref{F(U) end of induction}. In fact, the relations \eqref{interval decomposition} and \eqref{F(portion of U)} jointly imply
   \begin{align*} 
   	F(U) \cap J_{\ell} &= F(U) \cap [a, b) = F(U) \cap \bigcup_{r=0}^{n_0-1} \widetilde{I}_\ell(r) = F(U) \cap \Bigl[ \bigcup_{r=0}^{n_0-1} F \bigl({I}_\ell(r) \bigr) \Bigr] \\ & = \bigcup_{r=0}^{n_0-1} \Bigl[ F(U) \cap F \bigl({I}_\ell(r) \bigr) \Bigr]
   	= \bigcup_{r=0}^{n_0-1} F\bigl(U \cap {I}_\ell(r) \bigr) \subseteq \bigcup_{r=0}^{n_0-1} \bigcup_{i=1}^{C_{N-1}} \widetilde{U}_i(\ell, r).  \end{align*}
   Since $F(U)$ is covered by the union of the disjoint intervals $J_{\ell}$ over non-negative indices $\ell = 0$ mod $n_0$, we conclude that 
   \begin{align} F(U) = \bigcup_{\ell }^{\ast}  \Bigl[ F(U) \cap J_{\ell} \Bigr] &\subseteq \bigcup_{\ell }^{\ast}\bigcup_{r=0}^{n_0-1}  \bigcup_{i=1}^{C_{N-1}} \widetilde{U}_i(\ell, r) \nonumber \\ 
   	&=  \bigcup_{r=0}^{n_0-1}  \bigcup_{i=1}^{C_{N-1}} U^{\ast}_{ri}, \quad U^{\ast}_{ri} := \bigcup_{\ell} \widetilde{U}_i(\ell, r). \label{def: Uri} \end{align} 
   In the above display $\bigcup^{\ast}$ denotes a union ranging over positive integers $\ell = 0$ mod $n_0$. It follows from the definition \eqref{def: Uri} that for each $0 \leq r < n_0$ and $1 \leq i \leq C_{N-1}$, the set $U^{\ast}_{ri}$ is a disjoint union of sets, whose $\ell^{\text{th}}$ constituent $\widetilde{U_i}(\ell, r)$ lies in the $\ell^{\text{th}}$ gap interval $J_{\ell}$ of $B_0$. The property $\widetilde{U}_i(\ell, r) \in \Lambda(N-1, \lambda)$ from \eqref{F(portion of U)} leads to the conclusion that 
   \[U^{\ast}_{ri} \in \Lambda(N, \lambda)  \text{ with special sequence $B_0$}. \] In summary, $F(U)$ lies in the union of at most $n_0 C_{N-1}$ sets in $\Lambda(N, \lambda)$, completing the proof of \eqref{F(U) end of induction}.  A visual depiction of the gap intervals $J_{\ell}$ and their decomposition is given in Figure \ref{fig: bi-Lipschitz image}.   
		\vskip0.1in
		\noindent It therefore remains to show \eqref{F(portion of U)}.  
		Since $U\in \Lambda(N, \lambda)$ by assumption, it follows from Definition \ref{defn: Lacunary sets} that the portion of $U$ in any gap interval of $A$ must be of lower lacunary order. Since $I_{\ell}(r)$ is a gap interval of $A$ for every $\ell \geq 1$ and $0 \leq r < n_0$, we conclude that 
		\[ U \cap I_{\ell}(r) \in \Lambda(N-1, \lambda). \] The induction hypothesis applied to this portion of $U$ yields that $F(U \cap I_\ell(r))$ is contained in the union of at most $C_{N-1}$ lacunary sets of order $(N-1)$, which is the desired conclusion \eqref{F(portion of U)}. 
		\end{proof}
\begin{center} 
\begin{figure} 
\begin{tikzpicture}[scale=0.8, line cap=round, line join=round]

\usetikzlibrary{decorations.pathreplacing,patterns}

% Colors
\definecolor{myred}{RGB}{210,0,25}
\definecolor{myblue}{RGB}{0,100,180}

% Main line
\draw[black, very thick] (0,0) -- (12,0);

% Left dots
\foreach \x in {-1.0,-0.75,-0.5}
  \fill[black] (\x,0) circle (1pt);

% Main black points
\fill[black] (0,0) circle (3.5pt);
\fill[black] (2.5,0) circle (3.5pt);
\fill[black] (6.0,0) circle (3.5pt);
\fill[black] (12.0,0) circle (3.5pt);

% Smaller colored division points
\fill[myblue] (0.8,0) circle (1.6pt);
\fill[myblue] (4.0,0) circle (1.6pt);
\fill[myblue] (8.1,0) circle (1.6pt);

% Alternating blue/red small ticks
\foreach \x in {0.05,0.25,0.45,0.65}
  \draw[myblue, thick] (\x,-0.07) -- (\x,0.07);

\foreach \x in {0.95,1.15,1.35,1.55,1.75,1.95,2.15,2.35}
  \draw[myred, thick] (\x,-0.07) -- (\x,0.07);

\foreach \x in {2.85,3.1,3.35,3.6,3.85}
  \draw[myblue, thick] (\x,-0.07) -- (\x,0.07);

\foreach \x in {4.25,4.5,4.75,5.0,5.25,5.5,5.75}
  \draw[myred, thick] (\x,-0.07) -- (\x,0.07);

\foreach \x in {6.25,6.5,6.75,7.0,7.25,7.5,7.75,8.0}
  \draw[myblue, thick] (\x,-0.07) -- (\x,0.07);

\foreach \x in {8.35,8.6,8.85,9.1,9.35,9.6,9.85,10.1,10.35,10.6,10.85,11.1,11.35,11.6,11.85}
  \draw[myred, thick] (\x,-0.07) -- (\x,0.07);

% Top braces and labels
\draw[myblue, decorate, decoration={brace, amplitude=5pt}]
  (0.0,0.35) -- (0.85,0.35);
\node[myblue, above] at (0.43,0.55) {\footnotesize{$\widetilde{I}_{2n_0}(0)$}};

\draw[myred, decorate, decoration={brace, amplitude=5pt}]
  (0.95,0.35) -- (2.35,0.35);
\node[myred, above] at (1.65,0.55) {\footnotesize{$\widetilde{I}_{2n_0}(1)$}};

\draw[myblue, decorate, decoration={brace, amplitude=5pt}]
  (2.75,0.35) -- (4.0,0.35);
\node[myblue, above] at (3.35,0.55) {\footnotesize{$\widetilde{I}_{n_0}(0)$}};

\draw[myred, decorate, decoration={brace, amplitude=5pt}]
  (4.2,0.35) -- (5.95,0.35);
\node[myred, above] at (5.05,0.55) {\footnotesize{$\widetilde{I}_{n_0}(1)$}};

\draw[myblue, decorate, decoration={brace, amplitude=5pt}]
  (6.2,0.35) -- (8.1,0.35);
\node[myblue, above] at (7.15,0.55) {\footnotesize{$\widetilde{I}_{0}(0)$}};

\draw[myred, decorate, decoration={brace, amplitude=5pt}]
  (8.35,0.35) -- (11.95,0.35);
\node[myred, above] at (10.15,0.55) {\footnotesize{$\widetilde{I}_{0}(1)$}};

% Bottom interval braces
\draw[black, decorate, decoration={brace, amplitude=7pt, mirror}]
  (0.0,-0.7) -- (2.5,-0.7);
\node[below] at (1.25,-0.95) {$J_{2n_0}$};

\draw[black, decorate, decoration={brace, amplitude=7pt, mirror}]
  (2.5,-0.7) -- (6.0,-0.7);
\node[below] at (4.25,-0.95) {$J_{n_0}$};

\draw[black, decorate, decoration={brace, amplitude=7pt, mirror}]
  (6.0,-0.7) -- (12.0,-0.7);
\node[below] at (9.0,-0.95) {$J_0$};

% Bottom labels
\node[below] at (-1.5,-0.12) {\small{$F(\alpha)$}};
\node[below] at (0.0,-0.12) {\small{$F(a_{3n_0})$}};
\node[below] at (2.7,-0.12) {\small{$F(a_{2n_0})$}};
\node[below] at (6.25,-0.12) {\small{$F(a_{n_0})$}};
\node[below] at (12.1,-0.12) {\small{$F(a_0)$}};

\end{tikzpicture}
\caption{\small{An illustration of the proof of Proposition \ref{prop: AdFinLac under F} with $F$ decreasing and $n_0=2$. The black dots depict points of the special sequence $B_0 = F(A_0) \in {\tt{MonLac}}(\lambda)$. Each gap interval $J_{\ell}$ of $B_0$ consists of two subintervals $\tilde{I}_{\ell}(0)$ and $\tilde{I}_{\ell}(1)$. The proof shows that $F(U) \cap \tilde{I}_{\ell}(r) \in \Lambda(N-1, \lambda; C_{N-1})$, and hence $\bigcup_{\ell}^{\ast} \widetilde{I}_{\ell}(r) \cap F(U) \in \Lambda(N, \lambda; C_{N-1})$. }} \label{fig: bi-Lipschitz image}
\end{figure}
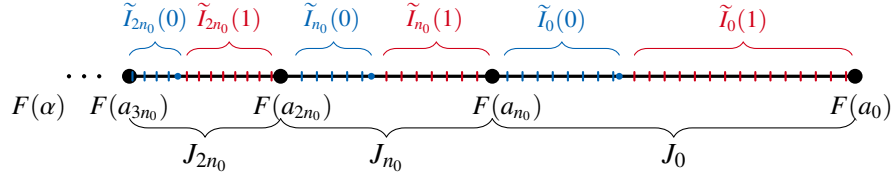
\end{center} 

\chapter{Appendix D: Lacunarity order and splitting number} \label{chapter: lac-split-appendix} 
In Proposition \ref{TREE-LACUNARY-TRADITIONAL} of Chapter \ref{section: split implies lacunarity}, we proved a key relation between the lacunarity order of a set and the splitting number of its tree in any integer base $M \geq 2$. Precisely, we showed that a set $U \subseteq [0,1)$ whose $M$-adic tree $\mathcal T(U; M)$ has finite splitting number must be admissible lacunary of finite order, in the sense of Definition \ref{defn: Admissible finite order lacunarity}. It is natural to ask whether the converse also holds, even though it is not directly used in any of the proofs in this monograph.  Establishing this reverse implication is the main objective of this chapter. 

\section{Converse of Proposition \ref{TREE-LACUNARY-TRADITIONAL}} \label{section: lacunarity implies split}
Let us recall a few definitions from previous chapters:
\vskip0.1in 
\begin{itemize} 
\item The class of sets $\Lambda(N, \lambda; R)$ from Definition \ref{defn: Admissible finite order lacunarity};
\vskip0.1in
\item The $M$-adic tree representation of a set from Section \ref{section: tree rep of sets}; 
\vskip0.1in
\item The splitting number of a tree from Section \ref{section: splitting number}. 
\end{itemize}  
\vskip0.1in
The main result of this chapter specifies the following connection between the three. 
\begin{proposition} \label{TREE-LACUNARY-TRADITIONAL-CONVERSE} 
For all integers $N, R \geq 1$, $M \geq 2$ and any $\lambda \in (0, 1)$, there exists a constant $C = C(N, \lambda, R; M) \geq 1$ with the following property. Suppose that  $U \subseteq [0,1)$ is a closed subset of $\mathbb R$. Then 
\begin{equation} \label{lacunarity implies finite split} 
U \in \Lambda(N, \lambda; R) \; \text{ must obey } \; {\text{split}}\bigl(\mathcal T(U; M) \bigr) \leq C(N, \lambda, R; M). 
\end{equation} 
In particular, this leads to the implication for every integer $M \geq 2$:
\[ U \text{ compact in } \mathbb R, \; U \in {\tt{AdFinLac}} \; \Longrightarrow \; \text{split}\bigl(\mathcal T(U; M)\bigr) < \infty. \]
This is the converse of \eqref{finite split implies finite order}.
\end{proposition} 
\noindent In accordance with the digit-expansion convention of Section \ref{section: rays as points},
all closed sets represented by $M$-adic trees will be regarded as subsets of $[0,1)$. This entails no loss of generality, since the endpoint $1$ is not represented by a maximal ray and therefore does not affect the associated tree. As explained in Section \ref{section: trees for bounded sets}, an arbitrary bounded set $V \subseteq [0, \infty)$ can be contained within a closed set $U$ of the type considered in Proposition \ref{TREE-LACUNARY-TRADITIONAL-CONVERSE} after an $M$-adic affine transformation, which leaves the splitting number of its tree unchanged. 
\subsection{Proof overview} 
The proof of Proposition \ref{TREE-LACUNARY-TRADITIONAL-CONVERSE} is organized around an analogy between finite-order lacunary sets and recursively structured trees. In Section \ref{section: Lacunary sets of order N}, a set in $\Lambda(N,\lambda)$ was assembled from lower-order lacunary subsets lying in the gaps of a special monotone lacunary sequence. Correspondingly, in Section \ref{section: TreeLac}, we will construct a tree of finite lacunarity order from lower-order subtrees rooted away from a distinguished ray, which serves as its structural spine. The main task is to make the correspondence between the set and the tree precise and then convert it into a bound for the splitting number.
\vskip0.1in
\noindent We begin in Section \ref{main ingredients section} with two facts about general rooted labelled trees. Lemma \ref{SPLITTING NUMBER PROPERTY} and Corollary \ref{Corollary: Finite unions of trees} control the splitting number of a finite union of trees in terms of the splitting numbers of its constituents. This reduces Proposition \ref{TREE-LACUNARY-TRADITIONAL-CONVERSE} from the class $\Lambda(N,\lambda;R)$ to $\Lambda(N,\lambda)$. Lemma \ref{LEMMA: SPLIT OFF DISTINGUISHED RAY}  gives a complementary (and much better) estimate for a tree that is assembled from a distinguished ray and a possibly infinite collection of subtrees rooted away from that ray.
\vskip0.1in
\noindent Motivated by the latter estimate, Section \ref{section: TreeLac} introduces the notion of {\em{finite-order tree lacunarity}}; see Definition \ref{defn: Tree Lacunary sets}. Lemmas  \ref{Lemma: Tree Lac Split} and \ref{lemma: NRM}, together with Corollary \ref{corollary: admissible treelac}, show that such a recursively constructed tree has controlled splitting number. This reduces the proof of Proposition \ref{TREE-LACUNARY-TRADITIONAL-CONVERSE} to the central intermediate statement, Proposition \ref{PROP: REDUCED CONVERSE}. Namely, the $M$-adic tree of every closed set in $\Lambda(N,\lambda)$ has uniformly bounded tree lacunarity order and is contained in the union of a uniformly bounded number of such trees.
\vskip0.1in
\noindent The proof of Proposition \ref{PROP: REDUCED CONVERSE} requires two further ingredients, stated in Section \ref{section: lacunary tree}. Proposition \ref{TREE STRUCTURE THEOREM} describes the structure of the $M$-adic tree associated with a rapidly converging lacunary monotone sequence; its consequence, Corollary \ref{COROLLARY: TREE LACUNARY SPLITTING NUMBER}, provides the base case $N=1$ of the induction-based proof of Proposition \ref{TREE-LACUNARY-TRADITIONAL-CONVERSE}. Lemma \ref{GAP LEMMA} identifies the tree vertices corresponding to the gaps of such a lacunary sequence. It thereby transfers the lower-order lacunary subsets occupying those gaps to controlled off-spine subtrees and supplies the inductive step.
\vskip0.1in
\noindent These ingredients are assembled in Section \ref{proof: converse}, where Proposition \ref{PROP: REDUCED CONVERSE} is proved by induction on $N$, and Proposition \ref{TREE-LACUNARY-TRADITIONAL-CONVERSE} follows from the resulting splitting-number bound together with Corollary \ref{Corollary: Finite unions of trees}. The general tree results from Section \ref{main ingredients section} are proved in Sections \ref{proof: splitting number of union} and \ref{section: split off distinguished ray}. The more specialized results stated in Sections \ref{section: TreeLac} and \ref{section: lacunary tree} are proved in Chapters \ref{Chapter: Appendix E} and \ref{section: proof of gap lemma}.
\section{Splitting number of tree unions} \label{main ingredients section}
\subsection{A finite union of trees} 
According to Definition \ref{defn: Admissible finite order lacunarity}, a set in $\Lambda(N, \lambda; R)$ is contained in the $R$-fold union of sets in $\Lambda(N, \lambda)$. It is therefore important to understand how splitting number behaves under finite unions. 
\begin{lemma}\label{SPLITTING NUMBER PROPERTY}
Given $L \geq 1$, let $\mathcal{T}^{[i]}$, $i=1, 2$ be two subtrees of a rooted labelled tree satisfying
\begin{equation} 
\text{split}\bigl(\mathcal T^{[i]} \bigr) \leq L, \quad \text{ for } i= 1, 2. 
\end{equation}
Suppose that $\mathcal T$ is the union of these two trees, i.e.~$\mathcal T$ is the tree whose collection of maximal rays is given by $\partial \mathcal T := \partial \mathcal{T}^{[1]} \cup \partial \mathcal T^{[2]}$. Then 
\begin{equation} \label{splitting number of two unions}  
\text{split}\bigl( \mathcal T \bigr) \leq 4(2L+3) \leq 20L. 
\end{equation} 
\end{lemma}
\vskip0.1in
\noindent Lemma \ref{SPLITTING NUMBER PROPERTY} is proved in Section \ref{proof: splitting number of union}. 
\vskip0.1in
\noindent Given a collection of $R$ trees of splitting number at most $L$, one may group the trees in batches of two, and iterate the relation \eqref{splitting number of two unions} $\lceil{\log_2 R} \rceil$ times to bound the splitting number of their union. The precise statement is recorded below. 
\begin{corollary} \label{Corollary: Finite unions of trees} 
For any choice of integers $L \geq 1$, $R \geq 2$, there exists a positive constant $C(L, R) \geq 1$, 
\[ C(L, R) \leq 20^{\lceil{\log_2 R} \rceil} L \leq R^{10}L, \] with the following property. 
\vskip0.1in
\noindent If $\{ \mathcal T^{[i]} : 1 \leq i \leq R \}$ is a finite collection of subtrees of a rooted, labelled tree, each with splitting number at most $L$, then 
\[ \text{split} \Bigl( \bigcup_{i=1}^{R} \mathcal T^{[i]}\Bigr) \leq  C(L, R). \] 
\end{corollary}
\vskip0.1in
\noindent {\em{Remark: }} By Lemma \ref{lemma: topological closure}, a closed set $U \subseteq [0, 1)$, $U \in \Lambda(N, \lambda; R)$ is contained in a bounded union of closed lacunary sets of the same order:  
\[ U \subseteq \bigcup_{i=1}^{RC_N} U_i, \quad U_i \in \Lambda(N, \lambda), \; U_i \text{ closed}. \] 
Intersecting with $U$ if needed, we may additionally assume $U_i \subseteq [0, 1)$. Consequently, 
\[ \mathcal T(U; M) \subseteq \bigcup_{i=1}^{R C_N} \mathcal T(U_i; M). \]
Corollary \ref{Corollary: Finite unions of trees} and the monotonicity of splitting numbers (Lemma \ref{monotone trees}) therefore reduce Proposition \ref{TREE-LACUNARY-TRADITIONAL-CONVERSE} to the case $R=1$. Henceforth, we focus on the following reduced statement. 
%The utility of Corollary \ref{Corollary: Finite unions of trees} is clear. A set in $\Lambda(N, \lambda; R)$ is, by definition, contained in the union of at most $R$ sets in $\Lambda(N, \lambda)$. Thus, by Corollary \ref{Corollary: Finite unions of trees} and the monotonicity of splitting numbers (Lemma \ref{monotone trees}), Proposition \ref{TREE-LACUNARY-TRADITIONAL-CONVERSE} reduces to the case $R=1$, so that $\Lambda(N, \lambda; R)$ in \eqref{lacunarity implies finite split} is replaced by $\Lambda(N, \lambda)$.   
\begin{lemma} \label{TREE-LACUNARY-TRADITIONAL-CONVERSE-reduction1} 
The conclusion \eqref{lacunarity implies finite split} of Proposition \ref{TREE-LACUNARY-TRADITIONAL-CONVERSE} follows from the existence of a constant $\Theta(N, \lambda; M) \geq 1$ such that
for all subsets $U \subseteq [0,1)$ that are closed in $\mathbb R$:  
\begin{equation} \label{R=1}
 U \in \Lambda(N, \lambda) \; \implies \; {\text{split}}\bigl(\mathcal T(U; M) \bigr) \leq \Theta(N, \lambda; M). 
\end{equation}  
\end{lemma} 
\subsection{A union of subtrees rooted away from a ray}  
Lemma \ref{SPLITTING NUMBER PROPERTY} and Corollary \ref{Corollary: Finite unions of trees} bound the splitting number of a finite, but otherwise arbitrary, union of trees in terms of the maximum splitting number of the constituents.  In the sequel, we will also need to control the splitting number of a possibly infinite union of subtrees that are in special arrangement within a larger tree. The following result gives an upper bound on the splitting number in this situation.  
\begin{lemma} \label{LEMMA: SPLIT OFF DISTINGUISHED RAY}
Let $\mathcal T$ be a tree rooted at $v_0$. Suppose there exist an integer $N \geq 0$ and a ray $\mathcal R^{\ast}$ rooted at $v_0$  such that 
\begin{equation} \label{split off special ray} 
{\text{split}}_{\mathcal T}(w) \leq N   \quad \text{ for every vertex } w \in \mathcal T, w \notin \mathcal R^{\ast}. 
\end{equation} 
 Then the following inequality holds:  
\begin{equation} \label{split of tree in terms of split off ray} 
{\text{split}}(\mathcal T) \leq N+1. 
\end{equation} 
Equality holds in \eqref{split of tree in terms of split off ray} if there are at least two vertices $w, w'$ with equality in \eqref{split off special ray}, $w \not\subseteq w'$, $w' \not\subseteq w$.  
\end{lemma} 
\vskip0.1in
\noindent {\em{Remarks: }} 
\begin{enumerate}[1.]
\item The hypothesis of Lemma \ref{LEMMA: SPLIT OFF DISTINGUISHED RAY} does not require $\mathcal R^{\ast}$ to be a ray in $\mathcal T$; in other words, not all vertices on $\mathcal R^{\ast}$ need to be in $\mathcal T$.   
\vskip0.1in
\item In view of the monotonicity of splitting numbers based on lineage as given by Lemma \ref{monotonicity}, the hypothesis \eqref{split off special ray} can be replaced by the equivalent, though superficially less restrictive, assumption 
\begin{equation} 
{\text{split}}_{\mathcal T}(w) \leq N \; \text{ for every } w \in \mathcal T \text{ with } w \notin \mathcal R^{\ast}, \text{ whose parent is in $\mathcal R^{\ast}$}. 
\end{equation}
\vskip0.1in
\item A ray $\mathcal R^{\ast}$ obeying the hypothesis of Lemma \ref{LEMMA: SPLIT OFF DISTINGUISHED RAY} serves as the ``spine'' of the tree $\mathcal T$, in the sense that all vertices away from $\mathcal R^{\ast}$ generate subtrees of $\mathcal T$ with lesser splitting number. This motivates the definition of a distinguished ray in the notion of tree lacunarity discussed in Section \ref{section: TreeLac}.   
\vskip0.1in 
\item Lemma \ref{LEMMA: SPLIT OFF DISTINGUISHED RAY} is proved in Section \ref{section: split off distinguished ray}.  
\end{enumerate} 

\section{Tree lacunarity of finite order} \label{section: TreeLac}
We start by developing some necessary terminology for an arbitrary $M$-adic tree 
\begin{align} \label{specify T} 
&\mathcal T \subseteq \mathcal T([0,1]; M) \text{ of infinite height rooted at a vertex $v_0$, and } \\ 
&\text{an infinite ray $\mathcal R^{\ast}$, not necessarily in $\partial \mathcal T$, also rooted at $v_0$. } \label{specify R-star} 
\end{align}  
The terminology we introduce mirrors the iterative construction of a finite-order lacunary set. Just as a monotone, lacunary sequence is the frame of reference for a finite order lacunary set, the distinguished ray $\mathcal R^{\ast}$ serves as the structural spine of a recursively constructed lacunary tree of finite order.
\vskip0.1in
\noindent For $\mathcal R^{\ast}$ as in \eqref{specify R-star}, let us define the collection of vertices located just off of $\mathcal R^{\ast}$:  
\begin{equation} 
\mathcal W^{\ast} = \mathcal W^{\ast}(\mathcal R^{\ast}) := \left\{w \in \mathcal T([0,1];M) \, \Bigl| \, \begin{aligned} &w \subsetneq v, v \in \mathcal R^{\ast}, w \notin \mathcal R^{\ast} \\ & h(w) = h(v) +1 \end{aligned} \right\}. \label{what is W-star}
\end{equation} 
In other words, $\mathcal W^{\ast}(\mathcal R^{\ast})$ consists of all vertices of the full $M$-adic tree $\mathcal T([0,1]; M)$ which are not on $\mathcal R^{\ast}$ but whose parents lie on $\mathcal R^{\ast}$.  
Let $\mathcal T_w([0,1]; M)$ denote the full $M$-adic subtree rooted at $w$. The collection of these subtrees, combined with the ray $\mathcal R^{\ast}$, partitions the full $M$-adic tree, with $\mathcal R^{\ast}$ acting as the ``spine'' for the decomposition: 
\begin{equation} \label{decomposition by spine} 
\mathcal T([0,1]; M) = \bigsqcup \bigl\{ \mathcal T_w([0,1]; M) : w \in \mathcal W^{\ast}(\mathcal R^{\ast}) \bigr\} \sqcup \{ \mathcal R^{\ast} \}. 
\end{equation} 
\noindent Among the vertices of $\mathcal W^{\ast}$, the ones that also belong to $\mathcal T$ comprise $\mathcal W(\mathcal T, \mathcal R^{\ast})$.   
\begin{equation}  
\mathcal W = \mathcal W(\mathcal T, \mathcal R^{\ast}) := \mathcal W^{\ast} \cap \mathcal T. \label{what is W(T)}
\end{equation} 
\noindent Let $\mathcal T_{w}$ denote the maximal subtree of $\mathcal T$ rooted at a vertex $w \in \mathcal W$. Intersecting both sides of \eqref{decomposition by spine} by $\mathcal T$, one sees that the subtrees $\mathcal T_w$ provide a decomposition for $\mathcal T$:
\begin{equation}  \label{tree decomp in w} 
\mathcal T \subseteq  \bigsqcup \bigl\{ \mathcal T_w : w \in \mathcal W(\mathcal T, \mathcal R^{\ast}) \bigr\} \sqcup \{ \mathcal R^{\ast}\}.  
\end{equation}  
As we will soon see, controlled splitting behaviour of subtrees of $\mathcal T$ rooted off of a single, specially chosen ray plays an important role in identifying the lacunarity order of the set encoded by $\mathcal T$.  
\vskip0.1in
\begin{definition} \label{defn: Tree Lacunary sets}
We declare the empty tree and every tree consisting of a single ray to be lacunary of order $0$. The collection of such trees is denoted 
\[ {\tt{TreeLac}}(0; M) \text{ for any integer } M \geq 2. \] 
\vskip0.1in
\noindent Suppose $N \geq 1$. A tree $\mathcal T \subseteq \mathcal T([0,1]; M)$  of the form \eqref{specify T} is said to be lacunary of order at most $N$ in base $M$, and denoted
\[ \mathcal T \in {\tt{TreeLac}}(N; M) \]  if there exists a ray $\mathcal R^{\ast}$ as in \eqref{specify R-star}, termed a distinguished ray of $\mathcal T$, such that 
\[ \mathcal T_w \in {\tt{TreeLac}}(N-1; M) \text{ for all } w \in \mathcal W(\mathcal T, \mathcal R^{\ast}).  \] 
\end{definition} 
\vskip0.1in
\noindent Since the splitting number of any tree of lacunarity order 0 is 0, a direct application of Lemma \ref{LEMMA: SPLIT OFF DISTINGUISHED RAY} yields
the following bound on the splitting number of ${\tt{TreeLac}}(N; M)$. We record this as a lemma. 
\begin{lemma}  \label{Lemma: Tree Lac Split}
Let $\mathcal T$ be as in \eqref{specify T}. Then 
\begin{equation} \label{splitting number for tree lac}
\text{split}(\mathcal T) \leq N \quad \text{ for any } \mathcal T \in {\tt{TreeLac}}(N; M).
\end{equation} 
\end{lemma} 
\begin{proof} 
The proof is by induction on $N$, with Lemma \ref{LEMMA: SPLIT OFF DISTINGUISHED RAY} at the inductive step. We leave this to the reader. 
\end{proof}
\vskip0.1in   
\noindent Incidentally, the proof of Lemma \ref{Lemma: Tree Lac Split} also shows that if a tree $\mathcal T \in {\tt{TreeLac}}(N; M)$ has splitting number $N$,  then its vertices of highest splitting number lie on $\mathcal R^{\ast}$. In this case, the special ray of $\mathcal T$ as identified in Lemma \ref{SplitsOnARay}, is a subset of the distinguished ray $\mathcal R^{\ast}$. Unlike a special ray (which could be finite),  a distinguished ray is always of infinite length and need not be unique.
\vskip0.1in
\noindent Finite unions of finite order lacunary sets $\Lambda(N, \lambda)$ led to the class $\Lambda(N, \lambda; R)$ in Definition \ref{defn: Admissible finite order lacunarity}. Carrying the analogy between finite order lacunary sets and trees in this direction leads to the following definition.  
\begin{definition} \label{defn: admissible tree lacunarity} 
Let $N \geq 0$ and $R \geq 1$ be integers. A tree $\mathcal T$ that is contained in the union of at most $R$ trees in  ${\tt{TreeLac}}(N; M)$ is denoted 
\[ \mathcal T \in {\tt{TreeLac}}(N, R; M). \] 
\end{definition} 
%\vskip1in
%\noindent An infinite ray $\mathcal R^{\ast}$ rooted at $v_0$ is said to be a {\em{distinguished ray}} for $\mathcal T$ if it contains all vertices of maximal splitting number in $\mathcal T$. In other words,   
%\[ \text{ all vertices $v \in \mathcal T$ with } \text{split}_{\mathcal T}(v) = \text{split}(\mathcal T) \text{ lie on $\mathcal R^{\ast}$}.  \]  
% It is not necessarily a member of $\partial \mathcal T$.   
%\vskip0.1in
%\noindent  For example, if $\mathcal T$ consists of finitely many rays corresponding to the set $\{2^{-j} : 1 \leq j \leq R\}$, then any ray $\mathcal R(e)$ that represents a point of the form 
%\[ e = \sum_{j=1}^{\infty} e_j 2^{-j}, \quad e_1 = \ldots = e_R = 0 \] 
%is distinguished for $\mathcal T$.
\noindent Let us record an elementary fact about ${\tt{TreeLac}}(N, R; M)$, whose proof follows a familiar argument for finite unions.   
\begin{lemma} \label{lemma: NRM} 
Let $N, R \geq 1$ be integers. Suppose that $\mathcal T$ is equipped with a distinguished ray $\mathcal R^{\ast}$ such that 
\begin{align} 
&\mathcal T_w \in {\tt{TreeLac}}(N-1, R; M) \text{ for all } w \in \mathcal W(\mathcal T, \mathcal R^{\ast}). \nonumber \\  
&\text{ Then } \; \mathcal T \in {\tt{TreeLac}}(N, R; M). \label{NRM-conclusion} 
\end{align}  
\end{lemma} 
\begin{proof}
For each $w \in \mathcal W(\mathcal T, \mathcal R^{\ast})$, let us decompose $\mathcal T_w$ as
\begin{align*} 
&\mathcal T_{w} \subseteq \bigcup_{i=1}^{R} \mathcal T_w^{[i]}, \;\text{ with } \mathcal T_w^{[i]} \in {\tt{TreeLac}}(N-1; M); \text{ this means } \mathcal T \subseteq \bigcup_{i=1}^{R} \mathcal T^{[i]}, \\ 
& \; \text{ where } \; \mathcal T^{[i]} := \bigcup \left\{\mathcal T_w^{[i]} \cup \{\mathcal R^{\ast} \} : w \in \mathcal W(\mathcal T, \mathcal R^{\ast}) \right\}.
\end{align*} 
It follows from Definition \ref{defn: Tree Lacunary sets} that $\mathcal T^{[i]} \in {\tt{TreeLac}}(N; M)$ for $1 \leq i \leq R$, establishing \eqref{NRM-conclusion}.
\end{proof} 
% \item Let $\mathcal N$ denote the collection of all non-splitting vertices of $\mathcal T$ all of whose descendants are also non-splitting, but whose parent is a splitting vertex. Namely, 
% \begin{equation} 
% \mathcal N := \left\{v \in  \mathcal T \; \Bigl| \;
% \begin{aligned} 
%  &v \text{ is a vertex of minimal height such that the maximal} \\  
%  & \text{subtree of $\mathcal T$ rooted at $v$ consists of a single ray} 
%  \end{aligned}\right\}. 
% \end{equation}       
% Clearly, no $v \in \mathcal N$ can lie on $\mathcal R^{\ast}$ and every $v \in \mathcal N$ contains a unique element of $\mathcal L$. Then, each $v \in \mathcal N$ satisfies 
% \begin{equation} 
%\# \Bigl\{I = [a, b] : I \cap v \neq \emptyset, \; a, b \in \mathfrak{L}, (a, b) \cap \mathfrak{L} = \emptyset \Bigr\} = 1.  
%\end{equation} 
%In other words, every $M$-adic interval representing $v \in \mathcal N$ intersects eaxctly one gap interval of $\mathfrak{L}$. 
\vskip0.1in
\noindent In light of Corollary \ref{Corollary: Finite unions of trees}, Lemma \ref{Lemma: Tree Lac Split} immediately yields the following consequence. 
\begin{corollary} \label{corollary: admissible treelac} 
Let $\mathcal T$ be as in \eqref{specify T}. Then for $N, R \geq 1$,  
\begin{equation} \label{splitting number for tree lac R}
\text{split}(\mathcal T) \leq R^{10}N \quad \text{ for any } \mathcal T \in {\tt{TreeLac}}(N, R; M).
\end{equation} 
\end{corollary} 

\subsection{The key reduction} The resemblance in the iterative definitions of a finite lacunary set $\Lambda(\cdot, \lambda)$ and a finite order lacunary tree ${\tt{TreeLac}}(\cdot; M)$ is not a coincidence. While the lacunarity order of the set  may not match the order of its corresponding tree, it will emerge from our analysis that one is finite if and only if the other one is. In fact, Proposition \ref{TREE-LACUNARY-TRADITIONAL-CONVERSE} is a by-product of the following structural statement. 
\begin{proposition} \label{PROP: REDUCED CONVERSE} 
    For $\lambda \in (0,1)$, integers $N \geq 1$, $M\geq 2$, there exist constants 
    \[\gamma_N = \gamma(N, \lambda; M) \in \mathbb N \; \text{ and } \; \rho_N = \rho(N, \lambda; M) \in \mathbb N\] such that the following implication holds for all closed sets $U$:
    %the $M$-adic tree of any set in $\Lambda(N, \lambda)$ is in ${\tt{TreeLac}}(\Theta_N; M)$: 
   \begin{align}  
 &U \subseteq [0,1), \; U \in \Lambda(N, \lambda) \; \implies \; \mathcal T\bigl(U; M\bigr) \in {\tt{TreeLac}}(\gamma_N, \rho_N; M). \label{Lac set & lac tree}\\   
 &{\text{As a result, by Corollary \ref{corollary: admissible treelac},  }} \; \text{split}\bigl( \mathcal T(U; M)\bigr) \leq  \rho_N^{10}\gamma_N . \label{splitting number of lac N}
    \end{align}
    \end{proposition}
\vskip0.1in
\noindent {\em{Remarks: }} 
\begin{enumerate}[1.]
\item Proposition \ref{PROP: REDUCED CONVERSE} establishes the conclusion \eqref{R=1}, with 
\[ \Theta(N, \lambda; M) = \rho_N^{10} \gamma_N, \] and thus proves Proposition \ref{TREE-LACUNARY-TRADITIONAL-CONVERSE}, via Lemma \ref{TREE-LACUNARY-TRADITIONAL-CONVERSE-reduction1}. 
\vskip0.1in
\item In Lemma \ref{lemma: topological closure}, we have shown that the topological closure of every set in $\Lambda(N, \lambda)$ can be covered by at most $C_N$ closed sets in the same family. Thus, for fixed parameters $\lambda, N$ and $M$, the implication \eqref{Lac set & lac tree} is equivalent to an analogous statement involving finite unions: 
\vskip0.1in
\noindent For all integers $R \geq 1$, there is a constant \[ \varrho_{N,R} := \varrho(N, R, \lambda; M) = R C_N \rho_N \] such that  
for all $U \subseteq [0,1)$ closed in $\mathbb R$, 
\begin{equation} \label{Lac set & lac tree with R}
 U \in \Lambda(N, \lambda; R) \; \implies \; \mathcal T\bigl(U; M\bigr) \in {\tt{TreeLac}}(\gamma_N, \varrho_{N, R}; M).
\end{equation} 
We will use the two conclusions \eqref{Lac set & lac tree} and \eqref{Lac set & lac tree with R} interchangeably.  
\vskip0.1in
\item Proposition \ref{PROP: REDUCED CONVERSE} is proved in Section \ref{proof: converse}. 
\end{enumerate} 
\section{Tree of a lacunary sequence} \label{section: lacunary tree}
As we have seen in Definition \ref{defn: Lacunary sets}, sequences lying in ${\tt{MonLac}}(\lambda)$ play a foundational role in the construction of finite-order lacunary sets. Arbitrary sets in $\Lambda(1, \lambda)$ are contained in the union of at most two such sequences. Additionally, sequences in ${\tt{MonLac}}(\lambda)$ constitute the structural backbone of sets in $\Lambda(N, \lambda)$ for $N \geq 2$. Not surprisingly, the structure and the splitting number of a tree representing a sequence in ${\tt{MonLac}}(\lambda)$ are central to our analysis. 
The next proposition summarizes the finer properties of such a tree that are relevant to the proof of Proposition \ref{TREE-LACUNARY-TRADITIONAL-CONVERSE}. 
\begin{proposition}
\label{TREE STRUCTURE THEOREM}
Fix an integer base $M \geq 2$ and any large absolute constant $C_0 \geq 3$. Then the $M$-adic tree 
\[ \mathcal T(\mathfrak A) := \mathcal T(\mathfrak A; M) \] 
for any infinite lacunary sequence $\mathfrak A \subseteq [0, 1)$, with 
\begin{equation} \label{rapidly converging sequence} 
\mathfrak A = \{a_j : j \geq 1\} \in {\tt{MonLac}}(M^{-C_0}), \quad \alpha = \lim_{j \rightarrow \infty} a_j  < 1
\end{equation} 
satisfies the following properties: 
\vskip0.1in
\begin{enumerate}[(a)]
\item (Lacunarity order at most 2) \label{TREE STRUCTURE THEOREM-lacunarity order} The tree $\mathcal T(\mathfrak A)$ satisfies 
\[ \mathcal T(\mathfrak A) \in {\tt{TreeLac}}(2; M), \text{ in the sense of Definition \ref{defn: Tree Lacunary sets}.} \]
Consequently, by Lemma \ref{Lemma: Tree Lac Split}, 
\begin{equation} \label{split 2} 
\text{split} \bigl( \mathcal T(\mathfrak A)\bigr) \leq 2 \text{ for any $\mathfrak A \in {\tt{MonLac}}(M^{-C_0})$.}  
 \end{equation}  
 \vskip0.1in
\item (A distinguished ray) \label{TREE STRUCTURE THEOREM-distinguished ray} The ray 
\begin{equation} \label{lacunary limit ray} 
\mathcal R^{\ast} = \mathcal R(\alpha) \; \text{ representing } \; \alpha = \lim \mathfrak A 
\end{equation}  belongs to $\partial \mathcal T(\mathfrak A)$; it contains infinitely many splitting vertices of $\mathcal T(\mathfrak A)$ and is a distinguished ray for $\mathcal T(\mathfrak A)$. \label{TREE STRUCTURE THEOREM - distinguished ray}
\vskip0.1in
%\item (Children of vertices of $\mathcal R^{\ast}$) Each vertex on $\mathcal R^{\ast}$ has at most two children in $\mathcal T(\mathfrak A)$; in other words, every splitting vertex of $\mathcal T(\mathfrak A)$ lying on $\mathcal R^{\ast}$  contributes exactly one child to the collection $\mathcal W = \mathcal W(\mathcal T(\mathfrak A), \mathcal R^{\ast})$ given by \eqref{what is W(T)}.  \label{part: number of children} 
\vskip0.1in
\item (Maximal subtrees rooted away from the distinguished ray) Let 
\[ \mathcal W= \mathcal W(\mathcal T(\mathfrak A), \mathcal R^{\ast}) \text{ be as in \eqref{what is W(T)}, with $\mathcal R^{\ast}$ given by \eqref{lacunary limit ray}.} \]  
If $\mathcal T_w = \mathcal T_w(\mathfrak A)$ is the maximal subtree of $\mathcal T(\mathfrak A)$ rooted at a vertex $w \in \mathcal  W$, then $\partial \mathcal T_w$ consists of finitely many rays. \label{tree structure theorem part c}
\end{enumerate}
\end{proposition} 
\vskip0.1in 
\noindent {\em{Remarks: }} 
\begin{enumerate}[1.]
\item It follows from Definition \ref{defn: Tree Lacunary sets} of tree lacunarity that 
\begin{equation}  
\mathcal T_w \in {\tt{TreeLac}}(1; M) \; \text{ for each $w \in \mathcal W = \mathcal W(\mathcal T(\mathfrak A), \mathcal R^{\ast})$}. 
\end{equation} 
As a result, by Lemma \ref{Lemma: Tree Lac Split}, one has 
\begin{equation} \sup_{w \in \mathcal W} \text{split}\bigl( \mathcal T_w \bigr) \leq 1. 
\end{equation}  
\vskip0.1in
\item Proposition \ref{TREE STRUCTURE THEOREM} is proved in Section \ref{Section: proof of TREE STRUCTURE THEOREM}.
\end{enumerate}

\subsection{Splitting number of a set in $\Lambda(1, \lambda)$} 
%The second main ingredient of the proof is the computation of the splitting number of a tree representing a sufficiently sparse lacunary sequence. Since lacunary sequences are the building blocks of $\Lambda(N, \lambda)$, it is not surprising that such a statement would be key for an induction-based proof on the lacunarity order. 
%\begin{lemma}\label{LEMMA: SPLITTING NUMBER LACUNARY SEQUENCE}
%There is an absolute constant $C_0 > 0$ such that the $M$-adic tree representing any infinite sequence 
%\begin{equation} 
%A \subseteq [0,1], \quad A = \{a_j : j \geq 1 \} \in {\tt{MonLac}}(M^{-C_0}) 
%\end{equation} 
%has the property 
%\begin{equation} 
%{\text{split}}\bigl( \mathcal T(A; M)\bigr) = 1. 
%\end{equation} 
%\end{lemma} 
We have seen in Lemma \ref{Lemma: Lacunary 1} \eqref{Lemma: Lacunary 1 (c)} that any set $U \in \Lambda(1, \lambda)$ can be covered by at most two lacunary sequences in ${\tt{MonLac}}(\lambda)$. In turn, any sequence $A = \{a_j : j \geq 1 \} \in {\tt{MonLac}}(\lambda)$ contains $R_0$ disjoint infinite subsequences of the form 
\[ A_{\ell} := \{a_{\ell + jR_0} : j \geq 0 \} \in  {\tt{MonLac}}(M^{-C_0}) \text{ for } 1 \leq \ell \leq R_0, \]
where the integer $R_0 = R_0(\lambda, C_0; M)$ satisfies 
\begin{equation} 
\lambda^{R_0} < M^{-C_0}; \text{ this ensures the decomposition } \;  A = \bigsqcup_{\ell=1}^{R_0} A_{\ell}.  \label{RMC_0}
\end{equation}    
%That means 
%\[ \mathcal T(U; M) \subseteq \mathcal T(A_1; M) \cup \mathcal T(A_2; M).  \] 
In view of Corollary \ref{Corollary: Finite unions of trees}, Proposition \ref{TREE STRUCTURE THEOREM} and \eqref{split 2},  one is therefore immediately led to the following corollary. 
\begin{corollary} \label{COROLLARY: TREE LACUNARY SPLITTING NUMBER} 
For every $\lambda \in (0, 1)$, $C_0 \geq 3$ and integer $M \geq 2$, there exist constants $R_0$ and $\gamma_1 = \gamma(1, \lambda; M)\geq 1$ depending only on these quantities with the following properties.
\vskip0.1in
\begin{enumerate}[(a)]
\item Any set $A \subseteq [0,1)$, $A \in {\tt{MonLac}}(\lambda)$ admits a cover of the form
\[  A \subseteq \bigsqcup_{\ell=1}^{R_0} A_{\ell}, \quad \text{ with }  A_{\ell} \in {\tt{MonLac}}(M^{-C_0}). \] The maximum number of sets in the cover $R_0 = R_0(\lambda, C_0; M)$ is the constant in \eqref{RMC_0}. 
\vskip0.1in
\item Any closed set $U \subseteq [0,1)$, $U \in \Lambda(1, \lambda)$ is contained in the union of at most $2R_0$ sequences lying in ${\tt{MonLac}}(M^{-C_0})$. 
\vskip0.1in
\item As a result, any closed set $U \subseteq [0,1)$ obeys  
\begin{align} \label{TreeLac of a Lac seq} 
&\mathcal T(U; M) \in {\tt{TreeLac}}(2, 2R_0; M) \text{ provided $U \in \Lambda(1, \lambda)$; } \\
&\text{this in turn implies } \; \text{split}\bigl(\mathcal T(U; M) \bigr) \leq 2 (2R_0)^{10}. \nonumber 
\end{align} 
\end{enumerate} 
\end{corollary} 
\subsection{Gaps of a lacunary sequence encoded in its tree} 
In order to prove Proposition \ref{PROP: REDUCED CONVERSE}, and therefore Proposition \ref{TREE-LACUNARY-TRADITIONAL-CONVERSE}, it is not enough to just identify the tree structure of a lacunary sequence $\mathfrak A$ of the form \eqref{rapidly converging sequence}. One also needs to understand where the gaps of $\mathfrak A$ lie in relation to its tree, since these intervals carry the lower order components of lacunary sets of arbitrary finite order.   
\vskip0.1in
\noindent For the purposes of the result stated in this section, its application in Section \ref{proof: converse} and its proof in Chapter \ref{section: proof of gap lemma}, it will be important to distinguish an
$M$-adic interval from the canonical half-open cell associated with it. If a vertex $u\in\mathcal T([0,1];M)$ represents the closed $M$-adic interval
\begin{equation} \label{defn: canonical cell}
u=[a,b], \; \text{ we define its {\em{canonical cell}} } \;  u^\sharp:=
\begin{cases}
[a,b), & b<1,\\
[a,1], & b=1.
\end{cases}
\end{equation} 
Thus the cells $u^\sharp$ at any fixed height form a disjoint partition
of $[0,1]$. This convention agrees with the choice in
Section \ref{section: rays as points} of the finitary base-$M$ expansion for an $M$-adic rational: a boundary point is assigned to the interval immediately to its right. While $u$ and $u^{\sharp}$ are interchangeable for measure-theoretic estimates, the distinction between them emerges whenever membership or intersection with a set is involved.
\vskip0.1in
\noindent 
For $\mathfrak A = \{ a_j: j \geq 1\} \in {\tt{MonLac}}(M^{-C_0})$ and $\mathcal T = \mathcal T(\mathfrak A)$ as in Proposition \ref{TREE STRUCTURE THEOREM}, we define a collection of vertices of $\mathcal T([0,1]; M)$ as follows:
\begin{equation} \label{def-UA}
\mathcal U(\mathfrak A) := \left\{ u \in \mathcal T[0,1] \Biggl| \;  \begin{aligned} 
&u \text{ is of minimal height along its lineage} \\ 
&\text{such that } u^{\sharp} \cap [\inf \mathfrak A, \sup \mathfrak A] \neq \emptyset \\  
&\text{and } \text{split}_{\mathcal T}(u) = 0
\end{aligned} 
\right\}.  
\end{equation}
Equivalently, $\mathcal U(\mathfrak A)$ consists of all the  vertices $u$ of the $M$-adic tree that satisfy: 
\begin{equation} \label{U properties} 
u^{\sharp} \cap [\inf \mathfrak A, \sup \mathfrak A] \neq \emptyset, \quad \#(\mathfrak A \cap u^{\sharp}) \leq 1, \quad \#(\mathfrak A \cap v^{\sharp}) \geq  2, 
\end{equation}
where $v$ is the parent of $u$. We call $u \in \mathcal U(\mathfrak A)$ a {\em{stopping vertex}} of $\mathfrak A$, since it is the first vertex on its lineage for which the maximal subtree $\mathcal T_u(\mathfrak A)$ has splitting number zero, or equivalently, at most one maximal ray. Let us note that a stopping vertex of $\mathfrak A$ need not lie in $\mathcal T(\mathfrak A)$.  
\vskip0.1in
%\[
%\left\{
%\begin{aligned}
%&\text{A splitting vertex on $\mathcal R_w^{\ast}$ that has a splitting descendant } \\ 
%&\text{in $\mathcal T_w$ contributes $(M-1)$ children  to $\mathcal U(\mathfrak A)$. }\\
%&\text{The youngest splitting vertex on $\mathcal R_w^{\ast}$ contributes all $M$.}
%\end{aligned} 
%\right\}
%\] 
\noindent The significance of the collection $\mathcal U(\mathfrak A)$ is that the canonical cell for each of its vertices meets at most a bounded number of gaps of $\mathfrak A$; additionally, these canonical cells cover the interval spanned by $\mathfrak A$. 
\begin{lemma} \label{GAP LEMMA} 
The relations between the vertices of $\mathcal U(\mathfrak A)$ given by \eqref{def-UA} and the gaps of $\mathfrak A$ are summarized below.
\begin{enumerate}[(a)]
\item \label{gap lemma parta} Every vertex $u \in \mathcal U(\mathfrak A)$ either obeys 
\begin{equation} \label{u-sharp intersection singleton}
u \cap \bigl[ \inf \mathfrak A, \sup \mathfrak A \bigr] = \{ \alpha \},   
\end{equation}  
or 
%is contained in the union of at most two adjacent gaps of $\mathfrak A$. Specifically, for every $u \in \mathcal U(\mathfrak A)$, 
there is an index $j \geq 1$ such that
\begin{equation} \label{u in 2 gaps}
u^\sharp\cap[\inf\mathfrak A,\sup\mathfrak A]
\subseteq
\begin{cases}
[a_{j+1},a_{j-1}), & \text{if $\mathfrak A$ is decreasing},\\
[a_{j-1},a_{j+1}), & \text{if $\mathfrak A$ is increasing}.
\end{cases}
\end{equation} 
%u^{\sharp} \subseteq \begin{cases} 
%\left\{
%\begin{aligned} 
%&\bigl[a_{j+1}, a_{j-1} \bigr) &\text{ if } j \geq 2, \\ 
%&[a_2, 1] &\text{ if } j = 1, \end{aligned}  
%\right\}
%&\text{ for $\mathfrak A$  decreasing}, \\ 
%\left\{
%\begin{aligned} 
%&\bigl[a_{j-1}, a_{j+1} \bigr) & \text{ if } j \geq 2, \\ 
%&[0, a_2) &\text{ if } j = 1, \end{aligned}  
%\right\}
%&\text{ for $\mathfrak A$ increasing}. \end{cases}
%\end{equation}
Here and below, we use the convention $a_0 = 0$ if $\mathfrak A$ is increasing, and $a_0 = 1$ if $\mathfrak A$ is decreasing. %{\footnote{We leave it to the reader to verify that \eqref{u-sharp intersection singleton} occurs only when $\mathfrak A$ is increasing, with $\alpha = \sup \mathfrak A$ a dyadic rational.}}
\vskip0.1in
\item \label{gap lemma partb} Conversely, the canonical cells associated with the vertices of $\mathcal U(\mathfrak A)$, together with the limit point $\alpha$, cover the closed interval spanned by the extreme values of $\mathfrak A$: 
\begin{equation} \label{extreme cover} 
[\inf \mathfrak A, \sup \mathfrak A] \subseteq  \bigcup \bigl\{u^{\sharp} : u \in \mathcal U(\mathfrak A) \bigr\} \cup \{\alpha \}. 
\end{equation}  
\vskip0.1in
\item \label{gap lemma - partc} Let $U \subseteq [0,1)$ be a closed set with the property
\[ U \subseteq [\inf \mathfrak A, \sup \mathfrak A]. \]
Suppose there exist integers $J\geq 0, R \geq 1$ such that the tree $\mathcal T(U;M)$ representing $U$ satisfies the condition
\begin{equation} \label{lac order of lowest level}
\mathcal T_u(U; M) \in {\tt{TreeLac}}(J, R; M) \; \text{ for all } u \in \mathcal U(\mathfrak A).  
%\sup \left\{ \text{split} \bigl( \mathcal T_u(U;M) : u \in \mathcal U(\mathfrak A)\bigr)\right\} \leq J.
\end{equation}
Here $\mathcal T_u(U;M)$ is the maximal subtree of $\mathcal T(U;M)$ rooted at $u$. If $u \notin \mathcal T(U; M)$, we understand $\mathcal T_u(U; M)$ to be the empty tree. 
Then 
%$\mathcal R^{\ast} = \mathcal R(\alpha)$ is a distinguished ray of $\mathcal T(U;M)$,  
\begin{align} 
&\mathcal T(U; M) \in {\tt{TreeLac}}(J+2, R; M); \label{order J+2} \\ 
&\text{As a result, }  \text{split} \bigl( \mathcal T(U;M) \bigr) \leq R^{10}(J+2). \label{Split TU}
\end{align} 
\end{enumerate}
\end{lemma} 
\vskip0.1in
\noindent Lemma \ref{GAP LEMMA} is proved in Chapter \ref{section: proof of gap lemma}, building on the material developed earlier in Chapter \ref{Chapter: Appendix E}. 

\begin{figure}[ht]
\centering
\begin{tikzpicture}[x=1.25cm,y=1cm,
    every node/.style={font=\small}]

% Baseline
\draw[->] (0.3,0) -- (7.8,0);

% Adjacent half-open gaps
\draw[line width=4pt,blue!35] (1,0) -- (4,0);
\draw[line width=4pt,green!40!black] (4,0) -- (7,0);

% Endpoint conventions for the gaps
\fill[blue!70!black] (1,0) circle (2pt);
\draw[blue!70!black,fill=white] (4,0) circle (2pt);

\fill[green!40!black] (4,0) circle (2pt);
\draw[green!40!black,fill=white] (7,0) circle (2pt);

% Sequence points
\fill (1,0) circle (1.5pt);
\fill (4,0) circle (1.5pt);
\fill (7,0) circle (1.5pt);

\node[below=5pt] at (1,0) {$a_{j+1}$};
\node[below=5pt] at (4,0) {$a_j$};
\node[below=5pt] at (7,0) {$a_{j-1}$};

% Gap labels
\node[above=5pt,blue!70!black] at (2.5,0)
    {$[a_{j+1},a_j)$};
\node[above=5pt,green!40!black] at (5.5,0)
    {$[a_j,a_{j-1})$};

% Canonical cell
\draw[line width=1.4pt,red] (2,-0.65) -- (7,-0.65);
\fill[red] (2,-0.65) circle (2.2pt);
\draw[red,fill=white,line width=.8pt]
    (7,-0.65) circle (2.2pt);

\node[left=3pt,red] at (2,-0.65) {$\ell$};
\node[below=4pt,red] at (4.6,-0.65)
    {$u^\sharp=[\ell,a_{j-1})$};

% Closed M-adic vertex
\draw[line width=.8pt,dashed,purple]
    (2,-1.35) -- (7,-1.35);
\fill[purple] (2,-1.35) circle (2pt);
\fill[purple] (7,-1.35) circle (2pt);

\node[below=4pt,purple] at (4.6,-1.35)
    {$u=[\ell,a_{j-1}]$};

% Endpoint guide
\draw[dotted] (7,0.15) -- (7,-1.55);

\end{tikzpicture}
\caption{\small{A stopping vertex $u$ relative to two adjacent gaps of a decreasing
sequence $\mathfrak A$. The canonical cell $u^\sharp$ is contained in
$[a_{j+1},a_{j-1})$ and may meet both adjacent half-open gaps, while
the corresponding closed $M$-adic vertex $u$ may additionally contain
the right endpoint $a_{j-1}$.}}
\label{fig:stopping-cell-gaps}
\end{figure}

\section{Proof of Propositions \ref{PROP: REDUCED CONVERSE} and \ref{TREE-LACUNARY-TRADITIONAL-CONVERSE}} 
\label{proof: converse}
Assuming the key ingredients stated in Sections \ref{main ingredients section}, \ref{section: TreeLac} and \ref{section: lacunary tree}, the proof of  Proposition \ref{PROP: REDUCED CONVERSE} is completed in this section. Proposition \ref{TREE-LACUNARY-TRADITIONAL-CONVERSE} follows from Proposition \ref{PROP: REDUCED CONVERSE}.  
%Lemma \ref{SPLITTING NUMBER PROPERTY}, Proposition \ref{TREE STRUCTURE THEOREM} and Corollary \ref{COROLLARY: TREE LACUNARY SPLITTING NUMBER}, 
\vskip0.1in
\noindent

\begin{proof}
    %By Lemma \ref{SPLITTING NUMBER PROPERTY}, we may assume that $\Omega\in\Lambda(N;\lambda)$ where $N=N_2$. 
    The primary objective is to prove \eqref{Lac set & lac tree} or its equivalent version \eqref{Lac set & lac tree with R}, since the other two claims \eqref{splitting number of lac N} and \eqref{lacunarity implies finite split} are consequences of it. The proof is by induction on $N$. 
\vskip0.1in
\noindent Applying Corollary \ref{COROLLARY: TREE LACUNARY SPLITTING NUMBER} with $C_0=10$, we obtain \eqref{Lac set & lac tree}  for the base case $N=1$ with $\gamma_1=2$ and $\rho_1 = 2R_0$. 
    \vskip0.1in
    \noindent Proceeding to the inductive step, let us fix a lacunarity order $N \geq 1$. By the induction hypothesis, we may assume that  for all sets $V \subseteq [0, 1)$ that are closed in $\mathbb R$, 
    \begin{align}  
&\mathcal T(V; M) \in {\tt{TreeLac}}(\gamma_N, \rho_N; M) \text{ provided }  V\in \Lambda(N, \lambda), \nonumber \\ 
&\text{or equivalently, } \mathcal T(V; M) \in {\tt{TreeLac}}(\gamma_N, \varrho_{N,R}; M) \text{ for }  V\in \Lambda(N, \lambda; R). \label{induction hypothesis: lac and split}
%&{\text{split}}\bigl(\mathcal T(V; M) \bigr) \leq C(N, \lambda, R; M)\leq R^{10} C(N, \lambda; M), \quad \text{ for all } R \geq 1. 
\end{align} 
Here $\varrho_{N, R} = R C_N\rho_N$, with the constant $C_N$ as in Lemma \ref{lemma: topological closure}.  
    %The corresponding estimate for $\Lambda(1, \lambda; R)$ follows from Corollary \ref{Corollary: Finite unions of trees} and Lemma \ref{monotone trees}.  
\vskip0.1in
\noindent Let us now consider a closed set $U$ contained in $[0, 1)$, $U \in \Lambda(N+1, \lambda)$, whose special sequence is given by 
\[ A = \{ \tilde{a}_j : j \geq 1 \} \in {\tt{MonLac}}(\lambda); \quad \alpha = \lim A. \] 
Since $U$ is closed in $\mathbb R$, and $U \subseteq [0, 1)$, it is compact and hence $\sup U < 1$. After applying, if necessary, an $M$-adic affine transformation to $U$ and its special sequence $A$, we may assume that the closure of $A$ is contained in $[0, 1)$ as well. Indeed, such a transformation preserves membership in
$\Lambda(N+1,\lambda)$, by Lemma \ref{lacunarity under linear operations}, and produces an isomorphic
$M$-adic tree, by Lemma \ref{lemma: trees under affine maps}.  In particular, \[ \alpha = \lim A < 1. \] 
Setting $C_0 = 10$, let us first choose a subsequence 
%is contained in the a bounded finite union of monotone tree-lacunary sequences 
\[ \mathfrak A \subseteq A, \quad \mathfrak A \in {\tt{MonLac}}(M^{-10}), \; \text{ with } \; \inf A = \inf \mathfrak A, \; \sup A = \sup \mathfrak A. \] 
This is possible by choosing 
\[ \mathfrak A = \{a_j = \tilde{a}_{1 + (j-1)R_0} : j \geq 1  \}, \text{ where } \lambda^{R_0} < M^{-10}.  \] 
Then, on one hand, $\mathfrak A$ obeys the hypotheses of Proposition \ref{TREE STRUCTURE THEOREM}. On the other, every gap interval of $\mathfrak A$ contains at most $R_0$ gaps of $A$. 
\vskip0.1in
\noindent The defining property of a special sequence implies, on one hand, that 
\begin{equation} \label{limits of U} 
U \subseteq [\inf \mathfrak A, \sup \mathfrak A]. 
\end{equation}  
On the other, the definition of $U \in \Lambda(N+1, \lambda)$ says that the intersection of $U$ with any gap interval of $A$ lies in $\Lambda(N, \lambda)$, which translates to the statement 
\begin{equation} \label{U in gap}
U \cap [a, b) \in \Lambda(N, \lambda; R_0) \text{ for any two consecutive elements $a, b$ of $\mathfrak A$}.
\end{equation} 
\vskip0.1in
\noindent 
%The $M$-adic tree $\mathcal T(\mathfrak A) = \mathcal T(\mathfrak A; M)$ obeys the properties stated in Proposition \ref{TREE STRUCTURE THEOREM}. In particular, this means that the ray $\mathcal R^{\ast}$ depicting $\alpha = \lim A$ is a distinguished ray for $\mathcal T(\mathfrak A)$.  Further, the containment \eqref{where is U} and the covering property \eqref{extreme cover} jointly imply that
%\begin{equation} \label{U containment} 
%U = U \cap [\inf \mathfrak A, \sup \mathfrak A] 
%\subseteq \bigcup \bigl\{u \cap U : u \in \mathcal U(\mathfrak A) \bigr\} \cup \{\alpha\} \\ 
%%&\subseteq \bigcup_{i=0}^{1} \bigcup_{w \in \mathcal W(\mathcal T, \mathcal R^{\ast})} \bigl\{u \cap U : u \in \mathcal U_w(\mathfrak A), \text{split}(\mathcal T_w) = i \bigr\} \cup \{\alpha\}. 
%\end{equation}
In view of the containment \eqref{limits of U} and  Lemma \ref{GAP LEMMA}\eqref{gap lemma parta}, we are led to one of the following scenarios: either 
\[ U \cap u \subseteq \{\alpha\}, \]
or $U \cap u$ is contained in the union of at most two adjacent half-open gaps of $\mathfrak A$, together with possibly the right endpoint of $u$. Combining this with \eqref{U in gap}, we obtain
\begin{equation} \label{u cap U}
U \cap u \in \Lambda(N, \lambda; 2R_0+1). 
\end{equation}  
Every ray of $\mathcal T_u(U;M)$ represents a point of
$U\cap u^\sharp\subseteq U\cap u$. Hence
$\mathcal T_u(U;M)$ is a subtree of $\mathcal T(U\cap u;M)$. 
\vskip0.1in
\noindent One may therefore apply the induction hypothesis \eqref{induction hypothesis: lac and split}, with $R = (2R_0+1)$, to the closed sets $U \cap u$ contained in $[0,1)$ appearing in \eqref{u cap U}.  
%Combined with Corollary \ref{Corollary: Finite unions of trees} and \eqref{splitting number for tree lac}, 
Since membership in ${\tt{TreeLac}}(J, R; M)$ is inherited by rooted subtrees, as follows directly from Definitions \ref{defn: Tree Lacunary sets} and \ref{defn: admissible tree lacunarity}, we obtain for every $u \in \mathcal U(\mathfrak A)$
\begin{equation}
\mathcal T_u\bigl( U ; M)  \in {\tt{TreeLac}}\bigl(\gamma_N, (2R_0+1) C_N \rho_N; M \bigr).
%&\text{split} \left(\mathcal T_u \bigl( U ; M) \right) = \text{split} \left(\mathcal T \bigl( U \cap u; M) \right) \leq (2R_0)^{10} \Theta_N. 
\end{equation} 
Applying Lemma \ref{GAP LEMMA} \eqref{gap lemma - partc} with $J = \gamma_N$ and $R = (2R_0+1) C_N \rho_N$ therefore yields
\[ 
\mathcal T\bigl( U ; M) \in {\tt{TreeLac}}(\gamma_{N+1}, \rho_{N+1}; M) \] 
where $\gamma_{N+1} =  \gamma_N +2$, $\rho_{N+1} = (2R_0+1) C_N \rho_N$.  
This completes the induction step, and establishes the desired conclusion \eqref{Lac set & lac tree} for all $N$. 
\vskip0.1in
\noindent The conclusion \eqref{Split TU} of Lemma \ref{GAP LEMMA} then implies 
\begin{equation} 
\text{split} \left(\mathcal T\bigl( U ; M) \right) \leq \rho_{N+1}^{10}\gamma_{N+1}.  
\end{equation} 
This establishes \eqref{R=1} with
\begin{equation*}
\Theta(N+1, \lambda; M) = \rho_{N+1}^{10}\gamma_{N+1}.
\end{equation*}
Lemma \ref{TREE-LACUNARY-TRADITIONAL-CONVERSE-reduction1}, together with Corollary \ref{Corollary: Finite unions of trees}, now yields \eqref{lacunarity implies finite split}; more precisely, one can take 
\begin{equation*}
C(N+1, \lambda, R; M) \leq R^{10} \Theta(N+1, \lambda; M). 
\end{equation*} 
This completes the proof of Propositions \ref{PROP: REDUCED CONVERSE} and \ref{TREE-LACUNARY-TRADITIONAL-CONVERSE}. 
\end{proof} 

\section{Splitting number of tree unions: Proof of Lemma \ref{SPLITTING NUMBER PROPERTY}} 
\label{proof: splitting number of union} 
\begin{proof}
Towards a contradiction, suppose if possible that 
\begin{equation} 
\text{split}(\mathcal T) > 4(2L+3) = 4(2L+2 + 1). 
\end{equation} 
Applying the pruning mechanism Proposition \ref{PRUNING STAGE 1} on $\mathcal T$ with $N =2L+2$ and $C_0=1$,  we obtain a model $(2L+2, 1)$ subtree $\mathcal U \subseteq \mathcal T$, in the sense of Section \ref{section: model trees}. For $j \geq 1$, let us define 
\[ \mathscr{V}_j := \bigl\{v : v \text{ is a $j^{\text{th}}$ splitting vertex of some maximal ray } \mathcal R \in \partial \mathcal U  \bigr\}. \] 
It follows from the definition of a model tree that every maximal ray of $\mathcal U$ contains exactly $(2L+2)$ splitting vertices. 
\vskip0.1in
\noindent Let $\mathcal{U}_j$ denote the subtree of $\mathcal{U}$, possibly with maximal rays of variable length, whose terminal vertices lie in $\mathscr{V}_j$. The model structure of $\mathcal U$ then implies that
every ray in $\partial \mathcal U_j$ has exactly $(j-1)$ splitting vertices, so that  
\begin{equation} 
\text{split}(\mathcal{U}_{L+2}) = L+1.  
\end{equation}   
\vskip0.1in
\noindent While every vertex of $\mathcal U$ belongs to either $\mathcal T^{[1]}$ or $\mathcal T^{[2]}$, we claim that:
\begin{equation} \label{distributed vertices} 
\mathscr{V}_{L+2} \not\subseteq \text{ the vertex set of }\mathcal T^{[i]} \text{ for either } i = 1 \text{ or } 2. 
\end{equation} 
Indeed, if this were the case, that would mean $\mathcal U_{L+2} \subseteq \mathcal T^{[i]}$ for some $i$. By the monotonicity of splitting numbers, one would obtain 
\[ L+1 = \text{split}\bigl(\mathcal{U}_{L+2} \bigr) \leq \text{split}\bigl(\mathcal T^{[i]} \bigr) \leq L, \text{ a contradiction}.  \] 
\vskip0.1in
\noindent The claim \eqref{distributed vertices} ensures a vertex $v \in \mathscr{V}_{L+2}$ such that 
\[ v \in \mathcal T^{[i]} \text{ but } v \notin \mathcal T^{[i']}; \; i \neq i'. \] Without loss of generality, let $i=1$, $i'=2$. Since $\mathcal T^{[2]}$ is a rooted subtree and $v \not\in \mathcal T^{[2]}$, the subtree $\widetilde{\mathcal U}(v) \subseteq \mathcal U$ rooted at $v$ does not contain any vertex of $\mathcal T^{[2]}$. Therefore, 
\[ \widetilde{\mathcal U}(v)  \subseteq \mathcal T^{[1]}. \]
Every ray of $\partial \widetilde{\mathcal U}(v)$ is a sub-ray of $\mathcal U$ originating from $v$, and therefore has exactly $(2L+2) - (L+1)$ splitting vertices. This means 
\[ L+1 = \text{split}\bigl( \widetilde{\mathcal U}(v)\bigr) \leq \text{split}\bigl(\mathcal T^{[1]} \bigr) \leq L, \]
again a contradiction. This completes the proof.    
%The model structure of $\mathcal U$ ensures each maximal ray of $\mathcal U$ has $2L+2$ splitting vertices. If the ray contains $v$, the sub-ray of $\mathcal U$ rooted at $v$ contains $(2L+2)-(L+1)= (L+1)$ splitting vertices of $\mathcal U$.  
\end{proof}

\section{Subtrees rooted off a ray: Proof of Lemma \ref{LEMMA: SPLIT OFF DISTINGUISHED RAY}} \label{section: split off distinguished ray} 
\begin{proof} 
In view of Corollary \ref{corollary-monotonicity} and the definition \eqref{splitting vertex} of the splitting number, the desired conclusion \eqref{split of tree in terms of split off ray} follows from the estimate
\begin{equation} \label{proof via contrapositive} 
\min_{\mathcal R \in \partial \mathcal S}\text{split}_{\mathcal S}(\mathcal R) \leq N+1
\end{equation}  
for every subtree $\mathcal S \subseteq \mathcal T$ rooted at $v_0$. We will prove \eqref{proof via contrapositive} via the contrapositive statement, considering two cases. 
\vskip0.1in 
\noindent Suppose there exists a subtree $\bar{\mathcal S} \subseteq \mathcal T$ rooted at $v_0$ for which \eqref{proof via contrapositive} fails, i.e.
\begin{equation} \label{will contradict this} 
\min_{\bar{\mathcal R} \in \partial \bar{\mathcal S}}\text{split}_{\bar{\mathcal S}}(\bar{\mathcal R}) \geq N+2. 
\end{equation}  Let $\bar{u}$ denote the (unique) splitting vertex of $\bar{\mathcal S}$ of minimal height. In other words, $\bar{\mathcal S}$ does not have any splitting vertex at height $h$, $0 \leq h < h(\bar{u})$.  
\vskip0.1in
\noindent {\em{Case 1: }} Suppose that $\bar{u} \notin \mathcal R^{\ast}$. 
%This means that $\bar{u}$ is either equal to, or is descended from, a vertex $w$ that is not on $\mathcal R^{\ast}$. 
Let $\bar{\mathcal S}_{\bar{u}}$ denote the maximal subtree of $\bar{\mathcal S}$ rooted at $\bar{u}$. Then there is a clear bijection between $\partial \bar{\mathcal S}_{\bar{u}}$ and $\partial \bar{\mathcal S}$. Namely, every ray $\bar{\mathcal R}_{\bar{u}} \in \partial \bar{\mathcal S}_{\bar{u}}$ is the restriction of a unique ray $\bar{\mathcal R} \in \partial \bar{\mathcal S}$; indeed, every $\bar{\mathcal R}_{\bar{u}} \in \partial \bar{\mathcal S}_{\bar{u}}$ extends uniquely towards $v_0$ to a ray $\bar{\mathcal R} \in \partial \bar{\mathcal S}$. Since $\bar u$ is the splitting vertex of $\bar{\mathcal S}$ of minimal height, no additional splitting vertex occurs between $v_0$ and $\bar u$. In other words, 
\[
{\text{split}}_{\bar{\mathcal S}_{\bar{u}}}(\mathcal R_{\bar{u}}) = {\text{split}}_{\bar{\mathcal S}}(\bar{\mathcal R}) \text{ for all rays } \mathcal R_{\bar{u}} \in \partial \bar{\mathcal S}_{\bar{u}}. \] 
The definition \eqref{splitting vertex} of the splitting number of a vertex then leads to the relation 
\[ \text{split}_{\mathcal T}({\bar{u}}) \geq \min_{\bar{\mathcal R}_{\bar{u}} \in \bar{\mathcal S}_{\bar{u}}} \text{split}_{\bar{\mathcal S}_{\bar{u}}}(\bar{\mathcal R}_{\bar{u}}) = \min_{\bar{\mathcal R} \in \partial \bar{\mathcal S}}\text{split}_{\bar{\mathcal S}}(\bar{\mathcal R}) \geq N+2. \]
This contradicts the assumption \eqref{split off special ray} for the vertex $w = \bar{u}$. 
\vskip0.1in
\noindent {\em{Case 2: }} Suppose next that the splitting vertex $\bar{u} \in \mathcal R^{\ast}$. 
%Our hypothesis dictates that $\bar{u}$ has exactly two children (in $\mathcal T$ and therefore in $\bar{\mathcal S}$), one of which lies in $\mathcal R^{\ast}$. 
Let $w \in \bar{\mathcal S}$ be a child of $\bar{u}$ not on $\mathcal R^{\ast}$. According to our assumption \eqref{will contradict this}, every ray $\bar{\mathcal R} \in \partial \bar{\mathcal S}$ contains at least $(N+2)$ splitting vertices, of which one is always $\bar{u}$. As before, let $\bar{\mathcal S}_w$ denote the maximal subtree of $\bar{\mathcal S}$ rooted at $w$. Then any ray $\bar{\mathcal R}_{w} \in \partial \bar{\mathcal S}_w$ is the restriction of a unique ray $\bar{\mathcal R} \in \partial \bar{\mathcal S}$. The ray $\bar{\mathcal R}_{w}$ contains all splitting vertices of $\bar{\mathcal R}$ except $\bar{u}$, and therefore splits at least $(N+1)$ times. In other words, 
\[ \text{split}_{\mathcal T}(w) \geq \min_{\bar{\mathcal R}_w \in \partial \bar{\mathcal S}_w} \text{split}_{\bar{\mathcal S}_w}(\bar{\mathcal R}_w)  \geq \min_{\bar{\mathcal R} \in \partial \bar{\mathcal S}}\Bigl[ \text{split}_{\bar{\mathcal S}}(\bar{\mathcal R}) -1 \bigr] \geq (N+1). \] 
This again contradicts  \eqref{split off special ray}, establishing \eqref{proof via contrapositive}, and hence \eqref{split of tree in terms of split off ray}.  
\vskip0.1in
\noindent Next, let us consider the equality assertion. Suppose there exist 
\[ w, w' \notin \mathcal R^{\ast}, \; w \subsetneq w', \; w' \subsetneq w, \; \text{ satisfying } \; 
%h(w) = h(v)+1, \; h(w') = h(v')+1,\\
%&v, v' \in \mathcal R^{\ast}, \;  h(v) < h(v'), \; 
{\text{split}}_{\mathcal T}(w) = {\text{split}}_{\mathcal T}(w') = N. 
\]
Let $v = D(w, w')$. Applying the second part of Lemma \ref{monotonicity} we obtain 
\[ \text{split}(\mathcal T) \geq \text{split}_{\mathcal T}(v) > \min \left\{ {\text{split}}_{\mathcal T}(w), \; {\text{split}}_{\mathcal T}(w') \right\} = N, \] 
i.e. $\text{split}(\mathcal T) \geq N+1$. Combining this with \eqref{split of tree in terms of split off ray} concludes the proof. 
\end{proof}

\chapter{Appendix E Part I: Tree of a lacunary sequence}\label{Chapter: Appendix E}
The purpose of this chapter is to analyse the $M$-adic tree associated with a rapidly converging lacunary sequence. We first relate $M$-adic digit expansions to Euclidean separation and record the resulting branching properties of finite ordered tuples. We then show that lacunarity forbids several tree configurations. These exclusions imply that subtrees rooted away from the limiting ray have lacunarity order at most one, leading to the proof of Proposition \ref{TREE STRUCTURE THEOREM}.  

\section{Digit expansions and Euclidean separation}
In Section \ref{tree encoding section}, we have seen how the digit sequence of a real number $x$ in base $M \geq 2$ carries precise geometric information about the point it represents. Given $x \in [0,1)$ with base-$M$ expansion 
\begin{equation} \label{expansion of x}  
x = \sum_{k=1}^{\infty} x_k M^{-k}, \quad \text{ with } x_k \in \mathbb Z_M \text{ for } k \geq 1, 
\end{equation}  
consider the corresponding ray $\mathcal R(x) \in \partial \mathcal T([0,1]; M)$  
\[ \mathcal R(x): \emptyset \mapsto \langle x_1 \rangle  \rightarrow \langle x_1, x_2 \rangle \rightarrow \cdots \rightarrow \langle x_1, x_2, \ldots, x_k \rangle \rightarrow \cdots. \] 
At height $k$, the ray follows the $x_k^{\text{th}}$ child of the vertex at height $(k-1)$. It therefore encodes $x$ through a nested sequence of closed $M$-adic intervals whose lengths tend to zero and whose intersection is $\{x\}$. Let us also recall that our numbering scheme for children follows the linear order on $\mathbb R$; the $0^{\text{th}}$ child is the leftmost subinterval of the parent, the $(M-1)^{\text{th}}$ child is the rightmost one.  
\vskip0.1in
\noindent As explained in Section \ref{section: rays as points}, the expansion \eqref{expansion of x} is unique as long as digit sequences with all but finitely many entries equal to $(M-1)$ are forbidden. We will continue to make this assumption for all $M$-adic expansions, so that the point-to-digit correspondence $x\mapsto \mathbf x = (x_1, x_2, \ldots)$ remains unambiguous. 
\vskip0.1in
\subsection{First separating height} 
Let us consider two points $x, y \in [0,1)$, $x \neq y$, with $M$-adic digit expansions 
\[ x = \sum_{k=1}^{\infty} x_k M^{-k}, \quad  y = \sum_{k=1}^{\infty} y_k M^{-k} \text{ with } x_k, y_k \in \mathbb Z_M \text{ for } k \geq 1.  \] 
The {\em{first separating height}} $\kappa(x, y)$ between $x$ and $y$ is the first position in the $M$-adic expansion, i.e., the first height of the $M$-adic tree, at which their digit sequences $\mathbf x = (x_1, x_2, \ldots)$ and $\mathbf y = (y_1, y_2, \ldots)$ differ. Namely,  
\begin{equation} \label{what is nu}
\kappa = \kappa(x, y) = \kappa(y, x) := \min \left\{k : x_k \neq y_k\right\}.
\end{equation}
Alternatively, it can also be described as 
 \[ \kappa(x, y) = \max \left\{k \geq 0 : x_{k'} = y_{k'} \text{ for all } 0 \leq k' \leq k \right\} + 1, \] with the convention $x_0 = y_0 = 0$ for $x, y \in [0, 1)$. The index $\kappa$ has a special significance in the context of the $M$-adic tree: the vertex $D(x, y)$, the youngest common ancestor of $x, y$ in $\mathcal T([0,1]; M)$, lies at height $(\kappa-1)$; in other words, 
\begin{equation}  \label{hD(x,y)} 
h(D(x, y)) = \kappa(x, y)-1 \quad \text{ or } \quad |D(x, y)| = M^{-\kappa+1}. 
\end{equation}   
It follows from the uniqueness of the $M$-adic expansions that two points $x$ and $y$ are distinct if and only if $\kappa(x, y)$ is  finite. In fact, the first separating height provides a generally non-optimal bound for the Euclidean distance between $x$ and $y$. We record this bound in the next lemma.  
\begin{lemma} \label{lemma: distance x y}
Let $x,y\in[0,1)$ with $x<y$, and set $\kappa=\kappa(x,y)$. Then 
\begin{equation}  x_{\kappa} < y_{\kappa}. \label{ineq at fsh} \end{equation} 
The converse is also true; namely \eqref{ineq at fsh} implies $x < y$.  Moreover, 
\begin{equation} \label{dist(x, y) upper}
y-x \leq M^{-\kappa + 1}. 
\end{equation} 
\end{lemma} 
\begin{proof} 
From the definition of $\kappa = \kappa(x, y)$, we have  
\[ x_k = y_k  \text{ for all } 1 \leq k < \kappa, \text{ and } x_{\kappa} \neq y_{\kappa}. \]  
Suppose $x_{\kappa} > y_{\kappa}$. Then 
\begin{align*} 
x - y &= (x_{\kappa} - y_{\kappa}) M^{-\kappa} + \sum_{k > \kappa} (x_k - y_k) M^{-k} \\
&\geq (x_{\kappa} - y_{\kappa}) M^{-\kappa} - \sum_{k > \kappa} (M - 1) M^{-k} \\
&\geq M^{-\kappa} - M^{-(\kappa + 1)} (M-1) \Bigl(1 - \frac{1}{M}\Bigr)^{-1} 
\geq  M^{- \kappa} - M^{-\kappa} = 0. 
\end{align*} 
This contradicts $x < y$. Since $x_{\kappa} \neq y_{\kappa}$, it follows that $x_{\kappa} < y_{\kappa}$.
\vskip0.1in 
\noindent The converse direction is essentially a retracing of steps. If $x_{\kappa} < y_{\kappa}$, then the same steps above lead to 
\begin{align*} y - x &= (y_{\kappa} - x_{\kappa}) M^{-\kappa} + \sum_{k > \kappa} (y_k - x_k) M^{-k} \\
&\geq M^{-\kappa} - M^{-\kappa} = 0. 
\end{align*}  
Equality could occur only if  
\[ y_{\kappa} - x_{\kappa} = 1, \quad  y_{k} = 0, \; x_k = M-1 \text{ for all $k > \kappa$}. \] The latter condition would make the digit expansion of $x$ eventually equal to $(M-1)$,
contrary to our convention. Hence $y - x > 0$, proving $x < y$.  
\vskip0.1in
\noindent The relation \eqref{dist(x, y) upper} is a direct consequence of \eqref{hD(x,y)}. Since $x, y$ both lie in $D(x, y)$, the distance between them is at most the length of $D(x,y)$. Thus 
\[ y - x \leq |D(x, y)| = M^{-\kappa+1}. \]
This completes the proof.  
\end{proof} 
The first separating height can be used in certain cases to determine an ordering among three points. Here is an example. 
\begin{lemma} \label{lemma: fsh transitive}
Let $x, y, z \in [0, 1)$ with $x < y$ and $x < z$. If 
\[ \kappa_1 = \kappa(x, y) < \kappa_2 = \kappa(x, z), \]
then $y < z$. 
\end{lemma} 
\begin{proof}
Using the definition \eqref{what is nu} of the first separating height, we may write $y - z$ as 
\begin{align*} y - z &= (y-x) - (z - x) = (y_{\kappa_1} - x_{\kappa_1}) M^{-\kappa}+ \sum_{k > \kappa_1} (y_{k} - z_k) M^{k} \\
&\leq - M^{-\kappa} + (M-1) \sum_{k > \kappa_1} M^{-k} = -M^{-\kappa} + M^{-\kappa} = 0.  
\end{align*}  
The estimate is an equality only if $y_{k} = M-1$ and $z_k = 0$ for all $k > \kappa_1$. This is not possible since digit expansions with all but finitely many digits equal to $(M-1)$ are not permitted. Thus $y < z$.  
\end{proof}  
\subsection{Euclidean distance and first separating height} 
The Euclidean distance between $x$ and $y$ may be much smaller than the upper bound in \eqref{dist(x, y) upper}. The relevant scale depends on what the rays $\mathcal R(x)$ and $\mathcal R(y)$ do after their first separation. If they enter nonadjacent children, the first separating height $\kappa$ controls their Euclidean distance. If they enter adjacent children, the rays may continue moving toward the common boundary of those children for many further generations, producing substantial cancellation. The following lemma makes this distinction precise.
%In general, the Euclidean distance between $x$ and $y$ could be much smaller than the upper bound specified by \eqref{dist(x, y) upper}. The $M$-adic scale of the Euclidean distance is, in fact, dictated by the behaviour of the rays depicting $x$ and $y$ right after they separate on the $M$-adic tree. The following lemma makes this connection precise. 
\begin{lemma} \label{lemma: Euclidean vs Madic distance} 
Let $x,y\in[0,1)$ with $x<y$, and let $\kappa=\kappa(x,y)$ be their first separating height.
\begin{enumerate}[(a)] 
\item \label{lemma: Euclidean vs Madic distance parta}  If $x$ and $y$ are separated by at least one $M$-adic interval at their first separating height, i.e. if  
\begin{equation} \label{diff at height nu} 
y_{\kappa} - x_{\kappa} \geq 2, 
\end{equation} 
then in addition to \eqref{dist(x, y) upper}, one has the following lower bound: 
\begin{equation} \label{dist(x, y) lower} 
y - x \geq M^{-\kappa}. 
\end{equation} 
\item  \label{lemma: Euclidean vs Madic distance partb} If $x$ and $y$ lie in adjacent $M$-adic intervals at their first separating height, i.e. if
\[ y_{\kappa} - x_{\kappa} =1, \] then the distance between $x$ and $y$ is given by 
\begin{equation} \label{dist(x, y) both sided} 
M^{-\mu} \leq y-x \leq 2M^{-\mu+1}, 
\end{equation} 
where $\mu$ denotes the first height beyond $\kappa$ where the difference of the digits of $y$ and $x$ exceeds its minimum possible value: 
\begin{equation} \label{def: mu>nu}
\mu = \mu(x, y) := \min \bigl\{k > \kappa: y_k - x_k > -(M-1) \bigr\}. 
\end{equation}   
\end{enumerate} 
\end{lemma} 
\vskip0.1in
\noindent {\em{Remarks: }} 
\begin{enumerate}[1.]
\item Combining parts \eqref{lemma: Euclidean vs Madic distance parta} and \eqref{lemma: Euclidean vs Madic distance partb} of Lemma \ref{lemma: Euclidean vs Madic distance} leads to a summary statement: any pair of points $x, y \in [0, 1)$ with $x < y$ obeys 
\begin{equation} \label{distance scales summary}
\left\{
\begin{aligned}  
& \hskip1in M^{-\mu} \leq y-x \leq M^{-\mu+2} \text{ where }\\ 
& \; \mu = \mu(x,y) := \begin{cases} \kappa(x, y) \text{ given by } \eqref{what is nu} &\text{ if } y_{\kappa} - x_{\kappa} \geq 2, \\
\mu(x, y) \text{ given by } \eqref{def: mu>nu} &\text{ if } y_{\kappa} - x_{\kappa} = 1.   
\end{cases}
\end{aligned} 
\right\}
\end{equation}  
The definition above applies when $x<y$. For $x>y$, we define \[\mu(x, y) = \mu(y, x), \]
so that $\mu$ is symmetric on pairs of distinct points. 
\vskip0.1in
\item The set of indices $k$ in \eqref{def: mu>nu} is non-empty, so $\mu$  is well-defined. Indeed, 
\[ y_k - x_k = -(M-1) \text{ holds if and only if } y_k = 0 \text{ and }x_k = (M-1). \] If this equality held for every $k>\kappa$, that would imply $x_k=M-1$ for every $k>\kappa$, contrary to our convention of excluding expansions that are eventually equal to $M-1$.
\vskip0.1in 
\item Suppose that $\mu > \kappa + 1$. By the definition of $\mu$, 
\[ x_k = M-1, \quad y_k = 0 \text{ for every } \kappa < k < \mu.  \] 
This has a geometric interpretation on the $M$-adic tree: after the rays $\mathcal R(x)$ and $\mathcal R(y)$ separate at height $\kappa$, and until they reach height $\mu-1$, the ray $\mathcal R(x)$ repeatedly follows the rightmost child of its current vertex, while $\mathcal R(y)$ repeatedly follows the leftmost child. This sustained movement toward the common boundary of the two adjacent $M$-adic intervals explains why $y-x$ is governed by the later scale $M^{-\mu}$, rather than by the first separating scale $M^{-\kappa}$.
% the height $\mu(x, y)$ depicting the Euclidean distance between $x, y$ admits a geometric interpretation on the $M$-adic tree: after separating from $\mathcal R(y)$ at height $\kappa$ and until the height $(\mu-1)$, the ray $\mathcal R(x)$ keeps moving continually and maximally to the right within its parent. In other words, every edge of $x$ between these heights joins a parent to its $(M-1)^{\text{th}}$, i.e. rightmost child. Similarly, the ray depicting $y$ keeps moving continually and maximally to the left between heights $\kappa$ and $(\mu-1)$, joining a parent to its $0^{\text{th}}$, or leftmost, child. 
\end{enumerate} 

\begin{figure}[t]
\centering

\resizebox{0.96\textwidth}{!}{%
\begin{tikzpicture}[
    x=1cm,
    y=1cm,
    every node/.style={font=\scriptsize},
    vertex/.style={
        circle,
        draw,
        fill=white,
        inner sep=1.35pt
    },
    xray/.style={
        blue!65!black,
        line width=0.9pt
    },
    yray/.style={
        red!70!black,
        line width=0.9pt
    },
    ordinary/.style={
        black!55,
        line width=0.55pt
    },
    guide/.style={
        black!40,
        densely dashed,
        line width=0.45pt
    },
    interval/.style={
        draw=black!45,
        line width=0.45pt
    }
]

% Optional divider between the two panels
\draw[black!18] (0,-0.05) -- (0,5.15);

% ============================================================
% Panel (a): nonadjacent children
% ============================================================

\begin{scope}[xshift=-3.75cm]

\node[font=\small] at (0,5.15)
    {(a) \(y_\kappa-x_\kappa\geq 2\)};

% Youngest common ancestor
\node[vertex] (Da) at (0,4.40) {};
\node[anchor=west] at (-1,4.6) {\(D(x,y)\),};
\node[anchor=west,text=black!55] at (-0.05,4.6)
    {\(h=\kappa-1\)};

% Children at the first separatinge height
\node[
    vertex,
    draw=blue!65!black,
    fill=blue!65!black
] (Xa) at (-1.65,3.35) {};

\node[vertex] (Ia) at (0,3.35) {};

\node[
    vertex,
    draw=red!70!black,
    fill=red!70!black
] (Ya) at (1.65,3.35) {};

\draw[xray]     (Da) -- (Xa);
\draw[ordinary] (Da) -- (Ia);
\draw[yray]     (Da) -- (Ya);

\node[
    anchor=east,
    text=blue!65!black
] at (-1.7,3.5) {\(x_\kappa\)};

\node[
    anchor=west,
    text=red!70!black
] at (1.7,3.5) {\(y_\kappa\)};

\draw[guide] (-2.75,3.35) -- (2.75,3.35);

\node[
    anchor=west,
    text=black!55
] at (2.82,3.35) {\(\kappa\)};

% A small amount of continuation below the separating height
\coordinate (XaL) at (-2.05,2.60);
\coordinate (XaR) at (-1.25,2.60);
\coordinate (YaL) at (1.25,2.60);
\coordinate (YaR) at (2.05,2.60);

\draw[xray] (Xa) -- (XaL);
\draw[xray] (Xa) -- (XaR);
\draw[yray] (Ya) -- (YaL);
\draw[yray] (Ya) -- (YaR);

\fill[blue!65!black] (XaL) circle (1.1pt);
\fill[blue!65!black] (XaR) circle (1.1pt);
\fill[red!70!black]  (YaL) circle (1.1pt);
\fill[red!70!black]  (YaR) circle (1.1pt);

% Corresponding intervals at height kappa
\draw[
    interval,
    fill=blue!8
] (-2.75,1.15) rectangle (-0.95,1.65);

\draw[
    interval,
    fill=black!8
] (-0.85,1.15) rectangle (0.85,1.65);

\draw[
    interval,
    fill=red!8
] (0.95,1.15) rectangle (2.75,1.65);

\fill[blue!65!black] (-1.55,1.40) circle (1.4pt);
\fill[red!70!black]  (1.55,1.40) circle (1.4pt);

\node[text=blue!65!black] at (-1.55,0.91) {\(x\)};
\node[text=red!70!black]  at (1.55,0.91) {\(y\)};
\node at (0,1.40) {\(M^{-\kappa}\)};

% Distance marker
\draw[line width=0.5pt]
    (-1.55,0.52) -- (1.55,0.52);

\draw[line width=0.5pt]
    (-1.55,0.43) -- (-1.55,0.61);

\draw[line width=0.5pt]
    (1.55,0.43) -- (1.55,0.61);

\node[
    fill=white,
    inner sep=1.2pt
] at (0,0.52) {\(y-x\geq M^{-\kappa}\)};

\end{scope}

% ============================================================
% Panel (b): adjacent children
% ============================================================

\begin{scope}[xshift=3.75cm]

\node[font=\small] at (0,5.15)
    {(b) \(y_\kappa-x_\kappa=1\)};

% Youngest common ancestor
\node[vertex] (Db) at (0,4.40) {};
\node[anchor=west] at (-1,4.6) {\(D(x,y)\), };
\node[anchor=west,text=black!55] at (-0.05,4.6)
    {\(h=\kappa-1\)};

% Adjacent children
\node[
    vertex,
    draw=blue!65!black,
    fill=blue!65!black
] (Xb) at (-1.20,3.35) {};

\node[
    vertex,
    draw=red!70!black,
    fill=red!70!black
] (Yb) at (1.20,3.35) {};

\draw[xray] (Db) -- (Xb);
\draw[yray] (Db) -- (Yb);

\node[
    anchor=east,
    text=blue!65!black
] at (-1.2,3.5) {\(x_\kappa\)};

\node[
    anchor=west,
    text=red!70!black
] at (1.2,3.5) {\(y_\kappa\)};

\draw[guide] (-2.75,3.35) -- (2.75,3.35);

\node[
    anchor=west,
    text=black!55
] at (2.82,3.35) {\(\kappa\)};

% Repeated movement toward the common boundary
\coordinate (XbOne)   at (-0.72,2.85);
\coordinate (XbTwo)   at (-0.40,2.35);
\coordinate (XbThree) at (-0.20,1.85);

\coordinate (YbOne)   at (0.72,2.85);
\coordinate (YbTwo)   at (0.40,2.35);
\coordinate (YbThree) at (0.20,1.85);

\draw[xray]
    (Xb) -- (XbOne) -- (XbTwo) -- (XbThree);

\draw[yray]
    (Yb) -- (YbOne) -- (YbTwo) -- (YbThree);

\foreach \P in {XbOne,XbTwo,XbThree}
    \fill[blue!65!black] (\P) circle (1.1pt);

\foreach \P in {YbOne,YbTwo,YbThree}
    \fill[red!70!black] (\P) circle (1.1pt);

\draw[guide] (-2.75,1.85) -- (2.75,1.85);

\node[
    anchor=west,
    text=black!55
] at (2.82,1.85) {\(\mu-1\)};

% Adjacent intervals at height kappa
\draw[
    interval,
    fill=blue!8
] (-2.75,1.15) rectangle (0,1.65);

\draw[
    interval,
    fill=red!8
] (0,1.15) rectangle (2.75,1.65);

% Points close to the common boundary
\fill[blue!65!black] (-0.30,1.40) circle (1.4pt);
\fill[red!70!black]  (0.30,1.40) circle (1.4pt);

\node[text=blue!65!black] at (-0.30,0.91) {\(x\)};
\node[text=red!70!black]  at (0.30,0.91) {\(y\)};

% Distance marker
\draw[line width=0.5pt]
    (-0.30,0.52) -- (0.30,0.52);

\draw[line width=0.5pt]
    (-0.30,0.43) -- (-0.30,0.61);

\draw[line width=0.5pt]
    (0.30,0.43) -- (0.30,0.61);

\node[
    fill=white,
    inner sep=1.2pt
] at (0,0.52)
    {\(M^{-\mu}\leq y-x\leq 2M^{-\mu+1}\)};

\end{scope}

\end{tikzpicture}%
}

\caption{\small{The two alternatives in Lemma \ref{lemma: Euclidean vs Madic distance}.
If the rays enter nonadjacent children at their first separating
height \(\kappa\), their Euclidean separation is bounded below at
scale \(M^{-\kappa}\).
If they enter adjacent children, continued movement toward their
common boundary can make the points much closer; the first height
\(\mu\) at which this maximal cancellation fails determines their
Euclidean separation.}}
\label{fig:first-separating-height}
\end{figure}

\begin{proof}[Proof of Lemma \ref{lemma: Euclidean vs Madic distance}] 
Let us start with part \eqref{lemma: Euclidean vs Madic distance parta}. 
Under the assumption \eqref{diff at height nu}, the quantity $(y-x)$ can be  
estimated from below, as follows: 
\begin{align*} 
y-x &= \sum_{k = \kappa}^{\infty} (y_k - x_k) M^{-k} \geq  2 M^{-\kappa} + \sum_{k = \kappa+1}^{\infty} (y_k - x_k) M^{-k} \\ 
& \geq 2M^{-\kappa} - (M-1) M^{-(\kappa + 1)} \frac{1}{1 - \frac{1}{M}} \geq 2M^{-\kappa} - M^{-\kappa} = M^{-\kappa}. 
\end{align*}
This is the claimed bound \eqref{dist(x, y) lower}. 
\vskip0.1in 
\noindent We turn now to part \eqref{lemma: Euclidean vs Madic distance partb}. If $\mu = \kappa + 1$, then the definition of $\mu$ gives 
\[y_{\kappa+1} - x_{\kappa+1} \geq -(M-2). \] This leads to the bound
\begin{align*} 
y - x &= M^{-\kappa} + \bigl(y_{\kappa+1} - x_{\kappa+1}\bigr) M^{-(\kappa + 1)} + \sum_{k = \kappa+2}^{\infty} (y_k - x_k) M^{-k}\\ 
&\geq M^{-\kappa} - (M-2) M^{-(\kappa + 1)} - (M-1) \sum_{k = \kappa+2}^{\infty} M^{-k} \\ 
&\geq 2M^{-(\kappa + 1)} - M^{-(\kappa + 1)} = M^{-(\kappa +1)} = M^{-\mu}, 
\end{align*} 
which agrees with the left inequality of \eqref{dist(x, y) both sided}. For the upper bound, the assumption $y_{\kappa} - x_{\kappa} = 1$ gives 
\begin{align*}
y - x &= M^{-\kappa} + \sum_{k > \kappa} (y_k-x_k) M^{-k} \\
&\leq M^{-\kappa} + (M-1) \sum_{k > \kappa}M^{-k} = 2M^{-\kappa} = 2M^{-\mu+1}. 
\end{align*} 
This establishes \eqref{dist(x, y) both sided} when $\mu = \kappa +1$. 
\vskip0.1in
\noindent Next, suppose that $\mu > \kappa + 1$. In this case,
\[  y_k - x_k = \begin{cases} 1 &\text{ for } k = \kappa, \\  -(M-1) &\text{ for } \kappa < k < \mu, \end{cases} \] and $y_k - x_k \geq -(M-2)$ for $k = \mu$.
These relations isolate the exact cancellation occurring between heights $\kappa$ and $\mu-1$. More precisely,
%Using these, one can compute certain terms and sub-sums occurring in the expansion of $(y-x)$. For instance,
\begin{equation} 
(y_{\kappa} - x_{\kappa}) M^{-\kappa} = M^{-\kappa}, \quad (y_{\mu} - x_{\mu}) M^{-\mu} \geq -(M-2) M^{-\mu}, \label{diff-1}
\end{equation}
whereas
\begin{align}
\sum_{k= \kappa+1}^{\mu - 1} (y_{k} - x_{k}) M^{- k} &= - (M-1) \sum_{k= \kappa+1}^{\mu - 1} M^{- k}  \nonumber \\
&= -(M-1) \frac{M^{-\kappa-1} - M^{-\mu}}{1 - \frac{1}{M}} \nonumber \\ 
&= - \bigl[ M^{-\kappa} - M^{-\mu +1} \bigr], \text{ and }  \label{diff-2} \\ 
\Bigl|\sum_{k = \mu + 1}^{\infty} (y_k - x_k) M^{-k} \Bigr| &\leq (M-1) \sum_{k=\mu+1}^{\infty} M^{-k} \leq M^{-\mu}. \label{diff-3}
\end{align} 
Combining \eqref{diff-1}--\eqref{diff-3} allows us to estimate $(y-x)$ both from above and from below. On one hand,   
\begin{align*}  
y - x & = (y_{\kappa} - x_{\kappa}) M^{-\kappa} + \sum_{k= \kappa+1}^{\mu - 1} (y_{k} - x_{k}) M^{- k}  \\
&\hskip1in + (y_{\mu} - x_{\mu}) M^{-\mu}  + \sum_{k = \mu + 1}^{\infty} (y_k - x_k) M^{-k} \\
&\geq M^{-\kappa} - \bigl[M^{-\kappa} - M^{-\mu+1} \bigr] - (M-2) M^{-\mu} - M^{-\mu} \geq M^{-\mu}.
\end{align*} 
This is the left inequality claimed in \eqref{dist(x, y) both sided}. The same relations also yield 
\begin{align*}
y - x & = (y_{\kappa} - x_{\kappa}) M^{-\kappa} + \sum_{k= \kappa+1}^{\mu - 1} (y_{k} - x_{k}) M^{- k} + \sum_{k = \mu}^{\infty} (y_k - x_k) M^{-k} \\
&\leq M^{-\kappa} - \bigl[M^{-\kappa} - M^{-\mu+1} \bigr] + (M-1) \frac{M^{-\mu}}{1 - \frac{1}{M}} \leq 2M^{-\mu+1}
\end{align*} 
This is the right inequality in \eqref{dist(x, y) both sided}.
\end{proof} 
\section{Splitting behaviour of an ordered tuple}
The results on $M$-adic and Euclidean distance obtained in the preceding section allow us to study the $M$-adic tree generated by a strictly ordered collection of points. After applying an $M$-adic affine transformation of the form \eqref{T: affine tf} if necessary, we may assume that the collection lies in $[0,1)$. By Lemma \ref{lemma: trees under affine maps}, the transformation induces a tree isomorphism, and therefore preserves all relevant tree properties. We enumerate the points in decreasing order. 
\vskip0.1in
\noindent For $p \geq 3$, let $\mathcal P_{p}$ denote a  $p$-tuple of real numbers:
\begin{equation} \label{ordered p tuple}
\mathcal P_p := \bigl\{ \varrho(i) : 1 \leq i \leq p \bigr\} \subseteq [0,1), \text{ with } 
\varrho(i+1) < \varrho(i), \; 1 \leq i < p. 
\end{equation}
The $M$-adic digit sequence of $\varrho(i)$ is denoted by ${\pmb{\rho}(i)} := (\rho_{1}(i), \rho_{2}(i), \ldots)$:
\begin{equation} \label{rhok}
 \varrho(i) = \sum_{k=1}^{\infty} \rho_k(i) M^{-k}, \quad 
 \rho_k(i) \in \mathbb Z_M := \{0, 1, \ldots, M-1\} \text{ for all } k \geq 1.
\end{equation}  
For brevity,  we write $\mathcal R[i]$ for the ray encoding $\varrho(i)$, and $D[i, i'], \kappa[i, i']$ for the youngest common ancestor and first separating height, respectively,  of the pair $(\varrho(i), \varrho(i'))$. Thus
\begin{align} 
 &\mathcal R[i] := \mathcal R \bigl( \varrho(i)\bigr), \; D[i, i'] = D(\varrho(i), \varrho(i')), \nonumber \\ 
&\kappa[i, i'] := \kappa \bigl( \varrho(i), \varrho(i')\bigr), \ \quad 1 \leq i < i' \leq p.  \label{nuii'}
\end{align} 
The following lemma records several consequences of the linear ordering of $\mathcal P_p$ for the tree $\mathcal T(\mathcal P_p; M)$. 
\begin{lemma} \label{lemma: order and tree}
For $p \geq 3$ and $\mathcal P_p$ as in \eqref{ordered p tuple}, the following conclusions hold: 
\begin{enumerate}[(a)]
\item \label{lemma: order and tree parta} For any $1 \leq i_1 < i_2 < i_3 \leq p$, the first separating heights obey 
\begin{equation} \label{nu for three}
\kappa[i_1, i_3] \leq \kappa[i_2, i_3]. 
\end{equation}
\vskip0.1in
\item \label{lemma: order and tree parta'} The inequality \eqref{nu for three} is strict if and only if 
\[ D[i_1, i_3] \text{ and } D[i_2, i_3] \text{ are distinct splitting vertices on $\mathcal R[i_3]$ } \] 
in the tree $\mathcal T(\mathcal P_p; M)$.
\vskip0.1in 
\item \label{lemma: order and tree partb} All the splitting vertices of $\mathcal T(\mathcal P_p; M)$ lie on $\mathcal R[p]$ if and only if 
\begin{equation} \label{arrangement of nu} 
\kappa[i, i'] = \kappa[i, p] \quad \text{ for every } 1 \leq i < i' \leq p-1. 
\end{equation} 
\end{enumerate}
\end{lemma} 
\noindent The configuration described in part \eqref{lemma: order and tree partb} is illustrated in
Figure~\ref{fig:ordered-tuple-spine}. In particular, several rays
may leave $\mathcal R[p]$ at the same splitting vertex, provided
they enter distinct children of that vertex; the corresponding
first separating heights are then equal.
\begin{figure}[t]
\centering

\resizebox{0.92\textwidth}{!}{%
\begin{tikzpicture}[
    x=1cm,
    y=0.9cm,
    every node/.style={font=\scriptsize},
    vertex/.style={
        circle,
        draw,
        fill=white,
        inner sep=1.4pt
    },
    childvertex/.style={
        circle,
        inner sep=1.5pt
    },
    rone/.style={
        blue!65!black,
        line width=0.9pt
    },
    rtwo/.style={
        red!70!black,
        line width=0.9pt
    },
    rthree/.style={
        green!45!black,
        line width=0.9pt
    },
    rfour/.style={
        black,
        line width=0.9pt
    },
    guide/.style={
        black!40,
        densely dashed,
        line width=0.45pt
    },
    continuation/.style={
        densely dotted,
        line width=0.75pt
    }
]

% Ordered tuple
%\node[font=\small] at (0,7.15)
%    {$\varrho(1)>\varrho(2)>\varrho(3)>\varrho(4)$};

% ------------------------------------------------------------
% The minimal ray R[4]
% ------------------------------------------------------------

\coordinate (Top) at (-2.40,6.65);
\coordinate (D12) at (-2.40,5.85);
\coordinate (C4a) at (-2.40,5.15);
\coordinate (D34) at (-2.40,3.55);
\coordinate (C4b) at (-2.40,2.85);
\coordinate (R4end) at (-2.40,1.10);

\draw[rfour]
    (Top) -- (D12) -- (C4a) -- (D34) -- (C4b) -- (R4end);

\node[anchor=west] at (-2.24,1.38) {$\mathcal R[4]$};
\node[anchor=west,text=black!55] at (-2.24,1.10)
    {minimal ray};

% ------------------------------------------------------------
% First splitting vertex:
% R[1] and R[2] leave through distinct children
% ------------------------------------------------------------

\node[vertex] at (D12) {};

\node[anchor=west] at (-2.2,6.4)
    {$D[1,2]=D[1,4]=D[2,4]$};

\node[anchor=west,text=black!55] at (-2.20,6)
    {height $\kappa[1,4]-1$};

% Children one generation below D12
\coordinate (C2) at (-0.80,5.15);
\coordinate (C1) at (1.15,5.15);

\draw[rtwo]
    (D12) -- (C2) -- (1.55,3.30);

\draw[rone]
    (D12) -- (C1) -- (3.05,3.00);

\node[
    childvertex,
    fill=red!70!black
] at (C2) {};

\node[
    childvertex,
    fill=blue!65!black
] at (C1) {};

\node[
    childvertex,
    fill=black
] at (C4a) {};

% First separating-height guide
\draw[guide]
    (-3.35,5.15) -- (3.75,5.15);

\node[
    anchor=west,
    text=black!55
] at (3.82,5.15)
    {$\kappa[1,4]=\kappa[2,4]$};

\node[
    anchor=west,
    text=red!70!black
] at (1.63,3.30)
    {$\mathcal R[2]$};

\node[
    anchor=west,
    text=blue!65!black
] at (3.13,3.00)
    {$\mathcal R[1]$};

% ------------------------------------------------------------
% Later splitting vertex:
% R[3] leaves R[4]
% ------------------------------------------------------------

\node[vertex] at (D34) {};

\node[anchor=west] at (-2.20,3.80)
    {$D[3,4]$};

\node[anchor=west,text=black!55] at (-2.20,3.49)
    {height $\kappa[3,4]-1$};

\coordinate (C3) at (-1.30,2.85);

\draw[rthree]
    (D34) -- (C3) -- (-0.35,1.10);

\node[
    childvertex,
    fill=green!45!black
] at (C3) {};

\node[
    childvertex,
    fill=black
] at (C4b) {};

% Second separating-height guide
\draw[guide]
    (-3.35,2.85) -- (3.75,2.85);

\node[
    anchor=west,
    text=black!55
] at (3.82,2.85)
    {$\kappa[3,4]$};

\node[
    anchor=west,
    text=green!45!black
] at (-0.27,1.10)
    {$\mathcal R[3]$};

% ------------------------------------------------------------
% Dotted continuations to the represented points
% ------------------------------------------------------------

\coordinate (P4) at (-2.40,0);
\coordinate (P3) at (-0.35,0);
\coordinate (P2) at (1.55,0);
\coordinate (P1) at (3.05,0);

\draw[
    rfour,
    continuation
] (R4end) -- (P4);

\draw[
    rthree,
    continuation
] (-0.35,1.10) -- (P3);

\draw[
    rtwo,
    continuation
] (1.55,3.30) -- (P2);

\draw[
    rone,
    continuation
] (3.05,3.00) -- (P1);

% ------------------------------------------------------------
% The ordered points on the real line
% ------------------------------------------------------------

\draw[black,line width=0.6pt]
    (-3.15,0) -- (3.75,0);

\fill[black]             (P4) circle (1.5pt);
\fill[green!45!black]    (P3) circle (1.5pt);
\fill[red!70!black]      (P2) circle (1.5pt);
\fill[blue!65!black]     (P1) circle (1.5pt);

\node[text=black] at (-2.40,-0.38)
    {$\varrho(4)$};

\node[text=green!45!black] at (-0.35,-0.38)
    {$\varrho(3)$};

\node[text=red!70!black] at (1.55,-0.38)
    {$\varrho(2)$};

\node[text=blue!65!black] at (3.05,-0.38)
    {$\varrho(1)$};

\node[
    anchor=west,
    text=black!55
] at (3.82,0)
    {$\mathbb R$};

% Monotonicity relation
\node[font=\small] at (0,-1.05)
    {$\kappa[1,4]=\kappa[2,4]<\kappa[3,4]$};

\end{tikzpicture}%
}

\caption{\small{The splitting structure of an ordered quadruple when every
splitting vertex lies on the minimal ray $\mathcal R[4]$.
The rays $\mathcal R[1]$ and $\mathcal R[2]$ leave
$\mathcal R[4]$ through distinct children of the same splitting
vertex, whereas $\mathcal R[3]$ leaves at a later height.
The dotted continuations indicate the points represented by the
four rays.}}
\label{fig:ordered-tuple-spine}
\end{figure}

\begin{proof} 
We begin with part \eqref{lemma: order and tree parta}. Let us assume for contradiction that
\[ m = \kappa[i_1, i_3] > n = \kappa[i_2, i_3]. \] 
Because $m > n$, the rays $\mathcal R[i_1]$ and $\mathcal R[i_3]$ still agree at height $n$, so
\[ \rho_n(i_1) = \rho_n(i_3). \]
Since $n = \kappa[i_2, i_3]$ and $\varrho(i_3) < \varrho(i_2)$, Lemma \ref{lemma: distance x y} gives 
\[ \rho_n(i_2) > \rho_n(i_3) = \rho_n(i_1). \] 
Moreover, all three digit sequences agree before height $n$. Consequently, 
 \[ \kappa[i_1, i_2] = n \; \text{ and } \; \rho_n(i_1) < \rho_n(i_2). \]
Applying Lemma \ref{lemma: distance x y} to the pair $\varrho(i_1), \varrho(i_2)$ then gives 
\[ \varrho(i_1) < \varrho(i_2), \] 
contradicting their prescribed ordering.    
\vskip0.1in
\noindent Next we consider part \eqref{lemma: order and tree parta'}. Suppose first that $\kappa[i_2, i_3] > \kappa[i_1, i_3]$. Let us recall from \eqref{hD(x,y)} that the first separating height of a pair is one higher than that of their youngest common ancestor. Therefore the assumed inequality means that 
\begin{equation} \label{different heights}
h \bigl(D[i_2, i_3] \bigr) = \kappa[i_2, i_3] - 1 > \kappa[i_1, i_3] - 1 = h \bigl(D[i_1, i_3] \bigr). 
\end{equation}  
Since both splitting vertices $D[i_1, i_3]$ and $D[i_2, i_3]$ lie on $\mathcal R[i_3]$ by definition, \eqref{different heights} implies they must be distinct.  
\vskip0.1in
\noindent Conversely, suppose that $D[i_1,i_3]$ and $D[i_2,i_3]$ are distinct. Since both vertices lie on the same ray $\mathcal R[i_3]$, they have different heights. Part \eqref{lemma: order and tree parta} gives 
\[ \kappa[i_1, i_3] \leq  \kappa[i_2, i_3], \text{ or equivalently }  h \bigl(D[i_1, i_3] \bigr) \leq h \bigl(D[i_2, i_3] \bigr). \] 
The heights cannot be equal, so both inequalities are strict. Therefore  $\kappa[i_1, i_3] < \kappa[i_2, i_3]$. 
\vskip0.1in
\noindent We turn to part \eqref{lemma: order and tree partb}. If \eqref{arrangement of nu} holds, then we deduce that
%$\bigl\{D[i,p] : 1 \leq i \leq p-1 \bigr\}$ is a collection of splitting vertices on $\mathcal R[p]$ of monotone non-decreasing heights. In particular, this means 
 \[ h \bigl(D[i, i'] \bigr) = h \bigl( D[i, p] \bigr) \text{ for all } 1 \leq i < i' \leq p-1. \] 
Since $D[i, i']$ and $D[i, p]$ both lie on $\mathcal R[i]$, this means   
\[D[i, i'] = D[i, p], \text{ placing each youngest common ancestor on } \mathcal R[p] .\] 
Every splitting vertex of $\mathcal T(\mathcal P_p;M)$ is the youngest common ancestor of two elements of $\mathcal P_p$, and hence equals $D[i,i']$ for some such pair.
It follows that every splitting vertex of $\mathcal T(\mathcal P_p; M)$ lies on $\mathcal R[p]$. 
%\noindent We turn now to part \eqref{lemma: order and tree partb}. 
\vskip0.1in
\noindent Conversely, suppose that all the splitting vertices of $\mathcal T(\mathcal P_p; M)$ lie on $\mathcal R[p]$. For each $1 \leq i \leq p-1$, let 
\[ v(i) := D[i, p], \] 
the last vertex of $\mathcal R[i]$ lying on $\mathcal R[p]$.
The child of $v(i)$ lying on $\mathcal R[i]$ does not belong to $\mathcal R[p]$. The maximal subtree rooted at this child contains
no splitting vertices and therefore consists solely of the
continuation of $\mathcal R[i]$.
%In other words, all the first separating heights 
%\[ \{\kappa[i, p]: 1 \leq i \leq p-1 \} \text{ are distinct}. \] 
In light of the linear ordering  \eqref{ordered p tuple} among $\{\varrho(i) : 1 \leq i \leq p\}$ and the inequality \eqref{nu for three} proved in part \eqref{lemma: order and tree parta}, we obtain 
\begin{align} 
&\kappa[1,p] \leq \kappa[2,p] \leq \cdots \leq \kappa[p-1,p], \text{ and therefore } \nonumber \\  
&D[1, p] \supseteq D[2, p] \supseteq \cdots \supseteq D[p-1, p] \; \text{ along } \mathcal R[p]. 
\end{align} 
We fix $1 \leq i < i' \leq p-1$ and consider two cases:
\begin{equation} \label{casification} 
\left\{
\begin{aligned} 
&{\text{\em{Case 1: }}} v(i) = D[i, p] \neq D[i', p] = v(i'), \;  \\ 
&{\text{\em{Case 2: }}} v(i) = D[i, p] = D[i', p] = v(i'). 
\end{aligned} 
\right\}
\end{equation} 
\vskip0.1in
\begin{itemize} 
\item In case 1 of \eqref{casification}, the ray $\mathcal R[i]$ leaves $\mathcal R[p]$ at $v(i)$, while $\mathcal R[i']$ continues along $\mathcal R[p]$ beyond $v(i)$. Hence 
\begin{equation} \label{ycm equality ii'p} D[i, i'] = v(i) = D[i, p] \text{ in this case.} 
\end{equation} 
\item In case 2 of \eqref{casification}, the rays $\mathcal R[i]$ and $\mathcal R[i']$ must enter distinct children of $v(i)$. Otherwise, they would have a later youngest
common ancestor away from $\mathcal R[p]$, producing a splitting vertex
off $\mathcal R[p]$. Thus \eqref{ycm equality ii'p} holds in this case as well, establishing it for all $1 \leq i < i' \leq p-1$.  
\end{itemize} 
\vskip0.1in
\noindent Taking heights of the vertices in \eqref{ycm equality ii'p} and using \eqref{hD(x,y)}, we obtain 
\[
\kappa[i,i']
 = h\bigl(D[i,i']\bigr)+1
 = h\bigl(D[i,p]\bigr)+1
 = \kappa[i,p],
\]
which is \eqref{arrangement of nu}.
\end{proof} 
%\vskip0.1in
%\noindent Conversely, suppose that $\mathcal T(\mathcal P_p; M)$ has a splitting vertex $v$ that does not lie on $\mathcal R(1)$. Let
%\[ v = D \bigl(\varrho(i_2), \varrho(i_3) \bigr) \text{ for some } 2 \leq i_2 < i_3 \leq p. \]
%This means that the youngest common ancestor $w$ of $\varrho(1)$ and $v$ must obey
%\[ w = D\bigl(\varrho(1), \varrho(i_2) \bigr) = D\bigl(\varrho(1), \varrho(i_3) \bigr); \text{ in other words } \kappa_{1, i_2} = \kappa_{1, i_3} = h(w)+1. \]
%This contradicts one of the strict inequalities in \eqref{arrangement of nu} and completes the proof.  

\section{A finite lacunary tuple} In this section and the next, we will consider a finite collection of points that occurs as a subset of a decreasing, infinite lacunary sequence with small lacunarity constant. The notion is introduced in this section, and the branching behaviour of rays depicting such a set of points is addressed in the next. This analysis forms the basis of the proof of Proposition \ref{TREE STRUCTURE THEOREM}, which prescribes the tree lacunarity of a monotone lacunary sequence.

\begin{definition}
 Given any constant $C_0 \geq 1$, an ordered $p$-tuple $\mathcal P_p$ given by \eqref{ordered p tuple} is called a {\em{$(M, C_0)$-lacunary tuple}}, and we write 
 \begin{align} 
 \mathcal P_p &\in {\tt{Lac}}_p(M, C_0), \text{ if } \\
 \bigl( \varrho(i+1) - \varrho(p) \bigr) &\leq M^{-C_0} \bigl( \varrho(i) - \varrho(p) \bigr) \text{ for } 1 \leq i < p.  
 \end{align} 
 \end{definition}
\vskip0.1in
\noindent The following properties are immediate consequences of the definition. Their verification is left to the reader.
\begin{itemize} 
\item The class ${\tt{Lac}}_p(M, \cdot)$ is monotone with respect to the lacunarity constant: 
\begin{equation} 
\quad {\tt{Lac}}_p(M, C_1) \subseteq {\tt{Lac}}_p(M, C_0) \text{ if } C_1 \geq C_0. 
\end{equation} 
\vskip0.1in 
\item The property of $(M, C_0)$-lacunarity is preserved under containment, provided the subset retains the minimal element of the original set.  
\begin{equation} \label{lacunarity under containment}
\left\{ 
\begin{aligned} 
&\text{if $\mathcal P_p \in {\tt{Lac}}_p(M, C_0)$ and $3 \leq q < p$, then for any } \mathcal P' \subseteq \mathcal P_p, \\  
&\mathcal P' \in {\tt{Lac}}_{q}(M, C_0) \text{ provided } \; \#(\mathcal P') = q, \text{ and } \varrho(p) \in \mathcal P'. 
\end{aligned}
\right\}
\end{equation}  
If $\varrho(p) \notin \mathcal P'$, then $\mathcal P'$ is still a lacunary tuple, though possibly with a smaller constant $C_0$. The exact statement appears in Lemma \ref{lemma: lacunary subtuple}.  
\end{itemize} 
\subsection{A lacunary tuple as a subset of a lacunary sequence} The definition of the class of lacunary tuples ${\tt{Lac}}_p(M, C_0)$ is superficially reminiscent of a finitary version of ${\tt{MonLac}}(M^{-C_0})$. Unlike the infinite setting, however, $\varrho(p)$ is the minimal element of the tuple rather than a distinguished limit point. Not surprisingly, $p$-element subsets of a sequence in ${\tt{MonLac}}(\cdot)$ belong to ${\tt{Lac}}_p(M,\cdot)$, after a suitable adjustment of the lacunarity constant. The next lemma, which makes this connection precise, is useful in the analysis of the branching structure of $\mathcal T(\mathcal P_p; M)$ in the next section, and ultimately in the proof of Proposition \ref{TREE STRUCTURE THEOREM}.
\begin{lemma} \label{lemma: lacunary sequence vs tuple}
Suppose that $\mathfrak A$ is a strictly decreasing sequence, 
\[ \mathfrak A \in {\tt{MonLac}}(M^{-C_0-1}), \quad \alpha = \lim \mathfrak A. \]
Then for every $p \geq 3$, any $p$-element subset $\mathcal P_p \subseteq \bar{\mathfrak A} = \mathfrak A \cup \{\alpha\}$,  
ordered as in \eqref{ordered p tuple}, belongs to ${\tt{Lac}}_p(M, C_0)$. 
\end{lemma} 
\begin{proof} 
For a collection 
\[ \mathcal P_p = \{ \varrho(p) < \varrho(p-1) < \cdots < \varrho(1) \} \subseteq \bar{\mathfrak A}, \] we first identify indices $1 \leq i_1  < i_2 < \cdots < i_{p-1} < i_p \leq \infty$ such that 
\[ \varrho(k) =  a_{i_k}, \; 1 \leq k \leq p, \; \text{ where } a_{\infty} := \alpha.  \] We consider two cases. 
\vskip0.1in
\noindent {\em{Case 1: $\alpha \in \mathcal P_p$, i.e., the limit point of $\mathfrak A$ is a member of $\mathcal P_p$. }}
 In this case, $i_p = \infty$ and $\varrho(p) = \alpha$. Then the lacunarity condition implies that for each $1 \leq k < p$, 
\begin{align*} 
\bigl[ \varrho(k+1) - \varrho(p) \bigr] = \bigl[ a_{i_{k+1}} - \alpha \bigr] & \leq M^{-(C_0+1)(i_{k+1} - i_k)} \bigl[a_{i_{k}} - \alpha \bigr] \\ 
&\leq M^{-C_0} \bigl[ \varrho(k) - \varrho(p) \bigr],  
\end{align*} 
since $i_{k+1} - i_k \geq 1$. This establishes that $\mathcal P_p \in {\tt{Lac}}_p(M, C_0)$.
\vskip0.1in
\noindent {\em{Case 2: $\alpha \notin \mathcal P_p$, i.e. $\mathcal P_p \subsetneq \mathfrak A$. }} In this case $i_p < \infty$. The monotonicity and lacunarity assumptions on $\mathfrak A$ then give 
\begin{align} 
 \left[ \varrho(k+1) - \varrho(p) \right] &\leq \left[ \varrho(k+1) - \alpha \right] \nonumber \\ 
 &\leq M^{-C_0-1} \left[ \varrho(k) - \alpha \right] \nonumber \\ 
&\leq M^{-C_0-1}\left[ \bigl(\varrho(k) - \varrho(p) \bigr) +  \bigl(\varrho(p) - \alpha \bigr) \right] \nonumber \\
 &\leq  M^{-C_0-1} \left[ \varrho(k) - \varrho(p) \right]  \left[ 1 + \frac{\varrho(p) - \alpha}{\varrho(k) - \varrho(p)} \right] \label{ratio-1} 
 \end{align} 
 We estimate the ratio in the second factor of the expression in  \eqref{ratio-1} as follows: 
 \begin{align} 
 \varrho(k) - \varrho(p) &= (\varrho(k) - \alpha) -  (\varrho(p) - \alpha) \nonumber \\
 &= \left[ \varrho(p) - \alpha \right] \times \left[ \frac{(\varrho(k) - \alpha)}{(\varrho(p) - \alpha)}-1 \right] \nonumber \\ 
 &\geq \left[ \varrho(p) - \alpha \right] \times \left[M^{(C_0+1)(p-k)}-1 \right]. \label{ratio-2}
 \end{align}  
 Inserting \eqref{ratio-2} into \eqref{ratio-1} leads to the estimate:
 \begin{align*} \frac{\left( \varrho(k+1) - \varrho(p) \right)}{\left( \varrho(k) - \varrho(p) \right)} &\leq M^{-C_0-1} \left[ 1 + \frac{1}{M^{(C_0+1)(p-k)}-1}\right] \\ 
 &= M^{-C_0-1} \frac{M^{(C_0+1)(p-k)}}{M^{(C_0+1)(p-k)}-1} \leq M^{-C_0},
 \end{align*}  
 where the last inequality follows from $M\geq 2$ and $(C_
0+1)(p-k)\geq1$. This establishes the desired conclusion $\mathcal P_p \in {\tt{Lac}}_p(M, C_0)$ in the second case. 
\end{proof} 
\begin{lemma} \label{lemma: lacunary subtuple} 
Suppose that $\mathcal P_p \in {\tt{Lac}}_p(M, C_0)$ for some $p > 3$ and $C_0 > 1$. Then for $3 \leq q < p$, 
\[ \mathcal P' \in {\tt{Lac}}_q(M, C_0-1) \text{ for all subsets } \mathcal P' \subsetneq \mathcal P_p, \; \#(\mathcal P') = q. \]
\end{lemma} 
\begin{proof} 
If $\varrho(p) \in \mathcal P'$, then \eqref{lacunarity under containment} yields the stronger conclusion $\mathcal P' \in {\tt{Lac}}_q(M, C_0)$. If $\varrho(p) \notin \mathcal P'$, let 
\[ \mathcal P' = \{ \varrho(i_q) < \varrho(i_{q-1}) < \cdots < \varrho(i_{1}) \}, \; 1 \leq i_1 < i_2 < \cdots < i_q < p. \]
The argument now follows the same lines as the proof of Case 2 in Lemma \ref{lemma: lacunary sequence vs tuple}. Specifically, for $1 \leq k < q-1$, 
\begin{equation} 
\rho(i_{k+1}) - \rho(i_q) \leq \rho(i_{k+1}) - \rho(p) \leq M^{-C_0(i_{k+1} - i_k)} \bigl[ \rho(i_{k}) - \rho(p) \bigr]. \label{kqp upper}
\end{equation}   
On the other hand, 
\begin{align}
\rho(i_{k}) - \rho(i_q) &= \bigl[ \rho(i_{k}) - \rho(p) \bigr] - \bigl[ \rho(i_{q}) - \rho(p) \bigr] \nonumber \\ 
& = \bigl[ \rho(i_{k}) - \rho(p) \bigr] \times  \left[ 1 - \frac{\rho(i_{q}) - \rho(p)}{\rho(i_{k}) - \rho(p)} \right] \nonumber \\
& \geq \bigl[ \rho(i_{k}) - \rho(p) \bigr] \left(1 - M^{-C_0(i_q - i_k)} \right). \label{kqp lower}
\end{align} 
Combining \eqref{kqp upper} and \eqref{kqp lower}, we arrive at 
\[ \frac{\rho(i_{k+1}) - \rho(i_q)}{\rho(i_{k}) - \rho(i_q)} \leq \frac{M^{-C_0(i_{k+1} - i_k)}}{\left(1 - M^{-C_0(i_q - i_k)} \right)} \leq \frac{M^{-C_0}}{1 - M^{-C_0}} \leq M^{-C_0+1}. \] 
Here the last inequality uses
\[ 1-M^{-C_0} \geq M^{-1}, \] which follows from $C_0>1$ and $M\geq 2$. This completes the proof. 
\end{proof} 
%The next lemma shows that any finite sub-collection of a sufficiently fast-decaying lacunary sequence $A$ is in ${\tt{Lac}}_p$ for some $p$.  
%\begin{lemma} 
%Given a strictly decreasing sequence 
%\[ A = \{a_j : j \geq 1\} \in {\tt{MonLac}}(M^{-C_0-1}), \quad A \subseteq [0, 1), \quad a = \lim_{j \rightarrow \infty} a_j,\]  
%
%\[ \] 
%\end{lemma} 
%\begin{proof} 
%
%\end{proof}
\section{Forbidden configurations of lacunary tuples} 
The assumption $\mathcal P_p \in {\tt{Lac}}_p(M, C_0)$ prevents the occurrence of certain branching configurations in the tree $\mathcal T(\mathcal P_p; M)$, thereby restricting its splitting number and the order of its tree lacunarity. This section identifies several such forbidden configurations and uses them to establish the structural conclusions about the tree $\mathcal T(\mathfrak A)$ stated in Proposition \ref{TREE STRUCTURE THEOREM}. 
\subsection{Tree for a lacunary tuple: progeny of a splitting vertex} 
\begin{lemma} \label{lemma: forbidden configuration 1} 
Suppose that $M \geq 4$ and the four rays corresponding to the points in 
\begin{equation} \label{linear ordering P4}  
\mathcal P_4 = \{\varrho(4) < \varrho(3) < \varrho(2) < \varrho(1)\} 
\end{equation}  
diverge from a single vertex; i.e. there exists $v \in \mathcal T(\mathcal P_4) = \mathcal T(\mathcal P_4; M)$ with
\begin{equation} \label{common ancestor of 4} 
v = D[i, i'] \; \text{ for } 1 \leq i < i' \leq 4.  
\end{equation} 
Then $\mathcal P_4 \notin {\tt{Lac}}_4(M, C_0)$ for any $C_0 > 1$.
\end{lemma} 
\begin{proof} 
As a consequence of the hypothesis \eqref{common ancestor of 4} of a common ancestor $v$, we obtain 
\[ \kappa = \kappa[i, i'] = h(v) + 1 \; \text{ for } 1 \leq i < i' \leq 4.\]
In view of the ordering \eqref{linear ordering P4}, Lemma \ref{lemma: distance x y} yields an ordering of the $\kappa^{\text{th}}$ digits of $\varrho(\cdot)$:  
\[ \rho_{\kappa}(4) < \rho_{\kappa}(3) < \rho_{\kappa}(2) < \rho_{\kappa}(1). \]   
Since each digit $\rho_{\kappa}(\cdot)$ lies in $\mathbb Z_{M}$, the above chain of inequalities means 
\[ \rho_{\kappa}(2) - \rho_{\kappa}(4) \geq 2, \quad \rho_{\kappa}(1) - \rho_{\kappa}(4) \geq 2. \] 
The relation \eqref{distance scales summary} therefore implies 
\[ M^{-\kappa} \leq \varrho(2) - \varrho(4) \leq M^{-\kappa+1}, \quad M^{-\kappa} \leq \varrho(1) - \varrho(4) \leq M^{-\kappa+1}. \]
Combining the two relations, we obtain 
\[ \frac{\varrho(2) - \varrho(4)}{\varrho(1) - \varrho(4)} \geq \frac{M^{-\kappa}}{M^{-\kappa+1}} = M^{-1}.\]  
This is contrary to the lacunarity requirement 
\[ \bigl[ \varrho(2) - \varrho(4) \bigr] \leq M^{-C_0} \bigl[ \varrho(1) - \varrho(4)\bigr] \text{ if } C_0 > 1, \] 
and completes the proof. 
\end{proof} 
\begin{corollary} 
Suppose that $M \geq 4$, $C_0 > 1$, and \[ \mathfrak A \in {\tt{MonLac}}(M^{-C_0-1})\] 
 is a sequence decreasing to $\alpha$. Then every vertex of $\mathcal T(\mathfrak A; M)$ lying on $\mathcal R^{\ast} = \mathcal R(\alpha)$ has at most three children. 
\end{corollary} 
\begin{proof} 
Suppose that a vertex $v$ has at least four children in $\mathcal T(\mathfrak A; M)$. Let us select four of these children and, from each one, choose an element of $\mathfrak A$ whose representing ray passes through that child. Indexing the resulting points as $\{\varrho(i) : 1 \leq i \leq 4 \}$ in decreasing order, we are reduced to the situation of Lemma \ref{lemma: forbidden configuration 1}. The lemma allows us to deduce that 
 \begin{equation}  \label{not in Lac}
 \mathcal P_4 \notin {\tt{Lac}}_4(M, C_0) \text{ for } C_0 > 1. 
 \end{equation}
 At the same time, $\mathcal P_4 \subseteq \mathfrak A$ while $\mathfrak A \in {\tt{MonLac}}(M^{-C_0-1})$. Lemma \ref{lemma: lacunary sequence vs tuple} therefore implies that 
 \begin{equation} \label{in Lac}  
 \mathcal P_4 \in {\tt{Lac}}_4(M, C_0). 
 \end{equation} 
 The contradiction of \eqref{not in Lac} and \eqref{in Lac} completes the proof.   
\end{proof} 
\subsection{Tree for a lacunary triple: digit following a separating height}
\begin{lemma}  \label{lemma: forbidden configurations - digit restrictions} 
Let us consider an ordered triple 
\[ \mathcal P_3 = \{\varrho(3) < \varrho(2) < \varrho(1)\} \text{ as in \eqref{ordered p tuple}}, \] with the following properties:
\begin{align}  
D[1, 2] &\subsetneq D[1, 3] = D[2, 3], 
\label{123 ancestors}  \\
\rho_{\kappa+1}(2) - \rho_{\kappa + 1}(3) &> -(M-1), \text{ where } \kappa = \kappa[1, 3] = \kappa[2, 3]. 
\label{<M+1}
\end{align}  
Here $\kappa[i, i']$ and $\rho_k(i)$ denote the first separating height and the $k^{\text{th}}$ $M$-adic digit, as defined in \eqref{nuii'} and \eqref{rhok}. Then 
\begin{equation}  \label{P3notlac}
\mathcal P_3 \not\in {\tt{Lac}}_3(M, C_0) \quad \text{ for any } C_0 > 2.
\end{equation}  
\end{lemma} 
\begin{proof} 
We will prove the assertion \eqref{P3notlac} by establishing the following distance relation: 
\begin{equation} \label{distance relation 13 23}
M^{-\kappa-1} \leq \varrho(i) - \varrho(3) \leq M^{-\kappa+1} \text{ for } i = 1,2. 
\end{equation} 
Indeed, \eqref{distance relation 13 23} implies 
\[ \frac{\varrho(2) - \varrho(3) }{\varrho(1) - \varrho(3) } \geq \frac{M^{- \kappa-1}}{M^{-\kappa+1}} = M^{-2},\]
which contradicts  the defining criterion of a $(M, C_0)$ lacunary triple, namely
\[ \frac{\varrho(2) - \varrho(3) }{\varrho(1) - \varrho(3) } \leq M^{-C_0}, \quad \text{ provided } C_0 > 2.  \]
\vskip0.1in
\noindent Henceforth, we will focus on proving \eqref{distance relation 13 23}. Fix an index $i\in\{1,2\}$, and consider two cases depending on the value of $\rho_{\kappa}(i) - \rho_{\kappa}(3)$. Lemma \ref{lemma: Euclidean vs Madic distance} is invoked in each case to estimate the size of $\varrho(i) - \varrho(3)$. 
\vskip0.1in
\noindent {\em{Case 1: }} First suppose $\rho_{\kappa}(i) - \rho_{\kappa}(3) \geq 2$. It then follows from  Lemma \ref{lemma: Euclidean vs Madic distance} \eqref{lemma: Euclidean vs Madic distance parta} that 
\begin{equation} \label{Euclidean distance in case 1}
M^{-\kappa} \leq \varrho(i) - \varrho(3) \leq M^{-\kappa+1}.
\end{equation} 
This is stronger than the claim \eqref{distance relation 13 23}.
\vskip0.1in
\noindent {\em{Case 2: }} Next let us consider the complementary case  
\begin{equation} \label{difference 1} 
\rho_{\kappa}(i) - \rho_{\kappa}(3) = 1.
\end{equation}  
If \eqref{difference 1} holds with $i=2$, then the assumption \eqref{<M+1}, combined with the definition \eqref{def: mu>nu} of $\mu$, leads to the relation  
\begin{equation} \label{mu23} 
\mu[2, 3] =  \kappa+1. 
\end{equation} 
The distance property given by \eqref{distance scales summary} then implies 
\[ M^{-\kappa-1} = M^{-\mu[2,3]}\leq \varrho(2) - \varrho(3) \leq M^{-\mu[2,3]+1} \leq M^{-\kappa}.   \]
This agrees with the stated bound \eqref{distance relation 13 23} for $i=2$. 
\vskip0.1in
\noindent On the other hand, if \eqref{difference 1} holds with $i=1$, the assumption \eqref{123 ancestors} implies 
\[ \kappa[1, 2] = h(D[1,2]) +1 > h(D[1, 3]) + 1 =  \kappa. \]
In other words, $\varrho(1)$ and $\varrho(2)$ separate at a height strictly larger than $\kappa$. In view of the ordering relation $\varrho(2) < \varrho(1)$, Lemma \ref{lemma: distance x y} yields $\rho_{\kappa+1}(1) \geq \rho_{\kappa+1}(2)$. 
Combining this with the assumption \eqref{<M+1}, we deduce  
\[\rho_{\kappa+1}(1) - \rho_{\kappa+1}(3) \geq \rho_{\kappa+1}(2) - \rho_{\kappa+1}(3) > -(M-1). \] 
In other words,  
\begin{equation} 
\mu[1, 3] = \kappa+1. \label{mu13}
\end{equation} 
Once again, the distance relation \eqref{distance scales summary} yields the desired estimate:
\begin{equation}  \label{Euclidean distance in case 2}
M^{-\kappa-1} = M^{-\mu[1,3]} \leq \varrho(1) - \varrho(3) \leq M^{-\mu[1,3]+1} = M^{-\kappa}.
\end{equation}  
This concludes the verification of \eqref{distance relation 13 23} in all cases. 
\end{proof} 

\subsection{Vertices of $\mathcal R^{\ast}$ that generate non-trivial subtrees of $\mathcal T(\mathfrak A)$} 
\begin{corollary} \label{corollary: digit after sep height} 
Suppose that 
\[ \mathfrak A \subseteq [0, 1), \quad \mathfrak A \in {\tt{MonLac}}(M^{-C_0-1}) \] is an infinite sequence decreasing to a limit $\alpha < 1$. Then $\mathcal T =  \mathcal T \bigl( \mathfrak A ; M\bigr)$ obeys the following property. 
\vskip0.1in
\noindent Let $\mathcal T_w = \mathcal T(\mathfrak A_w; M)$ denote the maximal subtree of $\mathcal T$ representing the set $\mathfrak A_w \subseteq \mathfrak A$, and rooted at a vertex 
\[ w \in \mathcal W = \mathcal W(\mathcal T, \mathcal R^{\ast}) \text{ where } \; \mathcal R^{\ast} = \mathcal R(\alpha) \text{ and } \mathcal W \text{ is defined in \eqref{what is W(T)}}.  \]
\begin{enumerate}[(a)]
\item \label{corollary: digit after sep height parta} If $\#(\partial \mathcal T_w) \geq 2$, then all rays of $\partial \mathcal T_w$, except possibly the one corresponding to $\max \mathfrak A_w$, have a common edge $e$ between heights $h(w)$ and $h(w)+1$.
\vskip0.1in 
\item \label{corollary: digit after sep height partb} The edge $e$ moves maximally to the left; i.e. $e$ joins $w$ to its leftmost child. 
\vskip0.1in
\item \label{corollary: digit after sep height partc} The corresponding edge $e'$ of $\mathcal R^{\ast}$ between the same heights $h(w)$ and $h(w)+1$ moves maximally to the right. Specifically, this means 
\begin{equation}  \label{<M-1, M-1} 
\alpha_{h(w)} < M-1, \quad \alpha_{h(w)+1} = M-1. 
\end{equation}  
\end{enumerate}
\end{corollary} 
\begin{proof} 
Given a vertex $w$ as in the hypothesis, let us consider the triple 
\[ \mathcal P_3 = \{ \varrho(3) < \varrho(2) < \varrho(1)\}, \quad \varrho(3) = \alpha, \; \varrho(2) = a_{i_2}, \; \varrho(1) = a_{i_1}.  \]
Here $a_{i_2} < a_{i_1}$ are any two elements of $\mathfrak A$ contained in $w$. Then, on one hand, Lemma \ref{lemma: lacunary sequence vs tuple} posits that 
\[ \mathcal P_3 \in {\tt{Lac}}_3(M, C_0). \]  
On the other, $\mathcal P_3$ satisfies the assumption \eqref{123 ancestors} of Lemma \ref{lemma: forbidden configurations - digit restrictions}: 
\[ w = D[1, 2] \subsetneq v = D[1, 3] = D[2, 3] \in \mathcal R^{\ast}, \; h(w) = h(v) + 1.\] 
This means that the second part of the assumption of Lemma \ref{lemma: forbidden configurations - digit restrictions}, namely \eqref{<M+1}, must be false. Hence, in the notation of that lemma, 
\[ \rho_{\kappa+1}(2) - \rho_{\kappa+1}(3) = -(M-1), \text{ with } \kappa = h(w).\]
Since $\rho_k(\cdot) \in \mathbb Z_M$, the equality above holds if and only if 
\begin{equation}  \label{digits - step1} \rho_{\kappa+1}(2) = 0 \; \text{ and } \; \rho_{\kappa+1}(3) = M-1. \end{equation} 
Moreover, since $\varrho(3)=\alpha<\varrho(2)$ and these two numbers first separate at height $\kappa=h(w)$, Lemma \ref{lemma: distance x y} gives
\begin{equation} \label{digits - step2}
\alpha_{h(w)} =\rho_{\kappa}(3)< \rho_{\kappa}(2) \leq M-1.
\end{equation} 
The second equality in \eqref{digits - step1}, together with \eqref{digits - step2}, proves \eqref{<M-1, M-1} and hence 
part \eqref{corollary: digit after sep height partc}. 
\vskip0.1in 
\noindent The first equality in \eqref{digits - step1} means that the edge $e$ on $\mathcal R(a_{i_2})$ originating at $w$ must move to the maximum allowable left. The second one means that  the corresponding edge $e'$ on $\mathcal R^{\ast}$ between the same heights moves to the maximum allowable right. The argument above holds for all elements $a_{i_2} \in \mathfrak A_w$ except the largest one, namely $\max \mathfrak A_w$. This shows that with the possible exception of this one element in $\mathfrak A_w$, the edge $e$ is common to all rays in $\partial \mathcal T_w$. This proves both parts \eqref{corollary: digit after sep height parta} and \eqref{corollary: digit after sep height partb}, and thus completes the proof.    
\end{proof} 
\noindent As a result of Corollary \ref{corollary: digit after sep height}, one can identify the vertices $v \in \mathcal R^{\ast}$ whose children $w \notin \mathcal R^{\ast}$ can generate trees $\mathcal T_w$ consisting of at least two rays.  
\begin{corollary} \label{corollary: nontrivial splitting heights} 
For $\mathfrak A$ and $\alpha = \sum_k \alpha_k M^{-k}$ as in Corollary \ref{corollary: digit after sep height}, let us define
\begin{align}  
\mathcal K^{\ast} &= \mathcal K^{\ast}(\mathcal R^{\ast}) := \left\{k \geq 1 : \alpha_{k} < M -1 \; \text{ but } \; \alpha_{k+1} = M -1 \right\}, \nonumber \\ 
 \mathcal W_{2+} &= \mathcal W_{2+}(\mathcal T, \mathcal R^{\ast}) := \left\{ w \in \mathcal W(\mathcal T, \mathcal R^{\ast}) : \# \bigl( \partial \mathcal T_w \bigr) \geq 2 \right\}. \label{def-W2+}
\end{align} 
Then 
\[ h(w)  \in \mathcal K^{\ast} \text{ for each } w \in \mathcal W_{2+}. \] In other words, only the vertices of $\mathcal R^{\ast}$ that occur at a height lying in $\mathcal K^{\ast}-1$ can produce a child $w \in \mathcal W^{\ast}$ containing more than one element of $\mathfrak A$. 
\end{corollary}  
\begin{proof}
This is a direct consequence of Corollary \ref{corollary: digit after sep height} \eqref{corollary: digit after sep height partc}. The details are left to the reader.  
\end{proof} 
%\begin{proof} 
%
%\end{proof} 
%\begin{corollary} 
%
%\end{corollary} 
\subsection{Subtrees rooted away from the minimal ray in a lacunary quadruple} 
So far, as a result of Corollaries \ref{corollary: digit after sep height}  and \ref{corollary: nontrivial splitting heights}, we have identified $\mathcal W^{\ast}$ as the vertices of $\mathcal T = \mathcal T(\mathfrak A)$ that can potentially generate non-trivial subtrees, i.e. subtrees with more than a single ray. In this section, we will discuss the structure of these subtrees. The following lemma distils this structure in the context of a lacunary quadruple. 
\vskip0.1in
\noindent Let $\mathcal P_4 \in {\tt{Lac}}_4(M, C_0)$ be an ordered quadruple as in \eqref{ordered p tuple} consisting of 
\begin{align} 
&\varrho(4) < \varrho(3) < \varrho(2) < \varrho(1);  \nonumber \\ 
&\text{suppose that } D[1, 4] = D[2, 4] = D[3, 4]. \label{quadruplet assumption 1} 
\end{align} 
 Let us denote by $v$ the vertex that simultaneously serves as the youngest common ancestor of the three pairs $(\varrho(i), \varrho(4))$, $i=1, 2,3$. Clearly, $v$ lies on $\mathcal R[4]$. The assumption \eqref{quadruplet assumption 1} of course means that the first separating heights of these three pairs are equal, i.e. 
\[ \kappa[1, 4] = \kappa[2, 4] = \kappa[3, 4], \; \text{ which we will denote as $\kappa$}.  \]
We will assume that $v$ has a child $w \notin \mathcal R[4]$ such that
\begin{equation}  \label{quadruplet assumption 2}
\mathcal P_3 \subseteq w, \; \text{ where } \; \mathcal P_3 := \{ \varrho(3) < \varrho(2) < \varrho(1) \}. 
\end{equation}  
This means that the maximal subtree $\mathcal T_w(\mathcal P_4) \subseteq \mathcal T(\mathcal P_4)$ rooted at $w$ consists of three rays, namely those identifying $\mathcal P_3$. Let us denote by $\mathcal R_w[i]$ the maximal ray of $\mathcal T_w(\mathcal P_4)$ representing $\varrho(i)$, $1 \leq i \leq 3$. 
\vskip0.1in
\noindent The next lemma asserts that for $C_0$ sufficiently large and for a lacunary quadruple $\mathcal P_4 \in {\tt{Lac}}_4(M, C_0)$, 
\[ 
\left\{ 
\begin{aligned} 
&\text{the ray $\mathcal R_w[3]$  corresponding to the smallest element $\varrho(3) \in \mathcal P_3$} \\ 
& \text{ is the special ray for $\mathcal T_w(\mathcal P_4)$ in the sense of Lemma \ref{SplitsOnARay}}. 
\end{aligned}
\right\} \] 
Namely, all splitting vertices of $\mathcal T_w(\mathcal P_4)$ necessarily lie on this ray.  
\begin{lemma} \label{lemma: behaviour of maximal subtrees} 
For $C_0 > 3$ and $\mathcal P_4 \in {\tt{Lac}}_4(M, C_0)$ obeying the above setup,  
\begin{equation} \label{triplet ancestors 0}
\text{every vertex of $\mathcal T_w(\mathcal P_4)$ not on } \mathcal R_w[3] \text{ is non-splitting. } 
\end{equation} 
In view of the ordering relation $\varrho(2) < \varrho(1)$, this is equivalent to 
\begin{equation} \label{triplet ancestors} 
D[2, 3] \subseteq D[1, 3] = D[1, 2]. 
\end{equation}   
\end{lemma} 
\begin{proof} 
Aiming for a contradiction, suppose that \eqref{triplet ancestors 0}, and equivalently \eqref{triplet ancestors}, does not hold. This means $\mathcal T_w(\mathcal P_4)$ must have at least one splitting vertex $w_0$ not on $\mathcal R_w[3]$. Since the tree $\mathcal T_w(\mathcal P_4)$ has three terminal vertices, and therefore can have at most two splitting vertices, the only possibility for $w_0$ is $D[1,2]$. The negation of \eqref{triplet ancestors} therefore implies 
\begin{equation}  \label{123 ancestors revisited} 
D[1, 2] \subsetneq D[1, 3] = D[2,3]. 
\end{equation}  
We will consider two cases, depending on the relative sizes of $\mu[3, 4]$ and $\kappa[2, 3]$. 
\vskip0.1in
\noindent {\em{Case 1: $\mu[3, 4] \leq \kappa[2, 3]+1$. }}
This relation along with \eqref{123 ancestors revisited} gives rise to the chain of inequalities:  
\[ \mu[3, 4] \leq \kappa[2, 3] + 1 \leq \kappa[1,2].\] 
On one hand, the ordering relation $\varrho(2) > \varrho(3) > \varrho(4)$, together with \eqref{distance scales summary}, provides a lower bound for $\varrho(2) - \varrho(4)$:
\begin{equation} 
\varrho(2) - \varrho(4) > \varrho(3) - \varrho(4) \geq M^{-\mu[3, 4]}. \label{24 lower}
\end{equation} 
On the other hand, Lemma \ref{lemma: distance x y} can be used together with \eqref{distance scales summary} to deduce an upper bound on $\varrho(1)-\varrho(4)$:
\begin{align}  
\varrho(1) - \varrho(4) &= \bigl[ \varrho(1) - \varrho(2)] +  \bigl[ \varrho(2) - \varrho(3)] + \bigl[ \varrho(3) - \varrho(4)] \nonumber \\ 
&\leq M^{-\kappa[1,2]+1} + M^{-\kappa[2,3]+1} + M^{-\mu[3, 4]+1} \nonumber \\
&\leq M^{-\mu[3,4]+1} \left[1 + M^{\mu[3,4] - \kappa[1,2]} + M^{\mu[3,4]-\kappa[2,3]} \right] \nonumber \\ 
&\leq   M^{-\mu[3,4]+1}(2+M) \leq M^{-\mu[3,4]+3}. \label{14 upper}
\end{align} 
Combining \eqref{24 lower} and \eqref{14 upper}, we obtain 
\begin{equation} \label{124 ratio lower}
\frac{\varrho(2) - \varrho(4)}{\varrho(1) - \varrho(4)}\geq M^{-3}. 
\end{equation} 
However, $\mathcal P_3' = \{\varrho(4) < \varrho(2) < \varrho(1) \} \subseteq \mathcal P_4$, and therefore, by virtue of \eqref{lacunarity under containment}, 
\begin{equation}  \label{124 ratio upper}
\mathcal P_3' \in  {\tt{Lac}}_3(M, C_0) \; \text{ i.e., } \; \frac{\varrho(2) - \varrho(4)}{\varrho(1) - \varrho(4)}\leq M^{-C_0}. 
\end{equation}    
The two conclusions \eqref{124 ratio lower} and \eqref{124 ratio upper} contradict each other for $C_0 > 3$, completing the proof for Case 1. 
\vskip0.1in
\noindent {\em{Case 2: $\mu[3, 4] \geq \kappa[2, 3]+2$. }}
In this case, we will argue that $\mathcal P_3$ given by \eqref{quadruplet assumption 2} obeys the requirements of Lemma \ref{lemma: forbidden configurations - digit restrictions}. Once this is established, the conclusion of the lemma says that $\mathcal P_3 \notin {\tt{Lac}}_3(M, C_0-1)$ for $C_0 > 3$. This would contradict 
\[ \mathcal P_3 \in {\tt{Lac}}_3(M, C_0-1) \text{ which follows from } \mathcal P_3 \subseteq \mathcal P_4 \text{ via } \eqref{lacunarity under containment}. \] 
The contradiction would then complete the argument in Case 2. 
\vskip0.1in
\noindent Thus it only remains to verify the hypotheses of Lemma \ref{lemma: forbidden configurations - digit restrictions} for $\mathcal P_3$. The first assumption \eqref{123 ancestors} of that lemma is precisely \eqref{123 ancestors revisited}. 
\vskip0.1in
\noindent Let us turn now to the verification of \eqref{<M+1}, which requires us to show that
\begin{equation*}  
\rho_{\kappa' + 1}(2) - \rho_{\kappa' + 1}(3) > -(M-1) \text{ with } \kappa' = \kappa[1, 3] = \kappa[2, 3]. 
\end{equation*}  
This is equivalent to showing 
\begin{equation} \label{pair equality} 
(\rho_{\kappa' + 1}(2), \rho_{\kappa' + 1}(3)) \neq (0, M-1).  
\end{equation} 
To confirm this, we appeal to the definition of $\mu$. Let us first observe that 
 \[ D[2, 3] \subseteq w = D[3, 4] \text{ implies } \kappa[2,3] \geq \kappa[3, 4]. \]
Therefore, in view of the defining constraint of Case 2,  
\[ \mu[3,4] \geq \kappa[2, 3]+2 > \kappa[3,4]. \]
The definition \eqref{distance scales summary} of $\mu[3, 4]$ then yields
\begin{equation}  \label{kappa34} 
\rho_{k}(3) = 0 \text{ and } \rho_k(4) = M-1 \text{ for all } \kappa[3, 4] +1 \leq k \leq \mu[3, 4]. 
\end{equation}  
Let us check that $k = \kappa'+1$ falls in the range of $k$ displayed in \eqref{kappa34}. For the lower bound, we note from \eqref{quadruplet assumption 1} and \eqref{quadruplet assumption 2} that  
\begin{align}  
&D[2, 3] \subseteq w \subsetneq v = D[3, 4]; \\
\label{check range 1} 
&\text{therefore } \kappa' + 1 = \kappa[2, 3] + 1 \geq \kappa[3, 4] + 1. 
\end{align}  
To find an upper bound on $\kappa'+1$, the defining constraint of this case says that 
\begin{equation} \kappa'+1 = \kappa[2, 3] + 1 \leq \mu[3, 4]-1.  \label{check range 2} 
\end{equation}  
Combining \eqref{check range 1} and \eqref{check range 2} with \eqref{kappa34}, we deduce that 
\[ \rho_{\kappa'+1}(3) = 0 \neq M-1, \text{ and therefore } (\rho_{\kappa' + 1}(2), \rho_{\kappa' + 1}(3)) \neq (0, M-1). \]
This verifies \eqref{pair equality}, and thus completes the proof. 
\end{proof}   
%Lemma \ref{lemma: forbidden configurations - digit restrictions} again, but this time on the triple $\mathcal P_3' := \{ \varrho(4) < \varrho(3) < \varrho(2) \}$. On one hand, 
%\[ \mathcal P_3' \subseteq \mathcal P_4 \text{ implies } \mathcal P_3' \in {\tt{Lac}}_3(M, C_0) \text{ by virtue of \eqref{lacunarity under containment}}.   \]
%On the other hand, the assumptions \eqref{quadruplet assumption 1} and \eqref{quadruplet assumption 2} on $\mathcal P_4$ imply
%\[ D[2, 3] \subseteq w \subsetneq v = D[2, 4] = D[3, 4]. \] 
%This is the analogue of the assumption \eqref{123 ancestors} in Lemma \ref{lemma: forbidden configurations - digit restrictions} for the triple $\mathcal P_3'$.
%\vskip0.1in   
%\noindent Applying the lemma to the triple $\mathcal P_3'$ then implies 
%\[ \rho_{\kappa+1}(3) = 0, \quad \rho_{\kappa+1}(4) = M-1, \text{ where } \kappa = \kappa[2, 4] = \kappa[3, 4]. \]
%From this, we claim that 
%\[ \rho_{\kappa'+1}(2) > 0. \] 
%Indeed 
%\[ \rho_r(2) = \rho_r(3) \text{ for } r < \kappa[2, 3], \quad \rho_r(2) = \rho_r(3) \text{ for } r = \kappa[2, 3]. \] 
%\[ \kappa[2, 3] > h(w) > \kappa[2, 4].  \] 

\begin{lemma} 
In the setup of Lemma \ref{lemma: behaviour of maximal subtrees}, the last splitting vertex $D[2,3]$ has to occur at or before height $\mu[3,4]$, i.e.,
\begin{equation} \label{last split}  
h\bigl( D[2, 3]\bigr) \leq \mu[3,4]. 
\end{equation}  
\end{lemma} 
\begin{proof} 
Aiming for a contradiction, suppose if possible that 
\begin{equation} \label{contradictory assumption}   
h\bigl( D[2, 3]\bigr) = \kappa[2, 3] - 1 > \mu[3,4].
\end{equation}
From \eqref{lemma: distance x y} and \eqref{distance scales summary}, we have the two distance estimates 
\[ \varrho(2) - \varrho(3) \leq M^{-\kappa[2, 3] + 1}, \quad \varrho(3) - \varrho(4) \geq M^{-\mu[3,4]}. \] 
These estimates permit the bound
\begin{align*} 
\varrho(2) - \varrho(4) &= \bigl[ \varrho(2) - \varrho(3) \bigr] + \bigl[ \varrho(3) - \varrho(4) \bigr] \\ 
& = \bigl[ \varrho(3) - \varrho(4) \bigr] \left[ 1 + \frac{\varrho(2) - \varrho(3)}{\varrho(3) - \varrho(4) } \right] \\ 
&\leq  \bigl[ \varrho(3) - \varrho(4) \bigr] \left[ 1 + M^{-\kappa[2, 3]+ \mu[3, 4] + 1} \right] \\
&\leq 2\bigl[ \varrho(3) - \varrho(4) \bigr] \leq M \bigl[ \varrho(3) - \varrho(4) \bigr], 
\end{align*} 
where the penultimate inequality uses the assumption \eqref{contradictory assumption}. Thus 
\[ M^{-1} \leq \frac{\varrho(3) - \varrho(4) }{\varrho(2) - \varrho(4)}. \] For $C_0 > 1$, this contradicts the property 
$\mathcal P_3' = \{ \varrho(4) < \varrho(3) < \varrho(2) \} \in {\tt{Lac}}_3(M, C_0)$, which is a consequence of \eqref{lacunarity under containment} and asserts that 
  \[  \frac{\varrho(3) - \varrho(4) }{\varrho(2) - \varrho(4)} \leq M^{-C_0}. \]
  This completes the proof.  
\end{proof} 
\subsection{Implications of forbidden configurations on a lacunary tree} 
\begin{corollary} \label{corollary: tree structure} 
Let $\mathfrak A$ and $\mathcal T = \mathcal T(\mathfrak A)$ be as in Corollary \ref{corollary: digit after sep height}, and let $\mathcal T_w$ denote the maximal subtree of $\mathcal T$ rooted at $w$. Then the following conclusions hold: 
\begin{enumerate}[(a)]
\item \label{TREE STRUCTURE THEOREM parta} For each $w \in \mathcal W = \mathcal W(\mathcal T, \mathcal R^{\ast})$ as defined in \eqref{what is W(T)}, the subtree $\mathcal T_{w} \subseteq \mathcal T$ consists of finitely many rays, i.e.,
\begin{equation} 
n_w := \# \bigl( \partial \mathcal T_w \bigr) < \infty \quad \text{ for all } w \in \mathcal W. 
\end{equation} 
The cardinality $\#(\partial \mathcal T_w)$ may grow without bound as $w$ ranges over $\mathcal W$. 
\vskip0.1in 
\item \label{TREE STRUCTURE THEOREM partd} (The subset of $\mathfrak A$ depicted by $\mathcal T_w$) The trees $\{\mathcal T_w: w \in \mathcal W\}$ represent disjoint blocks of points in $\mathfrak A$. 
\vskip0.1in
\noindent More precisely, let us define the finite sets $\mathfrak A_{w} \subseteq \mathfrak A$ and $\mathbb N_w \subseteq \mathbb N$ as follows,
\begin{equation} 
\mathcal T_{w} := \mathcal T(\mathfrak A_w; M), \quad \mathfrak A_{w} := \bigl\{a_j : j \in \mathbb N_w \bigr\}. 
\end{equation}
Then $\mathbb N_w$ is a block of contiguous indices of the form
\begin{equation} \label{contiguous}
\mathbb N_w = \{ j_{w} + 1, \ldots, j_w + n_w \} \; \text{ with } \; \mathbb N := \bigsqcup \bigl\{ \mathbb N_w : w \in \mathcal W(\mathcal T, \mathcal R^{\ast}) \bigr\}.   
\end{equation} 
This means $\mathfrak A_{w}$ is a finite block of consecutive elements of $\mathfrak A$.
\item \label{TREE STRUCTURE THEOREM partc} For each $w \in \mathcal W$ as defined in \eqref{what is W(T)}, one of the following conclusions holds: 
\begin{equation}  \label{two types of w} 
\mathcal T_w \in \begin{cases} {\tt{TreeLac}}(0; M) &\text{ if } w \notin \mathcal W_{2+}, \\ {\tt{TreeLac}}(1; M) \setminus {\tt{TreeLac}}(0; M) &\text{ if } w \in \mathcal W_{2+}. \end{cases} \end{equation} 
Here the class of trees {\tt{TreeLac}} is specified by Definition \ref{defn: Tree Lacunary sets} and the definition of $\mathcal W_{2+}$ appears in \eqref{def-W2+}.  
%As a result, 
%\[ \bigl[ \min \mathfrak A_w, \max \mathfrak A_w \bigr] \subseteq w.  \] 
\vskip0.1in
\item \label{TREE STRUCTURE THEOREM partc'}  The infinite ray $\mathcal R^{\ast}_w$ rooted at $w$ representing the extremal point $a_{j_w + n_w}$ of $\mathfrak A_w$ (according to the linear ordering on $\mathbb R$) is a distinguished ray of $\mathcal T_w$.  \label{second-order tree} 
\end{enumerate}  
\end{corollary}
\noindent Figure \ref{fig:two-level-spine} depicts the two-level spine structure of $\mathcal T(\mathfrak A)$, as described in Corollary \ref{corollary: tree structure}. 
\begin{figure}[ht]
\centering
\begin{tikzpicture}[
    x=1.05cm,
    y=1cm,
    every node/.style={font=\small},
    principal/.style={line width=1.4pt,blue!70!black},
    secondary/.style={line width=1.4pt,green!45!black},
    branch/.style={line width=.8pt},
    ray/.style={line width=.8pt,dashed},
    ordinary/.style={
        circle,draw,fill=white,
        inner sep=1.7pt,line width=.7pt
    },
    splitting/.style={
        circle,draw,fill=black,
        inner sep=2pt
    }
]

% Principal distinguished ray: the leftmost ray
\node[ordinary]  (r0) at (0,5.0) {};
\node[splitting] (r1) at (0,4.3) {};
\node[splitting] (r2) at (0,3.6) {};
\node[splitting] (r3) at (0,2.9) {};
\node[ordinary]  (r4) at (0,2.2) {};
\node[ordinary]  (r5) at (0,1.5) {};

\draw[principal]
    (r0)--(r1)--(r2)--(r3)--(r4)--(r5);
\draw[principal,dashed] (r5)--(0,0.65);

\node[blue!70!black,left=7pt] at (0,2.75)
    {$\mathcal R^\ast$};

% An earlier branch: farther from R^*
\node[ordinary] (wp0) at (4.8,3.6) {};
\node[ordinary] (wp1) at (4.8,2.9) {};
\node[ordinary] (wp2) at (4.8,2.2) {};

\draw[branch] (r1)--(wp0);
\draw[branch] (wp0)--(wp1)--(wp2);
\draw[ray] (wp2)--(4.8,1.35);

\node[above=4pt] at (wp0)
    {$w'\in\mathcal W\setminus\mathcal W_{2+}$};
\node[right=5pt] at (4.8,2.55)
    {$\mathcal T_{w'}$};

% A later branch: closer to R^*
\node[ordinary]  (w)  at (1.7,2.9) {};
\node[splitting] (s1) at (1.7,2.2) {};
\node[splitting] (s2) at (1.7,1.5) {};
\node[ordinary]  (s3) at (1.7,0.8) {};

\draw[branch] (r2)--(w);
\draw[secondary]
    (w)--(s1)--(s2)--(s3);
\draw[secondary,dashed] (s3)--(1.7,-0.05);

\node[above=4pt] at (w)
    {$w\in\mathcal W_{2+}$};

\node[green!45!black,left=7pt] at (1.7,1.75)
    {$\mathcal R_w^\ast$};

% Rays leaving R_w^*: later branches lie closer to R_w^*
\node[ordinary] (q1) at (3.9,1.5) {};
\node[ordinary] (q2) at (2.9,0.8) {};

\draw[branch] (s1)--(q1);
\draw[branch] (s2)--(q2);

\draw[ray] (q1)--(3.9,0.65);
\draw[ray] (q2)--(2.9,-0.05);

\node[right=5pt] at (2.2,1.4)
    {$\mathcal T_w$};

% Point represented by the secondary distinguished ray
\node[green!45!black,below=3pt,align=center]
    at (1.7,-0.05)
    {\scriptsize extremal point\\[-2pt]
     \scriptsize of $\mathfrak A_w$};

\end{tikzpicture}

\caption{\small{The ordered two-level tree structure described in
Corollary \ref{corollary: tree structure}. The principal distinguished ray $\mathcal R^\ast$, representing the limit point $\alpha=\inf\mathfrak A$, is the leftmost ray. Rays leaving $\mathcal R^\ast$ lie to its right and approach it as their splitting vertices move down $\mathcal R^\ast$. If $w\in\mathcal W_{2+}$, then $\mathcal T_w$ has the analogous structure: its distinguished ray $\mathcal R_w^\ast$ is leftmost, and all splitting vertices of $\mathcal T_w$ lie on $\mathcal R_w^\ast$. Filled vertices are splitting; open vertices are non-splitting.}}
\label{fig:two-level-spine}
\end{figure}

\begin{proof} 
For part \eqref{TREE STRUCTURE THEOREM parta}, it suffices to consider vertices $w \in \mathcal W_{2+}$ only, where $\mathcal W_{2+}$ is defined in \eqref{def-W2+}. Since $w$ is an $M$-adic interval not containing the sole limit point $\alpha = \lim \mathfrak A$, the number of points in $\mathfrak A \cap w$ is finite. In the notation of Corollary \ref{corollary: digit after sep height}, this means that the set $\mathfrak A_w$ is finite, which translates to its representative tree $\mathcal T_w$ having finitely many rays. 
\vskip0.1in
\noindent Turning to part \eqref{TREE STRUCTURE THEOREM partd}, suppose that $w, w'$ are two distinct elements of $\mathcal W$ such that
\[ w \subsetneq v, \; w' \subsetneq v', \; v, v' \in \mathcal R^{\ast}, h(w) = h(v) + 1 \leq h(v') +1 = h(w').  \]
Thus neither of the $M$-adic intervals $w, w'$ can contain the other, and therefore they must be disjoint. Being disjoint intervals, they are linearly ordered. 
%This means that every element of $\mathfrak A_{w}$ separates from $\alpha$ sooner than any element of $\mathfrak A_{w'}$; in other words,  
%\begin{align*} 
%&\kappa(a_{j}, \alpha)  < \kappa(a_{j'}, \alpha) \; \text{ for all $j, j'$ such that } a_j \in \mathfrak A_w, \; a_{j'} \in \mathfrak A_{w'};  \\  
%&\text{Lemma \ref{lemma: fsh transitive} then implies } a_{j} > a_{j'}.
%\end{align*}   
In other words, every element of $\mathfrak A_w = \{ a_j : j \in \mathbb N_w \} $ is strictly larger than that of $\mathfrak A_{w'}= \{ a_{j'} : j' \in \mathbb N_{w'}\} $, or conversely. Since $\mathfrak A$ is strictly decreasing, the associated index sets $\mathbb N_w$ and $\mathbb N_{w'}$ are correspondingly ordered: $\mathbb N_{w} \cap \mathbb N_{w'} = \emptyset$, with  
\begin{equation} \label{ordering of N_w and N_w'}
\; \text{ either } \max \mathbb N_w < \min \mathbb N_{w'} \; \; \text{ or } \; \; \max \mathbb N_{w'} < \min \mathbb N_{w}. 
\end{equation}
This means each $\mathbb N_w$ has to be of the form \eqref{contiguous}, which completes the proof of part \eqref{TREE STRUCTURE THEOREM partd}. It is worth mentioning in passing that if $h(w) < h(w')$, i.e. if $w, w'$ are not children of the same vertex $v \in \mathcal R^{\ast}$, then 
\begin{align*} 
&D(a_{j'}, \alpha) \subsetneq D(a_j, \alpha), \text{ and therefore } \\  
&\kappa(a_{j'}, \alpha) > \kappa(a_j, \alpha) \text{ for all } j \in \mathbb N_{w}, j' \in \mathbb N_{w'}. 
\end{align*}
Lemma \ref{lemma: fsh transitive} then gives $a_j > a_{j'}$, leading to the first alternative \eqref{ordering of N_w and N_w'}. 
\vskip0.1in
\noindent Next, we will prove parts \eqref{TREE STRUCTURE THEOREM partc} and \eqref{TREE STRUCTURE THEOREM partc'} together. If $w \notin \mathcal W_{2+}$, then $\mathcal T_w$ consists of a single ray; in other words,
\begin{equation}  \mathcal T_w \in {\tt{TreeLac}}(0; M) \text{ for } w \notin \mathcal W_{2+}. \label{not in W2+}
\end{equation} 
In this case, $\mathcal R_w^*$ is the unique ray of $\mathcal T_w$, so part \eqref{TREE STRUCTURE THEOREM partc'} also holds. 
\vskip0.1in
\noindent For $w \in \mathcal W_{2+}$, the tree $\mathcal T_w$ contains more than a single ray by definition of $\mathcal W_{2+}$, therefore 
$\mathcal T_w \notin {\tt{TreeLac}}(0;M)$. Let 
\[ \alpha_w := \min \mathfrak A_w \; \text{ and let } \; \mathcal R_w^{\ast} \text{ denote the sub-ray of $\mathcal R(\alpha_w)$ rooted at $w$.} \] 
We will show momentarily that 
\begin{equation} \label{show TreeLac1} 
\text{every vertex of $\mathcal T_{w}$ not lying on $\mathcal R_w^{\ast}$ is non-splitting.}
\end{equation} 
This will establish $\mathcal R_{w}^{\ast}$ as a distinguished ray of $\mathcal T_w$, with
\begin{equation} 
\mathcal T_w \in {\tt{TreeLac}}(1; M) \setminus {\tt{TreeLac}}(0;M) \text{ for } w \in \mathcal W_{2+},  \label{in W2+} 
\end{equation} 
in view of Definition \ref{defn: Tree Lacunary sets}. 
Jointly, \eqref{not in W2+} and \eqref{in W2+} prove \eqref{two types of w}, which is the conclusion of part \eqref{TREE STRUCTURE THEOREM partc}.
\vskip0.1in
\noindent Let us therefore focus on proving \eqref{show TreeLac1}. If $n_w=2$, then $\mathcal T_w$ has exactly one splitting vertex, namely the youngest common ancestor of its two rays. This vertex lies on $\mathcal R_w^*$, so \eqref{show TreeLac1} follows immediately. We may therefore assume that $n_w\geq 3$. Given any two elements  
\[ a_j, a_k \in \mathfrak A_{w} \setminus \{\alpha_w\}, \; j < k, \] we consider the lacunary quadruple 
\begin{align*} 
&\mathcal P_4 = \{ \varrho(4) < \varrho(3) < \varrho(2) < \varrho(1)\}, \; \text{ where } \\ 
&\varrho(4)= \alpha, \; \varrho(3) = \alpha_w, \; \varrho(2) = a_k, \; \varrho(1) = a_j.  
\end{align*}   
Then $\mathcal P_4$ satisfies the hypotheses of Lemma \ref{lemma: behaviour of maximal subtrees}, with 
\begin{align*} 
&v = \text{ parent of $w$ } = D[1, 4] = D[2, 4] = D[3, 4], \; \text{ and } \\ 
&\mathcal P_3 = \{\varrho(3) < \varrho(2) < \varrho(1) \} \subseteq w. 
\end{align*}  
According to the conclusions \eqref{triplet ancestors} and  \eqref{triplet ancestors 0} of the lemma, every vertex of $\mathcal T_w(\mathcal P_4) \setminus \mathcal R_w[3]$ is non-splitting, with 
\[ D \bigl(a_k, \alpha_w \bigr) = D[2, 3] \subseteq D[1,3] = D \bigl(a_j, \alpha_{w} \bigr)  = D[1, 2] = D \bigl(a_j, a_k \bigr). \] 
This in turn leads to the identity 
\[ \kappa(a_j, a_k) = \kappa(a_j, \alpha_w) \; \text{ for all } j_w+1 \leq j < k \leq j_w + n_w-1, \]
placing it in the set-up of Lemma \ref{lemma: order and tree} \eqref{lemma: order and tree partb}. By the conclusion of that lemma, all the splitting vertices of $\mathcal T_w$ lie on $\mathcal R_w^{\ast}$. This establishes \eqref{show TreeLac1} and concludes the proof. 
\end{proof} 
\section{Tree structure of lacunary sequences: Proof of Proposition \ref{TREE STRUCTURE THEOREM}} \label{Section: proof of TREE STRUCTURE THEOREM}  
\begin{proof}
Proving Proposition \ref{TREE STRUCTURE THEOREM} involves assembling the components that we have deduced thus far. Part \eqref{tree structure theorem part c} is the same as Corollary \ref{corollary: tree structure} \eqref{TREE STRUCTURE THEOREM parta}.   Corollary \ref{corollary: tree structure}  \eqref{TREE STRUCTURE THEOREM partc} says that for every vertex $w \in \mathcal W(\mathcal T(\mathfrak A), \mathcal R^{\ast})$,  
\[ \mathcal T_w \in {\tt{TreeLac}}(1; M) \] 
In view of Definition \ref{defn: Tree Lacunary sets}, this implies 
\[ \mathcal T(\mathfrak A) \in {\tt{TreeLac}}(2; M) \text{ with $\mathcal R^{\ast}$ as a distinguished ray.}  \]
This simultaneously establishes parts \eqref{TREE STRUCTURE THEOREM-lacunarity order} and \eqref{TREE STRUCTURE THEOREM-distinguished ray} of Proposition \ref{TREE STRUCTURE THEOREM}. The proof is therefore complete. 
\end{proof}

\chapter{Appendix E Part II: Gaps of a lacunary sequence via its tree} \label{section: proof of gap lemma}
\section{Chapter overview} 
In this chapter, we complete the proof of the gap lemma, namely Lemma \ref{GAP LEMMA}.
\vskip0.1in 
\noindent Given an infinite decreasing sequence
\[ \mathfrak A = \{a_j : j \geq 1\} \in {\tt{MonLac}}(M^{-C_0})  \]
as in Proposition \ref{TREE STRUCTURE THEOREM}, let us denote its $j^{\text{th}}$ gap by 
\begin{equation}
{\tt{Gap}}_j := \bigl[ a_{j+1}, a_j \bigr), \quad j \geq 1.
\end{equation}  
Whereas Proposition \ref{TREE STRUCTURE THEOREM} and Corollary \ref{corollary: tree structure} specify the tree structure of a monotone lacunary sequence $\mathfrak A$, our proof of Lemma \ref{GAP LEMMA} will identify the vertices of the full $M$-adic tree relative to $\mathcal T(\mathfrak A)$ that cover ${\tt{Gap}}_j$. The proof of Lemma \ref{GAP LEMMA} relies on two main pillars, which we summarize below as the main results of this chapter. 
\subsection{A stopping vertex lemma} The first one concerns the containment of gaps in certain vertices of the $M$-adic tree. 
\begin{lemma} \label{GAP LEMMA REVISITED} 
For every $j \geq 1$, 
\begin{equation}
{\tt{Gap}}_j \subseteq \bigcup \bigl\{ u^{\sharp} : u \in \mathcal U(\mathfrak A)\bigr\}, 
\end{equation} 
where $\mathcal U(\mathfrak A)$ is the collection of $M$-adic vertices defined in \eqref{def-UA}, and $u^{\sharp}$ is the canonical cell of $u$ given by \eqref{defn: canonical cell}.
\end{lemma}
\begin{proof}
Fix $j \geq 1$ and a point $x \in {\tt{Gap}}_j$. Since $\alpha = \lim \mathfrak A$ is the unique accumulation point of $\mathfrak A$ and $x \neq \alpha$, i.e., since every point of $\mathfrak A \setminus \{\alpha\}$ is isolated, there exists an open interval $I$ containing $x$ such that 
\[ \mathfrak A \cap I  = \left\{\begin{aligned} \{x\} &\text{ if } x \in \mathfrak A, \\ \emptyset &\text{ if } x \notin \mathfrak A \end{aligned} \right\}; \quad \text{in particular, } \#\bigl( \mathfrak A \cap I \bigr) \leq 1. \]
 Let us choose an $M$-adic ray representing $x$, and let 
 \[ v_0 \supset v_1 \supset v_2 \supset \cdots \]
 denote its successive vertices. Since the corresponding $M$-adic intervals shrink to $x$, there exists an integer $n$ such that 
 \[ v_n^{\sharp} \subset I; \; \text{ hence } \# \left( \mathfrak A \cap v_n^{\sharp} \right) \leq 1.  \]  
 Let $n_0$ be the smallest index for which this inequality holds. Since $\mathfrak A$ is infinite and $v_0 = [0,1]$, we deduce that $n_0 \geq 1$. Then 
 \[ \#\bigl(\mathfrak A \cap v_{n_0}^{\sharp} \bigr) \leq 1, \quad  \#\bigl(\mathfrak A \cap v_{n_0-1}^{\sharp} \bigr) \geq 2. \]  
 Thus $v_{n_0}$ is a vertex of minimal height along its lineage whose intersection with $\mathfrak A$ has cardinality at most one. Since 
 \[x \in v_{n_0}^{\sharp} \cap \bigl[ \inf \mathfrak A, \sup \mathfrak A\bigr], \text{ the characterization \eqref{U properties} gives } v_{n_0} \in \mathcal U(\mathfrak A). \] 
 We have therefore shown that every $x \in {\tt{Gap}}_j$ belongs to the canonical cell associated with some vertex of $\mathcal U(\mathfrak A)$. Consequently,
 \[ {\tt{Gap}}_j \subseteq \bigcup \bigl\{ u^{\sharp} : u \in \mathcal U(\mathfrak A)\bigr\}, \]
 as claimed.  
\end{proof} 
\subsection{Two-level decomposition of stopping vertices} 
Let us take a moment to recall the definition \eqref{def-UA} of the stopping vertex collection $\mathcal U(\mathfrak A)$, and its equivalent formulation \eqref{U properties}. In the next lemma, we find an explicit description of stopping vertices, grouping them according to the two levels of the spine structure. 
\vskip0.1in
\noindent Let us recall the definitions of $\mathcal W^{\ast}(\mathcal R^{\ast})$ and $\mathcal W_{2+}$ from \eqref{what is W-star} and \eqref{def-W2+} respectively. The members of $\mathcal W^{\ast}(\mathcal R^{\ast})$ are precisely the first-level off-spine vertices. Among them, some are roots of subtrees of $\mathcal T(\mathfrak A)$ that have a genuine secondary spine $\mathcal R_w^{\ast}$; these contribute to the collection $\mathcal W_{2+}$.  Next set
\begin{align}
&\mathcal U^{\ast}(\mathfrak A) := \left\{u \in \mathcal U(\mathfrak A) : u \in \mathcal W^{\ast}(\mathcal R^{\ast})  \right\}, \label{def-Ustar} \\ 
 &\mathcal U_w^{\ast}(\mathfrak A) := \left\{u \in \mathcal U(\mathfrak A) :  u \in \mathcal W^{\ast}(\mathcal R_w^{\ast})   \right\} \sqcup \{ u_w^{\ast}\}, \; w \in \mathcal W_{2+}, \label{def-Uwstar}
\end{align}
where $u_w^{\ast} \in \mathcal R_w^{\ast}$ is the child of the last splitting vertex on $\mathcal R_w^{\ast}$. A schematic diagram of the types of vertices defined in \eqref{def-Ustar} and \eqref{def-Uwstar} is given in Figure \ref{fig:stopping-vertex-types}.  
\begin{figure}[ht]
\centering
\begin{tikzpicture}[
    x=1cm,
    y=1cm,
    every node/.style={font=\small},
    principal/.style={line width=1.4pt,blue!70!black},
    secondary/.style={line width=1.4pt,green!45!black},
    branch/.style={line width=.8pt},
    continuation/.style={line width=.8pt,dashed},
    ordinary/.style={
        circle,draw,fill=white,
        inner sep=1.7pt,line width=.7pt
    },
    splitting/.style={
        circle,draw,fill=black,
        inner sep=2pt
    },
    stopping/.style={
        rectangle,draw=red!75!black,
        fill=red!15,inner sep=2.3pt,line width=.9pt
    }
]

% ------------------------------------------------------------
% Panel 1: stopping vertex immediately off R^*
% ------------------------------------------------------------
\begin{scope}[xshift=0cm]

\node[ordinary]  (a0) at (0,2.8) {};
\node[splitting] (a1) at (0,2.1) {};
\node[ordinary]  (a2) at (0,1.4) {};
\node[ordinary]  (a3) at (0,0.7) {};

\draw[principal] (a0)--(a1)--(a2)--(a3);
\draw[principal,continuation] (a3)--(0,0.15);

\node[stopping] (u1) at (1.25,1.4) {};
\draw[branch] (a1)--(u1);
%\draw[continuation] (u1)--(1.25,0.65);

\node[blue!70!black,left=6pt] at (0,1.7)
    {$\mathcal R^\ast$};
\node[red!75!black,right=5pt] at (u1)
    {$u$};

\node[align=center] at (0.6,-0.35)
    {\textbf{(i)} Immediately off $\mathcal R^\ast$\\
     $u\in\mathcal U^\ast(\mathfrak A)$};

\end{scope}

% ------------------------------------------------------------
% Panel 2: stopping vertex immediately off R_w^*
% ------------------------------------------------------------
\begin{scope}[xshift=4.35cm]

\node[ordinary]  (b0) at (0,2.8) {};
\node[splitting] (b1) at (0,2.1) {};
\node[ordinary]  (b2) at (0,1.4) {};
\node[ordinary]  (b3) at (0,0.7) {};

\draw[secondary] (b0)--(b1)--(b2)--(b3);
\draw[secondary,continuation] (b3)--(0,0.15);

\node[stopping] (u2) at (1.25,1.4) {};
\draw[branch] (b1)--(u2);
%\draw[continuation] (u2)--(1.25,0.65);

\node[green!45!black,left=6pt] at (0,1.7)
    {$\mathcal R_w^\ast$};
\node[red!75!black,right=5pt] at (u2)
    {$u$};

\node[align=center] at (0.6,-0.35)
    {\textbf{(ii)} Immediately off $\mathcal R_w^\ast$\\
     $u\in\mathcal U(\mathfrak A)
          \cap\mathcal W^\ast(\mathcal R_w^\ast)$};

\end{scope}

% ------------------------------------------------------------
% Panel 3: u_w^* follows the last splitting vertex
% ------------------------------------------------------------
\begin{scope}[xshift=8.7cm]

\node[ordinary]  (c0) at (0,2.8) {};
\node[ordinary]  (c1) at (0,2.1) {};
\node[splitting] (c2) at (0,1.4) {};
\node[stopping]  (uw) at (0,0.7) {};

\draw[secondary] (c0)--(c1)--(c2)--(uw);
\draw[secondary,continuation] (uw)--(0,0.15);

% The other child of the last splitting vertex
\node[ordinary] (cbranch) at (1.25,0.7) {};
\draw[branch] (c2)--(cbranch);
\draw[continuation] (cbranch)--(1.25,0.15);

\node[green!45!black,left=6pt] at (0,2.0)
    {$\mathcal R_w^\ast$};

\node[left=7pt,align=right] at (c2)
    {\scriptsize last splitting\\[-2pt]
     \scriptsize vertex};

\node[red!75!black,right=5pt] at (uw)
    {$u_w^\ast$};

\node[align=center] at (0.6,-0.35)
    {\textbf{(iii)} On $\mathcal R_w^\ast$\\
     $u=u_w^\ast$};

\end{scope}

\end{tikzpicture}

\caption{\small{The three local configurations of a stopping vertex in Lemma \ref{lemma: U decompostion}. A stopping vertex lies either immediately off the principal
distinguished ray $\mathcal R^\ast$, immediately off a secondary distinguished
ray $\mathcal R_w^\ast$, or on $\mathcal R_w^\ast$ as the child $u_w^\ast$ of its last
splitting vertex. Filled circles denote splitting vertices, open circles
denote non-splitting vertices, and red squares denote stopping vertices of the full $M$-adic tree; a stopping vertex need not itself belong to $\mathcal T(\mathfrak A)$.}}
\label{fig:stopping-vertex-types}
\end{figure}
\begin{lemma} \label{lemma: U decompostion}
The collection $\mathcal U(\mathfrak A)$ admits the following decomposition, in the notation of Corollary \ref{corollary: tree structure}:
\begin{equation} \label{U decomposition} 
\mathcal U(\mathfrak A) = \mathcal U^{\ast}(\mathfrak A) \bigsqcup \bigl\{ \mathcal U_w^{\ast}(\mathfrak A) : w \in \mathcal W_{2+}\bigr\}.
\end{equation} 
More precisely, every $u\in\mathcal U(\mathfrak A)$ satisfies exactly one of the following conditions:
\begin{enumerate}[1.]
\item The vertex $u$ lies immediately off the (principal) distinguished ray $\mathcal R^*$, so that $u\in \mathcal U^{\ast}(\mathfrak A)$.
\vskip0.1in
\item There is a unique $w\in\mathcal W_{2+}$ such that $u$ lies in the maximal subtree $\mathcal T_w$, and either $u \in \mathcal W^{\ast}(\mathcal R_w^{\ast})$ or $u=u_w^{\ast}$.
\end{enumerate} 
\end{lemma} 
\begin{proof}
To establish \eqref{U decomposition}, we will first show that the left side is contained in the right.  Let $u \in \mathcal{U}(\mathfrak A)$, and let $v$ be its parent. By the
characterization \eqref{U properties} of a stopping vertex,
\begin{equation}
    \#(\mathfrak A \cap u^{\sharp}) \leq 1,
    \qquad
    \#(\mathfrak A \cap v^{\sharp}) \geq 2.
    \label{eq:stopping-vertex-counts}
\end{equation}
We first observe that $u \notin \mathcal R^*$. Indeed, Proposition \ref{TREE STRUCTURE THEOREM} \eqref{TREE STRUCTURE THEOREM-distinguished ray} states that $\mathcal R^*$ contains infinitely many splitting vertices of
$\mathcal T= \mathcal T(\mathfrak A)$. Consequently, every vertex on $\mathcal R^*$ has splitting number at
least one, whereas every $u \in \mathcal{U}(\mathfrak A)$ has splitting number zero, by \eqref{def-UA}. 
\vskip0.1in
\noindent Let $w$ be the first vertex in the lineage of $u$ that lies off $\mathcal R^*$.
Thus
\[
    w \in \mathcal{W}^*(\mathcal R^*),
    \qquad
    u \subseteq w.
\]
We distinguish two cases.
\vskip0.1in
\noindent {\em{Case 1: $\#(\mathfrak A \cap w^{\sharp}) \leq 1$.}}
We claim that $u=w$. If $u \subsetneq w$, then the parent $v$ of $u$ is also contained in $w$, and hence
\[
    \#(\mathfrak A \cap v^{\sharp})
    \leq
    \#(\mathfrak A \cap w^{\sharp})
    \leq 1,
\]
contradicting \eqref{eq:stopping-vertex-counts}. Therefore $u=w$, and
hence
\[
    u \in
    \mathcal{U}(\mathfrak A) \cap \mathcal{W}^*(\mathcal R^*)
    =
    \mathcal{U}^*(\mathfrak A).
\]
\vskip0.1in
\noindent {\em{Case 2: $\#(\mathfrak A \cap w^{\sharp}) \geq 2$.}}
In this case $w \in \mathcal T(\mathfrak A)$, and the maximal subtree $\mathcal T_w = \mathcal T_w(\mathfrak A)$ contains at least two maximal rays. Hence $w \in \mathcal{W}_{2+}$. By
Corollary \ref{corollary: tree structure}, $\mathcal T_w$ has a distinguished ray $\mathcal R_w^*$, and every
splitting vertex of $\mathcal T_w$ lies on $\mathcal R_w^*$.
\vskip0.1in
\noindent Suppose first that $u \notin \mathcal R_w^*$. Let $u_0$ be the first vertex in
the lineage of $u$ that lies off $\mathcal R_w^*$. Then
\[
    u_0 \in \mathcal{W}^*(\mathcal R_w^*),
    \qquad
    u \subseteq u_0.
\]
Since every splitting vertex of $\mathcal T_w$ lies on $\mathcal R_w^*$, the subtree
rooted at $u_0$ has splitting number zero and therefore contains at
most one maximal ray. Equivalently,
\[
    \#\bigl(\mathfrak A \cap u_0^{\sharp} \bigr) \leq 1.
\]
If $u \subsetneq u_0$, then the parent $v$ of $u$ is also contained in
$u_0$, which gives
\[
    \#\bigl(\mathfrak A \cap v^{\sharp} \bigr) \leq 1,
\]
contrary to \eqref{eq:stopping-vertex-counts}. Thus $u=u_0$, and
consequently
\[
    u \in
    \mathcal{U}(\mathfrak A) \cap \mathcal{W}^*(\mathcal R_w^*)
    \subseteq
    \mathcal{U}_w^*(\mathfrak A).
\]
\vskip0.1in
\noindent It remains to consider the scenario under case 2 when $u \in \mathcal R_w^*$. By
\eqref{eq:stopping-vertex-counts}, the canonical cell of the parent $v$ of $u$ contains at
least two points of $\mathfrak A$, whereas that of $u$ contains at most one. Since $u$
lies on $\mathcal R_w^*$, its canonical cell contains the extremal point represented by this
ray, and therefore
\[
    \#(\mathfrak A \cap u^{\sharp})=1.
\]
It follows that $v$ is a splitting vertex of $T_w$. Moreover, no
descendant of $u$ on $\mathcal R_w^*$ can be splitting, since the canonical cell of every such descendant is contained in
$u^\sharp$, which contains only one point of $\mathfrak A$. Hence $v$ is the last splitting vertex on
$\mathcal R_w^*$, and $u$ is its child lying on that ray. In other words,
\[
    u=u_w^*
    \in
    \mathcal{U}_w^*(\mathfrak A).
\]
In summary, we have proved that every $u \in \mathcal{U}(\mathfrak A)$ belongs either to
$\mathcal{U}^*(\mathfrak A)$ or to $\mathcal{U}_w^*(\mathfrak A)$ for a unique $w \in \mathcal{W}_{2+}$. This establishes that the left of \eqref{U decomposition} is contained in the right. 
\vskip0.1in
\noindent Conversely, every vertex in $\mathcal{U}^*(\mathfrak A)$, and every vertex in
\[
    \mathcal{U}(\mathfrak A) \cap \mathcal{W}^*(\mathcal R_w^*),
\]
belongs to $\mathcal{U}(\mathfrak A)$ by definitions \eqref{def-Ustar} and \eqref{def-Uwstar}. It remains only to verify
this for $u_w^*$. Its parent is the last splitting vertex on $\mathcal R_w^*$,
and hence contains at least two points of $\mathfrak A$, whereas the subtree
rooted at $u_w^*$ consists of the single ray representing the
extremal point of $\mathfrak A_w$. Therefore
\[
    \#\bigl(\mathfrak A \cap (u_w^*)^{\sharp} \bigr)=1,
    \qquad
    \#\bigl(\mathfrak A \cap ({\text{parent of }} u_w^*)^{\sharp}\bigr)\geq 2.
\]
By the characterization \eqref{U properties}, we conclude that
$u_w^* \in \mathcal{U}(\mathfrak A)$. This proves that the two sides of \eqref{U decomposition} are equal. 
\vskip0.1in
\noindent Finally, the decomposition is disjoint. Indeed, distinct first-level
vertices $w, w' \in \mathcal{W}_{2+}$ are roots of disjoint subtrees, so 
\[
    \mathcal U_w^\ast(\mathfrak A)
    \cap
    \mathcal U_{w'}^\ast(\mathfrak A)
    =\emptyset.
\]
Moreover, none of
their descendants lies immediately off $\mathcal R^*$, and hence none belongs to $\mathcal U^\ast(\mathfrak A)$.  Finally,
for each $w\in\mathcal W_{2+}$, the vertex $u_w^\ast$ lies on
$\mathcal R_w^\ast$, whereas every member of $\mathcal W^\ast(\mathcal R_w^\ast)$
lies off that ray. Thus the union in \eqref{def-Uwstar} is also disjoint. Consequently,
\[
    \mathcal{U}(\mathfrak A)
    =
    \mathcal{U}^*(\mathfrak A)
    \bigsqcup
        \bigl\{ \mathcal{U}_w^*(\mathfrak A) : w \in \mathcal W_{2+} \bigr\}.
\]
This completes the proof.
\end{proof}
\section{Proof of Lemma \ref{GAP LEMMA}} \label{section: actual proof of gap lemma} 
\subsection{Proof of Lemma \ref{GAP LEMMA}\eqref{gap lemma parta}} 
\begin{proof}
We first give the proof when $\mathfrak A=\{a_j:j\geq 1\} \subseteq [0,1)$ is strictly decreasing. Recall in this case 
\[a_0=1 \text{ and } \alpha = \lim \mathfrak A = \inf \mathfrak A. \] 
Let us first show that the first alternative in \eqref{gap lemma parta}, namely \eqref{u-sharp intersection singleton} does not hold. Indeed, if
\[
    u\cap[\alpha,a_1]=\{\alpha\},
\]
then $\alpha$ must be the right endpoint of $u$. Since
$\alpha<1$, the canonical cell $u^\sharp$ does not contain its
right endpoint, and hence
\[
    u^\sharp\cap[\alpha,a_1]= \emptyset.
\]
This contradicts the non-empty intersection condition in \eqref{U properties}.
\vskip0.1in
\noindent We now prove \eqref{u in 2 gaps} for each $u\in\mathcal U(\mathfrak A)$. From the defining property \eqref{U properties} of $\mathcal U(\mathfrak A)$, we know that
\[ u^\sharp\cap[\inf\mathfrak A,\sup\mathfrak A]\neq\emptyset, \; \text{ and } \;  \#(\mathfrak A\cap u^\sharp)\leq 1. \]
Suppose first that $\#(\mathfrak A\cap u^\sharp) = 1$. Then
\begin{equation} \label{single intersection j} 
                    \mathfrak A\cap u^\sharp=\{a_j\} \text{ for some } j \geq 1.
\end{equation} 
If the index $j$ in \eqref{single intersection j} is at least 2, then we conclude that 
\[
         u^{\sharp} \cap \bigl[ \inf \mathfrak A, \sup \mathfrak A \bigr]  = u^\sharp\subseteq (a_{j+1},a_{j-1}) \subseteq \bigl[a_{j+1}, a_{j-1} \bigr).
\]
This is because $u^\sharp$ is an interval containing $a_j$ but no other point of
$\mathfrak A$; therefore it cannot extend past either neighbouring point
$a_{j+1}$ or $a_{j-1}$ if $j \geq 2$. 
\vskip0.1in 
\noindent If $j=1$, then the same argument yields 
\[ u^{\sharp} \subseteq (a_2, 1], \text{ therefore } u^{\sharp} \cap \bigl[ \inf \mathfrak A, \sup \mathfrak A\bigr] \subseteq (a_2, a_1] \subseteq \bigl[a_{j+1}, a_{j-1} \bigr) \] under our convention $a_0 =1$ and our assumption $a_1 < a_0$.
\vskip0.1in
\noindent Suppose next that $\mathfrak A\cap u^\sharp=\emptyset$. Since
$u^\sharp$ meets the interval spanned by $\mathfrak A$, it intersects one of
the closed intervals $[a_{j+1},a_j]$. On the other hand, $u^\sharp$ is connected
and contains no point of $\mathfrak A$ by the assumption of this case, so it cannot cross either endpoint
of that interval. Hence
\[
                    u^\sharp\subseteq \bigl(a_{j+1},a_j \bigr)
                    \subseteq \bigl[a_{j+1},a_{j-1} \bigr) \text{ for a suitable $j\geq 1$}.
\]
Thus, in both cases, there exists $j\geq 1$ such that
\[
                         u^\sharp \cap \bigl[ \inf \mathfrak A, \sup \mathfrak A \bigr] \subseteq \bigl[a_{j+1},a_{j-1} \bigr),
\]
which proves \eqref{u in 2 gaps} for decreasing $\mathfrak A$.
\vskip0.1in
\noindent Now suppose that $\mathfrak A=\{a_j:j\geq1\}$ is increasing to $\alpha=\sup \mathfrak A$, with $a_0=0$. Fix $u\in\mathcal U(\mathfrak A)$. If $u\cap[a_1,\alpha]=\{\alpha\}$, then \eqref{u-sharp intersection singleton} holds.
Otherwise, $u^\sharp$ cannot contain $\alpha$ or have $\alpha$
as its right endpoint: in either case \(u^\circ\) would contain
infinitely many points of $A$, contrary to \eqref{U properties}. If $\mathfrak A\cap u^\sharp=\{a_j\}$, the absence of the neighbouring points of $A$ from $u^\sharp$ gives
\[
u^\sharp \cap[a_1,\alpha]\subseteq[a_{j-1},a_{j+1}).
\]
If $\mathfrak A\cap u^\sharp=\emptyset$, connectedness places
$u^\sharp \cap[a_1,\alpha]$ within one gap $[a_j,a_{j+1})$, which is contained in $[a_{j-1},a_{j+1})$. Thus \eqref{u in 2 gaps} holds in the increasing case as well.
\end{proof} 
\vskip0.1in
\subsection{Proof of Lemma \ref{GAP LEMMA} \eqref{gap lemma partb}} 
\begin{proof}
Suppose $x \in [\inf \mathfrak A, \sup \mathfrak A]$. We aim to show that any such $x$ lies on the right hand side of \eqref{extreme cover}. We will first consider the case when $\mathfrak A$ is decreasing, then identify the changes required in the complementary case. 
\vskip0.1in
\noindent Let us settle the endpoints first. Clearly, if $x = \alpha = \inf \mathfrak A$, then $x$ already lies on the right hand side of \eqref{extreme cover}. Similarly, if $x = a_1$, then consider the $M$-adic ray representing $a_1$ under the convention of Section \ref{section: rays as points}.  Let
$u$ be the first vertex on this ray for which $\#(\mathfrak A\cap u^\sharp)\leq 1$. Then $u\in\mathcal U(\mathfrak A)$ by \eqref{U properties}, and the choice of the representing ray ensures that $a_1\in u^\sharp$.  Once again, $x$ lies on the right side of \eqref{extreme cover}.
\vskip0.1in
\noindent Once the endpoints $\inf \mathfrak A$ and $\sup \mathfrak A$ are removed, the remaining open interval can be realized as a disjoint union of gaps, i.e., 
\[ \bigl(\inf \mathfrak A, \sup \mathfrak A \bigr) = \bigsqcup_{j=1}^{\infty} \bigl[a_{j+1}, a_{j} \bigr). \] 
Therefore, in order to complete the proof of \eqref{extreme cover} it suffices to show that 
\begin{equation} \label{extreme cover minus endpoints} 
{\tt{Gap}}_j = \bigl[a_{j+1}, a_{j} \bigr) \subseteq \bigcup \bigl\{u^{\sharp} : u \in \mathcal U(\mathfrak A)\bigr\}. 
\end{equation}  
But this is precisely the conclusion of Lemma \ref{GAP LEMMA REVISITED} when $\mathfrak A$ is decreasing.
\vskip0.1in
\noindent If $\mathfrak A$ is increasing, the same argument applies, with $\inf \mathfrak A=a_1$ and $\sup \mathfrak A=\alpha$,
and with the interval $[a_1,\alpha)$ decomposed into the half-open
gaps $[a_j,a_{j+1})$. The remaining endpoint $\alpha$ is already included separately on the
right-hand side of \eqref{extreme cover}. This proves \eqref{extreme cover} in both cases.
\end{proof}
\subsection{Proof of Lemma \ref{GAP LEMMA} \eqref{gap lemma -  partc}} 
\begin{proof} 
Suppose that $U \subseteq [0,1)$ is a closed set contained between the extreme values of $\mathfrak A$. 
%In view of  part \eqref{gap lemma partb}, we obtain 
%\[ U \subseteq U \cap \bigl[ \inf \mathfrak A, \sup \mathfrak A \bigr] = \bigcup  \bigl\{U \cap u : u \in \mathcal U(\mathfrak A) \bigr\} \cup \{ \alpha \}. \] 
By the general decomposition \eqref{tree decomp in w} of trees, 
\begin{equation}
\mathcal T(U)  
%&= \mathcal T(U; M) \subseteq \bigcup  \bigl\{\mathcal T(U \cap u; M) : u \in \mathcal U(\mathfrak A) \bigr\} \cup \{ \mathcal R^{\ast} \} \nonumber \\ 
%&\subseteq \bigcup  \bigl\{\mathcal T_u(U; M) : u \in \mathcal U(\mathfrak A) \bigr\} \cup \{ \mathcal R^{\ast} \}  \nonumber \\
\subseteq \bigcup  \Bigl\{\mathcal T_w(U; M) : w \in \mathcal W\bigl(\mathcal T(U), \mathcal R^{\ast} \bigr)\Bigr\} \cup \{ \mathcal R^{\ast} \}, \label{classifying by u}
\end{equation}  
where $\mathcal W(\cdot, \cdot)$ is defined as in \eqref{what is W(T)}. 
\vskip0.1in
\noindent Fix $w\in\mathcal W\bigl(\mathcal T(U),R^\ast\bigr)$.
Then it follows from \eqref{what is W(T)} that $w \in \mathcal T(U)$, which in turn implies
\begin{equation} \label{intersection of w}
 w^{\sharp} \cap \bigl[ \inf U, \sup U \bigr] \neq \emptyset, \text{ and hence } w^{\sharp} \cap \bigl[ \inf \mathfrak A, \sup \mathfrak A\bigr] \neq \emptyset.  
\end{equation} 
Since the parent of $w$ lies on the principal distinguished ray
$\mathcal R^\ast$ of $\mathcal T(\mathfrak A)$, it contains infinitely many points of $\mathfrak A$; in particular,
\begin{equation} \#\bigl( \mathfrak A \cap (\text{parent of w})^{\sharp} \bigr) \geq 2. \label{parent failure} 
\end{equation}  Consequently,
exactly one of the following two cases occurs.
\vskip0.1in
\noindent {\em{Case 1: }}  If $\#(\mathfrak A\cap w^\sharp)\leq 1$, then \eqref{intersection of w} and \eqref{parent failure} combined with the characterization \eqref{U properties} gives $w\in\mathcal U(\mathfrak A)$. The hypothesis \eqref{lac order of lowest level} therefore yields
\begin{equation} \label{Case1-conclusion} 
        \mathcal T_w(U;M)\in {\tt{TreeLac}}(J,R;M)
        \subseteq {\tt{TreeLac}}(J+1,R;M).
\end{equation} 
\vskip0.1in
\noindent {\em{Case 2: }} If $\#(\mathfrak A\cap w^\sharp)\geq 2$, then $w$ belongs to the collection
$\mathcal W_{2+} = \mathcal W_{2+}(\mathcal T(\mathfrak A), \mathcal R^{\ast})$, defined in \eqref{def-W2+}. From Corollary \ref{corollary: tree structure}, the maximal subtree
of $\mathcal T(\mathfrak A)$ rooted at $w$ possesses the secondary distinguished ray $\mathcal R_w^\ast$. The hypothesis \eqref{lac order of lowest level} combined with Lemma \ref{lemma: U decompostion} then implies 
\[ \mathcal T_u (U;M)\in {\tt{TreeLac}}(J,R;M) \text{ for all } u \in \mathcal U_w^{\ast}(\mathfrak A). \]
We claim that every maximal subtree of $\mathcal T_w(U;M)$ rooted
immediately off $\mathcal R_w^\ast$ belongs to ${\tt{TreeLac}}(J,R;M)$. Indeed, let
\[
 z\in\mathcal W\bigl(\mathcal T_w(U;M), \mathcal R_w^\ast\bigr),
\]
i.e., $z$ is a vertex of $\mathcal T_w(U)$ which is not on $\mathcal R_w^{\ast}$, but whose parent is. In particular, this means
\begin{equation}  \label{intersection of z} 
z^{\sharp} \cap U \neq \emptyset, \text{ and hence } z^{\sharp} \cap \bigl[\inf \mathfrak A, \sup \mathfrak A \bigr] \neq \emptyset, 
\end{equation}  
by the same argument that led to \eqref{intersection of w}. 
\vskip0.1in 
\noindent According to Corollary \ref{corollary: tree structure}, every splitting vertex of $\mathcal T_w(\mathfrak A;M)$ lies on $\mathcal R_w^\ast$; therefore one has
\[
                         \#(\mathfrak A\cap z^\sharp)\leq 1.
\]
If the canonical cell of the parent of $z$ contains at least two points of $\mathfrak A$, then 
$z\in\mathcal U(\mathfrak A)$ by \eqref{intersection of z} and \eqref{U properties}, and hence
\[
             \mathcal T_z(U;M)\in {\tt{TreeLac}}(J,R;M)
\]
by the hypothesis \eqref{lac order of lowest level}.
\vskip0.1in
\noindent Otherwise, the parent of $z$ contains at most one point of $\mathfrak A$.
It then lies strictly beyond the last splitting vertex of
$\mathcal T_w(\mathfrak A;M)$ along $\mathcal R_w^\ast$. Consequently,
\[z\subseteq u_w^\ast, \]
where $u_w^\ast\in\mathcal U_w^\ast(\mathfrak A)$ is the child of the last
splitting vertex lying on $\mathcal R_w^\ast$. Therefore
\[
        \mathcal T_z(U;M)\subseteq\mathcal T_{u_w^\ast}(U;M).
\]
Since $u_w^\ast \in\mathcal U(A)$, the hypothesis \eqref{lac order of lowest level} and
inheritance of tree lacunarity by rooted subtrees give
\[
             \mathcal T_z(U;M)\in {\tt{TreeLac}}(J,R;M).
\]
Thus every subtree rooted immediately off $\mathcal R_w^\ast$ lies in
${\tt{TreeLac}}(J,R;M)$. Lemma \ref{lemma: NRM}, applied with distinguished ray $\mathcal R_w^\ast$, now gives
\begin{equation} \label{Case2-conclusion} 
             \mathcal T_w(U;M)\in {\tt{TreeLac}}(J+1,R;M).
\end{equation} 
Combining the conclusions \eqref{Case1-conclusion} and \eqref{Case2-conclusion} of Cases~1 and~2 respectively, we deduce that
\[
 \mathcal T_w(U;M)\in {\tt{TreeLac}}(J+1,R;M)
 \quad
 \text{for every }
 w\in\mathcal W\bigl(\mathcal T(U),\mathcal R^\ast\bigr).
\]
Applying Lemma \ref{lemma: NRM} once more, now with the principal distinguished ray $\mathcal R^\ast$, gives
\[\mathcal T(U;M)\in {\tt{TreeLac}}(J+2,R;M), \]
which proves \eqref{order J+2}. Finally, Corollary \ref{corollary: admissible treelac} yields
\[ \text{split}\bigl(\mathcal T(U;M)\bigr)
                  \leq R^{10(J+2)}, \]
which is \eqref{Split TU}.
\end{proof}

	\end{document}